\documentclass[11pt]{book}
\usepackage[T2A]{fontenc}
\usepackage[cp1251]{inputenc}
\usepackage[english, ukrainian, russian]{babel}
\usepackage{amsmath,amssymb,amsthm}
\usepackage{mathrsfs}
\usepackage{titlesec}
\usepackage{indentfirst}
\usepackage{graphicx}

\newtheorem{theorem}{\indent Теорема}[chapter]

\newtheorem{lemma}{\indent Лемма}[chapter]

\newtheorem{corollary}{\indent Следствие}[chapter]

\newtheorem{proposition}{\indent Предложение}[chapter]

\theoremstyle{definition}

\newtheorem{example}{\indent Пример}[chapter]

\newtheorem{definition}{\indent Определение}[chapter]

\newtheorem{remark}{\indent Замечание}[chapter]

\renewcommand{\chaptermark}[1]{\markboth{#1}{#1}}

\allowdisplaybreaks[3]

\titlelabel{\thetitle.\quad}

\begin{document}



\thispagestyle{empty}

\bf

\phantom{m}

\vspace{2cm}

\LARGE

\begin{flushright}

В. А. МИХАЙЛЕЦ \\ А. А. МУРАЧ

\end{flushright}

\vspace{0.5cm}

\huge

\begin{center}

ПРОСТРАНСТВА \medskip \\ ХЕРМАНДЕРА, \medskip \\  ИНТЕРПОЛЯЦИЯ И \medskip \\
ЭЛЛИПТИЧЕСКИЕ ЗАДАЧИ

\end{center}

\vspace{1cm}

\Large

\begin{center}

\textit{с предисловием академика \\ Ю. М. Березанского}

\end{center}

\vspace{2.5cm}

\Large

\begin{center}

Киев --- 2010

\end{center}


\normalsize

\rm

\newpage

\thispagestyle{empty}

\noindent УДК 517.95(2), 517.98(02), 517.956.223

\smallskip

Михайлец~В.~А., Мурач~А.~А. Пространства Хермандера, интерполяция и эллиптические
задачи.~--- Киев, 2010.

\bigskip

Монография дает первое систематическое изложение теории эллиптических (скалярных и
матричных) операторов и эллиптических краевых задач в шкалах гильбертовых
пространств Хермандера функций/распределений произвольной положительной или
отрицательной гладкости. Ее отличает использование метода интерполяции с
функциональным параметром абстрактных и соболевских гильбертовых пространств. Часть
результатов являются новыми и для соболевских шкал.

Книга предназначена для научных работников, преподавателей университетов и
аспирантов соответствующих специальностей.

\vspace{2cm}

Рецензенты:

\medskip

академик НАН Украины,

доктор физ.-мат. наук, профессор \phantom{мммм} Ю.~М.~Березанский

\medskip

доктор физ.-мат. наук, профессор \phantom{мммм}  С.~Д.~Ивасишен

\vspace{5cm}

 © В. А. Михайлец, А. А. Мурач, 2010

\newpage

\thispagestyle{empty}

\phantom{м}

\vspace{2cm}

\bf

\LARGE

\begin{flushright}

V. A. MIKHAILETS \\ A. A. MURACH

\end{flushright}

\vspace{0.5cm}

\huge

\begin{center}

H\"ORMANDER SPACES, \medskip \\  INTERPOLATION, \medskip \\
AND ELLIPTIC PROBLEMS

\end{center}

\vspace{1cm}

\Large

\begin{center}

\textit{with preface by academician \\ Yu. M. Berezansky}

\end{center}

\vspace{3.5cm}

\Large

\begin{center}

Kiev --- 2010

\end{center}

\normalsize

\rm

\newpage

\thispagestyle{empty}


\smallskip

Mikhailets~V.~A., Murach~A.~A. H\"ormander Spaces, Interpolation, and Elliptic
Problems.~--- Kiev, 2010.

\bigskip

The research monograph gives the first systematic exposition of the elliptic (scalar
and matrix) operators theory and elliptic boundary-value problems in the scales of
Hilbert spaces of H\"ormander of the functions/distributions of arbitrary positive
or negative smoothness. The book is based on the method of interpolation with a
functional parameter for the abstract and Sobolev unitary spaces. Some results are
also new for the Sobolev scales.

The monograph is intended for the researches, professors and PhD students.

\vspace{2cm}

Reviewers:

\medskip

academician of NAS of Ukraine,

Doctor Hab., Professor \phantom{mmmmmmm} Yu.~M.~Berezansky

\medskip

Doctor Hab., Professor \phantom{mmmmmmm}  S.~D.~Ivasyshen

\vspace{6cm}

 © V. A. Mikhailets, A. A. Murach, 2010

\normalsize



\tableofcontents

\chaptermark{}

\chapter*{\textbf{Предисловие}}




\addcontentsline{toc}{section}{\textbf{Предисловие Ю.~М.~Березанского}}


Хорошо известны фундаментальные применения пространств Соболева $W^{m}_{2}(G)$ к
изучению многомерных дифференциальных уравнений, в частности, эллиптического типа.
Без теории таких пространств изучение эллиптических задач по существу невозможно.
Вместе с тем уже примерно 40 лет тому назад была разработана теория более общих, чем
соболевские, пространств Хермандера. Существует значительное число работ,
посвященное применениям этих пространств к дифференциальным уравнениям.

Однако применения пространств Хермандера к краевым задачам для эллиптических
уравнений до последнего времени носили эпизодический характер. Основная часть этой
книги~--- достаточно систематическое применение пространств Хермандера именно к этим
задачам. Для этого авторам пришлось ввести и исследовать пространства Хермандера
„промежуточного типа”, функции из которых имеют промежуточную гладкость между
гладкостью функций из $W^{m}_{2}(G)$ и $W^{m+1}_{2}(G)$, где $m$~--- целое число.
Эти пространства вводятся и исследуются в книге. В~качестве $G$ может служить
область $n$-мерного пространства или компактное многообразие размерности~$n$.
Детальному введению и исследованию таких пространств посвящены первые две главы
книги.

В главах 3, 4 рассматриваются эллиптические уравнения и краевые задачи для них, как
однородные так и неоднородные. Для них в пространствах Хермандера устанавливается
ряд существенных результатов, подобных случаю соболевских пространств. Можно
сказать, что авторам удалось перенести классическую „соболевскую” теорию краевых
задач на случай пространств Хермандера. Отметим, что некоторые задачи,
сформулированные независимо от понятия пространства Хермандера, удалось решить с
помощью этих пространств.

Наконец, глава 5 книги посвящена перенесению упомянутых результатов на эллиптические
системы уравнений.

Книга мне представляется интересной и полезной. Я уверен, что ее следует перевести
на английский язык, тогда ее результаты были бы доступны значительно большему кругу
математиков. При переводе следует привести доказательства ряда вспомогательных
фактов, необходимых для изложения и принадлежащих другим авторам. Это сделало бы ее
доступной для более широкого круга читателей.

\bigskip

\large\emph{Академик НАН Украины Ю.~М.~Березанский}\normalsize

\bigskip\bigskip\bigskip

\addcontentsline{toc}{section}{Благодарности}

\noindent\Huge\textbf{Благодарности}

\normalsize

\bigskip\bigskip

Авторы особо признательны Ю.~М.~Березанскому за ценные советы и глубокое влияние,
во-многом определившее их научные интересы.

Мы искренне благодарны М.~С.~Аграновичу, Б.~П.~Панеяху, \fbox{И.~В.~Скрыпнику} \,и
\,\fbox{С.~Д.~Эйдельману} \,за общение, которое стимулировало нас.

Очень существенными для нас были: поддержка М.~Л.~Горбачука и А.~М.~Самойленко,
внимание Б.~Боярского, доброжелательное участие В.~П.~Бурского и С.~Д.~Ивасишена.

Мы также признательны всем коллегам, которые проявили живой интерес к новой теории и
ее применениям.


\newpage

\chaptermark{}

\chapter*{\textbf{Введение}}

\chaptermark{\emph Введение}



\addcontentsline{toc}{section}{\textbf{Введение}}

\markboth{\emph Введение}{Введение}

В теории уравнений с частными производными центральное место занимают вопросы о
существовании, единственности и регулярности решений. При этом свойства регулярности
исследуемых решений формулируются в терминах их принадлежности к эталонным классам
функциональных пространств. Чем тоньше градуирована шкала пространств, тем точнее и
содержательнее такие результаты.

В отличие от обыкновенных дифференциальных уравнений с гладкими коэффициентами эти
вопросы являются достаточно сложными. Так, известны линейные уравнения с частными
производными с гладкими коэффициентами и правыми частями, которые не имеют решений в
окрестности заданной точки даже в классе распределений [\ref{Lewy57};
\,\ref{Hermander65}, с.~210]. Далее, некоторые однородные уравнения, в частности,
эллиптического типа, с гладкими (но не аналитическими) коэффициентами имеют
нетривиальные решения с компактными носителями [\ref{Plis54}; \,\ref{Hermander86},
с.~257]. Поэтому нетривиальное нуль-пространство такого уравнения нельзя устранить
какими-либо однородными краевыми условиями, т.~е. оператор любой краевой задачи не
является инъективным. Наконец, вопрос о регулярности решения также весьма непрост.
Например, уже для оператора Лапласа известно [\ref{GilbargTrudinger89}, с.~73], что
$$
\triangle u=f\in\,C(\Omega)\nRightarrow\,u\in C^{\,2}(\Omega)
$$
для произвольной евклидовой области $\Omega$.

Наиболее полно указанные вопросы исследованы для линейных эллиптических уравнений,
систем и краевых задач. Это было сделано в 50 -- 60-х годах XX столетия в работах
С.~Агмона, А.~Дуглиса и Л.~Ниренберга [\ref{AgmonDouglisNirenberg62},
\ref{AgmonDouglisNirenberg64}, \ref{DouglisNirenberg55}], М.~C.~Аграновича и
А.~С.~Дынина [\ref{AgranovichDynin62}], Ю.~М.~Березанского, С.~Г.~Крейна и
Я.~А.~Ройтберга [\ref{BerezanskyKreinRoitberg63}, \ref{Berezansky65},
\ref{Roitberg64}, \ref{Roitberg65}, \ref{Roitberg96}], Ф.~Браудера [\ref{Browder56},
\ref{Browder59}], Л.~Р.~Волевича [\ref{Volevich63}, \ref{Volevich65}], Ж.-Л.~Лионса
и Э.~Мадженеса [\ref{LionsMagenes71}, \ref{Magenes66}], Л.~Н.~Слободецкого
[\ref{Slobodetsky58c}, \ref{Slobodetsky60}], В.~А.~Солонникова [\ref{Solonnikov64},
\ref{Solonnikov66}, \ref{Solonnikov67}], Л.~Хермандера [\ref{Hermander65}],
М.~Шехтера [\ref{Schechter60a}, \ref{Schechter61}, \ref{Schechter63}] и других. При
этом эллиптические уравнения были исследованы в классических шкалах пространств
Гельдера (нецелого порядка) и пространств Соболева (как положительного, так и
отрицательного порядков).

Фундаментальный результат теории эллиптических уравнений состоит в том, что они
порождают ограниченные нетеровы операторы (т.~е. операторы с конечным индексом) в
соответствующих парах пространств. Пусть, например, на замкнутом гладком
многообразии $\Gamma$ задано линейное эллиптическое дифференциальное уравнение
$Au=f$ порядка $m$. Тогда оператор
$$
A:\,H^{s+m}(\Gamma)\rightarrow H^{s}(\Gamma),\quad s\in\mathbb{R},
$$
ограничен и нетеров. При этом конечномерные пространства решений однородных
уравнений $Au=0$ и $A^{+}v=0$ лежат в $C^{\infty}(\Gamma)$. Здесь $A^{+}$~---
формально сопряженный оператор к $A$, а $H^{s+r}(\Gamma)$ и $H^{s}(\Gamma)$~---
гильбертовы пространства Соболева на $\Gamma$ соответственно порядков $s+r$ и $s$.
Отсюда следует важное свойство регулярности решения $u$ в соболевской шкале:
$$
(f\in H^{s}(\Gamma)\;\;\mbox{для некоторого}\;\;s\in\mathbb{R})\;\Rightarrow\; u\in
H^{s+m}(\Gamma).\eqno(1)
$$

Если многообразие имеет край, то нетеров оператор порождается эллиптической краевой
задачей для уравнения $Au=f$, например, краевой задачей Дирихле.

Некоторые из этих теорем были распространены Х.~Трибелем [\ref{Triebel80},
\ref{Triebel86}] и одним из авторов [\ref{94Dop12}, \ref{94UMJ12}] на более тонкие
шкалы функциональных пространств~--- Никольского--Бесова, Зигмунда и
Лизоркина--Трибеля.

Указанные результаты нашли различные приложения в теории дифференциальных уравнений,
математической физике, спектральной теории дифференциальных операторов (см.
монографии Ю.~М.~Березанского [\ref{Berezansky65}], Ю.~М.~Березанского, Г.~Ф.~Уса,
З.~Г.~Шефтеля [\ref{BerezanskyUsSheftel90}], О.~А.~Ладыженской и Н.~Н.~Уральцевой
[\ref{LadyszenskayaUraltseva64}], Ж.-Л.~Лионса [\ref{Lions72a}, \ref{Lions72b}],
Ж.-Л.~Лионса и Э.~Мадженеса [\ref{LionsMagenes71}], Я.~А.~Ройтберга
[\ref{Roitberg96}, \ref{Roitberg99}], И.~В.~Скрыпника [\ref{Skrypnik73}], Х.Трибеля
[\ref{Triebel80}], обзоры М.~С.~Аграновича [\ref{Agranovich65}, \ref{Agranovich90},
\ref{Agranovich97}] и приведенную там литературу).

С точки зрения приложений, особенно к спектральной теории, наиболее важен случай
гильбертовых пространств. Однако единственной шкалой гильбертовых пространств, в
которой были систематически исследованы свойства эллиптических операторов, до
последнего времени оставалась шкала пространств Соболева. Как выяснилось, для ряда
важных задач эта шкала не является достаточно тонкой.

Приведем два характерных примера. Первый касается исследования гладкости решения
эллиптического уравнения $Au=f$ на многообразии $\Gamma$. Согласно теореме вложения
Соболева
$$
H^{\sigma}(\Gamma)\subset C^{r}(\Gamma)\;\;\Leftrightarrow\;\;\sigma>r+n/2, \eqno(2)
$$
где целое $r\geq0$, а $n:=\dim\Gamma$. Этот факт вместе со свойством (1) позволяет
исследовать классическую гладкость решения $u$. Так, если $f\in H^{s}(\Gamma)$ для
некоторого $s>r-m+n/2$, то $u\in H^{s+m}(\Gamma)\subset C^{r}(\Gamma)$. Однако при
$s=r-m+n/2$ это не так, т.~е. в терминах соболевской шкалы нельзя выразить
неулучшаемые достаточные условия принадлежности решения $u$ к классу
$C^{r}(\Gamma)$. Аналогичная ситуация имеется и в теории эллиптических краевых
задач.

Второй показательный пример относится к спектральной теории. Предположим, что
дифференциальный оператор $A$ порядка $m>0$ эллиптичен на $\Gamma$ и самосопряжен в
пространстве $L_{2}(\Gamma)$. Рассмотрим разложение функции $f\in\nobreak
L_{2}(\Gamma)$ в ряд
$$
f=\sum_{j=1}^{\infty}\;c_{j}(f)\,h_{j} \eqno(3)
$$
по полной ортонормированной системе $(h_{j})_{j=1}^{\infty}$ собственных функций
оператора~$A$. Здесь $c_{j}(f)$~--- коэффициенты Фурье функции $f$ по функциям
$h_{j}$, которые занумерованы так, что модули соответствующих собственных чисел
образуют возрастающую (нестрого) последовательность. Согласно теореме
Меньшова--Радемахера (справедливой для общих ортогональных рядов) ряд (3) сходится
почти всюду на $\Gamma$, если
$$
\sum_{j=1}^{\infty}\,|c_{j}(f)|^{2}\,\log^{2}(j+1)<\infty. \eqno(4)
$$
В терминах принадлежности функции $f$ пространствам Соболева это условие нельзя
переформулировать эквивалентным образом, поскольку для произвольного $s>0$
$$
\|f\|_{H^{s}(\Gamma)}^{2}\,\asymp\,\sum_{j=1}^{\infty}\,|c_{j}(f)|^{2}\,j^{2s}.
$$
Можем только утверждать, что условие „$f\in H^{s}(\Gamma)$ для некоторого $s>0$”
влечет сходимость почти всюду на $\Gamma$ ряда $(3)$. Оно не выражает адекватным
образом условие $(4)$ теоремы Меньшова--Радемахера.

Широкое и содержательное обобщение пространств Соболева в категории гильбертовых
пространств было предложено Л.~Хермандером [\ref{Hermander65}, с.~54] в 1963 году.
Введенные им пространства параметризуются не числовым набором, а (достаточно общей)
весовой функцией. В частности, Хермандер рассмотрел гильбертовы пространства
\begin{gather}\label{1}
B_{2,\mu}(\mathbb{R}^{n}):=
\bigl\{u:\,\mu\,\mathcal{F}u\in L_{2}(\mathbb{R}^{n})\bigr\}, \tag{5}\\
\|u\|_{B_{2,\mu}(\mathbb{R}^{n})}:= \|\mu\,\mathcal{F}u\|_{L_{2}(\mathbb{R}^{n})}.
\notag
\end{gather}
Здесь $\mathcal{F}u$~--- преобразование Фурье медленно растущего распределения $u$,
заданного в $\mathbb{R}^{n}$, $\mu$~--- весовая функция $n$ переменных.

В~случае, когда
$$
\mu(\xi)=\langle\xi\rangle^{s},\quad
\langle\xi\rangle:=(1+|\xi|^{2})^{1/2},\quad\xi\in\mathbb{R}^{n},\quad
s\in\mathbb{R},
$$
мы получаем пространство Соболева $B_{2,\mu}(\mathbb{R}^{n})=H^{s}(\mathbb{R}^{n})$.

В 1965 году пространства (5) независимо ввели и исследовали также Л.~Р.~Волевич и
Б.~П.~Панеях [\ref{VolevichPaneah65}].

Пространства Хермандера занимают центральное место среди пространств обобщенной
(функциональной) гладкости. Они являются предметом глубоких исследований,
значительная часть которых выполнена в последние десятилетия: см. статью
Л.~Р.~Волевича и Б.~П.~Панеяха [\ref{VolevichPaneah65}], обзор Г.~А.~Калябина и
П.~И.~Лизоркина [\ref{KalyabinLizorkin87}], монографию Х.~Трибеля [\ref{Triebel06}],
недавние работы С.~Д.~Моуры [\ref{Moura01}], Б.~Опика и В.~Требелза
[\ref{OpicTrebels00}], В.~Фаркаса и Х.-Г.~Леопольда [\ref{FarkasLeopold06}],
В.~Фаркаса, Н.~Якоба и Р.~Л.~Шилинга [\ref{FarkasJacobScilling01b}], Д.~Д.~Хароске и
С.~Д.~Моуры [\ref{HaroskeMoura04}] и приведенную там литературу. Разные классы
пространств обобщенной гладкости естественно возникают в теоремах вложения
функциональных пространств, в теории интерполяции функциональных пространств, в
спектральной теории дифференциальных операторов, в теории случайных процессов.

Еще в 1963 году пространства (5) и более общие банаховы пространства
$B_{p,\mu}(\mathbb{R}^{n})$, где $1\leq p\leq\infty$,  были использованы Хермандером
[\ref{Hermander65}] для исследования свойств решений уравнений с частными
производными с постоянными коэффициентами, а также некоторых классов уравнений с
переменными коэффициентами, заданных в евклидовых областях. Однако, в отличие от
соболевских, пространства Хермандера не нашли широкого применения в теории общих
эллиптических уравнений на многообразиях и в теории эллиптических краевых задач. Это
связано с продолжительным отсутствием корректного определения таких пространств на
гладких многообразиях (оно не должно зависеть от выбора локальных карт, покрывающих
многообразие), а также отсутствием аналитического аппарата, пригодного для работы с
такими пространствами.

В случае пространств Соболева такой аппарат есть --- это интерполяция пространств.
Так, произвольное пространство Соболева дробного порядка можно получить путем
интерполяции некоторой пары соболевских пространств целых порядков. Последнее
существенно облегчает как исследование этих пространств, так и доказательство
различных теорем теории эллиптических уравнений, поскольку при интерполяции
сохраняются свойства ограниченности линейных операторов и их нетеровости (при
неизменном дефекте)

Поэтому представляется разумным выделить среди гильбертовых пространств Хермандера
те, которые получаются посредством интерполяции (теперь уже с функциональным
параметром) пар гильбертовых пространств Соболева. Для этого мы вводим класс
изотропных пространств
$$
H^{s,\varphi}(\mathbb{R}^{n}):=B_{2,\mu}(\mathbb{R}^{n})\quad\mbox{для}\quad
\mu(\xi)={\langle\xi\rangle^{s}\varphi(\langle\xi\rangle)}. \eqno(6)
$$
Здесь числовой параметр $s\in\mathbb{R}$, а функциональный параметр $\varphi$
медленно меняется на бесконечности по Карамата [\ref{Seneta85},
\ref{BinghamGoldieTeugels89}] (можно считать, что положительная функция $\varphi$
постоянна вне окрестности бесконечности). Например, в качестве $\varphi$ допускаются
логарифмическая функция, ее итерации, любые их степени, а также произведения таких
функций.

Класс пространств $(6)$ содержит соболевскую гильбертову шкалу
$\{H^{s}\}\equiv\{H^{s,1}\}$, привязан к ней числовым параметром, но градуирован
тоньше, чем соболевская шкала. Так,
$$
H^{s+\varepsilon}(\mathbb{R}^{n})\subset H^{s,\varphi}(\mathbb{R}^{n})\subset
H^{s-\varepsilon}(\mathbb{R}^{n})\quad\mbox{для любого}\quad\varepsilon>0.
$$
Поэтому числовой параметр $s$ определяет основную (степенную) гладкость, а
функциональный параметр $\varphi$ задает дополнительную (субстепенную) гладкость в
классе пространств (6). В частности, если $\varphi(t)\rightarrow\infty$ (или
$\varphi(t)\rightarrow0$) при $t\rightarrow\infty$, то $\varphi$ задает
дополнительную положительную (или отрицательную) гладкость. Т.~е. параметр $\varphi$
уточняет основную гладкость $s$. Поэтому класс пространств (6) естественно именовать
уточненной шкалой (по отношению к исходной соболевской шкале).

Уточненная шкала имеет важное свойство: каждое пространство
$H^{s,\varphi}(\mathbb{R}^{n})$ есть результат интерполяции с подходящим
функциональным параметром пары соболевских пространств
$H^{s-\varepsilon}(\mathbb{R}^{n})$ и $H^{s+\delta}(\mathbb{R}^{n})$, где
$\varepsilon,\delta>0$. Этот параметр является правильно меняющейся функцией на
бесконечности (по Карамата) порядка $\theta\in(0,\,1)$, где
$\theta:=\varepsilon/(\varepsilon+\delta)$. Более того, как оказалось, класс
пространств $(6)$ замкнут относительно такой интерполяции.

Таким образом, пространства Хермандера $H^{s,\varphi}(\mathbb{R}^{n})$ обладают
интерполяционным свойством по отношению к соболевской гильбертовой шкале. Это
означает, что всякий линейный оператор, ограниченный в каждом из пространств
$H^{s-\varepsilon}(\mathbb{R}^{n})$ и $H^{s+\delta}(\mathbb{R}^{n})$, является
ограниченным и в $H^{s,\varphi}(\mathbb{R}^{n})$. Интерполяционное свойство играет
тут решающую роль~--- оно дает возможность установить такие важные свойства
уточненной соболевской шкалы, которые позволяют применять ее в теории эллиптических
уравнений. Так, при помощи интерполяции доказывается, что пространства
$H^{s,\varphi}(\mathbb{R}^{n})$, как и соболевские пространства, инвариантны
относительно диффеоморфных преобразований $\mathbb{R}^{n}$. Это позволяет корректно
определить пространства $H^{s,\varphi}(\Gamma)$ на гладком замкнутом многообразии
$\Gamma$, поскольку запас распределений и топология в них не зависят от выбора
набора локальных карт, покрывающих $\Gamma$. Такие функциональные пространства
необходимы в теории эллиптических операторов на многообразиях, в теории
эллиптических краевых задач и неявным образом присутствуют в ряде проблем анализа.

Остановимся на некоторых результатах, демонстрирующих преимущества уточненной шкалы
по сравнению с соболевской. Они относятся к примерам, рассмотренным выше. Пусть, как
и прежде, $A$~--- эллиптический дифференциальный оператор на $\Gamma$ порядка $m$.
Тогда он задает ограниченные нетеровы операторы
$$
A:\,H^{s+m,\varphi}(\Gamma)\rightarrow H^{s,\varphi}(\Gamma)\quad\mbox{для
всех}\quad s\in\mathbb{R},\;\;\varphi\in\mathcal{M}.
$$
Здесь $\mathcal{M}$~--- класс медленно меняющихся функциональных параметров
$\varphi$, использованных в (6). Отметим, что оператор $A$ оставляет инвариантным
функциональный параметр $\varphi$, уточняющий основную гладкость $s$. Кроме того,
выполняется следующее свойство регулярности решения уравнения $Au=f$:
$$
(f\in H^{s,\varphi}(\Gamma)\;\;\mbox{для
некоторых}\;\;s\in\mathbb{R},\;\varphi\in\mathcal{M})\;\Rightarrow\; u\in
H^{s+m,\varphi}(\Gamma).
$$

Из него вытекает, что для уточненной шкалы справедливо следующее усиление теоремы
вложения Соболева: для целого числа $r\geq0$ и функции $\varphi\in\mathcal{M}$
вложение $H^{r+n/2,\varphi}(\Gamma)\subset C^{r}(\Gamma)$ равносильно тому, что
$$
\int\limits_{1}^{\infty}\frac{dt}{t\,\varphi^{2}(t)}<\infty. \eqno(7)
$$
Поэтому, если $f\in H^{r-m+n/2,\varphi}(\Gamma)$ для некоторого параметра
$\varphi\in\nobreak\mathcal{M}$, удовлетворяющего (7), то решение $u\in
C^{r}(\Gamma)$.

Подобные результаты справедливы также для эллиптических систем и для регулярных
эллиптических краевых задач.

Перейдем теперь к анализу сходимости спектрального разложения (3) (оператор $A$
порядка $m>0$ предполагается неограниченным самосопряженным в пространстве
$L_{2}(\Gamma)$). Условие (4) о сходимости ряда (3) почти всюду на $\Gamma$
эквивалентно включению
$$
f\in H^{0,\varphi}(\Gamma),\quad\mbox{где}\quad\varphi(t):=\max\{1,\log t\}.
$$
Оно существенно шире условия „$f\in H^{\varepsilon}(\Gamma)$ для некоторого
$\varepsilon>\nobreak0$”. Сходным образом можно выразить также условия
\emph{безусловной} сходимости ряда (3) почти всюду или в банаховом пространстве
$C^{r}(\Gamma)$ для целых $r\geq0$.

Эти и другие результаты показывают сколь полезной и удобной может быть уточненная
шкала. Ее можно использовать также и в других  областях современного анализа (см.,
например, работы П.~Матэ и У.~Тотенхана [\ref{MatheTautenhahn06}], М.~Хегланда
[\ref{Hegland95}, \ref{Hegland10}]).

Предлагаемая вниманию читателей монография дает первое систематическое изложение
теории эллиптических (скалярных и матричных) операторов и регулярных эллиптических
краевых задач в уточненных шкалах пространств, а также связанных с ней вопросов
интерполяции с функциональным параметром гильбертовых (абстрактных и соболевских)
пространств. Эта теория построена авторами в недавних работах [\ref{05UMJ5}~--
\ref{09Dop3}, \ref{07UMJ6}~-- \ref{07Dop6}, \ref{08MFAT1}~--
\ref{09OperatorTheory191}, \ref{08MFAT2}~-- \ref{09MFAT2}]. Содержание и структура
монографии достаточно полно отражены в ее оглавлении.


\newpage

\chapter{\textbf{Интерполяция и пространства Хермандера}}\label{ch2a}

\chaptermark{\emph Гл. \ref{ch2a}. Интерполяция и пространства Хермандера}



\section[Интерполяция с функциональным параметром]
{Интерполяция \\ с функциональным параметром}\label{sec2.1}

\markright{\emph \ref{sec2.1}. Интерполяция с функциональным параметром}

В этом параграфе мы рассмотрим интерполяцию с функциональным параметром пар
гильбертовых пространств. Она является естественным обобщением классического
интерполяционного метода Ж.-Л.~Лионса и С.~Г.~Крейна [\ref{LionsMagenes71}, с.~41;
\ref{FunctionalAnalysis72}, с.~253] на случай, когда в качестве параметра
интерполяции берется не число, а довольно общая функция. Интерполяция с
функциональным параметром будет одним из основных методов наших доказательств.

\subsection{Определение интерполяции}\label{sec2.1.1}

Приведем определение интерполяции с функциональным параметром пар гильбертовых
пространств и исследуем ряд ее свойств, необходимых в дальнейшем. Для наших целей
достаточно ограничиться случаем сепарабельных пространств.


\begin{definition}\label{def2.1} Упорядоченную пару $[X_{0},X_{1}]$
комплексных гильбертовых пространств $X_{0}$ и $X_{1}$ назовем \textit{допустимой},
если пространства $X_{0}$, $X_{1}$ сепарабельны и справедливо непрерывное плотное
вложение $X_{1}\hookrightarrow X_{0}$.
\end{definition}


Пусть задана допустимая пара $X=[X_{\,0},X_{1}]$ гильбертовых пространств. Как
известно [\ref{LionsMagenes71}, c.~22; \ref{FunctionalAnalysis72}, c.~251], для $X$
существует такой изометрический изоморфизм $J:X_{1}\leftrightarrow X_{\,0}$, что $J$
является самосопряженным положительно определенным оператором в пространстве
$X_{\,0}$ с областью определения $X_{1}$. Оператор $J$ называется
\textit{порождающим} для пары $X$. Этот оператор определяется парой $X$ однозначно.
В самом деле, пусть оператор $J_{1}$ также порождающий для пары $X$. Тогда операторы
$J$ и $J_{1}$ метрически равны:
$\|Ju\|_{X_{\,0}}=\|u\|_{X_{1}}=\|J_{1}u\|_{X_{\,0}}$ для любого $u\in X_{1}$. Кроме
того, эти операторы положительно определены. Отсюда следует, что они равны:
$J=J_{1}$.

Обозначим через $\mathcal{B}$ множество всех функций
$\psi:(0,\infty)\rightarrow(0,\infty)$ таких, что:

а)\;$\psi$ измерима по Борелю на полуоси $(0,\infty)$;

б) $\psi$ отделена от $0$ на каждом множестве $[r,\infty)$, где $r>0$;

в) $\psi$ ограничена на каждом отрезке $[a,b]$, где $0<a<b<\infty$.

Пусть $\psi\in\mathcal{B}$. В пространстве $X_{0}$ определен как функция от $J$
оператор $\psi(J)$. Обозначим через $[X_{0},X_{1}]_{\psi}$ или, короче, $X_{\psi}$
область определения оператора $\psi(J)$, наделенную скалярным произведением
$(u_{1},u_{2})_{X_{\psi}}:=(\psi(J)u_{1},\psi(J)u_{2})_{X_{0}}$ и соответствующей
нормой $\|\,u\,\|_{X_{\psi}}=(u,u)_{X_{\psi}}^{1/2}$.

Пространство $X_{\psi}$ сепарабельное гильбертово, причем справедливо непрерывное
плотное вложение $X_{\psi}\hookrightarrow X_{0}$. Действительно,
$\mathrm{Spec}\,J\subseteq[r,\infty)$ и $\psi(t)\geq c$ при $t\geq r$ для некоторых
положительных чисел $r,c$. (Как обычно, $\mathrm{Spec}\,J$ обозначает спектр
оператора $J$.) Поэтому $\mathrm{Spec}\,\psi(J)\subseteq[c,\infty)$, что влечет за
собой изометрический изоморфизм $\psi(J):X_{\psi}\leftrightarrow X_{0}$. Отсюда
вытекает полнота и сепарабельность пространства $X_{\psi}$ (а также то, что функция
$\|\cdot\|_{X_{\psi}}$ положительно определена и является нормой). Далее, поскольку
оператор $\psi^{-1}(J)$ ограничен в пространстве $X_{0}$, то оператор вложения
$I=\psi^{-1}(J)\psi(J):X_{\psi}\rightarrow X_{0}$ ограничен. Это вложение плотно,
так как область определения оператора $\psi(J)$ --- плотное линейное многообразие в
пространстве $X_{0}$.


\begin{remark}\label{rem2.1}
Пусть функции $\varphi,\psi\in\mathcal{B}$ такие, что $\varphi\asymp\psi$ в
окрестности $\infty$. Тогда в силу определения множества $\mathcal{B}$ справедливо
$\varphi\asymp\psi$ на $\mathrm{Spec}\,J$. Следовательно, $X_{\varphi}=X_{\psi}$ с
точностью до эквивалентности норм. (Как обычно, выражение $\varphi\asymp\psi$
означает, что обе функции $\varphi/\psi$ и $\psi/\varphi$ ограничены на
соответствующем множестве; при этом функции $\varphi$ и $\psi$ предполагаются
положительными.)
\end{remark}


\begin{definition}\label{def2.2}
Функцию $\psi\in\mathcal{B}$ назовем \textit{интерполяционным параметром}, если для
произвольных допустимых пар $X=[X_{0},X_{1}]$, $Y=[Y_{0},Y_{1}]$ гильбертовых
пространств и для любого линейного отображения $T$, заданного на $X_{0}$,
выполняется следующее условие. Если при $j\in\{0,1\}$ сужение отображения $T$ на
пространство $X_{j}$ является ограниченным оператором $T:X_{j}\rightarrow Y_{j}$, то
и сужение отображения $T$ на пространство $X_{\psi}$ является ограниченным
оператором $T:X_{\psi}\rightarrow Y_{\psi}$.
\end{definition}


Иначе говоря, параметр $\psi$ интерполяционный тогда и только тогда, когда
отображение $X\mapsto X_{\psi}$ является интерполяционным функтором, заданным на
категории допустимых пар $X$ гильбертовых пространств [\ref{BergLefstrem80}, с.~41;
\ref{Triebel80}, c.~18]. Если параметр $\psi$ интерполяционный, то будем говорить,
что пространство $X_{\psi}$ получено \textit{интерполяцией с функциональным
параметром} $\psi$ допустимой пары~$X$.

Далее мы исследуем основные свойства отображения $X\mapsto X_{\psi}$.

\subsection{Вложения пространств}\label{sec2.1.2}

Изучим некоторые свойства интерполяции, связанные с вложениями пространств.


\begin{theorem}\label{th2.1}
Пусть параметр $\psi\in\mathcal{B}$ интерполяционный, а $X=[X_{0},X_{1}]$~---
допустимая пара гильбертовых пространств. Тогда верны непрерывные плотные вложения
$X_{1}\hookrightarrow X_{\psi}\hookrightarrow X_{0}$.
\end{theorem}


\textbf{Доказательство.} В силу сказанного выше осталось доказать существование
непрерывного плотного вложения $X_{1}\hookrightarrow X_{\psi}$. Рассмотрим
ограниченные операторы вложения $I:X_{1}\rightarrow X_{0}$ и $I:X_{1}\rightarrow
X_{1}$. Поскольку параметр $\psi$ интерполяционный, они влекут ограниченность
оператора вложения $I:X_{1}\rightarrow X_{\psi}$, т. е. непрерывность вложения
$X_{1}\hookrightarrow X_{\psi}$.

Докажем его плотность. Возьмем произвольное $u\in X_{\psi}$. Тогда $v:=\psi(J)u\in
X_{0}$, где $J$ --- порождающий оператор для пары $X$. Поскольку $X_{1}$ плотно в
$X_{0}$, существует последовательность $(v_{k})\subset X_{1}$ такая, что
$v_{k}\rightarrow v$ в $X_{0}$ при $k\rightarrow\infty$. Отсюда вытекает сходимость
$$
u_{k}:=\psi^{-1}(J)v_{k}\rightarrow\psi^{-1}(J)v=u\quad\mbox{в}\quad
X_{\psi}\quad\mbox{при}\quad k\rightarrow\infty.
$$
Остается заметить, что
$$
u_{k}=\psi^{-1}(J)J^{-1}Jv_{k}=J^{-1}\psi^{-1}(J)Jv_{k} \in X_{1}.
$$
Теорема \ref{th2.1} доказана.


\begin{theorem}\label{th2.2}
Пусть функции $\psi,\chi\in\mathcal{B}$ такие, что функция $\psi/\chi$ ограничена в
окрестности $\infty$. Тогда для каждой допустимой пары $X=[X_{0},X_{1}]$
гильбертовых пространств справедливо непрерывное и плотное вложение
$X_{\chi}\hookrightarrow X_{\psi}$. Если вложение $X_{1}\hookrightarrow X_{0}$
компактно и $\psi(t)/\chi(t)\rightarrow0$ при $t\rightarrow\infty$, то вложение
$X_{\chi}\hookrightarrow X_{\psi}$ также компактно.
\end{theorem}


\textbf{Доказательство.} Пусть $J$ --- порождающий оператор для пары $X$. Заметим,
что $\mathrm{Spec}\,J\subseteq[r,\infty)$ для некоторого числа $r>0$. По условию
$\psi(t)/\chi(t)\leq c$ при $t\geq r$, следовательно,
$$
X_{\chi}=\mathrm{Dom}\,\chi(J)\subseteq\mathrm{Dom}\,\psi(J)=X_{\psi},
\quad\|\psi(J)\,u\|_{X_{0}}\leq c\,\|\chi(J)\,u\|_{X_{0}}.
$$
Отсюда на основании определения пространств $X_{\chi}$ и $X_{\psi}$ получаем
непрерывность вложения $X_{\chi}\hookrightarrow X_{\psi}$.

Докажем его плотность.

Рассмотрим изометрические изоморфизмы $\psi(J):X_{\psi}\leftrightarrow X_{0}$ и
$\chi(J):X_{\chi}\leftrightarrow X_{0}$. Пусть $u\in X_{\psi}$, тогда $\psi(J)\,u\in
X_{0}$. Поскольку пространство $X_{\chi}$ плотно вложено в $X_{0}$, существует
последовательность $(v_{k})\subset X_{\chi}$ такая, что $v_{k}\rightarrow\psi(J)\,u$
в $X_{0}$ при $k\rightarrow\infty$. Следовательно, $\psi^{-1}(J)\,v_{k}\rightarrow
u$ в $X_{\psi}$ при $k\rightarrow\infty$. Здесь
\begin{gather*}
\psi^{-1}(J)\,v_{k}=\psi^{-1}(J)\,\chi^{-1}(J)\,\chi(J)\,v_{k}=\\
\chi^{-1}(J)\,\psi^{-1}(J)\,\chi(J)\,v_{k} \in X_{\chi}.
\end{gather*}
Тем самым доказана плотность вложения $X_{\chi}\hookrightarrow X_{\psi}$.

Предположим теперь, что вложение $X_{1}\hookrightarrow X_{0}$ компактно и
$\psi(t)/\chi(t)\rightarrow0$ при $t\rightarrow\infty$. Докажем компактность
вложения $X_{\chi}\hookrightarrow X_{\psi}$. Пусть $(u_{k})$
--- произвольная ограниченная последовательность в
$X_{\chi}$. Поскольку последовательность элементов $w_{k}:=J^{-1}\,\chi(J)\,u_{k}$
ограничена в $X_{1}$, из нее можно выделить фундаментальную в $X_{0}$
подпоследовательность элементов $w_{k_{n}}=J^{-1}\,\chi(J)\,u_{k_{n}}$. Покажем, что
подпоследовательность $(u_{k_{n}})$ фундаментальна в $X_{\psi}$.

Пусть $E_{t}$, $t\geq r$,~--- разложение единицы в $X_{0}$, соответствующее
самосопряженному оператору\,$J$. Запишем:
\begin{gather}\notag 
\|u_{k_{n}}-u_{k_{m}}\|_{X_{\psi}}^{2}=
\|\psi(J)\,(u_{k_{n}}-u_{k_{m}})\|_{X_{0}}^{2}= \\ \notag
\|\psi(J)\,\chi^{-1}(J)\,J\,(w_{k_{n}}-w_{k_{m}})\|_{X_{0}}^{2}= \\
\int\limits_{r}^{\infty}\,\psi^{2}(t)\,\chi^{-2}(t)\,t^{2}\:
d\,\|E_{t}(w_{k_{n}}-w_{k_{m}})\|_{X_{0}}^{2}.\label{2.1}
\end{gather}
Выберем произвольное $\varepsilon>0$. Существует число
$\varrho=\varrho(\varepsilon)>r$ такое, что
$\psi(t)/\chi(t)\leq(2c_{0})^{-1}\varepsilon$ при $t\geq\varrho$ и
$$
c_{0}:=\sup\,\{\,\|w_{k}\|_{X_{1}}:\,k\in\mathbb{N}\,\}<\infty.
$$
Следовательно, для любых номеров $n,m$
\begin{gather} \notag 
\int\limits_{\varrho}^{\infty}\psi^{2}(t)\,\chi^{-2}(t)\,t^{2}\:
d\,\|E_{t}(w_{k_{n}}-w_{k_{m}})\|_{X_{0}}^{2}\,\leq \\ \notag
(2c_{0})^{-2}\,\varepsilon^{2}\,\int\limits_{\varrho}^{\infty}\,t^{2}\:
d\,\|E_{t}(w_{k_{n}}-w_{k_{m}})\|_{X_{0}}^{2}\,\leq \\ \notag
(2c_{0})^{-2}\,\varepsilon^{2}\,\|J\,(w_{k_{n}}-w_{k_{m}})\|_{X_{0}}^{2}\,=\\
\,(2c_{0})^{-2}\,\varepsilon^{2}\,\|w_{k_{n}}-w_{k_{m}}\|_{X_{1}}^{2}\,\leq
\,\varepsilon^{2}. \label{2.2}
\end{gather}
Кроме того, ввиду неравенства $\psi(t)/\chi(t)\leq c$ при $t\geq r$ имеем:
\begin{gather} \notag 
\int\limits_{r}^{\varrho}\,\psi^{2}(t)\,\chi^{-2}(t)\,t^{2}\:
d\,\|E_{t}(w_{k_{n}}-w_{k_{m}})\|_{X_{0}}^{2}\,\leq \\ \notag \,c^{2}\varrho^{2}\,
\int\limits_{r}^{\varrho}\,d\,\|E_{t}(w_{k_{n}}-w_{k_{m}})\|_{X_{0}}^{2}\,\leq \\
c^{2}\varrho^{2}\,\|w_{k_{n}}-w_{k_{m}}\|_{X_{0}}^{2}\rightarrow0
\quad\mbox{при}\quad n,m\rightarrow\infty. \label{2.3}
\end{gather}
Теперь соотношения \eqref{2.1} --- \eqref{2.3} влекут за собой неравенство
$\|u_{k_{n}}-u_{k_{m}}\|_{X_{\psi}}\leq2\varepsilon$ для достаточно больших номеров
$n,m$. Следовательно, подпоследовательность $(u_{k_{n}})$ фундаментальна в
пространстве $X_{\psi}$, что означает компактность вложения $X_{\chi}\hookrightarrow
X_{\psi}$.

Теорема \ref{th2.2} доказана.

\subsection{Свойство реитерации}\label{sec2.1.3}

Это свойство состоит в том, что повторное применение интерполяции с функциональным
параметром дает снова интерполяцию с некоторым функциональным параметром.


\begin{theorem}\label{th2.3}
Пусть $f,g,\psi\in\mathcal{B}$, причем функция $f/g$ ограничена в окрестности
$\infty$. Для каждой допустимой пары $X$ гильбертовых пространств справедливо
$[X_{f},X_{g}]_{\psi}=X_{\omega}$ с равенством норм. Здесь функция
$\omega\in\mathcal{B}$ определена по формуле $\omega(t):=f(t)\,\psi(g(t)/f(t))$ при
$t>0$. Если $f,g,\psi$ --- интерполяционные параметры, то $\omega$ --- тоже
интерполяционный параметр.
\end{theorem}


\textbf{Доказательство.} Поскольку функция $f/g$ ограничена в окрестности $\infty$,
пара $[X_{f},X_{g}]$ допустимая по теореме~\ref{th2.2}, и $\omega\in\mathcal{B}$.
Тем самым пространства $[X_{f},X_{g}]_{\psi}$ и $X_{\omega}$ определены. Докажем их
равенство.

Пусть оператор $J$~--- порождающий для пары $X=[X_{0},X_{1}]$, где
$\mathrm{Spec}\,J\subseteq[r,\infty)$ для некоторого числа $r>0$. Имеем
изометрические изоморфизмы
\begin{gather*}
f(J):X_{f}\leftrightarrow X_{0},\quad g(J):X_{g}\leftrightarrow X_{0},\\
B:=f^{-1}(J)\,g(J):X_{g}\leftrightarrow X_{f}.
\end{gather*}
Будем рассматривать $B$ как замкнутый оператор в пространстве $X_{f}$ с областью
определения $X_{g}$. Оператор $B$~--- порождающий для пары $[X_{f},X_{g}]$,
поскольку он положительно определенный и самосопряженный в $X_{f}$. Первое вытекает
из условия $f(t)/g(t)\leq c$ при $t\geq r$:
\begin{gather*}
(Bu,u)_{X_{f}}=(g(J)\,u,f(J)\,u)_{X_{0}}\geq \\
 c^{-1}\,(f(J)\,u,f(J)\,u)_{X_{0}}=c^{-1}\,\|u\|_{X_{f}}^{2},\quad u\in X_{f}.
\end{gather*}
Второе теперь следует из того, что $0$ --- регулярная точка оператора $B$.

При помощи спектральной теоремы приведем оператор $J$, самосопряженный в
пространстве $X_{0}$, к виду умножения на функцию: $J=I^{-1}(\alpha\cdot I)$. Здесь
$I:X_{0}\leftrightarrow L_{2}(U,d\mu)$  --- изометрический изоморфизм, $(U,\mu)$ ---
пространство с конечной мерой, $\alpha:U\rightarrow[r,\infty)$
--- измеримая функция. Изометрический изоморфизм $If(J):X_{f}\leftrightarrow
L_{2}(U,d\mu)$ приводит оператор $B$, самосопряженный в $X_{f}$, к виду умножения на
функцию $(g/f)\circ\alpha$:
\begin{gather*}
If(J)\,B\,u=Ig(J)\,u=(g\circ\alpha)\,Iu= \\
(g\circ\alpha)\,If^{-1}(J)f(J)\,u= ((g/f)\circ\alpha)\,If(J)\,u,\quad u\in X_{g}.
\end{gather*}
Следовательно, для произвольного $u\in X_{\psi}$ справедливо
\begin{gather*}
\|\psi(B)\,u\|_{X_{f}}=
\|(\psi\circ(g/f)\circ\alpha)\cdot(If(J)\,u)\|_{L_{2}(U,d\mu)}= \\
\|(\omega\circ\alpha)\cdot(Iu)\|_{L_{2}(U,d\mu)}=\|\omega(J)\,u\|_{X_{0}}.
\end{gather*}
Заметим здесь, что функция $f/\omega$ ограничена в окрестности $\infty$, откуда
$X_{\omega}\hookrightarrow X_{f}$ и поэтому $f(J)\,u$ определено. Таким образом,
доказано равенство $[X_{f},X_{g}]_{\psi}=X_{\omega}$.

Предположим теперь, что $f,g,\psi$ --- интерполяционные параметры. Покажем,что
$\omega$~--- интерполяционный параметр. Пусть допустимые пары $X=[X_{0},X_{1}]$,
$Y=[Y_{0},Y_{1}]$ и линейное отображение $T$~--- такие же как в
определении~\ref{def2.2}. Имеем ограниченные операторы $T:X_{f}\rightarrow Y_{f}$ и
$T:X_{g}\rightarrow Y_{g}$, что влечет за собой ограниченность оператора
$T:[X_{f},X_{g}]_{\psi}\rightarrow[Y_{f},Y_{g}]_{\psi}$. По уже доказанному,
$[X_{f},X_{g}]_{\psi}=X_{\omega}$ и $[Y_{f},Y_{g}]_{\psi}=Y_{\omega}$.
Следовательно, определен и ограничен оператор $T:X_{\omega}\rightarrow Y_{\omega}$,
т. е. $\omega$ --- интерполяционный параметр.

Теорема~\ref{th2.3} доказана.

\subsection[Интерполяция двойственных пространств]
{Интерполяция \\ двойственных пространств}\label{sec2.1.4}

Пусть $H$ --- ги\-льбертово пространство. Обозначим через $H'$ пространство,
двойственное к $H$, т. е. нормированное пространство всех линейных непрерывных
функционалов $l:H\rightarrow\mathbb{C}$. По теореме Ф.~Рисса, отображение
$S:v\mapsto(\,\cdot,v)_{H}$, где $v\in H$, определяет антилинейный изометрический
изоморфизм $S:H\leftrightarrow H'$. Отсюда следует, что пространство $H'$
гильбертово: норма в нем порождена скалярным произведением
$(l,m)_{H'}:=(S^{-1}l,S^{-1}m)_{H}$. Отметим, что мы не отождествляем $H$ и $H'$ с
помощью изоморфизма~$S$.


\begin{theorem}\label{th2.4}
Пусть функция $\psi\in\mathcal{B}$ такая, что функция $\psi(t)/t$ ограничена в
окрестности $\infty$. Для каждой допустимой пары $[X_{0},X_{1}]$ гильбертовых
пространств справедливо $[X_{1}',X_{0}']_{\psi}=[X_{0},X_{1}]_{\chi}'$ с равенством
норм. Здесь функция $\chi\in\mathcal{B}$ определена по формуле $\chi(t):=t/\psi(t)$
при $t>0$. Если $\psi$ --- интерполяционный параметр, то $\chi$ --- тоже
интерполяционный параметр.
\end{theorem}


\textbf{Доказательство.} Заметим, что пара $[X_{1}',X_{0}']$ допустимая при
естественном отождествлении функционалов из $X_{0}'$ с их сужениями на пространство
$X_{1}$. Из условия теоремы вытекает, что $\chi\in\mathcal{B}$. Тем самым
гильбертовы пространства $[X_{1}',X_{0}']_{\psi}$ и $[X_{0},X_{1}]_{\chi}'$
определены. Докажем их равенство.

Пусть оператор $J:X_{1}\leftrightarrow X_{0}$ порождающий для пары $[X_{0},X_{1}]$.
Рассмотрим изометрические изоморфизмы $S_{j}:X_{j}\leftrightarrow X_{j}'$,
$j=0,\,1$, фигурирующие в теореме Ф.~Рисса. Оператор $J'$, сопряженный с оператором
$J$, удовлетворяет равенству $J'=S_{1}J^{-1}S_{0}^{-1}$. Это вытекает из следующей
цепочки соотношений:
\begin{gather*}
(J'l)u=l(Ju)=(Ju,S_{0}^{-1}l)_{X_{0}}=(u,J^{-1}S_{0}^{-1}l)_{X_{1}}=\\
(S_{1}J^{-1}S_{0}^{-1}l)u\quad\mbox{для всех}\quad l\in X_{0}',\; u\in X_{1}.
\end{gather*}
Поэтому оператор $J'$ осуществляет изометрический изоморфизм
\begin{equation}\label{2.4}
J'=S_{1}J^{-1}S_{0}^{-1}:\,X_{0}'\leftrightarrow X_{1}'.
\end{equation}
Заметим здесь, что из равенств
\begin{gather*}
(u,JS_{1}^{-1}l)_{X_{0}}=(J^{-1}u,S_{1}^{-1}l)_{X_{1}}=l(J^{-1}u),\\
(u,J^{-1}S_{0}^{-1}l)_{X_{0}}=(J^{-1}u,S_{0}^{-1}l)_{X_{0}}=l(J^{-1}u),
\end{gather*}
справедливых для произвольных $l\in X_{0}'\hookrightarrow X_{1}'$, $u\in X_{0}$,
следует свойство
\begin{equation}\label{2.5}
JS_{1}^{-1}l=J^{-1}S_{0}^{-1}l\in X_{1}\quad\mbox{для всех}\quad l\in X_{0}'.
\end{equation}

Будем рассматривать $J'$ как замкнутый оператор в пространстве $X_{1}'$ с областью
определения $X_{0}'$. Оператор $J'$~--- порождающий для пары $[X_{1}',X_{0}']$,
поскольку он положительно определенный и самосопряженный в $X_{1}'$. Первое свойство
вытекает из положительной определенности оператора $J$ в пространстве $X_{0}$ и
свойства \eqref{2.5}:
\begin{gather*}
(J'l,l)_{X_{1}'}=(S_{1}J^{-1}S_{0}^{-1}l,l)_{X_{1}'}=
(J^{-1}S_{0}^{-1}l,S_{1}^{-1}l)_{X_{1}}=\\
(JJ^{-1}S_{0}^{-1}l,JS_{1}^{-1}l)_{X_{0}}=
(JJS_{1}^{-1}l,JS_{1}^{-1}l)_{X_{0}}\geq \\
c\,\|JS_{1}^{-1}l\|_{X_{0}}^{2}=c\,\|S_{1}^{-1}l\|_{X_{1}}^{2}=
c\,\|l\|_{X_{1}'}^{2}.
\end{gather*}
Здесь число $c>0$ не зависит от $l\in X_{0}'$. Второе свойство теперь следует из
того, что $0$~--- регулярная точка оператора $J'$ в силу \eqref{2.4}.

Воспользуемся приведением оператора $J$ к виду умножения на функцию
$J=I^{-1}(\alpha\cdot I)$ из доказательства теоремы~\ref{th2.3}. Изометрический
изоморфизм
\begin{equation}\label{2.6}
IJS_{1}^{-1}:\,X_{1}'\leftrightarrow L_{2}(U,d\mu)
\end{equation}
приводит оператор $J'$ к виду умножения на ту же функцию $\alpha$:
\begin{gather*}
(IJS_{1}^{-1})J'l=IS_{0}^{-1}l=IJJ^{-1}S_{0}^{-1}l=\\
\alpha\cdot IJ^{-1}S_{0}^{-1}l=\alpha\cdot IJS_{1}^{-1}l\quad\mbox{для любого}\quad
l\in X_{0}'
\end{gather*}
(мы воспользовались формулами \eqref{2.4} и \eqref{2.5}).

В силу теоремы~\ref{th2.2} выполняются непрерывные плотные вложения
$X_{0}'\hookrightarrow[X_{1}',X_{0}']_{\psi}$ и $[X_{0},X_{1}]_{\chi}\hookrightarrow
X_{0}$. Второе из них влечет непрерывность и плотность вложения
$X_{0}'\hookrightarrow[X_{0},X_{1}]_{\chi}'$. Покажем, что нормы в пространствах
$[X_{1}',X_{0}']_{\psi}$ и $[X_{0},X_{1}]_{\chi}'$ совпадают на плотном подмножестве
$X_{0}'$. Для произвольных $l\in X_{0}'$, $u\in[X_{0},X_{1}]_{\chi}$ запишем:
\begin{gather*}
l(u)=(u,S_{0}^{-1}l)_{X_{0}}=(\chi(J)u,\chi^{-1}(J)S_{0}^{-1}l)_{X_{0}}=\\
(v,\chi^{-1}(J)S_{0}^{-1}l)_{X_{0}},
\end{gather*}
где $v:=\chi(J)u\in X_{0}$. Отсюда
\begin{gather*}
\|l\|_{\,[X_{0},X_{1}]_{\chi}'}=
\sup\,\{\,|l(u)|\,/\,\|u\|_{\,[X_{0},X_{1}]_{\chi}}: u\in
[X_{0},X_{1}]_{\chi},\,u\neq0\,\}=\\
\sup\,\{\,|(v,\chi^{-1}(J)S_{0}^{-1}l)_{X_{0}}|\,/\,\|v\|_{X_{0}}: \,v\in
X_{0},\,v\neq0\,\}=\\
\|\chi^{-1}(J)S_{0}^{-1}l\|_{X_{0}}= \|I\chi^{-1}(J)S_{0}^{-1}l\|_{L_{2}(U,d\mu)}=\\
\|(\chi^{-1}\circ\alpha)\cdot IS_{0}^{-1}l\|_{L_{2}(U,d\mu)}.
\end{gather*}
С другой стороны, применяя изоморфизмы \eqref{2.6} и \eqref{2.4}, запишем:
\begin{gather*}
\|l\|_{\,[X_{1}',X_{0}']_{\psi}}=\|\psi(J')l\|_{X_{1}'}= \|\chi^{-1}(J')J'l\|_{X_{1}'}=\\
\|(IJS_{1}^{-1})\chi^{-1}(J')J'l\|_{L_{2}(U,d\mu)}=\\
\|(\chi^{-1}\circ\alpha)\cdot(IJS_{1}^{-1})J'l\|_{L_{2}(U,d\mu)}=\\
\|(\chi^{-1}\circ\alpha)\cdot IS_{0}^{-1}l\|_{L_{2}(U,d\mu)}.
\end{gather*}
Таким  образом, нормы в пространствах $[X_{1}',X_{0}']_{\psi}$ и
$[X_{0},X_{1}]_{\chi}'$ совпадают на плотном подмножестве $X_{0}'$. Поэтому
пространства совпадают.

Предположим теперь, что $\psi$ --- интерполяционный параметр. Покажем, что $\chi$
--- тоже интерполяционный параметр. Пусть допустимые пары $X=[X_{0},X_{1}]$,
$Y=[Y_{0},Y_{1}]$ и линейное отображение $T$~--- такие же как в
определении~\ref{def2.2}. Перейдя к сопряженному оператору $T'$, получим
ограниченные операторы $T':Y_{j}'\rightarrow X_{j}'$, где $j=0,\,1$. Поскольку
параметр $\psi$ интерполяционный, имеем ограниченный оператор
$T':[Y_{1}',Y_{0}']_{\psi}\rightarrow[X_{1}',X_{0}']_{\psi}$. По доказанному,
$[X_{1}',X_{0}']_{\psi}=[X_{0},X_{1}]_{\chi}'$ и
$[Y_{1}',Y_{0}']_{\psi}=[Y_{0},Y_{1}]_{\chi}'$ с равенством норм. Следовательно,
ограничен оператор $T':[Y_{0},Y_{1}]_{\chi}'\rightarrow[X_{0},X_{1}]_{\chi}'$.
Отсюда, перейдя ко второму сопряженному оператору $T''$, получим ограниченный
оператор $T'':[X_{0},X_{1}]_{\chi}''\rightarrow[Y_{0},Y_{1}]_{\chi}''$. Отождествляя
вторые сопряженные пространства с исходными гильбертовыми пространствами, будем
иметь ограниченный оператор $T:[X_{0},X_{1}]_{\chi}\rightarrow[Y_{0},Y_{1}]_{\chi}$.
Это означает, что $\chi$ --- интерполяционный параметр.

Теорема~\ref{th2.4} доказана.

\medskip

Отметим, что для некоторого достаточно широкого класса функциональных
интерполяционных параметров теорема~\ref{th2.4} доказана Г.~Шлензак
[\ref{Shlenzak74}, с.~52] (теорема~2).

В приложениях бывает удобно в качестве двойственного пространства $H'$ брать
пространство всех \emph{анти}линейных непрерывных функционалов на $H$. Для такого
пространства $H'$ теорема~\ref{th2.4}, очевидно, сохраняет силу. В дальнейшем из
контекста всегда будет ясно, из каких функционалов, линейных или антилинейных,
состоит пространство $H'$.

\subsection[Интерполяция ортогональных сумм прост\-ранств]
{Интерполяция ортогональных сумм \\ пространств}\label{sec2.1.5} Напомним, что
ортогональная сумма конечного или счетного множества сепарабельных гильбертовых
пространств сама является сепарабельным гильбертовым пространством.


\begin{theorem}\label{th2.5}
Пусть задано конечное или счетное множество допустимых пар гильбертовых пространств
$X^{(k)}:=[X_{0}^{(k)},X_{1}^{(k)}]$, $k\in\omega$. Предположим, что множество норм
операторов вложения $X_{1}^{(k)}\hookrightarrow X_{0}^{(k)}$, $k\in\omega$,
ограничено. Тогда для произвольной функции $\psi\in\mathcal{B}$ справедливо
$$
\Bigl[\:\bigoplus_{k\in\omega}X_{0}^{(k)},\,
\bigoplus_{k\in\omega}X_{1}^{(k)}\Bigr]_{\psi}=
\,\bigoplus_{k\in\omega}\,\bigl[X_{0}^{(k)},X_{1}^{(k)}\,\bigr]_{\psi}
$$
с равенством норм.
\end{theorem}


\textbf{Доказательство.} Пусть $\omega=\mathbb{N}$ (случай конечного множества
$\omega$ рассматривается аналогично и проще [\ref{Shlenzak74}, с.~53]). Пространства
$$
X_{0}:=\bigoplus_{k=1}^{\infty}X_{0}^{(k)},\quad
X_{1}:=\bigoplus_{k=1}^{\infty}X_{1}^{(k)}
$$
гильбертовы и сепарабельные. Непрерывность вложения $X_{1}\hookrightarrow X_{0}$
очевидна в силу условия теоремы. Покажем, что это вложение плотное. Пусть
$u:=(u_{1},u_{2},\ldots)\in X_{0}$. Для произвольных номеров $n,k$ существует
элемент $v_{n,k}\in X_{1}^{(k)}$ такой, что $\|u_{k}-v_{n,k}\|_{X_{0}^{(k)}}<1/n$.
Образуем последовательность векторов
$v^{(n)}:=(v_{n,1},\ldots,v_{n,n},0,0,\ldots)\in X_{1}$. Имеем:
\begin{gather*}
\|u-v^{(n)}\|_{X_{0}}^{2}=\sum_{k=1}^{n}\|u_{k}-v_{n,k}\|_{X_{0}}^{2}+
\sum_{k=n+1}^{\infty}\|u_{k}\|_{X_{0}}^{2}\leq \\
\frac{n}{n^{2}}+
\sum_{k=n+1}^{\infty}\|u_{k}\|_{X_{0}}^{2}\rightarrow0\;\;\mbox{при}
\;\;n\rightarrow\infty.
\end{gather*}
Значит, $X_{1}$ плотно в пространстве $X_{0}$, и пара $X:=[X_{0},X_{1}]$ допустимая.

Пусть оператор $J_{k}$ порождающий для пары $X^{(k)}$. Непосредственно проверяется,
что оператор $J:=(J_{1},J_{2},\ldots)$ порождающий для пары $X$. Покажем, что
$\psi(J)=(\psi(J_{1}),\psi(J_{2}),\ldots)$, где
$\mathrm{Dom}\,\psi(J)=\bigoplus_{k=1}^{\infty}X_{\psi}^{(k)}$. Приведем оператор
$J_{k}$ к виду умножения на функцию: $I_{k}J_{k}=\alpha_{k}\cdot I_{k}$. Здесь
$I_{k}:X_{0}^{(k)}\leftrightarrow L_{2}(V_{k},d\mu_{k})$
--- изометрический изоморфизм, $V_{k}$ --- пространство с конечной
мерой $\mu_{k}$, функция $\alpha_{k}:V_{k}\rightarrow(0,\infty)$ измеримая. Не
уменьшая общности, можно считать, что множества $V_{k}$ попарно не пересекаются.
Положим $V:=\bigcup_{k=1}^{\infty}V_{k}$. Подмножество $\Omega\subseteq V$ называем
измеримым, если для любого номера $k$ множество $\Omega\cap V_{k}$ измеримо
относительно меры $\mu_{k}$. На $\sigma$-алгебре измеримых множеств $\Omega\subseteq
V$ введем $\sigma$-конечную меру $\mu(\Omega):=\sum_{k=1}^{\infty}\mu_{k}(\Omega\cap
V_{k})$. Для произвольного вектора $u:=(u_{1},u_{2},\ldots)\in X_{0}$ определяем
измеримые функции $Iu$ и $\alpha$ на множестве $V$ по формулам
$(Iu)(\lambda):=(I_{k}u_{k})(\lambda)$ и $\alpha(\lambda):=\alpha_{k}(\lambda)$ при
$\lambda\in V_{k}$. Имеем изометрический изоморфизм $I:X_{0}\leftrightarrow
L_{2}(V,d\mu)$. Он приводит оператор $J$ к виду умножения на функцию $\alpha$,
поскольку
$$
(IJu)(\lambda)=(I_{k}J_{k}u_{k})(\lambda)=
\alpha_{k}(\lambda)(I_{k}u_{k})(\lambda)=\alpha(\lambda)(Iu)(\lambda)
$$
для любых $u\in X_{1}$ и $\lambda\in V_{k}$. Теперь можем записать:
\begin{gather*}
X_{\psi}=\mathrm{Dom}\,\psi(J)=\{\,u\in X_{0}:\,(\psi\circ\alpha)\cdot(Iu)\in
L_{2}(V,d\mu)\,\}=\\
\{\,u\in X_{0}:\,\sum_{k=1}^{\infty}\:\|(\psi\circ\alpha_{k})\cdot(I_{k}u_{k})\|
_{L_{2}(V_{k},d\mu_{k})}^{2}<\infty\,\}=\\
\{\,u:\,u_{k}\in\mathrm{Dom}\,\psi(J_{k}),\;
\sum_{k=1}^{\infty}\:\|\psi(J_{k})u_{k}\|
_{X_{0}^{(k)}}^{2}<\infty\,\}=\bigoplus_{k=1}^{\infty}X_{\psi}^{(k)}.
\end{gather*}
Кроме того, для любого $u\in \mathrm{Dom}\,\psi(J)$ справедливо
\begin{gather*}
(I\psi(J)u)(\lambda)=\psi(\alpha(\lambda))\,(Iu)(\lambda)=
\psi(\alpha_{k}(\lambda))\,(I_{k}u_{k})(\lambda)=\\
(I_{k}\psi(J_{k})u_{k})(\lambda)=
\bigl(I(\psi(J_{1})u_{1},\psi(J_{2})u_{2},\ldots)\bigr)(\lambda)
\quad\mbox{при}\;\;\lambda\in V_{k}.
\end{gather*}
Следовательно, $\psi(J)u=(\psi(J_{1})u_{1},\psi(J_{2})u_{2},\ldots)$, откуда
$$
\|u\|_{X_{\psi}}^{2}=\|\psi(J)u\|_{X_{0}}^{2}=
\sum_{k=1}^{\infty}\,\|\psi(J_{k})u_{k}\|_{X_{0}^{(k)}}^{2}=
\sum_{k=1}^{\infty}\,\|u_{k}\|_{X_{\psi}^{(k)}}^{2}.
$$

Теорема~\ref{th2.5} доказана.

\subsection[Интерполяция подпространств и факторпространств]
{Интерполяция подпространств \\ и факторпространств}\label{sec2.1.6} Напомним, что
по определению подпространство гильбертова пространства замкнуто. Далее будут
рассматриваться, вообще говоря, неортогональные проекторы на подпространство.


\begin{theorem}\label{th2.6}
Пусть $X=[X_{0},X_{1}]$ --- допустимая пара гильбертовых пространств, а $Y_{0}$
--- подпространство в $X_{0}$. Тогда $Y_{1}:=X_{1}\cap Y_{0}$ ---
подпространство в $X_{1}$. Предположим, что существует линейное отображение $P$,
которое для каждого $j\in\{0,1\}$ является проектором пространства $X_{j}$ на
подпространство $Y_{j}$. Тогда пары $[Y_{0},Y_{1}]$, $[X_{0}/Y_{0},X_{1}/Y_{1}]$
допустимые и для произвольного интерполяционного параметра $\psi\in\mathcal{B}$
справедливы следующие равенства пространств с эквивалентностью норм в них:
\begin{equation}\label{2.7}
[Y_{0},Y_{1}]_{\psi}=X_{\psi}\cap Y_{0}, \quad
[X_{0}/Y_{0},X_{1}/Y_{1}]_{\psi}=X_{\psi}/(X_{\psi}\cap Y_{0}).
\end{equation}
Здесь $X_{\psi}\cap Y_{0}$ --- подпространство в $X_{\psi}$.
\end{theorem}


\textbf{Доказательство.} Поскольку вложения $X_{1}\hookrightarrow X_{0}$ и
$X_{\psi}\hookrightarrow X_{0}$ непрерывны, линейные многообразия $Y_{1}=X_{1}\cap
Y_{0}$ и $X_{\psi}\cap Y_{0}$ замкнуты в пространствах $X_{1}$ и $X_{\psi}$
соответственно. Пары $[Y_{0},Y_{1}]$ и $[X_{0}/Y_{0},X_{1}/Y_{1}]$ допустимые, так
как $Y_{1}$ --- плотное подмножество в пространстве $Y_{0}$. Последнее доказывается
следующим образом. Для произвольного $u\in Y_{0}$ существует последовательность
элементов $u_{k}\in X_{1}$ такая, что $u_{k}\rightarrow u$ в $X_{0}$ при
$k\rightarrow\infty$. Отсюда в силу условия теоремы получаем $Pu_{k}\rightarrow
Pu=u$ в $X_{0}$ при $k\rightarrow\infty$, где $Pu_{k}\in Y_{1}$, т. е. упомянутую
плотность. Таким образом, левые и правые части равенств \eqref{2.7} определены
корректно. Доказательство этих равенств проводится аналогично доказательству теорем
1.17.1/1 и 1.17.2/1 из книги Х.~Трибеля [\ref{Triebel80}, c.~136, 138].
Теорема~\ref{th2.6} доказана.

\subsection{Интерполяция нетеровых операторов}\label{sec2.1.7}
Напомним следующее определение.


\begin{definition}\label{def2.3}
Линейный ограниченный оператор $T:X\rightarrow Y$, где $X$, $Y$~--- банаховы
пространства, называется \emph{нетеровым}, если его ядро $\ker T$ и коядро
$\mathrm{coker}\,T:=Y/\,T(X)$ конечномерные. \emph{Индексом} нетерового оператора
$T$ называется число $\mathrm{ind}\,T:=\dim\ker T-\dim(Y/\,T(X))$.
\end{definition}

Отметим, что область значений нетерового оператора замкнута [\ref{Hermander87},
с.~246] (лемма~19.1.1).


\begin{theorem}\label{th2.7}
Пусть задано две допустимые пары $X=[X_{\,0},X_{1}]$ и $Y=[Y_{\,0},Y_{1}]$
гильбертовых пространств. Пусть, кроме того, на $X_{\,0}$ задано линейное
отображение $T$, для которого существуют ограниченные нетеровы операторы
$T:X_{j}\rightarrow Y_{j}$, где $j\in\{0,1\}$, имеющие общее ядро $N$ и одинаковый
индекс $\varkappa$. Тогда для произвольного интерполяционного параметра
$\psi\in\mathcal{B}$ ограниченный оператор $T:X_{\psi}\rightarrow Y_{\psi}$ нетеров
с ядром $N$, областью значений $Y_{\psi}\cap T(X_{\,0})$ и тем же индексом
$\varkappa$.
\end{theorem}


\textbf{Доказательство.} По условию, имеем ограниченные операторы
$T_{j}:=T:X_{j}\rightarrow Y_{j}$ и $T_{\psi}:=T:X_{\psi}\rightarrow Y_{\psi}$.
Рассмотрим сопряженные к ним ограниченные операторы $T_{j}':Y_{j}'\rightarrow
X_{j}'$ и $T_{\psi}':Y_{\psi}'\rightarrow X_{\psi}'$. В силу теоремы~\ref{th2.1}
выполняются непрерывные плотные вложения $X_{1}\hookrightarrow
X_{\psi}\hookrightarrow X_{0}$ и $Y_{1}\hookrightarrow Y_{\psi}\hookrightarrow
Y_{0}$. Поэтому
$$
\ker T_{1}\subseteq \ker T_{\psi}\subseteq \ker T_{0}\quad\mbox{и}\quad\ker
T_{0}'\subseteq \ker T_{\psi}'\subseteq \ker T_{1}'.
$$
Но, по условию, $\ker T_{j}=N$ и
$$
\dim\ker T_{j}'=\dim Y_{j}/T(X_{j})=-\varkappa+\dim N.
$$
Следовательно, $\ker T_{\psi}=N$ и $\ker T_{\psi}'=\ker T_{j}'=:M$, где $\dim
M=-\varkappa+\dim N$. Итак, оператор $T:X_{\psi}\rightarrow Y_{\psi}$ имеет
конечномерные ядро $N$ и дефектное подпространство $M$, удовлетворяющие равенству
$\dim N-\dim M=\varkappa$. Остается показать, что $T(X_{\psi})=Y_{\psi}\cap
T(X_{0})$, поскольку тогда $T(X_{\psi})$ замкнуто в пространстве $Y_{\psi}$ и $\dim
Y_{\psi}/T(X_{\psi})=\dim M$.

Для этого рассмотрим гомеоморфизмы
\begin{equation}\label{2.8}
T:\,X_{j}/N\leftrightarrow T(X_{j})\quad\mbox{при}\quad j\in\{0,1\}.
\end{equation}
Они каноническим образом порождены ограниченными нетеровыми операторами $T_{j}$.
Напомним, что $T(X_{j})$ --- подпространство в $Y_{j}$. Применив к \eqref{2.8}
интерполяцию с параметром $\psi$, получим еще один гомеоморфизм
\begin{equation}\label{2.9}
T:\,[X_{0}/N,X_{1}/N]_{\psi}\leftrightarrow [T(X_{0}),T(X_{1})]_{\psi}.
\end{equation}
Заметим здесь, что пара $[X_{0}/N,X_{1}/N]$, очевидно, допустимая, откуда в силу
\eqref{2.8} вытекает допустимость пары $[T(X_{0}),T(X_{1})]$. Опишем
интерполяционные пространства, в которых действует гомеоморфизм \eqref{2.9}.

Рассмотрим ортогональную сумму $X_{0}=N\oplus E$. Ее сужение на пространство $X_{1}$
является прямой суммой подпространств $X_{1}=N\dotplus(E\cap X_{1})$. Обозначим
через $P$ ортопроектор пространства $X_{0}$ на $N$. Тогда его сужение на
пространство $X_{1}$ является проектором этого пространства на $N$, соответствуюшим
второй сумме. Следовательно, в силу теоремы~\ref{th2.6} и равенства $X_{\psi}\cap
N=N$ справедливо
\begin{equation}\label{2.10}
[X_{0}/N,X_{1}/N]_{\psi}=X_{\psi}/N
\end{equation}
с эквивалентностью норм.

Далее, рассмотрим ортогональную сумму $Y_{0}=T(X_{0})\oplus Z_{0}$, где $\dim
Z_{0}<\infty$. Отсюда, поскольку $Y_{1}$ плотно в пространстве $Y_{0}$, мы получаем
в силу леммы 2.1 из работы [\ref{GohbergKrein57}] разложение пространства $Y_{0}$ в
прямую сумму $Y_{0}=T(X_{0})\dotplus Z_{1}$, где $Z_{1}$
--- некоторое конечномерное подпространство в $Y_{1}$. Сужение этой суммы на
пространство $Y_{1}$ имеет вид $Y_{1}=(T(X_{0})\cap Y_{1})\dotplus Z_{1}$. Заметим
здесь, что $T(X_{0})\cap Y_{1}=T(X_{1})$. Это равенство вытекает из следующего
представления замкнутой области значений оператора $T_{j}$:
$$
T(X_{j})=\{f\in Y_{j}:\,l(f)=0\;\mbox{для любого}\;l\in M\}
$$
при $j\in\{0,1\}$. Таким образом, мы имеем прямые суммы $Y_{0}=T(X_{0})\dotplus
Z_{1}$ и $Y_{1}=T(X_{1})\dotplus Z_{1}$, причем вторая сумма является сужением
первой. Обозначим через $Q$ проектор пространства $Y_{0}$ на подпространство
$T(X_{0})$, соответствующим первой сумме. Тогда сужение оператора $Q$ на $Y_{1}$
является проектором пространства $Y_{1}$ на подпространство $T(X_{1})=Y_{1}\cap
T(X_{0})$, соответствующим второй сумме. Поэтому в силу теоремы~\ref{th2.4} можно
записать
\begin{equation}\label{2.11}
[T(X_{0}),T(X_{1})]_{\psi}=Y_{\psi}\cap T(X_{0})
\end{equation}
с эквивалентностью норм.

Теперь, применяя равенства \eqref{2.10} и \eqref{2.11} к гомеоморфизму \eqref{2.9},
получаем гомеоморфизм $T:X_{\psi}/N\leftrightarrow Y_{\psi}\cap T(X_{0})$. Отсюда
следует равенство
$$T(X_{\psi})=T(X_{\psi}/N)=Y_{\psi}\cap
T(X_{0}),
$$
на котором и заканчивается доказательство теоремы~\ref{th2.7}.

\medskip

Отметим, что аналог теоремы~\ref{th2.7} доказан Ж.~Жеймона [\ref{Geymonat65},
с.~281] для произвольного интерполяционого функтора, заданного на категории всех
совместимых пар банаховых пространств. Приведенное им доказательство схоже с нашим.

\subsection[Оценка нормы оператора в интерполяционных
пространствах]{Оценка нормы оператора \\ в интерполяционных
пространствах}\label{sec2.1.8}

Докажем следующую теорему.


\begin{theorem}\label{th2.8}
Для заданных интерполяционного параметра $\psi\in\mathcal{B}$ и числа $m>0$
существует число $c=c(\psi,m)>0$ такое, что
\begin{equation}\label{2.12}
\|T\|_{X_{\psi}\rightarrow Y_{\psi}}\leq c \max\,\bigl\{\,\|T\|_{X_{j}\rightarrow
Y_{j}}:\,j=0,\,1\,\bigr\}.
\end{equation}
Здесь $X=[X_{0},X_{1}]$ и $Y=[Y_{0},Y_{1}]$ --- произвольные допустимые пары
гильбертовых пространств, для которых нормы операторов вложений
$X_{1}\hookrightarrow X_{0}$ и $Y_{1}\hookrightarrow Y_{0}$ не превосходят число
$m$, а $T$~--- произвольное линейное отображение, заданное на пространстве $X_{0}$ и
определяющее ограниченные операторы $T:X_{j}\rightarrow Y_{j}$ при $j\in\{0,1\}$.
Постоянная $c$ не зависит от $X$, $Y$ и $T$.
\end{theorem}


\textbf{Доказательство.} Предположим, что теорема не верна. Тогда
\begin{equation}\label{2.13}
\|T_{k}\|_{X^{(k)}_{\psi}\rightarrow Y^{(k)}_{\psi}}>k\,m_{k} \quad\mbox{для
каждого}\;\;k\in\mathbb{N}.
\end{equation}
Здесь $X^{(k)}:=[X_{0}^{(k)},X_{1}^{(k)}]$ и
$Y^{(k)}:=[Y_{0}^{(k)},Y_{1}^{(k)}]$~--- некоторые допустимые пары гильбертовых
пространств, для которых нормы операторов вложений $X_{1}^{(k)}\hookrightarrow
X_{0}^{(k)}$ и $Y_{1}^{(k)}\hookrightarrow Y_{0}^{(k)}$ не превосходят число $m$, а
$T_{k}$~--- некоторое линейное отображение, заданное на пространстве $X_{0}^{(k)}$ и
определяющее ограниченные операторы $T_{k}:X_{j}^{(k)}\rightarrow Y_{j}^{(k)}$, где
$j\in\{0,1\}$. При этом используется обозначение
$$
m_{k}:=\max\bigl\{\,\|T_{k}\|_{X^{(k)}_{0}\rightarrow
Y^{(k)}_{0}},\;\|T_{k}\|_{X^{(k)}_{1}\rightarrow Y^{(k)}_{1}}\,\bigr\}>0.
$$
Рассмотрим ограниченные операторы
\begin{gather}
T:u=(u_{1},u_{2},\ldots)\mapsto
(m_{1}^{-1}\,T_{1}\,u_{1},\,m_{2}^{-1}\,T_{2}\,u_{2},\ldots),\notag \\
T:\,\bigoplus_{k=1}^{\infty}X_{j}^{(k)}\rightarrow\bigoplus_{k=1}^{\infty}Y_{j}^{(k)},
\quad j\in\{0,1\}. \label{2.14}
\end{gather}
Их ограниченность вытекает из следующих неравенств
\begin{gather*}
\sum_{k=1}^{\infty}\;\|m_{k}^{-1}\,T_{k}u_{k}\|_{Y_{j}^{(k)}}^{2}\leq \\
\sum_{k=1}^{\infty}\;m_{k}^{-2}\,\|T_{k}\|_{X_{j}^{(k)}\rightarrow
Y_{j}^{(k)}}^{2}\,\|u_{k}\|_{X_{j}^{(k)}}^{2}\leq
\sum_{k=1}^{\infty}\;\|u_{k}\|_{X_{j}^{(k)}}^{2}.
\end{gather*}
Поскольку параметр $\psi$ интерполяционный, то из ограниченности операторов
\eqref{2.14} вытекает, что следующий оператор определен и ограничен:
$$
T:\,\Bigl[\:\bigoplus_{k=1}^{\infty}X_{0}^{(k)},\,
\bigoplus_{k=1}^{\infty}X_{1}^{(k)}\Bigr]_{\psi}\rightarrow\;
\Bigl[\:\bigoplus_{k=1}^{\infty}Y_{0}^{(k)},\,
\bigoplus_{k=1}^{\infty}Y_{1}^{(k)}\Bigr]_{\psi}.
$$
Это в силу теоремы~\ref{th2.5} свидетельствует об ограниченности оператора
$$
T:\,\bigoplus_{k=1}^{\infty}X_{\psi}^{(k)}\rightarrow\;
\bigoplus_{k=1}^{\infty}Y_{\psi}^{(k)}.
$$
Пусть $c_{0}$ --- норма этого оператора. Для каждого номера $k$ рассмотрим вектор
$u^{(k)}:=(u_{1},\ldots,u_{k},\ldots)$ такой, что $u_{k}\in X_{\psi}^{(k)}$ и
$u_{j}=0$ при $j\neq k$. Имеем:
\begin{gather*}
\|T_{k}u_{k}\|_{Y_{\psi}^{(k)}}=
m_{k}\,\|Tu^{(k)}\|_{\bigoplus_{j=1}^{\infty}Y_{\psi}^{(j)}}\leq \\
m_{k}\,c_{0}\,\|u^{(k)}\|_{\bigoplus_{j=1}^{\infty}X_{\psi}^{(j)}}=
m_{k}\,c_{0}\,\|u_{k}\|_{X_{\psi}^{(k)}}
\end{gather*}
для всех $u_{k}\in X_{\psi}^{(k)}$. Следовательно,
$$
\|T_{k}\|_{X^{(k)}_{\psi}\rightarrow Y^{(k)}_{\psi}}\leq c_{0}\,m_{k} \quad\mbox{для
каждого номера}\;\;k,
$$
что противоречит условию \eqref{2.13}. Таким образом, наше предположение ложно, и
значит, теорема~\ref{th2.8} верна.

\medskip

Отметим, что неравенство \eqref{2.12}, где постоянная $c$ не зависит от $T$ (но
может зависеть от $X$ и $Y$), выполняется для любого интерполяционного функтора,
заданного на категории пар гильбертовых либо банаховых пространств
[\ref{KreinPetuninSemenov78}, с.~34]. Теорема~\ref{th2.8} усиливает этот факт
применительно к интерполяционному функтору $X\mapsto X_{\psi}$.

Допустимую пару $[X_{0},X_{1}]$ гильбертовых пространств называем \emph{нормальной},
если $\|u\|_{X_{0}}\leq\|u\|_{X_{1}}$ для любого $u\in X_{1}$. В силу теоремы
\ref{th2.8} постоянная $c$ в неравенстве \eqref{2.12} не зависит от допустимых пар
$X$, $Y$ и оператора $T$, если эти пары нормальные. Заметим, что всякую допустимую
пару $[X_{0},X_{1}]$ можно сделать нормальной, заменив, например, в пространстве
$X_{0}$ норму $\|u\|_{X_{0}}$ на пропорциональную норму $k\,\|u\|_{X_{0}}$, где
$0<k\leq m^{-1}$, а $m$
--- норма оператора вложения $X_{1}\hookrightarrow X_{0}$.

\subsection[Критерий интерполяционности функционального
параметра]{Критерий интерполяционности \\ функционального параметра}\label{2.1.9}
Опираясь на результаты Ж.~Петре [\ref{Peetre68}; \ref{BergLefstrem80}, с.~153],
докажем критерий того, что функция $\psi\in\mathcal{B}$ является интерполяционным
параметром.


\begin{definition}\label{def2.4}
Пусть заданы функция $\psi:(0,\infty)\rightarrow(0,\infty)$ и число $r\geq0$.
Функция $\psi$ называется \emph{псевдовогнутой} на полупрямой $(r,\infty)$, если
существует вогнутая функция $\psi_{1}:(r,\infty)\rightarrow(0,\infty)$ такая, что
$\psi(t)\asymp \psi_{1}(t)$ при $t>r$. Функция $\psi$ называется псевдовогнутой в
окрестности $\infty$, если она псевдовогнута на некоторой полупрямой $(r,\infty)$,
где $r$~--- достаточно большое число.
\end{definition}


\begin{theorem}\label{th2.9}
Функция $\psi\in\mathcal{B}$ является интерполяционным параметром тогда и только
тогда, когда она псевдовогнута в окрестности~$\infty$.
\end{theorem}


Доказательству этой теоремы предпошлем две леммы.


\begin{lemma}\label{lem2.1}
Пусть функция $\psi$ принадлежит множеству $\mathcal{B}$ и псевдовогнута в
окрестности $\infty$. Тогда существует вогнутая функция
$\psi_{0}:(0,\infty)\rightarrow(0,\infty)$ такая, что для каждого числа
$\varepsilon>0$ справедливо соотношение $\psi(t)\asymp \psi_{0}(t)$ при
$t\geq\varepsilon$.
\end{lemma}


\textbf{Доказательство леммы~\ref{lem2.1}.} Согласно условию существуют такие число
$r\gg1$ и вогнутая функция $\psi_{1}:(r,\infty)\rightarrow(0,\infty)$, что
$\psi(t)\asymp \psi_{1}(t)$ при $t>r$. Поскольку функция $\psi_{1}$ вогнута и
положительна на полуоси $(r,\infty)$, она там (нестрого) возрастает. Кроме того, для
каждой фиксированной точки $t_{0}\in(r,\infty)$ функция наклона
$(\psi_{1}(t)-\psi_{1}(t_{0}))/(t-t_{0})$, $t\in(r,\infty)\setminus \{t_{0}\}$,
(нестрого) убывает. Следовательно, в точке $r+1$ функция $\psi_{1}$ имеет правую
касательную, образующую острый или нулевой угол с осью абсцисс. Зададим на полуоси
$[0,\infty)$ функцию $\psi_{2}$ так, чтобы ее график совпадал на промежутке
$[0,r+1)$ с указанной касательной, а на полуоси $[r+1,\infty)$ --- с графиком
функции $\psi_{1}$. Функция $\psi_{2}$ возрастает на $[0,\infty)$ и вогнута на
$(0,\infty)$. Последнее вытекает из того, что для каждой фиксированной точки
$t_{0}\in(0,\infty)$ функция наклона $(\psi_{2}(t)-\psi_{2}(t_{0}))/(t-t_{0})$,
$t\in(0,\infty)\setminus \{t_{0}\}$, убывает. Положим
$\psi_{0}(t):=\psi_{2}(t)+|\psi_{2}(0)|$+1. Функция $\psi_{0}$ положительна,
возрастает и вогнута на полуоси $(0,\infty)$. Выберем произвольное число
$\varepsilon>0$. Заметим, что $\psi(t)\asymp1\asymp\psi_{0}(t)$ при
$t\in[\varepsilon,r+1+\varepsilon]$. Далее, поскольку функция $\psi_{2}$ возрастает
и положительна на полуоси $[r+1,\infty)$, имеем: $|\psi_{2}(0)|+1\leq
c\,\psi_{2}(t)$ при $t\geq r+1$, где $c:=(|\psi_{2}(0)|+1)/\psi_{2}(r+1)>0$. Отсюда
получаем соотношение $\psi(t)\asymp\psi_{1}(t)=\psi_{2}(t)\asymp\psi_{0}(t)$ при
$t\geq r+1$. Таким образом, $\psi(t)\asymp\psi_{0}(t)$ при $t\geq\varepsilon$, что и
требовалось доказать.


\begin{lemma}\label{lem2.2} Пусть заданы функция $\psi\in\mathcal{B}$
и число $r\geq0$. Функция $\psi$ псевдовогнута на полупрямой $(r,\infty)$ тогда и
только тогда, когда существует число $c>0$ такое, что
$$
\psi(t)/\psi(s)\leq c\,\max\{1,\,t/s\}\quad\mbox{для любых}\quad t,s>r.
$$
\end{lemma}


\textbf{Доказательство леммы~\ref{lem2.2}.} В случае $r=0$ эта лемма доказана
Ж.~Петре [\ref{Peetre68}; \ref{BergLefstrem80}, с.~152, 153] (при этом условие
$\psi\in\mathcal{B}$ излишнее). В~случае $r>0$ достаточность устанавливается
аналогично. Необходимость же сводится к случаю $r=0$ с помощью леммы~\ref{lem2.1}.
Действительно, предположив, что $\psi$ псевдовогнута на $(r,\infty)$, рассмотрим
функцию $\psi_{0}$ из этой леммы, где $\varepsilon=r$. Тогда
$$
\psi(t)/\psi(s)\asymp\psi_{0}(t)/\psi_{0}(s)\leq c_{0}\max\{1,\,t/s\}
$$
для любых $t,s>r$. (На самом деле $c_{0}=1$ для вогнутой функции $\psi_{0}$
[\ref{Peetre68}]). Лемма~\ref{lem2.2} доказана.


\medskip

\textbf{Доказательство теоремы~\ref{th2.9}.} \textit{Достаточность.} Предположим,
что функция $\psi\in\mathcal{B}$ псевдовогнута в окрестности $\infty$. Докажем, что
она является интерполяционным параметром.

Выберем произвольные допустимые пары $X=[X_{0},X_{1}]$, $Y=[Y_{0},Y_{1}]$
гильбертовых пространств и линейное отображение $T$ такие же как в
определении~\ref{def2.2}. Пусть операторы $J_{X}:X_{1}\leftrightarrow X_{0}$ и
$J_{Y}:Y_{1}\leftrightarrow Y_{0}$ порождающие для пар $X$ и $Y$ соответственно. При
помощи спектральной теоремы приведем эти операторы, самосопряженные в $X_{0}$ и
$Y_{0}$, к виду умножения на функцию:
\begin{equation}\label{2.15}
J_{X}=I_{X}^{-1}\,(\alpha\cdot I_{X})\quad\mbox{и}\quad
J_{Y}=I_{Y}^{-1}\,(\beta\cdot I_{Y}).
\end{equation}
Здесь $I_{X}:X_{0}\leftrightarrow L_{2}(U,d\mu)$ и $I_{Y}:Y_{0}\leftrightarrow
L_{2}(V,d\nu)$ --- изометрические изоморфизмы, $(U,\mu)$ и $(V,\nu)$ ---
пространства с конечными мерами, а $\alpha:U\rightarrow(0,\infty)$ и
$\beta:V\rightarrow(0,\infty)$ --- измеримые функции. Поскольку ограничены операторы
$T:X_{0}\rightarrow Y_{0}$ и $T:X_{1}\rightarrow Y_{1}$, то также ограничены и
операторы
\begin{gather}\label{2.16}
I_{Y}\,T\,I_{X}^{-1}:\,L_{2}(U,d\mu)\rightarrow L_{2}(V,d\nu),\\
I_{Y}\,J_{Y}\,T\,J_{X}^{-1}\,I_{X}^{-1}:\,L_{2}(U,d\mu)\rightarrow L_{2}(V,d\nu).
\label{2.17}
\end{gather}
В силу \eqref{2.15} имеем:
$$
I_{Y}\,J_{Y}\,T\,J_{X}^{-1}\,I_{X}^{-1}=(\beta\cdot
I_{Y})\,T\,I_{X}^{-1}(\alpha^{-1}\cdot).
$$
Следовательно, \eqref{2.17} влечет ограниченность оператора
\begin{gather}\notag
I_{Y}\,T\,I_{X}^{-1}=\\
\beta^{-1}\cdot(I_{Y}\,J_{Y}\,T\,J_{X}^{-1}\,I_{X}^{-1})
(\alpha\cdot):\,L_{2}(U,\alpha^{2}d\mu)\rightarrow
L_{2}(V,\beta^{2}d\nu).\label{2.18}
\end{gather}

Пусть вогнутая функция $\psi_{0}:(0,\infty)\rightarrow(0,\infty)$ такая же как в
лемме~\ref{lem2.1}. Отметим, что $\psi_{0}\in\mathcal{B}$ и (см.
замечание~\ref{sec2.1})
\begin{equation}\label{2.19}
X_{\psi}=X_{\psi_{0}},\quad Y_{\psi}=Y_{\psi_{0}}\quad\mbox{с эквивалентностью
норм}.
\end{equation}
Как установил Ж.~Петре [\ref{Peetre68}; \ref{BergLefstrem80}, с.~153],
псевдовогнутость положительной функции на полуоси $(0,\infty)$ равносильна ее
интерполяционности в смысле определения, данного в [\ref{BergLefstrem80}, с.~151].
Поэтому для функции $\psi_{0}$ ограниченность операторов \eqref{2.16}, \eqref{2.18}
влечет за собой ограниченность оператора
\begin{equation}\label{2.20}
I_{Y}\,T\,I_{X}^{-1}:\,L_{2}(U,(\psi_{0}\circ\alpha^{2})\,d\mu)\rightarrow
L_{2}(V,(\psi_{0}\circ\beta^{2})\,d\nu).
\end{equation}
Перейдем от \eqref{2.20} к оператору $T:X_{\psi_{0}}\rightarrow Y_{\psi_{0}}$ с
помощью изометрических изоморфизмов $\psi_{0}(J_{X}):\,X_{\psi_{0}}\leftrightarrow
X_{0}$ и $\psi_{0}(J_{Y}):\,Y_{\psi_{0}}\leftrightarrow Y_{0}$. Приведем их к виду
умножения на функцию:
$$
I_{X}\,\psi_{0}(J_{X})=(\psi_{0}\circ\alpha)\cdot I_{X}:\,
X_{\psi_{0}}\leftrightarrow L_{2}(U,d\mu),
$$
$$
I_{Y}\,\psi_{0}(J_{Y})=(\psi_{0}\circ\beta)\cdot I_{Y}:\,
Y_{\psi_{0}}\leftrightarrow L_{2}(V,d\nu).
$$
Получим изометрические изоморфизмы
$$
I_{X}=(\psi_{0}^{-1}\circ\alpha)\cdot(I_{X}\,\psi_{0}(J_{X})):\,
X_{\psi_{0}}\leftrightarrow L_{2}(U,(\psi^{2}\circ\alpha)\,d\mu),
$$
$$
I_{Y}=(\psi_{0}^{-1}\circ\beta)\cdot(I_{Y}\,\psi_{0}(J_{Y})):\,
Y_{\psi_{0}}\leftrightarrow L_{2}(V,(\psi^{2}\circ\beta)\,d\nu).
$$
Они вместе с \eqref{2.20} влекут за собой ограниченность оператора
$$
T=I_{Y}^{-1}(I_{Y}\,T\,I_{X}^{-1})I_{X}:\,X_{\psi_{0}}\rightarrow Y_{\psi_{0}}.
$$

Таким образом, ввиду равенств \eqref{2.19} имеем:
$$
(T:X_{j}\rightarrow Y_{j},\,j=0,1)\,\Rightarrow\,(T:X_{\psi_{0}}\rightarrow
Y_{\psi_{0}})\,\Rightarrow\,(T:X_{\psi}\rightarrow Y_{\psi}),
$$
где линейные операторы ограничены. Значит, по определению~\ref{def2.2}, функция
$\psi$ --- интерполяционный параметр. Достаточность доказана.

\medskip

\textit{Необходимость.} Предположим, что функция $\psi\in\mathcal{B}$ является
интерполяционным параметром. Покажем, что она псевдовогнута  в окрестности $\infty$.
Будем рассуждать подобно [\ref{Peetre68}; \ref{BergLefstrem80}, с.~152].

Рассмотрим пространство $L_{2}(U,d\mu)$, где $U=\{0,\,1\}$,
$\mu(\{0\})=\mu(\{1\})=1$, и определим на нем линейное отображение $T$ по формуле:
$(Tu)(0)=0$, $(Tu)(1)=u(0)$, где $u\in L_{2}(U,d\mu)$. Выберем произвольные числа
$s,t>1$ и положим $\omega(0):=s^{2}$, $\omega(1):=t^{2}$. Имеем допустимую пару
пространств $X:=[L_{2}(U,d\mu),L_{2}(U,\omega\,d\mu)]$ и ограниченные операторы
$$
T:\,L_{2}(U,d\mu)\rightarrow L_{2}(U,d\mu)\quad\mbox{и}\quad
T:\,L_{2}(U,\omega\,d\mu)\rightarrow L_{2}(U,\omega\,d\mu)
$$
с нормами $1$ и $t/s$ соответственно. Отсюда, поскольку $\psi$ --- интерполяционный
параметр, получаем ограниченный оператор $T:X_{\psi}\rightarrow X_{\psi}$, норма
которого в силу теоремы~\ref{th2.8}, где берем $Y=X$ и $m=1$, удовлетворяет
неравенству
\begin{equation}\label{2.21}
\|T\|_{X_{\psi}\rightarrow X_{\psi}}\leq\,c\,\max\{1,t/s\}.
\end{equation}
Здесь число $c>0$ не зависит от $t,s>1$.

Нетрудно вычислить норму в пространстве $X_{\psi}$. В самом деле, оператор $J$
умножения на функцию $\omega^{1/2}$ является порождающим для пары $X$.
Следовательно, поскольку $\psi(J)$~--- это оператор умножения на функцию
$\psi\circ\omega^{1/2}$, можем записать:
\begin{gather*}
\|u\|_{X_{\psi}}^{2}=\|(\psi\circ\omega^{1/2})\cdot u\|_{L_{2}(U,d\mu)}^{2}=
\psi^{2}(s)\,|u(0)|^{2}+\psi^{2}(t)\,|u(1)|^{2},\\ \quad
\|Tu\|_{X_{\psi}}^{2}=\psi^{2}(t)\,|u(0)|^{2}.
\end{gather*}
Отсюда получаем:
\begin{equation}\label{2.22}
\|T\|_{X_{\psi}\rightarrow X_{\psi}}=\psi(t)/\psi(s).
\end{equation}
Теперь соотношения \eqref{2.21}, \eqref{2.22} влекут за собой неравенство
$$
\psi(t)/\psi(s)\leq c\,\max\{1,t/s\}\quad\mbox{для любых}\quad t,s>1.
$$
Последнее в силу леммы~\ref{lem2.2} равносильно псевдовогнутости функции $\psi$ на
полуоси $(1,\infty)$. Необходимость доказана.

Теорема~\ref{th2.9} доказана.

\medskip

В заключение этого пункта приведем описание всех гильбертовых интерполяционных
пространств для заданной допустимой пары гильбертовых пространств.


\begin{definition}\label{def2.5}
Гильбертово пространство $H$ называется \emph{интерполяционным} для допустимой пары
гильбертовых пространств $[X_{0},X_{1}]$, если:
\begin{itemize}
\item [$\mathrm{(i)}$] справедливы непрерывные вложения $X_{1}\hookrightarrow
H\hookrightarrow X_{0}$;
\item [$\mathrm{(ii)}$] произвольный линейный оператор $T:X_{0}\rightarrow X_{0}$,
ограниченный на каждом из пространств $X_{0}$ и $X_{1}$, является также ограниченным
оператором на пространстве~$H$.
\end{itemize}
\end{definition}


Следующий важный результат принадлежит В.~И.~Овчинникову [\ref{Ovchinnikov84}, с.
511] (теорема 11.4.1).


\begin{proposition}\label{prop2.1}
Пусть $X=[X_{0},X_{1}]$ --- произвольная допустимая пара гильбертовых пространств.
Если гильбертово пространство $H$ является интерполяционным для этой пары, то
существует псевдовогнутая в окрестности $\infty$ функция $\psi\in\mathcal{B}$ такая,
что пространства $H$ и $X_{\psi}$ совпадают с точностью до эквивалентности норм.
\end{proposition}


Из предложения~\ref{prop2.1} и теоремы~\ref{th2.9} немедленно вытекает следствие.

\begin{corollary}\label{cor2.1}
Класс всех гильбертовых пространств, интерполяционных для заданной допустимой пары
$X$ гильбертовых пространств, совпадает (с точностью до эквивалентности норм) с
классом пространств $X_{\psi}$, где $\psi\in\mathcal{B}$~--- произвольная
псевдовогнутая функция в окрестности $\infty$.
\end{corollary}




\markright{\emph \ref{sec2.2}. Правильно меняющиеся функции и их обобщение}

\section[Правильно меняющиеся функции и их обобщение]
{Правильно меняющиеся функции \\ и их обобщение}\label{sec2.2}

\markright{\emph \ref{sec2.2}. Правильно меняющиеся функции и их обобщение}

Правильно меняющиеся функции (и слабо эквивалентные им) играют фундаментальную роль
в нашей работе. С их помощью параметризуются пространства Хермандера. Эти функции
также используются в качестве параметров интерполяции.

\subsection{Правильно меняющиеся функции}\label{sec2.2.1}

Приведем следующее важное для нас определение.


\begin{definition}\label{def2.6}
Положительная функция $\psi$, заданная на вещественной полуоси $[b,\infty)$,
называется \emph{правильно меняющейся} на $\infty$ функцией порядка
$\theta\in\mathbb{R}$, если $\psi$ измерима по Борелю на $[b_{\,0},\infty)$ для
некоторого числа $b_{\,0}\geq b$ и удовлетворяет условию
\begin{equation}\label{2.23}
\lim_{t\rightarrow\,\infty}\;\frac{\psi(\lambda\,t)}{\psi(t)}=
\lambda^{\theta}\quad\mbox{для любого}\quad \lambda>0.
\end{equation}
Правильно меняющаяся  на $\infty$ функция порядка $\theta=0$ называется
\emph{медленно меняющейся} на $\infty$.
\end{definition}


Понятие правильно меняющейся функции введено И.~Караматой в 1930~г. в работе
[\ref{Karamata30a}] (при этом рассматривались непрерывные функции). Правильно
меняющиеся функции примыкают к степенным, хорошо изучены и имеют многочисленные
приложения, в основном, благодаря их особой роли в теоремах тауберова типа (см.
монографии [\ref{Seneta85}, \ref{BinghamGoldieTeugels89}, \ref{Maric00},
\ref{Reshnick87}]).

Обозначим через $\mathrm{SV}$ множество всех медленно меняющихся на $\infty$
функций. Очевидно, что $\psi$~--- правильно меняющаяся на $\infty$ функция порядка
$\theta$ тогда и только тогда, когда $\psi(t)=t^{\theta}\varphi(t)$ при $t\gg1$ для
некоторой функции $\varphi\in\mathrm{SV}$. Поэтому при изучении правильно меняющихся
функций достаточно ограничиться медленно меняющимися функциями.

Приведем наиболее известный (эталонный) пример [\ref{Seneta85}, c.~48~-- 50]
медленно меняющейся функции.


\begin{example}\label{ex2.1}
Пусть заданы $k\in\mathbb{N}$ вещественных чисел $r_{1}$, $r_{2}$,~$\ldots$,
$r_{k}$. Положим
$$
\varphi(t)=(\log t)^{r_{1}}\,(\log\log t)^{r_{2}} \ldots (\log\ldots\log
t)^{r_{k}}\quad\mbox{при}\quad t\gg1.
$$
Тогда $\varphi\in\mathrm{SV}$.
\end{example}


Функции, рассмотренные в этом примере, образуют так называемую
\textit{логарифмическую мультишкалу}, имеющую ряд приложений в теории функциональных
пространств. Другие примеры функций класса $\mathrm{SV}$ будут приведены ниже.

Сформулируем два фундаментальных свойства медленно меняющихся функций. Они доказаны
И.~Караматой [\ref{Karamata30a}, \ref{Karamata33}] для непрерывных функций и
несколько позже рядом авторов для измеримых функций (см. монографии [\ref{Seneta85},
с.~10, 14; \ref{BinghamGoldieTeugels89}, с.~6, 12] и приведенные там ссылки).


\begin{proposition}[теорема о равномерной сходимости]\label{prop2.2}
Пусть $\varphi\in\mathrm{SV}$. Тогда для любого фиксированного отрезка $[a,b]$, где
$0<a<b<\infty$, отношение $\varphi(\lambda t)/\varphi(t)$ стремится к $1$ при
$t\rightarrow\infty$ равномерно относительно $\lambda\in[a,b]$.
\end{proposition}


\begin{proposition}[теорема о представлении]\label{prop2.3}
Пусть $\varphi\in\mathrm{SV}$., тогда
\begin{equation}\label{2.24}
\varphi(t)=\exp\Biggl(\beta(t)+\int\limits_{b}^{\:t}\frac{\alpha(\tau)}{\tau}\,d\tau\Biggr)
\quad\mbox{при}\quad  t\geq b
\end{equation}
для некоторых числа $b>0$, непрерывной функции $\alpha:
[b,\infty)\rightarrow\mathbb{R}$, стремящейся к нулю в $\infty$, и измеримой по
Борелю ограниченной функции $\beta:\nobreak[b,\infty)\rightarrow\mathbb{R}$, имеющей
конечный предел в $\infty$. Обратное также верно: всякая функция вида \eqref{2.24}
принадлежит классу $\mathrm{SV}$.
\end{proposition}


\begin{remark}\label{rem2.2}
Функция вида \eqref{2.24}, где $\beta(t)=\mathrm{const}$, называется
\emph{нормированной} медленно меняющейся функцией. Для нее
$\alpha(\tau)=\tau\varphi'(\tau)/\varphi(\tau)$ $[\ref{BinghamGoldieTeugels89},
c.~15]$.
\end{remark}

Из предложения~\ref{th2.3} вытекает полезное достаточное условие медленного
изменения функции (см., например, [\ref{Seneta85}, с.~15]).


\begin{proposition}\label{prop2.4}
Пусть дифференцируемая функция
$\varphi\nobreak:\nobreak(b,\infty)\rightarrow(0,\infty)$ удовлетворяет условию
$t\varphi\,'(t)/\varphi(t)\rightarrow0$ при $t\rightarrow \infty$. Тогда
$\varphi\in\mathrm{SV}$.
\end{proposition}


Предложение~\ref{prop2.4} приводит к следующим трем интересным примерам медленно
меняющихся функций.


\begin{example}\label{ex2.2}
Пусть $\psi(t):=\exp\varphi(t)$ при $t\gg1$, где функция $\varphi$ взята из
примера~\ref{ex2.1}, в котором полагаем $r_{1}<1$. Тогда $\psi\in\mathrm{SV}$.
\end{example}


\begin{example}\label{ex2.3}
Пусть $\alpha,\beta,\gamma\in\mathbb{R}$, причем $\beta\neq0$ и $0<\gamma<1$.
Положим $\omega(t)=\alpha+\beta\sin\,\ln^{\gamma}t$,
$\varphi(t)=(\ln\,t)^{\,\omega(t)}$ при $t>1$. Тогда $\varphi\in\mathrm{SV}$.
\end{example}


\begin{example}\label{ex2.4}
Пусть $\alpha,\beta,\gamma\in\mathbb{R}$, причем $\alpha\neq0$ и $0<\gamma<\beta<1$.
Положим $r(t)=\alpha(\ln\,t)^{-\beta}\sin\,\ln^{\gamma}t$ и $\varphi(t)=t^{\,r(t)}$
при $t>1$. Тогда $\varphi\in\mathrm{SV}$.
\end{example}


Пример~\ref{ex2.4} интересен тем, что поскольку
$$
\ln\varphi(t)=\alpha(\ln\,t)^{1-\beta}\,\sin\,\ln^{\gamma}t,
$$
то для функции $\varphi(t)$ множество предельных при $t\rightarrow\infty$ точек
заполняет всю полуось $(0,\infty)$.

\subsection{Квазиправильно меняющиеся функции}\label{sec2.2.2} Мы будем
использовать правильно меняющиеся функции порядка $\theta\in(0,1)$ в качестве
интерполяционных параметров. В силу замечания~\ref{rem2.1} свойство быть
интерполяционным параметром наследуется при переходе к (слабо) эквивалентной
функции. Поэтому полезно следующее обобщение понятия правильно меняющейся функции.


\begin{definition}\label{def2.7}
Положительную функцию $\psi$, заданную на вещественной полуоси $[b,\infty)$,
называем \emph{квазиправильно меняющейся} на $\infty$ функцией порядка
$\theta\in\mathbb{R}$, если существуют число $b_{1}\geq b$ и правильно меняющаяся на
$\infty$ функция $\psi_{1}:[b_{1},\infty)\rightarrow (0,\infty)$ порядка $\theta$
такие, что $\psi(t)\asymp\psi_{1}(t)$ при $t\geq b_{1}$. Функцию, квазиправильно
меняющуюся на $\infty$ порядка $\theta=0$, называем \emph{квазимедленно меняющейся}
на $\infty$.
\end{definition}


Обозначим через $\mathrm{QSV}$ множество всех квазимедленно меняющихся на $\infty$
функций. Воспользовавшись предложением~\ref{prop2.3}, получим следующее описание
класса $\mathrm{QSV}$.


\begin{theorem}\label{th2.10}
Класс $\mathrm{QSV}$ состоит из всех функций вида \eqref{2.24}, где число $b>0$,
функция $\alpha: [b,\infty)\rightarrow\mathbb{R}$ непрерывна и стремится к нулю в
$\infty$, а функция $\beta:[b,\infty)\rightarrow\mathbb{R}$ ограничена.
\end{theorem}


\textbf{Доказательство.} По определению~\ref{2.7}, $\varphi\in\mathrm{QSV}$ тогда и
только тогда, когда $\varphi(t)=\omega(t)\varphi_{1}(t)$ при $t\gg1$, где
$\varphi_{1}\in\mathrm{SV}$, а положительная функция $\omega$ ограничена и отделена
от нуля в окрестности $\infty$. Поэтому в силу предложения~\ref{prop2.3},
$\varphi\in\mathrm{QSV}\Leftrightarrow$
$$
\varphi(t)=
\exp\Biggl(\log\omega(t)+\beta(t)+\int\limits_{b}^{t}\frac{\alpha(\tau)}{\tau}\;d\tau\Biggr)
\quad\mbox{при}\;\;t\geq b,
$$
где функции  $\alpha$ и $\beta$ удовлетворяют условию этого предложения, а число
$b\gg1$. Отсюда следует теорема~\ref{th2.10}, поскольку функция $\log\omega+\beta$
ограничена на полуоси $[b,\infty)$. Теорема~\ref{th2.10} доказана.

\medskip

Далее фундаментальную роль будет играть следующее интерполяционное свойство
квазиправильно меняющихся функций.


\begin{theorem}\label{th2.11}
Пусть функция $\psi\in\mathcal{B}$ квазиправильно меняющаяся на $\infty$ порядка
$\theta$, где $0<\theta<1$. Тогда $\psi$ является интерполяционным параметром.
\end{theorem}


\textbf{Доказательство.} Запишем $\psi(t)=t^{\theta}\varphi(t)$ при $t>0$, где
$\varphi\in\mathrm{QSV}$. На основании теоремы~\ref{th2.10} представим функцию
$\varphi$ в виде \eqref{2.24}, где функции $\alpha$ и $\beta$ удовлетворяют условию
этой теоремы. Положим $\varepsilon:=\min\{\theta,1-\theta\}>0$ и выберем число
$b_{\varepsilon}\geq b$ такое, что $|\alpha(t)|<\varepsilon$ при $t>
b_{\varepsilon}$. Для произвольных $t,s>b_{\varepsilon}$ в силу \eqref{2.24} имеем:
\begin{gather*}
\frac{\varphi(t)}{\varphi(s)}=\exp\Biggl(\beta(t)-\beta(s)+
\int\limits_{s}^{t}\frac{\alpha(\tau)}{\tau}\;d\tau\Biggr)\leq \\
c\exp\Bigl|\int\limits_{s}^{t}\frac{\varepsilon}{\tau}\;d\tau\Bigr|=
c\max\left\{(t/s)^{\varepsilon},(s/t)^{\varepsilon}\right\}.
\end{gather*}
Здесь число $c>0$ не зависит от $t$ и $s$, так как функция $\beta$ ограничена.
Отсюда, поскольку $0\leq\theta\pm\varepsilon\leq1$, получаем:
\begin{gather*}
\psi(t)/\psi(s)=(t^{\theta}\varphi(t))/(s^{\theta}\varphi(s))\leq\\
c\max\left\{(t/s)^{\theta+\varepsilon},(t/s)^{\theta-\varepsilon}\right\}\leq
c\max\{1,t/s\},\quad t,s>b_{\varepsilon}.
\end{gather*}
Следовательно, в силу леммы~\ref{lem2.2} функция $\psi\in\mathcal{B}$ псевдовогнута
в окрестности $\infty$. Это согласно теореме~\ref{th2.9} равносильно тому, что
$\psi$ --- интерполяционный параметр. Теорема~\ref{th2.11} доказана.


\begin{remark}\label{rem2.3}
Прямое доказательство теоремы~\ref{th2.11} (без использования теоремы~\ref{th2.9})
приведено авторами [\ref{06UMJ2}] (п.~2).
\end{remark}


\begin{remark}\label{rem2.4}
Теорема~\ref{th2.11} не верна в предельных случаях $\theta=0$ и $\theta=\nobreak1$
даже, если дополнительно предположить, что при $t\rightarrow\infty$ выполняется
условие $\psi(t)\rightarrow\nobreak\infty$ в случае $\theta=0$ либо
$\psi(t)/t\rightarrow0$ в случае $\theta=1$. Приведем соответствующие примеры.
\end{remark}


\begin{example}\label{ex2.5}
Случай $\theta=0$. Пусть $h(t):=(\ln t)^{-1/2}\sin\ln^{1/4}t$ при $t>1$. Определим
функцию $\psi$ следующим образом: $\psi(t):=t^{h(t)}+\ln t$ при $t\geq3$ и
$\psi(t):=1$ при $0<t<3$. Заметим, что $\psi\in\mathcal{B}$ и
$\psi(t)\rightarrow\infty$ при $t\rightarrow\infty$. Непосредственно устанавливается
сходимость $t\psi'(t)/\psi(t)\rightarrow0$ при $t\rightarrow\infty$. Отсюда в силу
предложения~\ref{prop2.4} следует, что функция $\psi$ медленно меняющаяся на
$\infty$. Покажем, что она не является псевдовогнутой в окрестности $\infty$ и,
следовательно, в силу теоремы~\ref{th2.9} не является интерполяционным параметром.
Рассмотрим последовательности чисел $t_{k}:=\exp((2\pi k+\pi/2)^{4})$ и
$s_{k}:=\exp((2\pi k+\pi)^{4})$, где $k\in\mathbb{N}$. Вычисляем: $h(t_{k})=(2\pi
k+\pi/2)^{-2}$ и $h(s_{k})=0$, откуда
$$
\ln\psi(t_{k})\geq h(t_{k})\ln t_{k}=(2\pi k+\pi/2)^{2}
$$
и $\psi(s_{k})=1+(2\pi k+\pi)^{4}$. Следовательно,
$$
\frac{\psi(t_{k})}{\psi(s_{k})}\geq\frac{\exp((2\pi k+\pi/2)^{2})}{(1+(2\pi
k+\pi)^{4})}\rightarrow\infty
$$
при $k\rightarrow \infty$. Но $t_{k}<s_{k}$, значит, в силу леммы \ref{lem2.2}
функция $\psi$ не может быть псевдовогнутой в окрестности $\infty$.
\end{example}


\begin{example}\label{ex2.6}
Случай $\theta=1$. Воспользовавшись функцией $\psi$ из примера \ref{ex2.5}, положим
$\psi_{1}(t):=t/\psi(t)$ при $t>0$. По определению, функция $\psi_{1}$ правильно
меняющаяся на $\infty$ порядка $\theta=1$. Заметим, что $\psi_{1}\in\mathcal{B}$ и
$\psi_{1}(t)/t=1/\psi(t)\rightarrow0$ при $t\rightarrow\infty$. В силу
теоремы~\ref{th2.4} функция $\psi_{1}$ не является интерполяционным параметром,
поскольку в противном случае интерполяционным параметром была бы функция
$\psi(t)=t/\psi_{1}(t)$, что, как было показано в предыдущем примере, невозможно.
\end{example}


Далее будут использованы следующие свойства  класса $\mathrm{QSV}$.


\begin{theorem}\label{th2.12} Пусть $\varphi,\chi\in\mathrm{QSV}$.
Справедливы следующие утверждения.
\begin{itemize}
\item [$\mathrm{(i)}$] Существует положительная функция
$\varphi_{1}\in C^{\infty}((0;\infty))\cap\mathrm{SV}$ такая, что
$\varphi(t)\asymp\varphi_{1}(t)$ при $t\gg1$.
\item [$\mathrm{(ii)}$] Для произвольного числа $\theta>0$ справедливо
$t^{-\theta}\varphi(t)\rightarrow0$ и $t^{\theta}\varphi(t)\rightarrow\infty$ при
$t\rightarrow\infty$.
\item [$\mathrm{(iii)}$] Функции $\varphi+\chi$, $\varphi\,\chi$, $\varphi/\chi$ и
$\varphi^{\sigma}$, где $\sigma\in\mathbb{R}$, принадлежат классу $\mathrm{QSV}$
\item [$\mathrm{(iv)}$] Если число $\theta\geq0$, причем в случае $\theta=0$ справедливо
$\varphi(t)\rightarrow\infty$ при $t\rightarrow\infty$, то сложная функция
$\chi(t^{\theta}\varphi(t))$ аргумента $t$ принадлежит классу $\mathrm{QSV}$.
\end{itemize}
\end{theorem}


\textbf{Доказательство.} Если дополнительно предположить, что
$\varphi,\chi\in\mathrm{SV}$, то получим известные [\ref{Seneta85}, c. 23 -- 25]
свойства медленно меняющихся функций. (При этом в пункте (i) будет даже сильная
эквивалентность $\varphi(t)\sim\varphi_{1}(t)$ при $t\rightarrow\infty$.) Отсюда
немедленно вытекают пункты (i), (ii), (iii) для функций
$\varphi,\chi\in\mathrm{QSV}$.

Докажем пункт (iv). Пусть число $\lambda>0$. Так как $\varphi\in\mathrm{QSV}$, то
величина $\varphi(\lambda t)/\varphi(t)$ ограничена и отделена от нуля при $t\gg1$.
Поэтому в силу предложения~\ref{prop2.2}, для произвольной положительной функции
$\chi_{1}\in\mathrm{SV}$, удовлетворяющей условию $\chi_{1}(\tau)\asymp\chi(\tau)$,
$\tau\gg1$, мы имеем
\begin{gather*}
\chi_{1}\left((\lambda t)^{\theta}\varphi(\lambda t)\right)\big/
\chi_{1}\left(t^{\theta}\varphi(t)\right)=\\
\chi_{1}\left(\frac{\lambda^{\theta}\varphi(\lambda
t)}{\varphi(t)}\:t^{\theta}\varphi(t)\right)\Big/
\chi_{1}\left(t^{\theta}\varphi(t)\right)\rightarrow1\quad\mbox{при}\;\;t\rightarrow\infty.
\end{gather*}
Здесь мы воспользовались тем, что $t^{\theta}\varphi(t)\rightarrow\infty$ при
$t\rightarrow\infty$. Следовательно, функция $\chi_{1}(t^{\theta}\varphi(t))$
медленно меняющаяся на $\infty$. Но
$\chi(t^{\theta}\varphi(t))\asymp\chi_{1}(t^{\theta}\varphi(t))$ при $t\gg1$.
Значит, функция $\chi(t^{\theta}\varphi(t))$ принадлежат классу $\mathrm{QSV}$.
Пункт (iv), а с ним и теорема~\ref{th2.12} доказаны.

\medskip

Нам понадобится также следующее обобщение утверждения (iv) теоремы~\ref{th2.12} на
случай, когда $\theta=0$ и $\varphi(t)\nrightarrow\infty$ при $t\rightarrow\infty$.


\begin{theorem}\label{th2.13}Пусть заданы функции
$\varphi:[b,\infty)\rightarrow(0,\infty)$ и $\chi:(0,\infty)\rightarrow(0,\infty)$
класса $\mathrm{QSV}$. Предположим, что функция $1/\varphi$ ограничена на полуоси
$[b,\infty)$, а функции $\chi$ и $1/\chi$ ограничены на каждом отрезке
$[a_{1},a_{2}]$, где $0<a_{1}<a_{2}<\infty$. Тогда сложная функция
$\chi(\varphi(t))$ аргумента $t$ принадлежат классу $\mathrm{QSV}$.
\end{theorem}


\textbf{Доказательство.} По условию $\varphi(t)\geq r$ при $t\geq b$ для некоторого
числа $r>0$. Согласно теореме~\ref{th2.12} (i) существует число
$b_{1}\geq\max\{b,r\}$ и непрерывные медленно меняющиеся на $\infty$ функции
$\varphi_{1},\chi_{1}:[b_{1},\infty)\rightarrow(0,\infty)$ такие, что
$\varphi(t)\asymp\varphi_{1}(t)$ и $\chi(t)\asymp\chi_{1}(t)$ при $t\geq b_{1}$. При
этом можно считать, что $\varphi_{1}(t)\geq r$ при $t\geq b_{1}$. Доопределим:
$\varphi_{1}(t):=\varphi_{1}(b_{1})$ и $\chi_{1}(t):=\chi_{1}(b_{1})$ при
$0<t<b_{1}$. Тогда $\varphi_{1},\chi_{1}\in C((0,\infty))$ и $\varphi_{1}(t)\geq r$
при $t>0$, а также $\chi(t)\asymp\chi_{1}(t)$ при $t\geq r$ ввиду условия
доказываемой теоремы. Покажем, что $\chi(\varphi(t))\asymp\chi_{1}(\varphi_{1}(t))$
при $t\geq b_{1}$ и функция $\chi_{1}(\varphi_{1}(t))$ медленно меняющаяся на
$\infty$. Отсюда последует требуемое свойство $\chi\circ\varphi\in\mathrm{QSV}$.

Пусть $t\geq b_{1}$. Поскольку $c^{-1}\leq\varphi_{1}(t)/\varphi(t)\leq c$ для
некоторого числа $c\geq1$, в силу предложения~\ref{prop2.2} имеем для функции
$\chi_{1}\in\mathrm{SV}$
$$
\chi_{1}(\varphi_{1}(t))\,/\,\chi_{1}(\varphi(t))=
\chi_{1}\left(\frac{\varphi_{1}(t)}{\varphi(t)}\,\varphi(t)\right)
\Big/\chi_{1}(\varphi(t))\rightarrow1
$$
при $\varphi(t)\rightarrow\infty$. Поэтому существует число $\varrho\geq r$ такое,
что $\chi_{1}(\varphi(t))\asymp\chi_{1}(\varphi_{1}(t))$ при
$\varphi(t)\geq\varrho$. Кроме того,
$$
\chi_{1}(\varphi(t))\asymp\chi(\varphi(t))\asymp1\asymp\chi(\varphi_{1}(t))
\asymp\chi_{1}(\varphi_{1}(t))
$$
при $r\leq\varphi(t)\leq\varrho$. Следовательно,
$$
\chi(\varphi(t))\asymp\chi_{1}(\varphi(t))\asymp\chi_{1}(\varphi_{1}(t))
\quad\mbox{при}\;\;t\geq b_{1}.
$$

Фиксируем число $\lambda>0$. Докажем, что
$$
\chi_{1}(\varphi_{1}(\lambda t))/\chi_{1}(\varphi_{1}(t))\rightarrow1
\quad\mbox{при}\;\;t\rightarrow\infty.
$$
Выберем произвольно число $\varepsilon>0$. Так как $\varphi_{1}\in\mathrm{SV}$, то
$\beta_{\lambda}(t):=\varphi_{1}(\lambda t)/\varphi_{1}(t)\rightarrow1$ при
$t\rightarrow\infty$. В частности, $1/2\leq\beta_{\lambda}(t)\leq2$ при $t\geq
t_{\lambda}>0$. Отсюда в силу предложения~\ref{prop2.2} для функции
$\chi_{1}\in\mathrm{SV}$ найдется число $k=k(\varepsilon)>r$ такое, что
\begin{gather}\notag
|\chi_{1}(\varphi_{1}(\lambda t))/\chi_{1}(\varphi_{1}(t))-1|=\\
|\chi_{1}(\beta_{\lambda}(t)\varphi_{1}(t))/\chi_{1}(\varphi_{1}(t))-1|
<\varepsilon\label{2.25}
\end{gather}
при $t\geq t_{\lambda}$ и $\varphi_{1}(t)>k$. Кроме того, поскольку функция
$\chi_{1}>0$ равномерно непрерывна на отрезке $[r,k+1]$, существует число
$m=m(\varepsilon)>0$ такое, что
\begin{gather*}
|\chi_{1}(\varphi_{1}(\lambda t))-\chi_{1}(\varphi_{1}(t))|=
|\chi_{1}(\beta_{\lambda}(t)\varphi_{1}(t))-\chi_{1}(\varphi_{1}(t))|< \\
\varepsilon\min\{\chi_{1}(\tau):\,r\leq\tau\leq k\}\quad\mbox{при}\quad t\geq
m,\quad r\leq\varphi_{1}(t)\leq k.
\end{gather*}
Значит,
\begin{equation}\label{2.26}
|\chi_{1}(\varphi_{1}(\lambda t))/\chi_{1}(\varphi_{1}(t))-1|=
\frac{|\chi_{1}(\varphi_{1}(\lambda
t))-\chi_{1}(\varphi_{1}(t))|}{\chi_{1}(\varphi_{1}(t))}<\varepsilon
\end{equation}
при $t\geq m$ и $r\leq\varphi_{1}(t)\leq k$. Теперь соотношения \eqref{2.25} и
\eqref{2.26} влекут за собой неравенство
$$
\left|\frac{\chi_{1}(\varphi_{1}(\lambda
t))}{\chi_{1}(\varphi_{1}(t))}-1\right|<\varepsilon\quad\mbox{при}\quad
t\geq\max\{t_{\lambda},m\}.
$$
Последнее ввиду произвольности выбора числа $\varepsilon>0$ означает, что
$\chi_{1}(\varphi_{1}(\lambda t))/\chi_{1}(\varphi_{1}(t))\rightarrow1$ при
$t\rightarrow\infty$. Следовательно, функция $\chi_{1}(\varphi_{1}(t))$ медленно
меняющаяся на $\infty$. Теорема~\ref{th2.13} доказана.

\subsection{Вспомогательные результаты}\label{sec2.2.3}

В этом пункте мы докажем два вспомогательных утверждения о свойствах класса
$\mathrm{QSV}$, которые нам понадобятся впоследствии.


\begin{lemma}\label{lem2.3}
Пусть функция $\varphi\in\mathrm{QSV}$ положительна на полуоси $[1;\infty]$ и
ограничена вместе с функцией $1/\varphi$ на каждом отрезке $[1;b]$, где
$1<b<\infty$. Тогда для произвольного числа $\varepsilon>0$ существует число
$c(\varepsilon)>0$ такое, что для любых $t\geq1$ и $s\geq1$ выполняется неравенство
\begin{equation}\label{2.27}
\varphi(t)/\varphi(s)\leq c(\varepsilon)\,(1+|t-s|)^{\varepsilon}.
\end{equation}
\end{lemma}


\textbf{Доказательство} достаточно провести для случая, когда
$0\nobreak<\nobreak\varepsilon<1$. В~силу теоремы~\ref{th2.10} представим $\varphi$
по формуле \eqref{2.24}. Поскольку там $\alpha(\tau)\rightarrow0$ при
$\tau\rightarrow\infty$, то найдется число $b_{\varepsilon}\geq1$ такое, что
$|\alpha(\tau)|\leq\varepsilon$  при $\tau\geq b_{\varepsilon}$. Возьмем
произвольные числа $t\geq1$ и $s\geq1$. Далее в доказательстве будем обозначать
через $c_{1}$, $c_{2}$, $c_{3}$ конечные положительные постоянные, не зависящие от
$t$ и $s$. Докажем \eqref{2.27} отдельно для четырех возможных случаев расположения
чисел $t$ и $s$ относительно $b_{\varepsilon}$.

Первый случай: $t\geq b_{\varepsilon},\,s\geq b_{\varepsilon}$. В силу \eqref{2.24}
имеем:
\begin{gather*}
\frac{\varphi(t)}{\varphi(s)}=\exp\int\limits_{s}^{\,t}\frac{\alpha(\tau)}{\tau}\,d\tau
\leq\exp\Bigl|\int\limits_{s}^{\,t}\frac{\varepsilon}{\tau}\;d\tau\Bigr|=\\
\max\left\{\Bigl(\frac{t}{s}\Bigr)^{\,\varepsilon}\,,
\Bigl(\frac{s}{t}\Bigr)^{\,\varepsilon}\right\}=
\max\left\{\Bigl(1+\frac{t-s}{s}\Bigr)^{\,\varepsilon}\,,
\Bigl(1+\frac{s-t}{t}\Bigr)^{\,\varepsilon}\right\}\leq \\
(1+|t-s|)^{\,\varepsilon}.
\end{gather*}

Второй случай: $t\geq b_{\varepsilon}$ и $1\leq s\leq b_{\varepsilon}$. В силу
теоремы~\ref{th2.12} (ii) и условия этой леммы имеем: $\varphi(t)\leq
c_{1}\,t^{\,\varepsilon}$ и $1/\varphi(s)\leq c_{1}$. Отсюда
\begin{gather*}
\frac{\varphi(t)}{\varphi(s)}\leq
c_{1}^{\,2}\,t^{\varepsilon}=c_{1}^{\,2}\,(s+(t-s))^{\,\varepsilon}\leq\\
c_{1}^{\,2}\,(b_{\varepsilon}+|t-s|)^{\,\varepsilon}\leq
c_{1}^{\,2}\,b_{\varepsilon}(1+|t-s|)^{\,\varepsilon}.
\end{gather*}

Третий случай: $1\leq t\leq b_{\varepsilon}$ и $s\geq b_{\varepsilon}$. Аналогично
предыдущему случаю получаем: $1/\varphi(s)\leq c_{2}\,s^{\,\varepsilon}$,
$\varphi(t)\leq c_{2}$ и, следовательно,
\begin{gather*}
\frac{\varphi(t)}{\varphi(s)}\leq
c_{2}^{\,2}\,s^{\,\varepsilon}=c_{2}^{\,2}\,(t+(s-t))^{\,\varepsilon} \leq \\
c_{2}^{\,2}\,(b_{\varepsilon}+|s-t|)^{\,\varepsilon}\leq
c_{2}^{\,2}\,b_{\varepsilon}(1+|t-s|)^{\,\varepsilon}.
\end{gather*}

Четвертый случай: $1\leq t\leq b_{\varepsilon}$ и $1\leq s\leq b_{\varepsilon}$. Он
тривиален:  $\varphi(t)/\varphi(s)\leq c_{3}\leq c_{3}(1+|t-s|)^{\varepsilon}$.

Таким образом, неравенство \eqref{2.27} справедливо для произвольных $t\geq1$ и
$s\geq1$. Лемма~\ref{lem2.3} доказана.


\begin{lemma}\label{lem2.4}
Пусть функция $\psi_{1}\in\mathrm{QSV}$ положительна, непрерывна на полуоси
$[1;\infty)$ и удовлетворяет условию
\begin{equation}\label{2.28}
I_{1}:=\int_{1}\limits^{\infty}\frac{d\,t}{t\,\psi_{1}(t)}<\infty.
\end{equation}
Тогда существует такая функция $\psi_{0}\in\mathrm{SV}$, положительная и непрерывная
на $[1;\infty)$, что $\psi_{0}(t)/\psi_{1}(t)\rightarrow0$ при $t\rightarrow\infty$
и
\begin{equation}\label{2.29}
\int\limits_{1}^{\infty}\frac{d\,t}{t\,\psi_{0}(t)}<\infty.
\end{equation}
\end{lemma}


\textbf{Доказательство.} Ввиду теоремы~\ref{th2.12} (i) можно без потери общности
считать, что $\psi_{1}\in C((0,\infty))\cap\mathrm{SV}$. Положим
\begin{equation}\label{2.30}
\varphi(t):=\int\limits_{t}^{\infty}\frac{d\,t}{t\,\psi_{1}(t)}\quad\mbox{и}\quad\psi_{0}(t)
:=\psi_{1}(t)\,\sqrt{\varphi(t)}\quad\mbox{при}\;\;t\geq1.
\end{equation}
По условию, функция $\varphi$ конечная положительная на $[1;\infty)$, причем
$\varphi(t)\rightarrow0$ при $t\rightarrow\infty$. Кроме того, $\varphi$ имеет
непрерывную производную $\varphi\,'(t)=-(t\,\psi_{1}(t))^{-1}$,  $t\geq1$. Отсюда
ввиду включения $\psi_{1}\in\mathrm{SV}$ вычисляем по правилу Лопиталя при любом
$\lambda>0$:
$$
\lim_{t\rightarrow\infty}\frac{\varphi(\lambda\, t)}{\varphi(t)}=
\lim_{t\rightarrow\infty}\frac{\lambda\,(\lambda\,t\,\psi_{1}(\lambda\,
t))^{-1}}{(t\,\psi_{1}(t))^{-1}}=
\lim_{t\rightarrow\infty}\frac{\psi_{1}(t)}{\psi_{1}(\lambda\, t)}=1.
$$
Значит, $\varphi\in\mathrm{SV}$. Обратимся теперь к функции $\psi_{\,0}$. Она
положительна и непрерывна на $[1;\infty)$. Поскольку
$\psi_{1},\varphi\in\mathrm{SV}$, то $\psi_{0}\in\mathrm{SV}$ по определению
является медленно меняющейся функцией на $\infty$. Кроме того,
$\psi_{\,0}(t)/\psi_{1}(t)=\sqrt{\varphi(t)}\rightarrow0$ при $t\rightarrow\infty$,
и
$$
\int\limits_{1}^{\infty}\frac{d\,t}{t\,\psi_{0}(t)}=
\int\limits_{1}^{\infty}\frac{d\,t}{t\,\psi_{1}(t)\,\sqrt{\varphi(t)}}=
-\int\limits_{1}^{\,\infty}\frac{d\,\varphi(t)}{\sqrt{\varphi(t)}}=
-\int\limits_{I_{1}}^{0}\frac{d\,\tau}{\sqrt{\tau}}<\infty.
$$
Таким образом, $\psi_{0}$ удовлетворяет условиям леммы. Лемма~\ref{lem2.4} доказана.





\markright{\emph \ref{sec2.3}. Пространства Хермандера и уточненная шкала}

\section[Пространства Хермандера и уточненная шкала]
{Пространства Хермандера \\ и уточненная шкала}\label{sec2.3}

\markright{\emph \ref{sec2.3}. Пространства Хермандера и уточненная шкала}

В этом параграфе мы рассмотрим важный класс пространств Хермандера,
параметризованных с помощью правильно меняющихся функций. Это основной класс
пространств, в котором будут изучены эллиптические операторы. Мы называем его
уточненной шкалой.

\subsection{Предварительные сведения и обозначения}\label{sec2.3.0} 

Мы используем следующие общепринятые обозначения комплексных линейных топологических
пространств основных функций и распределений (обобщенных функций), заданных в
евклидовом пространстве $\mathbb{R}^{n}$ (см., например, [\ref{GelfandShilov58},
\ref{GelfandShilov59}, \ref{Schwartz50}, \ref{Schwartz51}]):

\begin{itemize}
\item [] $C^{\infty}_{0}(\mathbb{R}^{n}):=\mathcal{D}(\mathbb{R}^{n})$~--- пространство всех
бесконечно дифференцируемых функций $u:\mathbb{R}^{n}\rightarrow\mathbb{C}$ с
компактными носителями;

\item [] $\mathcal{S}(\mathbb{R}^{n})$~--- пространство Шварца всех
бесконечно дифференцируемых функций $u:\mathbb{R}^{n}\rightarrow\mathbb{C}$ c
конечными нормами
$$
\max\{\,(1+|x|)^{m}\,|\partial^{\alpha}_{x}u(x)|:\,x\in\mathbb{R}^{n},\;|\alpha|\leq
m\},\;\; m=0,1,2,\ldots;
$$
здесь и далее $\alpha=(\alpha_{1},\ldots,\alpha_{n})$~--- мультииндекс (вектор с
целыми неотрицательными координатами),
$|\alpha|:=\alpha_{1}+\ldots+\nobreak\alpha_{n}$,
$\partial^{\alpha}_{x}u(x):=\partial^{|\alpha|}u(x)/\partial
x_{1}^{\alpha_{1}}\ldots\partial x_{n}^{\alpha_{n}}$~--- частная производная,
соответствующая мультииндексу $\alpha$;

\item [] $\mathcal{S}'(\mathbb{R}^{n})$~--- двойственное к $\mathcal{S}(\mathbb{R}^{n})$
пространство распределений медленного роста, заданных в $\mathbb{R}^{n}$;

\item [] $\mathcal{D}'(\mathbb{R}^{n})$~--- двойственное к $\mathcal{D}(\mathbb{R}^{n})$
пространство всех распределений, заданных в $\mathbb{R}^{n}$.
\end{itemize}

С точки зрения приложений нам удобно трактовать распределения как антилинейные
функционалы. Поэтому в качестве элементов двойственных пространств
$\mathcal{S}'(\mathbb{R}^{n})$ и $\mathcal{D}'(\mathbb{R}^{n})$ мы берем
антилинейные непрерывные функционалы на $\mathcal{S}(\mathbb{R}^{n})$ и
$\mathcal{D}(\mathbb{R}^{n})$ соответственно. Взаимная двойственность пространств
основных функций и распределений в $\mathbb{R}^{n}$ рассматривается относительно
расширения по непрерывности полуторалинейной формы
$$
(u,v)_{\mathbb{R}^{n}}:=\int\limits_{\mathbb{R}^{n}}u(x)\,\overline{v(x)}\,dx.
$$
Это расширение мы обозначаем также через $(u,v)_{\mathbb{R}^{n}}$; оно равно
значению распределения $u$ на основной функции $v$.

Всякая локально суммируемая по Лебегу функция
$u:\mathbb{R}^{n}\rightarrow\mathbb{C}$ отождествляется с распределением~---
антилинейным функционалом $v\mapsto(u,v)_{\mathbb{R}^{n}}$, заданном на функциях
$v\in\mathcal{D}(\mathbb{R}^{n})$; такое распределение называется регулярным. В этом
смысле справедливы непрерывные и плотные вложения
$$
\mathcal{D}(\mathbb{R}^{n})\hookrightarrow\mathcal{S}(\mathbb{R}^{n})\hookrightarrow
\mathcal{S}'(\mathbb{R}^{n})\hookrightarrow\mathcal{D}'(\mathbb{R}^{n}).
$$

Мы обозначаем через $\mathcal{F}u$ или коротко $\widehat{u}$ преобразование Фурье
распределения $u\in\mathcal{S}'(\mathbb{R}^{n})$. Если
$u\in\mathcal{S}(\mathbb{R}^{n})$, то нами используется следующая формула
преобразования Фурье:
$$
(\mathcal{F}u)(\xi)=\widehat{u}(\xi):=
\int\limits_{\mathbb{R}^{n}}\,e^{ix\cdot\xi}\,u(x)\,dx,\quad\xi\in\mathbb{R}^{n}.
$$
Здесь, как обычно, $i$~--- мнимая единица, а
$x\cdot\xi:=x_{1}\xi_{1}+\ldots+x_{n}\xi_{n}$~--- скалярное произведение векторов
$x,\xi\in\mathbb{R}^{n}$. Преобразование Фурье является гомеоморфизмом пространства
$\mathcal{S}(\mathbb{R}^{n})$ на себя. Обратное преобразование Фурье вычисляется по
формуле
$$
u(x)=(\mathcal{F}^{-1}\widehat{u})(x)=(2\pi)^{-n}
\int\limits_{\mathbb{R}^{n}}\,e^{-ix\cdot\xi}\,\widehat{u}(\xi)\,d\xi,\quad
x\in\mathbb{R}^{n}.
$$
Преобразование Фурье продолжается по непрерывности до гомеоморфизма пространства
$\mathcal{S}'(\mathbb{R}^{n})$ на себя, при этом сохраняется равенство Парсеваля
$$
(\widehat{u},\widehat{v})_{\mathbb{R}^{n}}=(2\pi)^{n}(u,v)_{\mathbb{R}^{n}}\quad\mbox{для
всех}\;\;u\in\mathcal{S}'(\mathbb{R}^{n}),\;\;v\in\mathcal{S}(\mathbb{R}^{n}).
$$
Последнее может быть взято в качестве определения преобразования Фурье распределения
$u\in\mathcal{S}'(\mathbb{R}^{n})$.

Напомним еще некоторые стандартные обозначения функциональных пространств.
Пространство $L_{p}(\mathbb{R}^{n})$ состоит из всех измеримых по Лебегу функций
$u:\mathbb{R}^{n}\rightarrow\mathbb{C}$ таких, что
\begin{gather*}
\|u\|_{L_{p}(\mathbb{R}^{n})}^{p}:=\int\limits_{\mathbb{R}^{n}}\,|u(x)|^{p}\,dx<\infty
\quad\mbox{при}\;\;1<p<\infty,\\
\|u\|_{L_{\infty}(\mathbb{R}^{n})}:=\mathrm{ess\,sup}\{\,|u(x)|:\,x\in\mathbb{R}^{n}\}<\infty
\quad\mbox{при}\;\;p=\infty.
\end{gather*}
Пространство $L_{p}(\mathbb{R}^{n})$, где $1\leq p\leq\infty$, банахово относительно
нормы $\|u\|_{L_{p}(\mathbb{R}^{n})}$, при этом, естественно, отождествляют функции,
совпадающие в $\mathbb{R}^{n}$ почти всюду. В случае $p=2$ это пространство
становится гильбертовым: норма в нем порождена скалярным произведением
$(u,v)_{\mathbb{R}^{n}}$. Преобразование Фурье, умноженное на $(2\pi)^{-n/2}$,
является изометрическим изоморфизмом пространства $L_{2}(\mathbb{R}^{n})$ на себя.

Как обычно, $C^{k}(\mathbb{R}^{n})$, где целое $k\geq0$, обозначает пространство
всех функций $u:\mathbb{R}^{n}\rightarrow\mathbb{C}$, имеющих непрерывные частные
производные до порядка $k$ включительно. Подпространство
$C^{k}_{\mathrm{b}}(\mathbb{R}^{n})$ состоит из таких функций $u\in
C^{k}(\mathbb{R}^{n})$, у которых все частные производные до порядка $k$
включительно ограничены в $\mathbb{R}^{n}$. Пространство
$C^{k}_{\mathrm{b}}(\mathbb{R}^{n})$ банахово относительно нормы
$$
\sup\{\,|\partial^{\alpha}_{x}u(x)|:\,x\in\mathbb{R}^{n},\;|\alpha|\leq k\}.
$$
Нам понадобятся также пространства
$$
C^{\infty}(\mathbb{R}^{n}):=\bigcap_{k\geq0}C^{k}(\mathbb{R}^{n}),\quad
C^{\infty}_{\mathrm{b}}(\mathbb{R}^{n}):=\bigcap_{k\geq0}C^{k}_{\mathrm{b}}(\mathbb{R}^{n}).
$$

\subsection{Пространства Хермандера}\label{sec2.3.1}

Приведем определение функциональных пространств, введенных и изученных
Л.~Хермандером в монографии [\ref{Hermander65}, с.~54] (см. также
[\ref{Hermander86}, с.~13]). Эти пространства состоят из распределений, заданных в
$\mathbb{R}^{n}$, где $n\in\mathbb{N}$, и обозначаются через
$B_{p,\mu}(\mathbb{R}^{n})$. Здесь и далее в этом пункте число $p$ удовлетворяет
неравенству $1\leq p\leq\infty$, а функция $\mu=\mu(\xi)$ переменной
$\xi\in\mathbb{R}^{n}$ положительна, непрерывна и является весовой в смысле
следующего определения.

\begin{definition}\label{def2.8}
Функцию $\mu:\mathbb{R}^{n}\rightarrow(0,\infty)$ называют \emph{весовой}, если
существуют числа $c\geq1$ и $l>0$ такие, что
\begin{equation}\label{2.31}
\frac{\mu(\xi)}{\mu(\eta)}\leq c\,(1+|\xi-\eta|)^{l}\quad\mbox{для
любых}\quad\xi,\eta\in\mathbb{R}^{n}.
\end{equation}
\end{definition}


Следующее определение является базисным для нас.

\begin{definition}\label{def2.9}
\emph{Пространством Хермандера} $B_{p,\mu}(\mathbb{R}^{n})$ называют линейное
пространство всех распределений $u\in\mathcal{S}'(\mathbb{R}^{n})$ таких, что
преобразование Фурье $\widehat{u}$ является локально суммируемой по Лебегу в
$\mathbb{R}^{n}$ функцией, удовлетворяющей включению $\mu\,\widehat{u}\in
L_{p}(\mathbb{R}^{n})$. В пространстве $B_{p,\mu}(\mathbb{R}^{n})$ вводится норма по
формуле
$\|u\|_{B_{p,\mu}(\mathbb{R}^{n})}:=\|\mu\,\widehat{u}\|_{L_{p}(\mathbb{R}^{n})}$.
\end{definition}


Пространство $B_{p,\mu}(\mathbb{R}^{n})$ полно относительно введенной в нем нормы и
непрерывно вложено в $\mathcal{S}'(\mathbb{R}^{n})$. Если $1\leq p<\infty$, то это
пространство сепарабельно, и множество $C^{\infty}_{0}(\mathbb{R}^{n})$ плотно в нем
(см. [\ref{Hermander65}, с.~54], теорема 2.2.1). Особый интерес представляет случай
$p=2$, когда пространство $B_{p,\mu}(\mathbb{R}^{n})$ гильбертово.


\begin{remark}\label{rem2.5}
Л.~Хермандер [\ref{Hermander65}, с.~51] изначально предполагает, что функция $\mu$
удовлетворяет более сильному условию, чем \eqref{2.31}, а именно: существуют
положительные числа $c$ и $l$ такие, что
\begin{equation}\label{2.32}
\frac{\mu(\xi)}{\mu(\eta)}\leq(1+c\,|\xi-\eta|)^{l}\quad\mbox{для
любых}\quad\xi,\eta\in\mathbb{R}^{n}.
\end{equation}
Однако потом замечает [\ref{Hermander65}, с.~54], что множества функций,
удовлетворяющие условиям \eqref{2.31} и \eqref{2.32} приводят к одному и тому же
классу пространств $B_{p,\mu}(\mathbb{R}^{n})$.
\end{remark}


\begin{remark}\label{rem2.6}
Л.~Р.~Волевич и Б.~П.~Панеях [\ref{VolevichPaneah65}] ввели и изучили пространства
$H^{\mu}_{p}(\mathbb{R}^{n})$, $1<p<\infty$, тесно связанные с пространствами
Хермандера. А именно,
$$
H^{\mu}_{p}(\mathbb{R}^{n}):=\{u\in\mathcal{S}'(\mathbb{R}^{n}):
\mathcal{F}^{-1}(\mu\,\widehat{u})\in L_{p}(\mathbb{R}^{n})\}.
$$
При этом в случае $p\neq2$ на весовую функцию $\mu$ накладывались дополнительные
условия. В гильбертовом случае $p=2$ пространства $B_{p,\mu}(\mathbb{R}^{n})$ и
$H^{\mu}_{p}(\mathbb{R}^{n})$ совпадают.
\end{remark}

Среди свойств пространств Хермандера выделим следующую важную теорему вложения
[\ref{Hermander65}, с.~59].


\begin{proposition}\label{prop2.5}
Пусть $p,q\in[1,\infty]$, $1/p+1/q=1$, непрерывная функция
$\mu:\mathbb{R}^{n}\rightarrow(0,\infty)$ весовая, и целое число $k\geq0$. Тогда
условие
\begin{equation}\label{2.33}
(1+|\xi|)^{k}\,\mu^{-1}(\xi)\in L_{q}(\mathbb{R}^{n}_{\xi})
\end{equation}
влечет непрерывное вложение $B_{p,\mu}(\mathbb{R}^{n})\hookrightarrow
C^{k}_{\mathrm{b}}(\mathbb{R}^{n})$. Обратно, если
$$
\{u\in B_{p,\mu}(\mathbb{R}^{n}):\mathrm{supp}\,u\subset V\}\subset
C^{k}(\mathbb{R}^{n})
$$
для некоторого открытого непустого множества $V\subseteq\mathbb{R}^{n}$, то
выполняется условие \eqref{2.33}.
\end{proposition}


\subsection{Уточненная шкала}\label{sec2.3.2}

С точки зрения приложений к эллиптическим операторам нам целесообразно ограничиться
изотропными гильбертовыми пространствами Хермандера $B_{2,\mu}(\mathbb{R}^{n})$, где
$\mu(\xi)=\langle\xi\rangle^{s}\varphi(\langle\xi\rangle)$, $s\in\mathbb{R}$,
$\varphi\in\mathrm{QSV}$. Здесь и далее
$\langle\xi\rangle:=(1+\xi_{1}^{2}+\ldots+\xi_{n}^{2})^{1/2}$~--- сглаженный модуль
вектора $\xi=(\xi_{1},\ldots,\xi_{n})\in\mathbb{R}^{n}$. Это пространство мы
обозначим через $H^{s,\varphi}(\mathbb{R}^{n})$. Несколько детализируя выбор
функционального параметра $\varphi$, приведем определение пространства
$H^{s,\varphi}(\mathbb{R}^{n})$.

Обозначим через $\mathcal{M}$ множество всех функций
$\varphi:[1;\infty)\rightarrow(0;\infty)$ таких, что:
\begin{itemize}
\item [(i)] $\varphi$ измерима по Борелю на полуоси $[1;\infty)$;
\item [(ii)] функции $\varphi$ и $1/\varphi$ ограничены на каждом отрезке $[1;b]$, где
$1<b<\infty$;
\item [(iii)] $\varphi\in\mathrm{QSV}$.
\end{itemize}

Из теоремы \ref{th2.10} сразу получаем следующее описание класса~$\mathcal{M}$:
$$
\varphi\in\mathcal{M}\;\Leftrightarrow\;\varphi(t)=\exp\Biggl(\beta(t)+
\int\limits_{1}^{t}\frac{\alpha(\tau)}{\tau}\,d\tau\Biggr)\;\;\mbox{при}\;\;t\geq1;
$$
здесь функция $\alpha$ непрерывна и $\alpha(\tau)\rightarrow0$ при
$\tau\rightarrow\infty$, а функция $\beta$ измерима по Борелю и ограничена на
полуоси $[1,\infty)$.

Пусть $s\in\mathbb{R}$ и $\varphi\in\mathcal{M}$.


\begin{definition}\label{def2.10}
Линейное пространство $H^{s,\varphi}(\mathbb{R}^{n})$ состоит, по определению, из
всех распределений $u\in\mathcal{S}'(\mathbb{R}^{n})$ таких, что преобразование
Фурье $\widehat{u}$ является локально суммируемой по Лебегу в $\mathbb{R}^{n}$
функцией, удовлетворяющей условию
$$
\int\limits_{\mathbb{R}^{n}}
\langle\xi\rangle^{2s}\varphi^{2}(\langle\xi\rangle)\,|\widehat{u}(\xi)|^{2}\,
d\xi<\infty.
$$
В пространстве $H^{s,\varphi}(\mathbb{R}^{n})$ определено скалярное произведение по
формуле
$$
(u_{1},u_{2})_{H^{s,\varphi}(\mathbb{R}^{n})}:=
\int\limits_{\mathbb{R}^{n}}\langle\xi\rangle^{2s}\varphi^{2}(\langle\xi\rangle)
\,\widehat{u_{1}}(\xi)\,\overline{\widehat{u_{2}}(\xi)}\,d\xi.
$$
Оно стандартным образом порождает норму.
\end{definition}

В силу леммы~\ref{lem2.3} функция
$\mu(\xi)=\langle\xi\rangle^{s}\varphi(\langle\xi\rangle)$ является весовой. В~самом
деле, для любых $\xi,\eta\in\mathbb{R}^{n}$ имеем:
\begin{gather*}
\frac{\mu(\xi)}{\mu(\eta)}=\left(\frac{\langle\xi\rangle}{\langle\eta\rangle}\right)^{s}
\frac{\varphi(\langle\xi\rangle)}{\varphi(\langle\eta\rangle)}\leq \\
c(1+|\langle\xi\rangle -\langle\eta\rangle|)^{|s|+1}\leq c(1+|\xi-\eta|)^{|s|+1},
\end{gather*}
где постоянная $c>0$ не зависит от $\xi$ и $\eta$. Если функция $\varphi$
непрерывна, то пространство
$H^{s,\varphi}(\mathbb{R}^{n})=B_{2,\mu}(\mathbb{R}^{n})$ --- частный случай
пространства Хермандера. Отметим, что замена условия непрерывности на более слабое
условие измеримости по Борелю, сделанная при определении множества $\mathcal{M}$
функциональных параметров $\varphi$, не приводит к новым пространствам.
Действительно, ввиду теоремы~\ref{th2.12} (i), для каждого $\varphi\in\mathcal{M}$
существует функция $\varphi_{1}\in C^{\infty}([1,\infty))\cap\mathcal{M}$ такая, что
$\varphi\asymp\varphi_{1}$ на полуоси $[1,\infty)$, и значит, пространства
$H^{s,\varphi}(\mathbb{R}^{n})$ и $H^{s,\varphi_{1}}(\mathbb{R}^{n})$ равны с
точностью до эквивалентности норм.

В важном частном случае $\varphi\equiv1$ пространство
$H^{s,\varphi}(\mathbb{R}^{n})$ совпадает с гильбертовым пространством Соболева
$H^{s}(\mathbb{R}^{n})$ порядка~$s$. В общем случае справедливо следующее
утверждение.


\begin{lemma}\label{lem2.7}
Пусть $s\in\mathbb{R}$ и $\varphi\in\mathcal{M}$. Тогда выполняются непрерывные
вложения
\begin{equation}\label{2.34}
H^{s+\varepsilon}(\mathbb{R}^{n})\hookrightarrow
H^{s,\varphi}(\mathbb{R}^{n})\hookrightarrow
H^{s-\varepsilon}(\mathbb{R}^{n})\quad\mbox{для любого}\;\;\varepsilon>0.
\end{equation}
\end{lemma}


\textbf{Доказательство.} Пусть $\varepsilon>0$. Поскольку
$\varphi\in\mathcal{M}\subset\mathrm{QSV}$, то в силу теоремы~\ref{th2.12} (ii)
имеем: $t^{-\varepsilon}\leq\varphi(t)\leq t^{\,\varepsilon}$ при $t\gg1$. Отсюда и
из определения класса $\mathcal{M}$ следует существование числа $c>1$ такого, что
$c^{-1}t^{-\varepsilon}\leq\varphi(t)\leq c\,t^{\,\varepsilon}$ для всех $t\geq1$.
Поэтому
$$
c^{-1}\langle\xi\rangle^{s-\varepsilon}\leq
\langle\xi\rangle^{s}\varphi(\langle\xi\rangle)\leq
c^{}\langle\xi\rangle^{s+\varepsilon}\quad\mbox{для
любого}\quad\xi\in\mathbb{R}^{n}.
$$
Последнее сразу же влечет непрерывные вложения \eqref{2.34}. Лемма~\ref{lem2.7}
доказана.

\medskip

Вложения \eqref{2.34} полезно записать в такой форме:
\begin{equation}\label{2.35}
\bigcup_{\varepsilon>0}H^{s+\varepsilon}(\mathbb{R}^{n})=:H^{s+}(\mathbb{R}^{n})
\subset H^{s,\varphi}(\mathbb{R}^{n})\subset
H^{s-}(\mathbb{R}^{n}):=\bigcap_{\varepsilon>0}H^{s-\varepsilon}(\mathbb{R}^{n}).
\end{equation}
Как видим, в семействе пространств
\begin{equation}\label{2.36}
\{H^{s,\varphi}(\mathbb{R}^{n}): s\in\mathbb{R},\varphi\in\mathcal{M}\}
\end{equation}
числовой параметр $s$ задает основную (степенную) гладкость, а функциональный
параметр $\varphi$ определяет дополнительную гладкость, подчиненную основной. В
зависимости от того, будет ли $\varphi(t)\rightarrow\infty$ или
$\varphi(t)\rightarrow0$ при $t\rightarrow\infty$, параметр $\varphi$ задает
дополнительную либо положительную, либо отрицательную гладкость. Иначе говоря,
параметр $\varphi$ уточняет основную $s$-гладкость. Поэтому естественно дать
следующее определение.


\begin{definition}\label{def2.11}
Семейство функциональных пространств \eqref{2.36} называем \emph{уточненной шкалой}
в $\mathbb{R}^{n}$ (по отношению к соболевской шкале).
\end{definition}

\subsection{Свойства уточненной шкалы}\label{sec2.3.3}

Связь между уточненной и соболевской шкалами пространств не исчерпывается вложениями
\eqref{2.35}. Оказывается, каждое пространство уточненной шкалы можно получить
посредством интерполяции с подходящим функциональным параметром пары гильбертовых
пространств Соболева.


\begin{theorem}\label{th2.14}
Пусть заданы функция $\varphi\in\mathcal{M}$ и положительные числа
$\varepsilon,\delta$. Положим
$$
\psi(t):=
\begin{cases}
\;t^{\,\varepsilon/(\varepsilon+\delta)}\,
\varphi(t^{1/(\varepsilon+\delta)}) & \text{при\;\;\;$t\geq1$}, \\
\;\varphi(1) & \text{при\;\;\;$0<t<1$}.
\end{cases}
$$
Тогда:
\begin{itemize}
\item [$\mathrm{(i)}$] функция $\psi$ принадлежит множеству $\mathcal{B}$ и является
интерполяционным параметром;
\item [$\mathrm{(ii)}$] для произвольного $s\in\mathbb{R}$ справедливо
$$
\bigl[H^{s-\varepsilon}(\mathbb{R}^{n}),H^{s+\delta}(\mathbb{R}^{n})\bigr]_{\psi}
=H^{s,\varphi}(\mathbb{R}^{n})
$$
с равенством норм.
\end{itemize}
\end{theorem}


\textbf{Доказательство.} (i) В силу теоремы \ref{th2.12} (ii), (iv) функция $\psi$
принадлежит множеству $\mathcal{B}$ и является квазиправильно меняющаяся на $\infty$
порядка $\theta=\varepsilon/(\varepsilon+\delta)\in(0,\,1)$. Следовательно, $\psi$
--- интерполяционный параметр на основании теоремы \ref{th2.11}. Пункт (i) доказан.

(ii) Пусть $s\in\mathbb{R}$. Пара соболевских пространств
$H^{s-\varepsilon}(\mathbb{R}^{n})$, $H^{s+\delta}(\mathbb{R}^{n})$ допустимая,
причем псевдодифференциальный оператор с символом
$\langle\xi\rangle^{\varepsilon+\delta}$ является порождающим оператором $J$ для
этой пары. При помощи преобразования Фурье
$$
\mathcal{F}:H^{s-\varepsilon}(\mathbb{R}^{n})\leftrightarrow
L_{2}\bigl(\mathbb{R}^{n},\langle\xi\rangle^{2(s-\varepsilon)}d\xi\bigr)
$$
оператор $J$ приводится к виду умножения на функцию
$\langle\xi\rangle^{\varepsilon+\delta}$ аргумента $\xi\in\mathbb{R}^{n}$.
Следовательно, оператор $\psi(J)$ приводится к виду умножения на функцию
$\psi(\langle\xi\rangle^{\varepsilon+\delta})=
\langle\xi\rangle^{\varepsilon}\varphi(\langle\xi\rangle)$. Отсюда ввиду
\eqref{2.35} получаем:
\begin{gather*}
\bigl[H^{s-\varepsilon}(\mathbb{R}^{n}),H^{s+\delta}(\mathbb{R}^{n})\bigr]_{\psi}=\\
\left\{u\in H^{s-\varepsilon}(\mathbb{R}^{n}):
\langle\xi\rangle^{\varepsilon}\varphi(\langle\xi\rangle)\ \widehat{u}(\xi)\in
L_{2}\bigl(\mathbb{R}^{n},\langle\xi\rangle^{2(s-\varepsilon)}d\xi\bigr)\right\}=\\
\biggl\{u\in H^{s-\varepsilon}(\mathbb{R}^{n}):
\int\limits_{\mathbb{R}^{n}}\langle\xi\rangle^{2s}\varphi^{2}(\langle\xi\rangle)
\left|\widehat{u}(\xi)\right|^{2}d\xi<\infty\biggr\}=\\
H^{s-\varepsilon}(\mathbb{R}^{n})\cap
H^{s,\varphi}(\mathbb{R}^{n})=H^{s,\varphi}(\mathbb{R}^{n}).
\end{gather*}
Кроме того, норма в пространстве $\left[H^{s-\varepsilon}(\mathbb{R}^{n}),
H^{s+\delta}(\mathbb{R}^{n})\right]_{\psi}$ равна
\begin{gather*}
\|\psi(J)u\|_{H^{s-\varepsilon}(\mathbb{R}^{n})}= \\
\biggl(\;\,\int\limits_{\mathbb{R}^{n}}
|\langle\xi\rangle^{\varepsilon}\varphi(\langle\xi\rangle)\
\widehat{u}(\xi)|^{2}\,\langle\xi\rangle^{2(s-\varepsilon)}\ d\xi\biggr)^{1/2}=
\|u\|_{H^{s,\varphi}(\mathbb{R}^{n})}.
\end{gather*}
Пункт (ii) доказан.

Теорема \ref{th2.14} доказана.

\medskip

Отметим, что интерполяция с функциональным параметром пространств Соболева,
Хермандера и некоторых других рассмотрена в работах М.~Шехтера [\ref{Schechter67}],
Г.~Шлензак [\ref{Shlenzak74}], К.~Меруччи [\ref{Merucci84}] и Ф.~Кобоса,
Д.~Л.~Фернандеза [\ref{CobosFernandez88}].

В следующей теореме даны необходимые нам свойства уточненной шкалы в
$\mathbb{R}^{n}$.


\begin{theorem}\label{th2.15}
Пусть $s\in\mathbb{R}$ и $\varphi,\varphi_{1}\in\mathcal{M}$. Справедливы следующие
утверждения.
\begin{itemize}
\item [$\mathrm{(i)}$] Для любого $\varepsilon>0$ выполняется непрерывное и
плотное вложение $H^{s+\varepsilon,\varphi_{1}}(\mathbb{R}^{n})\hookrightarrow
H^{s,\varphi}(\mathbb{R}^{n})$.
\item [$\mathrm{(ii)}$] Функция
$\varphi(t)/\varphi_{1}(t)$ ограничена в окрестности бесконечности тогда и только
тогда, когда $H^{s,\varphi_{1}}(\mathbb{R}^{n})\hookrightarrow
H^{s,\varphi}(\mathbb{R}^{n})$. Это вложение непрерывное и плотное.
\item [$\mathrm{(iii)}$] Пусть задано целое $k\geq0$. Условие
\begin{equation}\label{2.37}
\int\limits_{1}^{\infty}\frac{dt}{t\,\varphi^{\,2}(t)}<\infty
\end{equation}
равносильно вложению
\begin{equation}\label{2.38}
H^{k+n/2,\varphi}(\mathbb{R}^{n})\hookrightarrow C^{k}_{\mathrm{b}}(\mathbb{R}^{n}).
\end{equation}
Это вложение непрерывно.
\item [$\mathrm{(iv)}$] Пространства $H^{s,\varphi}(\mathbb{R}^{n})$ и
$H^{-s,1/\varphi}(\mathbb{R}^{n})$ взаимно двойственны относительно расширения по
непрерывности скалярного произведения в $L_{2}(\mathbb{R}^{n})$.
\end{itemize}
\end{theorem}

\textbf{Доказательство.} (i) В силу леммы \ref{lem2.7} имеем непрерывные вложения
$$
H^{s+\varepsilon,\varphi_{1}}(\mathbb{R}^{n})\hookrightarrow
H^{s+\varepsilon/2}(\mathbb{R}^{n})\hookrightarrow
H^{s,\varphi}(\mathbb{R}^{n})\quad\mbox{для любого}\;\;\varepsilon>0.
$$
Они плотны, поскольку множество $C^{\infty}_{0}(\mathbb{R}^{n})$ плотно в каждом из
записанных пространств. Утверждение (i) доказано.

(ii) Поскольку $\varphi,\varphi_{1}\in\mathcal{M}$, то функция $\varphi/\varphi_{1}$
ограничена в окрестности бесконечности тогда и только тогда, когда $\mu(\xi)\leq c
\mu_{1}(\xi)$ для любого $\xi\in\mathbb{R}^{n}$. Здесь
$\mu(\xi):=\langle\xi\rangle^{s}\varphi(\langle\xi\rangle)$,
$\mu_{1}(\xi):=\langle\xi\rangle^{s}\varphi_{1}(\langle\xi\rangle)$, а постоянная
$c>0$ не зависит от $\xi$. Л.~Хермандер доказал [\ref{Hermander65}, с.~55] (теорема
2.2.2), что неравенство $\mu(\xi)\leq c \mu_{1}(\xi)$ равносильно непрерывному
вложению
$$
H^{s,\varphi_{1}}(\mathbb{R}^{n})=B_{2,\mu_{1}}(\mathbb{R}^{n})\hookrightarrow
B_{2,\mu}(\mathbb{R}^{n})=H^{s,\varphi}(\mathbb{R}^{n}).
$$
Это вложение плотно, так как множество $C^{\infty}_{0}(\mathbb{R}^{n})$ плотно в
записанных пространствах. Утверждение (ii) доказано.

(iii) Пусть целое $k\geq0$. Согласно предложению \ref{prop2.5}, где полагаем
$\mu(\xi):=\langle\xi\rangle^{k+n/2}\,\varphi(\langle\xi\rangle)$ и $p=2$, имеем
\begin{equation}\label{2.39}
\int\limits_{\mathbb{R}^{n}}\;
\frac{d\xi}{\langle\xi\rangle^{n}\,\varphi^{2}(\langle\xi\rangle)}<\infty\,\;
\Leftrightarrow\;H^{k+n/2,\varphi}(\mathbb{R}^{n})\hookrightarrow
C^{k}_{\mathrm{b}}(\mathbb{R}^{n}).
\end{equation}
При этом вложение непрерывно. Переходя в интеграле к сферическим координатам и делая
замену переменной $t=(1+r^{2})^{1/2}$, запишем:
\begin{gather}\notag
\int\limits_{\mathbb{R}^{n}}\;
\frac{d\xi}{\langle\xi\rangle^{n}\,\varphi^{2}(\langle\xi\rangle)}=
c\int\limits_{1}^{\infty}\frac{r^{n-1}\,dr}{(1+r^{2})^{n/2}\,
\varphi^{2}((1+r^{2})^{1/2})}=\\
c\int\limits_{1}^{\infty}\frac{(t^{2}-1)^{(n-1)/2}}{t^{n}\,\varphi^{2}(t)}\,
\frac{t\,dt}{(t^{2}-1)^{1/2}}=
c\int\limits_{1}^{\infty}\,\frac{\omega_{n}(t)\,dt}{t\,\varphi^{2}(t)}. \label{2.40}
\end{gather}
Здесь $c>0$~--- некоторое число, а $\omega_{n}(t):=(\sqrt{t^{2}-1}/t)^{n-2}$.
Отметим, что функции $\omega_{n}$ и $1/\omega_{n}$ ограничены на полуоси
$[2,\infty)$, и поскольку $n\geq1$, то
$$
\int\limits_{1}^{2}\omega_{n}(t)dt<\infty.
$$
Ввиду определения множества $\mathcal{M}$ функция $\varphi^{-2}(t)$ ограничена на
отрезке $[1,\,2]$. Следовательно,
\begin{gather*}
\int\limits_{1}^{\infty}\frac{dt}{t\varphi^{\,2}(t)}<\infty\Leftrightarrow
\int\limits_{2}^{\infty}\frac{dt}{t\varphi^{\,2}(t)}<\infty\Leftrightarrow \\
\int\limits_{2}^{\infty}\frac{\omega_{n}(t)dt}{t\varphi^{2}(t)}<\infty\Leftrightarrow
\int\limits_{1}^{\infty}\frac{\omega_{n}(t)dt}{t\varphi^{2}(t)}<\infty.
\end{gather*}
Отсюда в силу \eqref{2.40} получаем
\begin{equation}\label{2.41}
\int\limits_{1}^{\infty}\frac{dt}{t\varphi^{\,2}(t)}<\infty\;\Leftrightarrow\;
\int\limits_{\mathbb{R}^{n}}\;
\frac{d\xi}{\langle\xi\rangle^{n}\,\varphi^{2}(\langle\xi\rangle)}<\infty.
\end{equation}
Формулы \eqref{2.39} и \eqref{2.41} дают утверждение (iii).

(iv) Утверждение (iv) --- это частный случай теоремы 2.2.9 из монографии Хермандера
[\ref{Hermander65}, с.~61]. Отметим, что
$\varphi\in\mathcal{M}\Leftrightarrow1/\varphi\in\mathcal{M}$, и значит,
пространство $H^{-s,1/\varphi}(\mathbb{R}^{n})$ принадлежит уточненной шкале.

Теорема \ref{th2.15} доказана.




\newpage


\markright{\emph \ref{sec2.4b}. Равномерно эллиптические операторы}

\section[Равномерно эллиптические операторы в уточненной шкале]
{Равномерно эллиптические \\ операторы в уточненной шкале}\label{sec2.4b}

\markright{\emph\ref{sec2.4b}. Равномерно эллиптические операторы}

Здесь мы изучим равномерно эллиптические псевдодифференциальные операторы, заданные
в евклидовом пространстве $\mathbb{R}^{n}$. Мы получим априорную оценку решения
эллиптического уравнения в уточненной шкале пространств и исследуем его внутреннюю
гладкость.

\subsection{Псевдодифференциальные операторы}\label{sec2.4.1b}

Обстоятельное изложение теории псевдодифференциальных операторов (ПДО) приведено,
например, в монографиях [\ref{Taylor85}, \ref{Treves84a}, \ref{Hermander87},
\ref{Shubin78}] и в обзоре [\ref{Agranovich90}]. Для удобства читателя напомним
определение псевдодифференциального оператора в $\mathbb{R}^{n}$ и связанные с ним
понятия, необходимые нам. Будем пользоваться в основном терминологией и
обозначениями из обзора М.~С.~Аграновича [\ref{Agranovich90}] (\S~1, 3).

Пусть $m\in\mathbb{R}$. Обозначим через $S^{m}(\mathbb{R}^{2n})$ множество всех
функций $a\in C^{\infty}(\mathbb{R}^{2n})$, которые удовлетворяют следующему
условию: для произвольных мультииндексов $\alpha,\beta$ существует число
$c_{\alpha,\beta}>0$ такое, что
$$
|\,\partial_{x}^{\alpha}\,\partial_{\xi}^{\beta}\,a(x,\xi)\,|
\leq\,c_{\alpha,\beta}\,\langle\xi\rangle^{m-|\beta|}\quad\mbox{для любых}\quad
x,\xi\in\mathbb{R}^{n}.
$$

\begin{definition}\label{defPs1}
\emph{Псевдодифференциальный оператор} $A$ в $\mathbb{R}^{n}$ с \emph{символом}
$a\in S^{m}(\mathbb{R}^{2n})$ задается формулой
$$
(Au)(x):=(2\pi)^{-n}
\int\limits_{\mathbb{R}^{n}}\,e^{-ix\cdot\xi}\,a(x,\xi)\,\widehat{u}(\xi)\,d\xi,\quad
x\in\mathbb{R}^{n};
$$
здесь $u$~--- произвольная функция из пространства $\mathcal{S}(\mathbb{R}^{n})$.
\end{definition}

Обозначим через $\Psi^{m}(\mathbb{R}^{n})$ класс всех ПДО в $\mathbb{R}^{n}$ с
символами из $S^{m}(\mathbb{R}^{2n})$.

Важный пример ПДО класса $\Psi^{m}(\mathbb{R}^{n})$ дает линейный дифференциальный
оператор
\begin{equation}\label{Ps1}
A(x,D):=\sum_{|\mu|\leq m}\;a_{\mu}(x)\,D^{\mu}
\end{equation}
порядка $m$ и с коэффициентами $a_{\mu}\in C_{\mathrm{b}}(\mathbb{R}^{n})$. Его
символ равен
$$
a(x,\xi)=\sum_{|\mu|\leq m}\;a_{\mu}(x)\,\xi^{\mu},\quad x,\xi\in\mathbb{R}^{n}.
$$
Здесь $\mu=(\mu_{1},\ldots,\mu_{n})$~--- мультииндекс,
$D^{\mu}:=i^{|\mu|}\,\partial^{\mu}_{x}$, и как обычно,
$\xi^{\mu}:=\xi_{1}^{\mu_{1}}\ldots\xi_{n}^{\mu_{n}}$ для вектора
$\xi=(\xi_{1},\ldots,\xi_{n})$. (Преобразование Фурье переводит оператор
дифференцирования $D^{\mu}$ в оператор умножения на $\xi^{\mu}$.)

Заметим, что с увеличением параметра $m$ класс $\Psi^{m}(\mathbb{R}^{n})$
расширяется. Положим
$$
\Psi^{-\infty}(\mathbb{R}^{n}):=\bigcap_{m\in\mathbb{R}}\,\Psi^{m}(\mathbb{R}^{n}),\quad
\Psi^{\infty}(\mathbb{R}^{n}):=\bigcup_{m\in\mathbb{R}}\,\Psi^{m}(\mathbb{R}^{n}).
$$

Всякий ПДО $A\in\Psi^{\infty}(\mathbb{R}^{n})$ непрерывен в
$\mathcal{S}(\mathbb{R}^{n})$ и единственным образом продолжается до линейного
непрерывного оператора в $\mathcal{S}'(\mathbb{R}^{n})$. Этот оператор также
обозначается через $A$.

Далее перейдем к более узкому и наиболее важному для приложений семейству
полиоднородных (или классических) ПДО. Сначала введем множество однородных символов
порядка $m\in\mathbb{R}$. Обозначим через $S^{m}_{\mathrm{h}}(\mathbb{R}^{2n})$
множество всех бесконечно дифференцируемых комплекснозначных функций $b(x,\xi)$,
заданных при $x,\xi\in\mathbb{R}^{n}$ и $\xi\neq0$, которые удовлетворяют следующим
условиям:
\begin{itemize}
\item [(i)] $b(x,\lambda\xi)=\lambda^{m}b(x,\xi)$ для произвольного $\lambda>0$;
\item [(ii)] для произвольных мультииндексов $\alpha,\beta$ существует число
$c_{\alpha,\beta}>0$ такое, что
$$
|\,\partial_{x}^{\alpha}\,\partial_{\xi}^{\beta}\,b(x,\xi)\,|
\leq\,c_{\alpha,\beta}\,|\xi|^{m-|\beta|}\;\;\mbox{для любых}\;\;
x,\xi\in\mathbb{R}^{n},\;\xi\neq0.
$$
\end{itemize}

\begin{definition}\label{defPs2}
Пусть $A$~--- ПДО из $\Psi^{m}(\mathbb{R}^{n})$ с символом $a(x,\xi)$. $A$
называется \emph{полиоднородным} (или \emph{классическим}) ПДО порядка $m$, если
найдется такая убывающая (строго) последовательность вещественных чисел $m_{0}=m$,
$m_{j}$, $j=1,2,\ldots$, стремящаяся к $-\infty$, и существует такая
последовательность однородных символов $a_{j}\in
S^{m_{j}}_{\mathrm{h}}(\mathbb{R}^{2n})$, $j=0,1,2,\ldots$, что
$$
a(x,\xi)-\theta(\xi)\sum_{j=0}^{k}\,a_{j}(x,\xi)\in
S^{m_{k+1}}(\mathbb{R}^{2n}_{(x,\xi)})
$$
для каждого целого $k\geq0$. При этом функция $a_{0}(x,\xi)$ предполагается отличной
от тождественного нуля, а функция $\theta\in
C^{\infty}_{\mathrm{b}}(\mathbb{R}^{n})$ полагается равной $1$ вне некоторой
окрестности начала и равной $0$ в несколько меньшей его окрестности. Функция
$a_{0}(x,\xi)$ называется \emph{главным символом} полиоднородного ПДО $A$.
\end{definition}

Обозначим через $\Psi^{m}_{\mathrm{ph}}(\mathbb{R}^{n})$ класс всех полиоднородных
ПДО в $\mathbb{R}^{n}$ порядка $m$. Он не зависит от указанного выбора функции
$\theta$. Заметим, что поскольку главный символ полиоднородного ПДО предполагается
не равным тождественно нулю, то
$\Psi^{m}_{\mathrm{ph}}(\mathbb{R}^{n})\cap\Psi^{r}_{\mathrm{ph}}(\mathbb{R}^{n})=\varnothing$
при $m\neq r$.

Дифференциальный оператор \eqref{Ps1} принадлежит к классу
$\Psi^{m}_{\mathrm{ph}}(\mathbb{R}^{n})$, если отличен от тождественного нуля хотя
бы один коэффициент $a_{\mu}(x)$ с $|\mu|=m$. Главный символ этого оператора имеет
вид
$$
a_{0}(x,\xi)=\sum_{|\mu|=m}\;a_{\mu}(x)\,\xi^{\mu},\quad x,\xi\in\mathbb{R}^{n}.
$$

Нас интересуют полиоднородные ПДО равномерно эллиптические в $\mathbb{R}^{n}$.
Напомним их определение.

\begin{definition}\label{defPs3}
ПДО $A\in\Psi^{m}_{\mathrm{ph}}(\mathbb{R}^{n})$ и его главный символ $a_{0}(x,\xi)$
называются \emph{равномерно эллиптическими} в $\mathbb{R}^{n}$, если существует
число $c>0$ такое, что $|a_{0}(x,\xi)|\geq c$ для любых $x,\xi\in\mathbb{R}^{n}$ с
$|\xi|=1$.
\end{definition}

Важно, что равномерно эллиптический ПДО $A$ имеет параметрикс $B$, т. е. справедливо
следующее утверждение [\ref{Agranovich90}, с.~20] (теорема 1.8.3).


\begin{proposition}\label{prop2.6b}
Пусть ПДО $A\in\Psi^{m}_{\mathrm{ph}}(\mathbb{R}^{n})$ равномерно эллиптический в
$\mathbb{R}^{n}$. Тогда существует ПДО $B\in\Psi^{-m}_{\mathrm{ph}}(\mathbb{R}^{n})$
равномерно эллиптический в $\mathbb{R}^{n}$ такой, что
\begin{equation}\label{2.46b}
BA=I+T_{1},\quad AB=I+T_{2},
\end{equation}
где $T_{1}$ и $T_{2}$ --- некоторые ПДО из $\Psi^{-\infty}(\mathbb{R}^{n})$, а
$I$~--- тождественный оператор в $\mathcal{S}'(\mathbb{R}^{n})$.
\end{proposition}

\subsection{Априорная оценка решения}\label{sec2.4.2b}

Пусть задан ПДО $A\in\Psi^{m}_{\mathrm{ph}}(\mathbb{R}^{n})$, где $m\in\mathbb{R}$.
В этом и следующем пунктах предполагается, что ПДО $A$ равномерно эллиптический в
$\mathbb{R}^{n}$. Рассмотрим уравнение $Au=f$ в $\mathbb{R}^{n}$. Для его решения
$u$ установим априорную оценку в уточненной шкале.

Предварительно докажем следующую лемму о действии произвольного ПДО в уточненной
шкале.


\begin{lemma}\label{lem2.8b}\label{lem2.8}
Пусть $G$~--- ПДО класса $\Psi^{r}(\mathbb{R}^{n})$, где $r\in\mathbb{R}$. Тогда
сужение отображения $u\mapsto Gu$, $u\in\mathcal{S}'(\mathbb{R}^{n})$, является
линейным ограниченным оператором
\begin{equation}\label{2.44b}
G:H^{\sigma,\varphi}(\mathbb{R}^{n})\rightarrow H^{\sigma-r,\varphi}(\mathbb{R}^{n})
\end{equation}
для произвольных параметров $\sigma\in\mathbb{R}$ и $\varphi\in\mathcal{M}$.
\end{lemma}


\textbf{Доказательство.} В соболевском случае $\varphi\equiv1$ эта лемма известна
(см., например, [\ref{Agranovich90}, с.~10] (теорема 1.1.2) либо [\ref{Hermander87},
с.~107] (теорема 18.1.13)). Выберем произвольные $\sigma\in\mathbb{R}$ и
$\varphi\in\mathcal{M}$. Рассмотрим линейные ограниченные операторы
$$
G:H^{\sigma\mp1}(\mathbb{R}^{n})\rightarrow H^{\sigma\mp1-r}(\mathbb{R}^{n}).
$$
Применим интерполяцию с функциональным параметром $\psi$ из теоремы~\ref{th2.14},
где полагаем $\varepsilon=\delta=1$. В силу ее пункта (i) имеем ограниченный
оператор
$$
G:\bigl[H^{\sigma-1}(\mathbb{R}^{n}),H^{\sigma+1}(\mathbb{R}^{n})\bigr]_{\psi}\rightarrow
\bigl[H^{\sigma-r-1}(\mathbb{R}^{n}),H^{\sigma-r+1}(\mathbb{R}^{n})\bigr]_{\psi}.
$$
Отсюда в силу пункта (ii) вытекает, что ПДО $G$ определяет ограниченный оператор
\eqref{2.44b}. Лемма~\ref{lem2.8b} доказана.

\medskip

Как видим, ПДО класса $\Psi^{r}(\mathbb{R}^{n})$ понижает на $r$ основную гладкость
$\sigma$ и оставляет инвариантной дополнительную гладкость $\varphi$ пространства
$H^{\sigma,\varphi}(\mathbb{R}^{n})$.

В силу леммы~\ref{lem2.8b} имеем ограниченный оператор
$$
A:H^{s+m,\varphi}(\mathbb{R}^{n})\rightarrow H^{s,\varphi}(\mathbb{R}^{n})
$$
для произвольных параметров $s\in\mathbb{R}$ и $\varphi\in\mathcal{M}$.


\begin{theorem}\label{th2.16b}
Пусть $s\in\mathbb{R}$, $\sigma>0$ и $\varphi\in\mathcal{M}$. Существует число
$c=c(s,\sigma,\varphi)>0$ такое, что для произвольного распределения $u\in
H^{s+m,\varphi}(\mathbb{R}^{n})$ справедлива априорная оценка
\begin{equation}\label{2.49b}
\|u\|_{H^{s+m,\varphi}(\mathbb{R}^{n})}\leq
c\,\bigl(\,\|Au\|_{H^{s,\varphi}(\mathbb{R}^{n})}+
\|u\|_{H^{s+m-\sigma,\varphi}(\mathbb{R}^{n})}\,\bigr).
\end{equation}
\end{theorem}


\textbf{Доказательство.} Воспользуемся предложением \ref{prop2.6b}. В силу первого
равенства в \eqref{2.46b} запишем $u=BAu-T_{1}u$. Отсюда следует оценка
\eqref{2.49b}:
\begin{gather*}
\|u\|_{H^{s+m,\varphi}(\mathbb{R}^{n})}= \|BAu-T_{1}u\|_{H^{s+m,\varphi}(\mathbb{R}^{n})}\leq\\
\|BAu\|_{H^{s+m,\varphi}(\mathbb{R}^{n})}+\|T_{1}u\|_{H^{s+m,\varphi}(\mathbb{R}^{n})}\leq\\
c\,\|Au\|_{H^{s,\varphi}(\mathbb{R}^{n})}+c\,\|u\|_{H^{s+m-\sigma,\varphi}(\mathbb{R}^{n})}.
\end{gather*}
Здесь $c$ --- максимум норм операторов
\begin{gather}\label{2.50b}
B:\,H^{s,\varphi}(\mathbb{R}^{n})\rightarrow
H^{s+m,\varphi}(\mathbb{R}^{n}),\\
T_{1}:\,H^{s+m-\sigma,\varphi}(\mathbb{R}^{n}) \rightarrow
H^{s+m,\varphi}(\mathbb{R}^{n}).\label{2.51b}
\end{gather}
Эти операторы ограниченные в силу леммы \ref{lem2.8b} и предложения \ref{prop2.6b}.
Теорема~\ref{th2.16b} доказана.

Эта теорема уточняет известную априорную оценку решения равномерно эллиптического
уравнения в соболевской шкале (случай $\varphi\equiv1$) [\ref{Agranovich90}, с.~20]
(теорема 1.8.4).

\subsection{Гладкость решения}\label{sec2.4.3b}

Предположим, что правая часть уравнения  $Au=f$ имеет некоторую внутреннюю гладкость
в уточненной шкале на заданном открытом непустом множестве
$V\subseteq\mathbb{R}^{n}$. Изучим внутреннюю гладкости решения $u$ на этом
множестве. Рассмотрим сначала случай, когда $V=\mathbb{R}^{n}$. Обозначим через
$H^{-\infty}(\mathbb{R}^{n})$ объединение всех пространств
$H^{s,\varphi}(\mathbb{R}^{n})$, где $s\in\mathbb{R}$, $\varphi\in\mathcal{M}$. В
линейном пространстве $H^{-\infty}(\mathbb{R}^{n})$ вводится топология индуктивного
предела.


\begin{theorem}\label{th2.17b}
Предположим, что $u\in H^{-\infty}(\mathbb{R}^{n})$ является решением уравнения
$Au=f$ в $\mathbb{R}^{n}$, где $f\in H^{s,\varphi}(\mathbb{R}^{n})$ для некоторых
параметров $s\in\mathbb{R}$ и $\varphi\in\mathcal{M}$. Тогда $u\in
H^{s+m,\varphi}(\mathbb{R}^{n})$
\end{theorem}


\textbf{Доказательство.} В силу теоремы \ref{th2.15} (i) для распределения $u\in
H^{-\infty}(\mathbb{R}^{n})$ существует число $\sigma>0$ такое, что
\begin{equation}\label{2.52b}
u\in H^{s+m-\sigma,\varphi}(\mathbb{R}^{n}).
\end{equation}
Отсюда и из условия теоремы получаем на основании формул \eqref{2.46b},
\eqref{2.50b}, \eqref{2.51b} требуемое свойство:
$$
u=BAu-T_{1}u=Bf-T_{1}u \in H^{s+m,\varphi}(\mathbb{R}^{n}).
$$
Теорема~\ref{th2.17b} доказана.

\medskip

Рассмотрим теперь общий случай, когда $V$ --- произвольное открытое непустое
подмножество пространства $\mathbb{R}^{n}$. Положим
\begin{gather*}
H^{\sigma,\varphi}_{\mathrm{int}}(V):=\bigl\{w\in H^{-\infty}(\mathbb{R}^{n}):
\chi\,w\in H^{\sigma,\varphi}(\mathbb{R}^{n})\\ \mbox{для всех}\;\;\chi\in
C^{\infty}_{\mathrm{b}}(\mathbb{R}^{n}),\,\mathrm{supp}\,\chi\subset V,\;
\mathrm{dist}(\mathrm{supp}\,\chi,\partial V)>0\bigr\};
\end{gather*}
здесь $\sigma\in\mathbb{R}$ и $\varphi\in\mathcal{M}$. Топология в пространстве
$H^{\sigma,\varphi}_{\mathrm{int}}(V)$ задается полунормами
$w\mapsto\|\chi\,w\|_{H^{\sigma,\varphi}(\mathbb{R}^{n})}$, где функции $\chi$~---
те же, что и в определении этого пространства.


\begin{theorem}\label{th2.18b}
Предположим, что $u\in H^{-\infty}(\mathbb{R}^{n})$ является решением уравнения
$Au=f$ на множестве $V$, где $f\in H^{s,\varphi}_{\mathrm{int}}(V)$ для некоторых
параметров $s\in\mathbb{R}$ и $\varphi\in\mathcal{M}$. Тогда $u\in
H^{s+m,\varphi}_{\mathrm{int}}(V)$.
\end{theorem}


\textbf{Доказательство.} Покажем, что из условия $f\in\nobreak
H^{s,\varphi}_{\mathrm{int}}(V)$ вытекает следующее свойство повышения внутренней
гладкости решения уравнения $Au=f$: для каждого числа $r\geq1$ справедлива
импликация
\begin{equation}\label{2.55b}
u\in H^{s-r+m,\varphi}_{\mathrm{int}}(V)\;\Rightarrow\; u\in
H^{s-r+1+m,\varphi}_{\mathrm{int}}(V).
\end{equation}

Произвольно выберем функцию $\chi\in C^{\infty}_{\mathrm{b}}(\mathbb{R}^{n})$ такую,
что
\begin{equation}\label{2.56b}
\mathrm{supp}\,\chi\subset V\quad\mbox{и}\quad
\mathrm{dist}(\mathrm{supp}\,\chi,\partial V)>0.
\end{equation}
Для нее существует функция $\eta\in C^{\infty}_{\mathrm{b}}(\mathbb{R}^{n})$ такая,
что
\begin{equation}\label{2.57b}
\mathrm{supp}\,\eta\subset V,\;\;\mathrm{dist}(\mathrm{supp}\,\eta,\,\partial
V)>0,\;\;\eta=1\;\mbox{в окрестности}\;\mathrm{supp}\,\chi.
\end{equation}
Действительно, мы можем определить указанную функцию при помощи операции свертки
$\eta:=\chi_{2\varepsilon}\ast\omega_{\varepsilon}$, где
$\varepsilon:=\mathrm{dist}(\mathrm{supp}\,\chi,\partial V)/4$,
$\chi_{2\varepsilon}$ --- индикатор $2\varepsilon$-окрестности множества
$\mathrm{supp}\,\chi$, а функция $\omega_{\varepsilon}\in
C^{\infty}(\mathbb{R}^{n})$ удовлетворяет условиям
$$
\omega_{\varepsilon}\geq0,\;\;
\mathrm{supp}\,\omega_{\varepsilon}\subseteq\{x\in\mathbb{R}^{n}:\|x\|\leq
\varepsilon\},\;\;\int\limits_{\mathbb{R}^{n}}\omega_{\varepsilon}(x)\,dx=1.
$$
Непосредственно можно убедиться в том, что такая функция $\eta$ принадлежит классу
$C^{\infty}_{\mathrm{b}}(\mathbb{R}^{n})$ и имеет следующее свойство: $\eta\equiv1$
в $\varepsilon$-окрестности множества $\mathrm{supp}\,\chi$ и $\eta\equiv0$ вне
$3\varepsilon$-окрестности этого же множества, т.~е. $\eta$ удовлетворяет условиям
\eqref{2.57b}.

Переставив ПДО $A$ и оператор умножения на функцию $\chi$, запишем
\begin{gather}\notag
A\chi u=A\chi\eta u=\chi\,A\eta u+A'\eta u=\\ \notag
\chi\,Au+\chi\,A(\eta-1)u+A'\eta
u=\\
=\chi f+\chi\,A(\eta-1)u+A'\eta u\quad\mbox{в}\quad\mathbb{R}^{n}. \label{2.58b}
\end{gather}
Здесь ПДО $A'\in\Psi^{m-1}(\mathbb{R}^{n})$~--- коммутатор ПДО $A$ и оператора
умножения на функцию $\chi$. В силу леммы~\ref{lem2.8b} имеем ограниченный оператор
$$
A':\,H^{s-r+m,\varphi}(\mathbb{R}^{n})\rightarrow H^{s-r+1,\varphi}(\mathbb{R}^{n}).
$$
Следовательно,
\begin{equation}\label{2.59b}
u\in H^{s-r+m,\varphi}_{\mathrm{int}}(V)\;\Rightarrow\; A'\eta u\in
H^{s-r+1,\varphi}(\mathbb{R}^{n}).
\end{equation}
Далее, на основании условия $f\in\nobreak H^{s,\varphi}_{\mathrm{int}}(V)$ и ввиду
неравенства $r\geq1$ имеем
\begin{equation}\label{2.60b}
\chi f\in H^{s,\varphi}(\mathbb{R}^{n})\hookrightarrow
H^{s-r+1,\varphi}(\mathbb{R}^{n}).
\end{equation}
Кроме того, так как носители функций $\chi$ и $\eta-1$ не пересекаются, то ПДО
$$
\chi A(\eta-1)\in\Psi^{-\infty}(\mathbb{R}^{n}).
$$
Это сразу следует из формулы для символа композиции двух ПДО: $\chi A$ и оператора
умножения на функцию $\eta-1$ (см. [\ref{Agranovich90}, с.~13], теорема 1.2.4).
Отсюда, поскольку для $u\in H^{-\infty}(\mathbb{R}^{n})$ справедливо \eqref{2.52b}
при некотором $\sigma>0$, мы имеем в силу леммы~\ref{lem2.8b} включение
\begin{equation}\label{2.61b}
\chi\,A(\eta-1)u\in H^{s-r+1,\varphi}(\mathbb{R}^{n}).
\end{equation}

На основании формул \eqref{2.58b}~-- \eqref{2.61b} и теоремы~\ref{th2.17b} получаем,
что
\begin{gather*}
u\in H^{s-r+m,\varphi}_{\mathrm{int}}(V)\;\Rightarrow\; A\chi u\in
H^{s-r+1,\varphi}(\mathbb{R}^{n})\;\Rightarrow \\ \chi u\in
H^{s-r+1+m,\varphi}(\mathbb{R}^{n}).
\end{gather*}
Тем самым доказано \eqref{2.55b} ввиду произвольности выбора функции $\chi\in
C^{\infty}_{\mathrm{b}}(\mathbb{R}^{n})$, удовлетворяющей условию \eqref{2.56b}.

Теперь при помощи импликации \eqref{2.55b} легко вывести включение $u\in
H^{s,\varphi}_{\mathrm{int}}(V)$. Можно считать, что в формуле \eqref{2.52b} число
$\sigma>0$ целое. Значит, $u\in H^{s-\sigma+m,\varphi}_{\mathrm{int}}(V)$. Применив
\eqref{2.55b} последовательно для $r=\sigma,\,\sigma\nobreak-\nobreak1,\ldots,1$,
выводим требуемое включение:
\begin{gather*}
u\in H^{s-\sigma+m,\varphi}_{\mathrm{int}}(V) \;\Rightarrow\; u\in
H^{s-\sigma+1+m,\varphi}_{\mathrm{int}}(V)\;\Rightarrow\ldots\\ \Rightarrow\; u\in
H^{s+m,\varphi}_{\mathrm{int}}(V).
\end{gather*}

Теорема \ref{th2.18b} доказана.

\medskip

Теорема \ref{th2.18b} уточняет применительно к шкале пространств
$H^{s,\varphi}(\mathbb{R}^{n})$ известные утверждения о повышении внутренней
гладкости решений линейных эллиптических уравнений в соболевской шкале (см.,
например, [\ref{Berezansky65}, с.~189; \,\ref{Egorov84}, с.~77; \,\ref{Taylor85},
с.~72]). Как видим, уточненная гладкость $\varphi$ правой части эллиптического
уравнения наследуется его решением. Если оператор $A$ дифференциальный и множество
$V$ ограничено, то теорема~\ref{th2.18b} содержится в теореме Л.~Хермандера
[\ref{Hermander65}, с.~237] о регулярности решений гипоэллиптических уравнений,
выраженной в терминах пространств $B_{p,\mu}(\mathbb{R}^{n})$.


\begin{remark}\label{rem2.7b}
Следует различать \emph{внутреннюю} и \emph{локальную} гладкость на открытом
множестве $V\subset\mathbb{R}^{n}$. Пространство распределений, имеющих данную
локальную гладкость на этом множестве, определяется следующим образом:
\begin{gather*}
H^{\sigma,\varphi}_{\mathrm{loc}}(V):=\bigl\{w\in H^{-\infty}(\mathbb{R}^{n}):\\
\chi\,w\in H^{\sigma,\varphi}(\mathbb{R}^{n})\;\;\mbox{для всех}\;\;\chi\in
C^{\infty}_{0}(\mathbb{R}^{n}),\,\mathrm{supp}\,\chi\subset V\bigr\}.
\end{gather*}
В случае, когда множество $V$ ограничено, пространства
$H^{\sigma,\varphi}_{\mathrm{int}}(V)$ и $H^{\sigma,\varphi}_{\mathrm{loc}}(V)$
совпадают. Если же $V$ неограниченно, то может быть строгое включение
$H^{\sigma,\varphi}_{\mathrm{int}}(V)\subset H^{\sigma,\varphi}_{\mathrm{loc}}(V)$.
Для локальной уточненной гладкости справедлив аналог теоремы~\ref{th2.18b}; в ее
формулировке следует только заменить $\mathrm{int}$ на $\mathrm{loc}$ в обозначениях
пространств. Он легко выводится из теоремы~\ref{th2.18b}.
\end{remark}


Теорема~\ref{th2.18b} позволяет установить наличие непрерывных производных у решения
уравнения $Au=f$. При этом используется теорема~\ref{th2.15} (iii).


\begin{theorem}\label{th2.19b}
Пусть задано целое число $r\geq0$ и функция $\varphi\in\mathcal{M}$, удовлетворяющая
условию \eqref{2.37}. Предположим, что распределение $u\in
H^{-\infty}(\mathbb{R}^{n})$ является решением уравнения $Au=\nobreak f$ на открытом
множестве $V\subseteq\mathbb{R}^{n}$, где
\begin{equation}\label{2.62b}
f\in H^{r-m+n/2,\varphi}_{\mathrm{int}}(V).
\end{equation}
Тогда решение $u$ имеет на множестве $V$ непрерывные частные производные до порядка
$r$ включительно, причем эти производные ограничены на каждом множестве
$V_{0}\subset V$ таком, что $\mathrm{dist}(V_{0},\partial V)>0$. В~частности, если
$V=\mathbb{R}^{n}$, то $u\in C^{r}_{\mathrm{b}}(\mathbb{R}^{n})$.
\end{theorem}


\textbf{Доказательство.} В силу теоремы~\ref{th2.18b}, где полагаем $s:=r-m+n/2$,
справедливо включение $u\in H^{r+n/2,\varphi}_{\mathrm{int}}(V)$. Пусть функция
$\eta\in\nobreak C^{\infty}_{\mathrm{b}}(\mathbb{R}^{n})$ удовлетворяет условиям
$$
\mathrm{supp}\,\eta\subset V,\;\;\mathrm{dist}(\mathrm{supp}\,\eta,\,\partial
V)>0,\;\;\eta=1\;\mbox{в окрестности}\;V_{0}.
$$
Эта функция строится так же, как и в доказательстве теоремы~\ref{th2.18b}, если
заменить в нем множество $\mathrm{supp}\,\chi$ на $V_{0}$. Для распределения $\eta
u$ в силу теоремы~\ref{th2.15} (iii) имеем
$$
\eta u\in H^{r+n/2,\varphi}(\mathbb{R}^{n})\hookrightarrow
C^{r}_{\mathrm{b}}(\mathbb{R}^{n}).
$$
Отсюда вытекает, что все частные производные функции $u$ до порядка $r$ включительно
непрерывны и ограничены в некоторой окрестности множества $V_{0}$. Тогда эти
производные непрерывны и на множестве $V$, поскольку можно взять $V_{0}:=\{x_{0}\}$
для любой точки $x_{0}\in V$. Теорема~\ref{th2.19b} доказана.


\begin{remark}\label{rem2.8b}
Если использовать теорему~\ref{th2.19b} лишь для шкалы пространств Соболева, то
придется вместо \eqref{2.62b} потребовать, чтобы для некоторого числа
$\varepsilon>0$ выполнялось условие
$$
f\in H^{r-m+n/2+\varepsilon,1}_{\mathrm{int}}(V).
$$
Оно завышает основную гладкость правой части уравнения $Au=f$, что существенно
огрубляет результат.
\end{remark}

\begin{remark}\label{rem2.8c}
Условие \eqref{2.37} не только достаточное в теореме \ref{th2.19b}, но также и
необходимое на классе всех рассматриваемых решений уравнения $Au=f$. А именно,
\eqref{2.37} эквивалентно импликации
\begin{equation}\label{2.62c}
\bigl(u\in H^{-\infty}(\mathbb{R}^{n}),\:f:=Au\in
H^{r-m+n/2,\varphi}_{\mathrm{int}}(V)\bigr)\:\Rightarrow\:u\in C^{r}(V).
\end{equation}
В самом деле, если $u\in H^{r+n/2,\varphi}_{\mathrm{int}}(V)$, то $f=Au\in
H^{r-m+n/2,\varphi}_{\mathrm{int}}(V)$, откуда $u\in C^{r}(\Omega)$, если
\eqref{2.62c} выполняется. Поэтому \eqref{2.62c} влечет за собой \eqref{2.37} на
основании предложения \ref{2.5} и ввиду формулы \eqref{2.41}.
\end{remark}




\newpage

\markright{\emph \ref{notes1}. Примечания и комментарии}

\section{Примечания и комментарии}\label{notes1}

\markright{\emph \ref{notes1}. Примечания и комментарии}

\small

\textbf{К п. 1.1.} Первый метод интерполяции пространств был предложен независимо
Ж.-Л.~Лионсом [\ref{Lions58}] и С.~Г.~Крейном [\ref{Krein60a}]. Это~--- интерполяция
со степенным параметром пар гильбертовых пространств, в которой показатель степени
используется в качестве параметра интерполяции. Указанный метод был распространен на
пары нормированных пространств Ж.-Л.~Лионсом, Ж.~Петре [\ref{Lions59},
\ref{LionsPeetre61}, \ref{Peetre63}] (вещественная интерполяция) и С.~Г.~Крейном
[\ref{Krein60a}, \ref{Krein60b}], Ж.-Л.~Лионсом [\ref{Lions60}], А.~П.~Кальдероном
[\ref{Calderon64}], М.~Шехтером [\ref{Schechter67}](комплексная или голоморфная
интерполяция). Вообще говоря, комплексные и вещественные методы интерполяции
приводят к разным интерполяционным пространствам. Принципы построения общих
интерполяционных методов разработаны Э.~Гальярдо [\ref{Gagliardo68}].

В настоящее время известны различные методы вещественной и комплексной интерполяции
нормированных и более общих топологических пространств, в которых числовые наборы
выступают в качестве параметров интерполяции; см. монографии К.~Беннетта и Р.~Шапли
[\ref{BennetSharpley88}], Й.~Берга и Й.~Лёфстрёма [\ref{BergLefstrem80}],
Ю.~А.~Брудного и Н.~Я.~Кругляка [\ref{BrudnyiKrugljak91}], С.~Г.~Крейна,
Ю.~И.~Петунина и Е.~М.~Семёнова [\ref{KreinPetuninSemenov78}], Ж.-Л.~Лионса и
Э.~Мадженеса [\ref{LionsMagenes71}], В.~И.~Овчинникова [\ref{Ovchinnikov84}],
Х.~Трибеля [\ref{Triebel80}] и приведенную там обширную литературу.

Метод интерполяции нормированных пространств с функциональным параметром впервые
появился в статье К.~Фояша и Ж.-Л.~Лионса [\ref{FoiasLions61}], где также был
отдельно рассмотрен гильбертов случай. Интерполяции гильбертовых пространств с
функциональным параметром изучалась в работах В.~Ф.~Донохью [\ref{Donoghue67}],
Г.~Шлензак [\ref{Shlenzak74}], Е.~И.~Пустыльника [\ref{Pustylnik82}],
В.~И.~Овчинникова [\ref{Ovchinnikov84}] и авторов [\ref{06UMJ2}, \ref{06Dop6},
\ref{08MFAT1}]. В.~И.~Овчинников [\ref{Ovchinnikov84}] описал с помощью интерполяции
с функциональным параметром все гильбертовы пространства, интерполяционные
относительно заданной пары гильбертовых пространств. Отметим, что в ряде приложений
такую интерполяцию именуют методом переменных гильбертовых шкал; см., например,
работы М.~Хегланда [\ref{Hegland95}, \ref{Hegland10}], П.~Матэ и У.~Тотенхана
[\ref{MatheTautenhahn06}].

Различные методы интерполяции нормированных пространств с общим функциональным
параметром введены и изучены в работах Т.~Ф.~Калугиной [\ref{Kalugina75}],
Я.~Густавсcона [\ref{Gustavsson78}], С.~Янсона [\ref{Janson81}], К.~Меруччи
[\ref{Merucci82}], Л.-Э.~Персcона [\ref{Persson86}], Н.~Я.~Кругляка
[\ref{Krugljak93}] для вещественной интерполяции и в статье М.~И.~Карро,
Й.~Л.~Керд\'{а} [\ref{CarroCerda90}] для комплексной интерполяции. Отметим, что в
интерполяционных методах, предложенных С.~Янсоном [\ref{Janson81}], использован
очень широкий класс интерполяционных параметров~--- произвольные положительные
псевдовогнутые функции.

Все теоремы п. 1.1 установлены авторами: теоремы 1.1~-- 1.5 и 1.8, 1.9 доказаны в
[\ref{08MFAT1}] (п.~2), а теоремы 1.6 и 1.7 приведены в [\ref{06Dop6}] (п.~3).
Доказательство теоремы 1.4 аналогично доказательству, данному в книге Х.~Трибеля
[\ref{Triebel80}, с. 136, 138], а теоремы  1.5 близко к доказательству Ж.~Жеймона
[\ref{Geymonat65}, с. 133]. Х.~Трибель и Ж.~Жеймона рассматривали интерполяционные
функторы на категории всех соместимых пар банаховых пространств.

\medskip

\textbf{К п. 1.2.} Понятие правильно меняющейся функции введено И.~Караматой
[\ref{Karamata30a}] (в случае непрерывной функции). Им же [\ref{Karamata30b},
\ref{Karamata33}] установлены основные свойства правильно меняющихся функций. Теория
этих функций и различные ее приложения изложены в монографиях Н.~Х.~Бинхема,
Ч.~М.~Голди и Дж.~Л.~Тайгельза [\ref{BinghamGoldieTeugels89}], Г.~Л.~Гелук и
Л.~де~Хаан [\ref{GelukHaan87}], В.~Марика [\ref{Maric00}], С.~И.~Решника
[\ref{Reshnick87}], Е.~Сенеты [\ref{Seneta85}], Л.~де~Хаан [\ref{Haan70}].

Понятие квазиправильно меняющейся функции введено авторами в [\ref{06Dop6}, с.~15;
\ref{08MFAT1}, с.~90]; оно удобно в теории интерполяции пространств. Все теоремы п.
1.2 установлены в [\ref{08MFAT1}] (п.~3.1) за исключением теоремы 1.13.
Доказательство последней ранее не приводилось. Прямое доказательство важной
интерполяционной теоремы 1.11 дано авторами в [\ref{06UMJ2}] (п.~2). Оно использует
модификацию интерполяционного метода следов Ж.-Л.~Лионса [\ref{LionsMagenes71}]
(гл.~1, \S~3, 5), предложенную Г.~Шлензак [\ref{Shlenzak74}, с.~49]. Вспомогательные
леммы 1.3 и 1.4 установлены в [\ref{06UMJ2}] (п.~1).

\medskip

\textbf{К п. 1.3.} Теория распределений (обобщенных функций) берет начало от работ
С.~Л.~Соболева и Л.~Шварца. С.~Л.~Соболев ввел важные банаховы пространства
распределений, которые носят его имя и без которых невозможно представить
современную теорию дифференциальных уравнений. Изложение теории распределений и
пространств Соболева имеется, например, в монографиях С.~Л.~Соболева
[\ref{Sobolev50}, \ref{Sobolev74}], Л.~Шварца [\ref{Schwartz50}, \ref{Schwartz51}],
Р.~А.~Адамса [\ref{Adams75}], В.~С.~Владимирова [\ref{Vladimirov79}],
И.~М.~Гельфанда и Г.~Е.~Шилова [\ref{GelfandShilov58}, \ref{GelfandShilov59}],
С.~Г.~Михлина [\ref{Mihlin68}], Л.~Тартара [\ref{Tartar07}]. Разнообразные
приложения соболевских пространств стимулировали их глубокое изучение, что привело к
построению новых важных классов пространств распределений~--- пространства
Никольского, Бесова, шкала пространств Лизоркина--Трибеля, различные их весовые и
анизотропные аналоги; см. монографии О.~В.~Бесова, В.~П.~Ильина, С.~М.~Никольского
[\ref{BesovIlinNikolsky75}], С.~М.~Никольского [\ref{Nikolsky69}], Х.~Трибеля
[\ref{Triebel80}, \ref{Triebel86}, \ref{Triebel92}, \ref{Triebel06}] и приведенные
там литературные указания. Эти пространства параметризуются числовыми наборами.

Дальнейшее обобщение пространств Соболева получено на пути перехода от числовых к
функциональным параметрам. Последние описывают более тонко, чем числовые параметры,
свойства регулярности распределений, входящих в пространства. Это было сделано
Б.~Мальгранжем [\ref{Malgrange57}], а систематически Л. Хермандером в монографии
[\ref{Hermander65}] (гл.~II) и Л.~Р.~Волевичем, Б.~П.~Панеяхом в статье
[\ref{VolevichPaneah65}], которые ввели и изучили пространства, параметризуемые при
помощи весьма общего функционального параметра. Пространства Хермандера и
Волевича--Панеяха совпадают в гильбертовом случае. Некоторые приложения этих
пространств в теории уравнений с частными производными даны в монографиях
Л.~Хермандера [\ref{Hermander65}, \ref{Hermander86}] и Б.~П.~Панеяха
[\ref{Paneah00}].

Пространства Хермандера и Волевича--Панеяха занимают центральное место среди
пространств обобщенной гладкости. В~последние десятилетия такие пространства
являются предметом ряда глубоких исследований~--- см. обзоры П.~И.~Лизоркина
[\ref{Lizorkin86}], Г.~А.~Калябина и П.~И.~Лизоркина [\ref{KalyabinLizorkin87}],
монографии Х.~Трибеля [\ref{Triebel01}] (гл.~III), Н.~Якоба [\ref{Jacob010205}],
недавние работы В.~И.~Буренкова [\ref{Burenkov99}], Х.-Г.~Леопольда
[\ref{Leopold98}], С.~Д.~Моуры [\ref{Moura01}], Б.~Опика и В.~Требелза
[\ref{OpicTrebels00}], В.~Фаркаса и Х.-\nobreak Г.~Леопольда
[\ref{FarkasLeopold06}], В.~Фаркаса, Н.~Якоба и Р.~Л.~Шилинга
[\ref{FarkasJacobScilling01a}, \ref{FarkasJacobScilling01b}], Д.~Д.~Хароске и
С.~Д.~Моуры [\ref{HaroskeMoura04}], Д.~Э.~Эдмундса, П.~Гурки и Б.~Опика
[\ref{EdmundsGurkaOpic97}], Д.~Э.~Эдмундса и Д.~Д.~Хароске [\ref{EdmundsHaroske99}],
а также приведенную там литературу. Построены различные аналоги пространств
Никольского--Бесова и Лизоркина--Трибеля, параметризуемые посредством функциональных
параметров. Для них удалось получить точные теоремы вложения одних классов
пространств в другие, теоремы о продолжении, теоремы о следах и другие результаты. В
статьях М.~Шехтера [\ref{Schechter67}], С.~Меруччи [\ref{Merucci84}], Ф.~Кобоса и
Д.~Л.~Фернандеза [\ref{CobosFernandez88}] изучены интерполяционные свойства
различных классов пространств обобщенной гладкости.

Определение уточненной шкалы в $\mathbb{R}^{n}$ дано авторами в [\ref{05UMJ5}]; ее
свойства (теоремы 1.14 и 1.15) доказаны в [\ref{06UMJ3}] (п.~3).

Термин „уточненная шкала” использовался ранее Г.~Шлензак [\ref{Shlenzak74}] для
иного введенного ею класса гильбертовых пространств Хермандера. Этот класс не имеет
конструктивного описания и не привязан к соболевской шкале. Соответствующий ему
класс интерполяционных параметров довольно узок и подчинен излишнему условию на
поведение параметров вблизи нуля.

В книге Х.~Трибеля [\ref{Triebel01}] (гл.~III) и в статье Д.~Д.~Хароске и
С.~Д.~Моуры [\ref{HaroskeMoura04}] введены и изучены аналоги банаховых пространств
Никольского--Бесова и Лизоркина--Трибеля, регулярность распределений в которых
характеризуется, как и у нас, посредством двух параметров~--- основного числового и
дополнительного функционального. При этом Х.~Трибель использует логарифмические, а
Д.~Д.~Хароске и С.~Д.~Моура~--- более общие медленно меняющиеся функциональные
параметры.

\medskip

\textbf{К п. 1.4.} Алгебра псевдодифференциальных операторов (ПДО) построена и
изучена в основном в работах Дж.~Дж.~Кона, Л.~Ниренберга и Л.~Хермандера; см.
сборник статей [\ref{PsDO67}]. Среди ПДО важное место занимают эллиптические
операторы. Они нашли разнообразные приложения в теории эллиптических краевых задач
для дифференциальных уравнений, в спектральной теории, в теории функциональных
пространств и другие. Подробное изложение теории ПДО приведено, например, в
монографиях М.~Тейлора [\ref{Taylor85}], Ф.~Трева [\ref{Treves84a}], Л.~Хермандера
[\ref{Hermander87}], М.~А.~Шубина [\ref{Shubin78}] и обзорах М.~С.~Аграновича
[\ref{Agranovich65}, \ref{Agranovich90}]. Мы в основном следуем обозначениям и
терминологии, принятым в [\ref{Agranovich90}].

Для эллиптических дифференциальных уравнений известны внутренние априорные оценки
решений в подходящих парах пространств Гельдера (с нецелыми индексами) и пространств
Соболева, а также утверждения о локальной регулярности (или гладкости) решений в
таких пространствах. Изложение этих результатов и соответствующие литературные
указания приведены в книге Ю.~М.~Березанского [\ref{Berezansky65}] (гл.~III, \S~4).
В~пространствах Хермандера регулярность решений гипоэллиптических дифференциальных
уравнений изучена Л.~Хермандером в монографиях [\ref{Hermander65}, с. 139, 237] и
[\ref{Hermander86}, с. 79, 222]. Для эллиптических псевдодифференциальных уравнений
внутренние априорные оценки решений и утверждения о локальной регулярности решений в
соболевской шкале приведены в указанных выше работах по теории ПДО. Для ПДО,
равномерно эллиптических в $\mathbb{R}^{n}$, априорная оценка верна на всем
$\mathbb{R}^{n}$. Это верно и для утверждения о повышении регулярности решений; см.,
например, обзор М.~С.~Аграновича [\ref{Agranovich90}] (п. 1.8).

Все теоремы п. 1.4 установлены в [\ref{09Collection1}] применительно к более общим
равномерно эллиптическим по Петровскому матричным ПДО.

\normalsize



\chapter[Пространства Хермандера на замкнутом многообразии и их приложения]
{\textbf{Пространства Хермандера на замкнутом многообразии \\ и их приложения}}
\label{ch2b}

\chaptermark{\emph Гл. \ref{ch2b}. Пространства Хермандера на многообразии}



\section[Пространства Хермандера на замкнутом многообразии]
{Пространства Хермандера \\ на замкнутом многообразии}\label{sec2.5}

\markright{\emph \ref{sec2.5}. Пространства Хермандера на замкнутом многообразии}

В этом параграфе мы рассматриваем пространства Хермандера (уточненную шкалу) на
гладком замкнутом многообразии. Мы дадим эквивалентные определения этих пространств,
аналогичные тем, что используются для пространств Соболева.

\subsection{Эквивалентные определения}\label{sec2.5.1}

Далее в этой главе $\Gamma$ --- замкнутое (т. е. компактное и без края) бесконечно
гладкое ориентированное многоообразие размерности $n\geq1$. Предполагается, что на
$\Gamma$ задана некоторая $C^{\infty}$-плотность $dx$. Напомним, что
$\mathcal{D}'(\Gamma)$~--- топологическое линейное пространство всех распределений
на $\Gamma$, двойственное к пространству $C^{\infty}(\Gamma)$ относительно
расширения по непрерывности скалярного произведения в пространстве
$L_{2}(\Gamma,dx)=:L_{2}(\Gamma)$. Это расширение будем обозначать через
$(f,w)_{\Gamma}$, где $f\in\mathcal{D}'(\Gamma)$, $w\in C^{\infty}(\Gamma)$.

Пусть $s\in\mathbb{R}$ и $\varphi\in\mathcal{M}$. Дадим три эквивалентных
определения пространства $H^{s,\varphi}(\Gamma)$.

Первое определение характеризует $H^{s,\varphi}(\Gamma)$ исходя из локальных свойств
распределения $f\in\mathcal{D}'(\Gamma)$. Возьмем какой-либо конечный атлас из
$C^{\infty}$-структуры на $\Gamma$, образованный локальными картами
$\alpha_{j}:\mathbb{R}^{n}\leftrightarrow\Gamma_{j}$, где $j=1,\ldots,r$. Здесь
открытые множества $\Gamma_{j}$ составляют конечное покрытие многообразия $\Gamma$.
Пусть функции $\chi_{j}\in C^{\infty}(\Gamma)$, где $j=1,\ldots,r$, образуют
разбиение единицы на $\Gamma$, удовлетворяющее условию
$\mathrm{supp}\,\chi_{j}\subset\Gamma_{j}$.


\begin{definition}[локальное]\label{def2.13}
Линейное пространство $H^{s,\varphi}(\Gamma)$ состоит, по определению, из всех
распределений $f\in\mathcal{D}'(\Gamma)$ таких, что $(\chi_{j}f)\circ\alpha_{j}\in
H^{s,\varphi}(\mathbb{R}^{n})$ для каждого $j=1,\ldots,r$. Здесь
$(\chi_{j}f)\circ\alpha_{j}$
--- представление распределения $\chi_{j}f$ в локальной карте $\alpha_{j}$.
В пространстве $H^{s,\varphi}(\Gamma)$ определено скалярное произведение
распределений $f_{1}$, $f_{2}$ по формуле
$$
(f_{1},f_{2})_{H^{s,\varphi}(\Gamma)}:=\sum_{j=1}^{r}\,((\chi_{j}f_{1})\circ\alpha_{j},
(\chi_{j}\,f_{2})\circ\alpha_{j})_{H^{s,\varphi}(\mathbb{R}^{n})}.
$$
Оно стандартным образом задает норму.
\end{definition}


В важном частном случае $\varphi\equiv1$ пространство $H^{s,\varphi}(\Gamma)$
совпадает с пространством Соболева $H^{s}(\Gamma)$ порядка $s$. Последнее, как
известно [\ref{Hermander65}, с.~82; \ref{Shubin78}, с.~66],~--- полное и с точностью
до эквивалентных норм не зависит от выбора атласа и разбиения единицы на~$\Gamma$.
Следующее определение связывает пространство $H^{s,\varphi}(\Gamma)$ с
пространствами Соболева посредством интерполяции и показывает, что
$H^{s,\varphi}(\Gamma)$ не зависит (с точностью до эквивалентности норм) от выбора
атласа и разбиения единицы.


\begin{definition}[интерполяционное]\label{def2.14}
Пусть целые числа $k_{0}$ и $k_{1}$ такие, что $k_{0}<s<k_{1}$. По определению,
$$
H^{s,\varphi}(\Gamma):=\bigl[H^{k_{0}}(\Gamma),H^{k_{1}}(\Gamma)\bigr]_{\psi},
$$
где интерполяционный параметр $\psi$ задан формулой
$$
\psi(t)=
\begin{cases}
\;t^{\,(s-k_{0})/(k_{1}-k_{0})}\,
\varphi(t^{1/(k_{1}-k_{0})}) & \text{при\;\;\;$t\geq1$}, \\
\;\varphi(1) & \text{при\;\;\;$0<t<1$}.
\end{cases}
$$
\end{definition}


\begin{remark}
Функция $\psi$ из определения \ref{def2.14} является интерполяционным параметром в
силу теоремы \ref{th2.11}, поскольку $\psi$ --- правильно меняющаяся функция на
бесконечности порядка $\theta=(s-k_{0})/(k_{1}-k_{0})\in(0,\,1)$.
\end{remark}


В спектральной теории полезно определение пространства $H^{s,\varphi}(\Gamma)$,
связывающее норму в нем с некоторой функцией от оператора $1-\Delta_{\Gamma}$, где
$\Delta_{\Gamma}$~--- оператор Бельтрами--Лапласа на $\Gamma$. (При этом на
многообразии $\Gamma$ вводится риманова метрика; см., например, [\ref{Treves84b},
с.~336; \ref{Shubin78}, с.~163]).


\begin{definition}[операторное]\label{def2.15} Гильбертово пространство
$H^{s,\varphi}(\Gamma)$ определяется как пополнение линейного многообразия
$C^{\infty}(\Gamma)$ по норме
$$
f\,\mapsto\,
\|(1-\Delta_{\Gamma})^{s/2}\varphi((1-\Delta_{\Gamma})^{1/2})\,f\|_{L_{2}(\Gamma)},\quad
f\in C^{\infty}(\Gamma).
$$
\end{definition}


\begin{theorem}\label{th2.20}
Определения $\ref{def2.13}$, $\ref{def2.14}$ и $\ref{def2.15}$ эквивалентные: они
определяют одно и то же гильбертово пространство $H^{s,\varphi}(\Gamma)$ с точностью
до эквивалентности норм.
\end{theorem}


Мы докажем теорему \ref{th2.20} в пп. \ref{sec2.5.2} и \ref{sec2.5.3}. В связи с
этой теоремой уместно дать следующее определение.


\begin{definition}\label{def2.16} Семейство гильбертовых пространств
$$
\{H^{s,\varphi}(\Gamma):\,s\in\mathbb{R},\;\varphi\in\mathcal{M}\}
$$
называем уточненной шкалой на замкнутом многообразии $\Gamma$.
\end{definition}

\subsection{Интерполяционные свойства}\label{sec2.5.2}

В этом пункте мы изучим интерполяционные свойства уточненной шкалы на $\Gamma$ и,
как следствие, получим доказательство того, что определения $\ref{def2.13}$ и
$\ref{def2.14}$ равносильны. В качестве исходного определения пространства
$H^{s,\varphi}(\Gamma)$ используем локальное определение $\ref{def2.13}$.


\begin{theorem}\label{th2.21}
Пусть заданы функция $\varphi\in\mathcal{M}$ и положительные числа
$\varepsilon,\delta$. Тогда для произвольного $s\in\mathbb{R}$ справедливо
\begin{equation}\label{2.63}
\bigl[H^{s-\varepsilon}(\Gamma),\
H^{s+\delta}(\Gamma)\bigr]_{\psi}=H^{s,\varphi}(\Gamma)
\end{equation}
с эквивалентностью норм. Здесь $\psi$ --- интерполяционный параметр из теоремы
$\ref{th2.14}$.
\end{theorem}


\textbf{Доказательство.} Записанная слева в равенстве \eqref{2.63} пара соболевских
пространств допустимая [\ref{Shubin78}, с. 66, 68]. Мы выведем это равенство из
теоремы \ref{th2.14} при помощи известного приема „распрямления” и „склейки”
многообразия $\Gamma$. Согласно определению \ref{def2.13} линейное отображение
„распрямления”
$$
T:\,f\mapsto(\,(\chi_{1}f)\circ\alpha_{1},\ldots,(\chi_{r}f)\circ\alpha_{r}\,),
\quad f\in\mathcal{D}'(\Gamma),
$$
определяет изометрические операторы
\begin{gather}\label{2.64}
T:\,H^{\sigma}(\Gamma)\rightarrow(H^{\sigma}(\mathbb{R}^{n}))^{r},\quad
\sigma\in\mathbb{R},\\
T:\,H^{s,\varphi}(\Gamma)\rightarrow(H^{s,\varphi}(\mathbb{R}^{n}))^{r}.
\label{2.65}
\end{gather}
Поскольку параметр $\psi$ интерполяционный, то из ограниченности операторов
\eqref{2.64}, где $\sigma\in\{s-\varepsilon,s+\delta\}$, вытекает ограниченность
оператора
\begin{equation}\label{2.66}
T:\bigl[H^{s-\varepsilon}(\Gamma),H^{s+\delta}(\Gamma)\bigr]_{\psi}
\rightarrow\bigl[(H^{s-\varepsilon}(\mathbb{R}^{n}))^{r},
(H^{s+\delta}(\mathbb{R}^{n}))^{r}\bigr]_{\psi}.
\end{equation}
В силу теорем \ref{th2.5} и \ref{th2.14} справедливы следующие равенства пространств
и норм в них:
\begin{gather}\notag
\bigl[\,(H^{s-\varepsilon}(\mathbb{R}^{n}))^{r},\,
(H^{s+\delta}(\mathbb{R}^{n}))^{r}\,\bigr]_{\psi}=\\
\bigl(\,\bigl[H^{s-\varepsilon}(\mathbb{R}^{n}),
H^{s+\delta}(\mathbb{R}^{n})\bigr]_{\psi}\,\bigr)^{r}
=(H^{s,\varphi}(\mathbb{R}^{n}))^{r}. \label{2.67}
\end{gather}
Поэтому ограниченность оператора \eqref{2.66} означает ограниченность оператора
\begin{equation}\label{2.68}
T:\,\bigl[H^{s-\varepsilon}(\Gamma),H^{s+\delta}(\Gamma)\bigr]_{\psi}
\rightarrow(H^{s,\varphi}(\mathbb{R}^{n}))^{r}.
\end{equation}

Построим далее для $T$ левый обратный оператор „склейки”~$K$. Для каждого номера
$j=1,\ldots,r$ возьмем функцию $\eta_{j}\in C_{0}^{\infty}(\mathbb{R}^{n})$ такую,
что $\eta_{j}=1$ на множестве $\alpha_{j}^{-1}(\mathrm{supp}\,\chi_{j})$. Рассмотрим
линейное отображение
$$
K:\,(h_{1},\ldots,h_{r})\mapsto\sum_{j=1}^{r}\,
\Theta_{j}\bigl((\eta_{j}h_{j})\circ\alpha_{j}^{-1}\bigr),\quad
h_{1},\ldots,h_{r}\in\mathcal{S}'(\mathbb{R}^{n}).
$$
Здесь $(\eta_{j}h_{j})\circ\alpha_{j}^{-1}$~--- такое распределение, заданное в
открытом множестве $\Gamma_{j}\subseteq\Gamma$, что его представитель в локальной
карте $\alpha_{j}$ имеет вид $\eta_{j}h_{j}$. Кроме того, здесь также
$\Theta_{j}$~--- оператор продолжения нулем с $\Gamma_{j}$ на многообразие $\Gamma$.
Оператор $\Theta_{j}$ корректно определен на распределениях, носитель которых лежит
в $\Gamma_{j}$. В силу выбора функций $\chi_{j},\eta_{j}$ имеем:
\begin{gather*}
KTf=\sum_{j=1}^{r}\,\Theta_{j}\left(\bigl(\eta_{j}\,
((\chi_{j}f)\circ\alpha_{j})\bigr)\circ\alpha_{j}^{-1}\right)=\\
\sum_{j=1}^{r}\,\Theta_{j}\left((\chi_{j}f)\circ\alpha_{j}\circ\alpha_{j}^{-1}\right)=
\sum_{j=1}^{r}\,\chi_{j}f=f,
\end{gather*}
т. е.
\begin{equation}\label{2.69}
KTf=f\quad\mbox{для любого}\quad f\in\mathcal{D}'(\Gamma).
\end{equation}

Покажем, что линейное отображение $K$ определяет ограниченный оператор
\begin{equation}\label{2.70}
K:\,(H^{s,\varphi}(\mathbb{R}^{n}))^{r}\rightarrow H^{s,\varphi}(\Gamma).
\end{equation}
Для произвольного вектора $h=(h_{1},\ldots,h_{r})$ из пространства
$(H^{s,\varphi}(\mathbb{R}^{n}))^{r}$ можно записать:
\begin{gather}\notag
\bigl\|Kh\bigr\|^{2}_{H^{s,\varphi}(\Gamma)} =\sum_{l=1}^{r}\;\bigl\|(\chi_{l}\,Kh)
\circ\alpha_{\,l}\bigr\|_{H^{s,\varphi}(\mathbb{R}^{n})}^{2}=\\ \notag
\sum_{l=1}^{r}\,\Bigl\|\Bigl(\chi_{\,l}\,\sum_{j=1}^{r}\
\Theta_{j}\bigl((\eta_{j}h_{j})\circ\alpha_{j}^{-1}\bigr)\Bigr)
\circ\alpha_{\,l}\,\Bigr\|_{H^{\,s,\varphi}(\mathbb{R}^{n})}^{2}=\\ \notag
\sum_{l=1}^{r}\;\Bigl\|\,\sum_{j=1}^{r}(\eta_{j,l}\,h_{j})
\circ\beta_{\,j,l}\,\Bigr\|_{H^{s,\varphi}(\mathbb{R}^{n})}^{2}\leq \\
\sum_{l=1}^{r}\;\Bigl(\sum_{j=1}^{r}\,\bigl\|(\eta_{j,l}\,h_{j})\circ
\beta_{j,l}\,\bigr\|_{H^{\,s,\varphi}(\mathbb{R}^{n})}\Bigr)^{2}. \label{2.71}
\end{gather}
Здесь $\eta_{j,l}:=(\chi_{\,l}\circ\alpha_{j})\,\eta_{j}\in
C_{0}^{\infty}(\mathbb{R}^{n})$, а
$\beta_{j,l}:\,\mathbb{R}^{n}\,\leftrightarrow\,\mathbb{R}^{n}$ суть такой
$C^{\infty}$-диффеоморфизм, что $\beta_{j,l}\,=\alpha_{j}^{-1}\circ\alpha_{l}$ в
окрестности множества $\mathrm{supp}\,\eta_{j,l}$ и $\beta_{j,l}(x)=x$ для всех
$x\in\mathbb{R}^{n}$, достаточно больших по модулю. Как известно [\ref{Hermander87},
c.~625], оператор умножения на функцию класса $C_{0}^{\infty}(\mathbb{R}^{n})$ и
оператор замены переменных $u\mapsto u\circ\beta_{j,l}$ ограничены в каждом
пространстве $H^{\sigma}(\mathbb{R}^{n})$, где $\sigma\in\mathbb{R}$. Поэтому
линейный оператор $v\mapsto (\eta_{j,l}\,v)\circ\beta_{j,l}$ ограничен в
пространстве $H^{\sigma}(\mathbb{R}^{n})$. Отсюда в силу теоремы \ref{th2.14}
вытекает его ограниченность в пространстве $H^{s,\varphi}(\mathbb{R}^{n})$.
Следовательно, соотношения \eqref{2.71} влекут оценку
$$
\bigl\|Kh\bigr\|_{H^{s,\varphi}(\Gamma)}^{2}\leq c\,\sum_{j=1}^{r}\
\bigl\|h_{j}\bigr\|_{H^{s,\varphi}(\mathbb{R}^{n})}^{2},
$$
где число $c>0$ не зависит от вектора $h=(h_{1},\ldots,h_{r})$. Таким образом,
оператор \eqref{2.70} ограничен для любых $s\in\mathbb{R}$ и
$\varphi\in\mathcal{M}$.

В частности, ограничены операторы $K:\,(H^{\sigma}(\mathbb{R}^{n}))^{r}\rightarrow
H^{\sigma}(\Gamma)$, где $\sigma\in\mathbb{R}$. Взяв здесь значения
$\sigma\in\{s-\varepsilon,\,s+\delta\}$ и применив интерполяцию с параметром $\psi$,
ввиду равенства \eqref{2.67} получаем ограниченный оператор
\begin{equation}\label{2.72}
K:\,(H^{s,\varphi}(\mathbb{R}^{n}))^{r}\rightarrow\
\bigl[H^{s-\varepsilon}(\Gamma),\,H^{s+\delta}(\Gamma)\bigr]_{\psi}.
\end{equation}
Теперь из формул \eqref{2.65}, \eqref{2.72} и \eqref{2.69} следует, что
тождественный оператор $KT$ осуществляет непрерывное вложение пространства
$H^{s,\varphi}(\Gamma)$ в интерполяционное пространство
$[H^{s-\varepsilon}(\Gamma),\,H^{s+\delta}(\Gamma)]_{\psi}$. А~из формул
\eqref{2.68} и \eqref{2.70} вытекает, что тот же оператор $KT$ осуществляет обратное
непрерывное вложение.

Теорема \ref{th2.21} доказана.

\medskip

В случае, когда в теореме \ref{th2.21} числа $k_{0}:=s-\varepsilon$ и
$k_{1}:=s+\delta$ целые, получаем следующий важный результат.


\begin{corollary}\label{cor2.2}
Для произвольно заданных $s\in\mathbb{R}$ и $\varphi\in\mathcal{M}$ определения
$\ref{def2.13}$ и $\ref{def2.14}$ являются эквивалентными.
\end{corollary}


Из теоремы \ref{th2.21} и свойств интерполяции вытекают следующие свойства
уточненной шкалы на многообразии $\Gamma$.


\begin{theorem}\label{th2.22}
Пусть $s\in\mathbb{R}$ и $\varphi,\varphi_{1}\in\mathcal{M}$. Справедливы следующие
утверждения.
\begin{itemize}
\item [$\mathrm{(i)}$] Пространство $H^{s,\varphi}(\Gamma)$ полное (гильбертово)
и с точностью до эквивалентности норм не зависит от выбора атласа многообразия
$\Gamma$ и разбиения единицы, использованных в определении $\ref{def2.13}$.
\item [$\mathrm{(ii)}$] Множество $C^{\infty}(\Gamma)$ плотно в пространстве
$H^{s,\varphi}(\Gamma)$.
\item [$\mathrm{(iii)}$] Для произвольного числа $\varepsilon>0$ выполняется
компактное и плотное вложение $H^{s+\varepsilon,\varphi_{1}}(\Gamma)\hookrightarrow
H^{s,\varphi}(\Gamma)$.
\item [$\mathrm{(iv)}$] Функция $\varphi/\varphi_{1}$ ограничена в окрестности
$\infty$ тогда и только тогда, когда $H^{s,\varphi_{1}}(\Gamma)\hookrightarrow
H^{s,\varphi}(\Gamma)$. Это вложение плотно и непрерывно. Оно компактно тогда и
только тогда, когда $\varphi(t)/\varphi_{1}(t)\rightarrow0$ при
$t\rightarrow\infty$.
\item [$\mathrm{(v)}$] Пространства $H^{s,\varphi}(\Gamma)$ и
$H^{-s,1/\varphi}(\Gamma)$ взаимно двойственны с точностью до эквивалентности норм
относительно расширения по непрерывности скалярного произведения в $L_{2}(\Gamma)$.
\end{itemize}
\end{theorem}


\textbf{Доказательство.} (i) Пространство $H^{s,\varphi}(\Gamma)$ полное в силу
теоремы \ref{th2.21}~--- как результат интерполяции гильбертовых соболевских
пространств. Рассмотрим две пары $\mathcal{A}_{1}$ и $\mathcal{A}_{2}$, каждая из
которых образована атласом и соответствующим разбиением единицы на $\Gamma$,
использованным в определении \ref{def2.13}. Соответствующие этим парам пространства
уточненной шкалы и пространства Соболева обозначим соответственно через
$H^{s,\varphi}(\Gamma,\mathcal{A}_{j})$ и $H^{\sigma}(\Gamma,\mathcal{A}_{j})$, где
$j\in\{1,2\}$. Для соболевских пространств тождественное отображение определяет
гомеоморфизм $I:H^{\sigma}(\Gamma,\mathcal{A}_{1})\leftrightarrow
H^{\sigma}(\Gamma,\mathcal{A}_{2})$ при любом $\sigma\in\mathbb{R}$. Возьмем здесь
$\sigma=s\mp1$ и применим интерполяцию с параметром $\psi$ из теоремы \ref{th2.14}.
В силу теоремы \ref{th2.21} получаем гомеоморфизм
$I:H^{s,\varphi}(\Gamma,\mathcal{A}_{1})\leftrightarrow
H^{s,\varphi}(\Gamma,\mathcal{A}_{2})$. Он устанавливает указанную независимость
пространства $H^{s,\varphi}(\Gamma)$ от выбора атласа и разбиения единицы.
Утверждение (i) доказано.

(ii) В силу теорем \ref{th2.1} и \ref{th2.21} имеем непрерывное и плотное вложение
$H^{s+\delta}(\Gamma)\hookrightarrow H^{s,\varphi}(\Gamma)$. Как известно (см.,
например, [\ref{Shubin78}, с.~66], предложение~7.4), множество $C^{\infty}(\Gamma)$
плотно в соболевском пространстве $H^{s+\delta}(\Gamma)$. Следовательно, это
множество плотно в $H^{s,\varphi}(\Gamma)$. Утверждение (ii) доказано.

(iii) Пусть $\varepsilon>0$. Согласно теореме \ref{th2.21} существуют
интерполяционные параметры $\chi,\eta\in\mathcal{B}$ такие, что справедливы
следующие равенства пространств с точностью до эквивалентности норм в них:
\begin{gather*}
\bigl[H^{s+\varepsilon/2}(\Gamma),H^{s+2\varepsilon}(\Gamma)\bigr]_{\chi}=
H^{s+\varepsilon,\varphi_{1}}(\Gamma),\\
\bigl[H^{s-\varepsilon}(\Gamma),H^{s+\varepsilon/3}(\Gamma)\bigr]_{\eta}=
H^{s,\varphi}(\Gamma).
\end{gather*}
Отсюда в силу теоремы \ref{th2.1} получаем цепочку непрерывных вложений
$$
H^{s+\varepsilon,\varphi_{1}}(\Gamma)\hookrightarrow
H^{s+\varepsilon/2}(\Gamma)\hookrightarrow
H^{s+\varepsilon/3}(\Gamma)\hookrightarrow H^{s,\varphi}(\Gamma).
$$
Здесь центральное вложение соболевских пространств компактно [\ref{Shubin78}, с.~68]
(теорема~7.4). Следовательно, компактно вложение
$H^{s+\varepsilon,\varphi_{1}}(\Gamma)\hookrightarrow H^{s,\varphi}(\Gamma)$. Оно
плотно в силу пункта (ii). Утверждение (iii) доказано.

(iv) Пусть функция $\varphi/\varphi_{1}$ ограничена в окрестности $\infty$. На
основании теоремы \ref{th2.21} имеем следующие равенства пространств с точностью до
эквивалентности норм в них:
\begin{gather*}
\left[H^{s-1}(\Gamma),H^{s+1}(\Gamma)\right]_{\psi} =H^{s,\varphi}(\Gamma),\\
\left[H^{s-1}(\Gamma),H^{s+1}(\Gamma)\right]_{\psi_{1}} =H^{s,\varphi_{1}}(\Gamma).
\end{gather*}
Здесь интерполяционные параметры $\psi,\psi_{1}\in\mathcal{B}$ удовлетворяют
ус\-ло\-вию $\psi(t)/\psi_{1}(t)=\varphi(t^{1/2})/\varphi_{1}(t^{1/2})$ при
$t\geq1$. Следовательно, функция $\psi/\psi_{1}$ ограничена в окрестности $\infty$ и
по теореме \ref{th2.2} выполняется непрерывное и плотное вложение
$H^{s,\varphi_{1}}(\Gamma)\hookrightarrow H^{s,\varphi}(\Gamma)$. Далее, если
$\varphi(t)/\varphi_{1}(t)\rightarrow0$ при $t\rightarrow\infty$, то
$\psi(t)/\psi_{1}(t)\rightarrow0$ при $t\rightarrow\infty$. Отсюда и из компактности
вложения соболевских пространств $H^{s+1}(\Gamma)\hookrightarrow H^{s-1}(\Gamma)$
вытекает в силу теоремы \ref{th2.2} компактность вложения
$H^{s,\varphi_{1}}(\Gamma)\hookrightarrow H^{s,\varphi}(\Gamma)$.

Покажем теперь, что из вложения $H^{s,\varphi_{1}}(\Gamma)\hookrightarrow
H^{s,\varphi}(\Gamma)$ следует ограниченность функции $\varphi/\varphi_{1}$ на
полуоси $[1,\infty)$. Предположим, что это вложение выполняется. Воспользуемся
локальным определением \ref{def2.13}. Можно считать, что для некоторого открытого
непустого множества $U\subset\Gamma$ выполняются условия: $U\subset\Gamma_{1}$ и
$U\cap\Gamma_{j}=\varnothing$ при $j\neq1$. Для произвольного распределения $w\in
H^{s,\varphi_{1}}(\mathbb{R}^{n})$ такого, что
$\mathrm{supp}\,w\subset\alpha_{1}^{-1}(U)$, имеем
$$
\Theta(w\circ\alpha_{1}^{-1})\in H^{s,\varphi_{1}}(\Gamma)\hookrightarrow
H^{s,\varphi}(\Gamma).
$$
Здесь $\Theta$~--- оператор продолжения нулем распределения с множества $U$
на~$\Gamma$. Тогда
$$
w=\bigl(\chi_{1}(\Theta(w\circ\alpha_{1}^{-1}))\bigl)\circ\alpha_{1}\in
H^{s,\varphi}(\mathbb{R}^{n}).
$$
Следовательно, в силу теоремы 2.2.2 из монографии Л.~Хермандера [\ref{Hermander65},
с.~55] функция $\varphi(\langle\xi\rangle)/\varphi_{1}(\langle\xi\rangle)$ аргумента
$\xi\in\mathbb{R}^{n}$ ограничена. Значит, функция $\varphi/\varphi_{1}$ ограничена
на полуоси $[1,\infty)$.

Наконец, покажем, что из компактности вложения
$H^{s,\varphi_{1}}(\Gamma)\hookrightarrow H^{s,\varphi}(\Gamma)$ вытекает условие
$\varphi(t)/\varphi_{1}(t)\rightarrow0$ при $t\rightarrow\infty$. Если это вложение
компактно, то компактен и оператор вложения
$$
\{w\in
H^{s,\varphi_{1}}(\mathbb{R}^{n}):\,\mathrm{supp}\,w\subset\alpha_{1}^{-1}(U)\}
\hookrightarrow H^{s,\varphi}(\mathbb{R}^{n}).
$$
Отсюда на основании теоремы 2.2.3 из упомянутой монографии [\ref{Hermander65}, с.
56] получаем, что
$\varphi(\langle\xi\rangle)/\varphi_{1}(\langle\xi\rangle)\rightarrow0$ при
$|\xi|\rightarrow\infty$. Значит, $\varphi(t)/\varphi_{1}(t)\rightarrow0$ при
$t\rightarrow\infty$.

Утверждение (iv) доказано.

(v) В соболевском случае $\varphi\equiv1$ утверждение (v) известно [\ref{Shubin78},
с.~70] (теорема~7.7). Так, соболевские пространства $H^{s\pm1}(\Gamma)$ и
$H^{-s\mp1}(\Gamma)$ взаимно двойственны относительно расширения по непрерывности
скалярного произведения в $L_{2}(\Gamma)$. Это означает, что линейное отображение
$Q:w\mapsto(w,\,\cdot)_{\Gamma}$, где $w\in C^{\infty}(\Gamma)$, продолжается по
непрерывности до гомеоморфизмов
$Q:H^{s\mp1}(\Gamma)\leftrightarrow(H^{-s\pm1}(\Gamma))'$. Применив интерполяцию с
параметром $\psi$ из теоремы \ref{th2.21}, где $\varepsilon=\delta=1$, получаем еще
один гомеоморфизм
\begin{equation}\label{2.73}
Q:\left[H^{s-1}(\Gamma),H^{s+1}(\Gamma)\right]_{\psi}\leftrightarrow
\left[(H^{-s+1}(\Gamma))',(H^{-s-1}(\Gamma))'\right]_{\psi}.
\end{equation}
Здесь левое интерполяционное пространство равно $H^{s,\varphi}(\Gamma)$, а правое в
силу теоремы \ref{th2.4} представляется в виде
\begin{gather*}
\left[(H^{-s+1}(\Gamma))',(H^{-s-1}(\Gamma))'\right]_{\psi}=\\
\left[H^{-s-1}(\Gamma),H^{-s+1}(\Gamma)\right]_{\chi}'= (H^{-s,1/\varphi}(\Gamma))'.
\end{gather*}
Отметим, что последнее равенство справедливо, поскольку
$\chi(t):=t/\psi(t)=t^{1/2}/\varphi(t^{1/2})$ при $t\geq1$. Таким образом, формула
\eqref{2.73} влечет за собой гомеоморфизм
$Q:H^{s,\varphi}(\Gamma)\leftrightarrow(H^{-s,1/\varphi}(\Gamma))'$, т.~е. указанную
взаимную двойственность пространств $H^{s,\varphi}(\Gamma)$ и
$H^{-s,1/\varphi}(\Gamma)$. Пункт (v) доказан.

Теорема \ref{th2.22} доказана.

\medskip

Из теоремы \ref{th2.21} и свойства реитерации (теорема \ref{th2.3}) вытекает
следующее утверждение, показывающее, что уточненная шкала на $\Gamma$ замкнута
относительно интерполяции с функциональными параметрами, (квази)правильно
меняющимися на~$\infty$.


\begin{theorem}\label{th2.23}
Пусть $s_{0},s_{1}\in\mathbb{R}$, $s_{0}\leq s_{1}$ и
$\varphi_{0},\varphi_{1}\in\mathcal{M}$. В~случае $s_{0}=s_{1}$ предположим, что
функция $\varphi_{0}/\varphi_{1}$ ограничена в окрестности $\infty$. Пусть функция
$\psi\in\mathcal{B}$ квазиправильно меняющаяся на $\infty$ порядка $\theta$, где
$0<\theta<1$. В силу теоремы $\ref{th2.11}$, $\psi$ --- интерполяционный параметр.
Представим его в виде $\psi(t)=t^{\theta}\chi(t)$, где $\chi\in\mathrm{QSV}$.
Положим $s:=(1-\theta)s_{0}+\theta s_{1}$ и
\begin{equation}\label{2.74}
\varphi(t):=\varphi_{0}^{1-\theta}(t)\,\varphi_{1}^{\theta}(t)\,
\chi\left(t^{s_{1}-s_{0}}\varphi_{1}(t)/\varphi_{0}(t)\right)\quad\mbox{при}\quad
t\geq1.
\end{equation}
Тогда $\varphi\in\mathcal{M}$ и
\begin{equation}\label{2.75}
\left[\,H^{s_{0},\varphi_{0}}(\Gamma),
H^{s_{1},\varphi_{1}}(\Gamma)\,\right]_{\psi}= H^{s,\varphi}(\Gamma)
\end{equation}
с эквивалентностью норм.
\end{theorem}


\textbf{Доказательство.} Положительная функция $\varphi$ измерима по Борелю на
множестве $[1,\infty)$ и ограничена вместе с функцией $1/\varphi$ на каждом отрезке
$[1,b]$, где $1<b<\infty$, поскольку аналогичными свойствами обладают функции
$\varphi_{0},\varphi_{1}$ и $\chi$. Кроме того, условия
$\varphi_{0},\varphi_{1},\chi\in\mathrm{QSV}$ влекут за собой включение
$\varphi\in\mathrm{QSV}$ в силу теоремы \ref{th2.12} (iii), (iv) в случае
$s_{0}<s_{1}$, и теоремы \ref{th2.13}~--- в случае $s_{0}=s_{1}$. Таким образом,
$\varphi\in\mathcal{M}$. Далее, согласно теореме \ref{th2.22} (iii), (iv) пара
пространств $[H^{s_{0},\varphi_{0}}(\Gamma), H^{s_{1},\varphi_{1}}(\Gamma)]$
допустимая. Докажем равенство \eqref{2.75}.

Положим $\varrho:=s_{1}-s_{0}+1$ и $\varepsilon_{j}:=s_{j}-s+\varrho$,
$\delta_{j}:=s-s_{j}+\varrho$ при $j\in\{0,1\}$. Числа
$\varrho,\varepsilon_{j},\delta_{j}$ положительны, поскольку $s_{0}\leq s\leq
s_{1}$. Отметим следующие их свойства:
\begin{equation}\label{2.76}
\varepsilon_{j}+\delta_{j}=2\varrho,\quad
\varepsilon_{1}-\varepsilon_{0}=s_{1}-s_{0},
\quad(1-\theta)\varepsilon_{0}+\theta\varepsilon_{1}=\varrho.
\end{equation}
В силу теоремы \ref{th2.21} имеем при $j\in\{0,1\}$ равенства
$$
\bigl[H^{s_{j}-\varepsilon_{j}}(\Gamma),\
H^{s_{j}+\delta_{j}}(\Gamma)\bigr]_{\psi_{j}}= H^{s_{j},\varphi_{j}}(\Gamma)
$$
с эквивалентностью норм. Здесь интерполяционный параметр $\psi_{j}$ определяется по
формуле
\begin{equation}\label{2.77}
\psi_{j}(t):=
\begin{cases}
\;\;t^{\,\varepsilon_{j}/(\varepsilon_{j}+\delta_{j})}\,
\varphi_{j}(t^{1/(\varepsilon_{j}+\delta_{j})}) & \text{при}\;\;t\geq1, \\
\;\;\varphi_{j}(1) & \text{при}\;\;0<t<1.
\end{cases}
\end{equation}
Отсюда, поскольку $\psi$~--- интерполяционный параметр и
$s_{j}-\varepsilon_{j}=s-\varrho$, $s_{j}+\delta_{j}=s+\varrho$, можем записать
следующие равенства пространств с точностью до эквивалентности норм в них:
\begin{gather}\notag
\bigl[\,H^{s_{0},\varphi_{0}}(\Gamma),
H^{s_{1},\varphi_{1}}(\Gamma)\,\bigr]_{\psi}= \\
\left[\,\bigl[H^{s-\varrho}(\Gamma),H^{s+\varrho}(\Gamma)\bigr]_{\psi_{0}},\,
\bigl[H^{s-\varrho}(\Gamma),H^{s+\varrho}(\Gamma)\bigr]_{\psi_{1}}\, \right]_{\psi}.
\label{2.78}
\end{gather}
Заметим (см. формулу \eqref{2.76}), что функция
$$
\psi_{0}(t)/\psi_{1}(t)=t^{(s_{0}-s_{1})/(2\varrho)}\varphi_{0}(t^{1/(2\varrho)})/
\varphi_{1}(t^{1/(2\varrho)}),\quad t\geq1,
$$
ограничена в окрестности $\infty$ в силу теоремы \ref{th2.12} (ii) в случае
$s_{0}<s_{1}$, а также по условию в случае $s_{0}=s_{1}$. Применив к \eqref{2.78}
теорему \ref{th2.3} о свойстве реитерации, можем записать следующее равенство
пространств с точностью до эквивалентности норм в них:
\begin{equation}\label{2.79}
\left[\,H^{s_{0},\varphi_{0}}(\Gamma),
H^{s_{1},\varphi_{1}}(\Gamma)\,\right]_{\psi}= \left[\,H^{s-\varrho}(\Gamma),
H^{s+\varrho}(\Gamma)\,\right]_{\omega}.
\end{equation}
Здесь интерполяционный параметр $\omega$ определен по формуле
$$
\omega(t):=\psi_{0}(t)\,\psi(\psi_{1}(t)/\psi_{0}(t))\quad\mbox{при}\;\;t>0.
$$
В силу соотношений \eqref{2.76}, \eqref{2.77} и \eqref{2.74} после элементарных
вычислений получаем равенства: $\omega(t)=t^{1/2}\varphi(t^{1/(2\varrho)})$ при
$t\geq1$ и $\omega(t)=\varphi(1)$ при $0<t<1$. Поэтому согласно теореме \ref{th2.21}
имеем
\begin{equation}\label{2.80}
\bigl[\,H^{s-\varrho}(\Gamma), H^{s+\varrho}(\Gamma)\,\bigr]_{\omega}=
H^{s,\varphi}(\Gamma)
\end{equation}
с эквивалентностью норм. Из формул \eqref{2.79} и \eqref{2.80} следует \eqref{2.75}.

Теорема \ref{th2.23} доказана.


\begin{remark}
Теорема \ref{th2.23} верна в предельных случаях $\theta=0$ и $\theta=1$, если
дополнительно предположить, что функция $\psi$ псевдовогнута в окрестности~$\infty$.
Тогда в силу теоремы \ref{th2.9}, $\psi$ --- интерполяционный параметр, и проходит
приведенное выше доказательство. Например, теорема \ref{th2.23} остается
справедливой для каждой из функций $\psi(t):=\ln^{r}t$ и $\psi(t):=t/\ln^{r} t$, где
$t\gg1$ и $r>0$.
\end{remark}


\begin{remark}
Теорема \ref{th2.23} остается в силе, если в ее формулировке заменить $\Gamma$ на
$\mathbb{R}^{n}$. Доказательство аналогично приведенному выше. Таким образом,
уточненная шкала в $\mathbb{R}^{n}$ (как и на $\Gamma$) замкнута относительно
интерполяции с функциональными параметрами, квазиправильно меняющимися на~$\infty$.
\end{remark}

\subsection{Эквивалентные нормы}\label{sec2.5.3}

В этом пункте мы построим эквивалентные нормы в пространстве $H^{s,\varphi}(\Gamma)$
при помощи эллиптических положительно определенных ПДО. Как следствие, мы докажем,
что определения \ref{def2.13} и \ref{def2.15} эквивалентны.

Пусть на многообразии $\Gamma$ задан полиоднородный ПДО $A$ порядка $m>0$,
эллиптический на $\Gamma$. (Определение ПДО на замкнутом многообразии и связанные с
ним понятия приведены в п.~\ref{sec2.6.1b}.) Предположим также, что оператор
$A:C^{\infty}(\Gamma)\rightarrow C^{\infty}(\Gamma)$ положительно определенный в
пространстве $L_{2}(\Gamma)$, т. е. существует число $\varkappa>0$ такое, что
\begin{equation}\label{2.81}
(Au,u)_{\Gamma}\geq\varkappa\,(u,u)_{\Gamma}\quad\mbox{для любого}\quad u\in
C^{\infty}(\Gamma).
\end{equation}

Обозначим через $A_{0}$ замыкание в пространстве $L_{2}(\Gamma)$ оператора
$A:C^{\infty}(\Gamma)\rightarrow C^{\infty}(\Gamma)$. Поскольку ПДО $A$
эллиптический на $\Gamma$, это замыкание существует и определено на $H^{m}(\Gamma)$
[\ref{Agranovich90}, с.~28; \ref{Shubin78}, с.~77]. Из условия \eqref{2.81}
вытекает, что ПДО $A$ формально самосопряженный. Следовательно [\ref{Agranovich90},
с.~29; \ref{Shubin78}, с.~78], $A_{0}$ --- неограниченный самосопряженный оператор в
пространстве $L_{2}(\Gamma)$, причем
$\mathrm{Spec}\,A_{0}\subseteq[\varkappa,\infty)$.

Пусть $s\in\mathbb{R}$ и $\varphi\in\mathcal{M}$. Положим
\begin{equation}\label{2.82}
\varphi_{s,m}(t):=
\begin{cases}
\;t^{s/m}\varphi(t^{1/m}) & \text{при}\;\;t\geq1, \\
\;\varphi(1) & \text{при}\;\;0<t<1.
\end{cases}
\end{equation}
Поскольку функция $\varphi_{s,m}$ положительна и измерима по Борелю на полуоси
$(0,\infty)$, то в пространстве $L_{2}(\Gamma)$ определен, как функция от $A_{0}$,
самосопряженный оператор $\varphi_{s,m}(A_{0})$.


\begin{lemma}\label{lem2.9}
Утверждается следующее:
\begin{itemize}
\item [$\mathrm{(i)}$] область определения оператора $\varphi_{s,m}(A_{0})$ содержит
$C^{\infty}(\Gamma)$;
\item [$\mathrm{(ii)}$] отображение
\begin{equation}\label{2.83}
f\,\mapsto\,\|\varphi_{s,m}(A_{0})f\|_{L_{2}(\Gamma)},\quad f\in C^{\infty}(\Gamma),
\end{equation}
является нормой в пространстве $C^{\infty}(\Gamma)$.
\end{itemize}
\end{lemma}


\textbf{Доказательство.} (i) Выберем натуральное число $k>s/m$. Поскольку
$\varphi\in\mathcal{M}$, то функция $\varphi_{s,m}$ ограничена на каждом компактном
подмножестве полуоси $(0,\infty)$, и $t^{-k}\varphi_{s,m}(t)\rightarrow0$ при
$t\rightarrow\infty$ в силу теоремы \ref{th2.12} (ii), (iv). Следовательно,
существует число $c>0$ такое, что $\varphi_{s,m}(t)\leq c\,t^{k}$ при
$t\geq\nobreak\varkappa$. Рассмотрим неограниченный оператор $A_{0}^{k}$ в
пространстве $L_{2}(\Gamma)$. Поскольку $A:C^{\infty}(\Gamma)\rightarrow
C^{\infty}(\Gamma)$, то можем записать
$$
C^{\infty}(\Gamma)\subset\mathrm{Dom}\,A_{0}^{k}\subset
\mathrm{Dom}\,\varphi_{s,m}(A_{0}).
$$
Утверждение (i) доказано.

(ii) Согласно утверждению (i) отображение \eqref{2.83} корректно определено. Для
этого отображения все свойства нормы очевидны за исключением свойства положительной
определенности. Докажем его. Для произвольной функции $f\in C^{\infty}(\Gamma)$ на
основании спектральной теоремы запишем:
\begin{gather}\label{2.84}
\|\varphi_{s,m}(A_{0})f\|_{L_{2}(\Gamma)}^{2}=
\int\limits_{\varkappa}^{\infty}\varphi_{s,m}^{2}(t)\,d(E_{t}f,f)_{\Gamma}, \\
\|f\|_{L_{2}(\Gamma)}^{2}= \int\limits_{\varkappa}^{\infty}d(E_{t}f,f)_{\Gamma}.
\label{2.85}
\end{gather}
Здесь $E_{t}$, $t\geq\varkappa$,~--- разложение единицы в $L_{2}(\Gamma)$,
соответствующее самосопряженному оператору $A_{0}$. Теперь, если
$\|\varphi_{s,m}(A_{0})f\|_{L_{2}(\Gamma)}^{2}=\nobreak0$, то из \eqref{2.84} и
положительности функции $\varphi_{s,m}$ следует, что мера $(E(\cdot)f,f)_{\Gamma}$
множества $[\varkappa,\infty)$ равна~0. Отсюда в силу \eqref{2.85} получаем
равенство $f=0$ на $\Gamma$. Утверждение (ii) доказано.

Лемма~\ref{lem2.9} доказана.


\begin{theorem}\label{th2.24}
Для произвольных $s\in\mathbb{R}$ и $\varphi\in\mathcal{M}$ норма \eqref{2.83} и
норма в пространстве $H^{s,\varphi}(\Gamma)$ эквивалентны на $C^{\infty}(\Gamma)$.
Тем самым, пространство $H^{s,\varphi}(\Gamma)$ совпадает с точностью до
эквивалентности норм с пополнением линейной системы $C^{\infty}(\Gamma)$ по норме
\eqref{2.83}.
\end{theorem}


\textbf{Доказательство.} Сначала предположим, что $s>0$. Выберем $k\in\mathbb{N}$
такое, что $k\,m>s$. Поскольку оператор $A_{0}^{k}$ замкнут и положительно определен
в $L_{2}(\Gamma)$, его область определения $\mathrm{Dom}\,A_{0}^{k}$ является
гильбертовым пространством относительно скалярного произведения
$(A_{0}^{k}f,A_{0}^{k}g)_{\Gamma}$ функций $f,g$. При этом пара пространств
$[L_{2}(\Gamma),\mathrm{Dom}\,A_{0}^{k}]$ допустимая и оператор $A_{0}^{k}$
порождающий для нее. Кроме того, так как $A_{0}^{k}$
--- замыкание в $L_{2}(\Gamma)$ эллиптического ПДО $A^{k}$, пространства
$\mathrm{Dom}\,A_{0}^{k}$ и $H^{km}(\Gamma)$ равны с точностью до эквивалентности
норм. Пусть функция $\psi$ --- интерполяционный параметр из теорем \ref{th2.14},
\ref{th2.21}, где берем $\varepsilon=s$ и $\delta=k\,m-s$. Тогда
$\psi(t^{k})=\varphi_{s,m}(t)$ при $t>0$, и в силу теоремы \ref{th2.21} имеем:
\begin{gather*}
\|f\|_{H^{s,\varphi}(\Gamma)}\asymp
\|f\|_{[H^{0}(\Gamma),\,H^{km}(\Gamma)]_{\psi}}\asymp
\|f\|_{[L_{2}(\Gamma),\,\mathrm{Dom}\,A_{0}^{k\,}]_{\psi}}=\\
\|\psi(A_{0}^{k})f\|_{L_{2}(\Gamma)}= \|\varphi_{s,m}(A_{0})f\|_{L_{2}(\Gamma)},
\end{gather*}
где $f\in C^{\infty}(\Gamma)$.

Пусть теперь вещественное число $s$ произвольное. Выберем $k\in\mathbb{N}$ такое,
что $s+k\,m>0$. По доказанному
\begin{equation}\label{2.86}
\bigl\|g\bigr\|_{H^{s+km,\varphi}(\Gamma)}\asymp
\bigl\|\varphi_{s+km,m}(A_{0})\,g\bigr\|_{L_{2}(\Gamma)},\quad g\in
C^{\infty}(\Gamma).
\end{equation}
ПДО $A^{k}$ осуществляет гомеоморфизмы
\begin{gather}\label{2.87}
A^{k}:\,H^{\sigma+km}(\Gamma)\leftrightarrow H^{\sigma}(\Gamma),\quad\sigma\in\mathbb{R},\\
A^{k}:\,H^{s+km,\varphi}(\Gamma)\leftrightarrow H^{s,\varphi}(\Gamma).\label{2.88}
\end{gather}
Это мы докажем в следующем абзаце. Обозначим через $A^{-k}$ оператор, обратный к
$A^{k}$. Для произвольной функции $f\in C^{\infty}(\Gamma)$ имеем в силу
\eqref{2.87}:
$$
g:=A^{-k}f\in\bigcap_{\sigma\in\mathbb{R}}H^{\sigma+km}(\Gamma)=C^{\infty}(\Gamma)
$$
и $A_{0}^{k}A^{-k}f=f$. Отсюда в силу \eqref{2.88} и \eqref{2.86} получаем нужную
нам эквивалентность норм:
\begin{gather*}
\bigl\|f\bigr\|_{H^{s,\varphi}(\Gamma)}\asymp
\bigl\|A^{-k}f\bigr\|_{H^{s+km,\varphi}(\Gamma)}\asymp
\bigl\|\varphi_{s+km,m}(A_{0})A^{-k}f\bigr\|_{L_{2}(\Gamma)}=\\
\bigl\|\varphi_{s,m}(A_{0})A_{0}^{k}A^{-k}f\bigr\|_{L_{2}(\Gamma)}=
\bigl\|\varphi_{s,m}(A_{0})f\bigr\|_{L_{2}(\Gamma)},\quad f\in C^{\infty}(\Gamma).
\end{gather*}

Остается доказать, что $A^{k}$ осуществляет гомеоморфизмы \eqref{2.87} и
\eqref{2.88}. Поскольку ПДО $A$ эллиптический на $\Gamma$, он определяет
ограниченный нетеров оператор $A:H^{\sigma+m}(\Gamma)\rightarrow H^{\sigma}(\Gamma)$
при любом $\sigma\in\mathbb{R}$. Ядро этого оператора и индекс не зависят от
$\sigma$ (см., например, [\ref{Hermander87}, с.~262], теорема 19.2.1, либо
[\ref{Shubin78}, с.~76], теорема 8.1). Поскольку $0\notin\mathrm{Spec}\,A_{0}$,
самосопряженный оператор $A_{0}$ гомеоморфно отображает свою область определения
$H^{m}(\Gamma)$ на пространство $H^{0}(\Gamma)=L_{2}(\Gamma)$. Поэтому ядро
тривиально, а индекс равен нулю. Значит, ПДО $A$ устанавливает гомеоморфизмы
$$
A:H^{\sigma+m}(\Gamma)\leftrightarrow H^{\sigma}(\Gamma),\quad\sigma\in\mathbb{R}.
$$
Отсюда в результате $k$ итераций следуют гомеоморфизмы \eqref{2.87}. Наконец, взяв
$\sigma=s\mp1$ в \eqref{2.87} и воспользовавшись интерполяционной теоремой
\ref{th2.21}, получаем гомеоморфизм \eqref{2.88}.

Теорема~\ref{th2.24} доказана.

\medskip

Теорема~\ref{th2.24} в случае $A_{0}=1-\Delta_{\Gamma}$ дает следующий важный
результат.


\begin{corollary}\label{cor2.3}
Для произвольно заданных $s\in\mathbb{R}$ и $\varphi\in\mathcal{M}$ определения
$\ref{def2.13}$ и $\ref{def2.15}$ являются эквивалентными.
\end{corollary}


Выделим случай, когда пространство $H^{s,\varphi}(\Gamma)$ совпадает с областью
определения оператора $\varphi_{s,m}(A_{0})$.


\begin{theorem}\label{th2.25}
Пусть $s\geq0$ и $\varphi\in\mathcal{M}$. В случае $s=0$ предположим, что функция
$1/\varphi$ ограничена в окрестности $\infty$. Тогда пространство
$H^{s,\varphi}(\Gamma)$ совпадает с областью определения оператора
$\varphi_{s,m}(A_{0})$, причем норма в пространстве $H^{s,\varphi}(\Gamma)$
эквивалентна норме графика оператора $\varphi_{s,m}(A_{0})$.
\end{theorem}


\textbf{Доказательство.}  Область определения $\mathrm{Dom}\,\varphi_{s,m}(A_{0})$
замкнутого оператора $\varphi_{s,m}(A_{0})$ является гильбертовым пространством
относительно скалярного произведения графика этого оператора. Докажем, что норма
графика оператора $\varphi_{s,m}(A_{0})$ эквивалентна норме \eqref{2.83} на линейном
многообразии $C^{\infty}(\Gamma)$, и что это многообразие плотно в пространстве
$\mathrm{Dom}\,\varphi_{s}(A_{0})$. Отсюда в силу теоремы~\ref{th2.24} следует
теорема~\ref{th2.25}.

Согласно условию и теореме~\ref{th2.12} (ii) существует число $c>0$ такое, что
$\varphi_{s,m}(t)\geq c$ при $t>0$. Поэтому
$$
\bigl\|\varphi_{s,m}(A_{0})f\bigr\|_{L_{2}(\Gamma)}\geq c\,
\bigl\|f\bigr\|_{L_{2}(\Gamma)}\quad\mbox{для любого}\quad f\in C^{\infty}(\Gamma).
$$
Отсюда вытекает указанная эквивалентность норм. Остается доказать плотность
множества $C^{\infty}(\Gamma)$ в пространстве $\mathrm{Dom}\,\varphi_{s,m}(A_{0})$.

Пусть $f\in\mathrm{Dom}\,\varphi_{s,m}(A_{0})$. Поскольку $\varphi_{s,m}(A_{0})f\in
L_{2}(\Gamma)$, существует последовательность функций $h_{j}\in C^{\infty}(\Gamma)$
такая, что $h_{j}\rightarrow\varphi_{s,m}(A_{0})f$ в $L_{2}(\Gamma)$ при
$j\rightarrow\infty$. Так как $1/\varphi_{s,m}(t)\leq 1/c$ при $t>0$, оператор
$\varphi_{s,m}^{-1}(A_{0})$ ограничен в пространстве $L_{2}(\Gamma)$. Следовательно,
\begin{gather*}
f_{j}:=\varphi_{s,m}^{-1}(A_{0})h_{j}\rightarrow f,\\
\varphi_{s,m}(A_{0})f_{j}=h_{j}\rightarrow\varphi_{s,m}(A_{0})f\;\;\mbox{в}\;\;L_{2}(\Gamma)
\;\;\mbox{при}\;\;j\rightarrow\infty.
\end{gather*}
Иными словами, $f_{j}\rightarrow f$ по норме графика оператора
$\varphi_{s,m}(A_{0})$. Кроме того, поскольку $h_{j}\in C^{\infty}(\Gamma)$, то для
каждого $k\in\mathbb{N}$ справедливо
$$
f_{j}=A_{0}^{-k}\varphi_{s,m}^{-1}(A_{0})A_{0}^{k}\,h_{j}\in H^{km}(\Gamma).
$$
Следовательно, $f_{j}\in C^{\infty}(\Gamma)$ и плотность множества
$C^{\infty}(\Gamma)$ в пространстве $\mathrm{Dom}\,\varphi_{s,m}(A_{0})$
установлена.

Теорема \ref{th2.25} доказана.

\medskip

В заключение этого пункта дадим эквивалентную норму в пространстве
$H^{s,\varphi}(\Gamma)$, выраженную в терминах последовательностей.

Поскольку ПДО $A$ эллиптический на $\Gamma$ и положительно определенный в
$L_{2}(\Gamma)$, то он имеет следующие спектральные свойства (см., например,
[\ref{Agranovich90}, с. 95, 98], [\ref{Taylor85}, с. 309, 322] или [\ref{Shubin78},
с. 78, 132]). В пространстве $L_{2}(\Gamma)$ существует ортонормированный базис
$(h_{j})_{j=1}^{\infty}$, образованный собственными функциями $h_{j}$ оператора
$A_{0}$. При этом $h_{j}\in C^{\infty}(\Gamma)$, $Ah_{j}=\lambda_{j}h_{j}$, и
последовательность собственных значений $(\lambda_{j})_{j=1}^{\infty}$ положительна,
монотонно не убывает и стремится к $\infty$. Спектр оператора $A_{0}$ совпадает с
множеством всех его собственных значений $\{\lambda_{j}:j\in\mathbb{N}\}$.
Выполняется асимптотическая формула
\begin{equation}\label{2.89}
\lambda_{j}\sim c\,j^{\,m/n}\quad\mbox{при}\quad j\rightarrow\infty,
\end{equation}
где $c$ --- некоторое положительное число, зависящее от $A$ и $n=\dim\Gamma$.
Произвольное распределение $f\in\mathcal{D}'(\Gamma)$ разлагается в ряд Фурье
\begin{equation}\label{2.90}
f=\sum_{j=1}^{\infty}\;c_{j}(f)\,h_{j}\quad\mbox{(сходимость
в}\;\;\mathcal{D}'(\Gamma)),
\end{equation}
где $c_{j}(f):=(f,h_{j})_{\Gamma}$ --- коэффициенты Фурье распределения $f$ в базисе
$(h_{j})_{j=1}^{\infty}$.


\begin{theorem}\label{th2.26}
Пусть $s\in\mathbb{R}$ и $\varphi\in\mathcal{M}$. Тогда
\begin{gather}\label{2.91}
H^{s,\varphi}(\Gamma)=\Bigl\{f\in\mathcal{D}'(\Gamma):\,
\sum_{j=1}^{\infty}\;j^{\,2s/n}\,\varphi^{2}(j^{\,1/n})\,|c_{j}(f)|^{2}<\infty\Bigr\},\\
\|f\|_{H^{s,\varphi}(\Gamma)}\asymp\Bigl(\;\sum_{j=1}^{\infty}\;
j^{\,2s/n}\,\varphi^{2}(j^{\,1/n})\,|c_{j}(f)|^{2}\;\Bigr)^{1/2}. \label{2.92}
\end{gather}
\end{theorem}


Доказательству этой теоремы предпошлем три леммы. В них, как и прежде,
$s\in\mathbb{R}$ и $\varphi\in\mathcal{M}$.


\begin{lemma}\label{lem2.10}
Для любой функции $f\in\mathrm{Dom}\,(\varphi_{s,m}(A_{0}))$ выполняется равенство
\begin{equation}\label{2.93}
\|\varphi_{s,m}(A_{0})f\|_{L_{2}(\Gamma)}^{2}=
\sum_{j=1}^{\infty}\;\varphi_{s,m}^{2}(\lambda_{j})\,|c_{j}(f)|^{2}.
\end{equation}
\end{lemma}


\textbf{Доказательство.} Обозначим через $P_{j}$ ортопроектор в пространстве
$L_{2}(\Gamma)$ на одномерное подпространство $\{lh_{j}:l\in\mathbb{C}\}$. Тогда
$P_{j}f=c_{j}(f)\,h_{j}$ для любого $f\in L_{2}(\Gamma)$. Поскольку
$\mathrm{Spec}\,A_{0}=\{\lambda_{j}:j\in\mathbb{N}\}$, то получаем в силу
спектральной теоремы, что для любого $f\in \mathrm{Dom}\,(\varphi_{s,m}(A_{0}))$
справедливы равенства
$$
\varphi_{s,m}(A_{0})\,f=\sum_{j=1}^{\infty}\,\varphi_{s,m}(\lambda_{j})P_{j}f=
\sum_{j=1}^{\infty}\,\varphi_{s,m}(\lambda_{j})c_{j}(f)\,h_{j}.
$$
Здесь ряды сходятся в $L_{2}(\Gamma)$. Отсюда на основании равенства Парсеваля
получаем \eqref{2.93}. Лемма~\ref{lem2.10} доказана.


\begin{lemma}\label{lem2.11}
Для произвольного $f\in\mathcal{D}'(\Gamma)$ включение $f\in H^{s,\varphi}(\Gamma)$
равносильно неравенству
\begin{equation}\label{2.94}
\sum_{j=1}^{\infty}\;\varphi_{s,m}^{2}(\lambda_{j})\,|c_{j}(f)|^{2}<\infty.
\end{equation}
\end{lemma}


\textbf{Доказательство.} Пусть $f\in H^{s,\varphi}(\Gamma)$. Докажем неравенство
\eqref{2.94}. Выберем последовательность $(f_{k})\subset C^{\infty}(\Gamma)$ такую,
что $f_{k}\rightarrow f$ в $H^{s,\varphi}(\Gamma)$ при $k\rightarrow\infty$. На
основании леммы~\ref{lem2.10} и теоремы \ref{th2.24} имеем:
\begin{gather*}
\sum_{j=1}^{\infty}\;\varphi_{s,m}^{2}(\lambda_{j})\,|c_{j}(f_{k})|^{2}=
\|\varphi_{s,m}(A_{0})f_{k}\|_{L_{2}(\Gamma)}^{2}\asymp\\
\|f_{k}\|_{H^{s,\varphi}(\Gamma)}^{2}\leq\bigl(1+\|f\|_{H^{s,\varphi}(\Gamma)}^{2}\bigr),\quad
k\geq1.
\end{gather*}
Кроме того, из непрерывности вложения
$H^{s,\varphi}(\Gamma)\hookrightarrow\mathcal{D}'(\Gamma)$ вытекает, что
$$
c_{j}(f_{k})=(f_{k},h_{j})_{\Gamma}\rightarrow(f,h_{j})_{\Gamma}=c_{j}(f)
\quad\mbox{при}\;\;k\rightarrow\infty
$$
для каждого $j\in\mathbb{N}$. Поэтому в силу леммы Фату для положительных рядов
(см., например, [\ref{BerezanskyUsSheftel90}, с. 107])
$$
\sum_{j=1}^{\infty}\;\varphi_{s,m}^{2}(\lambda_{j})\,|c_{j}(f)|^{2}\leq
\varliminf_{k\rightarrow\infty}\,
\sum_{j=1}^{\infty}\;\varphi_{s,m}^{2}(\lambda_{j})\,|c_{j}(f_{k})|^{2}<\infty.
$$
Таким образом, включение $f\in H^{s,\varphi}(\Gamma)$ влечет неравенство
\eqref{2.94}.

Пусть теперь выполняется неравенство \eqref{2.94}. Докажем включение $f\in
H^{s,\varphi}(\Gamma)$. В силу \eqref{2.94} в пространстве $L_{2}(\Gamma)$ сходится
ортогональный ряд
\begin{equation}\label{2.95}
\sum_{j=1}^{\infty}\;\varphi_{s,m}(\lambda_{j})\,c_{j}(f)\,h_{j}=:h\in
L_{2}(\Gamma).
\end{equation}
Рассмотрим его частную сумму
$$
w_{k}:=\sum_{j=1}^{k}\;\varphi_{s,m}(\lambda_{j})\,c_{j}(f)\,h_{j}\in
C^{\infty}(\Gamma).
$$
Ввиду \eqref{2.95},
\begin{equation}\label{2.96}
w_{k}\rightarrow h\quad\mbox{в}\quad L_{2}(\Gamma)\quad\mbox{при}\quad
k\rightarrow\infty.
\end{equation}
Отсюда вытекает, что последовательность
$(\varphi_{s,m}^{-1}(A_{0})\,w_{k})_{k=1}^{\infty}$ фундаментальная в пространстве
$H^{s,\varphi}(\Gamma)$. В самом деле, поскольку $A_{0}h_{j}=\lambda_{j}h_{j}$, то
$\varphi_{s,m}^{-1}(A_{0})\,h_{j}=\varphi_{s,m}^{-1}(\lambda_{j})\,h_{j}$, и
\begin{gather*}
\varphi_{s,m}^{-1}(A_{0})\,w_{k}=
\sum_{j=1}^{k}\;\varphi_{s,m}(\lambda_{j})\,c_{j}(f)\,\varphi_{s,m}^{-1}(A_{0})h_{j}= \\
\sum_{j=1}^{k}\;c_{j}(f)\,h_{j}\in C^{\infty}(\Gamma).
\end{gather*}
Поэтому в силу теоремы \ref{th2.24} и формулы \eqref{2.96}
\begin{gather*}
\|\varphi_{s,m}^{-1}(A_{0})\,w_{k}-
\varphi_{s,m}^{-1}(A_{0})\,w_{p}\|_{H^{s,\varphi}(\Gamma)}\asymp\\
\|w_{k}-w_{p}\|_{L_{2}(\Gamma)}\rightarrow0\quad\mbox{при}\quad
k,p\rightarrow\infty,
\end{gather*}
т. е. последовательность $(\varphi_{s,m}^{-1}(A_{0})\,w_{k})_{k=1}^{\infty}$
фундаментальная в пространстве $H^{s,\varphi}(\Gamma)$. Обозначим ее предел через
$g$. Имеем:
\begin{equation}\label{2.97}
g=\lim_{k\rightarrow\infty}\;\varphi_{s,m}^{-1}(A_{0})\,w_{k}=
\sum_{j=1}^{\infty}\;c_{j}(f)\,h_{j}\quad\mbox{в}\quad H^{s,\varphi}(\Gamma).
\end{equation}
Отсюда ввиду равенства \eqref{2.90} получаем, что $f=g\in H^{s,\varphi}(\Gamma)$.
Таким образом, неравенство \eqref{2.94} влечет за собой включение $f\in\nobreak
H^{s,\varphi}(\Gamma)$.

Лемма \ref{lem2.11} доказана.


\begin{remark}\label{rem2.12}
Из формулы \eqref{2.97} следует, что для любого распределения $f\in
H^{s,\varphi}(\Gamma)$ ряд \eqref{2.90} сходится к $f$ в пространстве
$H^{s,\varphi}(\Gamma)$
\end{remark}


\begin{lemma}\label{lem2.12}
Справедлива эквивалентность норм
\begin{equation}\label{2.98}
\|f\|_{H^{s,\varphi}(\Gamma)}\asymp
\Bigl(\;\sum_{j=1}^{\infty}\;\varphi_{s,m}^{2}(\lambda_{j})\,|c_{j}(f)|^{2}\;\Bigr)^{1/2},
\quad f\in H^{s,\varphi}(\Gamma).
\end{equation}
\end{lemma}


\textbf{Доказательство.} Воспользуемся рассуждениями и обозначениями, сделанными в
доказательстве леммы~\ref{lem2.11}. Пусть $f\in H^{s,\varphi}(\Gamma)$, а $g$
определено по формуле \eqref{2.97}. Тогда $f=g$, и
\begin{equation}\label{2.99}
\|f\|_{H^{s,\varphi}(\Gamma)}=
\lim_{k\rightarrow\infty}\;\|\varphi_{s,m}^{-1}(A_{0})\,w_{k}\|_{H^{s,\varphi}(\Gamma)}.
\end{equation}
Здесь, напомним, $\varphi_{s,m}^{-1}(A_{0})\,w_{k}\in C^{\infty}(\Gamma)$. Поэтому
согласно теореме \ref{th2.24} существует такое число $c\geq1$, не зависящее от $f$,
что
\begin{gather}\notag
c^{-1}\,\|w_{k}\|_{L_{2}(\Gamma)}\leq
\|\varphi_{s,m}^{-1}(A_{0})\,w_{k}\|_{H^{s,\varphi}(\Gamma)}\leq \\
c\,\|w_{k}\|_{L_{2}(\Gamma)}\quad\mbox{при}\quad k\geq1.\label{2.100}
\end{gather}
В силу \eqref{2.95} и \eqref{2.96} имеем:
\begin{equation}\label{2.101}
\lim_{k\rightarrow\infty}\;\|w_{k}\|_{L_{2}(\Gamma)}=\|h\|_{L_{2}(\Gamma)}=
\Bigl(\;\sum_{j=1}^{\infty}\;\varphi_{s,m}^{2}(\lambda_{j})\,|c_{j}(f)|^{2}\;\Bigr)^{1/2}.
\end{equation}
Остается перейти в неравенстве \eqref{2.100} к пределу при $k\rightarrow\infty$ и
воспользоваться равенствами \eqref{2.99}, \eqref{2.101}. Получаем эквивалентность
норм \eqref{2.98}. Лемма~\ref{lem2.12} доказана.

\medskip

\textbf{Доказательство теоремы~\ref{th2.26}.} Из асимптотической формулы
\eqref{2.89} и включения $\varphi\in\mathcal{M}$ вытекает, что
\begin{equation}\label{2.102}
\varphi_{s,m}(\lambda_{j})\asymp
j^{\,s/n}\varphi(j^{\,1/n})\quad\mbox{при}\;\;j\geq1.
\end{equation}
В самом деле, по определению функции $\varphi_{s,m}$, данному в формуле
\eqref{2.82},
$$
\varphi_{s,m}(\lambda_{j})=\lambda_{j}^{s/m}\,\varphi(\lambda_{j}^{1/m})
\quad\mbox{при}\;\;\lambda_{j}\geq1.
$$
Поскольку $\varphi\in\mathcal{M}\subset\mathrm{QSV}$, существует положительная
функция $\varphi_{1}\in\mathrm{SV}$ такая, что $\varphi(t)\asymp\varphi_{1}(t)$ при
$t\gg1$. Значит,
\begin{equation}\label{2.103}
\varphi_{s,m}(\lambda_{j})\asymp\lambda_{j}^{s/m}\,\varphi_{1}(\lambda_{j}^{1/m})
\quad\mbox{при}\;\;j\gg1.
\end{equation}
Для функции $\varphi_{1}\in\mathrm{SV}$ в силу предложения \ref{prop2.2} (теорема о
равномерной сходимости) и формулы \eqref{2.89} имеем:
$$
\lim_{j\rightarrow\infty}\;\frac{\varphi_{1}(\lambda_{j}^{1/m})}{\varphi_{1}(j^{\,1/n})}=
\lim_{j\rightarrow\infty}\;
\frac{\varphi_{1}((\lambda_{j}^{1/m}j^{\,-1/n})\,j^{\,1/n})}{\varphi_{1}(j^{\,1/n})}=1.
$$
Отсюда и из \eqref{2.103} и \eqref{2.89} следует, что
\begin{equation}\label{2.104}
\varphi_{s,m}(\lambda_{j})\asymp j^{\,s/n}\varphi_{1}(j^{\,1/n})\asymp
j^{\,s/n}\varphi(j^{\,1/n})\quad\mbox{при}\;\;j\gg1.
\end{equation}
Поскольку функции $\varphi$ и $1/\varphi$ ограничены на каждом отрезке $[1,b]$, где
$1<b<\infty$, то \eqref{2.104} влечет за собой \eqref{2.102}.

Остается отметить, что теорема~\ref{th2.26} немедленно следует из лемм
\ref{lem2.11}, \ref{lem2.12} и формулы \eqref{2.102}.

Теорема~\ref{th2.26} доказана.


\begin{example}\label{ex2.7}
Пусть $\Gamma$~--- окружность радиуса 1 и $A:=1-d^{2}/dt^{2}$, где $t$ задает
натуральную параметризацию на $\Gamma$. Собственные функции
$h_{j}(t):=(2\pi)^{-1}\mathrm{e}^{ijt}$, $j\in\mathbb{Z}$, оператора $A$ образуют
ортонормированный базис в $L_{2}(\Gamma)$; они соответствуют собственным числам
$\lambda_{j}=1+j^{\,2}$. Пусть $s\in\mathbb{R}$ и $\varphi\in\mathcal{M}$. В силу
теоремы \ref{th2.24}
\begin{gather*}
\|f\|_{H^{s,\varphi}(\Gamma)}^{2}\asymp\|\varphi_{s,2}(A_{0})f\|_{L_{2}(\Gamma)}^{2}=\\
\sum_{j=-\infty}^{\infty}\;(1+j^{\,2})^{s}
\,\varphi^{2}((1+j^{\,2})^{1/2})\,|c_{j}(f)|^{2}.
\end{gather*}
Отметим, что здесь можно взять базис, образованный вещественными собственными
функциями $h_{0}(t):=(2\pi)^{-1}$, $h_{j}(t):=\pi^{-1}\cos jt$ и
$h_{-j}(t):=\pi^{-1}\sin jt$ для всех целых $j\geq\nobreak1$. Тогда
$$
\|f\|_{H^{s,\varphi}(\Gamma)}^{2}\asymp|a_{0}(f)|^{2}+
\sum_{j=1}^{\infty}\;j^{\,2s}\,\varphi^{2}(j)\,\bigl(|a_{j}(f)|^{2}+|b_{j}(f)|^{2}\bigr),
$$
где $a_{0}(f)$, $a_{j}(f)$ и $b_{j}(f)$~--- коэффициенты Фурье распределения $f$ по
указанным собственным функциям. В~этом случае пространство $H^{s,\varphi}(\Gamma)$
тесно связано с пространствами периодических вешественных функций, рассмотренными
А.~И.~Степанцом [\ref{Stepanets87}] (гл.~1, \S~7), [\ref{Stepanets02}] (ч.~1, гл.~3,
п.~7.1).
\end{example}

\subsection{Теорема вложения}\label{sec2.5.4}
Докажем важную теорему о вложении пространства $H^{s,\varphi}(\Gamma)$ в
пространство $C^{k}(\Gamma)$.


\begin{theorem}\label{th2.27}
Пусть заданы функция $\varphi\in\mathcal{M}$ и целое число $k\geq0$. Тогда условие
\eqref{2.37} равносильно вложению
\begin{equation}\label{2.105}
H^{k+n/2,\varphi}(\Gamma)\hookrightarrow C^{k}(\Gamma).
\end{equation}
Это вложение компактно.
\end{theorem}


\textbf{Доказательство.} Предположим, что верно \eqref{2.37}. Тогда по теореме
\ref{th2.15} (iii) справедливо непрерывное вложение
$H^{k+n/2,\varphi}(\mathbb{R}^{n})\hookrightarrow
C^{k}_{\mathrm{b}}(\mathbb{R}^{n})$. Отсюда нетрудно вывести непрерывное вложение
\eqref{2.105}, если воспользоваться локальным определением \ref{def2.13}
пространства $H^{k+n/2,\varphi}(\Gamma)$. Действительно, для любого $f\in
H^{k+n/2,\varphi}(\Gamma)$ имеем
$$
(\chi_{j}f)\circ\alpha_{j}\in H^{k+n/2,\varphi}(\mathbb{R}^{n})\hookrightarrow
C^{k}_{\mathrm{b}}(\mathbb{R}^{n})\quad\mbox{при}\;\;j=1,\ldots,r.
$$
Следовательно,
$$
f=\sum_{j=1}^{r}\;\chi_{j}f\in C^{k}(\Gamma).
$$
Кроме того,
\begin{gather*}
\|f\|_{C^{k}(\Gamma)}\leq\sum_{j=1}^{r}\;\|\chi_{j}f\|_{C^{k}(\Gamma)}\leq
c_{1}\sum_{j=1}^{r}\;
\|(\chi_{j}f)\circ\alpha_{j}\|_{C^{k}_{\mathrm{b}}(\mathbb{R}^{n})}\leq \\
c_{2}\sum_{j=1}^{r}\;\|(\chi_{j}f)\circ\alpha_{j}\|_{H^{k+n/2,\varphi}(\mathbb{R}^{n})}
\leq c_{3}\,\|f\|_{H^{k+n/2,\varphi}(\Gamma)}.
\end{gather*}
Здесь положительные числа $c_{1}$, $c_{2}$ и $c_{3}$ не зависят от $f$. Вложение
\eqref{2.105} и его непрерывность доказаны.

Покажем, что это вложение компактное. На основании теоремы \ref{th2.12} (i) можно
без ограничения общности считать функцию $\varphi\in\mathcal{M}$ непрерывной. В~силу
\eqref{2.37} функция $\psi_{1}:=\varphi^{2}$ удовлетворяет условию леммы
\ref{lem2.4}. Пусть $\psi_{0}$~--- функция из формулировки этой леммы. Тогда функция
$\varphi_{0}=\sqrt{\psi_{0}}\in\mathcal{M}$ удовлетворяет условию
$\varphi_{0}(t)/\varphi(t)\rightarrow0$ при $t\rightarrow\infty$ и неравенству
\eqref{2.37} с $\varphi_{0}$ вместо $\varphi$. Следовательно, по уже доказанному и в
силу теоремы \ref{th2.22} (iv)
$$
H^{k+n/2,\varphi}(\Gamma)\hookrightarrow
H^{k+n/2,\varphi_{0}}(\Gamma)\hookrightarrow C^{k}(\Gamma),
$$
причем первое вложение компактно, а второе непрерывно. Значит, вложение
\eqref{2.105} компактно.

Остается доказать, что из включения $H^{k+n/2,\varphi}(\Gamma)\subseteq
C^{k}(\Gamma)$ следует условие \eqref{2.37}. Снова воспользуемся локальным
определением \ref{def2.13} и будем использовать построения из доказательства
утверждения (iv) теоремы \ref{th2.22}. Для произвольного распределения $u\in
H^{k+n/2,\varphi}(\mathbb{R}^{n})$ такого, что
$\mathrm{supp}\,u\subset\alpha_{1}^{-1}(U)$, имеем
$$
\Theta(u\circ\alpha_{1}^{-1})\in H^{k+n/2,\varphi}(\Gamma)\subset C^{k}(\Gamma).
$$
Значит, $u\in C^{k}(\mathbb{R}^{n})$. Следовательно, в силу предложения
\ref{prop2.5} справедливо условие $\langle\xi\rangle^{k}\,\mu^{-1}(\xi)\in
L_{2}(\mathbb{R}^{n}_{\xi})$, где
$\mu(\xi):=\langle\xi\rangle^{k+n/2}\,\varphi(\langle\xi\rangle)$. Ввиду
\eqref{2.41} оно равносильно условию \eqref{2.37}.

Теорема \ref{th2.27} доказана.




\newpage


\markright{\emph \ref{sec2.6b}. Эллиптические операторы на замкнутом многообразии}

\section[Эллиптические операторы на замкнутом многообразии]
{Эллиптические операторы \\ на замкнутом многообразии}\label{sec2.6b}

\markright{\emph \ref{sec2.6b}. Эллиптические операторы на замкнутом многообразии}

В этом параграфе мы изучаем эллиптические ПДО, заданные на бесконечно гладком
замкнутом (компактном) многообразии $\Gamma$. Мы докажем, что эти ПДО ограничены и
нетеровы в подходящих парах пространств уточненной шкалы. Будет исследована также
глобальная и локальная гладкость решения эллиптического уравнения. Кроме того, мы
рассмотрим один класс эллиптических ПДО, зависящих от комплексного параметра, для
которых указанный оператор является гомеоморфизмом при больших по модулю значениях
параметра.

\subsection{ПДО на замкнутом многообразии}\label{sec2.6.1b}

Здесь для удобства читателя напомним определение ПДО на замкнутом многообразии
$\Gamma$ и связанные с ним понятия, необходимые нам. Как и прежде, будем
пользоваться терминологией и обозначениями из обзора М.~С.~Аграновича
[\ref{Agranovich90}] (\S~2).

\begin{definition}\label{defPsMan1}
Линейный оператор $A$, действующий в $C^{\infty}(\Gamma)$, называется ПДО (на
$\Gamma$) класса $\Psi^{m}(\Gamma)$, где $m\in\mathbb{R}$, если выполняются
следующие условия.
\begin{itemize}
\item [(i)] Для произвольных функций $\varphi,\psi\in C^{\infty}(\Gamma)$, носители
которых не пересекаются, отображение $u\mapsto\varphi A(\psi u)$, $u\in
C^{\infty}(\Gamma)$, продолжается по непрерывности до оператора порядка $-\infty$ в
соболевской шкале, т.~е. до ограниченного оператора класса $H^{s}(\Gamma)\rightarrow
H^{s+r}(\Gamma)$ для всех $s\in\mathbb{R}$ и $r>0$;
\item [(ii)] Пусть на $\Gamma$ произвольно выбраны локальная карта
$\alpha:\mathbb{R}^{n}\rightarrow\Gamma_{\alpha}$ класса $C^{\infty}$ и открытое
множество $\Omega$, замыкание которого содержится в координатной окрестности
$\Gamma_{\alpha}\subset\Gamma$. Тогда существует ПДО
$A_{\Omega}\in\Psi^{m}(\mathbb{R}^{n})$ такой, что для любых функций
$\varphi,\psi\in C^{\infty}(\Gamma)$, носители которых лежат в $\Omega$, выполняется
равенство
$$
(\varphi A(\psi u))\circ\alpha=(\varphi\circ\alpha)A_{\Omega}((\psi
u)\circ\alpha),\quad u\in C^{\infty}(\Gamma).
$$
\end{itemize}
\end{definition}

В связи с пунктом (i) этого определения отметим, что используемое в нем условие
„быть оператором порядка $-\infty$ в соболевской шкале” равносильно условию „быть
интегральным оператором
\begin{equation}\label{PsMan1}
u(y)\mapsto\int\limits_{\Gamma}\,K(x,y)\,u(y)\,dy,\quad u\in C^{\infty}(\Gamma),
\end{equation}
на $\Gamma$ с бесконечно гладким ядром $K(x,y)$”. Поэтому введенное понятие ПДО на
$\Gamma$ не связано на самом деле с той или иной шкалой гильбертовых пространств.

Положим
$$
\Psi^{-\infty}(\Gamma):=\bigcap_{m\in\mathbb{R}}\,\Psi^{m}(\Gamma),\quad
\Psi^{\infty}(\Gamma):=\bigcup_{m\in\mathbb{R}}\,\Psi^{m}(\Gamma).
$$
Класс $\Psi^{-\infty}(\Gamma)$ совпадает c классом всех интегральных операторов
\eqref{PsMan1} с бесконечно гладкими ядрами. Всякий ПДО $A\in\Psi^{\infty}(\Gamma)$
непрерывен в $C^{\infty}(\Gamma)$ и единственным образом продолжается до линейного
непрерывного оператора в $\mathcal{D}'(\Gamma)$. Этот оператор также обозначается
через $A$.

\begin{definition}\label{defPsMan2}
ПДО $A\in\Psi^{m}(\Gamma)$, где $m\in\mathbb{R}$, называется полиоднородным (или
классическим) на $\Gamma$ порядка $m$, если в определении~\ref{defPsMan1} все
$A_{\Omega}\in\Psi^{m}_{\mathrm{ph}}(\mathbb{R}^{n})$.
\end{definition}

Обозначим через $\Psi^{m}_{\mathrm{ph}}(\Gamma)$ класс всех полиоднородных ПДО на
$\Gamma$ порядка $m$.

\begin{definition}\label{defPsMan3}
Пусть задан ПДО $A\in\Psi^{m}_{\mathrm{ph}}(\Gamma)$. Его главный символ
$a_{0}(x,\xi)$ определяется как функция аргументов $x\in\Gamma$ и $\xi\in
T^{\ast}_{x}\Gamma$, $\xi\neq0$, совпадающая локально по $x$ с главным символом
соответствующего ПДО $A_{\Omega}\in\Psi^{m}_{\mathrm{ph}}(\mathbb{R}^{n})$.
\end{definition}

Здесь, как обычно, $T^{\ast}_{x}\Gamma$ обозначает кокасательное пространство к
многообразию $\Gamma$ в точке $x\in\Gamma$. Важно, что главный символ
полиоднородного ПДО на $\Gamma$ не зависит от выбора локальных карт на $\Gamma$. Он
является бесконечно гладкой функцией аргументов $x$ и $\xi$, положительно однородной
по $\xi$ порядка $m$.

\begin{definition}\label{defPsMan4}
Линейный оператор $A^{+}$, заданный в $C^{\infty}(\Gamma)$ и \emph{формально
сопряженный} к ПДО $A\in\Psi^{\infty}(\Gamma)$, определяется по формуле
$$
(Au,v)_{\Gamma}=(u,A^{+}v)_{\Gamma}\quad\mbox{для всех}\quad u,v\in
C^{\infty}(\Gamma).
$$
Если $A=A^{+}$, то оператор $A$ называется \emph{формально самосопряженным}.
\end{definition}

Отметим, что понятие формально сопряженного ПДО введено относительно заданной на
$\Gamma$ плотности $dx$, поскольку $(\cdot,\cdot)_{\Gamma}$~--- скалярное
произведение в $L_{2}(\Gamma,dx)$. Если $A\in\Psi^{m}(\Gamma)$ для некоторого
$m\in\mathbb{R}$, то и $A^{+}\in\Psi^{m}(\Gamma)$. Если, кроме того, ПДО $A$
полиоднородный с главным символом $a_{0}$, то ПДО $A^{+}$ также полиоднородный с
комплексно-сопряженным главным символом~$\overline{a_{0}}$.

\begin{definition}\label{def2.17b}
ПДО $A\in\Psi^{m}_{\mathrm{ph}}(\Gamma)$ и его главный символ $a_{0}(x,\xi)$
называются \emph{эллиптическими} на $\Gamma$, если $a_{0}(x,\xi)\neq0$ для
произвольных точки $x\in\Gamma$ и ковектора $\xi\in
T^{\ast}_{x}\Gamma\setminus\{0\}$.
\end{definition}

ПДО $A\in\Psi^{m}_{\mathrm{ph}}(\Gamma)$, эллиптический на $\Gamma$, имеет
параметрикс $B\in\Psi^{-m}_{\mathrm{ph}}(\Gamma)$, также эллиптический на $\Gamma$,
т.~е. справедлив аналог предложения \ref{prop2.6b}. Правда, он нам не понадобится.

Важным примером эллиптического ПДО на $\Gamma$ является оператор Бельтрами--Лапласа
$\Delta_{\Gamma}$. Напомним его определение. Пусть на многообразии $\Gamma$ введена
риманова метрика, т.~е. задано бесконечно гладкое ковариантное вещественное
тензорное поле $g(x)$, $x\in\Gamma$, где $g(x)=(g_{j,k}(x))_{j,k=1}^{n}$~---
симметричная положительно определенная матрица. Риманова метрика определяет на
$\Gamma$ плотность $dx:=(\det g(x))^{1/2}dx_{1}\ldots dx_{1}$ в локальных
координатах $x=(x_{1},\ldots,x_{n})$. По определению, действие оператора
Бельтрами--Лапласа на функцию $u\in C^{2}(\Gamma)$ задается формулой
$$
(\Delta_{\Gamma}u)(x):=(\det g(x))^{-1/2}\sum_{j,k=1}^{n}\partial_{x_{j}}\bigl((\det
g(x))^{1/2}\,g^{j,k}(x)\,\partial_{x_{k}}u(x)\bigr),
$$
где $g^{-1}(x)=(g^{j,k}(x))_{j,k=1}^{n}$~--- матрица, обратная к $g(x)$. Функция
$\Delta_{\Gamma}u$ не зависит от выбора локальных координат на $\Gamma$. Главный
символ дифференциального оператора $\Delta_{\Gamma}$ равен
$$
\sum_{j,k=1}^{n}g^{j,k}(x)\,\xi_{j}\,\xi_{k},
$$
и поэтому эллиптичность $\Delta_{\Gamma}$ следует из положительной определенности
матрицы $g^{-1}(x)$. Оператор $\Delta_{\Gamma}$ формально самосопряженный
относительно плотности $dx$.

\subsection[Эллиптический оператор в уточненной шкале]
{Эллиптический оператор в уточненной шкале}\label{sec2.6.2b}

Пусть задан ПДО $A\in\Psi^{m}_{\mathrm{ph}}(\Gamma)$ порядка $m\in\mathbb{R}$. В~п.
\ref{sec2.6.2b} и \ref{sec2.6.3b} мы предполагаем, что $A$ еллиптичен на $\Gamma$.

Положим
\begin{gather}\label{PsMan2}
\mathcal{N}:=\{\,u\in
C^{\infty}(\Gamma):\,Au=0\;\;\mbox{на}\;\;\Gamma\,\},\\
\mathcal{N}^{+}:=\{\,v\in C^{\infty}(\Gamma):\,A^{+}v=0\;\;\mbox{на}\;\;\Gamma\,\}.
\label{PsMan3}
\end{gather}
Поскольку оба ПДО $A$ и $A^{+}$ эллиптические на  $\Gamma$, пространства
$\mathcal{N}$ и $\mathcal{N}^{+}$ конечномерные [\ref{Agranovich90}, с. 28] (теорема
2.3.3).

Изучим свойства ПДО $A$ в уточненной шкале на $\Gamma$. Предварительно докажем
следующее утверждение о действии произвольного ПДО в этой шкале (аналог
леммы~\ref{lem2.8b} ).


\begin{lemma}\label{lem2.13b}\label{lem2.13}
Пусть $G$ --- ПДО класса $\Psi^{r}(\Gamma)$, где $r\in\mathbb{R}$. Тогда сужение
отображения $u\mapsto Gu$, $u\in\mathcal{D}'(\Gamma)$, на пространство
$H^{\sigma,\varphi}(\Gamma)$ является линейным ограниченным оператором
\begin{equation}\label{2.107b}
G:H^{\sigma,\varphi}(\Gamma)\rightarrow H^{\sigma-r,\,\varphi}(\Gamma)
\end{equation}
для любых $\sigma\in\mathbb{R}$ и $\varphi\in\mathcal{M}$.
\end{lemma}


\textbf{Доказательство.} В случае $\varphi\equiv1$ эта лемма известна
[\ref{Agranovich90}, с.~23]. Выберем произвольные $\sigma\in\mathbb{R}$ и
$\varphi\in\mathcal{M}$. Рассмотрим линейные ограниченные операторы
$$
G:H^{\sigma\mp1}(\Gamma)\rightarrow H^{\sigma\mp1-r}(\Gamma).
$$
Применим интерполяцию с функциональным параметром $\psi$ из теоремы \ref{th2.21},
где $\varepsilon=\delta=1$. Получим ограниченный оператор
$$
G:\bigl[\,H^{\sigma-1}(\Gamma),\,H^{\sigma+1}(\Gamma)\,\bigr]_{\psi}\rightarrow
\bigl[\,H^{\sigma-r-1}(\Gamma),\,H^{\sigma-r+1}(\Gamma)\,\bigr]_{\psi}.
$$
Отсюда в силу теоремы \ref{th2.21} получаем ограниченность оператора \eqref{2.107b}.
Лемма~\ref{lem2.13b} доказана.

\medskip

В силу леммы~\ref{lem2.8b} имеем ограниченный оператор
\begin{equation}\label{2.108b}
A:H^{s+m,\varphi}(\Gamma)\rightarrow H^{s,\varphi}(\Gamma)
\end{equation}
для произвольных параметров $s\in\mathbb{R}$ и $\varphi\in\mathcal{M}$. Изучим
свойства этого оператора.


\begin{theorem}\label{th2.28b}
Для произвольных $s\in\mathbb{R}$ и $\varphi\in\mathcal{M}$ ограниченный оператор
\eqref{2.108b} нетеров. Его ядро совпадает с $\mathcal{N}$, а область значений равна
\begin{equation}\label{2.110b}
\bigl\{f\in H^{s,\varphi}(\Gamma):\,(f,w)_{\Gamma}=0\;\mbox{для
всех}\;w\in\mathcal{N}^{+}\bigr\}.
\end{equation}
Индекс оператора \eqref{2.108b} равен $\dim\mathcal{N}-\dim\mathcal{N}^{+}$ и не
зависит от $s$, $\varphi$.
\end{theorem}


\textbf{Доказательство.} В случае $\varphi\equiv1$ (соболевская шкала) эта теорема
известна (см., например, [\ref{Agranovich90}], \S~2, теоремы 2.2.6, 2.3.3 и 3.3.12,
либо [\ref{Hermander87}, с.~262], теорема 19.2.1). Отсюда общий случай
$\varphi\in\mathcal{M}$ получается с помощью интерполяции с функциональным
параметром. А именно, пусть $s\in\mathbb{R}$. Имеем ограниченные нетеровы операторы
\begin{equation}\label{2.111b}
A:H^{s\mp1+m}(\Gamma)\rightarrow H^{s\mp1}(\Gamma)
\end{equation}
с общим ядром $\mathcal{N}$ и одинаковым индексом
$\varkappa:=\dim\mathcal{N}-\dim\mathcal{N}^{+}$. При этом
\begin{gather}\notag
A\bigl(H^{s\mp1+m}(\Gamma)\bigr)=\\
\bigl\{f\in H^{s\mp1}(\Gamma):(f,w)_{\Gamma}=0\;\mbox{для
всех}\;w\in\mathcal{N}^{+}\bigr\}. \label{2.112b}
\end{gather}
Применим к \eqref{2.111b} интерполяцию с функциональным параметром $\psi$ из теоремы
\ref{th2.21}, где $\varepsilon=\delta=1$. Получим ограниченный оператор
\begin{equation*}
A:\,\bigl[\,H^{s-1+m}(\Gamma),\,H^{s+1+m}(\Gamma)\,\bigr]_{\psi}\rightarrow
\bigl[\,H^{s-1}(\Gamma),\,H^{s+1}(\Gamma)\,\bigr]_{\psi},
\end{equation*}
который в силу теоремы \ref{th2.21} совпадает с оператором \eqref{2.108b}.
Следовательно, согласно теореме \ref{th2.7} оператор \eqref{2.108b} нетеров с ядром
$\mathcal{N}$ и индексом $\varkappa=\dim\mathcal{N}-\dim\mathcal{N}^{+}$. Область
значений этого оператора равна
$$
H^{s,\varphi}(\Gamma)\cap A\bigl(H^{s-1+m}(\Gamma)\bigr).
$$
Отсюда, в силу \eqref{2.112b} получаем, что она равна \eqref{2.110b}. Теорема
\ref{th2.28b} доказана.

\medskip

Согласно этой теореме, $\mathcal{N}^{+}$ --- дефектное подпространство оператора
\eqref{2.108b}. Заметим, что в виду теоремы \ref{th2.22} (v) оператор
\begin{equation}\label{2.113b}
A^{+}:H^{-s,1/\varphi}(\Gamma)\rightarrow H^{-s-m,1/\varphi}(\Gamma)
\end{equation}
является сопряженным к оператору \eqref{2.108b}. Поскольку ПДО $A^{+}$ эллиптический
на $\Gamma$, то в силу теоремы \ref{th2.28b} ограниченный оператор \eqref{2.113b}
нетеров и имеет ядро $\mathcal{N}^{+}$ и дефектное подпространство~$\mathcal{N}$.
Отметим также [\ref{AtiyahSinger63}; \ref{Agranovich90}, с.~32], что индексы
операторов \eqref{2.108b} и \eqref{2.113b} равны 0, если только $\dim\Gamma\geq2$.

Если пространства $\mathcal{N}$ и $\mathcal{N}^{+}$ тривиальны, то из теоремы
\ref{th2.28b} и теоремы Банаха об обратном операторе вытекает, что оператор
\eqref{2.108b} является гомеоморфизмом. В общем случае гомеоморфизм удобно задавать
с помощью следующих проекторов.

Представим пространства, в которых действует оператор \eqref{2.108b}, в виде прямых
сумм замкнутых подпространств:
\begin{gather*}
H^{s+m,\varphi}(\Gamma)=\mathcal{N}\dotplus\bigl\{u\in H^{s+m,\varphi}(\Gamma):
(u,v)_{\Gamma}=0\;\mbox{для всех}\;v\in\mathcal{N}\bigr\},\\
H^{s,\varphi}(\Gamma)=\mathcal{N}^{+}\dotplus\bigl\{f\in H^{s,\varphi}(\Gamma):
(f,w)_{\Gamma}=0\;\mbox{для всех}\;w\in\mathcal{N}^{+}\bigr\}.
\end{gather*}
Указанные разложения в прямые суммы существуют, поскольку в них слагаемые имеют
тривиальное пересечение, и конечная размерность первого из них равна коразмерности
второго. (Действительно, например, в первой сумме факторпространство пространства
$H^{s+m,\varphi}(\Gamma)$ по второму слагаемому является двойственным пространством
к подпространству $\mathcal{N}$ пространства $H^{-s-m,1/\varphi}(\Gamma)$).
Обозначим через $\mathcal{P}$ и $\mathcal{P}^{+}$ косые проекторы соответственно
пространств $H^{s+m,\varphi}(\Gamma)$ и $H^{s,\varphi}(\Gamma)$ на вторые слагаемые
в указанных суммах параллельно первым слагаемым. Эти проекторы не зависят от $s$ и
$\varphi$.


\begin{theorem}\label{th2.29b}
Для произвольных $s\in\mathbb{R}$ и $\varphi\in\mathcal{M}$ сужение оператора
\eqref{2.108b} на подпространство $\mathcal{P}(H^{s+m,\varphi}(\Gamma))$ является
гомеоморфизмом
\begin{equation}\label{2.114b}
A:\mathcal{P}(H^{s+m,\varphi}(\Gamma))\leftrightarrow
\mathcal{P}^{+}(H^{s,\varphi}(\Gamma)).
\end{equation}
\end{theorem}


\textbf{Доказательство.} По теореме \ref{th2.28b}, $\mathcal{N}$~--- ядро, а
$\mathcal{P}^{+}(H^{s,\varphi}(\Gamma))$~--- область значений оператора
\eqref{2.108b}. Следовательно, оператор \eqref{2.114b}~--- биекция. Кроме того, этот
оператор ограничен. Значит, он является гомеоморфизмом в силу теоремы Банаха об
обратном операторе. Теорема~\ref{th2.29b} доказана.

\medskip

Из теоремы \ref{th2.29b} вытекает следующая априорная оценка решения эллиптического
уравнения $Au=f$ на $\Gamma$ (аналог теоремы теоремы \ref{th2.16b} для замкнутого
многообразия).


\begin{theorem}\label{th2.30b}
Пусть $s\in\mathbb{R}$, $\sigma>0$ и $\varphi\in\mathcal{M}$. Существует число
$c=c(s,\sigma,\varphi)>0$ такое, что для произвольного распределения $u\in
H^{s+m,\varphi}(\Gamma)$ справедлива априорная оценка
\begin{equation}\label{2.116b}
\|u\|_{H^{s+m,\varphi}(\Gamma)}\leq c\,\bigl(\,\|Au\|_{H^{s,\varphi}(\Gamma)}+
\|u\|_{H^{s+m-\sigma,\varphi}(\Gamma)}\,\bigr).
\end{equation}
\end{theorem}


\textbf{Доказательство.} Так как $\mathcal{N}$~--- конечномерное подпространство в
пространствах $H^{s+m,\varphi}(\Gamma)$ и $H^{s+m-\sigma,\varphi}(\Gamma)$, то нормы
в них эквивалентны на $\mathcal{N}$. В частности, для распределения
$u-\mathcal{P}u\in\mathcal{N}$ справедливо
$$
\|u-\mathcal{P}u\|_{H^{s+m,\varphi}(\Gamma)}\leq\,
c_{1}\,\|u-\mathcal{P}u\|_{H^{s+m-\sigma,\varphi}(\Gamma)}
$$
с постоянной $c_{1}>0$, не зависящей от $u\in H^{s+m,\varphi}(\Gamma)$. Отсюда
получаем
\begin{gather*}
\|u\|_{H^{s+m,\varphi}(\Gamma)}\leq\,
\|u-\mathcal{P}u\|_{H^{s+m,\varphi}(\Gamma)}+\|\mathcal{P}u\|_{H^{s+m,\varphi}(\Gamma)}\leq\\
c_{1}\,\|u-\mathcal{P}u\|_{H^{s+m-\sigma,\varphi}(\Gamma)}+
\|\mathcal{P}u\|_{H^{s+m,\varphi}(\Gamma)}\leq\,\\
c_{1}\,c_{2}\,\|u\|_{H^{s+m-\sigma,\varphi}(\Gamma)}+
\|\mathcal{P}u\|_{H^{s+m,\varphi}(\Gamma)},
\end{gather*}
где $c_{2}$ --- норма проектора $1-\mathcal{P}$, действующего в пространстве
$H^{s+m-\sigma,\varphi}(\Gamma)$. Итак,
\begin{equation}\label{2.117b}
\|u\|_{H^{s+m,\varphi}(\Gamma)}\leq
c_{1}\,c_{2}\,\|u\|_{H^{s+m-\sigma,\varphi}(\Gamma)}+
\|\mathcal{P}u\|_{H^{s+m,\varphi}(\Gamma)}.
\end{equation}
Воспользуемся теперь условием $Au=f$. Поскольку $\mathcal{N}$~--- ядро оператора
\eqref{2.108b} и $u-\mathcal{P}u\in\mathcal{N}$, то $A\mathcal{P}u=f$. Таким
образом, $\mathcal{P}u$~--- прообраз распределения $f$ относительно гомеоморфизма
\eqref{2.114b}. Следовательно,
$$
\|\mathcal{P}u\|_{H^{s+m,\varphi}(\Gamma)}\leq\,
c_{3}\,\|f\|_{H^{s,\varphi}(\Gamma)},
$$
где $c_{3}$ --- норма оператора, обратного к \eqref{2.114b}. Отсюда и из неравенства
\eqref{2.117b} немедленно вытекает оценка \eqref{2.116b}. Теорема \ref{th2.30b}
доказана.

\medskip

Отметим, что, если $\mathcal{N}=\{0\}$, т. е.  уравнение $Au=f$ имеет не более
одного решения, то норма $\|u\|_{H^{s+m-\sigma}(\Gamma)}$ в правой части оценки
\eqref{2.116b} можно опустить. Если же $\mathcal{N}\neq\{0\}$, то для каждого
распределения $u$ эту величину можно сделать как угодно малой за счет выбора
достаточно большого числа $\sigma$.

\subsection[Гладкость решения эллиптического уравне-\break ния]
{Гладкость решения \\ эллиптического уравнения}\label{sec2.6.3b}

Пусть $\Gamma_{0}$~--- непустое открытое подмножество многообразия $\Gamma$.
Исследуем локальную гладкость решения эллиптического уравнения $Au=f$ на
$\Gamma_{0}$ в уточненной шкале. Рассмотрим сначала случай, когда
$\Gamma_{0}=\Gamma$.


\begin{theorem}\label{th2.31b}
Предположим, что распределение $u\in\mathcal{D}'(\Gamma)$ является решением
уравнения $Au=f$ на $\Gamma$, где $f\in H^{s,\varphi}(\Gamma)$ для некоторых
параметров $s\in\mathbb{R}$ и $\varphi\in\mathcal{M}$. Тогда $u\in
H^{s+m,\varphi}(\Gamma)$.
\end{theorem}


\textbf{Доказательство.} Поскольку многообразие $\Gamma$ компактно, то пространство
$\mathcal{D}'(\Gamma)$ является объединением соболевских пространств
$H^{\sigma}(\Gamma)$, где $\sigma\in\mathbb{R}$. Следовательно, для распределения
$u\in\mathcal{D}'(\Gamma)$ существует число $\sigma<s$, такое что $u\in
H^{\sigma+m}(\Gamma)$. В~силу теорем \ref{th2.28b} и \ref{th2.22} (iii) справедливо
равенство
$$
(H^{s,\varphi}(\Gamma))\cap A(H^{\sigma+m}(\Gamma))= A(H^{s+m,\varphi}(\Gamma)).
$$
Поэтому, из условия $f\in H^{s,\varphi}(\Gamma)$ вытекает, что
$$
f=Au\in A(H^{s+m,\varphi}(\Gamma)).
$$
Таким образом, на многообразии $\Gamma$ наряду с равенством $Au=f$ выполняется также
равенство $Av=f$ для некоторого распределения $v\in H^{s+m,\varphi}(\Gamma)$.
Следовательно, $A(u-v)=0$ на $\Gamma$ и, согласно теореме \ref{th2.28b} справедливо
включение
$$
w:=u-v\in\mathcal{N}\subset C^{\infty}(\Gamma)\subset H^{s+m,\varphi}(\Gamma)
$$
Значит, $u=v+w\in H^{s+m,\varphi}(\Gamma)$. Теорема \ref{th2.31b} доказана.

\medskip

Рассмотрим теперь случай произвольного $\Gamma_{0}$. Обозначим
\begin{gather*}
H^{\sigma,\varphi}_{\mathrm{loc}}(\Gamma_{0}):=\{h\in\mathcal{D}'(\Gamma): \\
\chi\,h\in H^{\sigma,\varphi}(\Gamma)\;\mbox{для всех}\;\chi\in
C^{\infty}(\Gamma),\,\mathrm{supp}\,\chi\subset \Gamma_{0}\};
\end{gather*}
здесь $\sigma\in\mathbb{R}$, $\varphi\in\mathcal{M}$.


\begin{theorem}\label{th2.32b}
Предположим, что распределение $u\in\mathcal{D}'(\Gamma)$ является решением
уравнения $Au=f$ на $\Gamma_{0}$, где $f\in
H^{s,\varphi}_{\mathrm{loc}}(\Gamma_{0})$ для некоторых параметров $s\in\mathbb{R}$
и $\varphi\in\mathcal{M}$. Тогда $u\in H^{s+m,\varphi}_{\mathrm{loc}}(\Gamma_{0})$.
\end{theorem}


Теорема \ref{th2.32b} доказывается аналогично теореме \ref{th2.18b}. При этом надо
вместо теоремы \ref{th2.17b} использовать теорему \ref{th2.31b}.

Теорема \ref{th2.32b} уточняет применительно к шкале пространств
$H^{s,\varphi}(\Gamma)$ известные утверждения о повышении локальной гладкости
решения эллиптического уравнения на многообразии в соболевской шкале
[\ref{Egorov84}, с.~77;\,\ref{Shubin78}, с.~61]. Как видим, локальная уточненная
гладкость $\varphi$ правой части эллиптического уравнения наследуется его решением.

Из теорем \ref{th2.32b} и \ref{th2.27} немедленно вытекает следующее достаточное
условие непрерывности производных у решения $u$ уравнения $Au=f$.


\begin{corollary}\label{cor2.4b}
Пусть заданы целое число $r\geq0$ и функция $\varphi\in\mathcal{M}$, удовлетворяющая
условию \eqref{2.37}. Предположим, что распределение $u\in\mathcal{D}'(\Gamma)$
является решением уравнения $Au=f$ на множестве $\Gamma_{0}$, где $f\in
H^{r-m+n/2,\varphi}_{\mathrm{loc}}(\Gamma_{0})$. Тогда $u\in C^{r}(\Gamma_{0})$.
\end{corollary}

Отметим, что здесь условие \eqref{2.37} не только достаточное для справедливости
включения $u\in C^{r}(\Gamma_{0})$, но также и необходимое на классе всех
рассматриваемых решений уравнения $Au=f$; см. замечание~\ref{rem2.8c}.

\subsection{Эллиптический оператор с параметром}\label{sec2.6.4b} Различные
эллиптические операторы с параметром изучались в работах Ш.~Агмона, Л.~Ниренберга
[\ref{AgmonNirenberg63}], М.~С.~Аграновича, М.~И.~Вишика [\ref{AgranovichVishik64}],
А.~Н.~Кожевникова [\ref{Kozhevnikov73}] и их последователей (см. обзор
[\ref{Agranovich90}] и цитируемую там литературу). Ими было установлено, что при
достаточно больших по модулю значениях комплексного параметра эллиптический оператор
является гомеоморфизмом в подходящих парах соболевских пространств, причем норма
оператора допускает некоторую двустороннюю оценку с постоянными, не зависящими от
параметра. Мы рассмотрим широкий класс эллиптических ПДО с параметром на
многообразии $\Gamma$ и распространим указанный результат на уточненную шкалу
пространств.

Определении эллиптического ПДО с параметром мы берем из обзора [\ref{Agranovich90},
с.~57]. Произвольно зафиксируем числа $q\in\mathbb{N}$ и $m>\nobreak0$. Рассмотрим
ПДО $A(\lambda)\in\Psi^{mq}(\Gamma)$, зависящий от комплексного параметра $\lambda$
следующим образом:
\begin{equation}\label{2.120b}
A(\lambda):=\,\sum_{r\,=\,0}^{q}\,\lambda^{q-r}\,A^{(r)}.
\end{equation}
Здесь $A^{(r)}$, $r=0,\ldots,q$,~--- полиоднородные ПДО на $\Gamma$ порядков
$\mathrm{ord}\,A^{(r)}\leq mr$, причем  причем  $A_{0}$ --- это оператор умножения
на некоторую функцию $a_{0}\in C^{\infty}(\Gamma)$.

Пусть $K$~--- фиксированный замкнутый угол на комплексной плоскости с вершиной в
начале координат (не исключается случай, когда $K$ вырождается в луч).


\begin{definition}\label{def2.18b}
ПДО $A(\lambda)$ называется \emph{эллиптическим с параметром} в угле $K$, если
\begin{equation}\label{2.122b}
\sum_{r\,=\,0}^{q}\,\lambda^{q-r}\,a^{r,0}(x,\xi)\neq0
\end{equation}
для любых $x\in\Gamma$, $\xi\in T_{x}^{\ast}\Gamma$ и $\lambda\in K$ таких, что
$(\xi,\lambda)\neq0$. Здесь $a^{r,0}(x,\xi)$~--- главный символ ПДО $A^{(r)}$ в
случае $\mathrm{ord}\,A^{(r)}=mr$, либо $a^{r,0}(x,\xi)\equiv0$ в случае
$\mathrm{ord}\,A^{(r)}<mr$. При этом $a^{0,0}(x,\xi)\equiv a_{0}(\xi)$, а функции
$a^{1,0}(x,\xi),\,a^{2,0}(x,\xi),\ldots$ считаются равными 0 при $\xi=0$ (такое
допущение обусловлено тем, что главные символы изначально не определены при
$\xi=0$).
\end{definition}


\begin{example}\label{exPsMan1}
Рассмотрим ПДО $A-\lambda I$, где $A\in\Psi_{\mathrm{\mathrm{ph}}}^{m}(\Gamma)$, а
$I$~--- тождественный оператор. Для этого ПДО условие эллиптичности с параметром в
угле $K$ означает просто, что $a_{0}(x,\xi)\notin K$ при $\xi\neq0$. Здесь, как и
прежде, $a_{0}(x,\xi)$~--- главный символ ПДО $A$. Этот пример важен в спектральной
теории эллиптических операторов.
\end{example}

Далее в п. \ref{sec2.6.4b} мы предполагаем, что ПДО $A(\lambda)$ удовлетворяет
определению \ref{def2.18b}.

Из этого определения следует, что ПДО $A(\lambda)$ эллиптический на $\Gamma$ при
каждом фиксированном $\lambda\in\mathbb{C}$. В самом деле, главный символ ПДО
$A(\lambda)$ равен $a^{q,0}(x,\xi)$ для каждого $\lambda$ и, как это следует из
\eqref{2.122b} при $\lambda=0$, удовлетворяет неравенству $a^{q,0}(x,\xi)\neq0$ для
любых $x\in\Gamma$ и $\xi\in T^{\ast}_{x}\Gamma\setminus\{0\}$. Последнее означает,
что ПДО $A(\lambda)$ эллиптический на $\Gamma$. Cогласно теореме \ref{th2.28b} мы
имеем ограниченный и нетеров оператор
\begin{equation}\label{2.124b}
A(\lambda):H^{s+mq,\varphi}(\Gamma)\rightarrow H^{s,\varphi}(\Gamma)
\end{equation}
для произвольных $s\in\mathbb{R}$, $\varphi\in\mathcal{M}$ и $\lambda\in\mathbb{C}$.
Более того, поскольку $A(\lambda)$~--- эллиптический с параметром ПДО в угле $K$, он
имеет важные дополнительные свойства.


\begin{theorem}\label{th2.33b}
Справедливы следующие утверждения.
\begin{itemize}
\item [$\mathrm{(i)}$] Существует число $\lambda_{0}>0$ такое, что для
каждого значения параметра $\lambda\in K$, удовлетворяющего условию $|\lambda|\geq
\lambda_{0}$, ПДО $A(\lambda)$ устанавливает при любых $s\in\mathbb{R}$ и
$\varphi\in\mathcal{M}$ гомеоморфизм
\begin{equation}\label{2.125b}
A(\lambda):H^{s+mq,\varphi}(\Gamma)\leftrightarrow H^{s,\varphi}(\Gamma).
\end{equation}
\item [$\mathrm{(ii)}$] Для произвольных параметров
$s\in\mathbb{R}$ и $\varphi\in\mathcal{M}$ найдется число $c=c(s,\varphi)\geq1$
такое, что выполняется двусторонняя оценка
\begin{gather}\notag
c^{-1}\,\|A(\lambda)u\|_{H^{s,\varphi}(\Gamma)}\leq\\ \notag
\bigl(\,\|u\|_{H^{s+mq,\varphi}(\Gamma)}+
|\lambda|^{q}\,\|u\|_{H^{s,\varphi}(\Gamma)}\,\bigr)\leq\\
c\,\bigl\|A(\lambda)u\|_{H^{s,\varphi}(\Gamma)}\label{2.126b}
\end{gather}
для любого $\lambda\in K$, $|\lambda|\geq\lambda_{0}$, и произвольного $u\in
H^{s+mq,\varphi}(\Gamma)$.
\end{itemize}
\end{theorem}


В случае $\varphi\equiv1$ (соболевские пространства) эта теорема известна
[\ref{Agranovich90}, с.~58] (теорема 4.1.2); при этом левое неравенство в
двусторонней оценке \eqref{2.126b} справедливо без предположения об эллиптичности с
параметром ПДО $A(\lambda)$.

Докажем отдельно утверждения (i) и (ii) теоремы \ref{th2.33b}. Мы выведем общий
случай $\varphi\in\mathcal{M}$ из соболевского случая $\varphi\equiv1$.

\medskip

\textbf{Доказательство теоремы~\ref{th2.33b} (i).} Пусть $s\in\mathbb{R}$ и
$\varphi\in\mathcal{M}$. Поскольку для каждого $\lambda\in\mathbb{C}$ ПДО
$A(\lambda)$ эллиптический на $\Gamma$, то в силу теоремы \ref{th2.28b} ограниченный
нетеров оператор \eqref{2.124b} имеет независящие от $s$ и $\varphi$ конечномерные
ядро $\mathcal{N}(\lambda)$ и дефектное подпространство $\mathcal{N}^{+}(\lambda)$.
Далее воспользуемся тем, что теорема \ref{th2.33b} верна в случае $\varphi\equiv1$.
Существует число $\lambda_{0}>0$ такое, что для каждого $\lambda\in K$,
удовлетворяющего условию $|\lambda|\geq\lambda_{0}$, ПДО $A(\lambda)$ устанавливает
гомеоморфизм
$$
A(\lambda):H^{s+mq}(\Gamma)\leftrightarrow H^{s+m}(\Gamma).
$$
Следовательно, для указанного $\lambda$ пространства $\mathcal{N}(\lambda)$ и
$\mathcal{N}^{+}(\lambda)$ тривиальны, т.~е. линейный ограниченный оператор
\eqref{2.124b} является биекцией. Отсюда по теореме Банаха об обратном операторе
получаем гомеоморфизм \eqref{2.125b}. Утверждение (i) теоремы~\ref{th2.33b}
доказано.

\medskip

Утверждение (ii) теоремы~\ref{th2.33b} мы докажем с помощью интерполяции с
функциональным параметром. Для этого нам понадобится следующее пространство.

Пусть задано функцию $\varphi\in\mathcal{M}$ и числа $\sigma\in\mathbb{R}$,
$\varrho>0$, $\theta>0$. Обозначим через $H^{\sigma,\varphi}(\Gamma,\varrho,\theta)$
пространство $H^{\sigma,\varphi}(\Gamma)$, наделенное нормой, зависящей от числовых
параметров $\varrho$ и $\theta$ следующим образом:
$$
\|h\|_{H^{\sigma,\varphi}(\Gamma,\varrho,\theta)}:=
\left(\;\|h\|_{H^{\sigma,\varphi}(\Gamma)}^{2}+\varrho^{2}\,
\|h\|_{H^{\sigma-\theta,\varphi}(\Gamma)}^{2}\:\right)^{1/2}.
$$
Это определение корректно в силу непрерывного вложения
$H^{\sigma,\varphi}(\Gamma)\hookrightarrow H^{\sigma-\theta,\varphi}(\Gamma)$.
Отсюда вытекает, что нормы в пространствах
$H^{\sigma,\varphi}(\Gamma,\varrho,\theta)$ и $H^{\sigma,\varphi}(\Gamma)$
эквивалентны. Норма в пространстве $H^{\sigma,\varphi}(\Gamma,\varrho,\theta)$
порождена скалярным произведением
$$
(h_{1},h_{2})_{H^{\sigma,\varphi}(\Gamma,\varrho,\theta)}:=
(h_{1},h_{2})_{H^{\sigma,\varphi}(\Gamma)}+
\varrho^{2}\,(h_{1},h_{2})_{H^{\sigma-\theta,\varphi}(\Gamma)}.
$$
Следовательно, это пространство гильбертово. Как и прежде, в соболевском случае
$\varphi\equiv1$ индекс $\varphi$ в обозначениях опускаем. Возвращаясь к
формулировке утверждения (ii) теоремы \ref{th2.33b}, заметим, что
\begin{gather*}
\|u\|_{H^{s+mq,\varphi}(\Gamma,|\lambda|^{q},mq)}\leq \\
\|u\|_{H^{s+mq,\varphi}(\Gamma)}+
|\lambda|^{q}\,\|u\|_{H^{s,\varphi}(\Gamma)}\leq\\
\sqrt{2}\;\|u\| _{H^{s+mq,\varphi}(\Gamma,|\lambda|^{q},mq)}.
\end{gather*}

В силу теоремы \ref{th2.21} пространства
$$
\bigl[\,H^{\sigma-\varepsilon}(\Gamma,\varrho,\theta),\,
H^{\sigma+\delta}(\Gamma,\varrho,\theta)\,\bigr]_{\psi}\quad\mbox{и}\quad
H^{\sigma,\varphi}(\Gamma,\varrho,\theta)
$$
равны с точностью до эквивалентности норм. Здесь числа $\varepsilon$ и $\delta$
положительные, а функциональный параметр $\psi$ такой, как в теореме \ref{th2.21}.
Оказывается, что в двусторонних оценках норм этих пространств постоянные можно
выбрать так, чтобы они не зависели от параметра~$\varrho$.


\begin{lemma}\label{lem2.14}
Пусть $\sigma\in\mathbb{R}$, $\varphi\in\mathcal{M}$, и заданы положительные числа
$\theta$, $\varepsilon$, $\delta$. Тогда существует число $c_{0}\geq1$ такое, что
для произвольных $\varrho>0$ и $h\in H^{\sigma,\varphi}(\Gamma)$ верна двусторонняя
оценка норм
\begin{gather}\notag
c_{0}^{-1}\,\|h\|_{H^{\sigma,\varphi}(\Gamma,\varrho,\theta)}\leq
\|h\|_{[\,H^{\sigma-\varepsilon}(\Gamma,\varrho,\theta),\,
H^{\sigma+\delta}(\Gamma,\varrho,\theta)\,]_{\psi}}\leq \\
c_{0}\,\|h\|_{H^{\sigma,\varphi}(\Gamma,\varrho,\theta)}.\label{2.129}
\end{gather}
Здесь $\psi$ --- интерполяционный параметр из теоремы $\ref{th2.21}$, а число
$c_{0}$ не зависит от $\varrho$ и $h$.
\end{lemma}


\textbf{Доказательство.} Пусть параметр $\varrho>0$. Мы сначала докажем аналог
оценки \eqref{2.129} для пространств распределений в $\mathbb{R}^{n}$, а затем с
помощью операторов „распрямления” и „склейки” перейдем к пространствам распределений
на многообразии $\Gamma$ (сравнить с доказательством теоремы \ref{th2.21}).
Обозначим через $H^{\sigma,\varphi}(\mathbb{R}^{n},\varrho,\theta)$ пространство
$H^{\sigma,\varphi}(\mathbb{R}^{n})$, наделенное гильбертовой нормой
\begin{gather}\notag
\|w\|_{H^{\sigma,\varphi}(\mathbb{R}^{n},\varrho,\theta)}:=
\left(\;\|w\|_{H^{\sigma,\varphi}(\mathbb{R}^{n})}^{2}+ \varrho^{2}\,\|w\|
_{H^{\sigma-\theta,\varphi}(\mathbb{R}^{n})}^{2}\:\right)^{1/2}=\\
\biggl(\;\;\int\limits_{\mathbb{R}^{n}}\,\langle\xi\rangle^{2\sigma}\,
\bigl(1+\varrho^{2}\,\langle\xi\rangle^{-2\theta}\bigr)\,
\varphi^{2}(\langle\xi\rangle)\,|\widehat{w}(\xi)|^{2}\,d\xi\,\biggr)^{1/2}.
\label{2.130}
\end{gather}
Для каждого фиксированного $\varrho>0$ эта норма эквивалентна норме в пространстве
$H^{\sigma,\varphi}(\mathbb{R}^{n})$. Следовательно, пространство
$H^{\sigma,\varphi}(\mathbb{R}^{n},\varrho,\theta)$ гильбертово. Аналогично
определяются пространства $H^{\sigma-\varepsilon}(\mathbb{R}^{n},\varrho,\theta)$ и
$H^{\sigma+\delta}(\mathbb{R}^{n},\varrho,\theta)$. В силу теоремы \ref{th2.14}
пространства
\begin{equation}\label{2.131}
\bigl[\,H^{\sigma-\varepsilon}(\mathbb{R}^{n},\varrho,\theta),\,
H^{\sigma+\delta}(\mathbb{R}^{n},\varrho,\theta)\,\bigr]_{\psi}
\end{equation}
и $H^{\sigma,\varphi}(\mathbb{R}^{n},\varrho,\theta)$ равны с точностью до
эквивалентности норм при каждом фиксированном $\varrho>0$. Покажем, что в
двусторонних оценках норм этих пространств можно взять постоянные, не зависящие от
параметра~$\varrho$.

Вычислим норму в пространстве \eqref{2.131}. Обозначим через $J$
псевдодифференциальный оператор в пространстве $\mathbb{R}^{n}$ с символом
$\langle\xi\rangle^{\varepsilon+\delta}$, где аргумент $\xi\in\mathbb{R}^{n}$.
Непосредственно проверяется, что оператор $J$ порождающий для пары пространств,
записанной в формуле \eqref{2.131}. С помощью изометрического изоморфизма
$$
\mathcal{F}:H^{\sigma-\varepsilon}(\mathbb{R}^{n},\varrho,\theta)\leftrightarrow
L_{2}\bigl(\mathbb{R}^{n},\langle\xi\rangle^{2(\sigma-\varepsilon)}
(1+\varrho^{2}\,\langle\xi\rangle^{-2\theta})\,d\xi\bigr),
$$
где $\mathcal{F}$ --- преобразование Фурье, оператор $J$ приводится к виду умножения
на функцию $\langle\xi\rangle^{\varepsilon+\delta}$, а оператор $\psi(J)$ --- к виду
умножения на функцию $\psi(\langle\xi\rangle^{\varepsilon+\delta})=
\langle\xi\rangle^{\varepsilon}\varphi(\langle\xi\rangle)$. Следовательно,
\begin{gather*}
\|w\|^{2} _{[H^{\sigma-\varepsilon}(\mathbb{R}^{n},\varrho,\theta),\,
H^{\sigma+\delta}(\mathbb{R}^{n},\varrho,\theta)]_{\psi}}=
\|\psi(J)\,w\|_{H^{\sigma-\varepsilon}(\mathbb{R}^{n},\varrho,\theta)}^{2}=\\
\int\limits_{\mathbb{R}^{n}}\,\langle\xi\rangle^{2(\sigma-\varepsilon)}
\bigl(1+\varrho^{2}\,\langle\xi\rangle^{-2\theta}\bigr)\:
\bigl|\langle\xi\rangle^{\varepsilon}\varphi\bigl(\langle\xi\rangle\bigr)\,
\widehat{w}(\xi)\bigr|^{2}\,d\xi=
\|w\|_{H^{\sigma,\varphi}(\mathbb{R}^{n},\varrho,\theta)}^{2}.
\end{gather*}
Таким образом,
\begin{equation}\label{2.132}
\|w\|_{[H^{\sigma-\varepsilon}(\mathbb{R}^{n},\varrho,\theta),\,
H^{\sigma+\delta}(\mathbb{R}^{n},\varrho,\theta)]_{\psi}}=
\|w\|_{H^{\sigma,\varphi}(\mathbb{R}^{n},\varrho,\theta)}
\end{equation}
для произвольных $w\in H^{\sigma,\varphi}(\mathbb{R}^{n})$ и $\varrho>0$.

Теперь выведем неравенство \eqref{2.129} из \eqref{2.132}. Мы используем локальное
определение \ref{def2.13} пространств $H^{s,\varphi}(\Gamma)$, $s\in\mathbb{R}$,
$\varphi\in\mathcal{M}$, для фиксированных конечного атласа и разбиения единицы на
$\Gamma$. Рассмотрим линейное отображение „распрямления” многообразия~$\Gamma$:
$$
T:\,h\mapsto\bigl(\,(\chi_{1}h)\circ\alpha_{1},\ldots,
(\chi_{r}h)\circ\alpha_{r}\,\bigr),\quad h\in\mathcal{D}'(\Gamma).
$$
Непосредственно проверяется, что это отображение задает изометрические операторы
\begin{gather}\label{2.133}
T:\,H^{\sigma,\varphi}(\Gamma,\varrho,\theta)\rightarrow
(H^{\sigma,\varphi}(\mathbb{R}^{n},\varrho,\theta))^{r}, \\
T:\,H^{s}(\Gamma,\varrho,\theta)\rightarrow(H^{s}(\mathbb{R}^{n},\varrho,\theta))^{r},
\quad s\in\{\sigma-\varepsilon,\sigma+\delta\}. \label{2.134}
\end{gather}
Применив к операторам \eqref{2.134} интерполяцию с параметром $\psi$, получим
ограниченный оператор
\begin{gather}\notag
T:\,\bigl[\,H^{\sigma-\varepsilon}(\Gamma,\varrho,\theta),\,
H^{\sigma+\delta}(\Gamma,\varrho,\theta)\,\bigr]_{\psi}\,\rightarrow \\
\bigl[\;(H^{\sigma-\varepsilon}(\mathbb{R}^{n},\varrho,\theta))^{r},\,
(H^{\sigma+\delta}(\mathbb{R}^{n},\varrho,\theta))^{r}\,\bigr]_{\psi}. \label{2.135}
\end{gather}
Поскольку записанные здесь пары пространств нормальные, то в силу теоремы
\ref{th2.8} норма оператора \eqref{2.135} не превышает некоторого числа $c$, не
зависящего от параметра $\varrho$. Отсюда на основании теоремы \ref{th2.5} и
равенства \eqref{2.132} получаем ограниченный оператор
\begin{gather}\notag
T:\,\bigl[\,H^{\sigma-\varepsilon}(\Gamma,\varrho,\theta),\,
H^{\sigma+\delta}(\Gamma,\varrho,\theta)\,\bigr]_{\psi}\,\rightarrow \\
(H^{\sigma,\varphi}(\mathbb{R}^{n},\varrho,\theta))^{r} \quad\mbox{с нормой}\;\leq
c.\label{2.136}
\end{gather}

Далее, наряду с отображением $T$ рассмотрим линейное отображение „склейки”
$$
K:\,(w_{1},\ldots,w_{r})\mapsto\sum_{j=1}^{r}\Theta_{j}\bigl(\,
(\eta_{j}w_{j})\circ\alpha_{j}^{-1}\bigr),
$$
где $w_{1},\ldots,w_{r}$ --- распределения в $\mathbb{R}^{n}$. Здесь функция
$\eta_{j}\in C^{\infty}(\mathbb{R}^{n})$ финитная и равная 1 на множестве
$\alpha_{j}^{-1}(\mathrm{supp}\,\chi_{j})$, а $\Theta_{j}$ --- оператор продолжения
нулем на $\Gamma$. В силу \eqref{2.70} имеем ограниченные операторы
\begin{gather}\label{2.137}
K:\,(H^{s}(\mathbb{R}^{n}))^{r}\rightarrow H^{s}(\Gamma),\quad
s\in\mathbb{R},\\
K:\,(H^{s,\varphi}(\mathbb{R}^{n}))^{r}\rightarrow H^{s,\varphi}(\Gamma),\quad
s\in\mathbb{R}. \label{2.138}
\end{gather}
Пусть $c_{2}$ --- максимум норм операторов \eqref{2.137}, где
$$
s\in\{\sigma-\varepsilon,\,\sigma-\varepsilon-\theta,\,\sigma+\delta,\,
\sigma+\delta-\theta\},
$$
и операторов \eqref{2.138}, где $s\in\{\sigma,\,\sigma-\theta$\}. Число $c_{2}$ не
зависит от параметра $\varrho$. Непосредственно проверяется, что нормы операторов
\begin{gather}\label{2.139}
K:\,(H^{\sigma,\varphi}(\mathbb{R}^{n},\varrho,\theta))^{r} \rightarrow
H^{\sigma,\varphi}(\Gamma,\varrho,\theta),\\
K:\,(H^{s}(\mathbb{R}^{n},\varrho,\theta))^{r}\rightarrow
H^{s}(\Gamma,\varrho,\theta), \quad s\in\{\sigma-\varepsilon,\sigma+\delta\},
\label{2.140}
\end{gather}
не превышают числа $c_{2}$. Применив к \eqref{2.140} интерполяцию с параметром
$\psi$, получим ограниченный оператор
\begin{gather*}
K:\,\bigl[\,(H^{\sigma-\varepsilon}(\mathbb{R}^{n},\varrho,\theta))^{r},\,
(H^{\sigma+\delta}(\mathbb{R}^{n},\varrho,\theta))^{r}\,\bigr]_{\psi}
\,\rightarrow\\
\bigl[\,H^{\sigma-\varepsilon}(\Gamma,\varrho,\theta),\,
H^{\sigma+\delta}(\Gamma,\varrho,\theta)\,\bigr]_{\psi},
\end{gather*}
норма которого не превышает числа $c\,c_{2}$ в силу теоремы \ref{th2.8}. Отсюда на
основании теоремы \ref{th2.5} и равенства \eqref{2.132} получаем ограниченный
оператор
\begin{gather}\notag
K:\,(H^{\sigma,\varphi}(\mathbb{R}^{n},\varrho,\theta))^{r}\,\rightarrow \\
\bigl[\,H^{\sigma-\varepsilon}(\Gamma,\varrho,\theta),\,
H^{\sigma+\delta}(\Gamma,\varrho,\theta)\,\bigr]_{\psi}\quad\mbox{с нормой}\;\leq
c_{3},\label{2.141}
\end{gather}
где число $c_{3}:=\sqrt{c_{1}}\,c\,c_{2}$ не зависит от параметра $\varrho$.

В силу \eqref{2.69}, произведение $KT=I$ --- тождественный оператор. Поэтому на
основании формул \eqref{2.133} (изометрический оператор) и  \eqref{2.141} имеем
ограниченный оператор
$$
I=KT:\,H^{\sigma,\varphi}(\Gamma,\varrho,\theta)\rightarrow
\bigl[\,H^{\sigma-\varepsilon}(\Gamma,\varrho,\theta),\,
H^{\sigma+\delta}(\Gamma,\varrho,\theta)\,\bigr]_{\psi}
$$
с нормой $\leq c_{3}$. Кроме того, в силу \eqref{2.136} и \eqref{2.139} (норма
второго оператора не превышает числа $c_{2}$) имеем еще один ограниченный оператор
$$
I=KT:\,\bigl[\,H^{\sigma-\varepsilon}(\Gamma,\varrho,\theta),\,
H^{\sigma+\delta}(\Gamma,\varrho,\theta)\,\bigr]_{\psi} \rightarrow
H^{\sigma,\varphi}(\Gamma,\varrho,\theta)
$$
с нормой $\leq c_{2}\,c$. Отсюда немедленно вытекает двусторонняя оценка
\eqref{2.129}, где число $c_{0}:=\max\{1,\,c_{3},\,c_{2},\,c\}$ не зависит от
параметра $\varrho$.

Лемма \ref{lem2.14} доказана.

\medskip

\textbf{Доказательство теоремы~\ref{th2.33b} (ii).} Пусть $s\in\mathbb{R}$ и
$\varphi\in\mathcal{M}$. Напомним, что теорема \ref{th2.33b} справедлива в
соболевском случае $\varphi\equiv1$. Поэтому существует такое число $\lambda_{0}>0$,
что для каждого значения параметра $\lambda\in K$, удовлетворяющего условию
$|\lambda|\geq \lambda_{0}$, справедливы гомеоморфизмы
\begin{equation}\label{2.142b}
A(\lambda):H^{s\mp1+mq}(\Gamma,|\lambda|^{q},mq)\leftrightarrow H^{s\mp1}(\Gamma),
\end{equation}
причем норма оператора \eqref{2.142b} вместе с нормой обратного оператора ограничены
равномерно по параметру $\lambda$. Пусть $\psi$~--- интерполяционный параметр из
теоремы \ref{th2.21}, где берем $\varepsilon=\delta=1$. Применив интерполяцию с этим
параметром к \eqref{2.142b}, получим гомеоморфизм
\begin{gather}\notag
A(\lambda):\,\bigl[\,H^{s-1+mq}(\Gamma,|\lambda|^{q},mq),\,
H^{s+1+mq}(\Gamma,|\lambda|^{q},mq)\,\bigr]_{\psi}\,\leftrightarrow\\
\bigl[\,H^{s-1}(\Gamma),\,H^{s+1}(\Gamma)\,\bigr]_{\psi}. \label{2.143b}
\end{gather}
При этом в силу теоремы \ref{th2.8} норма оператора \eqref{2.143b} и норма обратного
к нему оператора ограничены равномерно по параметру~$\lambda$. (Записанные в формуле
\eqref{2.143b} допустимые пары пространств нормальные.) Остается воспользоваться
леммой \ref{lem2.14}, где полагаем
$$
\sigma:=s+mq,\;\;\varrho:=|\lambda|^{q},\;\;\theta:=mq,\;\;\varepsilon=\delta=1,
$$
и теоремой \ref{th2.21}. Согласно ним \eqref{2.143b} влечет гомеоморфизм
\begin{equation}\label{2.145b}
A(\lambda):H^{s+mq,\varphi}(\Gamma,|\lambda|^{q},mq)\leftrightarrow
H^{s,\varphi}(\Gamma)
\end{equation}
такой, что норма оператора \eqref{2.145b} вместе с нормой обратного оператора
ограничены равномерно по параметру $\lambda$. Это означает двустороннюю оценку
\eqref{2.126b}, где число $c$ не зависит от параметра $\lambda$ и от распределения
$u\in H^{s+mq,\varphi}(\Gamma)$. Утверждение (ii) теоремы~\ref{th2.33b} доказано.

\medskip

Из теоремы \ref{th2.33b} (i) вытекает следующее утверждение об индексе
эллиптического с параметром ПДО (сравнить с [\ref{AgranovichVishik64}, с.~96]).


\begin{corollary}\label{cor2.5b}
Предположим, что ПДО $A(\lambda)$ эллиптический с параметром на некотором замкнутом
луче $K:=\{\lambda\in\mathbb{C}:\,\arg\lambda =\mathrm{const}\}$. Тогда оператор
\eqref{2.124b} имеет нулевой индекс при любом $\lambda\in\mathbb{C}$.
\end{corollary}


\textbf{Доказательство.} Для каждого фиксированного $\lambda\in\mathbb{C}$ ПДО
$A(\lambda)$ эллиптический на $\Gamma$. Поэтому в силу теоремы \ref{th2.28b}
оператор \eqref{2.124b} имеет конечный индекс, не зависящий от $s\in\mathbb{R}$ и
$\varphi\in\mathcal{M}$. Кроме того, этот индекс не зависит и от параметра
$\lambda$. В самом деле, $\lambda$ влияет лишь на младшие члены в сумме
\eqref{2.120b}: $A(\lambda)-A(0)\in\Psi^{m(q-1)}$. Поэтому ввиду леммы
\ref{lem2.13b} имеем ограниченный оператор
$$
A(\lambda)-A(0):H^{s+mq,\varphi}(\Gamma)\rightarrow H^{s+m,\varphi}(\Gamma).
$$
Но в силу теоремы \ref{th2.22} (iii) и условия $m>0$ справедливо компактное вложение
$H^{s+m,\varphi}(\Gamma)\hookrightarrow H^{s,\varphi}(\Gamma)$. Следовательно,
компактен оператор
$$
A(\lambda)-A(0):H^{s+mq,\varphi}(\Gamma)\rightarrow H^{s,\varphi}(\Gamma).
$$
Отсюда вытекает (см., например, [\ref{Hermander87}, с.~249]), что нетеровы операторы
$A(\lambda)$ и $A(0)$ имеют одинаковый индекс, т.~е. последний не зависит от
параметра $\lambda$. Далее, согласно теореме \ref{th2.33b} (i) при достаточно
больших по модулю значениях параметра $\lambda\in K$, справедлив гомеоморфизм
\eqref{2.125b}. Следовательно, индекс оператора $A(\lambda)$ равен нулю при
$\lambda\in K$, $|\lambda|\gg1$, а значит, и при любом $\lambda\in\mathbb{C}$.
Следствие \ref{cor2.5b} доказано.




\markright{\emph \ref{sec2.7}. Сходимость спектральных разложений}

\section{Сходимость спектральных разложений}\label{sec2.7}

\markright{\emph \ref{sec2.7}. Сходимость спектральных разложений}

В этом пункте мы дадим приложение введенных нами пространств Хермандера к
исследованию сходимости спектральных разложений почти всюду и в метриках пространств
$C^{k}$, где целое $k\geq0$.

\subsection[Сходимость почти всюду произвольных ортогональных рядов]
{Сходимость почти всюду \\ произвольных ортогональных рядов}\label{sec2.7.1}

Сформулируем сначала необходимые нам классические результаты о сходимости почти
всюду произвольных ортогональных рядов. Пусть $\Gamma$~--- измеримое пространство с
конечной мерой $\mu$, и $(h_{j})_{j=1}^{\infty}$ --- ортонормированная система
функций в $L_{2}(\Gamma,\mu)$. (Функции $h_{j}$, вообще говоря, комплекснозначные.)
Мы используем обозначение $\log$ для логарифмической функции при произвольном
фиксированном основании $a>1$.


\begin{proposition}[общая теорема Меньшова-Ра\-де\-ма\-хе\-ра]\label{prop2.7}
Пусть последовательность $(a_{j})_{j=1}^{\infty}$ комплексных чисел такая, что
\begin{equation}\label{2.146}
L:=\sum_{j=1}^{\infty}\;|a_{j}|^{2}\,\log^{2}\,(j+1)<\infty.
\end{equation}
Тогда ортогональный ряд
\begin{equation}\label{2.147}
\sum_{j=1}^{\infty}\;a_{j}\,h_{j}(x)
\end{equation}
сходится $\mu$-почти всюду на $\Gamma$. Кроме того, если
$$
S^{\ast}(x):=\sup_{1\leq k<\infty}\;\Bigl|\,\sum_{j=1}^{k}\;a_{j}\,h_{j}(x)\,\Bigr|,
$$
--- мажоранта частных сумм ряда \eqref{2.147}, то
\begin{equation}\label{2.148}
\|S^{\ast}\|_{L_{2}(\Gamma)}\leq C\,\sqrt{L},
\end{equation}
где число $C>0$ зависит лишь от $\Gamma$ и $\mu$.
\end{proposition}


Этот результат доказан независимо Д.~Е.~Меньшовым [\ref{Menschoff23}] и
Г.~Радемахером [\ref{Rademacher22}] в случае, когда $\Gamma$ --- конечный интервал
на оси $\mathbb{R}$, $\mu$~--- мера Лебега, и функции $h_{j}$ вещественные.
Приведенное в монографии [\ref{KashinSaakyan84}, с. 292 -- 294] доказательство
теоремы Меньшова-Радемахера проходит в рассматриваемой нами общей ситуации.

Заметим, что теорема Меньшова-Радемахера точна. Д.~Е.~Мень\-шов [\ref{Menschoff23}]
построил пример ортонормированной системы $(h_{j})_{j=1}^{\infty}$ в $L_{2}((0,1))$
такой, что для каждой числовой последовательности $(\omega_{j})_{j=1}^{\infty}$,
удовлетворяющей условиям
$$
1=\omega_{1}\leq\omega_{2}\leq\omega_{3}\leq\ldots\quad\mbox{и}\quad
\lim_{j\rightarrow\infty}\,\frac{\omega_{j}}{\log^{2}j}=0,
$$
найдется расходящийся почти всюду ряд вида \eqref{2.147}, для коэффициентов которого
выполняется неравенство
$$
\sum_{j=1}^{\infty}\;|a_{j}|^{2}\,\omega_{j}<\infty.
$$
(Изложение этого результата имеется, например, в [\ref{KashinSaakyan84}, с. 295]).

Отметим [\ref{KashinSaakyan84}, с. 307], что из сходимости ряда \eqref{2.147}
$\mu$-почти всюду не вытекает, вообще говоря, что этот ряд безусловно сходится
$\mu$-почти всюду. Напомним, что ряд \eqref{2.147} называется  безусловно сходящимся
$\mu$-почти всюду на $\Gamma$, если для любой перестановки натурального ряда
$\sigma=(\sigma(j))_{j=1}^{\infty}$ сходится $\mu$-почти всюду ряд
\begin{equation}\label{2.149}
\sum_{j=1}^{\infty}\;a_{\sigma(j)}\,h_{\sigma(j)}.
\end{equation}
(При этом множество меры нуль расходимости ряда \eqref{2.149} может зависить от
перестановки $\sigma$.)


\begin{proposition}[общая теорема Орлича-Ульянова]\label{prop2.8}
Пусть последовательность $(a_{j})_{j=1}^{\infty}$ комплексных чисел и возрастающая
(нестрого) последовательность $(\omega_{j})_{j=1}^{\infty}$ положительных чисел
удовлетворяют условиям
\begin{gather}\label{2.150}
\sum_{j=2}^{\infty}\;|a_{j}|^{2}\,(\log^{2}j)\,\omega_{j}<\infty,\\
\sum_{j=2}^{\infty}\;\frac{1}{j\,(\log j)\,\omega_{j}}<\infty. \label{2.151}
\end{gather}
Тогда ряд \eqref{2.147} безусловно сходится $\mu$-почти всюду на $\Gamma$.
\end{proposition}


Это --- эквивалентная формулировка теоремы Орлича [\ref{Orlicz27}], полученная
П.~Л.~Ульяновым [\ref{Ulyanov63}, с. 53] (см. также [\ref{Ulyanov64}, с. 53]). Как
показал К.~Тандори [\ref{Tandori62}], теорема Орлича окончательная в том смысле, что
условие \eqref{2.150} на последовательность $(\omega_{j})_{j=1}^{\infty}$ нельзя
ослабить.

В.~Орлич и П.~Л.~Ульянов ограничились случаем, когда $\Gamma$~--- конечный интервал
на оси $\mathbb{R}$, $\mu$~--- мера Лебега, а функции $h_{j}$ вещественные. В
рассматриваемом нами общем случае предложение \ref{prop2.8} остается верным. Это
следует из более общей теоремы Тандори [\ref{Tandori62}], доказательство которой,
данное в монографии [\ref{KashinSaakyan84}, с. 308], проходит в рассматриваемой нами
общей ситуации.

\subsection[Сходимость почти всюду спектральных разложений]
{Сходимость почти всюду \\ спектральных разложений}\label{sec2.7.2}

Изучим теперь сходимость почти всюду спектральных разложений по собственным функциям
эллиптических ПДО. Далее в этом пункте $\Gamma$ --- бесконечно гладкое замкнутое
(компактное) многообразие размерности $n\geq1$, на котором фиксирована плотность
$dx$. Пусть $A$ --- классический эллиптический ПДО на $\Gamma$ положительного
порядка. Предположим, что $A$ является (неограниченным) нормальным оператором в
гильбертовом пространстве $L_{2}(\Gamma,dx)$. Пусть $(h_{j})_{j=1}^{\infty}$ ---
полная ортонормированная система собственных функций этого оператора. Они
занумерованы так, что модули соответствующих собственных чисел образуют монотонно
неубывающую последовательность. Для произвольной функции $f\in L_{2}(\Gamma,dx)$
имеем
\begin{equation}\label{2.152}
f=\sum_{j=1}^{\infty}\;c_{j}(f)\,h_{j}\quad\mbox{в}\;\;L_{2}(\Gamma,dx),
\end{equation}
где $c_{j}(f):=(f,h_{j})_{\Gamma}$ --- коэффициенты Фурье распределения $f$ в базисе
$(h_{j})_{j=1}^{\infty}$.

Ряд \eqref{2.152} называется сходящимся в указанном смысле на некотором
функциональном классе $X(\Gamma)$, если для любой функции $f\in X(\Gamma)$ этот ряд
сходится к $f$ надлежащим образом.

Изучим сходимость почти всюду ряда \eqref{2.152} на классах~--- пространствах
Хермандера. Заметим, что, если для некоторой функции $f\in L_{2}(\Gamma,dx)$ ряд
\eqref{2.152} сходится почти всюду на $\Gamma$, то $f$ будет суммой этого ряда почти
всюду на $\Gamma$. Это вытекает из того, что как сходимость в $L_{2}(\Gamma,dx)$,
так и сходимость почти всюду на $\Gamma$ влекут за собой сходимость по мере,
порожденной плотностью $dx$ (см., например, [\ref{BerezanskyUsSheftel90}, с. 70,
186]).

Рассмотрим мажоранту частных сумм ряда \eqref{2.152} для функции $f\in
L_{2}(\Gamma,dx)$:
$$
S^{\ast}(f,x):=\sup_{1\leq
k<\infty}\;\Bigl|\,\sum_{j=1}^{k}\;c_{j}(f)\,h_{j}(x)\,\Bigr|.
$$

В статье К.~Мини [\ref{Meaney82}] доказано для случая $A=\Delta_{\Gamma}$, что ряд
\eqref{2.152} сходится почти всюду на $\Gamma$ на функциональном классе
$$
H^{0+}(\Gamma):= \bigcup_{\varepsilon>0}H^{\varepsilon}(\Gamma).
$$
При этом для произвольных $\varepsilon>0$ и $f\in H^{\varepsilon}(\Gamma)$
выполняется оценка
$$
\|S^{\ast}(f,\cdot)\|_{L_{2}(\Gamma,dx)}\leq
C_{\varepsilon}\|f\|_{H^{\varepsilon}(\Gamma)},
$$
где число $C_{\varepsilon}$ не зависит от $f$.

Мы обобщим и уточним этот результат, используя пространства Хермандера. Обозначим
$\log^{\ast}t:=\max\{1,\log t\}$.


\begin{theorem}\label{th2.34}
На функциональном классе $H^{0,\log^{\ast}}(\Gamma)$ ряд \eqref{2.152} сходится
почти всюду на $\Gamma$. При этом для любой функции $f\in H^{0,\log^{\ast}}(\Gamma)$
выполняется оценка
\begin{equation}\label{2.153}
\|S^{\ast}(f,\cdot)\|_{L_{2}(\Gamma,dx)}\leq C\,\|f\|_{H^{0,\log^{\ast}}(\Gamma)},
\end{equation}
где число $C>0$ не зависит от $f$.
\end{theorem}


\textbf{Доказательство.} Пусть $f\in H^{0,\log^{\ast}}(\Gamma)$. Если оператор $A$
самосопряженный и положительно определенный в $L_{2}(\Gamma,dx)$, то в силу
теоремы~\ref{th2.26}
$$
\|f\|_{H^{0,\log^{\ast}}(\Gamma)}^{2}\asymp\sum_{j=1}^{\infty}\;
(\log^{\ast}(j^{\,1/n}))^{2}\,|c_{j}(f)|^{2},
$$
где $\log^{\ast}(j^{\,1/n})\asymp\log(j+1)$ при $j\geq1$. Следовательно,
$$
\sum_{j=1}^{\infty}\;|c_{j}(f)|^{2}\,\log^{2}(j+1)\asymp
\|f\|_{H^{0,\log^{\ast}}(\Gamma)}^{2}<\infty.
$$
Поэтому согласно общей теореме Меньшова-Радемахера (предложение \ref{prop2.7}) ряд
\eqref{2.152} сходится почти всюду на $\Gamma$ (к $f$) и справедлива оценка
\eqref{2.153}.

Общая ситуация, когда оператор $A$ нормальный, сводится к рассмотренному случаю,
переходом к самосопряженному и положительно определенному оператору
$B:=1+A^{\ast}A$. При этом следует учесть тот факт, что $(h_{j})_{j=1}^{\infty}$~---
полная в $L_{2}(\Gamma,dx)$ система собственных функций как нормального оператора
$A$, так и самосопряженного оператора $B$, а их нумерация согласована.

Теорема~\ref{th2.34} доказана.

\medskip

Уточнением этого результата является следующая теорема о безусловной сходимости.

\begin{theorem}\label{th2.35}
Пусть возрастающая функция $\varphi\in\mathcal{M}$ удовлетворяет условию
\begin{equation}\label{2.154}
\int\limits_{2}^{\infty}\frac{dt}{t\,(\log t)\,\varphi^{2}(t)}<\infty.
\end{equation}
Тогда на функциональном классе $H^{0,\varphi\log^{\ast}}(\Gamma)$ ряд \eqref{2.152}
безусловно сходится почти всюду на $\Gamma$.
\end{theorem}


\textbf{Доказательство.} Пусть $f\in H^{0,\varphi\log^{\ast}}(\Gamma)$. Если
оператор $A$ самосопряженный и положительно определенный в $L_{2}(\Gamma,dx)$, то по
теореме \ref{th2.26} можем записать
\begin{equation}\label{2.155}
\|f\|_{H^{0,\varphi\log^{\ast}}(\Gamma)}^{2}\asymp\sum_{j=1}^{\infty}\;
\varphi^{2}(j^{\,1/n})(\log^{\ast}(j^{\,1/n}))^{2}\,|c_{j}(f)|^{2}.
\end{equation}
Рассмотрим возрастающую последовательность чисел
$\omega_{j}:=\varphi^{2}(j^{\,1/n})$, $j\in\mathbb{N}$. В силу \eqref{2.155} имеем
\begin{equation}\label{2.156}
\sum_{j=2}^{\infty}\;|c_{j}(f)|^{2}\,(\log^{2}j)\,\omega_{j}<\infty.
\end{equation}
Кроме того, согласно условию \eqref{2.154},
\begin{gather}\notag
\sum_{j=3}^{\infty}\;\frac{1}{j\,(\log j)\,\omega_{j}}\leq
\int\limits_{2}^{\infty}\frac{d\tau}{\tau\,(\log\tau)\,\varphi^{2}(\tau^{1/n})}= \\
\int\limits_{2^{1/n}}^{\infty}\frac{n\,t^{n-1}\,dt}{t^{n}\,n\,(\log
t)\,\varphi^{2}(t)}<\infty.\label{2.157}
\end{gather}
В силу общей теоремы Орлича-Ульянова (предложение \ref{prop2.8}), неравенства
\eqref{2.156} и \eqref{2.157} влекут за собой безусловную сходимость почти всюду на
$\Gamma$ ряда \eqref{2.152}.

Как и в доказательстве теоремы \ref{th2.34}, общая ситуация, когда оператор $A$
нормальный, сводится к рассмотренному выше случаю, переходом к самосопряженному
оператору $B:=1+A^{\ast}A$.

Теорема \ref{th2.35} доказана.

\medskip

Теоремы \ref{th2.34} и \ref{th2.35} полностью вкладываются в теорию общих
ортогональных рядов. Однако в них даны условия на функцию $f$ в конструктивных
терминах гладкости функций.

\subsection[Сходимость спектральных разложений в метрике пространства $C^{k}$]
{Сходимость спектральных разложений \\ в метрике пространства
$C^{k}$}\label{sec2.7.3}

В заключение п.~\ref{sec2.7} докажем критерий сходимости ряда \eqref{2.152} в
метрике пространств $C^{k}(\Gamma)$, где целое $k\geq0$, на классах
$H^{s,\varphi}(\Gamma)$.


\begin{theorem}\label{th2.36}
Пусть заданы целое число $k\geq0$ и функция $\varphi\in\mathcal{M}$. Ряд
\eqref{2.152} сходится в метрике пространства $C^{k}(\Gamma)$ на классе
$H^{k+n/2,\varphi}(\Gamma)$ тогда и только тогда, когда функция $\varphi$
удовлетворяет условию \eqref{2.37}.
\end{theorem}


\textbf{Доказательство.} \textit{Достаточность.} Предположим, что $\varphi$
удовлетворяет условию \eqref{2.37}. Пусть $f\in H^{k+n/2,\varphi}(\Gamma)$. Как
отмечалось в замечании \ref{rem2.12}, ряд \eqref{2.152} сходится к $f$ в метрике
пространства $H^{k+n/2,\varphi}(\Gamma)$. Согласно теореме \ref{th2.27}, неравенство
\eqref{2.37} влечет непрерывность вложения $H^{k+n/2,\varphi}(\Gamma)\hookrightarrow
C^{k}(\Gamma)$. Следовательно, ряд \eqref{2.152} сходится к $f$ в пространстве
$C^{k}(\Gamma)$. Достаточность доказана.

\textit{Необходимость.} Предположим, что для каждой функции $f\in
H^{k+n/2,\varphi}(\Gamma)$ ряд \eqref{2.152} сходится (к $f$) в метрике пространства
$C^{k}(\Gamma)$. Тогда $H^{k+n/2,\varphi}(\Gamma)\subseteq C^{k}(\Gamma)$, что
влечет условие \eqref{2.37} в силу теоремы \ref{th2.27}. Необходимость доказана.

Теорема \ref{th2.36} доказана.




\markright{\emph \ref{sec2.8}. RO-меняющиеся функции и пространства Хермандера}

\section[RO-меняющиеся функции и пространства Хермандера]
{RO-меняющиеся функции \\ и пространства Хермандера}\label{sec2.8}

\markright{\emph \ref{sec2.8}. RO-меняющиеся функции и пространства Хермандера}

В этом пункте мы опишем \emph{все} (с точностью до эквивалентности норм)
интерполяционные гильбертовы пространства для пар гильбертовых пространств Соболева
$H^{s_{0}}(\mathbb{R}^{n})$ и $H^{s_{1}}(\mathbb{R}^{n})$, где
$s_{0},s_{1}\in\mathbb{R}$ и $s_{0}<s_{1}$. Как выяснится, класс этих
интерполяционных пространств в точности состоит из изотропных пространств Хермандера
с весовыми функциями, RО-меняющимися на бесконечности по Авакумовичу. Мы рассмотрим
указанные пространства также на гладком замкнутом многообразии и дадим некоторые их
приложения.

\subsection{RO-меняющиеся функции по Авакумовичу}\label{sec2.8.1}

Приведем базисное для нас определение.


\begin{definition}\label{def2.19}
Пусть $\mathrm{RO}$ --- множество всех измеримых по Борелю функций
$\varphi:\nobreak[1,\infty)\rightarrow(0,\infty)$, для которых существуют числа
$a>1$ и $c\geq1$ такие, что
\begin{equation}\label{2.158}
c^{-1}\leq\frac{\varphi(\lambda t)}{\varphi(t)}\leq c\quad\mbox{для любых}\quad
t\geq1,\;\lambda\in[1,a]
\end{equation}
(постоянные $a$ и $c$ зависят от $\varphi$). Такие функции называют
\emph{RO-меняющимися} на бесконечности.
\end{definition}


Класс RO-меняющихся функций введен В.~Авакумовичем [\ref{Avakumovic36}] в 1936~г. и
достаточно полно изучен (см., например, [\ref{Seneta85}, с. 86 -- 99]).

Отметим некоторые известные свойства функций класса $\mathrm{RO}$ [\ref{Seneta85},
с. 87 -- 89]. Ясно, что для функции $\varphi\in\mathrm{RO}$ число $a$ в формуле
\eqref{2.158} можно брать сколь угодно большим. Отсюда вытекает следующее свойство.


\begin{proposition}\label{prop2.9}
Всякая функция $\varphi\in\mathrm{RO}$ ограничена и отделена от нуля на каждом
отрезке $[1,b]$, где $1<b<\infty$.
\end{proposition}


\begin{proposition}\label{prop2.10}
Справедливо следующее описание класса $\mathrm{RO}$:
$$
\varphi\in\mathrm{RO}\;\Leftrightarrow\;\varphi(t)=
\exp\Biggl(\beta(t)+\int\limits_{1}^{t}\frac{\varepsilon(\tau)}{\tau}\;d\tau\Biggr),
\;\;t\geq1,
$$
где вещественные функции $\beta$ и $\varepsilon$ измеримы по Борелю и ограничены на
полуоси $[1,\infty)$.
\end{proposition}


\begin{proposition}\label{prop2.11}
Для произвольной фунции $\varphi:[1,\infty)\rightarrow(0,\infty)$ условие
\eqref{2.158} равносильно следующему: существуют числа $s_{0},s_{1}\in\mathbb{R}$,
$s_{0}\leq s_{1}$, и $c\geq1$ такие, что
\begin{equation}\label{2.159}
t^{-s_{0}}\varphi(t)\leq
c\,\tau^{-s_{0}}\varphi(\tau),\;\;\tau^{-s_{1}}\varphi(\tau)\leq
c\,t^{-s_{1}}\varphi(t)
\end{equation}
для всех $t\geq1$ и $\tau\geq t$.
\end{proposition}


Условие \eqref{2.159} показывает, что функция $t^{-s_{0}}\varphi(t)$ эквивалентна
некоторой возрастающей функции, а функция $t^{-s_{1}}\varphi(t)$ эквивалентна
некоторой убывающей функции на полуоси $[1,\infty)$. При этом эквивалетность
понимается в слабом смысле, а возрастание и убывание~--- в нестрогом смысле.

Положив $\lambda:=\tau/t$ в условии \eqref{2.159}, перепишем его в эквивалентном
виде
\begin{equation}\label{2.160}
c^{-1}\lambda^{s_{0}}\leq\frac{\varphi(\lambda t)}{\varphi(t)}\leq
c\lambda^{s_{1}}\quad\mbox{для всех}\quad t\geq1,\;\lambda\geq1.
\end{equation}

Для произвольной функции $\varphi\in\mathrm{RO}$ обозначим:
\begin{gather*}
\sigma_{0}(\varphi):=\sup\,\{s_{0}\in\mathbb{R}:\,\mbox{выполняется
\eqref{2.160}}\},\\
\sigma_{1}(\varphi):=\inf\,\{s_{1}\in\mathbb{R}:\,\mbox{выполняется
\eqref{2.160}}\}.
\end{gather*}
Очевидно, что $-\infty<\sigma_{0}(\varphi)\leq\sigma_{1}(\varphi)<\infty$.

Доопределив функцию $\varphi(t):=\varphi^{2}(1)/\varphi(t^{-1})$ при $0<t<1$,
получим положительную на полуоси $(0,\infty)$ функцию $\varphi$, для которой
\eqref{2.160} влечет условие
$$
c^{-2}\lambda^{s_{0}}\leq\frac{\varphi(\lambda t)}{\varphi(t)}\leq
c^{2}\lambda^{s_{1}} \quad\mbox{для любых}\quad t>0,\;\lambda\geq1.
$$
Отсюда следует, что числа $\sigma_{0}(\varphi)$ и $\sigma_{1}(\varphi)$ равны
соответственно нижнему и верхнему показателям растяжения функции $\varphi$
[\ref{KreinPetuninSemenov78}, с.~76]. Их можно вычислить по формулам
[\ref{KreinPetuninSemenov78}, с.~75]
\begin{equation}\label{2.161}
\sigma_{0}(\varphi)= \lim_{\lambda\rightarrow0+}\,\frac{\log
m_{\varphi}(\lambda)}{\log\lambda},\quad \sigma_{1}(\varphi)=
\lim_{\lambda\rightarrow\infty}\,\frac{\log m_{\varphi}(\lambda)}{\log\lambda},
\end{equation}
где
$$
m_{\varphi}(\lambda):=\sup_{t>0}\,\frac{\varphi(\lambda t)}{\varphi(t)},\quad
\lambda>0,
$$
--- функция растяжения функции $\varphi$. Отметим [\ref{Boyd67}], что правые части
формул \eqref{2.161} равны по определению нижнему и верхнему индексам Бойда
функции~$m_{\varphi}$.

В важном случае $\sigma_{0}(\varphi)=\sigma_{1}(\varphi)=:\sigma$ число $\sigma$
называют \emph{порядком изменения} функции~$\varphi$. Отметим, что всякая измеримая
по Борелю функция $\varphi:[1,\infty)\rightarrow(0,\infty)$, квазиправильно
меняющаяся порядка $\sigma$ на $\infty$, а также ограниченная и отделенная от нуля
на каждом отрезке $[1,b]$, где $1<b<\infty$, принадлежит классу $\mathrm{RO}$ и
имеет порядок изменения $\sigma$. Это вытекает из [\ref{Seneta85}, с. 26, 27] и
предложения~\ref{prop2.11}. В частности, $\mathcal{M}\subset\mathrm{RO}$.

\subsection[Интерполяционные пространства пары соболевских пространств]
{Интерполяционные пространства \\ пары соболевских пространств}\label{sec2.8.2}

В этом пункте мы опишем все интерполяционные гильбертовы пространства для пар
гильбертовых пространств Соболева. Для этого выделим сначала следующий класс
гильбертовых пространств Хермандера.


\begin{definition}\label{def2.20}
Пусть $\varphi\in\mathrm{RO}$. Линейное пространство $H^{\varphi}(\mathbb{R}^{n})$
состоит, по определению, из всех распределений
$u\in\nobreak\mathcal{S}'(\mathbb{R}^{n})$ таких, что их преобразование Фурье
$\widehat{u}$ локально суммируемо по Лебегу в $\mathbb{R}^{n}$ и удовлетворяет
неравенству
$$
\int\limits_{\mathbb{R}^{n}}\varphi^{2}(\langle\xi\rangle)
\:|\widehat{u}(\xi)|^{2}\,d\xi<\infty.
$$
В пространстве $H^{\varphi}(\mathbb{R}^{n})$ определено скалярное произведение
распределений $u_{1}$, $u_{2}$ по формуле
$$
(u_{1},u_{2})_{H^{\varphi}(\mathbb{R}^{n})}:=
\int\limits_{\mathbb{R}^{n}}\varphi^{2}(\langle\xi\rangle)
\,\widehat{u_{1}}(\xi)\,\overline{\widehat{u_{2}}(\xi)}\,d\xi.
$$
Оно стандартным образом задает норму.
\end{definition}


Для каждого функционального параметра $\varphi\in\mathrm{RO}$ пространство
$H^{\varphi}(\mathbb{R}^{n})$ является частным (изотропным) гильбертовым случаем
пространств Хермандера: $H^{\varphi}(\mathbb{R}^{n})=B_{2,\mu}(\mathbb{R}^{n})$, где
$\mu(\xi):=\varphi(\langle\xi\rangle)$, $\xi\in\mathbb{R}^{n}$. Здесь, как мы сейчас
покажем, функция $\mu$ является весовой.


\begin{lemma}\label{lem2.15}
Пусть $\varphi\in\mathrm{RO}$. Тогда функция $\mu(\xi):=\varphi(\langle\xi\rangle)$
аргумента $\xi\in\mathbb{R}^{n}$ является весовой в смысле
определения~$\ref{def2.8}$.
\end{lemma}


\textbf{Доказательство.} Пусть $\xi,\,\eta\in\mathbb{R}^{n}$. Возведением в квадрат
проверяется неравенство
$|\langle\xi\rangle-\nobreak\langle\eta\rangle|\leq|\,|\xi|-|\eta|\,|$. Отсюда в
случае $\langle\xi\rangle\geq\langle\eta\rangle$ получаем, что
$$
\frac{\langle\xi\rangle}{\langle\eta\rangle}=
1+\frac{\langle\xi\rangle-\langle\eta\rangle}{\langle\eta\rangle}\leq
1+|\xi|-|\eta|\leq1+|\xi-\eta|.
$$
Поэтому в силу предложения~\ref{prop2.11}
$$
\frac{\varphi(\langle\xi\rangle)}{\varphi(\langle\eta\rangle)}\leq
c\,\left(\frac{\langle\xi\rangle}{\langle\eta\rangle}\right)^{s_{1}}\leq
c\,(1+|\xi-\eta|)^{\max\{0,s_{1}\}}.
$$
Если же $\langle\eta\rangle\geq\langle\xi\rangle$, то
$$
\frac{\varphi(\langle\xi\rangle)}{\varphi(\langle\eta\rangle)}\leq
c\,\left(\frac{\langle\xi\rangle}{\langle\eta\rangle}\right)^{s_{0}}=
c\,\left(\frac{\langle\eta\rangle}{\langle\xi\rangle}\right)^{-s_{0}}\leq
c\,(1+|\xi-\eta|)^{\max\{0,-s_{0}\}}.
$$
Таким образом, для любых $\xi,\,\eta\in\mathbb{R}^{n}$ выполняется оценка
$$
\frac{\mu(\xi)}{\mu(\eta)}=
\frac{\varphi(\langle\xi\rangle)}{\varphi(\langle\eta\rangle)}\leq
c\,(1+|\xi-\eta|)^{l},
$$
где числа $c\geq1$ и $l:=\max\{0,-s_{0},s_{1}\}$ не зависят от $\xi$ и $\eta$. Это
означает, что функция $\mu$ весовая в смысле определения~$\ref{def2.8}$. Лемма
\ref{lem2.15} доказана.

\medskip

Заметим, что если $\varphi_{0}\in\mathcal{M}$ и $s\in\mathbb{R}$, то функция
$\varphi_{s}(t):=t^{s}\varphi_{0}(t)$ принадлежит классу $\mathrm{RO}$.
Следовательно, класс пространств
$\{H^{\varphi}(\mathbb{R}^{n}):\varphi\in\mathrm{RO}\}$ содержит в себе уточненную
шкалу пространств Хермандера.

Укажем свойства пространства $H^{\varphi}(\mathbb{R}^{n})$, вытекающие из
соответствующих свойств пространств Хермандера [\ref{Hermander65}, с. 54 -- 62].


\begin{proposition}\label{prop2.12}
Пусть $\varphi,\varphi_{1}\in\mathrm{RO}$. Тогда:
\begin{itemize}
\item [$\mathrm{(i)}$] Пространство $H^{\varphi}(\mathbb{R}^{n})$ сепарабельное
гильбертово.
\item [$\mathrm{(ii)}$] Множество $C^{\infty}_{0}(\mathbb{R}^{n})$ плотно в
$H^{\varphi}(\mathbb{R}^{n})$.
\item [$\mathrm{(iii)}$] Функция
$\varphi(t)/\varphi_{1}(t)$ ограничена в окрестности $\infty$ тогда и только тогда,
когда $H^{\varphi_{1}}(\mathbb{R}^{n})\hookrightarrow H^{\varphi}(\mathbb{R}^{n})$.
Это вложение непрерывно и плотно.
\item [$\mathrm{(iv)}$] Для любых вещественных чисел $s_{0}<\sigma_{0}(\varphi)$ и
$s_{1}>\sigma_{1}(\varphi)$ справедливы непрерывные и плотные вложения
$H^{s_{1}}(\mathbb{R}^{n})\hookrightarrow H^{\varphi}(\mathbb{R}^{n})\hookrightarrow
H^{s_{0}}(\mathbb{R}^{n})$.
\item [$\mathrm{(v)}$] Пространства $H^{\varphi}(\mathbb{R}^{n})$ и
$H^{1/\varphi}(\mathbb{R}^{n})$ взаимно двойственные относительно расширения по
непрерывности скалярного произведения в пространстве $L_{2}(\mathbb{R}^{n})$.
\item [$\mathrm{(vi)}$] Для каждого фиксированного целого числа $k\geq0$ условие
\begin{equation}\label{2.162}
\int\limits_{1}^{\,\infty}\;t^{2k+n-1}\,\varphi^{-2}(t)\,dt<\infty
\end{equation}
равносильно вложению $H^{\varphi}(\mathbb{R}^{n})\hookrightarrow
C^{k}_{\mathrm{b}}(\mathbb{R}^{n})$. Это вложение непрерывно.
\end{itemize}
\end{proposition}


Сделаем ряд пояснений к предложению \ref{prop2.12}. Утверждение (iv) следует из
(iii) и неравенства \eqref{2.160}, в котором полагаем $t:=1$. Поскольку
$\varphi\in\mathrm{RO}\Leftrightarrow1/\varphi\in\mathrm{RO}$, то пространство
$H^{1/\varphi}(\mathbb{R}^{n})$ из утверждения (v) определено корректно. Утверждение
(vi) следует из теоремы вложения Хермандера (предложение \ref{prop2.5}). При этом
следует учесть, что $\eqref{2.33}\Leftrightarrow\eqref{2.162}$ для функции
$\mu(\xi):=\varphi(\langle\xi\rangle)$, $\xi\in\mathbb{R}^{n}$, если перейти к
сферическим координатам.

Изучим интерполяционные свойства пространств Хермандера
$H^{\varphi}(\mathbb{R}^{n})$, где функция $\varphi\in\mathrm{RO}$.


\begin{theorem}\label{th2.37}
Пусть заданы функции $\varphi_{0},\varphi_{1}\in\mathrm{RO}$ и $\psi\in\mathcal{B}$.
Предположим, что функция $\varphi_{0}/\varphi_{1}$ ограничена в окрестности
$\infty$, а $\psi$ --- интерполяционный параметр. Положим
$$
\varphi(t):=\varphi_{0}(t)\,\psi(\varphi_{1}(t)/\varphi_{0}(t))\quad\mbox{при} \quad
t\geq1.
$$
Тогда $\varphi\in\mathrm{RO}$ и
\begin{equation}\label{2.163}
[H^{\varphi_{0}}(\mathbb{R}^{n}),H^{\varphi_{1}}(\mathbb{R}^{n})]_{\psi}=
H^{\varphi}(\mathbb{R}^{n})\quad\mbox{с равенством норм}.
\end{equation}
\end{theorem}


\textbf{Доказательство.} Сначала покажем, что $\varphi\in\mathrm{RO}$. В силу
определения, функция $\varphi$ измерима по Борелю на полуоси $[1,\infty)$. Проверим,
что она удовлетворяет условию \eqref{2.158}. Поскольку
$\varphi_{0},\varphi_{1}\in\mathrm{RO}$, существуют числа $a>1$ и $c>1$ такие, что
\begin{equation}\label{2.164}
c^{-1}\leq\frac{\varphi_{j}(\lambda t)}{\varphi_{j}(t)}\leq c \;\;\mbox{для
всех}\;\;t\geq1,\;\;\lambda\in[1,a],\;\;j\in\{0,\,1\}.
\end{equation}
Из ограниченности функции $\varphi_{0}/\varphi_{1}$ в окрестности $\infty$ вытекает
ввиду предложения~\ref{prop2.9}, что существует число $\varepsilon>0$, при котором
\begin{equation}\label{2.165}
\frac{\varphi_{1}(t)}{\varphi_{0}(t)}>\varepsilon\quad\mbox{для любого}\quad t\geq1.
\end{equation}
Далее, так как $\psi$ --- интерполяционный параметр, то по теореме~\ref{th2.9}
функция $\psi$ псевдовогнута в окрестности $\infty$. Тогда в силу лемм \ref{lem2.1}
и \ref{lem2.2} функция $\psi$ слабо эквивалентна вогнутой функции на полуоси
$(\varepsilon,\infty)$, а это равносильно следующему условию: существует число
$c_{0}>1$ такое, что
\begin{equation}\label{2.166}
\frac{\psi(\tau)}{\psi(t)}\leq
c_{0}\max\left\{1,\frac{\tau}{t}\right\}\quad\mbox{для любых}\quad
\tau>\varepsilon,\;\;t>\varepsilon.
\end{equation}
Отсюда вытекает неравенство
\begin{equation}\label{2.167}
\frac{\psi(t)}{\psi(\tau)}\geq c_{0}^{-1}\min\left\{1,\frac{t}{\tau}\right\}
\quad\mbox{для любых}\quad\tau>\varepsilon,\;\;t>\varepsilon.
\end{equation}

Теперь из формул \eqref{2.164}, \eqref{2.165} и \eqref{2.166} следует, что для любых
$t\geq1$ и $\lambda\in[1,a]$
\begin{gather*}
\frac{\varphi(\lambda t)}{\varphi(t)}= \frac{\varphi_{0}(\lambda
t)}{\varphi_{0}(t)}\cdot\frac{\psi(\varphi_{1}(\lambda t)/\varphi_{0}(\lambda
t))}{\psi(\varphi_{1}(t)/\varphi_{0}(t))}\leq \\
c\cdot c_{0}\max\Bigl\{1,\frac{\varphi_{1}(\lambda t)/\varphi_{0}(\lambda
t)}{\varphi_{1}(t)/\varphi_{0}(t)}\Bigr\}\leq c^{3}c_{0}.
\end{gather*}
Аналогично из формул \eqref{2.164}, \eqref{2.165} и \eqref{2.167} вытекает, что
$$
\frac{\varphi(\lambda t)}{\varphi(t)}\geq c^{-1}
c_{0}^{-1}\min\Bigl\{1,\frac{\varphi_{1}(\lambda t)/\varphi_{0}(\lambda
t)}{\varphi_{1}(t)/\varphi_{0}(t)}\Bigr\}\geq c^{-3}c_{0}^{-1}.
$$
Таким образом, функция $\varphi$ удовлетворяет условию \eqref{2.158} и поэтому
$\varphi\in\mathrm{RO}$.

Докажем теперь равенство \eqref{2.163}. В силу предложения \ref{prop2.12} (iii) пара
$[H^{\varphi_{0}}(\mathbb{R}^{n}),H^{\varphi_{1}}(\mathbb{R}^{n})]$ является
допустимой. Псевдодифференциальный оператор с символом
$\varphi_{1}(\langle\xi\rangle)/\varphi_{0}(\langle\xi\rangle)$, где
$\xi\in\mathbb{R}^{n}$, является порождающим оператором $J$ для этой пары. С помощью
преобразования Фурье
$$
\mathcal{F}:H^{\varphi_{0}}(\mathbb{R}^{n})\leftrightarrow
L_{2}(\mathbb{R}^{n},\varphi_{0}^{2}(\langle\xi\rangle)\,d\xi)
$$
порождающий оператор $J$ приводится к виду умножения на функцию
$\varphi_{1}(\langle\xi\rangle)/\varphi_{0}(\langle\xi\rangle)$. Следовательно,
оператор $\psi(J)$ приводится к виду умножения на функцию
$$
\psi\left(\frac{\varphi_{1}(\langle\xi\rangle)}{\varphi_{0}(\langle\xi\rangle)}\right)=
\frac{\varphi(\langle\xi\rangle)}{\varphi_{0}(\langle\xi\rangle)}.
$$
Поэтому для любой функции $u\in C^{\infty}_{0}(\mathbb{R}^{n})$ можем записать:
\begin{align*}
\|u\|_{[H^{\varphi_{0}}(\mathbb{R}^{n}),H^{\varphi_{1}}(\mathbb{R}^{n})]_{\psi}}
^{2}&=\|\psi(J)u\|_{H^{\varphi_{0}}(\mathbb{R}^{n})}^{2}= \\
=\int\limits_{\mathbb{R}^{n}}\varphi_{0}^{2}(\langle\xi\rangle)\:
|(\widehat{\psi(J)u})(\xi)|^{2}\,d\xi&=
\int\limits_{\mathbb{R}^{n}}\varphi^{2}(\langle\xi\rangle)\:|\widehat{u}(\xi)|^{2}\,
d\xi=\|u\|_{H^{\varphi}(\mathbb{R}^{n})}^{2}.
\end{align*}
Отсюда следует равенство пространств \eqref{2.163}, поскольку множество
$C^{\infty}_{0}(\mathbb{R}^{n})$ плотно в каждом из них. Это вытекает из
предложения~\ref{prop2.12} (ii) и теоремы~\ref{th2.1}, согласно которой пространство
$H^{\varphi_{1}}(\mathbb{R}^{n})$ непрерывно и плотно вложено в
$[H^{\varphi_{0}}(\mathbb{R}^{n}),H^{\varphi_{1}}(\mathbb{R}^{n})]_{\psi}$.

Теорема \ref{th2.37} доказана.


\begin{theorem}\label{th2.38}
Пусть заданы функция $\varphi\in\mathrm{RO}$ и вещественные числа $s_{0}$, $s_{1}$
такие, что $s_{0}<\sigma_{0}(\varphi)$ и $s_{1}>\sigma_{1}(\varphi)$. Положим
$$
\psi(t):=
\begin{cases}
\;t^{-s_{0}/(s_{1}-s_{0})}\,\varphi(t^{1/(s_{1}-s_{0})}) & \text{при}\;\;t\geq1, \\
\;\varphi(1) & \text{при}\;\;0<t<1.
\end{cases}
$$
Тогда функция $\psi\in\mathcal{B}$ является интерполяционным параметром~и
\begin{equation}\label{2.168}
[H^{s_{0}}(\mathbb{R}^{n}),H^{s_{1}}(\mathbb{R}^{n})]_{\psi}=
H^{\varphi}(\mathbb{R}^{n})\quad\mbox{с равенством норм}.
\end{equation}
\end{theorem}


\textbf{Доказательство.} Так как $\varphi(t)=t^{s_{0}}\,\psi(t^{s_{1}}/t^{s_{0}})$
при $t\geq1$, то теорема~\ref{th2.38} следует из теорем \ref{th2.9} и \ref{th2.37},
если мы докажем, что функция $\psi$ принадлежит множеству $\mathcal{B}$ и является
псевдовогнутой в окрестности $\infty$. В силу предложения~\ref{prop2.11}, функция
$\psi$ удовлетворяет условию \eqref{2.166} для $\varepsilon=1$. В самом деле, если
$t\geq1$ и $\tau\geq1$, то
\begin{gather*}
\frac{\psi(\tau)}{\psi(t)}=\left(\frac{\tau}{t}\right)^{-s_{0}/(s_{1}-s_{0})}\,
\frac{\varphi(\tau^{1/(s_{1}-s_{0})})}{\varphi(t^{1/(s_{1}-s_{0})})}\leq\\
\left(\frac{\tau}{t}\right)^{-s_{0}/(s_{1}-s_{0})}\:c\,
\max\left\{\left(\frac{\tau}{t}\right)^{s_{1}/(s_{1}-s_{0})},\,
\left(\frac{\tau}{t}\right)^{s_{0}/(s_{1}-s_{0})}\right\}= \\
c\,\max\left\{\frac{\tau}{t},1\right\}.
\end{gather*}
Отсюда, в частности, вытекает, что функция $\psi$ отделена от нуля на полуоси
$[1,\infty)$ и поэтому $\psi\in\mathcal{B}$ ввиду предложения~\ref{prop2.9}.
Остается сослаться на лемму~\ref{lem2.2}, согласно которой условие \eqref{2.166}
равносильно тому, что функция $\psi$ псевдовогнута на полуоси
$(\varepsilon,\infty)$. Теорема~\ref{th2.38} доказана.

\medskip

Теперь мы в состоянии доказать следующее фундаментальное свойство класса пространств
Хермандера
\begin{equation}\label{2.169}
\{H^{\varphi}(\mathbb{R}^{n}):\,\varphi\in\mathrm{RO}\}.
\end{equation}


\begin{theorem}\label{th2.39}
Класс пространств \eqref{2.169} совпадает (с точностью до эквивалентности норм) с
множеством всех гильбертовых интерполяционных пространств для пар гильбертовых
пространств Соболева
\begin{equation}\label{2.170}
[H^{s_{0}}(\mathbb{R}^{n}),H^{s_{1}}(\mathbb{R}^{n})],
\end{equation}
где $s_{0},s_{1}\in\mathbb{R}$ и $s_{0}<s_{1}$.
\end{theorem}


\textbf{Доказательство.} В силу теоремы \ref{th2.38} каждое пространство
$H^{\varphi}(\mathbb{R}^{n})$, где $\varphi\in\mathrm{RO}$, является
интерполяционным для пары пространств \eqref{2.170} при условии, что
$s_{0}<\sigma_{0}(\varphi)$ и $s_{1}>\sigma_{1}(\varphi)$. Обратно, если гильбертово
пространство $H$ является интерполяционным для пары пространств \eqref{2.170}, где
$s_{0},s_{1}\in\mathbb{R}$ и $s_{0}<s_{1}$, то согласно предложению \ref{prop2.1}
$$
H=[H^{s_{0}}(\mathbb{R}^{n}),H^{s_{1}}(\mathbb{R}^{n})]_{\psi}\quad\mbox{с
эквивалентностью норм}
$$
для некоторой псевдовогнутой в окрестности $\infty$ функции $\psi\in\mathcal{B}$. По
теореме \ref{th2.9} эта функция является интерполяционным параметром. Отсюда на
основании теоремы \ref{th2.37} получаем, что $H=H^{\varphi}(\mathbb{R}^{n})$, где
функция $\varphi\in\mathrm{RO}$ определена по формуле
$\varphi(t):=t^{s_{0}}\psi(t^{s_{1}-s_{0}})$ при $t\geq1$. Теорема~\ref{th2.39}
доказана.


\begin{remark}\label{rem2.13}
В силу теорем \ref{th2.37}, \ref{th2.9} и предложения \ref{prop2.1} класс
пространств \eqref{2.169} замкнут относительно интерполяции, результатом которой
есть гильбертово пространство.
\end{remark}


Рассмотрим теперь интерполяционные пространства Хермандера на бесконечно гладком
замкнутом многообразии $\Gamma$ размерности $n\geq1$. Следующая теорема дает
эквивалентные определения пространства $H^{\varphi}(\Gamma)$,
$\varphi\in\mathrm{RO}$ (сравнить с п.~\ref{sec2.5.1}). Мы используем обозначения из
п.~\ref{sec2.5.1}.


\begin{theorem}\label{th2.40}
Пусть $\varphi\in\mathrm{RO}$. Нижеследующие определения дают одно и тоже
гильбертово $H^{\varphi}(\Gamma)$ с точностью до эквивалентности норм.
\begin{itemize}
\item [$\mathrm{(i)}$] Линейное пространство $H^{\varphi}(\Gamma)$
состоит, по определению, из всех распределений $f\in\mathcal{D}'(\Gamma)$ таких, что
$(\chi_{j}f)\circ\alpha_{j}\in H^{\varphi}(\mathbb{R}^{n})$ для каждого
$j=1,\ldots,r$. В~пространстве $H^{\varphi}(\Gamma)$ задано скалярное произведение
по формуле
$$
(f_{1},f_{2})_{H^{\varphi}(\Gamma)}:=\sum_{j=1}^{r}\,((\chi_{j}f_{1})\circ\alpha_{j},
(\chi_{j}f_{2})\circ\alpha_{j})_{H^{\varphi}(\mathbb{R}^{n})}
$$
и соответствующая норма.
\item [$\mathrm{(ii)}$] Пусть целые числа $k_{0}$ и $k_{1}$ такие, что
$k_{0}<\sigma_{0}(\varphi)$ и $k_{1}>\sigma_{1}(\varphi)$. По определению,
$$
H^{\varphi}(\Gamma):=[H^{k_{0}}(\Gamma),H^{k_{1}}(\Gamma)]_{\psi},
$$
где интерполяционный параметр $\psi$ задан формулой
$$
\psi(t)=
\begin{cases}
\;t^{-k_{0}/(k_{1}-k_{0})}\,
\varphi(t^{1/(k_{1}-k_{0})}) & \text{при\;\;\;$t\geq1$}, \\
\;\varphi(1) & \text{при\;\;\;$0<t<1$}.
\end{cases}
$$
\item [$\mathrm{(iii)}$] По определению, пространство
$H^{\varphi}(\Gamma)$ --- пополнение линейного многообразия $C^{\infty}(\Gamma)$ по
гильбертовой норме
$$
f\,\mapsto\,\|\varphi((1-\Delta_{\Gamma})^{1/2})\,f\|_{L_{2}(\Gamma)},\quad f\in
C^{\infty}(\Gamma).
$$
\end{itemize}
\end{theorem}


Эта теорема является частным случаем следующих двух теорем. За исходное определение
пространства $H^{\varphi}(\Gamma)$ мы берем п. (i) теоремы~\ref{th2.40}.


\begin{theorem}\label{th2.41}
Интерполяционные теоремы $\ref{th2.37}$ и $\ref{th2.38}$ сохраняют силу, если в их
формулировках заменить $\mathbb{R}^{n}$ на $\Gamma$ и равенство норм на их
эквивалентность.
\end{theorem}


Теорема~\ref{th2.41} доказывается аналогично теореме~\ref{th2.21}.

Пусть $A$ --- эллиптический ПДО порядка $m>0$. Предположим, что $A$ является
(неограниченным) самосопряженным положительно определенным оператором в гильбертовом
пространстве $L_{2}(\Gamma)$. Доопределим функцию $\varphi\in\mathrm{RO}$ равенством
$\varphi(t):=\varphi(1)$ при $0<t<1$.


\begin{theorem}\label{th2.42}
Для произвольного $\varphi\in\mathrm{RO}$ норма в пространстве $H^{\varphi}(\Gamma)$
эквивалентна норме
\begin{equation}\label{2.171}
f\,\mapsto\,\|\varphi(A^{1/m})\,f\|_{L_{2}(\Gamma)}.
\end{equation}
на плотном множестве $C^{\infty}(\Gamma)$. Тем самым пространство
$H^{\varphi}(\Gamma)$ совпадает с пополнением множества $C^{\infty}(\Gamma)$ по
норме \eqref{2.171}.
\end{theorem}


Теорема~\ref{th2.42} доказывается аналогично теореме~\ref{th2.24}.

В заключение этого пункта приведем аналог теоремы~\ref{th2.39} для для класса
пространств
\begin{equation}\label{2.172}
\{H^{\varphi}(\Gamma):\,\varphi\in\mathrm{RO}\}.
\end{equation}


\begin{theorem}\label{th2.43}
Класс пространств \eqref{2.172} совпадает (с точностью до эквивалентности норм) с
множеством всех гильбертовых интерполяционных пространств для пар гильбертовых
пространств Соболева $[H^{s_{0}}(\Gamma),H^{s_{1}}(\Gamma)]$, где
$s_{0},s_{1}\in\mathbb{R}$ и $s_{0}<s_{1}$.
\end{theorem}


Эта теорема доказывается аналогично теореме~\ref{th2.39}.

\subsection{Приложения к эллиптическим операторам}\label{sec2.8.3}


Изучим некоторые приложения классов пространств \eqref{2.169} и \eqref{2.172} к
эллиптическим ПДО.

Сначала рассмотрим ПДО в $\mathbb{R}^{n}$. Изложенные ниже результаты полезно
сравнить с п.~\ref{sec2.4b}. Предварительно изучим действие ПДО в классе пространств
\eqref{2.169}. Положим $\varrho(t):=t$ при $t\geq1$.

\begin{lemma}\label{lem2.16}
Пусть $A$~--- ПДО класса $\Psi^{m}(\mathbb{R}^{n})$, где $m\in\mathbb{R}$. Тогда
сужение отображения $u\mapsto Au$, $u\in\mathcal{S}'(\mathbb{R}^{n})$, на
пространство $H^{\varphi\varrho^{m}}(\mathbb{R}^{n})$ является линейным ограниченным
оператором
\begin{equation}\label{2.173}
A:H^{\varphi\varrho^{m}}(\mathbb{R}^{n})\rightarrow H^{\varphi}(\mathbb{R}^{n})
\quad\mbox{для любого}\quad\varphi\in\mathrm{RO}.
\end{equation}
\end{lemma}

Эта лемма доказывается с помощью интерполяционной теоремы~\ref{th2.38} аналогично
лемме~\ref{lem2.8b}.

Предположим, что ПДО $A\in\Psi^{m}_{\mathrm{ph}}(\mathbb{R}^{n})$ равномерно
эллиптический в $\mathbb{R}^{n}$. Тогда отображение \eqref{2.173} имеет следующие
свойства.


\begin{theorem}\label{th2.44}
Пусть заданы функция $\varphi\in\mathrm{RO}$ и число $\sigma>0$. Существует число
$c=c(\varphi,\sigma)>0$ такое, что для произвольного распределения $u\in
H^{\varphi\varrho^{m}}(\mathbb{R}^{n})$ справедлива априорная оценка
$$
\|u\|_{H^{\varphi\varrho^{m}}(\mathbb{R}^{n})}\leq
c\,\bigl(\,\|Au\|_{H^{\varphi}(\mathbb{R}^{n})}+
\|u\|_{H^{\varphi\varrho^{m-\sigma}}(\mathbb{R}^{n})}\,\bigr).
$$
\end{theorem}


Эта теорема доказывается аналогично теореме~\ref{th2.16b}.

Пусть $V$ --- произвольное непустое открытое подмножество пространства
$\mathbb{R}^{n}$. Для $\varphi\in\mathrm{RO}$ положим
\begin{gather*}
H^{\varphi}_{\mathrm{int}}(V):=\bigl\{w\in H^{-\infty}(\mathbb{R}^{n}): \chi\,w\in
H^{\varphi}(\mathbb{R}^{n})\\ \mbox{для всех}\;\chi\in
C^{\infty}_{\mathrm{b}}(\mathbb{R}^{n}),\,\mathrm{supp}\,\chi\subset V,\;
\mathrm{dist}(\mathrm{supp}\,\chi,\partial V)>0\bigr\}.
\end{gather*}


\begin{theorem}\label{th2.45}
Предположим, что $u\in H^{-\infty}(\mathbb{R}^{n})$ является решением уравнения
$Au=f$ на множестве $V$, где $f\in H^{\varphi}_{\mathrm{int}}(V)$ для некоторого
параметра $\varphi\in\mathrm{RO}$. Тогда $u\in
H^{\varphi\varrho^{m}}_{\mathrm{int}}(V)$.
\end{theorem}


Эта теорема доказывается аналогично теореме~\ref{th2.18b}.

Из теоремы \ref{th2.45} и предложения \ref{prop2.12} (vi) вытекает следующее (более
тонкое, чем теорема \ref{th2.19b}) достаточное условие существования непрерывных
производных у решения уравнения $Au=f$.


\begin{theorem}\label{th2.45b}
Пусть заданы целое число $r\geq0$ и функция $\varphi\in\mathcal{M}$, удовлетворяющая
условию
$$
\int\limits_{1}^{\,\infty}\;t^{2(r-m)+n-1}\varphi^{-2}(t)\,dt<\infty.
$$
Предположим, что распределение $u\in H^{-\infty}(\mathbb{R}^{n})$ является решением
уравнения $Au=\nobreak f$ на открытом множестве $V\subseteq\mathbb{R}^{n}$, где
$f\in H^{\varphi}_{\mathrm{int}}(V)$. Тогда верно заключение теоремы
$\ref{th2.19b}$.
\end{theorem}


Далее кратко остановимся на приложении пространств $H^{\varphi}(\Gamma)$,
$\varphi\in\mathrm{RO}$, к эллиптическим ПДО на бесконечно гладком замкнутом
многообразии $\Gamma$.

Пусть ПДО $A\in\Psi^{m}_{\mathrm{ph}}(\Gamma)$, где $m\in\mathbb{R}$, является
эллиптическим на $\Gamma$. Для $A$ конечномерные пространства $\mathcal{N}$ и
$\mathcal{N}^{+}$ определены по формулам \eqref{PsMan2} и \eqref{PsMan3}.


\begin{theorem}\label{th2.46}
Для произвольного параметра $\varphi\in\mathrm{RO}$ сужение отображения $u\mapsto
Au$, $u\in\mathcal{D}'(\Gamma)$, на пространство $H^{\varphi\varrho^{m}}(\Gamma)$
является ограниченным оператором
\begin{equation}\label{2.175}
A: H^{\varphi\varrho^{m}}(\Gamma)\rightarrow H^{\varphi}(\Gamma).
\end{equation}
Этот оператор нетеров. Его ядро совпадает с $\mathcal{N}$, а область значений равна
$$
\bigl\{f\in H^{\varphi}(\Gamma):\,(f,w)_{\Gamma}=0\;\;\mbox{для
всех}\;\;w\in\mathcal{N}^{+}\bigr\}.
$$
Индекс оператора \eqref{2.175} равен $\dim\mathcal{N}-\dim\mathcal{N}^{+}$ и не
зависит от~$\varphi$.
\end{theorem}


Эта теорема доказывается аналогично теореме~\ref{th2.28b}. Из теоремы~\ref{th2.46}
следует аналоги теорем \ref{th2.29b} -- \ref{th2.32b} для шкалы пространств
$H^{\varphi}(\Gamma)$, $\varphi\in\mathrm{RO}$. Это устанавливается с помощью тех же
рассуждений, что и в п. \ref{sec2.6.2b}, \ref{sec2.6.3b}. Формулировки этих аналогов
мы приводить не будем.





\markright{\emph \ref{notes2}. Примечания и комментарии}

\section{Примечания и комментарии}\label{notes2}

\markright{\emph \ref{notes2}. Примечания и комментарии}

\small

\textbf{К п. 2.1.} Установленные нами эквивалентные определения пространств
Хермандера на замкнутом (компактном) гладком многообразии аналогичны тем, которые
используются для пространств Соболева; см., например, монографии Ж.-Л.~Лионса и
Э.~Мадженеса [\ref{LionsMagenes71}] (п.~7.3), М.~Тейлора [\ref{Taylor85}] (гл.~I,
\S~5).

В статье Г.Шлензак [\ref{Shlenzak74}, с.~56] используется введенный локально
некоторый класс пространств Хермандера на границе области. Независимость этих
пространств и топологии в них от выбора локальных карт в [\ref{Shlenzak74}] не
доказывается.

Вещественные пространства Хермандера на окружности, где функции задаются
тригонометрическими рядами Фурье, использовались в работах Й.~Пёшеля
[\ref{Poschel08}], П.~Б.~Джакова и Б.~С.~Митягина [\ref{DjakovMityagin06},
\ref{DjakovMityagin09}], В.~А.~Михайлеца и В.~Н.~Молибоги
[\ref{MikhailetsMolyboga09a}, \ref{MikhailetsMolyboga10}]. Эти пространства тесно
связаны с пространствами периодических вещественных функций, введенных
А.~И.~Степанцом [\ref{Stepanets87}, \ref{Stepanets02}]. Они широко используются в
теории приближений.

Все результаты п. 2.1 установлены нами в статье [\ref{08MFAT1}] (п. 3.3, 3.4), за
исключением теоремы 2.7, доказательство которой не было опубликовано. Частично эти
результаты анонсированы ранее в [\ref{06Dop10}]. Ряд результатов~--- теоремы 2.2,
2.3 и 2.8 справедливы также и для компактного гладкого многообразия с краем
[\ref{06UMJ3}] (п.~3).

\medskip

\textbf{К п. 2.2.} Исчисление классических ПДО на гладких многообразиях без края
построено Л.~Хермандером [\ref{Hermander67b}]. Систематическое изложение теории
эллиптических ПДО на таких многообразиях имеется, например, в его монографии
[\ref{Hermander87}] (гл.~19) и в обзоре М.~С.~Аграновича [\ref{Agranovich90}]
(\S~2). Приведенные там классические результаты о нетеровости эллиптических ПДО и
регулярности решений в соболевской шкале нашли важные приложения в теории
эллиптических краевых задач для дифференциальных уравнений, в спектральной теории, в
теории функциональных пространств и др., см. также монографии М.~Тейлора
[\ref{Taylor85}], Ф.~Трева [\ref{Treves84a}, \ref{Treves84b}], М.~А.~Шубина
[\ref{Shubin78}].

Отдельный интерес представляют эллиптические ПДО, осуществляющие гомеоморфизмы в
соболевской шкале. Широкий класс таких операторов~--- эллиптические операторы с
параметром выделен и изучен в работах Ш.~Агмона и Л.~Ниренберга [\ref{Agmon62},
\ref{AgmonNirenberg63}], М.~С.~Аграновича и М.~И.~Вишика [\ref{AgranovichVishik64}];
см. также обзор М.~С.~Аграновича [\ref{Agranovich90}] (\S~4). Эти операторы имеют
важные приложения в спектральной теории и теории параболических уравнений.

Все теоремы п. 2.2. установлены в [\ref{07UMJ6}], часть из них анонсирована в
[\ref{06Dop10}] применительно к дифференциальным операторам.

\medskip

\textbf{К п. 2.3.} Для произвольных ортогональных рядов классическая теорема о
сходимости почти всюду доказана независимо Д.~Е.~Меньшовым [\ref{Menschoff23}] и
Г.~Радемахером [\ref{Rademacher22}]. Теоремы о безусловной сходимости установлены
В.~Орличем [\ref{Orlicz27}] и К.~Тандори [\ref{Tandori62}]. Мы используем
эквивалентную формулировку теоремы Орлича, полученную П.~Л.~Ульяновым
[\ref{Ulyanov63}, с. 53] (см. также [\ref{Ulyanov64}, с. 53]). Все эти теоремы
доказаны упомянутыми авторами для случая вещественных рядов, заданных на интервале
оси $\mathbb{R}$, и меры Лебега. Изложение этих результатов имеется, например, в
монографиях Г.~Алексича [\ref{Alexits63}], Б.~С.~Кашина и А.~А.~Саакяна
[\ref{KashinSaakyan84}].

Теоремы 2.15 и 2.16 анонсированы нами в [\ref{09Dop3}] (п.~5). Они обобщают и
существенно уточняют результат К.~Мини [\ref{Meaney82}], который использует
соболевскую шкалу и рассматривает разложения по собственным функциям оператора
Бельтрами-Лапласа. Теорема 2.17 ранее не приводилась.

\medskip

\textbf{К п. 2.4.} Класс RO-меняющихся функций введен В.~Авакумовичем
[\ref{Avakumovic36}] в 1936~г. и достаточно полно изучен: см., например, монографии
Н.~Х.~Бинхема, Ч.~М.~Голди и Дж.~Л.~Тайгельза [\ref{BinghamGoldieTeugels89}, с.~65],
Е.~Сенеты [\ref{Seneta85}, с.~86].

Пространства Хермандера, параметризуемые RO-ме\-ня\-ю\-щи\-ми\-ся функциями изучены
нами в [\ref{08Collection1}, \ref{09Dop3}], где установлены все теоремы п. 2.4. А
именно, в статье [\ref{08Collection1}] исследованы пространства в $\mathbb{R}^{n}$ и
даны их приложения к равномерно эллиптическим ПДО. В~заметке [\ref{09Dop3}] введены
и изучены пространства на замкнутом гладком многообразии, приведены их приложения к
эллиптическим ПДО и исследованию сходимости спектральных разложений.

\normalsize



\chapter[Полуоднородные эллиптические краевые задачи]
{\textbf{Полуоднородные \\ эллиптические \\ краевые задачи}}\label{ch3}

\chaptermark{\emph Гл. \ref{ch3}. Полуоднородные эллиптические краевые
задачи}



\section[Регулярная эллиптическая краевая задача]
{Регулярная эллиптическая \\ краевая задача}

\markright{\emph \ref{sec3.1}. Регулярная эллиптическая краевая
задача}\label{sec3.1}

Здесь мы приведем определение регулярной эллиптической краевой задачи в ограниченной
евклидовой области, а также сформулируем некоторые понятия, относящиеся к этой
задаче.

\subsection{Определение задачи}\label{sec3.1.1}

Далее в главах \ref{ch3} и \ref{ch4} предполагается, что $\Omega$ --- произвольная
ограниченная область в евклидовом пространстве $\mathbb{R}^{n}$ $(n\geq2)$ с
границей $\Gamma$, которая является бесконечно гладким многообразием без края
размерности $n-1$. Область $\Omega$ локально лежит по одну сторону от $\Gamma$.
Положим $\overline{\Omega}:=\Omega\cup\Gamma$ и
$\widehat{\Omega}:=\mathbb{R}^{n}\setminus\Omega$. Обозначим через $\nu(x)$ орт
внутренней нормали к границе $\Gamma$ в точке $x\in\Gamma$, а через $\nu$ бесконечно
гладкое векторное поле этих ортов.

В области $\Omega$ рассматривается краевая задача
\begin{gather}\label{3.1}
Lu\equiv\sum_{|\mu|\leq\,2q}l_{\mu}(x)\,D^{\mu}u=f\;\;\mbox{в}\;\;\Omega, \\
B_{j}\,u\equiv\sum_{|\mu|\leq\,m_{j}}b_{j,\mu}(x)\,D^{\mu}u=g_{j}\;\;\mbox{на}
\;\;\Gamma,\;\;j=1,\ldots,q. \label{3.2}
\end{gather}
Здесь $L=L(x,D)$ --- линейное дифференциальное выражение на $\overline{\Omega}$
четного порядка $2q\geq2$, а $B_{j}=B_{j}(x,D)$, $j=1,\ldots,q$, --- граничные
линейные дифференциальные выражения на $\Gamma$ порядков
$\mathrm{ord}\,B_{j}=m_{j}\leq2q-1$. Все коэффициенты дифференциальных выражений $L$
и $B_{j}$ предполагаются бесконечно гладкими комплекснозначными функциями:
$l_{\mu}\in C^{\infty}(\,\overline{\Omega}\,)$ и $b_{j,\mu}\in C^{\infty}(\Gamma)$.
Положим $B:=(B_{1},\ldots,B_{q})$.

В формулах \eqref{3.1} и \eqref{3.2} и далее используются стандартные обозначения:
$\mu:=(\mu_{1},\ldots,\mu_{n})$~--- мультииндекс, $|\mu|:=\mu_{1}+\ldots+\mu_{n}$,
$D^{\mu}:=D_{1}^{\mu_{1}}\ldots D_{n}^{\mu_{n}}$, $D_{k}:=i\partial/\partial x_{k}$
при $k=1,\ldots,n$, $i$ --- мнимая единица, $x=(x_{1},\ldots,x_{n})$ --- точка
пространства $\mathbb{R}^n$.

Всюду далее предполагается, что краевая задача \eqref{3.1}, \eqref{3.2} регулярная
эллиптическая в области $\Omega$.

Напомним соответствующее определение (см., например, [\ref{LionsMagenes71}, с.~137,
138], [\ref{FunctionalAnalysis72}, с.~167]).

Дифференциальным выражениям $L(x,D)$ при фиксированном $x\in\overline{\Omega}$ и
$B_{j}(x,D)$ при фиксированном $x\in\Gamma$  поставим в соответствие однородные
характеристические полиномы
$$
L^{(0)}(x,\xi):=\sum_{|\mu|=2q}l_{\mu}(x)\,\xi^{\mu}\quad\mbox{и}\quad
B_{j}^{(0)}(x,\xi):=\sum_{|\mu|=\,m_{j}}b_{j,\mu}(x)\,\xi^{\mu}.
$$
Они называются также главными символами выражений $L$ и $B_{j}$. Здесь переменное
$\xi=(\xi_{1},\ldots,\xi_{n})\in\mathbb{C}^n$, и
$\xi^{\mu}:=\xi_{1}^{\mu_{1}}\ldots\xi_{n}^{\mu_{n}}$.


\begin{definition}\label{def3.1}
Краевая задача \eqref{3.1}, \eqref{3.2} называется \emph{регулярной эллиптической} в
области $\Omega$, если выполняются следующие условия:
\begin{itemize}
\item [(i)] Дифференциальное выражение $L$ является \emph{правильно эллиптическим}
на $\overline{\Omega}$, т.~е. для произвольной точки $x\in\overline{\Omega}$ и любых
линейно независимых векторов $\xi',\xi''\in\mathbb{R}^n$ многочлен
$L^{(0)}(x,\xi'+\tau\xi'')$ переменного $\tau$ имеет ровно $q$ корней
$\tau^{+}_{j}(x;\xi',\xi'')$, $j=1,\ldots,q$, с положительной мнимой частью и
столько же корней с отрицательной мнимой частью (с учетом их кратности).
\item [(ii)] Система граничных выражений$\{B_{1},\ldots,B_{q}\}$ удовлетворяет
условию \emph{дополнительности} по отношению к $L$ на $\Gamma$, т.~е. для
произвольной точки $x\in\Gamma$ и любого вектора $\xi\neq0$ касательного к $\Gamma$
в точке $x$, многочлены $B_{j}^{(0)}(x,\xi+\tau\nu(x))$, $j=1,\ldots,q$, переменного
$\tau$ линейно независимы по модулю многочлена
$\prod_{j=1}^{q}(\tau-\tau^{+}_{j}(x;\xi,\nu(x)))$.
\item [(iii)] Система выражений $\{B_{1},\ldots,B_{q}\}$ является \emph{нормальной}, т.~е.
их порядки $m_{j}$, $j=1,\ldots,q$, попарно различны, и
$B_{j}^{(0)}(x,\nu(x))\neq\nobreak0$ для любого $x\in\Gamma$.
\end{itemize}
\end{definition}


Всюду далее (за исключением п.~4.1.3) предполагается, что краевая задача
\eqref{3.1}, \eqref{3.2} является регулярной эллиптической в области $\Omega$.


\begin{remark}\label{rem3.1}
Из условия (i) определения \ref{def3.1} следует, что
\begin{equation}\label{3.3}
L^{(0)}(x,\xi)\neq0\quad\mbox{для любых}\quad
x\in\overline{\Omega},\;\;\xi\in\mathbb{R}^{n}\setminus\{0\},
\end{equation}
т. е. дифференциальное выражение $L$ является \emph{эллиптическим} на
$\overline{\Omega}$. Если $n\geq3$, то условия (i) и \eqref{3.3} эквивалентные
[\ref{FunctionalAnalysis72}, с.~166]. Если все коэффициенты выражения $L$
вещественные, то эта эквивалентность имеет место и при $n=2$.
\end{remark}


\begin{remark}\label{rem3.2}
Условие дополнительности (определение \ref{def3.1} (ii)) впервые сформулировано
Я.~Б.~Лопатинским [\ref{Lopatinsky53}, \ref{Lopatinsky84}] и, в частных случаях,
З.~Я.~Шапиро [\ref{Shapiro53}]. Известны иные эквивалентные формы этого условия
[\ref{Agranovich97}, с.~7--9]. Условие нормальности системы граничных выражений
(определение \ref{def3.1} (iii)) введено независимо Н.~Ароншайном, А.~Мильграмом
[\ref{AronszajnMilgram52}] и М.~Шехтером [\ref{Schechter60a}].
\end{remark}


\begin{example}\label{ex3.1}
Пусть $k\in\mathbb{Z}$ и $0\leq k\leq q$. Система граничных выражений
$B_{j}u:=\partial^{k+j-1}u/\partial\nu^{k+j-1}$, $j=1,\ldots,q$, является нормальной
и удовлетворяет условию дополнительности на $\Gamma$ по отношению к любому
дифференциальному выражению $L$, правильно эллиптическому на $\overline{\Omega}$
(см., например, [\ref{Triebel80}, с.~454]). В случае $k=0$ мы имеем краевую задачу
Дирихле для уравнения $Lu=f$.
\end{example}

\subsection{Формально сопряженная задача}\label{sec3.1.2}

Наряду с задачей \eqref{3.1}, \eqref{3.2} рассмотрим краевую задачу
\begin{gather}\label{3.4}
L^{+}v\equiv\sum_{|\mu|\leq\,2q}D^{\mu}(\overline{l_{\mu}(x)}\,v)=\omega
\;\;\mbox{в}\;\;\Omega,\\
B^{+}_{j}v=h_{j}\;\;\mbox{на} \;\;\Gamma,\;\;j=1,\ldots,q. \label{3.5}
\end{gather}
Она формально сопряжена к задаче \eqref{3.1}, \eqref{3.2} относительно формулы
Грина:
\begin{equation}\label{3.6}
(Lu,v)_{\Omega}+\sum_{j=1}^{q}\;(B_{j}u,\,C_{j}^{+}v)_{\Gamma}
=(u,L^{+}v)_{\Omega}+\sum_{j=1}^{q}\;(C_{j}u,\,B_{j}^{+}v)_{\Gamma},
\end{equation}
справедливой для любых функций $u,v\in C^{\infty}(\,\overline{\Omega}\,)$
[\ref{FunctionalAnalysis72}, с.~168]. Здесь $\{B^{+}_{j}\}$, $\{C_{j}\}$,
$\{C^{+}_{j}\}$~--- некоторые нормальные системы граничных линейных дифференциальных
выражений с коэффициентами класса $C^{\infty}(\Gamma)$. Их порядки удовлетворяют
условию
\begin{equation}\label{3.7}
\mathrm{ord}\,B_{j}+\mathrm{ord}\,C^{+}_{j}=
\mathrm{ord}\,C_{j}+\mathrm{ord}\,B^{+}_{j}=2q-1.
\end{equation}
Кроме того, здесь и далее через $(\cdot,\cdot)_{\Omega}$ и $(\cdot,\cdot)_{\Gamma}$
обозначены скалярные произведения в пространствах $L_{2}(\Omega)$ и $L_{2}(\Gamma)$,
а также естественные расширения по непрерывности этих скалярных произведений.

Дифференциальное выражение $L^{+}$ называется формально сопряженным к выражению $L$,
а система граничных выражений $\{B^{+}_{1},\ldots,B^{+}_{q}\}$ называется
сопряженной к системе $\{B_{1},\ldots,B_{q}\}$ относительно выражения $L$.
Сопряженная система определяется неоднозначно, но все сопряженные системы
эквивалентны в следующем смысле [\ref{FunctionalAnalysis72}, с.~168;
\ref{Schechter77}, с.~232]: если
$\{\widetilde{B}^{+}_{1},\ldots,\widetilde{B}^{+}_{q}\}$
--- еще одна система, сопряженная  к $\{B_{1},\ldots,B_{q}\}$ относительно $L$, то
\begin{gather}\notag
\{v\in C^{\infty}(\,\overline{\Omega}\,):\;
B^{+}_{j}v=0\;\;\mbox{на}\;\;\Gamma,\;\;j=1,\ldots,q\}=\\
\{v\in C^{\infty}(\,\overline{\Omega}\,):\;
\widetilde{B}^{+}_{j}v=0\;\;\mbox{на}\;\;\Gamma,\;\;j=1,\ldots,q\}.\label{3.8}
\end{gather}

Известно [\ref{FunctionalAnalysis72}, с.~168; \ref{Schechter77}, с.~228], что
краевая задача регулярная эллиптическая тогда и только тогда, когда формально
сопряженная к ней задача является регулярной эллиптической.

Положим
\begin{gather*}
N:=\{u\in C^{\infty}(\,\overline{\Omega}\,):\;Lu=0\;\,\mbox{в}\;\,\Omega,\;\,
B_{j}u=0\;\,\mbox{на}\;\,\Gamma,\;\,j=1,\ldots,q\},\\
N^{+}:=\{v\in C^{\infty}(\,\overline{\Omega}\,):\\
L^{+}v=0\;\,\mbox{в}\;\,\Omega,\;\,
B^{+}_{j}v=0\;\,\mbox{на}\;\,\Gamma,\;\,j=1,\ldots,q\}.
\end{gather*}
В силу \eqref{3.8} множество $N^{+}$ не зависит от выбора сопряженной системы
граничных выражений $\{B^{+}_{1},\ldots,B^{+}_{q}\}$. Поскольку задачи \eqref{3.1},
\eqref{3.2} и \eqref{3.4}, \eqref{3.5} являются регулярными эллиптическими,
пространства $N$ и $N^{+}$ конечномерны [\ref{LionsMagenes71}, с.~191;
\ref{FunctionalAnalysis72}, с.~168].


\begin{example}\label{ex3.2}
Рассмотрим задачу Дирихле для дифференциального уравнения $Lu=f$, где выражение $L$
правильно эллиптическое на $\overline{\Omega}$. Для этой краевой задачи сопряженной
будет задача Дирихле для уравнения $L^{+}v=\omega$ [\ref{FunctionalAnalysis72},
с.~168]. При этом $\dim N=\dim N^{+}$ [\ref{LionsMagenes71}, с.~227].
\end{example}


В этой главе изучаются \emph{полуоднородные} регулярные эллиптические краевые
задачи, т. е. предполагается, что в \eqref{3.1}, \eqref{3.2} либо $f=0$ в области
$\Omega$, либо все $g_{j}=0$ на границе $\Gamma$. Эти важные классы
задач будут рассмотрены отдельно в п. 3.3 и 3.4.




\markright{\emph \ref{sec3.2}. Пространства Хермандера для евклидовых областей}

\section[Пространства Хермандера для евклидовых облас-\break тей]
{Пространства Хермандера \\ для евклидовых областей}\label{sec3.2}

\markright{\emph \ref{sec3.2}. Пространства Хермандера для евклидовых областей}

Здесь мы изучим отдельно классы пространств Хермандера для открытых и замкнутых
евклидовых областей. Первые состоят из распределений, \emph{заданных в} открытой
области, а вторые образованы распределениями в $\mathbb{R}^{n}$,
\emph{сосредоточенными на} замкнутой области [\ref{VolevichPaneah65}, с. 20, 22]. В
качестве открытой области мы берем $\Omega$, а в качестве замкнутой области ---
$\overline{\Omega}$. Указанные пространства необходимы нам для исследования краевой
задачи \eqref{3.1}, \eqref{3.2}. Мы изучим связь между ними, а также пространствами
на границе $\Gamma$. Как обычно, $\mathcal{D}'(\Omega)$~--- топологическое линейное
пространство распределений, заданных в области $\Omega$.

\subsection{Пространства для открытых областей}\label{sec3.2.1}

Пусть $s\in\mathbb{R}$ и $\varphi\in\nobreak\mathcal{M}$.


\begin{definition}\label{def3.2}
Линейное пространство $H^{s,\varphi}(\Omega)$ состоит, по определению, из сужений в
область $\Omega$ всех распределений $w\in H^{s,\varphi}(\mathbb{R}^{n})$. В
пространстве $H^{s,\varphi}(\Omega)$ определена норма распределения $u\in
H^{s,\varphi}(\Omega)$ по формуле
\begin{gather}\notag
\|u\|_{H^{s,\varphi}(\Omega)}:=\\ \inf\left\{\,
\|w\|_{H^{s,\varphi}(\mathbb{R}^{n})}:\,w\in H^{s,\varphi}(\mathbb{R}^{n}),\;\;
w=u\;\;\mbox{на}\;\;\Omega\,\right\}.\label{3.9}
\end{gather}
\end{definition}


\begin{theorem}\label{th3.1}
Пространство $H^{s,\varphi}(\Omega)$ сепарабельно и гильбертово относительно нормы
\eqref{3.9}.
\end{theorem}


\textbf{Доказательство.} Согласно определению, $H^{s,\varphi}(\Omega)$ является
факторпространством гильбертова пространства $H^{s,\varphi}(\mathbb{R}^{n})$ по
подпространству
\begin{equation}\label{3.10}
\bigl\{\,w\in
H^{s,\varphi}(\mathbb{R}^{n}):\,\mathrm{supp}\,w\subseteq\widehat{\Omega}\,\bigr\}.
\end{equation}
(Линейное многообразие \eqref{3.10} замкнуто в топологии пространства
$H^{s,\varphi}(\mathbb{R}^{n})$ в силу непрерывного вложения
$H^{s,\varphi}(\mathbb{R}^{n})\hookrightarrow\mathcal{D}'(\mathbb{R}^{n})$, см. ниже
теорему \ref{th3.6}). Следовательно, пространство $H^{s,\varphi}(\Omega)$
гильбертово относительно скалярного произведения
\begin{equation}\label{3.11}
(u_{1},u_{2})_{H^{s,\varphi}(\Omega)}:=(w_{1}-\Pi w_{1},w_{2}-\Pi
w_{2})_{H^{s,\varphi}(\mathbb{R}^{n})}.
\end{equation}
Здесь $u_{j}\in H^{s,\varphi}(\Omega)$, $w_{j}\in H^{s,\varphi}(\mathbb{R}^{n})$,
$u_{j}=w_{j}$ в $\Omega$ при $j\in\{1,2\}$, а $\Pi$~--- ортопроектор в
$H^{s,\varphi}(\mathbb{R}^{n})$ на подпространство \eqref{3.10}. (Правая часть
равенства \eqref{3.11} не зависит от выбора распределений $w_{1}$ и $w_{2}$). Норма
\eqref{3.9} порождена этим скалярным произведением. Cепарабельность пространства
$H^{s,\varphi}(\Omega)$ следует в силу определения \ref{def3.2} из сепарабельности
пространства $H^{s,\varphi}(\mathbb{R}^{n})$. Теорема~\ref{th3.1} доказана.

\medskip

В частном случае $\varphi\equiv1$ пространство $H^{s,\varphi}(\Omega)$ мы обозначаем
также через $H^{s}(\Omega)$. Это --- пространство Соболева порядка $s\in\mathbb{R}$
в области $\Omega$ [\ref{Triebel80}, c. 384].

Класс пространств
\begin{equation}\label{3.12}
\{H^{s,\varphi}(\Omega):\,s\in\mathbb{R},\,\varphi\in\mathcal{M}\}
\end{equation}
мы называем уточненной шкалой в области $\Omega$. Изучим свойства этой шкалы.


\begin{theorem}\label{th3.2}
Пусть заданы функция $\varphi\in\mathcal{M}$ и положительные числа
$\varepsilon,\delta$. Тогда для произвольного $s\in\mathbb{R}$
\begin{equation}\label{3.13}
[H^{s-\varepsilon}(\Omega),H^{s+\delta}(\Omega)]_{\psi}= H^{s,\varphi}(\Omega)
\end{equation}
с эквивалентностью норм. Здесь $\psi$ --- интерполяционный параметр из теоремы
$\ref{th2.14}$.
\end{theorem}


\textbf{Доказательство.} Написанная слева в \eqref{3.13} пара пространств Соболева
допустимая. Рассмотрим оператор $R_{\Omega}$ сужения распределения
$u\in\mathcal{D}'(\mathbb{R}^{n})$ в область $\Omega$. Имеем линейные ограниченные
сюрьективные операторы
\begin{gather}\label{3.14}
R_{\Omega}:\,H^{s-\varepsilon}(\mathbb{R}^{n})\rightarrow
H^{s-\varepsilon}(\Omega),\;\;
R_{\Omega}:\,H^{s+\delta}(\mathbb{R}^{n})\rightarrow H^{s+\delta}(\Omega),\\
R_{\Omega}:\,H^{s,\varphi}(\mathbb{R}^{n})\rightarrow H^{s,\varphi}(\Omega).
\label{3.15}
\end{gather}
Применив интерполяцию с параметром $\psi$ в формуле \eqref{3.14}, получим
ограниченный оператор
$$
R_{\Omega}:\,
[H^{s-\varepsilon}(\mathbb{R}^{n}),H^{s+\delta}(\mathbb{R}^{n})]_{\psi}\rightarrow
[H^{s-\varepsilon}(\Omega),H^{s+\delta}(\Omega)]_{\psi}.
$$
Это в силу теоремы \ref{th2.14} влечет за собой ограниченность оператора
$$
R_{\Omega}:\,H^{s,\varphi}(\mathbb{R}^{n})\rightarrow
[H^{s-\varepsilon}(\Omega),H^{s+\delta}(\Omega)]_{\psi}.
$$
Отсюда, ввиду сюрьективности оператора \eqref{3.15}, получаем включение
\begin{equation}\label{3.16}
H^{s,\varphi}(\Omega)\subseteq
[H^{s-\varepsilon}(\Omega),H^{s+\delta}(\Omega)]_{\psi}.
\end{equation}

Докажем обратное включение и его непрерывность. В монографии [\ref{Triebel80},
c.~386] (теорема 4.2.2) для каждого $k\in\mathbb{N}$ построено линейное отображение
$T_{k}$, продолжающее произвольное распределение $u\in H^{-k}(\Omega)$ в
пространство $\mathbb{R}^{n}$ и являющееся ограниченным оператором
\begin{equation}\label{3.17}
T_{k}:\ H^{\sigma}(\Omega)\rightarrow
H^{\sigma}(\mathbb{R}^{n})\quad\mbox{при}\;\;|\sigma|<k.
\end{equation}
Возьмем $k\in\mathbb{N}$ такое, что $|s-\varepsilon|<k$, $|s+\delta|<k$ и рассмотрим
ограниченные операторы \eqref{3.17} для $\sigma=s-\varepsilon$ и для
$\sigma=s+\nobreak\delta$. Поскольку $\psi$~--- интерполяционный параметр, это
влечет ограниченность оператора
$$
T_{k}:\,[H^{s-\varepsilon}(\Omega),H^{s+\delta}(\Omega)]_{\psi}\rightarrow
[H^{s-\varepsilon}(\mathbb{R}^{n}),H^{s+\delta}(\mathbb{R}^{n})]_{\psi}.
$$
Отсюда в силу теоремы \ref{th2.14} получаем ограниченный оператор
\begin{equation}\label{3.18}
T_{k}:\,[H^{s-\varepsilon}(\Omega),H^{s+\delta}(\Omega)]_{\psi}\rightarrow
H^{s,\varphi}(\mathbb{R}^{n}).
\end{equation}
Произведение ограниченных операторов \eqref{3.15} и \eqref{3.18} дает ограниченный
тождественный оператор
$$
I=R_{\Omega}T_{k}:\,[H^{s-\varepsilon}(\Omega),H^{s+\delta}(\Omega)]_{\psi}\rightarrow
H^{s,\varphi}(\Omega).
$$
Таким образом, наряду с включением \eqref{3.16} выполняется обратное ему непрерывное
вложение. Следовательно, справедливо равенство пространств \eqref{3.13}, причем, по
теореме Банаха об обратном операторе, нормы в этих пространствах эквивалентны.

Теорема \ref{th3.2} доказана.

\medskip

Напомним, что
$$
C^{\infty}_{0}(\Omega):=\{u\in C^{\infty}(\Omega):\,\mathrm{supp}\,u\subset\Omega\}.
$$
Мы отождествляем функции $u\in C^{\infty}_{0}(\Omega)$ с их продолжением нулем в
$\mathbb{R}^{n}$. Из контекста будет ясно на каком множестве~--- $\Omega$ или
$\mathbb{R}^{n}$~--- задана функция $u\in C^{\infty}_{0}(\Omega)$.


\begin{theorem}\label{th3.3}
Пусть $s\in\mathbb{R}$ и $\varphi,\varphi_{1}\in\mathcal{M}$. Справедливы следующие
утверждения.
\begin{itemize}
\item [$\mathrm{(i)}$] Множество $C^{\infty}(\,\overline{\Omega}\,)$ плотно в
пространстве $H^{s,\varphi}(\Omega)$.
\item [$\mathrm{(ii)}$] Если $s<1/2$, то множество $C^{\infty}_{0}(\Omega)$
плотно в пространстве $H^{s,\varphi}(\Omega)$.
\item [$\mathrm{(iii)}$] Для произвольного числа $\varepsilon>0$ выполняется
компактное и плотное вложение $H^{s+\varepsilon,\varphi_{1}}(\Omega)\hookrightarrow
H^{s,\varphi}(\Omega)$.
\item [$\mathrm{(iv)}$] Функция $\varphi/\varphi_{1}$ ограничена в окрестности
$\infty$ тогда и только тогда, когда $H^{s,\varphi_{1}}(\Omega)\hookrightarrow
H^{s,\varphi}(\Omega)$. Это вложение плотно и непрерывно. Оно компактно тогда и
только тогда, когда $\varphi(t)/\varphi_{1}(t)\rightarrow0$ при
$t\rightarrow\infty$.
\end{itemize}
\end{theorem}


\textbf{Доказательство.} Утверждение (i) непосредственно следует из плотности
множества $C^{\infty}_{0}(\mathbb{R}^{n})$ в пространстве
$H^{s,\varphi}(\mathbb{R}^{n})$ (см. лемму 3.1 статьи [\ref{VolevichPaneah65},
с.~20]).

Утверждение (ii) верно в соболевском случае $\varphi\equiv1$ (см., например,
[\ref{Triebel80}, с.~412], теорема 4.7.1 (d)). Отсюда для произвольного
$\varphi\in\mathcal{M}$ оно выводится с помощью теорем \ref{th3.2} и \ref{th2.1}.

Пусть $s<1/2$. Согласно эти теоремам выполняется непрерывное и плотное вложение
$H^{s+\delta}(\Omega)\hookrightarrow H^{s,\varphi}(\Omega)$ при любом $\delta>0$.
Кроме того, если $s+\delta<1/2$, то множество $C^{\infty}_{0}(\Omega)$ плотно в
соболевском пространстве $H^{s+\delta}(\Omega)$. Следовательно это множество плотно
и в $H^{s,\varphi}(\Omega)$. Утверждение (ii) доказано.

Утверждения (iii) и (iv) содержится в теоремах 7.4 и 8.1 статьи
[\ref{VolevichPaneah65}, с. 44, 48]. Плотность вложений следует из утверждения (i).

Теорема \ref{th3.3} доказана.


\begin{theorem}\label{th3.4}
Пусть заданы функция $\varphi\in\mathcal{M}$ и целое число $k\geq0$. Тогда условие
\eqref{2.37} равносильно вложению
\begin{equation}\label{3.19}
H^{k+n/2,\varphi}(\Omega)\hookrightarrow C^{k}(\,\overline{\Omega}\,).
\end{equation}
Оно компактно.
\end{theorem}


\textbf{Доказательство.} Предположим, что функция $\varphi$ удовлетворяет условию
\eqref{2.37}. Тогда по теореме \ref{th2.15} (iii) справедливо непрерывное вложение
\begin{equation}\label{3.20}
H^{k+n/2,\varphi}(\mathbb{R}^{n})\hookrightarrow C^{k}_{\mathrm{b}}(\mathbb{R}^{n}).
\end{equation}
Отсюда следует непрерывное вложение \eqref{3.19}. Действительно, для произвольных
функций $u\in H^{k+n/2,\varphi}(\Omega)$ и $w\in H^{k+n/2,\varphi}(\mathbb{R}^{n})$
таких, что $u=w$ в области $\Omega$, имеем: $w\in
C^{k}_{\mathrm{b}}(\mathbb{R}^{n})$, $u\in C^{k}(\,\overline{\Omega}\,)$ и
$$
\|u\|_{C^{k}(\,\overline{\Omega}\,)}\leq\|w\|_{C^{k}_{\mathrm{b}}(\mathbb{R}^{n})}\leq
c\,\|w\|_{H^{k+n/2,\varphi}(\mathbb{R}^{n})},
$$
где $c$ --- норма оператора вложения \eqref{3.20}. Переходя в этом неравенстве к
инфимуму по указанным функциям $w$, получаем оценку
$$
\|u\|_{C^{k}(\,\overline{\Omega}\,)}\leq c\,\|u\|_{H^{k+n/2,\varphi}(\Omega)}.
$$
Непрерывное вложение \eqref{3.19} доказано. Его компактность устанавливается так же,
как и компактность вложения \eqref{2.105} из доказательства теоремы~\ref{th2.27}.
Остается отметить, что включение \eqref{2.105} влечет свойство \eqref{2.37} в силу
предложения \ref{prop2.5} (для $p=q=2$, $V:=\Omega$) и эквивалентности \eqref{2.41}.
Теорема \ref{th3.4} доказана.


\begin{remark}\label{rem3.3}
В случае $k=0$ эта теорема содержится в теореме 7.5 статьи Л.~Р.~Волевича и
Б.~П.~Панеяха [\ref{VolevichPaneah65}, с. 45].
\end{remark}


Изучим вопрос о существовании и свойствах следов на границе $\Gamma$ произвольных
распределений $u\in H^{s,\varphi}(\Omega)$. Поскольку $\Gamma$~--- замкнутое
компактное бесконечно гладкое многообразие размерности $n-1$, то на $\Gamma$
определена уточненная шкала пространств $H^{\sigma,\varphi}(\Gamma)$,
$\sigma\in\mathbb{R}$, $\varphi\in\mathcal{M}$.


\begin{theorem}\label{th3.5}
Для произвольных параметров $s>1/2$ и $\varphi\in\nobreak\mathcal{M}$ линейное
отображение
\begin{equation}\label{3.21}
u\,\mapsto\,u\upharpoonright\Gamma\;\mbox{--- след
функции}\;\;u\;\;\mbox{на}\;\;\Gamma,\quad u\in C^{\infty}(\,\overline{\Omega}\,),
\end{equation}
продолжается по непрерывности до ограниченного сюрьективного оператора
\begin{equation}\label{3.22}
R_{\Gamma}:\,H^{s,\varphi}(\Omega)\rightarrow H^{s-1/2,\varphi}(\Gamma).
\end{equation}
Он имеет линейный ограниченный правый обратный оператор
\begin{equation}\label{3.23}
S_{\Gamma}:\,H^{s-1/2,\varphi}(\Gamma)\rightarrow H^{s,\varphi}(\Omega),
\end{equation}
такой, что отображение $S_{\Gamma}$ не зависит от $s$ и $\varphi$.
\end{theorem}


\textbf{Доказательство.} В случае $\varphi\equiv1$ (пространства Соболева) эта
теорема известна: см., например, монографию Х.~Трибеля [\ref{Triebel80}, с.~484]
(лемма 5.4.4). Отсюда общий случай $\varphi\in\mathcal{M}$ мы выведем с помощью
интерполяции. Выберем число $\varepsilon>0$ так, чтобы $s-\varepsilon>1/2$. Имеем
ограниченные линейные операторы
\begin{gather*}
R_{\Gamma}:\,H^{s\mp\varepsilon}(\Omega)\rightarrow
H^{s\mp\varepsilon-1/2,\varphi}(\Gamma),\\
S_{\Gamma}:\,H^{s\mp\varepsilon-1/2,\varphi}(\Gamma)\rightarrow
H^{s\mp\varepsilon,\varphi}(\Omega).
\end{gather*}
Применив здесь интерполяцию с параметром $\psi$ из теоремы \ref{th2.14}, где берем
$\delta:=\varepsilon$, получим в силу теорем \ref{th2.21} и \ref{th3.2} ограниченные
операторы \eqref{3.22} и \eqref{3.23}. Поскольку $R_{\Gamma}S_{\Gamma}h=h$ для
любого $h\in H^{s-\varepsilon-1/2}(\Gamma)$, то оператор \eqref{3.23} правый
обратный к оператору \eqref{3.22}, и последний сюрьективен. Теорема \ref{th3.5}
доказана.


\begin{remark}\label{rem3.4}
Теоремы о следах на гиперплоскостях и плоских кусках границы доказаны для общих
пространств Хермандера в [\ref{Hermander65}, с. 60] и [\ref{VolevichPaneah65}, с.
36, 38, 46]
\end{remark}


\begin{corollary}\label{cor3.1}
Пусть $\sigma>0$ и $\varphi\in\mathcal{M}$. Тогда
\begin{gather}\label{3.24}
H^{\sigma,\varphi}(\Gamma)=\{R_{\Gamma}f:\,f\in H^{\sigma+1/2,\varphi}(\Omega)\},\\
\|h\|_{H^{\sigma,\varphi}(\Gamma)}\asymp \notag\\
\inf\,\bigl\{\,\|f\|_{H^{\sigma+1/2,\varphi}(\Omega)}:\,f\in
H^{\sigma+1/2,\varphi}(\Omega),\;R_{\Gamma}f=h\,\bigr\}. \label{3.25}
\end{gather}
\end{corollary}


\textbf{Доказательство.} В теореме \ref{th3.5} положим $s:=\sigma+1/2$. Равенство
\eqref{3.24} следует из сюрьективности оператора \eqref{3.22}. Докажем
эквивалентность норм \eqref{3.25}. Для произвольных функций $h\in
H^{\sigma,\varphi}(\Gamma)$ и $f\in H^{\sigma+1/2,\varphi}(\Omega)$, удовлетворяющих
равенству $R_{\Gamma}f=h$, имеем:
$$
\|h\|_{H^{\sigma,\varphi}(\Gamma)}\leq
c_{1}\,\|f\|_{H^{\sigma+1/2,\varphi}(\Omega)},
$$
где $c_{1}$ --- норма оператора \eqref{3.22}. Переходя в этом неравенстве к инфимуму
по указанным функциям $f$, получаем оценку
$$
\|h\|_{H^{\sigma,\varphi}(\Gamma)}\leq
c_{1}\inf\,\bigl\{\|f\|_{H^{\sigma+1/2,\varphi}(\Omega)}:\,f\in
H^{\sigma+1/2,\varphi}(\Omega),\;R_{\Gamma}f=h\bigr\}.
$$
Обратная оценка следует из того, что для функции $f:=S_{\Gamma}h$
$$
\|f\|_{H^{\sigma+1/2,\varphi}(\Omega)}\leq
c_{2}\,\|h\|_{H^{\sigma,\varphi}(\Gamma)},
$$
где $c_{2}$ --- норма оператора \eqref{3.23}. Следствие \ref{cor3.1} доказано.


\begin{remark}\label{rem3.5}
Если $s<1/2$, то отображение \eqref{3.21} нельзя продолжить до непрерывного
оператора $R_{\Gamma}:H^{s,\varphi}(\Omega)\rightarrow\mathcal{D}'(\Gamma)$.
Действительно, если бы это было возможно, то в силу теоремы \ref{th3.3} (ii) для
каждого распределения $f\in H^{s,\varphi}(\Omega)$ его след $R_{\Gamma}f=0$ на
$\Gamma$, что неверно, например, для функции $f:=1$ на $\overline{\Omega}$. Это
замечание сохраняет силу для пространств Соболева и в случае $s=1/2$, поскольку
множество $C^{\infty}_{0}(\Omega)$ плотно в $H^{1/2}(\Omega)$ [\ref{Triebel80},
с.~412] (теорема 4.7.1 (d)).
\end{remark}

\subsection{Пространства для замкнутых областей}\label{sec3.2.2}

Пусть $s\in\mathbb{R}$, $\varphi\in\nobreak\mathcal{M}$, а $Q$ --- произвольное
непустое замкнутое подмножество пространства $\mathbb{R}^{n}$.


\begin{definition}\label{def3.3}
Линейное пространство $H^{s,\varphi}_{Q}(\mathbb{R}^{n})$ состоит, по определению,
из всех распределений $w\in H^{s,\varphi}(\mathbb{R}^{n})$ таких, что их носители
$\mathrm{supp}\,w\subseteq Q$. Пространство $H^{s,\varphi}_{Q}(\mathbb{R}^{n})$
наделяется скалярным произведением и нормой из пространства
$H^{s,\varphi}(\mathbb{R}^{n})$.
\end{definition}


\begin{theorem}\label{th3.6}
Пространство $H^{s,\varphi}_{Q}(\mathbb{R}^{n})$ сепарабельно и гильбертово.
\end{theorem}


\textbf{Доказательство.} Пусть последовательность $(w_{j})$ фундаментальна в
$H^{s,\varphi}_{Q}(\mathbb{R}^{n})$. Поскольку пространство
$H^{s,\varphi}(\mathbb{R}^{n})$ полное, то она имеет в нем некоторый предел $w$.
Ввиду непрерывности вложения
$H^{s,\varphi}(\mathbb{R}^{n})\hookrightarrow\mathcal{D}'(\mathbb{R}^{n})$ получаем
сходимость $w_{j}\rightarrow w$ в $\mathcal{D}'(\mathbb{R}^{n})$ при
$j\rightarrow\infty$. Отсюда и из включений $\mathrm{supp}\,w_{j}\subseteq Q$
следует, что $\mathrm{supp}\,w\subseteq Q$. Таким образом, последовательность
$(w_{j})$ имеет предел $w$ в пространстве $H^{s,\varphi}_{Q}(\mathbb{R}^{n})$.
Полнота этого пространства доказана. Оно сепарабельно как подпространство
сепарабельного пространства $H^{s,\varphi}(\mathbb{R}^{n})$. Теорема \ref{th3.6}
доказана.

\medskip

В соболевском случае $\varphi\equiv1$ мы будем опускать индекс $\varphi$ в
обозначении пространства $H^{s,\varphi}_{Q}(\mathbb{R}^{n})$ и других пространств,
вводимых в этой главе.

Для нас важен случай, когда $Q:=\overline{\Omega}$. Изучим свойства пространства
$H^{s,\varphi}_{\overline{\Omega}}(\mathbb{R}^{n})$ и его связь с уточненной шкалой
в области~$\Omega$.


\begin{theorem}\label{th3.7}
Пусть заданы функция $\varphi\in\mathcal{M}$ и положительные числа
$\varepsilon,\delta$. Тогда для произвольного $s\in\mathbb{R}$
\begin{equation}\label{3.26}
[H^{s-\varepsilon}_{\overline{\Omega}}(\mathbb{R}^{n}),
H^{s+\delta}_{\overline{\Omega}}(\mathbb{R}^{n})]_{\psi}
=H^{s,\varphi}_{\overline{\Omega}}(\mathbb{R}^{n})
\end{equation}
с эквивалентностью норм. Здесь $\psi$ --- интерполяционный параметр из теоремы
$\ref{th2.14}$.
\end{theorem}


\textbf{Доказательство.} Поскольку $\Omega$ --- ограниченная область с бесконечно
гладкой границей, то (см., например, [\ref{Triebel80}, c.~395], теорема 4.3.2/1 (b))
множество $C^{\infty}_{0}(\Omega)$ плотно в пространстве
$H^{\sigma}_{\overline{\Omega}}(\mathbb{R}^{n})$ при любом $\sigma\in\mathbb{R}$.
Следовательно, непрерывное вложение
$H^{s+\delta}_{\overline{\Omega}}(\mathbb{R}^{n})\hookrightarrow
H^{s-\varepsilon}_{\overline{\Omega}}(\mathbb{R}^{n})$ плотно, пара пространств
$H^{s-\varepsilon}_{\overline{\Omega}}(\mathbb{R}^{n})$,
$H^{s+\delta}_{\overline{\Omega}}(\mathbb{R}^{n})$ допустимая, и левая часть формулы
\eqref{3.26} определена.

Мы выведем эту формулу из теоремы \ref{th2.14} с помощью теоремы \ref{th2.6}
(интерполяция подпространств). Нам понадобится отображение, являющееся проектором
каждого пространства $H^{\sigma}(\mathbb{R}^{n})$, где $s-\varepsilon\leq\sigma\leq
s+\delta$, на подпространство $H^{\sigma}_{\overline{\Omega}}(\mathbb{R}^{n})$.
Построим это отображение следующим образом. Возьмем такое число $r>0$, что $|x|<r$
для всех $x\in\Omega$, и положим
$$
G=\{x\in\mathbb{R}^{n}:\,|x|<4r\}\setminus\overline{\Omega}.
$$
Обозначим через $R$ отображение, которое каждому распределению  в $\mathbb{R}^{n}$
ставит в соответствие его сужение в область $G$. Получим линейный ограниченный
оператор
\begin{equation}\label{3.27}
R:\,H^{\sigma}(\mathbb{R}^{n})\rightarrow H^{\sigma}(G)\quad\mbox{для
каждого}\quad\sigma\in\mathbb{R}.
\end{equation}
Заметим здесь, что $G$ --- ограниченная открытая область с бесконечно гладкой
границей. Тогда (см., например, [\ref{Triebel80}, c.~386], теорема 4.2.2) для любого
компакта $K\subset\mathbb{R}$ cуществует линейное отображениe $T$, которое является
ограниченным оператором
\begin{equation}\label{3.28}
T:\,H^{\sigma}(G)\rightarrow H^{\sigma}(\mathbb{R}^{n})\;\;\mbox{для
каждого}\;\;\sigma\in K
\end{equation}
и продолжает распределение с области $G$ в $\mathbb{R}^{n}$. Это означает, что
оператор \eqref{3.28} правый обратный к оператору \eqref{3.27}. Возьмем
$K=[s-\varepsilon,\,s+\delta]$ и рассмотрим отображение
$$
P_{0}:\,u\mapsto u-TR\,u,\quad u\in H^{s-\varepsilon}(\mathbb{R}^{n}).
$$
В силу ограниченности операторов \eqref{3.27} и \eqref{3.28} имеем линейный
ограниченный оператор
$$
P_{0}:\,H^{\sigma}(\mathbb{R}^{n})\rightarrow
H^{\sigma}(\mathbb{R}^{n})\quad\mbox{для каждого}\;\;\sigma\in
[s-\varepsilon,\,s+\delta].
$$
Он проектирует пространство $H^{\sigma}(\mathbb{R}^{n})$ на подпространство
$H^{\sigma}_{\widehat{G}}(\mathbb{R}^{n})$ при каждом $\sigma\in
[s-\varepsilon,\,s+\delta]$. В самом деле,
\begin{gather*}
u\in H^{\sigma}(\mathbb{R}^{n})\,\Rightarrow \\
RP_{0}\,u=R(u-TR\,u)=Ru-RTR\,u=Ru-Ru=0\,\Rightarrow\\ P_{0}\,u\in
H^{\sigma}_{\widehat{G}}(\mathbb{R}^{n}).
\end{gather*}
Кроме того,
$$
u\in H^{\sigma}_{\widehat{G}}(\mathbb{R}^{n})\,\Rightarrow\,Ru=0\,\Rightarrow\
P_{0}u=u-TR\,u=u.
$$

Теперь выберем функцию $\chi\in C^{\infty}(\mathbb{R}^{n})$ такую, что $\chi(x)=1$
при $|x|\leq2r$ и $\chi(x)=0$ при $|x|\geq3r$, и рассмотрим отображение
$$
P:\,u\mapsto\chi\cdot P_{0}\,u,\quad u\in H^{s-\varepsilon}(\mathbb{R}^{n}).
$$
Поскольку оператор умножения на функцию $\chi$ ограничен в каждом пространстве
Соболева в $\mathbb{R}^{n}$, то $P$
--- проектор пространства $H^{\sigma}(\mathbb{R}^{n})$ на подпространство
$H^{\sigma}_{\overline{\Omega}}(\mathbb{R}^{n})$ при любом $\sigma\in
[s-\varepsilon,\,s+\delta]$. Следовательно, в силу теорем \ref{th2.6} и \ref{th2.14}
имеем следующие равенства пространств с точностью до эквивалентности норм в них:
\begin{gather*}
[H^{s-\varepsilon}_{\overline{\Omega}}(\mathbb{R}^{n}),
H^{s+\delta}_{\overline{\Omega}}(\mathbb{R}^{n})]_{\psi}=\\
[H^{s-\varepsilon}(\mathbb{R}^{n}),H^{s+\delta}(\mathbb{R}^{n})]_{\psi}\,
\cap\,H^{s-\varepsilon}_{\overline{\Omega}}(\mathbb{R}^{n})=\\
H^{s,\varphi}(\mathbb{R}^{n})\,
\cap\,H^{s-\varepsilon}_{\overline{\Omega}}(\mathbb{R}^{n})=\,
H^{s,\varphi}_{\overline{\Omega}}(\mathbb{R}^{n}).
\end{gather*}

Теорема \ref{th3.7} доказана.


\begin{theorem}\label{th3.8}
Пусть $s\in\mathbb{R}$ и $\varphi\in\mathcal{M}$. Справедливы следующие утверждения.
\begin{itemize}
\item [$\mathrm{(i)}$] Множество $C^{\infty}_{0}(\Omega)$ плотно в
$H^{s,\varphi}_{\overline{\Omega}}(\mathbb{R}^{n})$.
\item [$\mathrm{(ii)}$] Если $|s|<1/2$, то пространства $H^{s,\varphi}(\Omega)$
и $H^{s,\varphi}_{\overline{\Omega}}(\mathbb{R}^{n})$ равны как пополнения линейного
многообразия $C^{\infty}_{0}(\Omega)$ по эквивалентным нормам.
\item [$\mathrm{(iii)}$] Пространства $H^{s,\varphi}_{\overline{\Omega}}(\mathbb{R}^{n})$ и
$H^{-s,1/\varphi}(\Omega)$ взаимно двойственны относительно расширения по
непрерывности скалярного произведения в $L_{2}(\Omega)$.
\item [$\mathrm{(iv)}$] Утверждения $\mathrm{(iii)}$ и $\mathrm{(iv)}$ теоремы
$\ref{th3.3}$ сохраняют силу, если в их формулировках заменить пространства
$H^{\cdot,\cdot}(\Omega)$ на пространства
$H^{\cdot,\cdot}_{\overline{\Omega}}(\mathbb{R}^{n})$ с теми же верхними индексами.
\end{itemize}
\end{theorem}


\textbf{Доказательство.} Утверждения (i), (iii) и (iv) содержатся в соответственно в
лемме 3.3, п. 3.4 и теоремах 7.1, 8.1 статьи [\ref{VolevichPaneah65}, с. 23, 25, 42,
48].

Докажем утверждение (ii). Предположим, что $|s|<1/2$. Утверждение (ii) в соболевском
случае $\varphi\equiv1$ доказано, например, в [\ref{Triebel80}, c. 395] (теорема
4.3.2/1 (a), (с)). Отсюда общий случай $\varphi\in\mathcal{M}$ выведем с помощью
интерполяции. Выберем число $\varepsilon>\nobreak0$ так, чтобы
$|s\mp\varepsilon|<1/2$. Тождественное отображение, заданное на множестве
$C^{\infty}_{0}(\Omega)$, продолжается по непрерывности до гомеоморфизмов
$$
I:\,H^{s\mp\varepsilon}_{\overline{\Omega}}(\mathbb{R}^{n})\leftrightarrow
H^{s\mp\varepsilon}(\Omega).
$$
Применив здесь интерполяцию с параметром $\psi$ из теоремы \ref{th2.14}, где берем
$\delta:=\varepsilon$, получим в силу теорем \ref{th3.2} и \ref{th3.7} еще один
гомеоморфизм
$$
I:\,H^{s,\varphi}_{\overline{\Omega}}(\mathbb{R}^{n})\leftrightarrow
H^{s,\varphi}(\Omega),
$$
что и доказывает утверждение (ii).

Теорема \ref{th3.8} доказана.

\subsection[Оснащение $L_{2}(\Omega)$ пространствами Херман-\break дера]
{Оснащение $L_{2}(\Omega)$ пространствами \\ Хермандера}\label{sec3.2.3}

При изучении эллиптической задачи с однородными краевыми условиями будет полезна
следующая шкала пространств Хермандера.


\begin{definition}\label{def3.4}
Пусть $s\in\mathbb{R}$ и $\varphi\in\mathcal{M}$. В случае $s\geq0$ обозначим через
$H^{s,\varphi,(0)}(\Omega)$ гильбертово пространство $H^{s,\varphi}(\Omega)$.
В~случае $s<0$ обозначим через $H^{s,\varphi,(0)}(\Omega)$ гильбертово пространство
$H^{s,\varphi}_{\overline{\Omega}}(\mathbb{R}^{n})$, двойственное в силу теоремы
\ref{th3.8} (iii) к пространству $H^{-s,1/\varphi}(\Omega)$ относительно расширения
по непрерывности скалярного произведения в $L_{2}(\Omega)$.
\end{definition}


Шкала пространств
\begin{equation}\label{3.29}
\{H^{s,\varphi,(0)}(\Omega):\,s\in\mathbb{R},\;\varphi\in\mathcal{M}\}
\end{equation}
двусторонняя по параметру $s$ и уточненная по параметру $\varphi$. При $s\geq0$
имеем позитивную часть шкалы, состоящую из пространств
$H^{s,\varphi,(0)}(\Omega)=H^{s,\varphi}(\Omega)$ распределений, заданных в области
$\Omega$. При $s<0$ имеем негативную часть шкалы, состоящую из пространств
$H^{s,\varphi,(0)}(\Omega)=H^{s,\varphi}_{\overline{\Omega}}(\mathbb{R}^{n})$
распределений, сосредоточенных в замыкании области $\Omega$. Таким образом, шкала
\eqref{3.29} образована пространствами распределений различной природы.

Мы отождествляем (что естественно) функции из пространства
$L_{2}(\Omega)=H^{0}(\Omega)$ с их продолжениями нулем в $\mathbb{R}^{n}$. В этом
смысле $L_{2}(\Omega)=H^{0}_{\overline{\Omega}}(\mathbb{R}^{n})$. В~силу теорем
\ref{th3.3} (i), (iii) и \ref{th3.8} (i) имеем непрерывные и плотные вложения
\begin{equation}\label{3.30}
H^{s,\varphi,(0)}(\Omega)\hookleftarrow L_{2}(\Omega)\hookleftarrow
H^{-s,1/\varphi,(0)}(\Omega),\quad s<0,\;\;\varphi\in\mathcal{M}.
\end{equation}
Здесь крайние пространства взаимно двойственны относительно расширения по
непрерывности скалярного произведения в $L_{2}(\Omega)$. Тем самым мы получили
\emph{оснащение} гильбертова пространства $L_{2}(\Omega)$ пространствами шкалы
\eqref{3.29}. (Определение гильбертова оснащения и связанные с ним понятия см. в
[\ref{Berezansky65}; \ref{BerezanskyUsSheftel90}, с. 462].)

В соболевском случае $\varphi\equiv1$ оснащение \eqref{3.30} введено и изучено
Ю.~М.~Березанским [\ref{Berezansky65}; \ref{BerezanskyUsSheftel90}, с. 475]. В этом
случае пространство $H^{s,\varphi,(0)}(\Omega)$ будем обозначать также через
$H^{s,(0)}(\Omega)$.

Из плотности непрерывных вложений \eqref{3.30} следует, что функции класса
$C^{\infty}(\,\overline{\Omega}\,)$ (продолженные нулем в $\mathbb{R}^{n}$) образуют
плотное подмножество в каждом негативном пространстве $H^{s,\varphi,(0)}(\Omega)$,
$s<0$. Они же плотны и в каждом позитивном пространстве $H^{s,\varphi,(0)}(\Omega)$,
$s>0$. Это позволяет рассматривать непрерывное плотное вложение вида
$H^{r,\chi,(0)}(\Omega)\hookrightarrow H^{s,\varphi,(0)}(\Omega)$, где $s,r\in
\mathbb{R}$ и $\varphi,\chi\in\mathcal{M}.$ Оно естественным образом означает, что
$$
\|u\|_{H^{s,\varphi,(0)}(\Omega)}\leq\mathrm{const}\,
\|u\|_{H^{r,\chi,(0)}(\Omega)}\quad\mbox{для любого}\;\;u\in
C^{\infty}(\,\overline{\Omega}\,),
$$
причем тождественное отображение, заданное на множестве
$C^{\infty}(\,\overline{\Omega}\,)$, продолжается по непрерывности до ограниченного
\textit{инъективного} оператора $I:\nobreak H^{r,\chi,(0)}(\Omega)\rightarrow
H^{s,\varphi,(0)}(\Omega)$. Его называют оператором вложения пространства
$H^{r,\chi,(0)}(\Omega)$ в пространство $H^{s,\varphi,(0)}(\Omega)$ (по этому поводу
см. [\ref{BerezanskyUsSheftel90}, с. 504]).

Изучим свойства шкалы \eqref{3.29}.


\begin{theorem}\label{th3.9}
Пусть $s\in\mathbb{R}$ и $\varphi,\varphi_{1}\in\mathcal{M}$. Справедливы следующие
утверждения.
\begin{itemize}
\item [$\mathrm{(i)}$] С точностью до эквивалентности норм
\begin{gather}\label{3.31}
H^{s,\varphi,(0)}(\Omega)=H^{s,\varphi}_{\overline{\Omega}}(\mathbb{R}^{n})
\quad\mbox{при}\;\;s<1/2 \\
H^{s,\varphi,(0)}(\Omega)=H^{s,\varphi}(\Omega)\quad\mbox{при}\;\;s>-1/2.
\label{3.32}
\end{gather}
\item [$\mathrm{(ii)}$] Если $s<1/2$, то множество $C^{\infty}_{0}(\Omega)$ плотно
в $H^{s,\varphi,(0)}(\Omega)$.
\item [$\mathrm{(iii)}$] Пространства $H^{s,\varphi,(0)}(\Omega)$ и
$H^{-s,\,1/\varphi,(0)}(\Omega)$ взаимно двойственны (при $s\neq0$ с равенством
норм, а при $s=0$ с эквивалентностью норм) относительно расширения по непрерывности
скалярного произведения в $L_{2}(\Omega)$.
\item [$\mathrm{(iv)}$] Для произвольного числа $\varepsilon>0$ выполняется
компактное и плотное вложение
$H^{s+\varepsilon,\varphi_{1},(0)}(\Omega)\hookrightarrow
H^{s,\varphi,(0)}(\Omega)$.
\end{itemize}
\end{theorem}


\textbf{Доказательство.} (i) Равенства \eqref{3.31} и \eqref{3.32} сразу следуют из
теоремы \ref{th3.8} (ii) и определения пространства $H^{s,\varphi,(0)}(\Omega)$.

(ii) Это утверждение является следствием равенства \eqref{3.31} и теоремы
\ref{th3.8}~(i).

(iii) В случае $s\neq0$ утверждение (iii) содержится в определении \ref{def3.4}.
В~случае $s=0$ оно следует из теоремы \ref{th3.8} (iii) и формулы \eqref{3.31}:
$$
(H^{0,\varphi,(0)}(\Omega))'=(H^{0,\varphi}(\Omega))'=
H^{0,1/\varphi}_{\overline{\Omega}}(\mathbb{R}^{n})=H^{0,1/\varphi,(0)}(\Omega)
$$
(последнее равенство выполняется с точностью до эквивалентности норм).

(iv) В случае $s\geq0$ утверждение (iv) совпадает с теоремой \ref{th3.3} (iii).
Случай $s+\varepsilon\leq0$ содержится в теореме \ref{th3.8} (iv); при этом мы
используем \eqref{3.31}, если $s+\varepsilon=0$. Теперь оставшийся случай, когда
$s<0<s+\varepsilon$ является следствием компактных плотных вложений
$$
H^{s+\varepsilon,\varphi_{1},(0)}(\Omega)\hookrightarrow
H^{0,\varphi_{1},(0)}(\Omega)\hookrightarrow H^{s,\varphi,(0)}(\Omega).
$$

Теорема \ref{th3.9} доказана.



\begin{theorem}\label{th3.10}
Пусть заданы функция $\varphi\in\mathcal{M}$ и положительные числа
$\varepsilon,\delta$. Тогда для произвольного $s\in\mathbb{R}$
\begin{equation}\label{3.33}
[H^{s-\varepsilon,(0)}(\Omega),H^{s+\delta,(0)}(\Omega)]_{\psi}=
H^{s,\varphi,(0)}(\Omega)
\end{equation}
с эквивалентностью норм. Здесь $\psi$ --- интерполяционный параметр из теоремы
$\ref{th2.14}$.
\end{theorem}


\textbf{Доказательство.} В случае, когда $s-\varepsilon\geq0$ либо $s+\delta\leq0$
эта теорема непосредственно следует из интерполяционных теорем \ref{th3.2} и
\ref{th3.7}. Рассмотрим оставшийся случай, когда $s-\varepsilon<0<s+\delta$. Положим
$\lambda:=\min\{\varepsilon/2,\delta/2,1/4\}$. Поскольку $s\mp\lambda<1/2$ либо
$s\mp\lambda>-1/2$, то в силу теоремы \ref{th3.9} (i) и упомянутых интерполяционных
теорем
\begin{equation}\label{3.34}
[H^{s-\lambda,(0)}(\Omega),H^{s+\lambda,(0)}(\Omega)]_{\eta}=
H^{s,\varphi,(0)}(\Omega).
\end{equation}
Здесь интерполяционный параметр $\eta$ задан по формуле
$$
\eta(t):=
\begin{cases}
\;t^{1/2}\,
\varphi(t^{1/(2\lambda)}) & \text{при\;\;\;$t\geq1$}, \\
\;\varphi(1) & \text{при\;\;\;$0<t<1$}.
\end{cases}
$$

Кроме того, так как $s-\varepsilon<s\mp\lambda<s+\delta$, то по теореме 12.5 из
гл.~1 монографии Ж.-Л.~Лионса и Э.~Мадженеса [\ref{LionsMagenes71}, с.~97] мы имеем:
\begin{gather}\notag
[H^{s-\varepsilon,(0)}(\Omega),H^{s+\delta,(0)}(\Omega)]_{\psi_{\mp}}=\\
[(H^{-s+\varepsilon}(\Omega))',H^{s+\delta}(\Omega)]_{\psi_{\mp}}=
H^{s\mp\lambda,(0)}(\Omega).\label{3.35}
\end{gather}
Здесь $\psi_{\mp}(t):=t^{\theta_{\mp}}$, где число $\theta_{\mp}\in(0,\,1)$
определяется из условия
$s\mp\lambda=(1-\theta_{\mp})(s-\varepsilon)+\theta_{\mp}(s+\delta)$. Отсюда
получаем, что $\theta_{\mp}=(\varepsilon\mp\lambda)/(\varepsilon+\delta)$.

Теперь равенства \eqref{3.34} и \eqref{3.35} (с эквивалентностью норм) влекут в силу
теоремы реитерации \ref{th2.3} равенство
\begin{gather*}
H^{s,\varphi,(0)}(\Omega)=[H^{s-\lambda,(0)}(\Omega),H^{s+\lambda,(0)}(\Omega)]_{\eta}=\\
[\,[H^{s-\varepsilon,(0)}(\Omega),H^{s+\delta,(0)}(\Omega)]_{\psi_{-}},
[H^{s-\varepsilon,(0)}(\Omega),H^{s+\delta,(0)}(\Omega)]_{\psi_{+}}]_{\eta}=\\
[H^{s-\varepsilon,(0)}(\Omega),H^{s+\delta,(0)}(\Omega)]_{\omega}.
\end{gather*}
Здесь
\begin{gather*}
\omega(t):=\psi_{-}(t)\,\eta\Bigl(\frac{\psi_{+}(t)}{\psi_{-}(t)}\Bigr)=\\
t^{\theta_{-}}\,\eta(t^{\theta_{+}-\theta_{-}})=
t^{(\varepsilon-\lambda)/(\varepsilon+\delta)}\,\eta(t^{2\lambda/(\varepsilon+\delta)})=\\
t^{(\varepsilon-\lambda)/(\varepsilon+\delta)}\,t^{\lambda/(\varepsilon+\delta)}\,
\varphi(t^{1/(\varepsilon+\delta)})=\psi(t)\quad\mbox{при}\quad t\geq1.
\end{gather*}
Следовательно (см. замечание \ref{rem2.1}), справедливо равенство \eqref{3.34} с
точностью до эквивалентности норм.

Теорема \ref{th3.10} доказана.





\section[Краевая задача для однородного эллиптического уравнения]
{Краевая задача для однородного \\ эллиптического уравнения}\label{sec3.3}

\markright{\emph \ref{sec3.3}. Краевая задача для однородного уравнения}

Рассмотрим регулярную эллиптическую краевую задачу \eqref{3.1}, \eqref{3.2} в
случае, когда эллиптическое уравнение \eqref{3.1} однородно:
\begin{equation}\label{3.36}
L\,u=0\;\;\text{в}\;\;\Omega,\quad
B_{j}\,u=g_{j}\;\;\text{на}\;\;\Gamma\;\;\text{при}\;\; j=1,\ldots,q.
\end{equation}
Изучим оператор $B=(B_{1},\ldots,B_{q})$, отвечающий этой задаче, в уточненной шкале
пространств Хермандера.

\subsection[Основной результат: ограниченность и нетеровость оператора]
{Основной результат: \\ ограниченность и нетеровость оператора}\label{sec3.3.1}

Положим
\begin{equation}\label{3.37}
K_{L}^{\infty}(\Omega):=\{\,u\in
C^{\infty}(\,\overline{\Omega}\,):\;L\,u=0\quad\mbox{в}\quad\Omega\,\}.
\end{equation}
Cвяжем с задачей \eqref{3.36} линейное отображение
\begin{equation}\label{3.38}
u\mapsto Bu=(B_{1}\,u,\ldots,B_{q}\,u),\;\;\mbox{где}\;\;u\in
K_{L}^{\infty}(\Omega).
\end{equation}

Пусть $s\in\mathbb{R}$ и $\varphi\in\mathcal{M}$. Обозначим
\begin{equation}\label{3.39}
K_{L}^{s,\varphi}(\Omega):=\{\,u\in
H^{s,\varphi}(\Omega):\;L\,u=0\;\;\mbox{в}\;\;\Omega\,\}.
\end{equation}
Поскольку вложение $H^{s,\varphi}(\Omega)\hookrightarrow\mathcal{D}'(\Omega)$
непрерывно, то $K_{L}^{s,\varphi}(\Omega)$~--- замкнутое подпространство в
$H^{s,\varphi}(\Omega)$.

Действительно, пусть $u\in H^{s,\varphi}(\Omega)$ и последовательность
$(u_{j})\subset K_{L}^{s,\varphi}(\Omega)$ такие, что $u_{j}\rightarrow u$ в
$H^{s,\varphi}(\Omega)$ при $j\rightarrow\infty$. Тогда $u_{j}\rightarrow u$ в
$\mathcal{D}'(\Omega)$, откуда $0=Lu_{j}\rightarrow Lu$ в $\mathcal{D}'(\Omega)$ при
$j\rightarrow\infty$. Следовательно, $Lu=0$ в области $\Omega$, т.~е. $u\in
K_{L}^{s,\varphi}(\Omega)$.

Мы рассматриваем $K_{L}^{s,\varphi}(\Omega)$ как гильбертово пространство
относительно скалярного произведения, индуцированного из $H^{s,\varphi}(\Omega)$.

\medskip

Сформулируем основной результат п. \ref{sec3.3}


\begin{theorem}\label{th3.11} Для произвольных параметров
$s\in\mathbb{R}$ и $\varphi\in\mathcal{M}$ множество $K_{L}^{\infty}(\Omega)$ плотно
в $K_{L}^{s,\varphi}(\Omega)$, а отображение $\eqref{3.38}$ продолжается по
непрерывности до ограниченного линейного оператора
\begin{equation}\label{3.40}
B:\,K_{L}^{s,\varphi}(\Omega)\rightarrow
\bigoplus_{j=1}^{q}\,H^{s-m_{j}-1/2,\,\varphi}(\Gamma)=:
\mathcal{H}_{s,\varphi}(\Gamma).
\end{equation}
Этот оператор нетеров, имеет ядро $N$ и область значений
\begin{equation}\label{3.41}
\Bigl\{(g_{1},\ldots,g_{q})\in\mathcal{H}_{s,\varphi}(\Gamma):\,
\sum_{j=1}^{q}\,(g_{j},C^{+}_{j}\,v)_{\Gamma}=0\;\mbox{для всех}\;v\in N^{+}\Bigl\}.
\end{equation}
Индекс оператора \eqref{3.40} не зависит от $s$ и $\varphi$.
\end{theorem}


Отметим, что в формуле \eqref{3.41} величина $(g_{j},C^{+}_{j}\,v)_{\Gamma}$~--- это
значение антилинейного функционала $g_{j}$ на функции $C^{+}_{j}\,v\in
C^{\infty}(\,\overline{\Omega}\,)$. Следовательно, множество \eqref{3.41} замкнуто в
пространстве $\mathcal{H}_{s,\varphi}(\Gamma)$. Далее, согласно теореме \ref{th3.11}
множество
\begin{equation}\label{3.42}
G:=\left\{\,\bigl(C_{1}^{+}v,\ldots,C_{q}^{+}v\bigr):\,v\in N^{+}\,\right\}.
\end{equation}
является дефектным подпространством оператора \eqref{3.40}: оно ортогонально области
значений этого оператора относительно расширения по непрерывности скалярного
произведения в пространстве $(L_{2}(\Gamma))^{q}$. Индекс оператора \eqref{3.40}
равен $\dim N-\dim G$. Ясно, что $\dim G\leq\dim N^{+}$, где возможно и строгое
неравенство. Это вытекает из [\ref{Hermander86}, с.~257].

В случае, когда $\varphi\equiv1$ и $s\notin\{1/2-k:\,k\in\mathbb{N}\}$, утверждение
теоремы \ref{th3.11} содержится в теореме Ж.-Л.~Лионса и Э.~Мадженеса
[\ref{LionsMagenes71}, с. 216, 217] о разрешимости неоднородной задачи \eqref{3.1},
\eqref{3.2} в двусторонней шкале пространств Соболева. Общий случай
$\varphi\in\mathcal{M}$ и $s\in\mathbb{R}$  мы выведим из теоремы Лионса и Мадженеса
с помощью интерполяции с подходящим функциональным параметром и последующего сужения
оператора задачи на пространство $K_{L}^{s,\varphi}(\Omega)$. Это будет сделано в
п.~\ref{sec3.3.4}. В~п.~\ref{sec3.3.2} и \ref{sec3.3.3} приведены необходимые для
этого утверждения из теории интерполяции и эллиптических краевых задач.

Отметим, что теорема \ref{th3.11} является новой даже в соболевском случае
$\varphi\equiv1$, если число $s<0$ полуцелое (см. ниже в п.~\ref{sec3.3.3}
замечание~\ref{rem3.6}). Эту теорему можно считать некоторым аналогом теоремы
Харнака о сходимости последовательностей гармонических функций (см., например,
[\ref{Mihlin68}, с. 247]). При этом вместо равномерной метрики используется метрика
пространства $H^{s,\varphi}(\Omega)$.

В связи с теоремой \ref{th3.11} отметим также исследование Р.~Сили [\ref{Seeley66};
\ref{Agranovich97}, с. 48 -- 52] данных Коши решений однородного эллиптического
уравнения в двусторонней шкале пространств Соболева.

\subsection[Одна теорема об интерполяции подпрост\-ранств]
{Одна теорема об интерполяции \\ подпространств}\label{sec3.3.2}

Здесь мы сформулируем и докажем одно (несколько громоздкое по формулировке)
утверждение об интерполяции подпространств, связанных с линейным оператором. Оно
сыграет важную роль в доказательстве основного результата. Для случая голоморфной
(комплексной) интерполяции это утверждение доказано в монографии Лионса и Мадженеса
[\ref{LionsMagenes71}, с. 119 -- 121] (теорема 14.3). Мы покажем, что оно
справедливо и для интерполяции гильбертовых пространств с функциональным параметром.
При этом наше доказательство в отличие от цитированной работы не будет использовать
конструкцию интерполяционного функтора.

Предварительно примем следующее обозначение. Пусть $H$, $\Phi$ и $\Psi$
--- гильбертовы пространства, причем непрерывно
$\Phi\hookrightarrow\Psi$. Пусть также задано линейный ограниченный оператор
$T:H\rightarrow\Psi$. Обозначим
$$
(H)_{T,\Phi}=\{u\in H:\,Tu\in\Phi\}.
$$
Пространство $(H)_{T,\Phi}$ гильбертово относительно скалярного произведения графика
$$
(u_{1},u_{2})_{(H)_{T,\Phi}}=(u_{1},u_{2})_{H}+(Tu_{1},Tu_{2})_{\Phi}
$$
и не зависит от $\Psi$.


\begin{theorem}\label{th3.12} Пусть задано шесть гильбертовых пространств $X_{0}$,
$Y_{0}$, $Z_{0}$, $X_{1}$, $Y_{1}$, $Z_{1}$ и три линейных отображения $T$, $R$,
$S$, которые удовлетворяют следующим условиям:
\begin{itemize}
\item [$\mathrm{(i)}$] пары $X=[X_{0},X_{1}]$ и $Y=[Y_{0},Y_{1}]$ допустимые;
\item [$\mathrm{(ii)}$] $Z_{0}$ и $Z_{1}$ являются подпространствами некоторого линейного
пространства $E$;
\item [$\mathrm{(iii)}$] справедливы непрерывные вложения $Y_{j}\hookrightarrow Z_{j}$
при $j\in\{0,1\}$;
\item [$\mathrm{(iv)}$] отображение $T$ задано на $X_{0}$ и определяет ограниченные операторы
$T:X_{j}\rightarrow Z_{j}$ при $j\in\{0,1\}$;
\item [$\mathrm{(v)}$] отображение $R$ задано на $E$ и определяет ограниченные операторы
$R:Z_{j}\rightarrow X_{j}$ при $j\in\{0,1\}$;
\item [$\mathrm{(vi)}$] отображение $S$ задано на $E$ и определяет ограниченные операторы
$S:Z_{j}\rightarrow Y_{j}$ при $j\in\{0,1\}$;
\item [$\mathrm{(vii)}$] для любого $\omega\in E$ справедливо $TR\,\omega=\omega+S\omega$.
\end{itemize}
Тогда пара пространств $[\,(X_{0})_{T,Y_{0}},\,(X_{1})_{T,Y_{1}}\,]$ допустимая и
для произвольного интерполяционного параметра $\psi\in\mathcal{B}$ справедливо
следующее равенство пространств с точностью до эквивалентности норм:
\begin{equation}\label{3.43}
[\,(X_{0})_{T,Y_{0}},\,(X_{1})_{T,Y_{1}}\,]_{\psi} =(X_{\psi})_{T,Y_{\psi}}.
\end{equation}
\end{theorem}


\textbf{Доказательство.} Ввиду условий (iii) и (iv) пространства
$(X_{j})_{T,Y_{j}}$, где $j\in\{0,1\}$, определены корректно. Покажем, что
пространство в правой части равенства \eqref{3.43} также определено корректно. В
силу условия (i) определены пространства $X_{\psi}$ и $Y_{\psi}$; для них
справедливы непрерывные вложения $X_{\psi}\hookrightarrow X_{0}$ и
$Y_{\psi}\hookrightarrow Y_{0}$. Теперь первое вложение и условие (iv) при $j=0$
влекут ограниченность оператора $T:X_{\psi}\rightarrow Z_{0}$. Кроме того, второе
вложение и условие (iii) влекут непрерывность вложения $Y_{\psi}\hookrightarrow
Z_{0}$. Таким образом, пространство в правой части равенства \eqref{3.43} определено
корректно и является гильбертовым пространством, как и пространства
$(X_{j})_{T,Y_{j}}$ при $j\in\{0,1\}$.

Далее в доказательстве нам понадобится отображение
\begin{equation}\label{3.44}
Pu=-RTu+u,\quad u\in X_{0}.
\end{equation}
В силу условий (iv) и (v), оператор $P:X_{j}\rightarrow X_{j}$ ограничен при любом
$j\in\{0,1\}$. Более того, условия (vi) и (vii) влекут для произвольного $u\in
X_{j}$ следующее:
$$
TPu=-TRTu+Tu=-(Tu+STu)+Tu=-STu\in Y_{j},
$$
т. е. $Pu\in(X_{j})_{T,Y_{j}}$. Кроме того, из ограниченности оператора
$P:X_{j}\rightarrow X_{j}$ и условий (iv), (vi) вытекает оценка
\begin{gather*}
\|Pu\|_{(X_{j})_{T,Y_{j}}}^{2}=\|Pu\|^{2}_{X_{j}}+ \|TPu\|^{2}_{Y_{j}}=\\
\|Pu\|^{2}_{X_{j}}+ \|-STu\|^{2}_{Y_{j}}\leq c_{1}\|u\|_{X_{j}}^{2},
\end{gather*}
в которой число $c_{1}>0$ не зависит от $u$. Таким образом, отображение \eqref{3.44}
задает ограниченные операторы
\begin{equation}\label{3.45}
P:X_{j}\rightarrow(X_{j})_{T,Y_{j}}\quad\mbox{при каждом}\;\;j\in\{0,1\}.
\end{equation}

Рассмотрим также сужения отображения $R$ на $Y_{j}$ при $j\in\{0,1\}$. В силу
условий (iii) и (v) мы имеем ограниченный оператор $R:Y_{j}\rightarrow X_{j}$. Более
того, из условий (vi) и (vii) вытекает, что для любого $\omega\in Y_{j}$ выполняется
равенство $TR\,\omega=\omega+S\omega\in Y_{j}$, т.~е.
$R\,\omega\in(X_{j})_{T,Y_{j}}$. Кроме того, из условий (iii), (vi), (vii) и
ограниченности оператора $R:Y_{j}\rightarrow X_{j}$ следует оценка
\begin{gather*}
\|R\,\omega\|_{(X_{j})_{T,Y_{j}}}^{2}= \|R\,\omega\|^{2}_{X_{j}}+
\|TR\,\omega\|^{2}_{Y_{j}}= \\ \|R\,\omega\|^{2}_{X_{j}}+
\|\omega+S\omega\|^{2}_{Y_{j}}\leq\\
\|R\,\omega\|^{2}_{X_{j}}+ \bigl(\|\omega\bigr\|_{Y_{j}}+
\bigl\|S\omega\|_{Y_{j}}\bigr)^{2}\leq\\
c_{2}\,\|\omega\|_{Y_{j}}^{2}+\bigl(\|\omega\|_{Y_{j}}+
c_{3}\,\|\omega\|_{Z_{j}}\bigr)^{2}\leq c_{4}\,\|\omega\|_{Y_{j}}^{2}.
\end{gather*}
Здесь постоянными  $c_{2}$, $c_{3}$ и $c_{4}$ не зависят от $\omega$. Таким образом,
ограничены операторы
\begin{equation}\label{3.46}
R:Y_{j}\rightarrow(X_{j})_{T,Y_{j}}\quad\mbox{при}\quad j\in\{0,1\}.
\end{equation}

Теперь с помощью операторов \eqref{3.45} и \eqref{3.46} покажем, что пара
пространств
\begin{equation}\label{3.47}
[\,(X_{\,0})_{T,Y_{0}},\,(X_{1})_{T,Y_{1}}\,]
\end{equation}
допустимая. Установим сначала сепарабельность пространства $(X_{j})_{T,Y_{j}}$ при
каждом $j\in\{0,1\}$. В силу условия (i) пространства $X_{j}$ и $Y_{j}$
сепарабельные. Возьмем какие-либо счетные множества $X_{j}^{0}$ и $Y_{j}^{0}$,
лежащие и плотные в $X_{j}$ и $Y_{j}$ соответственно. Построим по ним счетное
множество
$$
Q=\{Pu_{0}+Rv_{0}:\,u_{0}\in X_{j}^{0},\,v_{0}\in Y_{j}^{0}\}
$$
и аппроксимируем его элементами произвольное $u\in(X_{j})_{T,Y_{j}}$. Поскольку
$u\in\nobreak X_{j}$ и $Tu\in Y_{j}$, то найдутся такие последовательности элементов
$u_{k}\in X_{j}^{0}$ и $v_{k}\in Y_{j}^{0}$, что $u_{k}\rightarrow u$ в $X_{j}$ и
$v_{k}\rightarrow Tu$ в $Y_{j}$ при $k\rightarrow\infty$. Отсюда с помощью
ограниченных операторов \eqref{3.45}, \eqref{3.46} и равенства \eqref{3.44} имеем
\begin{equation}\label{3.48}
w_{k}:=Pu_{k}+Rv_{k}\rightarrow
Pu+RTu=u\;\,\mbox{в}\;\,(X_{j})_{T,Y_{j}}\;\,\mbox{при}\;\,k\rightarrow\infty,
\end{equation}
причем $w_{k}\in Q$. Значит, счетное множество $Q$ плотно в пространстве
$(X_{j})_{T,Y_{j}}$, т. е. последнее сепарабельно.

Для доказательства того, что эта пара \eqref{3.47} допустимая остается установить
плотность непрерывного вложения $(X_{1})_{T,Y_{1}}\hookrightarrow
(X_{0})_{T,Y_{0}}$. Возьмем произвольное $u\in (X_{0})_{T,Y_{0}}$; тогда $u\in
X_{0}$ и $Tu\in Y_{0}$. В силу условия (i) пространство $X_{1}$ плотно в $X_{0}$, а
пространство $Y_{1}$ плотно в $Y_{0}$. Следовательно, существуют такие
последовательности элементов $u_{k}\in X_{1}$ и $v_{k}\in Y_{1}$, что
$u_{k}\rightarrow u$ в $X_{\,0}$ и $v_{k}\rightarrow Tu$ в $Y_{0}$ при
$k\rightarrow\infty$. Отсюда с помощью операторов \eqref{3.45}, \eqref{3.46} и
равенства \eqref{3.44} имеем \eqref{3.48} для $j=0$ и $w_{k}\in (X_{1})_{T,Y_{1}}$.
Таким образом, $(X_{1})_{T,Y_{1}}$ плотно в $(X_{0})_{T,Y_{0}}$.

Перейдем к доказательству формулы \eqref{3.43}. Установим сначала непрерывное
вложение
\begin{equation}\label{3.49}
[\,(X_{0})_{T,Y_{0}},\,(X_{1})_{T,Y_{1}}\,]_{\psi}\hookrightarrow
(X_{\psi})_{T,Y_{\psi}}.
\end{equation}
В силу определения пространства $(X_{j})_{T,Y_{j}}$ ограничены операторы
$$
I:(X_{j})_{T,Y_{j}}\rightarrow X_{j}\quad\mbox{и}\quad
T:(X_{j})_{T,Y_{j}}\rightarrow Y_{j}\quad\mbox{при}\quad j\in\{0,1\}.
$$
Здесь, как обычно, через $I$ обозначено тождественное отображение. Отсюда, поскольку
параметр $\psi$ интерполяционный, получаем ограниченность операторов
\begin{gather*}
I:[\,(X_{0})_{T,Y_{0}},\,(X_{1})_{T,Y_{1}}\,]_{\psi}\rightarrow X_{\psi},\\
T:[\,(X_{0})_{T,Y_{0}},\,(X_{1})_{T,Y_{1}}\,]_{\psi}\rightarrow Y_{\psi}.
\end{gather*}
Значит, если $u\in[\,(X_{0})_{T,Y_{0}},\,(X_{1})_{T,Y_{1}}\,]_{\psi}$, то $u\in
X_{\psi}$, $Tu\in\nobreak Y_{\psi}$, причем
$$
\|u\|_{X_{\psi}}^{2}+\|Tu\|_{Y_{\psi}}^{2} \leq
c\,\bigl\|u\bigr\|_{[(X_{0})_{T,Y_{0}},\,(X_{1})_{T,Y_{1}}]_{\psi}}^{2}
$$
с некоторой постоянной с, не зависящей от $u$. Иными словами, справедливо
непрерывное вложение \eqref{3.49}.

Теперь в силу теоремы Банаха об обратном операторе остается доказать включение,
обратное к \eqref{3.49}. Для этого применим к \eqref{3.45} и \eqref{3.46}
интерполяцию с параметром $\psi$. Получим ограниченные операторы
\begin{gather*}
P:X_{\psi}\rightarrow
[\,(X_{0})_{T,Y_{0}},\,(X_{1})_{T,Y_{1}}\,]_{\psi},\\
R:Y_{\psi}\rightarrow [\,(X_{0})_{T,Y_{0}},\,(X_{1})_{T,Y_{1}}\,]_{\psi}.
\end{gather*}
Следовательно, если $u\in (X_{\psi})_{T,Y_{\psi}}$, т. е. $u\in X_{\psi}$, $Tu\in
Y_{\psi}\,$, то ввиду \eqref{3.44}
$$
u=Pu+RTu\in[\,(X_{\,0})_{T,Y_{0}},\,(X_{1})_{T,Y_{1}}\,]_{\psi}.
$$
Тем самым справедливо включение, обратное к \eqref{3.49}.

Теорема \ref{th3.12} доказана.

\subsection[Эллиптическая краевая задача в пространствах Соболева]
{Эллиптическая краевая задача \\ в пространствах Соболева}\label{sec3.3.3}

Для доказательства основного результата нам понадобится классическая теорема о
разрешимости неоднородной краевой задачи \eqref{3.1}, \eqref{3.2} в позитивных
пространствах Соболева (см., например, [\ref{Berezansky65}; \ref{LionsMagenes71},
с.~191; \ref{FunctionalAnalysis72}, с.~169]).


\begin{proposition}\label{prop3.1}
Отображение
\begin{equation}\label{3.50}
u\mapsto(Lu,B_{1}\,u,\ldots,B_{q}\,u),\quad u\in C^{\infty}(\,\overline{\Omega}\,),
\end{equation}
продолжается по непрерывности до ограниченного нетерового оператора
\begin{equation}\label{3.51}
(L,B):H^{s}(\Omega)\rightarrow H^{s-2q}(\Omega)\oplus\mathcal{H}_{s}(\Gamma)
\end{equation}
при любом вещественном $s\geq2q$. Ядро этого оператора совпадает с $N$, а область
значений состоит из всех векторов
$$
(f,g_{1},\ldots,g_{q})\in H^{s-2q}(\Omega)\oplus\mathcal{H}_{s}(\Gamma),
$$
которые удовлетворяют условию
\begin{equation}\label{3.52}
(f,v)_{\Omega}+\sum_{j=1}^{q}\,(g_{j},C^{+}_{j}\,v)_{\Gamma}=0\quad\mbox{для
всех}\;\;v\in N^{+}.
\end{equation}
Индекс оператора \eqref{3.51} равен $\dim N-\dim N^{+}$ и не зависит от~$s$.
\end{proposition}


Здесь и далее обозначено (ср. с формулой \eqref{3.40})
$$
\mathcal{H}_{s}(\Gamma):=\bigoplus_{j=1}^{q}H^{s-m_{j}-1/2}(\Gamma),\quad
s\in\mathbb{R} .
$$

Предложение \ref{prop3.1} распространено на случай произвольного вещественного $s$ в
работах Ж.-Л.~Лионса, Э.~Мадженеса [\ref{LionsMagenes62}, \ref{LionsMagenes63},
\ref{LionsMagenes71}] и Я.~А.~Ройтберга [\ref{Roitberg64}, \ref{Roitberg65},
\ref{Roitberg96}] (см. также изложение результатов Ройтберга в [\ref{Berezansky65}]
(гл. III, \S~6). При этом оператор $(L,B)$ был изучен в пространствах, построенных
различным образом с помощью пространств Соболева соответствующих порядков. Здесь нам
понадобится конструкция, изложенная в монографии Ж.-Л.~Лионса и Э.~Мадженеса
[\ref{LionsMagenes71}] (гл. 2, \S~6, 7), которая в отличие от цитированных работ
Я.~А.~Ройтберга остается в рамках пространств распределений в области $\Omega$. Для
простоты формулировок ограничимся случаем целого $s$ (этого нам будет достаточно).
Общий случай будет рассмотрен ниже в п.~\ref{sec4.4.1}.

Пусть функция $\varrho_{1}\in C^{\infty}(\,\overline{\Omega}\,)$ положительна в
$\Omega$ и равна расстоянию до границы $\Gamma$ в некоторой ее окрестности.

Для произвольного целого $\sigma\geq0$ положим
\begin{equation}\label{3.53}
\Xi^{\sigma}(\Omega):=\{u\in\mathcal{D}'(\Omega):\,\varrho_{1}^{|\mu|}D^{\mu}u\in
L_{2}(\Omega),\;\;|\mu|\leq\sigma\}.
\end{equation}
Здесь $\mu$ --- $n$-мерный мультииндекс. Пространство $\Xi^{\sigma}(\Omega)$
гильбертово относительно скалярного произведения
$$
(u_{1},u_{2})_{\Xi^{\sigma}(\Omega)}:=\sum_{|\mu|\leq\,\sigma}\,
(\varrho_{1}^{|\mu|}D^{\mu}u_{1},\,\varrho_{1}^{|\mu|}D^{\mu}u_{2})_{L_{2}(\Omega)}.
$$
Справедливы непрерывные плотные вложения
\begin{equation}\label{3.54}
H^{\sigma}_{0}(\Omega)\hookrightarrow\Xi^{\sigma}(\Omega)\hookrightarrow
L_{2}(\Omega).
\end{equation}
Здесь $H^{\sigma}_{0}(\Omega)$ --- замыкание множества $C^{\infty}_{0}(\Omega)$ в
топологии пространства $H^{\sigma}(\Omega)$.

Обозначим через $\Xi^{-\sigma}(\Omega)$ гильбертово пространство, двойственное к
$\Xi^{\sigma}(\Omega)$ относительно скалярного произведения в $L_{2}(\Omega)$.
Поскольку [\ref{Triebel80}, с. 414] (теорема 4.8.2 (a)) пространства
$H^{\sigma}_{0}(\Omega)$ и $H^{-\sigma}(\Omega)$ взаимно двойственны относительно
этого же скалярного произведения, то \eqref{3.54} влечет непрерывность плотных
вложений
\begin{equation}\label{3.55}
L_{2}(\Omega)\hookrightarrow\Xi^{-\sigma}(\Omega)\hookrightarrow
H^{-\sigma}(\Omega)\quad\mbox{для целого}\;\;\sigma>0.
\end{equation}
Из правого вложения вытекает, что пространство $\Xi^{-\sigma}(\Omega)$ состоит из
распределений в области $\Omega$.

Для произвольного целого $s<2q$ определим линейное пространство
$$
D^{s}_{L}(\Omega):=\{u\in H^{s}(\Omega):\,Lu\in\Xi^{s-2q}(\Omega)\}.
$$
Введем в нем скалярное произведение графика
$$
(u_{1},u_{2})_{D^{s}_{L}(\Omega)}=(u_{1},u_{2})_{H^{s}(\Omega)}+
(Lu_{1},Lu_{2})_{\Xi^{s-2q}(\Omega)}.
$$
Пространство $D^{s}_{L}(\Omega)$ полно относительно этого скалярного произведения.

Множество $C^{\infty}(\,\overline{\Omega}\,)$ плотно в пространстве
$D^{s}_{L}(\Omega)$. В силу \eqref{3.54} и ограниченности оператора
$L:H^{2q}(\Omega)\rightarrow L_{2}(\Omega)$ выполняются непрерывные плотные вложения
\begin{equation}\label{3.56}
H^{2q}(\Omega)\hookrightarrow D^{s}_{L}(\Omega)\hookrightarrow H^{s}(\Omega)
\quad\mbox{при целом }s<2q.
\end{equation}


\begin{remark}\label{rem3.6}
В цитируемой выше монографии Лионса и Мадженеса [\ref{LionsMagenes71}] через
$H^{s}(\Omega)$ при $s<0$ обозначено пространство, двойственное к
$H^{-s}_{0}(\Omega)$ относительно скалярного произведения в $L_{2}(\Omega)$. Это
двойственное пространство совпадает с используемым нами пространством
$H^{s}(\Omega)$ для всех неполуцелых $s<0$ [\ref{Triebel80}, с. 414] (теорема 4.8.2
(a)) и не совпадает с $H^{s}(\Omega)$ для полуцелых $s<0$, см. также
п.~\ref{sec4.4.1}.
\end{remark}


Следующий результат доказан Лионсом и Мадженесом [\ref{LionsMagenes71}, с. 206, 216]
(теоремы 6.7, 7.4).


\begin{proposition}\label{prop3.2}
Отображение \eqref{3.50} продолжается по непрерывности до ограниченного нетерового
оператора
\begin{equation}\label{3.57}
(L,B):D^{s}_{L}(\Omega)\rightarrow\Xi^{s-2q}(\Omega)\oplus\mathcal{H}_{s}(\Gamma)
\end{equation}
при любом целом $s<2q$. Ядро этого оператора совпадает с $N$, а область значений
состоит из всех векторов
$$
(f,g_{1},\ldots,g_{q})\in\Xi^{s-2q}(\Omega)\oplus\mathcal{H}_{s}(\Gamma),
$$
которые удовлетворяют условию \eqref{3.52}. Индекс оператора \eqref{3.51} равен
$\dim N-\dim N^{+}$ и не зависит от $s$.
\end{proposition}


Нам также понадобится одно утверждение о гомеоморфизме, который осуществляет
оператор, отвечающий некоторой однородной краевой задаче Дирихле. Зафиксируем
произвольное целое число $r\geq1$ и возьмем $r$-ую итерацию $L^{r}$ выражения $L$.
Пусть $L^{r+}$ --- выражение, формально сопряженное к $L^{r}$. Рассмотрим линейное
дифференциальное выражение $L^{r}L^{r+}+{1}$ порядка $4qr$ с коэффициентами из
$C^{\infty}(\,\overline{\Omega}\,)$. Для произвольного целого $\sigma\geq2qr$
положим
$$
H^{\sigma}_{D}(\Omega):=\{u\in H^{\sigma}(\Omega):
\,\gamma_{j}u=0\;\;\mbox{на}\;\;\Gamma,\;\;j=0,\ldots,2qr-1\}.
$$
Здесь и далее $\gamma_{j}u:=(\partial^{j}u/\partial^{j}_{\nu})\upharpoonright\Gamma$
--- оператор следа на $\Gamma$ нормальной производной порядка $j$. При этом след
понимается в смысле теоремы \ref{th3.5}, согласно которой мы имеем ограниченный
оператор $\gamma_{j}:H^{\sigma}(\Omega)\rightarrow H^{\sigma-j-1/2}(\Gamma)$.
Поэтому линейное пространство $H^{\sigma}_{D}(\Omega)$ полно относительно скалярного
произведения в пространстве $H^{\sigma}(\Omega)$.


\begin{lemma}\label{lem3.1}
Пусть целое $r\geq1$. Сужение отображения $u\mapsto L^{r}L^{r+}u+u$, где $u\in
\mathcal{D}'(\Omega)$, задает гомеоморфизм
\begin{equation}\label{3.58}
L^{r}L^{r+}+1\,:\;H^{\sigma}_{D}(\Omega)\leftrightarrow H^{\sigma-4qr}(\Omega)
\end{equation}
для произвольного целого $\sigma\geq2qr$.
\end{lemma}


\textbf{Доказательство.} Дифференциальное выражение $L^{r}L^{r+}+{1}$ правильно
эллиптическое в $\overline{\Omega}$, поскольку таковым является выражение $L$.
Рассмотрим неоднородную краевую задачу Дирихле
$$
L^{r}L^{r+}\,u+u=f\;\;\text{в}\;\;\Omega,\quad
\gamma_{j}\,u=g_{j}\;\;\text{на}\;\;\Gamma\;\;\mbox{при}\;\;j=0,\ldots,2qr-1.
$$
Она регулярная эллиптическая и, как установлено в [\ref{LionsMagenes71}, с. 223,
227], оператор этой задачи является ограниченным и нетеровым с нулевым индексом в
паре пространств
\begin{gather}\notag
(L^{r}L^{r+}+{1};\gamma_{0},\ldots,\gamma_{2qr-1}): \\ H^{\sigma}(\Omega)\rightarrow
H^{\sigma-4qr}(\Omega)\oplus\bigoplus_{j=0}^{2qr-1}H^{\sigma-j-1/2}(\Gamma)\label{3.59}
\end{gather}
для любого целого $\sigma\geq2qr$. Ядро $N_{D}$ оператора \eqref{3.59} лежит в
$C^{\infty}(\,\overline{\Omega})$. С помощью интегрирования по частям нетрудно
вывести, что оно является тривиальным: $u\in N_{D}\Rightarrow$
$$
(u,u)_{\Omega}=-(L^{r}L^{r+}u,u)_{\Omega}
=-(L^{r+}u,L^{r+}u)_{\Omega}\leq0\;\;\Rightarrow\;\;u=0.
$$
Заметим, что при перебрасывании дифференциального выражения $L^{r}$ порядка $2qr$ с
помощью интегрирования по частям появятся выражения вида
$(\,\cdot\,,\gamma_{j}u)_{\Gamma}$, где $j=0,\ldots,2qr-1$; а они равны нулю для
$u\in N_{D}$. Следовательно, оператор \eqref{3.59}~--- гомеоморфизм. Поэтому его
сужение на подпространство $H^{\sigma}_{D}(\Omega)$ определяет гомеоморфизм
\eqref{3.58}. Лемма \ref{lem3.1} доказана.


\begin{remark}\label{}
Утверждение подобное лемме~\ref{lem3.1} имеется в монографии Х.~Трибеля
[\ref{Triebel80}, с. 506].
\end{remark}

\subsection{Доказательство основного результата}\label{sec3.3.4}

Докажем основной результат п.~\ref{sec3.1} --- теорему~\ref{th3.11}.

\medskip

\textbf{Доказательство теоремы~\ref{th3.11}.} Пусть $s\in\mathbb{R}$ и
$\varphi\in\mathcal{M}$. Выберем целое число $r\geq1$ такое, что
\begin{equation}\label{3.60}
2q(1-r)<s<2qr,
\end{equation}
и воспользуемся предложением \ref{prop3.2} для целого $s=2q(1-r)\leq0$ и
предложением \ref{prop3.1} для $s=2qr\geq2q$. Мы получим, что отображение
\eqref{3.50} продолжается по непрерывности до ограниченных нетеровых операторов
\begin{gather}\label{3.61}
(L,B):D^{2q(1-r)}_{L}(\Omega)\rightarrow\Xi^{-2qr}(\Omega)
\oplus\mathcal{H}_{2q(1-r)}(\Gamma),\\
(L,B):H^{2qr}(\Omega)\rightarrow H^{2q(r-1)}(\Omega)\oplus\mathcal{H}_{2qr}(\Gamma),
\label{3.62}
\end{gather}
имеющих общее ядро $N$ и одинаковый индекс.

Отметим, что пары пространств
\begin{equation}\label{3.63}
[D^{\,2q(1-r)}_{L}(\Omega),H^{2qr}(\Omega)]\quad\mbox{и}\quad
[\Xi^{-2qr}(\Omega),H^{2q(r-1)}(\Omega)]
\end{equation}
допустимые. В самом деле, в силу \eqref{3.55} и \eqref{3.56} выполняются непрерывные
плотные вложения
\begin{gather*}
 H^{2qr}(\Omega)\hookrightarrow H^{2q}(\Omega)\hookrightarrow
D^{2q(1-r)}_{L}(\Omega), \\ H^{2q(r-1)}(\Omega)\hookrightarrow
L_{2}(\Omega)\hookrightarrow\Xi^{-2qr}(\Omega).
\end{gather*}
Следовательно, правые пространства пар \eqref{3.63} непрерывно и плотно вкладываются
в соответствующие левые пространства. Отсюда, поскольку правые пространства Соболева
сепарабельные, вытекает, что и левые пространства сепарабельные. Значит, пары
\eqref{3.63} допустимые.

Ввиду \eqref{3.60} положим
\begin{equation}\label{3.64}
\varepsilon=s-2q(1-r)>0,\quad\delta=2qr-s>0.
\end{equation}
Пусть $\psi$ --- интерполяционный параметр из теорем \ref{th3.2} и \ref{th2.21},
соответствующий выбранным нами параметрам $\varphi$, $\varepsilon$ и $\delta$.
Применив интерполяцию с параметром $\psi$ к пространствам, в которых действуют
ограниченные нетеровы операторы \eqref{3.61} и \eqref{3.62}, получим согласно
теоремам \ref{th2.7} и \ref{th2.5} ограниченный нетеров оператор
\begin{gather}\notag
(L,B):\,[D^{2q(1-r)}_{L}(\Omega),H^{2qr}(\Omega)]_{\psi}\rightarrow\\
[\Xi^{-2qr}(\Omega),H^{2q(r-1)}(\Omega)]_{\psi}\oplus
[\mathcal{H}_{2q(1-r)}(\Gamma),\mathcal{H}_{2qr}(\Gamma)]_{\psi}. \label{3.65}
\end{gather}

Здесь в силу теорем \ref{th2.21}, \ref{th2.5} и ввиду \eqref{3.64} мы имеем:
\begin{gather*}
[\mathcal{H}_{2q(1-r)}(\Gamma),\mathcal{H}_{2qr}(\Gamma)]_{\psi}= \\
\biggl[\bigoplus_{j=1}^{q}H^{2q(1-r)-m_{j}-1/2}(\Gamma),
\bigoplus_{j=1}^{q}H^{2qr-m_{j}-1/2}(\Gamma)\biggr]_{\psi}=\\
\bigoplus_{j=1}^{q}[H^{2q(1-r)-m_{j}-1/2}(\Gamma),H^{2qr-m_{j}-1/2}(\Gamma)]_{\psi}=\\
\bigoplus_{j=1}^{q}
[H^{s-m_{j}-1/2-\varepsilon}(\Gamma),H^{s-m_{j}-1/2+\delta}(\Gamma)]_{\psi}= \\
\bigoplus_{j=1}^{q}H^{s-m_{j}-1/2,\varphi}(\Gamma)=\mathcal{H}_{s,\varphi}(\Gamma)
\end{gather*}
с точностью до эквивалентности норм. Обозначив
\begin{equation}\label{3.66}
Z(\Omega):=[\Xi^{-2qr}(\Omega),H^{2q(r-1)}(\Omega)]_{\psi},
\end{equation}
получим, что оператор \eqref{3.65} ограничен и нетеров в паре пространств
\begin{equation}\label{3.67}
(L,B):\,[D^{2q(1-r)}_{L}(\Omega),H^{2qr}(\Omega)]_{\psi}\rightarrow
Z(\Omega)\oplus\mathcal{H}_{s,\varphi}(\Gamma).
\end{equation}
Этот оператор имеет то же ядро $N$ и тот же индекс, что и операторы \eqref{3.61},
\eqref{3.62}.

Опишем при помощи $Z(\Omega)$ область определения оператора \eqref{3.67}. Это будет
сделано на основании теоремы \ref{th3.12}, в условии которой мы полагаем
\begin{gather*}
X_{0}=H^{2q(1-r)}(\Omega),\quad Y_{0}=\Xi^{-2qr}(\Omega),\quad
Z_{0}=E=H^{-2qr}(\Omega),\\
X_{1}=H^{2qr}(\Omega),\quad Y_{1}=Z_{1}=H^{2q(r-1)}(\Omega),\quad T=L.
\end{gather*}
Из того, что вторая пара в формуле \eqref{3.63} допустимая, а также из правого
вложения \eqref{3.55} вытекает, что условия (i), (ii) и (iii) теоремы \ref{th3.12}
выполняются. Условие (iv) этой теоремы также выполняется, поскольку ограничен
оператор $L:H^{\sigma}(\Omega)\rightarrow H^{\sigma-2q}(\Omega)$ для произвольного
$\sigma\in\mathbb{R}$. Нам, кроме того, нужны линейные отображения $R$ и $S$,
удовлетворяющие условиям (v), (vi) и (vii). Определим их следующим образом.
Воспользуемся леммой \ref{lem3.1} и рассмотрим отображение $(L^{r}L^{r+}+1)^{-1}$,
обратное к гомеоморфизму \eqref{3.58}. Получим линейный ограниченный оператор
\begin{equation}\label{3.68}
(L^{r}L^{r+}+1)^{-1}:H^{\sigma-4qr}(\Omega)\rightarrow H^{\sigma}(\Omega)
\end{equation}
для произвольного целого $\sigma\geq2qr$. Положим
$$
R=L^{r-1}L^{r+}(L^{r}L^{r+}+1)^{-1}\quad\mbox{и}\quad S=-(L^{r}L^{r+}+1)^{-1}.
$$
В силу \eqref{3.68} при $\sigma=2qr$ и при $\sigma=2q(3r-1)$ получаем ограниченные
операторы
\begin{gather*}
R:Z_{0}=H^{-2qr}(\Omega)\rightarrow H^{2qr-2q(2r-1)}(\Omega)=X_{0},\\
R:Z_{1}=H^{2q(r-1)}(\Omega)\rightarrow H^{2q(3r-1)-2q(2r-1)}(\Omega)=X_{1},\\
S:Z_{0}=H^{-2qr}(\Omega)\rightarrow H^{2qr}(\Omega)\hookrightarrow
H^{0}(\Omega) \hookrightarrow\Xi^{-2qr}(\Omega)=Y_{0},\\
S:Z_{1}=H^{2q(r-1)}(\Omega)\rightarrow H^{2q(3r-1)}(\Omega)\hookrightarrow
H^{2qr}(\Omega)=X_{1}.
\end{gather*}
Кроме того, на множестве $E=H^{-2qr}(\Omega)$ справедливы равенства
\begin{gather*}
TR=LL^{r-1}L^{r+}(L^{r}L^{r+}+1)^{-1}= \\ (L^{r}L^{r+}+1-1)(L^{r}L^{r+}+1)^{-1}=1-S.
\end{gather*}
Таким образом, все условия теоремы \ref{th3.12} выполняются. Согласно этой теореме
для интерполяционного параметра $\psi$ имеем равенство пространств с точностью до
эквивалентности норм
\begin{equation}\label{3.69}
[(X_{0})_{L,Y_{0}},(X_{1})_{L,Y_{1}}]_{\psi}=(X_{\psi})_{L,Y_{\psi}}.
\end{equation}
Здесь
$$
(X_{0})_{L,Y_{0}}=\{u\in
H^{2q(1-r)}(\Omega):\,Lu\in\Xi^{-2qr}(\Omega)\}=D_{L}^{2q(1-r)}(\Omega),
$$
причем нормы в крайних пространствах равны. Далее, в силу ограниченности оператора
$L:H^{2qr}(\Omega)\rightarrow H^{2q(r-1)}(\Omega)$ верно равенство
$$
(X_{1})_{L,Y_{1}}=\{u\in H^{2qr}(\Omega):\,Lu\in
H^{2q(r-1)}(\Omega)\}=H^{2qr}(\Omega)
$$
с эквивалентностью норм в крайних пространствах. Кроме того, согласно теореме
\ref{th3.2} и ввиду \eqref{3.64} имеем
\begin{gather*}
X_{\psi}=[H^{2q(1-r)}(\Omega),H^{2qr}(\Omega)]_{\psi}= \\
[H^{s-\varepsilon}(\Omega),H^{s+\delta}(\Omega)]_{\psi}=H^{s,\varphi}(\Omega).
\end{gather*}
Таким образом, соотношение \eqref{3.69} принимает вид
\begin{equation}\label{3.70}
[D_{L}^{2q(1-r)}(\Omega),H^{2qr}(\Omega)]_{\psi}=\{u\in
H^{s,\varphi}(\Omega):\,Lu\in Z(\Omega)\},
\end{equation}
причем в последнем пространстве введено скалярное произведение графика, относительно
которого это пространство полно (мы также воспользовались обозначением \eqref{3.66},
согласно которому $Y_{\psi}=Z(\Omega)$).

Подставив теперь равенство \eqref{3.70} в \eqref{3.67}, получим ограниченный
оператор
\begin{equation}\label{3.71}
(L,B):\,\{u\in H^{s,\varphi}(\Omega):\,Lu\in Z(\Omega)\}\rightarrow
Z(\Omega)\oplus\mathcal{H}_{s,\varphi}(\Gamma).
\end{equation}
По уже доказанному, он нетеров и имеет ядро $N$. Кроме того, поскольку оператор
\eqref{3.71} получен в результате интерполяции, примененной к нетеровым операторам
\eqref{3.61} и \eqref{3.62}, то на основании теоремы \ref{th2.7} и предложения
\ref{prop3.2} область значений оператора \eqref{3.71} равна
$$
Z(\Omega)\oplus\mathcal{H}_{s,\varphi}(\Gamma)\cap (L,B)(D^{2q(1-r)}_{L}(\Omega)),
$$
т. е. состоит из всех векторов
$$
(f,g_{1},\ldots,g_{q})\in Z(\Omega)\oplus\mathcal{H}_{s,\varphi}(\Gamma),
$$
которые удовлетворяют условию \eqref{3.52}. Сужение оператора \eqref{3.71} на
подпространство
$$
K_{L}^{s,\varphi}(\Omega)=\{u\in
H^{s,\varphi}(\Omega):\;L\,u=0\;\;\mbox{в}\;\;\Omega\,\}
$$
определяет ограниченный оператор
\begin{equation}\label{3.72}
B:\,K_{L}^{s,\varphi}(\Omega)\rightarrow\mathcal{H}_{s,\varphi}(\Gamma).
\end{equation}
Его ядро равно $N\cap\,K_{L}^{s,\varphi}(\Omega)=N$ и, значит, конечномерно, а
область значений совпадает с подпространством \eqref{3.41} и, следовательно,
замкнута и имеет конечную коразмерность, равную размерности пространства $G$,
определенного по формуле \eqref{3.42}. Таким образом, оператор \eqref{3.72} нетеров
с ядром $N$, областью значений \eqref{3.41} и конечным индексом $\dim N-\dim G$, не
зависящим от $s$, $\varphi$.

Остается показать, что множество $K_{L}^{\infty}(\Omega)$ плотно в
$K_{L}^{s,\varphi}(\Omega)$, и что оператор \eqref{3.72} является продолжением по
непрерывности отображения \eqref{3.38}. В~связи с этим заметим следующее: поскольку
оператор \eqref{3.71} является продолжением отображения \eqref{3.50}, то в силу
определения оператора \eqref{3.72} последний является продолжением отображения
\eqref{3.38}. Поэтому для того, чтобы закончить доказательство, надо установить
плотность множества $K_{L}^{\infty}(\Omega)$ в $K_{L}^{s,\varphi}(\Omega)$. Сделаем
это при помощи гомеоморфизма
\begin{equation}\label{3.73}
B:\,K_{L}^{s,\varphi}(\Omega)/\,N\leftrightarrow\,\mathcal{R}_{s,\varphi}(\Gamma),
\end{equation}
который порожден нетеровым оператором \eqref{3.72}. Здесь через
$\mathcal{R}_{s,\varphi}(\Gamma)$ обозначена область значений \eqref{3.41} оператора
\eqref{3.72}. Рассмотрим гомеоморфизм $B^{-1}$, обратный к \eqref{3.73}. Он каждому
вектору $g=(g_{1},\ldots,g_{q})\in\mathcal{R}_{s,\varphi}(\Gamma)$ ставит в
соответствие класс смежности
$$
B^{-1}g=[\,u\,]=\{u+w:\,w\in N\}
$$
элемента $u\in K_{L}^{s,\varphi}(\Omega)$, такого, что $Bu=g$.

Покажем предварительно, что отображение \eqref{3.73} обладает следующим свойством
повышения гладкости:
\begin{gather}\notag
g\in\mathcal{R}_{s,\varphi}(\Gamma)\cap(C^{\infty}(\Gamma))^{q}\,\Rightarrow \\
B^{-1}g=[\,u\,]\;\;\mbox{для некоторого}\;\;u\in K_{L}^{\infty}(\Omega).\label{3.74}
\end{gather}
Пусть
$$
g=(g_{1},\ldots,g_{q})\in\mathcal{R}_{s,\varphi}(\Gamma)\cap(C^{\infty}(\Gamma))^{q}.
$$
Поскольку $\mathcal{R}_{s,\varphi}(\Gamma)$ --- это множество \eqref{3.41}, то в
силу предложения \ref{th3.1} эллиптическая краевая задача \eqref{3.36} имеет решение
$u\in H^{2q}(\Omega)$. Правые части этой задачи бесконечно гладкие; следовательно
[\ref{LionsMagenes71}, с. 191], выполняется включение $u\in
C^{\infty}(\,\overline{\Omega}\,)$. Таким образом, $u\in K_{L}^{\infty}(\Omega)$ и
$Bu=g$ для оператора \eqref{3.72}, что и доказывает \eqref{3.74}.

Теперь нетрудно установить упомянутую плотность. Возьмем произвольное распределение
$u\in K_{L}^{s,\varphi}(\Omega)$ и по нему образуем вектор
\begin{equation}\label{3.75}
g=Bu\in\mathcal{R}_{s,\varphi}(\Gamma)\subset\mathcal{H}_{s,\varphi}(\Gamma).
\end{equation}
Поскольку множество $C^{\infty}(\Gamma)$ плотно в $H^{\sigma,\varphi}(\Gamma)$ при
каждом $\sigma\in\mathbb{R}$, то для $g$ существует последовательность векторов
$g^{(k)}$ такая, что
\begin{equation}\label{3.76}
g^{(k)}\in(C^{\infty}(\Gamma))^{q}\quad\mbox{и}\quad g^{(k)}\rightarrow g
\;\;\mbox{в}\;\;\mathcal{H}_{s,\varphi}(\Gamma) \;\;\mbox{при}\;\;k\rightarrow
\infty.
\end{equation}
Заметим далее следующее: так как $\mathcal{R}_{s,\varphi}(\Gamma)$ и $G$ ---
замкнутые подпространства в $\mathcal{H}_{s,\varphi}(\Gamma)$, удовлетворяющие
условиям $\mathcal{R}_{s,\varphi}(\Gamma)\cap G=\{0\}$ и
$\mathrm{codim}\,\mathcal{R}_{s,\varphi}(\Gamma)=\dim G$, то
$\mathcal{H}_{s,\varphi}(\Gamma)$ является прямой суммой этих подпространств.
Воспользовавшись этой суммой, запишем: $g=g+0$ и $g^{(k)}=h^{(k)}+\omega^{(k)}$, где
$h^{(k)}\in\mathcal{R}_{s,\varphi}(\Gamma)$ и $\omega^{(k)}\in G$. Отсюда и из
\eqref{3.76} получаем следующие два утверждения:
\begin{gather*}
h^{(k)}=g^{(k)}-\omega^{(k)}\in
\mathcal{R}_{s,\varphi}(\Gamma)\cap(C^{\infty}(\Gamma))^{q},\\
h^{(k)}\rightarrow g\;\;\mbox{в}\;\;\mathcal{R}_{s,\varphi}(\Gamma)\;\;\mbox{(т. е.
в $\mathcal{H}_{s,\varphi}(\Gamma)\,)$}\;\;\mbox{при}\;\;k\rightarrow\infty.
\end{gather*}
Первое в силу \eqref{3.74} влечет $B^{-1}h^{(k)}=[\,u_{k}\,]$ для некоторого
$u_{k}\in K_{L}^{\infty}(\Omega)$. Из второго в виду \eqref{3.73} и \eqref{3.75}
вытекает
$$
[\,u_{k}\,]=B^{-1}h^{(k)}\rightarrow B^{-1}g=[\,u\,],
$$
т. е.
$$
[\,u_{k}-u\,]\rightarrow0\;\;\mbox{в}\;\;K_{L}^{s,\varphi}(\Omega)/\,N
\;\;\mbox{при}\;\;k\rightarrow\infty.
$$
Последнее означает, что
$$
u_{k}-u+w_{k}\rightarrow0\;\;\mbox{в}\;\;K_{L}^{s,\varphi}(\Omega)\;\;
\mbox{при}\;\;k\rightarrow\infty
$$
для некоторой последовательности функций $w_{k}\in N\subset K_{L}^{\infty}(\Omega)$.
Таким образом, произвольное распределение $u\in K_{L}^{s,\varphi}(\Omega)$
аппроксимировано в пространстве $K_{L}^{s,\varphi}(\Omega)$ последовательностью
функций $u_{k}+w_{k}\in K_{L}^{\infty}(\Omega)$. Значит, множество
$K_{L}^{\infty}(\Omega)$ плотно в $K_{L}^{s,\varphi}(\Omega)$.

Теорема \ref{th3.11} доказана.

\subsection[Свойства решений однородного эллиптического уравнения]
{Свойства решений однородного \\ эллиптического уравнения}\label{sec3.3.5}

Из теоремы \ref{th3.11} следует, что в случае тривиальных ядра $N$ и дефектного
подпространства $G$ оператор \eqref{3.40}, соответствующий задаче \eqref{3.36},
является гомеоморфизмом. В общем случае этот оператор определяет гомеоморфизм
\begin{equation}\label{3.77}
B:\,K_{L}^{s,\varphi}(\Omega)/\,N\leftrightarrow\,\mathcal{R}_{s,\varphi}(\Gamma)
\;\;\mbox{для любых}\;\;s\in\mathbb{R},\;\;\varphi\in\mathcal{M}.
\end{equation}
Здесь, напомним, $\mathcal{R}_{s,\varphi}(\Gamma)$ --- подпространство \eqref{3.41}.
(Заметим, что оператор, обратный к \eqref{3.77}, ограничен по теореме Банаха об
обратном операторе.) Набор гомеоморфизмов \eqref{3.77} дает решение задачи
\eqref{3.36} для произвольных распределений $g_{1},\ldots,g_{q}\in
\mathcal{D}'(\Gamma)$, удовлетворяющих условию
$$
(g_{1},\,C^{+}_{1}\,v)_{\Gamma}+\ldots+
(g_{q},\,C^{+}_{q}\,v)_{\Gamma}=0\quad\mbox{для любого}\quad v\in N^{+}.
$$
При этом выполняется следующая априорная оценка решения $u$.


\begin{theorem}\label{th3.13}
Пусть заданы параметры $s\in\mathbb{R}$, $\varphi\in\mathcal{M}$ и $\varepsilon>0$.
Существует число $c=c(s,\varphi,\varepsilon)>0$ такое, что для любого $u\in
K_{L}^{s,\varphi}(\Omega)$ верна оценка
\begin{equation}\label{3.78}
\|u\|_{H^{s,\varphi}(\Omega)}\leq
c\,\bigl(\,\|Bu\|_{\mathcal{H}_{s,\varphi}(\Gamma)}+
\|u\|_{H^{s-\varepsilon}(\Omega)}\,\bigr).
\end{equation}
\end{theorem}


\textbf{Доказательство.} Для любого распределения $u\in K_{L}^{s,\varphi}(\Omega)$ в
силу гомеоморфизма \eqref{3.77} имеем
\begin{equation}\label{3.79}
\inf\bigl\{\,\|u+w\|_{H^{s,\varphi}(\Omega)}:\,w\in N\,\bigr\} \leq
c_{0}\,\|Bu\|_{\mathcal{H}_{s,\varphi}(\Gamma)}.
\end{equation}
Здесь $c_{0}$ --- норма оператора, обратного к \eqref{3.77}. Далее, поскольку $N$
--- конечномерное подпространство в $H^{s,\varphi}(\Omega)$ и в
$H^{s-\varepsilon}(\Omega)$, нормы в этих двух пространствах эквивалентны на $N$. В
частности, для любого $w\in N$
$$
\|w\|_{H^{s,\varphi}(\Omega)}\leq c_{1}\,\|w\|_{H^{s-\varepsilon}(\Omega)}
$$
с постоянной $c_{1}$, не зависящей от $u$ и $w$. Кроме того,
\begin{gather*}
\|w\|_{H^{s-\varepsilon}(\Omega)}\leq \|u+w\|_{H^{s-\varepsilon}(\Omega)}+
\|u\|_{H^{s-\varepsilon}(\Omega)}\leq\\
c_{2}\,\|u+w\|_{H^{s,\varphi}(\Omega)}+\|u\|_{H^{s-\varepsilon}(\Omega)}.
\end{gather*}
Здесь $c_{2}$ --- норма оператора вложения $H^{s,\varphi}(\Omega)\hookrightarrow
H^{s-\varepsilon}(\Omega)$. Следовательно,
\begin{gather*}
\|u\|_{H^{s,\varphi}(\Omega)} \leq \|u+w\|_{H^{s,\varphi}(\Omega)}+
\|w\|_{H^{s,\varphi}(\Omega)}\leq\\
\|u+w\|_{H^{s,\varphi}(\Omega)}+c_{1}\,\|w\|_{H^{s-\varepsilon}(\Omega)}\leq\\
\leq(1+c_{1}c_{2})\,\|u+w\|_{H^{s,\varphi}(\Omega)}+
c_{1}\,\|u\|_{H^{s-\varepsilon}(\Omega)}.
\end{gather*}
Перейдем теперь к инфимуму по $w\in N$ и воспользуемся неравенством \eqref{3.79}.
Получим
$$
\|u\|_{H^{s,\varphi}(\Omega)}\leq (1+c_{1}c_{2})\,c_{0}\,
\|Bu\|_{\mathcal{H}_{s,\varphi}(\Gamma)}+ c_{1}\,\|u\|_{H^{s-\varepsilon}(\Omega)},
$$
т. е. оценку \eqref{3.78}, если положить $c=\max\{(1+c_{1}c_{2})c_{0},\,c_{1}\}$.
Теорема \ref{th3.13} доказана.

\medskip

Если в неравенстве \eqref{3.78} правая часть конечна, то конечна и левая часть.


\begin{theorem}\label{th3.14}
Пусть $s\in\mathbb{R}$, $\varphi\in\mathcal{M}$ и $\varepsilon>0$. Предположим, что
распределение $u\in H^{s-\varepsilon}(\Omega)$ является решением задачи
$\eqref{3.36}$, где
\begin{equation}\label{3.80}
B_{j}u=g_{j}\in H^{s-m_{j}-1/2,\,\varphi}(\Gamma)\;\;\mbox{для
каждого}\;\;j=1,\ldots,q.
\end{equation}
Тогда $u\in H^{s,\varphi}(\Omega)$.
\end{theorem}


\textbf{Доказательство.} По условию, $u\in K_{L}^{s-\varepsilon,\,1}(\Omega)$, и
$Bu=g$, где вектор $g=(g_{1},\ldots,g_{q})$. Следовательно, в силу \eqref{3.80} и
теоремы \ref{th3.11} (описание области значений оператора \eqref{3.40}) имеем
$$
g\in B(K_{L}^{s-\varepsilon,\,1}(\Omega))\cap\mathcal{H}_{s,\varphi}(\Gamma)=
B(K_{L}^{s,\varphi}(\Omega)).
$$
Поэтому существует распределение $u_{0}\in K_{L}^{s,\varphi}(\Omega)$ такое, что
$Bu_{0}=g$. Отсюда с учетом теоремы \ref{th3.11} (описание ядра оператора
\eqref{3.40}) последовательно получаем: $B(u-u_{0})=0$,
$$
w=u-u_{0}\in N\subset C^{\infty}(\,\overline{\Omega}\,),\quad u=u_{0}+w\in
H^{s,\varphi}(\Omega).
$$
Теорема \ref{th3.14} доказана.

\medskip

Теорема \ref{th3.14} --- это утверждение о повышении гладкости решения $u$ задачи
\eqref{3.36} вплоть до границы $\Gamma$. При этом, как видим, уточненная гладкость
$\varphi$ правых частей задачи наследуется решением. Отметим (см., например,
[\ref{Hermander65}, с. 237], теорема 7.4.1), что всякое решение однородного
эллиптического уравнения $Lu=0$ в области $\Omega$ обладает свойством $u\in
C^{\infty}(\Omega)$. Поэтому в теореме \ref{th3.14} существенно, что гладкость
решения $u$ повышается вплоть до границы области $\Omega$.


\begin{corollary}\label{cor3.2}
Пусть $\sigma\in\mathbb{R}$. Предположим, что распределение $u\in
H^{\sigma}(\Omega)$ является решением задачи $\eqref{3.36}$, в которой
\begin{equation}\label{3.81}
g_{j}\in H^{m-m_{j}+(n-1)/2,\,\varphi}(\Gamma)\;\;\mbox{для
каждого}\;\;j=1,\ldots,q,
\end{equation}
где $m:=\max\{m_{1},\ldots,m_{q}\}$, а функция $\varphi\in\mathcal{M}$ удовлетворяет
условию $\eqref{2.37}$. Тогда $u\in C^{m}(\,\overline{\Omega}\,)$, и поскольку также
$u\in C^{\infty}(\Omega)$, то $u$ является классическим решением задачи
$\eqref{3.36}$.
\end{corollary}


\textbf{Доказательство.} Условие \eqref{3.81} совпадает с \eqref{3.80}, если
положить в нем $s=m+n/2$. Поэтому в силу теорем \ref{th3.14} и \ref{th3.4}
$$
u\in H^{m+n/2\,,\varphi}(\Omega)\hookrightarrow C^{m}(\,\overline{\Omega}\,),
$$
что и требовалось доказать.

\medskip

Отметим, что для классического решения $u$ левые части задачи \eqref{3.36}
вычисляются с помощью классических производных, при этом $B_{j}u\in C(\Gamma)$.






\section[Эллиптическая задача с однородными краевыми условиями]
{Эллиптическая задача с \\ однородными краевыми условиями}\label{sec3.4}

\markright{\emph \ref{sec3.4}. Задача с однородными краевыми условиями}

Рассмотрим регулярную эллиптическую краевую задачу \eqref{3.1}, \eqref{3.2} в
случае, когда краевые условия \eqref{3.2} однородные:
\begin{gather}\label{3.82}
L\,u=f\;\;\text{в}\;\;\Omega,\\
B_{j}\,u=0\;\;\text{на}\;\;\Gamma\;\;\text{при}\;\; j=1,\ldots,q.\label{3.83}
\end{gather}

Изучим свойства отображения $u\mapsto Lu$, где $u$ удовлетворяет равенствам
\eqref{3.83}, в шкале пространств Хермандера.

\subsection[Теорема о гомеоморфизмах для эллиптического оператора]
{Теорема о гомеоморфизмах \\ для эллиптического оператора}\label{sec3.4.1}

Введем необходимые нам пространства распределений, удовлетворяющих однородным
краевым условиям. Пусть $s\in\mathbb{R}$ и $\varphi\in\mathcal{M}$. Для краткости
в~п.~\ref{sec3.4} будем обозначать гильбертово пространство
$H^{s,\varphi,(0)}(\Omega)$ через $H^{s,\varphi}$.

Обозначим через $(\mathrm{b.c.})$ однородные краевые условия \eqref{3.83}. Положим
$$
C^{\infty}(\mathrm{b.c.}):=\{u\in C^{\infty}(\,\overline{\Omega}\,):
\;B_{j}u=0\;\;\mbox{на}\;\;\Gamma,\;\;j=1,\ldots,q\}.
$$
Обозначим через $H^{s,\varphi}(\mathrm{b.c.})$ замыкание множества
$C^{\infty}(\mathrm{b.c.})$ в топологии гильбертова пространства $H^{s,\varphi}$.

Наряду с задачей \eqref{3.82}, \eqref{3.83} рассмотрим формально сопряженную задачу
с однородными краевыми условиями:
\begin{gather}\label{3.84}
L^{+}\,v=g\;\;\text{в}\;\;\Omega,\\
B^{+}_{j}\,v=0\;\;\text{на}\;\;\Gamma\;\;\text{при}\;\;j=1,\ldots,q. \label{3.85}
\end{gather}
Обозначим через $(\mathrm{b.c.})^{+}$ однородные краевые условия \eqref{3.85} и
положим
$$
C^{\infty}(\mathrm{b.c.})^{+}:=\{v\in
C^{\infty}(\,\overline{\Omega}\,):\;B_{j}^{+}v=0\;\;\mbox{на}\;\;\Gamma,\;\;j=1,\ldots,q\}.
$$
Обозначим через $H^{s,\varphi}(\mathrm{b.c.})^{+}$ замыкание множества
$C^{\infty}(\mathrm{b.c.})^{+}$ в топологии гильбертова пространства
$H^{s,\varphi}$.

Линейные пространства $H^{s,\varphi}(\mathrm{b.c.})$ и
$H^{s,\varphi}(\mathrm{b.c.})^{+}$ гильбертовы относительно скалярного произведения
в пространстве $H^{s,\varphi}$.

В силу равенства \eqref{3.8} множество $C^{\infty}(\mathrm{b.c.})^{+}$, а значит, и
пространство $H^{s,\varphi}(\mathrm{b.c.})^{+}$ не зависят от выбора системы
граничных выражений $\{B^{+}_{1},\ldots,B^{+}_{q}\}$, сопряженной к системе
$\{B_{1},\ldots,B_{q}\}$ относительно дифференциального выражения~$L$.

Согласно формуле Грина \eqref{3.6}
\begin{gather}\label{3.86}
(Lu,v)_{\Omega}=(u,L^{+}v)_{\Omega} \\ \mbox{для всех}\;\;u\in
C^{\infty}(\mathrm{b.c.}),\;\;v\in C^{\infty}(\mathrm{b.c.})^{+}. \notag
\end{gather}
Поэтому образ $Lu$ произвольной функции $u\in C^{\infty}(\mathrm{b.c.})$ естественно
интерпретировать как антилинейный ограниченный функционал $(Lu,\cdot)_{\Omega}$ на
пространстве $H^{2q-s,1/\varphi}(\mathrm{b.c.})^{+}$. Кроме того, поскольку
пространства $H^{2q-s,\,1/\varphi}$ и $H^{s-2q,\,\varphi}$ взаимно двойственны
относительно формы $(\cdot,\cdot)_{\Omega}$, относительно нее также взаимно
двойственны подпространство $H^{2q-s,\,1/\varphi}(\mathrm{b.c.})^{+}$ и
факторпространство $H^{s-2q,\,\varphi}/M_{s-2q,\,\varphi}$\,, где
$$
M_{s-2q,\,\varphi}:=\{h\in H^{s-2q,\,\varphi}:\,(h,w)_{\Omega}=0\;\mbox{для
всех}\;w\in C^{\infty}(\mathrm{b.c.})^{+}\}.
$$
Значит, $Lu$ можно интерпретировать как класс смежности
$$
\{Lu+h:\,h\in M_{s-2q,\,\varphi}\},
$$
принадлежащий факторпространству $H^{s-2q,\,\varphi}/M_{s-2q,\,\varphi}$.

Мы рассматриваем линейное отображение $u\mapsto Lu$ как оператор, действующий из
пространства $H^{s,\varphi}(\mathrm{b.c.})$ в пространство
$H^{s-2q,\,\varphi}/M_{s-2q,\,\varphi}$, которое отождествляется с двойственным
пространством $(H^{2q-s,\,1/\varphi}(\mathrm{b.c.})^{+})'$ (оно состоит из
антилинейных функционалов). Для того, чтобы сформулировать теорему о гомеоморфизмах
для такого оператора $L$, нам понадобятся проекторы пространства $H^{s,\varphi}$ на
подпространства, ортогональные соответственно к $N$ и $N^{+}$ относительно
билинейной формы $(\,\cdot\,,\,\cdot\,)_{\Omega}$. Эти проекторы существуют,
поскольку пространства $N$ и $N^{+}$ конечномерны.


\begin{lemma}\label{lem3.2} Пусть $s\in\mathbb{R}$ и $\varphi\in\mathcal{M}$.
Произвольный элемент $u\in H^{s,\varphi}$ допускает единственное представление вида
$u=u_{0}+u_{1}$, где $u_{0}\in N$, а $\,u_{1}\in H^{s,\varphi}$ удовлетворяет
условию $(u_{1},w)_{\Omega}=0$ для любого $w\in N$. При этом отображение $P:
u\mapsto u_{1}$ является проектором пространства $H^{s,\varphi}$ на подпространство
\begin{equation}\label{3.87}
\{u_{1}\in H^{s,\varphi}:\,(u_{1},w)_{\Omega}=0\;\,\mbox{для всех}\;\,w\in N\},
\end{equation}
и образ $Pu$ не зависит от $s,\varphi$. Сужение отображения $P$ на
$H^{s,\varphi}(\mathrm{b.c.})$ является проектором пространства
$H^{s,\varphi}(\mathrm{b.c.})$ на подпространство
\begin{equation}\label{3.88}
\{u_{1}\in H^{s,\varphi}(\mathrm{b.c.}):\,(u_{1},w)_{\Omega}=0\;\,\mbox{для
всех}\;\,w\in N\}.
\end{equation}
Эта лемма сохраняет силу, если в ее формулировке заменить $N$ на $N^{+}$, $P$ на
$P^{+}$ и $(\mathrm{b.c.})$ на $(\mathrm{b.c.})^{+}$ соответственно. При этом
$M_{s,\varphi}$ --- подпространство в $P^{+}(H^{s,\varphi})$.
\end{lemma}


\textbf{Доказательство.} Как отмечалось выше, $N$ --- конечномерное подпространство
в $H^{s,\varphi}$. Ясно, что $\dim N$ совпадает с коразмерностью подпространства
\eqref{3.87}, причем $N$ и \eqref{3.87} имеют тривиальное пересечение. Следовательно
[\ref{FunctionalAnalysis72}, с. 56], пространство $H^{s,\varphi}$ разлагается в
прямую сумму подпространств $N$ и \eqref{3.87} с проектором $P$ на подпространство
\eqref{3.87}, который, очевидно, не зависит от $s$ и $\varphi$. Отсюда, поскольку
$N\subset H^{s,\varphi}(\mathrm{b.c.})$, вытекает, что пространство
$H^{s,\varphi}(\mathrm{b.c.})$ разлагается в прямую сумму подпространств $N$ и
\eqref{3.88} с тем же проектором $P$ на подпространство \eqref{3.88}. Как видим,
существование проектора $P$ следует из двух условий: пространство $N$ конечномерно и
$N\subset H^{s,\varphi}(\mathrm{b.c.})$. Значит, поскольку пространство $N^{+}$
конечномерно и $N^{+}\subset H^{s,\varphi}(\mathrm{b.c.})^{+}$, то лемма остается
верной и для проектора $P^{+}$ при указанных заменах обозначений в ее формулировке.
Наконец, заметим, что последнее предложение леммы вытекает из включения
$N^{+}\subset C^{\infty}(\mathrm{b.c.})^{+}$. Лемма \ref{lem3.2} доказана.

\medskip

Сформулируем основной результат п.~\ref{sec3.4} --- теорему о гомеоморфизмах,
осуществляемых эллиптическим оператором $L$ в уточненной шкале пространств
Хермандера.


\begin{theorem}\label{th3.15} Пусть $s\in\mathbb{R}$ и $\varphi\in\mathcal{M}$,
причем
\begin{equation}\label{3.89}
s\neq j+1/2\;\;\mbox{для каждого}\;\;j=0,\,1,\ldots,2q-1.
\end{equation}
Отображение $u\mapsto Lu$, где $u\in C^{\infty}(\mathrm{b.c.})$, а $Lu$
интерпретируется либо как класс смежности $\{Lu+h:\,h\in M_{s-2q,\,\varphi}\}$, либо
как функционал $(Lu,\,\cdot\,)_{\Omega}$, продолжается по непрерывности до
ограниченного линейного оператора
\begin{equation}\label{3.90}
L:H^{s,\varphi}(\mathrm{b.c.})\rightarrow H^{s-2q,\,\varphi}/M_{s-2q,\,\varphi}=
\bigl(\,H^{2q-s,\,1/\varphi}(\mathrm{b.c.})^{+}\,\bigr)'.
\end{equation}
Сужение оператора $\eqref{3.90}$ на подпространство $\eqref{3.87}$ является
гомеоморфизмом
\begin{equation}\label{3.91}
L:P(H^{s,\varphi}(\mathrm{b.c.}))\leftrightarrow
P^{+}(H^{s-2q,\,\varphi})/M_{s-2q,\,\varphi}.
\end{equation}
\end{theorem}


\begin{remark}\label{rem3.8}
В соболевском случае $\varphi\equiv1$ теорема \ref{th3.15} доказана
Ю.~М.~Березанским, С.~Г.~Крейном и Я.~А.~Ройтбергом
[\ref{BerezanskyKreinRoitberg63}, с. 747; \ref{Berezansky65}, с.~262] для целых $s$.
Для всех вещественных $s$ она доказана в монографии Ройтберга [\ref{Roitberg96},
с.~168] (см. также обзор [\ref{Agranovich97}, с.~85, 86]). При этом для полуцелых
$s\in\{j+1/2:j=0,\,1,\ldots,2q-1\}$ пространства, в которых действует оператор
\eqref{3.90}, определялись с помощью интерполяции.
\end{remark}


Отметим, что, если пространства $N$ и $N^{+}$ тривиальны, то оператор \eqref{3.90}
становится гомеоморфизмом.

Из теоремы \ref{th3.15} следует свойство нетеровости оператора \eqref{3.90} и
априорная оценка решения уравнения $Lu=f$.


\begin{theorem}\label{th3.16} Пусть $s\in\mathbb{R}$ и $\varphi\in\mathcal{M}$,
причем выполняется условие $\eqref{3.89}$. Тогда справедливы следующие утверждения.
\begin{itemize}
\item [$\mathrm{(i)}$] Ограниченный оператор $\eqref{3.90}$ нетеров с ядром $N$ и
областью значений
\begin{equation}\label{3.92}
\left\{[f]\in
H^{s-2q,\,\varphi}/\,M_{s-2q,\,\varphi}:\,(f,w)_{\Omega}=0\;\,\mbox{для
всех}\;\,w\in N^{+}\right\}.
\end{equation}
Индекс этого оператора равен $\dim N-\dim N^{+}$ и не зависит от $s,\varphi$. Здесь
$[f]=\{f+h:\,h\in M_{s-2q,\,\varphi}\}$ --- класс смежности элемента $f\in
H^{s-2q,\,\varphi}$.
\item [$\mathrm{(ii)}$] Верна следующая априорная оценка решения $u\in
H^{s,\varphi}(\mathrm{b.c.})$ уравнения $Lu=[f]$: для произвольного числа
$\varepsilon>0$ существует такое число $c=c(s,\varphi,\varepsilon)>0$, не зависящее
от $u$, что
\begin{equation}\label{3.93}
\|u\|_{H^{s,\varphi}}\leq\,c\,\bigl(\,
\|f\|_{H^{s-2q,\varphi}}+\|u\|_{H^{s-\varepsilon}}\,\bigr).
\end{equation}
\end{itemize}
\end{theorem}


Мы докажем теоремы \ref{th3.15} и \ref{th3.16} ниже в п.~\ref{sec3.4.3}.

Из теоремы \ref{th3.16} следует, что $N^{+}$ является дефектным подпространством
оператора \eqref{3.90}. Если $N=\{0\}$, то в оценке \eqref{3.93} можно опустить
слагаемое $\|u\|_{H^{s-\varepsilon}}$.

Отметим также следующее. Поскольку формально сопряженная краевая задача
\eqref{3.84}, \eqref{3.85} регулярная эллиптическая, то теоремы \ref{th3.15} и
\ref{th3.16} сохраняют силу для оператора $L^{+}$ вместо $L$ (с очевидными
изменениями в формулировке). Так, линейное отображение $v\mapsto L^{+}v$, где $v\in
C^{\infty}(\mathrm{b.c.})^{+}$, продолжается по непрерывности до ограниченного
нетерового оператора
\begin{equation}\label{3.94}
L^{+}:H^{2q-s,\,1/\varphi}(\mathrm{b.c.})^{+}\rightarrow
H^{-s,\,1/\varphi}/\,M_{-s,\,1/\varphi}^{+}= (H^{s,\varphi}(\mathrm{b.c.}))'
\end{equation}
с ядром $N^{+}$ и дефектным подпространством $N$. Здесь обозначено
$$
M_{-s,\,1/\varphi}^{+}=\{h\in H^{-s,\,1/\varphi}:\,(h,w)_{\Omega}=0\;\,\mbox{для
всех}\;\,w\in C^{\infty}(\mathrm{b.c.})\},
$$
а предположения относительно параметров $s$ и $\varphi$ такие как в теоремах
\ref{th3.15} и \ref{th3.16}. Сужение оператора \eqref{3.94} на подпространство
\begin{gather*}
P^{+}(H^{2q-s,\,1/\varphi}(\mathrm{b.c.})^{+})=\\ \{v\in
H^{2q-s,\,1/\varphi}(\mathrm{b.c.})^{+}:\,(v,w)_{\Omega}=0\;\,\mbox{для
всех}\;\,w\in N^{+}\}
\end{gather*}
является гомеоморфизмом
$$
L^{+}:P^{+}(H^{2q-s,\,1/\varphi}(\mathrm{b.c.})^{+})\leftrightarrow
P(H^{-s,\,1/\varphi})/M_{-s,\,1/\varphi}^{+}.
$$

В силу ограниченности операторов \eqref{3.90} и \eqref{3.94} равенство \eqref{3.86}
продолжается по непрерывности до соотношения
\begin{gather*}
(Lu,v)_{\Omega}=(u,L^{+}v)_{\Omega} \\ \mbox{для всех}\;\;u\in
H^{s,\varphi}(\mathrm{b.c.}),\;\;v\in H^{2q-s,\,1/\varphi}(\mathrm{b.c.})^{+}.
\end{gather*}
Это означает, что операторы \eqref{3.90} и \eqref{3.94} взаимно сопряжены
относительно полуторалинейной формы $(\,\cdot\,,\,\cdot\,)_{\Omega}$~--- расширения
по непрерывности скалярного произведения в $L_{2}(\Omega)$.

\subsection[Интерполяция и однородные краевые усло-\break вия]
{Интерполяция и однородные \\ краевые условия}\label{sec3.4.2}

Изучим интерполяционные свойства пространств, в которых действует оператор
$\eqref{3.90}$. Эти свойства сыграют важную роль в доказательствах теорем
\ref{th3.15} и \ref{th3.16}. Напомним, что в соболевском случае $\varphi\equiv1$ мы
опускаем индекс $\varphi$ в обозначения пространств, введенных в предыдущем
п.~\ref{sec3.4.1}.


\begin{theorem}\label{th3.17} Пусть $s\in\mathbb{R}$ и $\varphi\in\mathcal{M}$,
причем
\begin{equation}\label{3.95}
s\neq m_{j}+1/2\;\;\mbox{для каждого}\;\;j=1,\ldots,q.
\end{equation}
Тогда справедливы следующие утверждения.
\begin{itemize}
\item [$\mathrm{(i)}$] Для указанного $s$ найдется такое число $\varrho=\varrho(s)>0$, не
зависящее от $\varphi$, что при любом $\varepsilon\in(0,\varrho)$ справедливо
равенство пространств с точностью до эквивалентности норм:
$$
[H^{s-\varepsilon}(\mathrm{b.c.}),H^{s+\varepsilon}(\mathrm{b.c.})]_{\psi}=
H^{s,\varphi}(\mathrm{b.c.}),
$$
где $\psi$ --- не зависящий от $s$ интерполяционный параметр из
теоремы~$\ref{th2.14}$, в которой берем $\varepsilon=\delta$.
\item [$\mathrm{(ii)}$] Если число $s$, удовлетворяющее условию $\eqref{3.95}$,
положительно, то пространство $H^{s,\varphi}(\mathrm{b.c.})$ состоит из всех
распределений $u\in H^{s,\varphi}(\Omega)$, которые удовлетворяют условию
$B_{j}\,u=0$ на $\Gamma$ для всех $j=1,\ldots,q$ таких, что $s>m_{j}+1/2$.
\item [$\mathrm{(iii)}$] Если $s<1/2$, то $H^{s,\varphi}(\mathrm{b.c.})=H^{s,\varphi}$.
\item [$\mathrm{(iv)}$] Утверждения $\mathrm{(i)}$, $\mathrm{(ii)}$ и $\mathrm{(iii)}$
сохраняют силу, если в их формулировках заменить $m_{j}$ на $m_{j}^{+}$,
$(\mathrm{b.c.})$ на $(\mathrm{b.c.})^{+}$ и $B_{j}$ на $B_{j}^{+}$ соответственно.
\end{itemize}
\end{theorem}


\begin{remark}\label{rem3.9} В связи с утверждением (i) теоремы \ref{th3.17} отметим, что
интерполяция со степенным параметром пространств Соболева, удовлетворяющих
однородным краевым условиям, рассмотрена П.~Гриваром [\ref{Grisvard67}] и Р.~Сили
[\ref{Seeley72}] (см. также [\ref{Triebel80}, с.~400]). Из результатов этих работ
вытекает, что условие \eqref{3.95} отбросить нельзя. Далее, в утверждении (ii)
теоремы \ref{th3.17} выражение $B_{j}\,u$ понимается в смысле теоремы \ref{th3.5} о
следах. Условие \eqref{3.95} для этого утверждения обусловлено сказанным в замечании
\ref{rem3.5}.
\end{remark}


\textbf{Доказательство теоремы \ref{th3.17}} проведем отдельно для случаев $s>0$ и
$s<1/2$.

\emph{Случай } $s>0$. Изменив при необходимости нумерацию операторов системы
$\{B_{j}:\,j=1,\ldots,q\}$, получим, что
$$
0\leq m_{1}<m_{2}<\ldots<m_{q}\leq2q-1.
$$
Дополнительно положим $m_{0}:=-1/2$\, и $m_{q+1}:=+\infty$. В силу условия
\eqref{3.95} найдется номер $r\in\{0,1,\ldots,q\}$ такой, что
$m_{r}+1/2<s<m_{r+1}+1/2$. Обозначим через $\varrho=\varrho(s)$ положительное
расстояние от точки $s$ до множества
$$
\{j+1/2:\,j=-1/2,\,0,\,1,\ldots,2q-1\}\setminus\{s\}
$$
и возьмем произвольное $\varepsilon\in(0,\varrho)$. Тогда
\begin{gather}\label{3.96}
s\mp\varepsilon\neq j+1/2\;\;\mbox{для каждого}\;\;j=0,\,1,\ldots,2q-1,\\
0\leq m_{r}+1/2<s-\varepsilon<s<s+\varepsilon<m_{r+1}+1/2. \label{3.97}
\end{gather}

Пусть $\psi$ --- интерполяционный параметр из формулировки теоремы \ref{th2.14}, в
которой $\varepsilon=\delta$; этот параметр не зависит от~$s$. Проинтерполируем с
параметром $\psi$ пару пространств $H^{s\mp\varepsilon}(\mathrm{b.c.})$ на основании
теорем \ref{th3.10} и \ref{th2.6} (интерполяция подпространств). Для этого необходим
некоторый проектор $\mathcal{P}$ каждого пространства
$H^{s\mp\varepsilon}(\Omega)=H^{s\mp\varepsilon}$ на подпространство
$H^{s\mp\varepsilon}(\mathrm{b.c.})$. Он существует в силу следующих соображений.

Предположим сначала, что $r\neq0$ и рассмотрим набор $\{B_{j}:\,j=1,\ldots,r\}$.
Поскольку система $\{B_{j}:\,j=1,\ldots,q\}$ нормальна, то и рассмотренный набор,
как ее часть, является нормальной системой граничных выражений. Теперь сошлемся на
монографию Х.~Трибеля [\ref{Triebel80}, с.~484, 485] (лемма 5.4.4), в которой
построено линейное отображение $\mathcal{P}$, являющееся проектором каждого
пространства $H^{\sigma}(\Omega)=H^{\sigma}$, где $\sigma>m_{r}+1/2$, на
подпространство
\begin{equation}\label{3.98}
\{\,u\in
H^{\sigma}(\Omega):\;B_{j}\,u=0\;\;\mbox{на}\;\;\Gamma,\;\;j=1,\ldots,r\,\}.
\end{equation}
Возьмем здесь $\sigma=s\mp\varepsilon$; тогда в силу \eqref{3.97} подпространство
\eqref{3.98} допускает следующее описание:
\begin{gather}\notag
\{u\in H^{\sigma}(\Omega):\,B_{j}\,u=0\;\;\mbox{на}\;\;\Gamma\\ \mbox{для
всех}\;\;j=1,\ldots,q\;\;\mbox{таких, что}\;\;\sigma>m_{j}+1/2\}.\label{3.99}
\end{gather}
Но, как показано в монографии Я.~А.~Ройтберга [\ref{Roitberg96}, с.~167], из
регулярной эллиптичности задачи \eqref{3.82}, \eqref{3.83} вытекает в виду условия
\eqref{3.96} плотность множества $C^{\infty}(\mathrm{b.c.})$ в подпространстве
\eqref{3.99} пространства $H^{\sigma}(\Omega)$. Значит, множество
$C^{\infty}(\mathrm{b.c.})$ плотно в пространстве \eqref{3.98}, т.~е. последнее
совпадает с $H^{\sigma}(\mathrm{b.c.})$. Таким образом, отображение $\mathcal{P}$
является проектором пространства $H^{s\mp\varepsilon}(\Omega)$ на подпространство
$H^{s\mp\varepsilon}(\mathrm{b.c.})$.

Если $r=0$, то $0<s\mp\varepsilon<m_{j}+1/2$ для каждого номера $j=1,\ldots,q$.
Следовательно [\ref{Roitberg96}, с.~167], множество $C^{\infty}(\mathrm{b.c.})$
плотно в пространстве $H^{s\mp\varepsilon}(\Omega)=H^{s\mp\varepsilon}$ и поэтому
$H^{s\mp\varepsilon}(\mathrm{b.c.})=H^{s\mp\varepsilon}(\Omega)$. Значит, при $r=0$
в качестве $\mathcal{P}$ надо взять тождественное отображение. Итак, проектор
$\mathcal{P}$ построен.

Это позволяет на основании теорем \ref{th2.6} (интерполяция подпространств) и
\ref{th3.10} написать следующие равенства пространств с точностью до эквивалентности
норм в них:
\begin{gather}\notag
[H^{s-\varepsilon}(\mathrm{b.c.}), H^{s+\varepsilon}(\mathrm{b.c.})]_{\psi}=\\
[H^{s-\varepsilon}(\Omega),H^{s+\varepsilon}(\Omega)]_{\psi}
\,\cap\,H^{s-\varepsilon}(\mathrm{b.c.})=H^{s,\varphi}(\Omega)\,\cap\,
H^{s-\varepsilon}(\mathrm{b.c.})=\notag\\
H^{s,\varphi}(\Omega)\,\cap\,\{u\in
H^{s-\varepsilon}(\Omega):\;B_{j}\,u=0\;\;\mbox{на}\;\;\Gamma,\;\;j=1,\ldots,r\}=\notag\\
\{u\in H^{s,\varphi}(\Omega):\,B_{j}\,u=0\;\;\mbox{на}\;\;\Gamma\notag\\
\mbox{для всех}\;\;j=1,\ldots,q\;\;\mbox{таких, что}\;\;s>m_{j}+1/2\,\bigr\}.
\label{3.100}
\end{gather}
Заметим здесь, что последнее равенство вытекает из условия \eqref{3.97}. Таким
образом, интерполяционное пространство
$$
[H^{s-\varepsilon}(\mathrm{b.c.}),H^{s+\varepsilon}(\mathrm{b.c.})]_{\psi}
$$
совпадает (с точностью до эквивалентности норм) с подпространством \eqref{3.100}
пространства $H^{s,\varphi}(\Omega)=H^{s,\varphi}$. Значит, ввиду теоремы
\ref{th2.1} пространство $H^{s+\varepsilon}(\mathrm{b.c.})$ непрерывно и плотно
вложено в \eqref{3.100}. Отсюда вытекает, что множество $C^{\infty}(\mathrm{b.c.})$
плотно в \eqref{3.100}, т.~е. пространство \eqref{3.100} совпадает с
$H^{s,\varphi}(\mathrm{b.c.})$. Тем самым утверждения (i) (для случая $s>0$) и (ii)
доказаны.

\emph{Случай } $s<1/2$. Положим $\varrho=1/2-s>0$ и возьмем произвольное число
$\varepsilon\in(0,\varrho)$. Поскольку $s-\varepsilon<s<s+\varepsilon<1/2$, то в
силу \ref{th3.9} (ii) множество $C^{\infty}_{0}(\Omega)$, а значит, и более широкое
множество $C^{\infty}(\mathrm{b.c.})$ плотны в пространствах $H^{s,\varphi}$ и
$H^{s\mp\varepsilon}$. Поэтому $H^{s,\varphi}(\mathrm{b.c.})= H^{s,\varphi}$ и
$H^{s\mp\varepsilon}(\mathrm{b.c.})=H^{s\mp\varepsilon}$. Отсюда в силу теоремы
\ref{th3.10} непосредственно следует утверждение (i) доказываемой теоремы в случае
$s<1/2$.

Утверждение (iii) этой теоремы является следствием плотности множества
$C^{\infty}(\mathrm{b.c.})$ в пространстве $H^{s,\varphi}$.

Доказательство теоремы \ref{th3.17} для условия $(\mathrm{b.c.})$ закончено.
Заметим, что оно опирается лишь на регулярную эллиптичность краевой задачи
\eqref{3.82}, \eqref{3.83}. Следовательно, поскольку формально сопряженная задача
\eqref{3.84}, \eqref{3.85} также регулярная эллиптическая, то верно утверждение
(iv).

Теорема \ref{th3.17} доказана.

\medskip

Пусть $\sigma\in\mathbb{R}$ и $\varphi\in\mathcal{M}$. Изучим далее свойства
пространства $(H^{-\sigma,\,1/\varphi}(\mathrm{b.c.})^{+})'$, антидвойственного к
пространству $H^{-\sigma,\,1/\varphi}(\mathrm{b.c.})^{+}$.

Напомним, что
$$
M_{\sigma,\varphi}=\{h\in H^{\sigma,\varphi}:\,(h,w)_{\Omega}=0\;\,\mbox{для
всех}\;\,w\in C^{\infty}(\mathrm{b.c.})^{+}\}.
$$
Множество $M_{\sigma,\varphi}$ замкнуто в пространстве $H^{\sigma,\varphi}$. В самом
деле, в силу теоремы \ref{th3.9} (iii) функция $w\in
C^{\infty}(\mathrm{b.c.})^{+}\subset H^{-\sigma,\,1/\varphi}$ порождает линейный
непрерывный функционал $(\,\cdot\,,w)_{\Omega}$ на пространстве
$H^{\sigma,\varphi}$. Следовательно, если последовательность распределений $h_{j}\in
M_{\sigma,\varphi}$ такая, что $h_{j}\rightarrow h$ в $H^{\sigma,\varphi}$ при
$j\rightarrow \infty$, то
$$
(h,w)_{\Omega}=\lim_{j\rightarrow \infty}\,(h_{j},w)_{\Omega}=0\;\;\mbox{для
любого}\;\;w\in C^{\infty}(\mathrm{b.c.})^{+},
$$
т. е. $h\in M_{\sigma,\varphi}$. Итак, $M_{\sigma,\varphi}$ --- подпространство в
гильбертовом пространстве $H^{\sigma,\varphi}$. Поэтому факторпространство
$H^{\sigma,\varphi}/\,M_{\sigma,\varphi}$ тоже гильбертово.


\begin{theorem}\label{th3.18} Пусть $\sigma\in\mathbb{R}$ и $\varphi\in\mathcal{M}$.
Справедливы следующие утверждения.
\begin{itemize}
\item [$\mathrm{(i)}$] Если $\sigma>-1/2$, то $M_{\sigma,\varphi}=\{0\}$.
\item [$\mathrm{(ii)}$] Факторпространство $H^{\sigma,\varphi}/\,M_{\sigma,\varphi}$ и
подпространство $H^{-\sigma,\,1/\varphi}(\mathrm{b.c.})^{+}$ взаимно двойственны
(при $s\neq0$ с равенством норм, а при $s=0$ с эквивалентностью норм) относительно
расширения по непрерывности скалярного произведения в $L_{2}(\Omega)$. Точнее, ---
относительно билинейной формы
$\bigl(\,[u]\,,v\bigr)_{\Omega}:=\bigl(u,v\bigr)_{\Omega}$, где $[u]=\{u+h:\,h\in
M_{\sigma,\varphi}\}$ --- класс смежности элемента $u\in H^{\sigma,\varphi}$, а
$v\in H^{-\sigma,\,1/\varphi}(\mathrm{b.c.})^{+}$.
\item [$\mathrm{(iii)}$] Верны компактные плотные вложения
$$
H^{\sigma+\varepsilon}/\,M_{\sigma+\varepsilon}\hookrightarrow
H^{\sigma,\varphi}/\,M_{\sigma,\varphi}\hookrightarrow
H^{\sigma-\varepsilon}/\,M_{\sigma-\varepsilon}
$$
для любого $\varepsilon>0$;
\item [$\mathrm{(iv)}$] Если число $\sigma$ удовлетворяет условию
\begin{equation}\label{3.101}
-\sigma\neq m_{j}^{+}+1/2\quad\mbox{для каждого}\quad j=1,\ldots,q,
\end{equation}
то для $\sigma$ найдется такое число $\varrho=\varrho(\sigma)>0$, не зависящее от
$\varphi$, что при любом $\varepsilon\in(0,\varrho)$ справедливо равенство
пространств с точностью до эквивалентности норм
$$
[\,H^{\sigma-\varepsilon}/\,M_{\sigma-\varepsilon},
\,H^{\sigma+\varepsilon}/\,M_{\sigma+\varepsilon}\,]_{\psi}=
H^{\sigma,\varphi}/\,M_{\sigma,\varphi}.
$$
Здесь $\psi$ --- не зависящий от $\sigma$ интерполяционный параметр из теоремы
$\ref{th3.10}$, в которой берем $\varepsilon=\delta$.
\end{itemize}
\end{theorem}


\textbf{Доказательство.} (i) Предположим, что $\sigma>-1/2$ и $h\in
M_{\sigma,\varphi}$. В силу теоремы \ref{th3.9} (iii) функционал
$(h,\,\cdot\,)_{\Omega}$ ограничен на пространстве $H^{-\sigma,\,1/\varphi}$.
Согласно предположению, он равен нулю на множестве $C^{\infty}(\mathrm{b.c.})^{+}$,
которое в виду теоремы \ref{th3.9} (ii) плотно в этом пространстве. Значит указанный
функционал равен нулю как элемент сопряженного пространства
$(H^{-\sigma,\,1/\varphi})'$. Отсюда на основании теоремы \ref{th3.9} (iii) вытекает
равенство $h=0$. Утверждение (i) доказано.

(ii) Поскольку множество $C^{\infty}(\mathrm{b.c.})^{+}$ плотно в пространстве
$H^{-\sigma,\,1/\varphi}(\mathrm{b.c.})^{+}$, то
$$
M_{\sigma,\varphi}=\{h\in H^{\sigma,\varphi}:\,(h,w)_{\Omega}=0\;\,\mbox{для
всех}\;\,w\in H^{-\sigma,\,1/\varphi}(\mathrm{b.c.})^{+}\}
$$
и, следовательно,
$$
H^{-\sigma,\,1/\varphi}(\mathrm{b.c.})^{+}=\{w\in
H^{-\sigma,\,1/\varphi}:\,(h,w)_{\Omega}=0\;\,\mbox{для всех}\;\,h\in
M_{\sigma,\varphi}\}.
$$
Теперь утверждение (ii) является следствием теоремы \ref{th3.9} (iii) и известной
теоремы [\ref{FunctionalAnalysis72}, с.~47] об пространствах, сопряженных к
подпространству и факторпространству данного пространства.

(iii) Согласно теореме \ref{th3.9} (iv) верны компактные плотные вложения
$$
H^{\sigma+\varepsilon}\hookrightarrow H^{\sigma,\varphi}\hookrightarrow
H^{\sigma-\varepsilon}\quad\mbox{для любого}\;\;\varepsilon>0.
$$
Отсюда вытекает, что отображения
$$
\{u+h:\,h\in M_{\sigma+\varepsilon}\}\mapsto \{u+h:\,h\in
M_{\sigma,\varphi}\}\mapsto \{u+h:\,h\in M_{\sigma-\varepsilon}\},
$$
где $u\in H^{\sigma+\varepsilon}$ для первого отображения и $u\in
H^{\sigma,\varphi}$ для второго отображения, определяют компактные плотные вложения,
сформулированные в утверждении (iii).

(iv) Предположим, что число $s=-\sigma$ удовлетворяет условию \eqref{3.101}.
Рассмотрим сначала случай $s>0$. Обратимся к доказательству теоремы \ref{th3.17}, в
котором вместо задачи \eqref{3.82}, \eqref{3.83} возьмем формально сопряженную
задачу \eqref{3.84}, \eqref{3.85}. При этом в доказательстве следует заменить
обозначения $m_{j}$, $B_{j}$, $(\mathrm{b.c.})$ и $P$ на $m_{j}^{+}$, $B_{j}^{+}$,
$(\mathrm{b.c.})^{+}$ и $P^{+}$ соответственно. По доказанному, $\mathcal{P}^{+}$
--- проектор каждого пространства $H^{s\mp\varepsilon}$ на подпространство
$H^{s\mp\varepsilon}(\mathrm{b.c.})^{+}$. Здесь $\varepsilon\in(0,\varrho)$, где
$\varrho$ --- некоторое достаточно малое положительное число, причем
$s\mp\varepsilon>0$. Пусть $\Pi^{+}$ --- оператор, сопряженный к $\mathcal{P}^{+}$
относительно полуторалинейной формы $(\,\cdot\,,\,\cdot\,)_{\Omega}$. В силу теоремы
\ref{th3.9} (iii) и уже доказанного утверждения (ii) настоящей теоремы мы имеем
линейный ограниченный оператор
\begin{equation}\label{3.102}
\Pi^{+}:\,H^{\sigma\pm\varepsilon}/\,M_{\sigma\pm\varepsilon}\rightarrow
H^{\sigma\pm\varepsilon}.
\end{equation}
Он обладает следующим свойством:
\begin{equation}\label{3.103}
u\in H^{\sigma\pm\varepsilon}\;\Rightarrow\;u-\Pi^{+}[u]\in
M_{\sigma\pm\varepsilon}.
\end{equation}
Здесь $[u]=\{u+h:\,h\in M_{\sigma\pm\varepsilon}\}$ --- класс смежности элемента
$u$. Это вытекает из того, что для любых $u\in H^{\sigma\pm\varepsilon}$ и $w\in
C^{\infty}(\mathrm{b.c.})^{+}$
$$
(\Pi^{+}[u],w)_{\Omega}=(\,[u],P^{+}w)_{\Omega}= (\,[u],w)_{\Omega}=(u,w)_{\Omega},
$$
т. е. $(u-\Pi^{+}[u],w)_{\Omega}=0$. Теперь из свойств \eqref{3.102}, \eqref{3.103}
следует, что отображение
$$
u\mapsto u-\Pi^{+}[u],\quad u\in H^{\sigma\pm\varepsilon},
$$
является проектором пространства $H^{\sigma\pm\varepsilon}=
H^{\sigma\pm\varepsilon}_{\overline{\Omega}}(\mathbb{R}^{n})$ на подпространство
$M_{\sigma\pm\varepsilon}$ (последнее равенство вытекает из условия
$\sigma\pm\varepsilon=-(s\mp\varepsilon)<0)$. Это позволяет на основании теорем
\ref{th2.6} (интерполяция факторпространств) и \ref{th3.10} записать следующие
равенства пространств с точностью до эквивалентности норм:
\begin{gather*}
[\,H^{\sigma-\varepsilon}/\,M_{\sigma-\varepsilon}\,,
\,H^{\sigma+\varepsilon}/\,M_{\sigma+\varepsilon}\,]_{\psi}=\\
[\,H^{\sigma-\varepsilon}_{\overline{\Omega}}(\mathbb{R}^{n})
/\,M_{\sigma-\varepsilon}\,,
\,H^{\sigma+\varepsilon}_{\overline{\Omega}}(\mathbb{R}^{n})
/\,M_{\sigma+\varepsilon}\,]_{\psi}=\\
[\,H^{\sigma-\varepsilon}_{\overline{\Omega}}(\mathbb{R}^{n}),
\,H^{\sigma+\varepsilon}_{\overline{\Omega}}(\mathbb{R}^{n})\,]_{\psi}\,\big/
\bigl(\,[\,H^{\sigma-\varepsilon}_{\overline{\Omega}}(\mathbb{R}^{n}),
\,H^{\sigma+\varepsilon}_{\overline{\Omega}}(\mathbb{R}^{n})\,]_{\psi}
\cap\,M_{\sigma-\varepsilon}\,\bigr)=\\
H^{\sigma,\varphi}_{\overline{\Omega}}(\mathbb{R}^{n})\,\big/
\bigl(\,H^{\sigma,\varphi}_{\overline{\Omega}}(\mathbb{R}^{n})
\cap\,M_{\sigma-\varepsilon}\,\bigr)=\\
H^{\sigma,\varphi}\,\big/
\bigl(\,H^{\sigma,\varphi}\cap\,M_{\sigma-\varepsilon}\,\bigr)=
H^{\sigma,\varphi}\,/\,M_{\sigma,\varphi}.
\end{gather*}
Здесь $\psi$ --- интерполяционный параметр из теоремы \ref{th3.10}, где берем
$\varepsilon=\delta$. Тем самым утверждение (iv) доказано в случае $s=-\sigma>0$.

Противоположный случай $\sigma\geq0$ тривиален в силу уже доказанного
утверждения~(i). В самом деле, положив $\varrho:=1/2$, мы можем записать для
произвольного $\varepsilon\in(0,\,1/2)$ на основании теорем \ref{th3.9} (i) и
\ref{th3.10} следующее:
\begin{gather*}
[\,H^{\sigma-\varepsilon}/\,M_{\sigma-\varepsilon}\,,
\,H^{\sigma+\varepsilon}/\,M_{\sigma+\varepsilon}\,]_{\psi}=
[\,H^{\sigma-\varepsilon},\,H^{\sigma+\varepsilon}\,]_{\psi}=\\
[\,H^{\sigma-\varepsilon}(\Omega),\,H^{\sigma+\varepsilon}(\Omega)\,]_{\psi}=
H^{\sigma,\varphi}(\Omega)=H^{\sigma,\varphi}=
H^{\sigma,\varphi}\,/\,M_{\sigma,\varphi}.
\end{gather*}
Здесь, как и прежде, равенства пространств выполняются с точностью до
эквивалентности норм. Утверждение (iv) доказано.

Теорема \ref{th3.18} доказана.

\medskip

Как отмечалось выше, мы отождествляем в смысле теоремы \ref{th3.18} (ii)
факторпространство $H^{\sigma,\varphi}/\,M_{\sigma,\varphi}$ и сопряженное
пространство $(H^{-\sigma,\,1/\varphi}(\mathrm{b.c.})^{+})'$. При этом для каждого
распределения $u\in H^{\sigma,\varphi}$ его класс смежности $[u]=\{u+h:\,h\in
M_{\sigma,\varphi}\}$ отождествляется с антилинейным ограниченным функционалом
$(u,\,\cdot\,)_{\Omega}$ на пространстве
$H^{-\sigma,\,1/\varphi}(\mathrm{b.c.})^{+}$. Соответствующие нормы класса смежности
и функционала равны при $s\neq0$ и эквивалентны при $s=0$. Для краткости это
отождествление пространств мы записываем (несколько условно) как их равенство:
\begin{equation}\label{3.104}
H^{\sigma,\varphi}/\,M_{\sigma,\varphi}=
\bigl(\,H^{-\sigma,\,1/\varphi}(\mathrm{b.c.})^{+}\,\bigr)'
\end{equation}
для любых $\sigma\in\mathbb{R}$ и $\varphi\in\mathcal{M}$.

\subsection[Доказательство теорем о гомеоморфизмах и нетеровости]
{Доказательство теорем \\ о гомеоморфизмах и нетеровости}\label{sec3.4.3}

Докажем теоремы \ref{th3.15} и  \ref{th3.16}, сформулированные в п.~\ref{sec3.4.1}.

\textbf{Доказательство теоремы \ref{th3.15}.} В случае $\varphi\equiv1$
(пространства Соболева) эта теорема доказана в монографии Я.~А.~Ройтберга
[\ref{Roitberg96}, с.~168] (теорема 5.5.2). Случай произвольного
$\varphi\in\mathcal{M}$ мы выведем из этого результата с помощью интерполяции с
функциональным параметром.

В силу условия \eqref{3.89} имеем
\begin{equation}\label{3.105}
s\neq m_{j}+1/2,\;\;2q-s\neq m_{j}^{+}+1/2\;\;\mbox{для всех}\;\;j=1,\ldots,q.
\end{equation}
Следовательно, согласно теоремам \ref{th3.17} (i) и \ref{th3.18} (iv) (в последней
берем $\sigma=s-2q$) существует достаточно малое число $\varepsilon>0$ такое, что
\begin{gather}\label{3.106}
[\,H^{s-\varepsilon}(\mathrm{b.c.}),
\,H^{s+\varepsilon}(\mathrm{b.c.})\,]_{\psi}=H^{s,\varphi}(\mathrm{b.c.}),\\
[\,H^{s-2q-\varepsilon}/\,M_{s-2q-\varepsilon}\,,
\,H^{s-2q+\varepsilon}/\,M_{s-2q+\varepsilon}\,]_{\psi}=\notag\\
H^{s-2q,\,\varphi}/\,M_{s-2q,\,\varphi}. \label{3.107}
\end{gather}
Здесь $s\mp\varepsilon\neq j+1/2$ для каждого $j=0,\,1,\ldots,2q-1$, а $\psi$~---
интерполяционный параметр. Теперь сошлемся на теорему 5.5.2 монографии
Я.~А.~Ройтберга [\ref{Roitberg96}, с.~168], согласно которой линейное отображение
$u\mapsto Lu$, где $u\in C^{\infty}(\mathrm{b.c.})$, продолжается по непрерывности
до ограниченных операторов
$$
L:H^{s\mp\varepsilon}(\mathrm{b.c.})\rightarrow
H^{s\mp\varepsilon-2q}/M_{s\mp\varepsilon-2q}
$$
и соответствующих гомеоморфизмов
$$
L:P(H^{s\mp\varepsilon}(\mathrm{b.c.})\,)\leftrightarrow
P^{+}(H^{s\mp\varepsilon-2q})/M_{s\mp\varepsilon-2q}.
$$
Применив к ним интерполяцию с параметром $\psi$, получим ограниченный оператор
\begin{gather}\notag
L:[\,H^{s-\varepsilon}(\mathrm{b.c.}),
\,H^{s+\varepsilon}(\mathrm{b.c.})\,]_{\psi}\rightarrow\\
[\,H^{s-\varepsilon-2q}/\,M_{s-\varepsilon-2q}\,,
\,H^{s+\varepsilon-2q}/\,M_{s+\varepsilon-2q}\,]_{\psi} \label{3.108}
\end{gather}
и гомеоморфизм
\begin{gather}\notag
L:\,[\,P(H^{s-\varepsilon}(\mathrm{b.c.})),
\,P(H^{s+\varepsilon}(\mathrm{b.c.}))\,]_{\psi}\leftrightarrow\\
[\,P^{+}(H^{s-\varepsilon-2q})/\,M_{s-\varepsilon-2q}\,,
\,P^{+}(H^{s+\varepsilon-2q})/\,M_{s+\varepsilon-2q}\,]_{\psi}. \label{3.109}
\end{gather}
(Пары пространств, написанные в \eqref{3.109}, допустимые, что вытекает из теоремы
\ref{th2.6}; см. ниже.) В силу интерполяционных формул \eqref{3.106} и \eqref{3.107}
(а также равенства \eqref{3.104}) оператор \eqref{3.108} становится ограниченным
оператором \eqref{3.90} из формулировки доказываемой теоремы. Остается показать, что
\eqref{3.109} есть гомеоморфизм \eqref{3.91}.

Докажем сначала, что область определения гомеоморфизма \eqref{3.109} совпадает с
подпространством $P(H^{s,\varphi}(\mathrm{b.c.}))$. Согласно лемме \ref{lem3.2}
отображение $P$ является проектором пространства
$H^{s\mp\varepsilon}(\mathrm{b.c.})$ на подпространство
$$
P(H^{s\mp\varepsilon}(\mathrm{b.c.}))= \{u\in
H^{s\mp\varepsilon}(\mathrm{b.c.}):\,(u,w)_{\Omega}=0\;\,\mbox{для всех}\,\;w\in
N\}.
$$
Это позволяет на основании теоремы \ref{th2.6} (интерполяция подпространств) и
формулы \eqref{3.106} записать следующие равенства пространств с точностью до
эквивалентности норм:
\begin{gather*}
[\,P(H^{s-\varepsilon}(\mathrm{b.c.})),
\,P(H^{s+\varepsilon}(\mathrm{b.c.}))\,]_{\psi}=\\
[\,H^{s-\varepsilon}(\mathrm{b.c.}), \,H^{s+\varepsilon}(\mathrm{b.c.})\,]_{\psi}\,
\cap\,P(H^{s-\varepsilon}(\mathrm{b.c.}))=\\
H^{s,\varphi}(\mathrm{b.c.})\,\cap\,P(H^{s-\varepsilon}(\mathrm{b.c.}))
=P(H^{s,\varphi}(\mathrm{b.c.})).
\end{gather*}
Итак,
\begin{equation}\label{3.110}
[\,P(H^{s-\varepsilon}(\mathrm{b.c.})),
\,P(H^{s+\varepsilon}(\mathrm{b.c.}))\,]_{\psi}=P(H^{s,\varphi}(\mathrm{b.c.})).
\end{equation}

Покажем далее, что область значений гомеоморфизма \eqref{3.109} равна
$P^{+}(H^{s-2q,\,\varphi})/M_{s-2q,\,\varphi}$.  Согласно лемме \ref{lem3.2}
отображение $P^{+}$ является проектором пространства $H^{s\mp\varepsilon-2q}$ на
подпространство
$$
P^{+}(H^{s\mp\varepsilon-2q})=\{f\in
H^{s\mp\varepsilon-2q}:\;(f,w)_{\Omega}=0\;\,\mbox{для всех}\,\;w\in N^{+}\}.
$$
Рассмотрим линейное отображение
\begin{gather}\notag
[f]=\{f+h:\,h\in M_{s\mp\varepsilon-2q}\}\mapsto \\ [P^{+}f]=\{P^{+}f+h:\,h\in
M_{s\mp\varepsilon-2q}\}.\label{3.111}
\end{gather}
где $f\in H^{s\mp\varepsilon-2q}$. Оно определено корректно. Действительно,
поскольку $M_{s\mp\varepsilon-2q}\subset P^{+}(H^{s\mp\varepsilon-2q})$, то
$P^{+}h=h$ для каждого $h\in M_{s\mp\varepsilon-2q}$. Следовательно, класс смежности
$[P^{+}f]$ не зависит от выбора представителя $f$ класса смежности $[f]$. Из
сказанного выше об отображении $P^{+}$ вытекает, что отображение \eqref{3.111}
является проектором пространства $H^{s\mp\varepsilon-2q}/M_{s\mp\varepsilon-2q}$ на
подпространство $P^{+}(H^{s\mp\varepsilon-2q})/M_{s\mp\varepsilon-2q}$. Это
позволяет на основании теоремы \ref{th2.6} (интерполяция подпространств) и формулы
\eqref{3.107} записать следующие равенства пространств с точностью до
эквивалентности норм в них:
\begin{gather*}
[\,P^{+}(H^{s-\varepsilon-2q})/\,M_{s-\varepsilon-2q}\,,
\,P^{+}(H^{s+\varepsilon-2q})/\,M_{s+\varepsilon-2q}\,]_{\psi}=\\
[\,H^{s-\varepsilon-2q}/\,M_{s-\varepsilon-2q}\,,
\,H^{s+\varepsilon-2q}/\,M_{s+\varepsilon-2q}\,]_{\psi}\,\cap \\
(\,P^{+}(H^{s-\varepsilon-2q})/\,M_{s-\varepsilon-2q}\,)=\\
(\,H^{s-2q,\,\varphi}/\,M_{s-2q,\,\varphi})\,\,\cap\,
(\,P^{+}(H^{s-\varepsilon-2q})/\,M_{s-\varepsilon-2q}\,)=\\
(\,H^{s-2q,\,\varphi}\,\cap\,P^{+}(H^{s-\varepsilon-2q})\,)\big/
\,M_{s-2q,\,\varphi}\,=
\\ P^{+}(H^{s-2q,\,\varphi})/M_{s-2q,\,\varphi}.
\end{gather*}
Итак,
\begin{gather}\notag
[\,P^{+}(H^{s-\varepsilon-2q})/\,M_{s-\varepsilon-2q}\,,
\,P^{+}(H^{s+\varepsilon-2q})/\,M_{s+\varepsilon-2q}\,]_{\psi}=\\
P^{+}(H^{s-2q,\,\varphi})/M_{s-2q,\,\varphi}.\label{3.112}
\end{gather}

Теперь из интерполяционных формул \eqref{3.110} и \eqref{3.112} следует, что
отображение \eqref{3.109} осуществляет гомеоморфизм \eqref{3.91}.

Теорема \ref{th3.15} доказана.


\textbf{Доказательство теоремы \ref{th3.16}.} Напомним, что пространства $N$ и
$N^{+}$ конечномерны. Покажем сначала, что $N$ --- ядро оператора \eqref{3.90}.
Поскольку $N\subset C^{\infty}(\mathrm{b.c.})$, образом элемента $u\in N$ при
отображении \eqref{3.90} является класс смежности
$$
\{Lu+h:\,h\in M_{s-2q,\,\varphi}\}=\{0+h:\,h\in M_{s-2q,\,\varphi}\},
$$
т. е. нулевой элемент факторпространства $H^{s-2q,\,\varphi}/M_{s-2q,\,\varphi}$.
Верно и обратное: если элемент $u\in H^{s,\varphi}$ удовлетворяет условию $Lu=0$, то
написав в силу леммы \ref{lem3.2} разложение $u=u_{0}+Pu$, где $u_{0}\in N$,
получим, по уже доказанному, равенство
$$
0=Lu=Lu_{0}+LPu=LPu.
$$
Отсюда в силу гомеоморфизма \eqref{3.91} вытекает, что $Pu=0$, т.~е. $u=u_{0}\in N$.
Итак, $N$ --- конечномерное ядро оператора \eqref{3.90}.

Из гомеоморфизма \eqref{3.91} также следует, что область значений оператора
\eqref{3.90} совпадает с факторпространством
$P^{+}(H^{s-2q,\,\varphi})/M_{s-2q,\,\varphi}$, т.~е. с \eqref{3.92} и поэтому
замкнута в пространстве $H^{s-2q,\,\varphi}/M_{s-2q,\,\varphi}$. Теперь
коразмерность области значений равна размерности факторпространства
\begin{gather*} (\,H^{s-2q,\,\varphi}/M_{s-2q,\,\varphi}\,)\big/
(\,P^{+}(H^{s-2q,\,\varphi})/M_{s-2q,\,\varphi}\,)= \\
H^{s-2q,\,\varphi}/P^{+}(H^{s-2q,\,\varphi}).
\end{gather*}
Последняя в силу леммы \ref{lem3.2} совпадает с $\dim N^{+}$ и поэтому конечна.
Таким образом, оператор \eqref{3.90} нетеров и имеет индекс, равный $\dim N-\dim
N^{+}$.

Остается установить априорную оценку \eqref{3.93}. Возьмем произвольное
распределение $u\in H^{s,\varphi}(\mathrm{b.c.})$ и число $\varepsilon>0$. Согласно
лемме \ref{3.2} справедливо $u-Pu\in N$. Но $N$~--- конечномерное подпространство в
пространствах $H^{s,\varphi}$ и $H^{s-\varepsilon}$; следовательно, нормы в них
эквивалентны на $N$. В частности,
$$
\|u-Pu\|_{H^{s,\varphi}}\leq\,c_{1}\,\|u-Pu\|_{H^{s-\varepsilon}}
$$
с некоторой постоянной $c_{1}>0$, не зависящей от $u$. Отсюда получаем
\begin{gather*}
\|u\|_{H^{s,\varphi}}\leq\,\|u-Pu\|_{H^{s,\varphi}}+ \|Pu\|_{H^{s,\varphi}}\leq \\
c_{1}\,\|u-Pu\|_{H^{s-\varepsilon}}+
\|Pu\|_{H^{s,\varphi}}\leq\\
c_{1}\,\|u\|_{H^{s-\varepsilon}}+c_{1}\,\|Pu\|_{H^{s-\varepsilon}}+
\|Pu\|_{H^{s,\varphi}}\leq \\ c_{1}\,\|u\|_{H^{s-\varepsilon}}+
(c_{1}c_{2}+1)\,\|Pu\|_{H^{s,\varphi}}.
\end{gather*}
Здесь $c_{2}$~--- норма оператора вложения $H^{s,\varphi}\hookrightarrow
H^{s-\varepsilon}$ (см. теорему \ref{th3.9} (iv)). Итак,
\begin{equation}\label{3.113}
\|u\|_{H^{s,\varphi}}\leq\,c_{1}\,\|u\|_{H^{s-\varepsilon}}+
(c_{1}c_{2}+1)\,\|Pu\|_{H^{s,\varphi}}.
\end{equation}
Пусть теперь $Lu=[f]$. Поскольку, по уже доказанному, $N$ --- ядро оператора
\eqref{3.90} и $u-Pu\in N$, то также $LPu=[f]$. Таким образом, $Pu$ --- прообраз
класса смежности $[f]$ при гомеоморфизме \eqref{3.91}. Следовательно,
$$
\|Pu\|_{H^{s,\varphi}}\leq\,c_{0}\,\|f\|_{H^{s-2q,\varphi}},
$$
где $c_{0}$ --- норма оператора, обратного к \eqref{3.91}. Отсюда и из неравенства
\eqref{3.113} немедленно вытекает априорная оценка \eqref{3.93}.

Теорема \ref{th3.16} доказана.

\subsection[Локальное повышении гладкости решения вплоть до границы]
{Локальное повышении гладкости \\ решения вплоть до границы}\label{sec3.4.4}

Обозначим через $H^{-\infty}$ объединение всех пространств $\mathrm{H}^{s,\varphi}$,
где $s\in\mathbb{R}$ и $\varphi\in\mathcal{M}$. Положим
$$
M_{-\infty}:=\{h\in H^{-\infty}:\,(h,w)_{\Omega}=0\;\,\mbox{для всех}\,\;w\in
C(\mathrm{b.c.})^{+}\}.
$$
Операторы \eqref{3.90} при $s\in\mathbb{R}$ и $\varphi\in\mathcal{M}$, определяют
линейное отображение
$$
L:\,H^{-\infty}\rightarrow H^{-\infty}/M_{-\infty}
$$
и изоморфизм
\begin{equation}\label{3.114}
L:\,P(H^{-\infty})\leftrightarrow P^{+}(H^{-\infty})/M_{-\infty}
\end{equation}

Зададимся следующим вопросом. Пусть распределение $u\in H^{-\infty}$ удовлетворяет
уравнению
\begin{equation}\label{3.115}
Lu=\{f+h: h\in M_{-\infty}\},
\end{equation}
причем $f$ имеет заданную гладкость в уточненной шкале пространств на некотором
открытом в $\overline{\Omega}$ множестве $U$. Что тогда можно сказать о гладкости
решения и на этом множестве? Дадим ответ на этот вопрос.

Пусть $U$ --- открытое множество в пространстве $\mathbb{R}^{n}$, причем
$\Omega_{0}:=\overline{\Omega}\cap U\neq\varnothing$. Положим
$\Gamma_{0}:=\Gamma\cap U$ (возможен случай $\Gamma_{0}=\varnothing$). Введем
следующие пространства уточненной гладкости на $\Omega_{0}$.

Пусть $s\in\mathbb{R}$ и $\varphi\in\mathcal{M}$. Обозначим
\begin{gather*}
H_{\mathrm{loc}}^{s,\varphi,(0)}(\Omega_{0}):=\\ \{u\in H^{-\infty}:\;\chi u\in
H^{s,\varphi}\;\;\mbox{для всех}\;\;\chi\in C^{\infty}
(\overline{\Omega}),\;\mathrm{supp}\,\chi\subset\Omega_{0}\}.
\end{gather*}
Отметим здесь, что оператор умножения на функцию $\chi\in
C^{\infty}(\overline{\Omega})$ ограничен в пространстве $H^{s,\varphi}$. В
соболевском случае $\varphi\equiv1$ это хорошо известно (см., например,
[\ref{Roitberg96}, с.~57], п.~1.12), а для произвольного $\varphi\in\mathcal{M}$
следует из интерполяционной теоремы \ref{th3.10}.

Далее, ввиду теоремы \ref{th3.17} (ii) положим
\begin{gather*}
H_{\mathrm{loc}}^{s,\varphi}(\Omega_{0},\mathrm{b.c.},\Gamma_{0}):=\{u\in
H_{\mathrm{loc}}^{s,\varphi,(0)}(\Omega_{0}):\, B_{j}\,u=0\;\;\text{на}\;\Gamma_{0} \\
\mbox{для всех}\;\;j=1,\ldots,q\;\;\mbox{таких, что}\;\;s>m_{j}+1/2\}.
\end{gather*}
Здесь предполагается, что $s\neq m_{j}+1/2$ для всех $j=1,\ldots,q$.


\begin{theorem}\label{th3.19}
Пусть $u\in \mathrm{H}^{-\infty}$ является решением уравнения \eqref{3.115}, в
котором $f\in H_{\mathrm{loc}}^{s-2q,\varphi,(0)}(\Omega_{0})$ для некоторого
$s\in\mathbb{R}$, удовлетворяющего \eqref{3.89}, и некоторого
$\varphi\in\mathcal{M}$. Тогда $u\in
H_{\mathrm{loc}}^{s,\varphi}(\Omega_{0},\mathrm{b.c.},\Gamma_{0})$.
\end{theorem}


Теорема \ref{th3.19} --- это утверждение о локальном повышении уточненной гладкости
решения $u$ краевой задачи \eqref{3.82}, \eqref{3.83} вплоть до границы области
$\Omega$. Отметим, что в случае $\overline{\Omega}\subset U$, т.~е. когда
$\Omega_{0}=\overline{\Omega}$ и $\Gamma_{0}=\Gamma$, „локальные” пространства
совпадают с „глобальными”:
$$
H_{\mathrm{loc}}^{s,\varphi,(0)}(\Omega_{0})=H^{s,\varphi}\quad\mbox{и}\quad
H_{\mathrm{loc}}^{s,\varphi}(\Omega_{0},\mathrm{b.c.},\Gamma_{0})=
H^{s,\varphi}(\mathrm{b.c.}).
$$
Поэтому теорема \ref{th3.19} содержит в себе также утверждение о глобальном
повышении гладкости, т.е. во всей замкнутой области $\overline{\Omega}$. Наконец,
отметим еще случай $U\subset\Omega$, т.~е. $\Gamma_{0}=\varnothing$, который
приводит к утверждению о повышении гладкости внутри области $\Omega$.


\begin{remark}\label{rem3.10}
В соболевском случае $\varphi\equiv1$ теорема \ref{th3.19} доказана
Ю.~М.~Березанским, С.~Г.~Крейном и Я.~А.~Ройтбергом
[\ref{BerezanskyKreinRoitberg63}; \ref{Roitberg63}; \ref{Berezansky65};
\ref{Roitberg96}, с.~221].
\end{remark}


\textbf{Доказательство теоремы \ref{th3.19}.} Как отмечалось в замечании
\ref{rem3.9}, эта теорема известна в случае $\varphi\equiv1$. Отсюда мы выведем ее
для произвольного $\varphi\in\mathcal{M}$ с помощью теоремы \ref{th3.15} о
гомеоморфизме.

Для произвольного числа $\varepsilon>0$ положим
$$
U_{\varepsilon}:=\{x\in U:\,\mathrm{dist}(x,\partial U)>\varepsilon\},\quad
\Omega_{\varepsilon}:=\overline{\Omega}\cap U_{\varepsilon},\quad
\Gamma_{\varepsilon}:=\Gamma\cap U_{\varepsilon}.
$$
Здесь, как обычно, $\partial U$ --- граница множества $U$. Выберем такую функцию
$\chi_{\varepsilon}\in C^{\infty}(\overline{\Omega})$, что
$\mathrm{supp}\,\chi_{\varepsilon}\subset\Omega_{0}$, и $\chi_{\varepsilon}\equiv1$
на $\Omega_{\varepsilon}$.

Представим распределение $f$ в виде $f=\chi_{\varepsilon}f+(1-\chi_{\varepsilon})f$.
Поскольку $f\in H_{\mathrm{loc}}^{s-2q,\varphi,(0)}(\Omega_{0})$, и
$\mathrm{supp}\,\chi_{\varepsilon}\subset\Omega_{0}$, то $\chi_{\varepsilon}f\in
H^{s-2q,\varphi}$. Отсюда в силу гомеоморфизма \eqref{3.91} из теоремы \ref{th3.15},
существует распределение $u_{\varepsilon}\in P(H^{s,\varphi}(\mathrm{b.c.}))$ такое,
что
\begin{equation}\label{3.116}
Lu_{\varepsilon}=[\,P^{+}\chi_{\varepsilon}f\,].
\end{equation}

Далее, поскольку $1-\chi_{\varepsilon}=0$ на множестве $\Omega_{\varepsilon}$, то
$\chi(1-\chi_{\varepsilon})f\equiv0$ для любой функции $\chi\in
C^{\infty}(\overline{\Omega})$, у которой
$\mathrm{supp}\,\chi\subset\Omega_{\varepsilon}$. Отсюда, в частности, следует что
$(1-\chi_{\varepsilon})f\in H_{\mathrm{loc}}^{s+1-2q,(0)}(\Omega_{\varepsilon})$.
Согласно лемме \ref{3.2} запишем:
$$
(1-\chi_{\varepsilon})f=g_{\varepsilon}+P^{+}(1-\chi_{\varepsilon})f,
$$
где $g_{\varepsilon}\in N^{+}\subset C^{\infty}(\mathrm{b.c.})$. Значит,
$$
P^{+}(1-\chi_{\varepsilon})f\in H_{\mathrm{loc}}^{s+1-2q,(0)}(\Omega_{\varepsilon}).
$$
Кроме того, в силу \eqref{3.114} существует $v_{\varepsilon}\in P(H^{-\infty})$
такое, что
\begin{equation}\label{3.117}
Lv_{\varepsilon}=[\,P^{+}(1-\chi_{\varepsilon})f\,].
\end{equation}
Следовательно, согласно теореме 7.3.1 из монографии Я.~А.~Ройтберга
[\ref{Roitberg96}, с.~221] выполняется включение $v_{\varepsilon}\in
H_{\mathrm{loc}}^{s+1}(\Omega_{\varepsilon},\mathrm{b.c.},\Gamma_{\varepsilon})$.
Отсюда в силу \eqref{3.116}, \eqref{3.117} и условия $u_{\varepsilon}\in
H^{s,\varphi}(\mathrm{b.c.})$ имеем равенство
$$
L(u_{\varepsilon}+v_{\varepsilon})=
[\,P^{+}\chi_{\varepsilon}f\,]+[\,P^{+}(1-\chi_{\varepsilon})f\,]=[\,P^{+}f\,],
$$
в котором
$$
u_{\varepsilon}+v_{\varepsilon}\in H^{s,\varphi}(\mathrm{b.c.})\cup
H_{\mathrm{loc}}^{s+1}(\Omega_{\varepsilon},\mathrm{b.c.},\Gamma_{\varepsilon})\subset
H_{\mathrm{loc}}^{s,\varphi}(\Omega_{\varepsilon},\mathrm{b.c.},\Gamma_{\varepsilon}).
$$
Здесь мы воспользовались вложением $H^{s+1}\hookrightarrow H^{s,\varphi}$ и
определением пространства
$H_{\mathrm{loc}}^{s,\varphi}(\Omega_{\varepsilon},\mathrm{b.c.},\Gamma_{\varepsilon})$.

Но $Lu=[f]$, значит, в силу \eqref{3.114} верно $f=P^{+}f$. Следовательно,
$$
L(u-u_{\varepsilon}-v_{\varepsilon})=[f]-[\,P^{+}f\,]=0,
$$
где $u_{\varepsilon}+v_{\varepsilon}\in
H_{\mathrm{loc}}^{s,\varphi}(\Omega_{\varepsilon},\mathrm{b.c.},\Gamma_{\varepsilon})$.
Отсюда в силу теоремы \ref{th3.15} имеем:
$$
v:=u-u_{\varepsilon}-v_{\varepsilon}\in N\subset C^{\infty}(\mathrm{b.c.}).
$$
Таким образом,
$$
u=u_{\varepsilon}+v_{\varepsilon}+v\in
H_{\mathrm{loc}}^{s,\varphi}(\Omega_{\varepsilon},\mathrm{b.c.},\Gamma_{\varepsilon}).
$$

Напомним, что в доказательстве число $\varepsilon>0$ произвольно. Очевидно, что
$$
H_{\mathrm{loc}}^{s,\varphi}(\Omega_{0},\mathrm{b.c.},\Gamma_{0})=
\bigcap_{\varepsilon>0}
H_{\mathrm{loc}}^{s,\varphi}(\Omega_{\varepsilon},\mathrm{b.c.},\Gamma_{\varepsilon}).
$$
Поэтому $u\in H_{\mathrm{loc}}^{s,\varphi}(\Omega_{0},\mathrm{b.c.},\Gamma_{0})$.

Теорема \ref{th3.19} доказана.




\markright{\emph \ref{sec3.5}. Некоторые свойства пространств Хермандера}

\section[Некоторые свойства пространств Хермандера]
{Некоторые свойства \\ пространств Хермандера}\label{sec3.5}

\markright{\emph \ref{sec3.5}. Некоторые свойства пространств Хермандера}

В этом пункте мы изучим некоторые подпространства пространства
$H^{s,\varphi}(\Omega)$, а также дадим одно эквивалентное описание этого
пространства. Приведенные здесь результаты понадобятся в дальнейшем, в частности, в
п.~\ref{sec4.5}.

\subsection{Пространство $H^{s,\varphi}_{0}(\Omega)$ и его свойства}\label{sec3.5.1}

Пусть $s\in\mathbb{R}$ и $\varphi\in\mathcal{M}$. Обозначим через замыкание
множества $C^{\infty}_{0}(\Omega)$ в пространстве $H^{s,\varphi}(\Omega)$. Мы
рассматриваем $H^{s,\varphi}_{0}(\Omega)$ как гильбертово пространство относительно
скалярного произведения в $H^{s,\varphi}(\Omega)$. В этом пункте мы изучим свойства
пространства $H^{s,\varphi}_{0}(\Omega)$, что приведет нас к различным эквивалентным
описаниям пространства $H^{s,\varphi}(\Omega)$ для неполуцелых отрицательных
индексов~$s$. Эти важные результаты будут нужны нам в дальнейшем. Мы отложили их
доказательство до настоящего пункта, поскольку для этого нам понадобится один
результат п.~\ref{sec3.4.2}.

Для $k\in\mathbb{N}$ положим
\begin{equation}\label{3.118}
C^{\infty}_{\nu,k}(\,\overline{\Omega}\,):=\{u\in
C^{\infty}(\,\overline{\Omega}\,):\,D^{j-1}_{\nu}u=0\;\mbox{на}\;\Gamma,\;j=1,\ldots,k\}.
\end{equation}
Здесь и далее нам удобно использовать обозначение
$D_{\nu}:=i\,\partial/\partial\nu$. Обозначим также через
$H^{s,\varphi}_{\nu,k}(\Omega)$ замыкание множества
$C^{\infty}_{\nu,k}(\,\overline{\Omega}\,)$ в пространстве $H^{s,\varphi}(\Omega)$.
Мы рассматриваем $H^{s,\varphi}_{\nu,k}(\Omega)$ как гильбертово пространство
относительно скалярного произведения в $H^{s,\varphi}(\Omega)$.

Если $s\geq0$, то $H^{s,\varphi}_{\nu,k}(\Omega)$ --- это пространство
$H^{s,\varphi}(\mathrm{b.c.})$, введенное в п.~\ref{sec3.4.1}. Здесь
$(\mathrm{b.c.})$ обозначает однородные краевые условия Дирихле, фигурирующие в
\eqref{3.118}. Поэтому согласно теореме \ref{th3.17}~(ii)
\begin{gather}\notag
H^{s,\varphi}_{\nu,k}(\Omega)=\{u\in
H^{s,\varphi}(\Omega):\;D^{j-1}_{\nu}u=0\;\;\mbox{на}\;\;\Gamma\phantom{mmm} \\
\mbox{для всех}\;\;j=1,\ldots,k\}\quad\mbox{при}\;\;s>k-1/2.\label{3.119}
\end{gather}

В силу теоремы \ref{th3.3} (ii)
\begin{equation}\label{3.120}
H^{s,\varphi}_{0}(\Omega)=H^{s,\varphi}(\Omega)\quad\mbox{при}\;\;s<1/2.
\end{equation}
Если $s>1/2$, то включение $H^{s,\varphi}_{0}(\Omega)\subset H^{s,\varphi}(\Omega)$
строгое. Это вытекает из теоремы \ref{th3.5} о следах.


\begin{theorem} \label{th3.20}
Пусть $s>1/2$, $s-1/2\notin\mathbb{Z}$, и $\varphi\in\mathcal{M}$. Справедливы
следующие утверждения.
\begin{itemize}
\item [$\mathrm{(i)}$] $H^{s,\varphi}_{0}(\Omega)$ совпадает с пространством
\eqref{3.119}, где $k:=[s+1/2]$.
\item [$\mathrm{(ii)}$] Нормы в пространствах $H^{s,\varphi}_{0}(\Omega)$ и
$H^{s,\varphi}_{\overline{\Omega}}(\mathbb{R}^{n})$ эквивалентны на плотном
множестве $C^{\infty}_{0}(\Omega)$, и поэтому
$H^{s,\varphi}_{0}(\Omega)=H^{s,\varphi}_{\overline{\Omega}}(\mathbb{R}^{n})$ с
точностью до эквивалентности норм.
\item [$\mathrm{(iii)}$] Пространства $H^{s,\varphi}_{0}(\Omega)$ и
$H^{-s,1/\varphi}(\Omega)$ взаимно двойственные (с точностью до эквивалентности
норм) относительно скалярного произведения в $L_{2}(\Omega)$.
\end{itemize}
\end{theorem}


\textbf{Доказательство.} В случае $\varphi\equiv1$ (пространства Соболева) эта
теорема известна. Ее доказательство имеется, например, в монографии Х.~Трибеля
[\ref{Triebel80}, с. 396, 412, 414]. А именно, утверждение (i) содержится в теореме
4.7.1 (a), утверждение (ii) --- в теореме 4.3.2/1 (c) и утверждение (iii) --- в
теореме 4.8.2 (a) этой монографии. Выведем отсюда последовательно утверждения (i)~--
(iii) для общего случая.

(i) Пусть $k:=[s+1/2]$, тогда $k-1/2<s<k+1/2$. Выберем число $\varepsilon>0$ так,
чтобы
\begin{equation}\label{3.121}
k-1/2<s\mp\varepsilon<k+1/2.
\end{equation}
Имеем непрерывное и плотное вложение
$$
H^{s+\varepsilon}_{0}(\Omega)= H^{s+\varepsilon}_{\nu,k}(\Omega)\hookrightarrow
H^{s,\varphi}_{\nu,k}(\Omega).
$$
Как указывалось выше, записанное здесь равенство пространств Соболева верно, так как
$[s+\varepsilon+1/2]=k$. Поэтому множество $C^{\infty}_{0}(\Omega)$ плотно в
$H^{s,\varphi}_{\nu,k}(\Omega)$. Следовательно,
$H^{s,\varphi}_{0}(\Omega)=H^{s,\varphi}_{\nu,k}(\Omega)$ и утверждение~(i)
доказано.

(ii) Мы выведем утверждение (ii) из соболевского случая $\varphi\equiv1$ с помощью
интерполяции. Пусть числа $k$ и $\varepsilon$ такие как и в предыдущем абзаце.
Обозначим через $\mathcal{O}$ оператор продолжения нулем функции с области
$\overline{\Omega}$ в пространство $\mathbb{R}^{n}$. Отображение
$u\mapsto\mathcal{O}u$, где $u\in C^{\infty}_{0}(\Omega)$, продолжается по
непрерывности до гомеоморфизмов
$$
\mathcal{O}:\,H^{s\mp\varepsilon}_{0}(\Omega)\leftrightarrow
H^{s\mp\varepsilon}_{\overline{\Omega}}(\mathbb{R}^{n}).
$$
Применим здесь интерполяцию с функциональным параметром $\psi$ из теоремы
\ref{th2.14}, где $\delta=\varepsilon$. Ввиду теоремы \ref{th3.7} получим еще один
гомеоморфизм
\begin{equation}\label{3.122}
\mathcal{O}:\,[H^{s-\varepsilon}_{0}(\Omega),H^{s+\varepsilon}_{0}(\Omega)\bigr]_{\psi}
\leftrightarrow H^{s,\varphi}_{\overline{\Omega}}(\mathbb{R}^{n}).
\end{equation}

Опишем область определения оператора \eqref{3.122}. В силу неравенства \eqref{3.121}
и утверждения (i) имеем:
$$
H^{s\mp\varepsilon}_{0}(\Omega)=H^{s\mp\varepsilon}_{\nu,k}(\Omega)\quad\mbox{и}
\quad H^{s,\varphi}_{\nu,k}(\Omega)=H^{s,\varphi}_{0}(\Omega).
$$
Отсюда на основании теоремы \ref{th3.17} (i) получаем:
\begin{gather*}
[H^{s-\varepsilon}_{0}(\Omega),H^{s+\varepsilon}_{0}(\Omega)]_{\psi}=
[H^{s-\varepsilon}_{\nu,k}(\Omega), H^{s+\varepsilon}_{\nu,k}(\Omega)\bigr]_{\psi}= \\
H^{s,\varphi}_{\nu,k}(\Omega)=H^{s,\varphi}_{0}(\Omega).
\end{gather*}
Здесь второе равенство пространств выполняется с точностью до эквивалентности норм.
Следовательно, \eqref{3.122} означает гомеоморфизм
$$
\mathcal{O}:\,H^{s,\varphi}_{0}(\Omega)\leftrightarrow
H^{s,\varphi}_{\overline{\Omega}}(\mathbb{R}^{n}).
$$
Тем самым доказано утверждение (ii).

Утверждение (iii) сразу следует из теоремы \ref{th3.8} (iii) и утверждения~(ii).

Теорема \ref{th3.20} доказана.

\subsection{Эквивалентное описание $H^{s,\varphi}(\Omega)$}\label{sec3.5.2}

Согласно теореме \ref{th3.20} (iii) пространство $H^{s,\varphi}(\Omega)$ для
неполуцелых отрицательных индексов $s$ можно определить как двойственное
пространство к $H^{-s,1/\varphi}_{0}(\Omega)$ относительно скалярного произведения в
$L_{2}(\Omega)$. В соболевском случае такое определение использовали Ж.-Л.~Лионс и
Э.~Мадженес [\ref{LionsMagenes71}, \ref{Magenes66}, \ref{LionsMagenes62},
\ref{LionsMagenes63}] при изучении краевых задач (см. замечание \ref{rem3.6}).

Отсюда следует другое (эквивалентное) описание пространства $H^{s,\varphi}(\Omega)$.
Напомним, что, по определению,
$$
H^{s,\varphi}(\Omega)=\{w\!\upharpoonright\!\Omega:\,w\in
H^{s,\varphi}(\mathbb{R}^{n})\}.
$$
Оказывается, здесь можно ограничиться распределениями $w$, сосредоточенными
на~$\overline{\Omega}$.


\begin{theorem} \label{th3.21}
Пусть $s<1/2$, $s-1/2\notin\mathbb{Z}$, и $\varphi\in\mathcal{M}$. Тогда
\begin{gather}\label{3.123}
H^{s,\varphi}(\Omega)= H^{s,\varphi}_{\overline{\Omega}}(\mathbb{R}^{n})/
H^{s,\varphi}_{\Gamma}(\mathbb{R}^{n})=\{w\!\upharpoonright\!\Omega:\,w\in
H^{s,\varphi}_{\overline{\Omega}}(\mathbb{R}^{n})\}, \\
\|u\|_{H^{s,\varphi}(\Omega)}\asymp
\inf\,\{\,\|w\|_{H^{s,\varphi}(\mathbb{R}^{n})}:\, w\in
H^{s,\varphi}_{\overline{\Omega}}(\mathbb{R}^{n}),\;w=u\;\mbox{в}\;\,\Omega\,\}.
\label{3.124}
\end{gather}
\end{theorem}


\textbf{Доказательство.} В силу теоремы \ref{th3.20} (iii) имеем для неполуцелых
$s<-1/2$ равенство
\begin{equation}\label{3.125}
H^{s,\varphi}(\Omega)=(H^{-s,1/\varphi}_{0}(\Omega))'\quad\mbox{с эквивалентностью
норм}.
\end{equation}
Двойственность пространств мы рассматриваем относительно скалярного произведения в
$L_{2}(\Omega)$. Если $-1/2<s<1/2$, то \eqref{3.125} также верно в силу теоремы
\ref{th3.8} (ii), (iii) и равенства \eqref{3.120}:
$$
H^{s,\varphi}(\Omega)=(H^{-s,1/\varphi}_{\overline{\Omega}}(\mathbb{R}^{n}))'=
(H^{-s,1/\varphi}(\Omega))'=(H^{-s,1/\varphi}_{0}(\Omega))'.
$$
Опишем $(H^{-s,1/\varphi}_{0}(\Omega))'$ как двойственное пространство к
подпространству $H^{-s,1/\varphi}_{0}(\Omega)$ пространства
$H^{-s,1/\varphi}(\Omega)$. Ввиду теоремы \ref{th3.8} (iii) имеем:
$$
(H^{-s,1/\varphi}_{0}(\Omega))'=(H^{-s,1/\varphi}(\Omega))'/G_{s,\varphi}=
H^{s,\varphi}_{\overline{\Omega}}(\mathbb{R}^{n})/G_{s,\varphi},
$$
Здесь
\begin{gather*}
G_{s,\varphi}:=\{w\in
H^{s,\varphi}_{\overline{\Omega}}(\mathbb{R}^{n}):\,(w,v)_{\Omega}=0
\;\,\mbox{для всех}\,\;v\in H^{-s,1/\varphi}_{0}(\Omega)\}= \\
\{w\in
H^{s,\varphi}_{\overline{\Omega}}(\mathbb{R}^{n}):\,(w,v)_{\Omega}=0\;\,\mbox{для
всех}\,\;v\in C^{\infty}_{0}(\Omega)\}= H^{s,\varphi}_{\Gamma}(\mathbb{R}^{n}).
\end{gather*}
Таким образом, в силу \eqref{3.125}
$$
H^{s,\varphi}(\Omega)= H^{s,\varphi}_{\overline{\Omega}}(\mathbb{R}^{n})/
H^{s,\varphi}_{\Gamma}(\mathbb{R}^{n})
$$
с точностью до эквивалентных норм. Отсюда немедленно следуют формулы \eqref{3.123} и
\eqref{3.124}.

Теорема \ref{th3.21} доказана.




\markright{\emph \ref{notes3}. Примечания и комментарии}

\section{Примечания и комментарии}\label{notes3}

\markright{\emph \ref{notes3}. Примечания и комментарии}

\small

\textbf{К п. 3.1.} Понятие общей эллиптической краевой задачи впервые сформулировано
Я.~Б.~Лопатинским [\ref{Lopatinsky53}, \ref{Lopatinsky84}] и, в частных случаях,
З.~Я.~Шапиро [\ref{Shapiro53}]. Н.~Ароншайн, А.~Мильграм [\ref{AronszajnMilgram52}]
и независимо М.~Шехтер [\ref{Schechter60a}] ввели важное условие нормальности
системы граничных выражений, обеспечивающее существование формально сопряженной
краевой задачи в классе дифференциальных операторов.

Систематическое изложение теории общих эллиптических краевых задач имеется,
например, в монографиях Ю.~М.~Березанского [\ref{Berezansky65}], Ж.-Л.~Лионса и
Э.~Мадженеса [\ref{LionsMagenes71}], О.~И.~Панича [\ref{Panich86}], Я.~А.~Ройтберга
[\ref{Roitberg96}], Х.~Трибеля [\ref{Triebel80}], Л.~Хермандера [\ref{Hermander65},
\ref{Hermander87}], М.~Шехтера [\ref{Schechter77}], в обзоре М.~С.~Аграновича
[\ref{Agranovich97}].

\medskip

\textbf{К п. 3.2.} Для открытых и замкнутых евклидовых областей гильбертовы
пространства Хермандера введены и изучены Л.~Р.~Волевичем и Б.~П.~Панеяхом
[\ref{VolevichPaneah65}]. Ими использованы стандартные определения, обычно
используемые в теории функциональных пространств (см., например, монографию
Х.~Трибеля [\ref{Triebel80}, с.~384, 395]).

Уточненные шкалы пространств для евклидовых областей введены авторами в
[\ref{05UMJ5}]. В~отличие от Л.~Р.~Волевича и Б.~П.~Панеяха мы устанавливаем
свойства пространств Хермандера на областях, используя интерполяционные формулы,
связывающие эти пространства с соболевскими шкалами. Теоремы 3.1~-- 3.5 о свойствах
уточненной шкалы в евклидовой области доказаны нами в [\ref{06UMJ3}] (п.~3) в более
общей ситуации уточненной шкалы на многообразии с краем. Теоремы 3.6~-- 3.8 о
свойствах уточненной шкалы на замкнутой евклидовой области установлены в
[\ref{06UMB4}] (п.~4). Там же доказана теорема 3.9, а доказательство
интерполяционной теоремы 3.10 ранее не приводилось. Последние две теоремы описывают
свойства оснащения пространства $L_{2}(\Omega)$ пространствами Хермандера.

Определение гильбертова оснащения и связанные с ним понятия см., например, в
монографиях Ю.~М.~Березанского [\ref{Berezansky65}, с.~47], Ю.~М.~Березанского,
Г.~Ф.~Уса, З.~Г.~Шефтеля [\ref{BerezanskyUsSheftel90}, с.~462]. Оснащение
соболевскими пространствами введено и изучено Ю.~М.~Березанским [\ref{Berezansky65},
с.~66; \ref{BerezanskyUsSheftel90}, с.~475].

В~статье Г.~Шлензак [\ref{Shlenzak74}, с.~55] доказана интерполяционная формула,
связывающая соболевскую шкалу с некоторыми пространствами Хермандера, заданными на
евклидовой области с гладкой границей. Последние применены Г.~Шлензак к теории общих
эллиптических краевых задач.

\medskip

\textbf{К п. 3.3.} Все теоремы этого пункта доказаны авторами в [\ref{06UMJ11}].
Основной результат (теорема 3.11) является новым даже в соболевском случае для
полуцелых $s<0$. Для остальных значений параметра $s$ он содержится в теореме
Ж.-Л.~Лионса и Э.~Мадженеса [\ref{LionsMagenes71}, с. 216, 217] о разрешимости
неоднородной регулярной эллиптической краевой задачи в двусторонней шкале
пространств Соболева (см. также их статьи [\ref{LionsMagenes62},
\ref{LionsMagenes63}]).

\medskip

\textbf{К п. 3.4.} Теорема о гомеоморфизмах, осуществляемых эллиптическим оператором
в двусторонней соболевской шкале в случае однородных нормальных краевых условий,
доказана Ю.~М.~Березанским, С.~Г.~Крейном и Я.~А.~Ройтбергом в статье
[\ref{BerezanskyKreinRoitberg63}, с. 747]. Там же установлена теорема о локальном
повышении гладкости решения эллиптического уравнения вплоть до границы области.
Независимо М.~Шехтер [\ref{Schechter63}] установил соответствующую априорную оценку
решения эллиптического уравнения. Доказательство этих результатов приведено в
монографиях Ю.~М.~Березанского [\ref{Berezansky65}] (гл.~III, \S~6, п. 6, 12) и
Я.~А.~Ройтберга [\ref{Roitberg96}] (п. 5.5, 7.3). См. также обзор М.~С.~Аграновича
[\ref{Agranovich97} с.~85, 86] и учебное пособие Ю.~М.~Березанского, Г.~Ф.~Уса,
З.~Г.~Шефтеля [\ref{BerezanskyUsSheftel90}, c.~538, 553] (в последнем рассмотрен
случай, когда эллиптическое выражение имеет второй порядок и вещественные
коэффициенты).

Обобщение этих результатов на случай краевых условий, не являющихся нормальными,
получено Ю.~В.~Костарчуком и Я.~А.~Ройтбергом [\ref{KostarchukRoitberg73}, с.~274],
а на случай эллиптических систем~--- И.~Я.~Ройтберг и Я.~А.~Ройтбергом
[\ref{RoitbergRoitberg00}, с.~298]; см. также монографию Я.~А.~Ройтберга
[\ref{Roitberg99}] (п. 1.3.7).

В связи с п. 3.4.2 отметим, что интерполяция с числовыми параметрами пространств
Соболева, удовлетворяющих однородным краевым условиям, исследована П.~Гриваром
[\ref{Grisvard67}], Р.~Сили [\ref{Seeley72}], Й.~Лёфстрёмом [\ref{Lefstrem92}]
(результаты первых двух авторов изложены также в монографии Х.~Трибеля
[\ref{Triebel80}, с.~400]).

Все теоремы п 3.4 доказаны нами в [\ref{06UMB4}], за исключением последней теоремы
3.19, которая анонсирована вместе с теоремами 3.15 и 3.16 в [\ref{05UMJ5}].

\medskip

\textbf{К п. 3.5.} Доказанные здесь результаты дают другие эквивалентные определения
пространств Хермандера $H^{s,\varphi}(\Omega)$ в евклидовой области для неполуцелых
$s<0$. В соболевском случае эти определения использовали Ж.-Л.~Лионс и Э.~Мадженес
[\ref{LionsMagenes71}, \ref{Magenes66}, \ref{LionsMagenes62}, \ref{LionsMagenes63}]
при изучении краевых задач. Теоремы этого пункта ранее не публиковались.

\normalsize



\chapter[Неоднородные эллиптические краевые задачи]
{\textbf{Неоднородные \\ эллиптические \\ краевые задачи}}\label{ch4}

\chaptermark{\emph Гл. \ref{ch4}. Неоднородные эллиптические краевые задачи}



\section[Эллиптические краевые задачи в шкале позитивных пространств]
{Эллиптические краевые задачи \\ в шкале позитивных пространств}\label{sec4.1}

\markright{\emph \ref{sec4.1}. Краевые задачи в шкале позитивных пространств}

В этом параграфе мы изучим неоднородные эллиптические краевые задачи в шкалах
позитивных пространств Хермандера. Для этих пространств числовой индекс, задающий
основную гладкость, положителен. В п. \ref{sec4.1.1} и \ref{sec4.1.2}
рассматривается регулярная эллиптическая краевая задача, а в остальных пунктах ---
другие важные классы эллиптических краевых задач.

\subsection[Теоремы о нетеровости и гомеоморфизмах]
{Теоремы о нетеровости \\ и гомеоморфизмах}\label{sec4.1.1}

Рассмотрим неоднородную регулярную эллиптическую краевую задачу \eqref{3.1},
\eqref{3.2}:
\begin{equation}\label{4.1}
L\,u=f\;\;\mbox{в}\;\;\Omega,\quad
B_{j}\,u=g_{j}\;\;\mbox{на}\;\;\Gamma\;\;\mbox{при}\;\; j=1,\ldots,q.
\end{equation}
Свяжем с ней линейное отображение
\begin{equation}\label{4.2}
u\mapsto(Lu,Bu):=(Lu,B_{1}u,\ldots,B_{q}u),\quad u\in
C^{\infty}(\,\overline{\Omega}\,).
\end{equation}
Изучим свойства оператора $(L,B)$, являющегося продолжением по непрерывности этого
отображения в соответствующих парах позитивных пространств Хермандера. Напомним, что
конечномерные пространства $N, N^{+}\subset C^{\infty}(\,\overline{\Omega}\,)$
определены в п. 3.1.


\begin{theorem}\label{th4.1} Для произвольных значений параметров
$s>2q$ и $\varphi\in\mathcal{M}$ отображение \eqref{4.2} продолжается по
непрерывности до ограниченного линейного оператора
\begin{gather}\label{4.3}
(L,B):H^{s,\varphi}(\Omega)\rightarrow
H^{s-2q,\varphi}(\Omega)\oplus\bigoplus_{j=1}^{q}H^{s-m_{j}-1/2,\varphi}(\Gamma)=:\\
\mathcal{H}_{s,\varphi}(\Omega,\Gamma). \notag
\end{gather}
Этот оператор нетеров. Его ядро равно $N$, а область значений состоит из всех
векторов $(f,g_{1},\ldots,g_{q})\in\mathcal{H}_{s,\varphi}(\Omega,\Gamma)$ таких,
что
\begin{equation}\label{4.4}
(f,v)_{\Omega}+\sum_{j=1}^{q}\;(g_{j},C^{+}_{j}v)_{\Gamma}=0\quad\mbox{для
всех}\quad v\in N^{+}.
\end{equation}
Индекс оператора \eqref{4.3} равен $\dim N-\dim N^{+}$ и не зависит от
$s$,~$\varphi$.
\end{theorem}


\textbf{Доказательство.} В соболевском случае $\varphi\equiv1$ и $s\geq2q$ эта
теорема~--- классический результат о разрешимости регулярной эллиптической краевой
задачи [\ref{FunctionalAnalysis72}, с. 169]. Он доказан, например, в монографиях
Ю.~М.~Березанского [\ref{Berezansky65}] (гл.~3, \S\,6) и Ж.-Л.~Лионса, Э.~Мадженеса
[\ref{LionsMagenes71}] (ч.~2, п.~5.2). Отсюда общий случай $\varphi\in\mathcal{M}$
получается с помощью интерполяции с функциональным параметром.

А именно, пусть $s>2q$ и $\varepsilon:=(s-2q)/2>0$. Отображение \eqref{4.2}
продолжается по непрерывности до ограниченных нетеровых операторов
\begin{gather}\label{4.5}
(L,B):H^{s\mp\varepsilon}(\Omega)\rightarrow H^{s\mp\varepsilon-2q}(\Omega)\oplus
\bigoplus_{j=1}^{q}H^{s\mp\varepsilon-m_{j}-1/2}(\Gamma)=:\\
\mathcal{H}_{s\mp\varepsilon}(\Omega,\Gamma). \notag
\end{gather}
Они имеют общее ядро $N$, одинаковый индекс $\varkappa:=\dim N-\dim N^{+}$ и области
значений
\begin{gather}\notag
(L,B)(H^{s\mp\varepsilon}(\Omega))=\\
\{\,(f,g_{1},\ldots,g_{q})\in\mathcal{H}_{s\mp\varepsilon}(\Omega,\Gamma):\,\mbox{верно
\eqref{4.4}\,}\}.\label{4.6}
\end{gather}

Применим к \eqref{4.5} интерполяцию с функциональным параметром $\psi$ из теоремы
\ref{th2.14}, где полагаем $\varepsilon=\delta$. Получим ограниченный оператор
\begin{equation}\label{4.7}
(L,B):\,[H^{s-\varepsilon}(\Omega),H^{s+\varepsilon}(\Omega)]_{\psi}\rightarrow
[\mathcal{H}_{s-\varepsilon}(\Omega,\Gamma),\mathcal{H}_{s+\varepsilon}(\Omega,\Gamma)]_{\psi},
\end{equation}
продолжающий по непрерывности отображение \eqref{4.2}. Это в силу интерполяционных
теорем \ref{th2.5}, \ref{th2.21} и \ref{th3.2} означает, что отображение \eqref{4.2}
продолжается по непрерывности до ограниченного оператора \eqref{4.3}, равного
\eqref{4.7}. Согласно теореме \ref{th2.7} нетеровость операторов \eqref{4.5} влечет
нетеровость оператора \eqref{4.3}, который наследует их ядро $N$ и индекс
$\varkappa=\dim N-\dim N^{+}$. Кроме того, область значений оператора \eqref{4.3}
равна $\mathcal{H}_{s,\varphi}(\Omega,\Gamma)\cap(L,B)(H^{s-\varepsilon}(\Omega))$.
Отсюда, в силу \eqref{4.6} получаем, что она такая как в формулировке настоящей
теоремы.

Теорема \ref{th4.1} доказана.


В силу теоремы \ref{th4.1} для произвольной функции $u\in H^{s,\varphi}(\Omega)$,
где $s>2q$ и $\varphi\in\mathcal{M}$, определены по замыканию правые части $f\in
H^{s-2q,\varphi}(\Omega)$ и $g_{j}\in H^{s-m_{j}-1/2,\varphi}(\Gamma)$ краевой
задачи \eqref{4.1}.

Если $N=N^{+}=\{0\}$ (дефект краевой задачи отсутствует), то оператор \eqref{4.3}
является гомеоморфизмом пространства $H^{s,\varphi}(\Omega)$ на
$\mathcal{H}_{s,\varphi}(\Omega,\Gamma)$. Это следует из теоремы \ref{th4.1} и
теоремы Банаха об обратном операторе. В общем случае гомеоморфизм удобно задавать с
помощью следующих проекторов.

Пусть $s>2q$ и $\varphi\in\mathcal{M}$. Представим пространства, в которых действует
оператор \eqref{4.3}, в виде прямых сумм (замкнутых) подпространств:
\begin{gather}\label{4.8}
H^{s,\varphi}(\Omega)=N\dotplus\{u\in
H^{s,\varphi}(\Omega):(u,w)_{\Omega}=0\;\mbox{для всех}\;w\in N\},\\
\mathcal{H}_{s,\varphi}(\Omega,\Gamma)=\{(v,0,\ldots,0):\,v\in N^{+}\}\dotplus
(L,B)(H^{s,\varphi}(\Omega)). \label{4.9}
\end{gather}
Такие разложения в прямые суммы существуют. В самом деле, равенство \eqref{4.8}
является сужением разложения пространства $L_{2}(\Omega)$ в ортогональную сумму
подпространства $N$ и его дополнения. Равенство \eqref{4.9} следует из теоремы
\ref{th4.1}, согласно которой подпространства в его правой части имеют тривиальное
пересечение, и (конечная) размерность первого из них совпадает с коразмерностью
второго.

Обозначим через $P$ и $Q^{+}$ косые проекторы соответственно пространств
$H^{s,\varphi}(\Omega)$ и $\mathcal{H}_{s,\varphi}(\Omega,\Gamma)$ на вторые
слагаемые в суммах \eqref{4.8} и \eqref{4.9} параллельно первым слагаемым. Эти
проекторы не зависят от $s$ и $\varphi$.


\begin{theorem}\label{th4.2} Для произвольных параметров
$s>2q$ и $\varphi\in\mathcal{M}$ сужение отображения $\eqref{4.3}$ на
подпространство $P(H^{s,\varphi}(\Omega))$ является гомеоморфизмом
\begin{equation}\label{4.10}
(L,B):\,P(H^{s,\varphi}(\Omega))\leftrightarrow
Q^{+}(\mathcal{H}_{s,\varphi}(\Omega,\Gamma)).
\end{equation}
\end{theorem}


\textbf{Доказательство.} Согласно теореме \ref{th4.1}, $N$ --- ядро, а
$Q^{+}(\mathcal{H}_{s,\varphi}(\Omega,\Gamma))$
--- область значений оператора \eqref{4.3}. Следовательно, ограниченный оператор
\eqref{4.10} --- биекция. Значит, он гомеоморфизмом в силу теоремы Банаха об
обратном операторе. Теорема~\ref{th4.2} доказана.


Из теоремы \ref{th4.2} вытекает следующая априорная оценка решения эллиптической
краевой задачи \eqref{4.1}.


\begin{theorem}\label{th4.3} Пусть $s>2q$ и $\varphi\in\mathcal{M}$.
Предположим, что функция $u\in H^{s,\varphi}(\Omega)$ является решением краевой
задачи \eqref{4.1}, где
$(f,g_{1},\ldots,g_{q})\in\mathcal{H}_{s,\varphi}(\Omega,\Gamma)$. Тогда выполняется
оценка
\begin{equation}\label{4.11}
\|u\|_{H^{s,\varphi}(\Omega)}\leq
c\,\bigl(\,\|(f,g_{1},\ldots,g_{q})\|_{\mathcal{H}_{s,\varphi}(\Omega,\Gamma)}+
\|u\|_{L_{2}(\Omega)}\,\bigr),
\end{equation}
в которой число $c=c(s,\varphi)>0$ не зависит от $u$ и $(f,g_{1},\ldots,g_{q})$.
\end{theorem}


\textbf{Доказательство.} Воспользуемся разложением \eqref{4.8} и запишем функцию
$u\in H^{s,\,\varphi}(\Omega)$ в виде $u=u_{0}+u_{1}$, где $u_{0}:=(1-P)u\in N$ и
$u_{1}:=Pu\in P(H^{s,\varphi}(\Omega))$. В силу теоремы \ref{th4.2}
\begin{gather*}
\|u_{1}\|_{H^{s,\varphi}(\Omega)}\leq
c_{1}\|(L,B)u_{1}\|_{\mathcal{H}_{s,\varphi}(\Omega,\Gamma)}=\\
c_{1}\|(L,B)u\|_{\mathcal{H}_{s,\varphi}(\Omega,\Gamma)}=
c_{1}\|(f,g_{1},\ldots,g_{q})\|_{\mathcal{H}_{s,\varphi}(\Omega,\Gamma)}.
\end{gather*}
Здесь $c_{1}$ --- норма оператора, обратного к \eqref{4.10}. Кроме того, так как $N$
конечномерно и $1-P$ суть ортопроектор на $N$ в $L_{2}(\Omega)$, то
$$
\|u_{0}\|_{H^{s,\varphi}(\Omega)}\leq c_{0}\|u_{0}\|_{L_{2}(\Omega)}\leq
c_{0}\|u\|_{L_{2}(\Omega)}.
$$
Здесь число $c_{0}>0$ не зависит от $u$ и $(f,g_{1},\ldots,g_{q})$. Сложив эти
неравенства, получаем \eqref{4.11}. Теорема \ref{th4.3} доказана.


Если $N=\{0\}$, т.~е.  краевая задача \eqref{4.1} имеет не более одного решения, то
слагаемое $\|u\|_{L_{2}(\Omega)}$ в правой части априорной оценки \eqref{4.11} можно
опустить.

\subsection{Гладкость решения вплоть до границы}\label{sec4.1.2}

Предположим, что правая часть эллиптической краевой задачи \eqref{4.1} имеет на
открытом в $\overline{\Omega}$ множестве некоторую гладкость в позитивной уточненной
шкале. Изучим гладкость решения $u$ на этом множестве. Рассмотрим сначала случай
глобальной гладкости на всей замкнутой области $\overline{\Omega}$.


\begin{theorem}\label{th4.4} Предположим, что
функция $u\in H^{2q}(\Omega)$ является решением краевой задачи \eqref{4.1}, где
$$
f\in H^{s-2q,\varphi}(\Omega)\quad\mbox{и}\quad g_{j}\in
H^{s-m_{j}-1/2,\varphi}(\Gamma)\quad\mbox{при}\quad j=1,\ldots,q
$$
для некоторых параметров $s>2q$ и $\varphi\in\mathcal{M}$. Тогда $u\in
H^{s,\varphi}(\Omega)$.
\end{theorem}


\textbf{Доказательство.} Согласно теореме \ref{th4.1}, справедливой в соболевском
случае $s=2q$ и $\varphi\equiv1$, вектор $F:=(f,g_{1},\ldots,g_{q})\in
\mathcal{H}_{s,\varphi}(\Omega,\Gamma)$ удовлетворяет условию \eqref{4.4}.
Следовательно, в силу той же теоремы, $F\in(L,B)(H^{s,\varphi}(\Omega))$. Поэтому
наряду с условием $(L,B)u=F$ справедливо равенство $(L,B)v=F$ для некоторого $v\in
H^{s,\varphi}(\Omega)$. Отсюда, $(L,B)(u-v)=0$, что ввиду теоремы \ref{th4.1} влечет
включение
$$
w:=u-v\in N\subset C^{\infty}(\,\overline{\Omega}\,)\subset H^{s,\varphi}(\Omega).
$$
Таким образом, $u=v+w\in H^{s,\varphi}(\Omega)$. Теорема \ref{th4.4} доказана.


Рассмотрим теперь случай локальной гладкости. Пусть $U$ --- открытое множество в
$\mathbb{R}^{n}$, имеющее непустое пересечение с областью $\Omega$. Положим
$\Omega_{0}:=U\cap\Omega$ и $\Gamma_{0}:=U\cap\Gamma$ (возможен случай
$\Gamma_{0}=\varnothing$). Для произвольных параметров $\sigma\in\mathbb{R}$ и
$\varphi\in\mathcal{M}$ введем локальный аналог пространства Хермандера в евклидовой
области:
\begin{gather*}
H^{\sigma,\varphi}_{\mathrm{loc}}(\Omega_{0},\Gamma_{0}):=\{u\in\mathcal{D}'(\Omega): \\
\chi\,u\in H^{\sigma,\varphi}(\Omega)\;\;\mbox{для всех}\;\;\chi\in
C^{\infty}(\,\overline{\Omega}\,),\;\mathrm{supp}\,\chi\subset\Omega_{0}\cup\Gamma_{0}\}.
\end{gather*}
Напомним, что $\mathcal{D}'(\Omega)$~--- топологическое линейное пространство всех
распределений, заданных в области $\Omega$. Нам понадобится также локальное
пространство $H^{\sigma,\varphi}_{\mathrm{loc}}(\Gamma_{0})$, введенное ранее в
п.~\ref{sec2.6.3b}. Как и прежде, в случае $\varphi\equiv1$ индекс $\varphi$ в
обозначении этих пространств будем опускать.

Отметим включения
\begin{equation}\label{4.12}
H^{\sigma,\varphi}(\Omega)\subset
H^{\sigma,\varphi}_{\mathrm{loc}}(\Omega_{0},\Gamma_{0}) \quad\mbox{и}\quad
H^{\sigma,\varphi}(\Gamma)\subset H^{\sigma,\varphi}_{\mathrm{loc}}(\Gamma_{0}).
\end{equation}
Они следует из того, что умножение на функцию из класса
$C^{\infty}(\,\overline{\Omega}\,)$ (соответственно, из $C^{\infty}(\Gamma)$)
является ограниченным оператором в пространстве $H^{\sigma,\varphi}(\Omega)$ (в
$H^{\sigma,\varphi}(\Gamma)$). В соболевском случае $\varphi\equiv1$ этот факт
известен [\ref{Agranovich97}, с. 13]; отсюда случай произвольного
$\varphi\in\mathcal{M}$ вытекает в силу интерполяционных теорем \ref{th2.21} и
\ref{th3.2}.


\begin{theorem}\label{th4.5} Предположим, что
функция $u\in H^{2q}(\Omega)$ является решением краевой задачи \eqref{4.1} задачи,
где
\begin{equation}\label{4.13}
f\in H^{s-2q,\varphi}_{\mathrm{loc}}(\Omega_{0},\Gamma_{0}),\;\;g_{j}\in
H^{s-m_{j}-1/2,\varphi}_{\mathrm{loc}}(\Gamma_{0}),\;\,j=1,\ldots,q,
\end{equation}
для некоторых параметров $s>2q$ и $\varphi\in\mathcal{M}$. Тогда решение $u\in
H^{s,\varphi}_{\mathrm{loc}}(\Omega_{0},\Gamma_{0})$.
\end{theorem}


\textbf{Доказательство.} В соболевском случае $\varphi\equiv1$ эта теорема известна;
см., например, доказательство в [\ref{Berezansky65}], гл.~III, \S~4. Следовательно,
так как
\begin{gather*}
f\in H^{s-2q,\varphi}_{\mathrm{loc}}(\Omega_{0},\Gamma_{0})\subset
H^{s-\varepsilon-2q}_{\mathrm{loc}}(\Omega_{0},\Gamma_{0}),\\
g_{j}\in H^{s-m_{j}-1/2,\varphi}_{\mathrm{loc}}(\Gamma_{0})\subset
H^{s-\varepsilon-m_{j}-1/2}_{\mathrm{loc}}(\Gamma_{0}),
\end{gather*}
то $u\in H^{s-\varepsilon}_{\mathrm{loc}}(\Omega_{0},\Gamma_{0})$, где
$\varepsilon:=\min\{(s-2q),\,1/2\}$. Выведим из последнего включения и условия
\eqref{4.13} требуемое свойство $u\in
H^{s,\varphi}_{\mathrm{loc}}(\Omega_{0},\Gamma_{0})$.

Пусть функции $\chi,\eta\in C^{\infty}(\,\overline{\Omega}\,)$ такие, что
$\mathrm{supp}\,\chi,\,\mathrm{supp}\,\eta\subset\Omega_{0}\cup\Gamma_{0}$ и
$\eta=\nobreak1$ в окрестности $\mathrm{supp}\,\chi$. Переставив оператор умножения
на функцию $\chi$ с дифференциальными операторами $L$ и $B_{j}$, где $j=1,\ldots,q$,
можем записать для произвольного $v\in C^{\infty}(\,\overline{\Omega}\,)$ следующие
равенства:
\begin{gather} \label{4.14}
L(\chi v)=L(\chi\eta v)=\chi\,L(\eta v)+L'(\eta v)= \chi\,Lv+L'(\eta v),\\
B_{j}(\chi v)=B_{j}(\chi\eta v)=\chi\,B_{j}(\eta v)+B'_{j}(\eta
v)=\chi\,B_{j}v+B'_{j}(\eta v). \label{4.15}
\end{gather}
Здесь $L'$ --- некоторое линейное дифференциальное выражение на $\overline{\Omega}$,
а $B'_{j}$~--- некоторое граничное линейное дифференциальное выражение на $\Gamma$.
Коэффициенты этих выражений бесконечно гладкие, а порядки удовлетворяют условиям
$\mathrm{ord}\,L'\leq2q-1$ и $\mathrm{ord}\,B'_{j}\leq m_{j}-1$. Следовательно
[\ref{Agranovich97}, c.~16], отображения $v\mapsto L'v$ и $v\mapsto B'_{j}v$, где
$v\in C^{\infty}(\overline{\Omega})$, продолжаются по непрерывности до ограниченных
операторов
\begin{gather} \label{4.16}
L':H^{\lambda}(\Omega)\rightarrow H^{\lambda-2q+1}(\Omega)\;\,\mbox{для
любого}\;\,\lambda\in\mathbb{R}, \\
B'_{j}:H^{\lambda}(\Omega)\rightarrow H^{\lambda-m_{j}+1/2}(\Gamma)\;\,\mbox{для
любого}\;\,\lambda>m_{j}-1/2. \label{4.17}
\end{gather}
Отсюда при $\lambda:=2q$ вытекает, что равенства \eqref{4.14}, \eqref{4.15}
продолжаются по непрерывности с класса функций $v\in C^{\infty}(\overline{\Omega})$
на класс функций $v\in H^{2q}(\Omega)$. Возьмем в этих равенствах $v:=u\in
H^{2q}(\Omega)$, где $u$ --- решение задачи \eqref{4.1}. Получим соотношения
\begin{gather} \label{4.18}
L(\chi u)=\chi\,f+L'(\eta u)\;\;\mbox{в}\;\;\Omega, \\
B_{j}(\chi u)=\chi\,g_{j}+B'_{j}(\eta u)\;\;\mbox{на}\;\;\Gamma\;\;\mbox{для
любого}\;\;j=1,\ldots,q. \label{4.19}
\end{gather}

Так как $\eta u\in H^{s-\varepsilon}(\Omega)$, то в силу \eqref{4.16}, \eqref{4.17},
где $\lambda:=s-\varepsilon$, и определения числа $\varepsilon$ имеем:
\begin{gather*}
L'(\eta u)\in H^{s-\varepsilon-2q+1}(\Omega)\subset H^{s-2q,\varphi}(\Omega),\\
B'_{j}(\eta u)\in H^{s-\varepsilon-m_{j}+1/2}(\Gamma)\subset
H^{s-m_{j}-1/2,\varphi}(\Gamma).
\end{gather*}
Отсюда и из условия \eqref{4.13} вытекает (см. также \eqref{4.12}), что функция
$\chi u\in H^{2q}(\Omega)$ является решением эллиптической краевой задачи
\eqref{4.18}, \eqref{4.19}, правые части которой имеют следующую гладкость:
\begin{gather*}
\chi\,f+L'(\eta u)\in H^{s-2q,\varphi}(\Omega),\\
\chi\,g_{j}+B'_{j}(\eta u)\in H^{s-m_{j}-1/2,\varphi}(\Gamma).
\end{gather*}
Поэтому согласно теореме \ref{th4.4} верно включение $\chi u\in
H^{s,\varphi}(\Omega)$. Оно ввиду произвольности функции $\chi\in
C^{\infty}(\,\overline{\Omega}\,)$, удовлетворяющей условию
$\mathrm{supp}\,\chi\subset\Omega_{0}\cup\Gamma_{0}$, означает, что $u\in
H^{s,\varphi}_{\mathrm{loc}}(\Omega_{0},\Gamma_{0})$.

Теорема  \ref{th4.5} доказана.


В теореме \ref{th4.5} отметим случай $\Gamma_{0}=\varnothing$, который приводит к
утверждению о повышении локальной гладкости решения в окрестностях внутренних точек
области $\Omega$.

В качестве приложения теорем \ref{th4.4} и \ref{th4.5} установим одно достаточное
условие того, что решение $u$ эллиптической краевой задачи \eqref{4.1} является
\emph{классическим}, т.~е. принадлежит классу $C^{2q}(\Omega)\cap
C^{m}(\,\overline{\Omega}\,)$, где $m:=\max\{m_{1},\ldots,m_{q}\}$. Если $u$ ---
такое решение, то левые части равенств в \eqref{4.1} вычисляются с помощью
классических производных, а сами равенства выполняются в каждой точке множества
$\Omega$ либо $\Gamma$ соответственно. При этом их правые части имеют следующую
гладкость:
\begin{equation}\label{4.20}
f\in C(\Omega)\;\;\mbox{и}\;\;g_{j}\in C^{m-m_{j}}(\Gamma)\;\;\mbox{для
каждого}\;\;j=1,\ldots,q.
\end{equation}
Обратное, вообще говоря, не верно: из условия \eqref{4.20} не следует, что решение
$u$ является классическим [\ref{GilbargTrudinger89}, с.~73]. С помощью пространств
уточненной шкалы мы так усилим это условие, что оно станет достаточным для
классичности решения~$u$.


\begin{theorem}\label{th4.6}
Предположим, что функция $u\in H^{2q}(\Omega)\cap H^{2q,\varphi}(\Omega)$ является
решением задачи \eqref{4.1}, где
\begin{gather}\label{4.21}
f\in H^{n/2,\varphi}_{\mathrm{loc}}(\Omega,\varnothing)\cap
H^{m-2q+n/2,\varphi}(\Omega), \\
g_{j}\in H^{m-m_{j}+(n-1)/2,\varphi}(\Gamma)\quad\mbox{для всех}\quad j=1,\ldots,q,
\label{4.22}
\end{gather}
а функциональный параметр $\varphi\in\mathcal{M}$ удовлетворяет условию
\eqref{2.37}. Тогда решение $u$ классическое, т.е. $u\in C^{2q}(\Omega)\cap
C^{m}(\,\overline{\Omega}\,)$.
\end{theorem}


\begin{remark}\label{rem4.1}
Условия \eqref{4.21} и \eqref{4.22} влекут свойство \eqref{4.20}. Это следует из
теорем \ref{th2.27} и \ref{th3.4}. Отметим также, что в теореме \ref{th4.6} условие
\eqref{2.37} не только достаточное для классичности решения $u$, но также и
необходимое на классе всех рассматриваемых решений эллиптической краевой задачи; см.
также замечание~\ref{rem2.8c}.
\end{remark}


\textbf{Доказательство теоремы 4.6.} Ввиду условия \eqref{4.21} имеем для
$s:=2q+n/2$ включения
\begin{gather*}
f\in H^{s-2q,\varphi}_{\mathrm{loc}}(\Omega,\varnothing),\\
g_{j}\in\mathcal{D}'(\Gamma)=
H^{s-m_{j}-1/2,\varphi}_{\mathrm{loc}}(\varnothing)\quad\mbox{при}\quad
j=1,\ldots,q.
\end{gather*}
Отсюда на основании теорем \ref{th4.5} и \ref{th3.4} получаем:
$$
u\in H^{s,\varphi}_{\mathrm{loc}}(\Omega,\varnothing)=
H^{2q+n/2,\varphi}_{\mathrm{loc}}(\Omega,\varnothing)\subset C^{2q}(\Omega).
$$

Докажем также, что $u\in C^{m}(\,\overline{\Omega}\,)$. Если $s:=m+n/2>2q$, то ввиду
\eqref{4.21}, \eqref{4.22} выполняется условие теоремы \ref{th4.4}. Следовательно, в
силу теорем \ref{th4.4} и \ref{th3.4}
$$
u\in H^{s,\varphi}(\Omega)=H^{m+n/2,\varphi}(\Omega)\subset
C^{m}(\,\overline{\Omega}\,).
$$
Если же $m+n/2\leq2q$, то, по условию,
$$
u\in H^{2q,\varphi}(\Omega)\subseteq H^{m+n/2,\varphi}(\Omega)\subset
C^{m}(\,\overline{\Omega}\,).
$$

Таким образом, $u\in C^{2q}(\Omega)\cap C^{m}(\,\overline{\Omega}\,)$. Теорема
\ref{th4.6} доказана.


\begin{remark}\label{rem4.2}
Если использовать теорему \ref{th4.6} лишь для соболевской шкалы пространств, то
придется вместо \eqref{4.21}, \eqref{4.22} потребовать, чтобы для некоторого числа
$\varepsilon>0$ выполнялись условия
\begin{gather*}
f\in H^{n/2+\varepsilon}_{\mathrm{loc}}(\Omega,\varnothing)\cap
H^{m-2q+n/2+\varepsilon}(\Omega), \\
g_{j}\in H^{m-m_{j}+(n-1)/2+\varepsilon}(\Gamma)\quad\mbox{при}\quad j=1,\ldots,q.
\end{gather*}
Они завышают основную гладкость правых частей краевой задачи \eqref{4.1}, что
огрубляет результат.
\end{remark}

\subsection[Нерегулярная эллиптическая краевая за-\break дача]
{Нерегулярная эллиптическая \\ краевая задача}\label{sec4.1.3}

В этом пункте мы предполагаем, что краевая задача \eqref{4.1} является эллиптической
в области $\Omega$, но нерегулярной. Это означает, что она удовлетворяет п. (i) и
(ii) определения \ref{3.1}, но не удовлетворяет п. (iii) этого определения. Для
такой задачи остаются верными (с некоторыми изменениями) все результаты предыдущих
п. \ref{sec4.1.1}, \ref{sec4.1.2}. Изменения касаются лишь теорем \ref{th4.1},
\ref{th4.2} и вызваны тем, что нерегулярная эллиптическая краевая задача не имеет
формально сопряженной краевой задачи в классе дифференциальных. Приведем аналоги
этих теорем.


\begin{theorem}\label{th4.7}
Для произвольных параметров $s>2q$ и $\varphi\in\mathcal{M}$ отображение \eqref{4.2}
продолжается по непрерывности до ограниченного нетерового оператора \eqref{4.3}.
Ядро этого оператора равно $N$, а область значений состоит из всех векторов
$(f,g_{1},\ldots,g_{q})\in\mathcal{H}_{s,\varphi}(\Omega,\Gamma)$ таких, что
\begin{equation}\label{4.23}
(f,v)_{\Omega}+\sum_{j=1}^{q}\;(g_{j},v_{j})_{\Gamma}=0\quad\mbox{для
всех}\quad(v,v_{1},\ldots,v_{q})\in W.
\end{equation}
Здесь $W$ --- некоторое конечномерное пространство, лежащее в
$C^{\infty}(\,\overline{\Omega}\,)\times(C^{\infty}(\Gamma))^{q}$ и не зависящее от
$s$, $\varphi$. Индекс оператора \eqref{4.3} равен $\dim N-\dim W$ и не зависит от
$s$, $\varphi$.
\end{theorem}


В соболевском случае $\varphi\equiv1$ эта теорема известна (см., например,
[\ref{Hermander87}, с.~315, 324] или [\ref{Agranovich97}, с.~20, 21]). Отсюда общий
случай $\varphi\in\mathcal{M}$ выводится с помощью интерполяции аналогично
доказательству теоремы~\ref{th4.1}.

В силу теоремы \ref{th4.7} можно записать для произвольных $s>2q$ и
$\varphi\in\mathcal{M}$
$$
\mathcal{H}_{s,\varphi}(\Omega,\Gamma)=W\dotplus (L,B)(H^{s,\varphi}(\Omega)).
$$
Обозначим через $Q$ косой проектор пространства
$\mathcal{H}_{s,\varphi}(\Omega,\Gamma)$ на подпространство
$(L,B)(H^{s,\varphi}(\Omega))$ параллельно $W$. Этот проектор не зависит от $s$ и
$\varphi$.


\begin{theorem}\label{th4.8}
Для произвольных параметров $s>2q$ и $\varphi\in\nobreak\mathcal{M}$ сужение
отображения $\eqref{4.3}$ на подпространство $P(H^{s,\varphi}(\Omega))$ является
гомеоморфизмом
$$
(L,B):\,P(H^{s,\varphi}(\Omega))\leftrightarrow
Q(\mathcal{H}_{s,\varphi}(\Omega,\Gamma)).
$$
\end{theorem}


Теорема \ref{th4.8} непосредственно следует из теоремы \ref{th4.7} и теоремы Банаха
об обратном операторе.


\begin{example}\label{ex4.1}
Краевая задача с косой производной для уравнения Лапласа:
$$
\Delta u=f\;\;\mbox{в}\;\;\Omega,\quad\frac{\partial
u}{\partial\tau}=g\;\;\mbox{на}\;\;\Gamma.
$$
Здесь $\tau$ --- бесконечно гладкое поле ортов, касательных к границе~$\Gamma$. Если
$\dim\Omega=2$, то эта задача эллиптическая, но не регулярная. Если
$\dim\Omega\geq3$, то она не является и эллиптической [\ref{Agranovich97}, с.~10].
\end{example}

\begin{example}\label{ex4.1bb}
Пусть $\tau_{1},\ldots,\tau_{n-1}$~--- линейно независимая система бесконечно
гладких полей ненулевых векторов, касательных к границе $\Gamma$. (Напомним, что
$n=\dim\Omega$.) Следующая краевая задача эллиптическая, но не регулярная:
$$
\Delta^{n-1}u =f\;\;\mbox{в}\;\;\Omega,\quad\;\frac{\partial
u}{\partial\tau_{j}}=g_{j}\;\;\mbox{на}\;\;\Gamma\;\;\mbox{при}\;\;j=1,\ldots,n-1.
$$
\end{example}

\begin{example}\label{ex4.1cc}
Пусть $\dim\Omega=2$. Рассмотрим краевую задачу
$$
\Delta^{2}u =f\;\;\mbox{в}\;\;\Omega,\quad\;\;\frac{\partial
u}{\partial\nu}+\frac{\partial u}{\partial\tau}=g_{1},\;\;\frac{\partial
u}{\partial\nu}-\frac{\partial u}{\partial\tau}=g_{2}\;\;\mbox{на}\;\;\Gamma.
$$
Здесь $\tau$~--- бесконечно гладкое поле ортов, касательных к границе~$\Gamma$. Эта
задача эллиптическая, но не регулярная.
\end{example}

Другие примеры эллиптических, но нерегулярных краевых задач приведены, например, в
работе Я.~А.~Ройтберга [\ref{Roitberg69}].


В заключение этого пункта напомним важный факт [\ref{Agranovich97}, с. 20]: если для
краевой задачи \eqref{4.1} ограниченный оператор \eqref{4.3} нетеров для некоторого
$s>2q$ при условии $\varphi\equiv1$, то эта задача эллиптична в области $\Omega$,
т.~е. удовлетворяет п. (i) и (ii) определения \ref{3.1}.

\subsection[Эллиптическая краевая задача с парамет-\break ром]
{Эллиптическая краевая задача \\ с параметром}\label{sec4.1.4}

В работах Ш.~Агмона, Л.~Ниренберга [\ref{Agmon62}, \ref{AgmonNirenberg63}] и
М.~С.~Аграновича, М.~И.~Вишика [\ref{AgranovichVishik64}] (см. также обзор
[\ref{Agranovich97}], \S~4) выделен подкласс эллиптических краевых задач~---
эллиптические краевые задачи с параметром, обладающие следующим важным свойством.
При достаточно больших по модулю значениях комплексного параметра оператор,
соответствующий задаче, является гомеоморфизмом в подходящих парах соболевских
пространств, причем норма оператора допускает двустороннюю оценку с постоянными, не
зависящими от параметра. В этом пункте будет доказано, что для пространств
уточненной шкалы справедлив аналог этого результата.

Приведем определение эллиптической краевой задачи с параметром. Рассмотрим
неоднородную краевую задачу
\begin{gather} \label{4.24}
L(\lambda)\,u=f\quad\mbox{в}\quad\Omega,\\
B_{j}(\lambda)\,u=g_{j}\quad\mbox{на}\quad\Gamma\quad\mbox{при}\quad j=1,\ldots,q,
\label{4.25}
\end{gather}
зависящую от комплексного параметра $\lambda$ следующим образом:
\begin{equation}\label{4.26}
L(\lambda):=\sum_{r=0}^{2q}\,\lambda^{2q-r}L_{r}\quad\;\mbox{и}\quad\;
B_{j}(\lambda):=\sum_{r=0}^{m_{j}}\,\lambda^{m_{j}-r}B_{j,r}.
\end{equation}
Здесь $L_{r}$ --- линейное дифференциальное выражение на $\overline{\Omega}$, а
$B_{j,r}$ --- граничное линейное дифференциальное на $\Gamma$; коэффициенты этих
выражений --- бесконечно гладкие комплекснозначные функции, а порядки не превышают
числа $r$. Как и прежде, фиксированные целые числа $q$ и $m_{j}$ удовлетворяют
условиям $q\geq1$ и $0\leq m_{j}\leq 2q-1$. Отметим, что $L(0)=L_{2q}$ и
$B_{j}(0)=B_{j,m_{j}}$.

Сопоставим дифференциальным выражениям \eqref{4.26} некоторые однородные полиномы от
$(\xi,\lambda)\in\mathbb{C}^{n+1}$. Положим
$$
L^{(0)}(x;\xi,\lambda):=\sum_{r=0}^{2q}\,\lambda^{2q-r}L^{(0)}_{r}(x,\xi)
\quad\mbox{при}\quad
x\in\overline{\Omega},\;\xi\in\mathbb{C}^{n},\;\lambda\in\mathbb{C}.
$$
Здесь $L^{(0)}_{r}(x,\xi)$ --- главный символ выражения $L_{r}$ в случае, когда
$\mathrm{ord}\,L_{r}=r$, либо $L^{(0)}_{r}(x,\xi)\equiv0$ в случае, когда
$\mathrm{ord}\,L_{r}<r$. Аналогично для каждого $j=1,\ldots,q$ положим
$$
B^{(0)}_{j}(x;\xi,\lambda):=\sum_{r=0}^{m_{j}}\,\lambda^{m_{j}-r}B^{(0)}_{j,r}(x,\xi)
\quad\mbox{при}\quad x\in\Gamma,\;\xi\in \mathbb{C}^{n},\;\lambda\in\mathbb{C}.
$$
Здесь $B^{(0)}_{j,r}(x,\xi)$ --- главный символ граничного дифференциального
выражения $B_{j,r}$ в случае $\mathrm{ord}\,B_{j,r}=r$, либо
$B^{(0)}_{j,r}(x,\xi)\equiv0$ в случае $\mathrm{ord}\,B_{j,r}<r$. Отметим, что
$L^{(0)}(x;\xi,\lambda)$ и $B^{(0)}_{j}(x;\xi,\lambda)$ --- однородные полиномы от
$(\xi,\lambda)\in\mathbb{C}^{n+1}$ степеней $2q$ и $m_{j}$ соответственно.

Пусть $K$ --- фиксированный замкнутый угол на комплексной плоскости с вершиной в
начале координат (не исключается случай, когда $K$ вырождается в луч).


\begin{definition}\label{def4.1}
Краевая задача \eqref{4.24}, \eqref{4.25} называется \emph{эллиптической с
параметром} в угле $K$, если выполняются следующие условия:
\begin{itemize}
\item [(i)] Для любых $x\in\overline{\Omega}$,
$\xi\in\mathbb{R}^{n}$ и $\lambda\in K$ таких, что $|\xi|+|\lambda|\neq0$,
справедливо $L^{(0)}(x;\xi,\lambda)\neq0$.
\item [(ii)] Для любых фиксированных точки $x\in\Gamma$, вектора
$\xi\in\nobreak\mathbb{R}^{n}$, касательного к границе $\Gamma$ в точке $x$, и
параметра $\lambda\in K$ таких, что $|\xi|+|\lambda|\neq0$, многочлены
$B^{(0)}_{j}(x;\xi+\tau\nu(x),\lambda)$, $j=1,\ldots,q$, переменного $\tau$ линейно
не\-зависимы по модулю многочлена
$\prod_{j=1}^{q}(\tau-\tau^{+}_{j}(x;\xi;\lambda))$. Здесь
$\tau^{+}_{1}(x;\xi;\lambda),\ldots,\tau^{+}_{q}(x;\xi;\lambda)$~--- все
$\tau$-корни полинома $L^{(0)}(x;\xi+\tau\nu(x),\lambda)$, имеющие положительную
мнимую часть и записанные с учетом их кратности.
\end{itemize}
\end{definition}


\begin{remark}\label{rem4.3}
Условие (ii) определения \ref{def4.1} сформулировано корректно в том смысле, что
полином $L^{(0)}(x;\xi+\tau\nu(x),\lambda)$ имеет ровно $q$ $\tau$-корней с
положительной мнимой частью и столько же корней с отрицательной мнимой частью (с
учетом их кратности). В самом деле, из условия (i) следует, что для каждой точки
$x\in\overline{\Omega}$ дифференциальное выражение
$$
L(x;D,D_{t}):=\sum_{r=0}^{2q}\,D_{t}^{2q-r}L_{r}(x,D)
$$
является эллиптическим. Поскольку оно содержит операторы дифференцирования по
$n+1\geq3$ вещественной переменной $x_{1},\ldots,x_{n},t$, из его эллиптичности
следует правильная эллиптичность [\ref{FunctionalAnalysis72}, с.~166]. Поэтому
$\tau$-корни полинома $L^{(0)}(x;\xi+\tau\nu(x),\lambda)$ имеют указанное свойство.
\end{remark}


Приведем некоторые примеры эллиптических с параметром краевых задач
[\ref{Agranovich97}, с. 24].


\begin{example}\label{ex4.2}
Пусть дифференциальное выражение $L(\lambda)$ удовлетворяет условию (i) определения
\ref{def4.1}. Тогда краевая задача Дирихле для уравнения $L(\lambda)=f$ является
эллиптической с параметром в угле $K$. Здесь краевые условия не зависят от параметра
$\lambda$.
\end{example}


\begin{example}\label{ex4.3}
Краевая задача
$$
\Delta u+\lambda^{2}u=f\;\;\mbox{в}\;\;\Omega,\quad\frac{\partial
u}{\partial\nu}-\lambda u=g\;\;\mbox{на}\;\;\Gamma
$$
является эллиптической с параметром в каждом угле
$$
K_{\varepsilon}:=
\{\lambda\in\mathbb{C}:\,\varepsilon\leq|\mathrm{\arg\lambda}|\leq\pi-\varepsilon\},
$$
где $0<\varepsilon<\pi/2$, а комплексная плоскость разрезана вдоль отрицательной
полуоси.
\end{example}


Далее в этом пункте предполагается, что краевая задача \eqref{4.24}, \eqref{4.25}
эллиптическая с параметром в угле $K$.

При $\lambda=0$ условие (i) определения \ref{def4.1} влечет эллиптичность на
$\overline{\Omega}$ и правильную эллиптичность на $\Gamma$ дифференциального
выражения $L(0)$ (см. замечание \ref{rem4.3}). Условие (ii) означает, что набор
$\{B_{1}(0),\ldots,B_{q}(0)\}$ граничных дифференциальных выражений удовлетворяет
условию дополнительности по отношению к $L(0)$ на $\Gamma$. Следовательно, при
$\lambda=0$ краевая задача \eqref{4.24}, \eqref{4.25} является эллиптической (не
обязательно регулярной) в области $\overline{\Omega}$. Поскольку параметр $\lambda$
влияет только на младшие члены дифференциальных выражений $L(\lambda)$ и
$B_{j}(\lambda)$, эта задача будет эллиптической при всех $\lambda\in\mathbb{C}$.
Согласно теореме~\ref{4.7} отображение
$$
u\mapsto(L(\lambda)u,B(\lambda)u):=(L(\lambda)u,B_{1}(\lambda)u,\ldots,B_{q}(\lambda)u),\quad
u\in C^{\infty}(\,\overline{\Omega}\,),
$$
продолжается по непрерывности до ограниченного нетерового оператора
\begin{equation}\label{4.27}
(L(\lambda),B(\lambda)):\,H^{s,\varphi}(\Omega)\rightarrow
\mathcal{H}_{s,\varphi}(\Omega,\Gamma)
\end{equation}
для произвольных параметров $s>2q$, $\varphi\in\mathcal{M}$ и
$\lambda\in\mathbb{C}$. Его индекс не зависит как от $s$, $\varphi$, так и от
$\lambda$, поскольку $\lambda$ влияет только на младшие члены [\ref{Hermander87},
с.~324] (теорема 20.1.8).

Поскольку краевая задача \eqref{4.24}, \eqref{4.25} эллиптическая с параметром в
угле $K$, оператор \eqref{4.27} имеет следующие важные свойства.


\begin{theorem}\label{th4.9}
Справедливы следующие утверждения.
\begin{itemize}
\item [$\mathrm{(i)}$] Существует число $\lambda_{0}>0$ такое, что для каждого
значения параметра $\lambda\in K$, удовлетворяющего условию
$|\lambda|\geq\nobreak\lambda_{0}$, и произвольных $s>2q$, $\varphi\in\mathcal{M}$
оператор \eqref{4.27} является гомеоморфизмом пространства $H^{s,\varphi}(\Omega)$
на пространство $\mathcal{H}_{s,\varphi}(\Omega,\Gamma)$.
\item [$\mathrm{(ii)}$] Для произвольных фиксированных параметров
$s>2q$ и $\varphi\in\mathcal{M}$ найдется число $c=c(s,\varphi)\geq1$ такое, что для
каждого значения параметра $\lambda\in K$, $|\lambda|\geq\max\{\lambda_{0},1\}$, и
любой функции $u\in H^{s,\varphi}(\Omega)$ выполняется двусторонняя оценка
\begin{gather}
c^{-1}\bigl(\,\|u\|_{H^{s,\varphi}(\Omega)}+
|\lambda|^{s}\varphi(|\lambda|)\,\|u\|_{L_{2}(\Omega)}\,\bigr)\leq \notag \\
\|L(\lambda)u\|_{H^{s-2q,\varphi}(\Omega)}+
|\lambda|^{s-2q}\varphi(|\lambda|)\,\|L(\lambda)u\|_{L_{2}(\Omega)}+ \notag
\\ \notag
\sum_{j=1}^{q}\,\bigl(\,\|B_{j}(\lambda)u\|_{H^{s-m_{j}-1/2,\varphi}(\Gamma)}+ \\
\notag
|\lambda|^{s-m_{j}-1/2}\varphi(|\lambda|)\,\|B_{j}(\lambda)u\|_{L_{2}(\Gamma)}\,\bigr)\leq \\
c\,\bigl(\,\|u\|_{H^{s,\varphi}(\Omega)}+
|\lambda|^{s}\varphi(|\lambda|)\,\|u\|_{L_{2}(\Omega)}\,\bigr). \label{4.28}
\end{gather}
Здесь число $c$ не зависит от $u$ и $\lambda$.
\end{itemize}
\end{theorem}


\begin{remark}\label{rem4.4}
Утверждение (ii) теоремы \ref{th4.9} нуждается в комментарии. При фиксированном
$\lambda$ оценка \eqref{4.28} записана для норм, эквивалентных нормам
$\|u\|_{H^{s,\varphi}(\Omega)}$ и
$\|(L(\lambda),B(\lambda))u\|_{\mathcal{H}_{s,\varphi}(\Omega,\Gamma)}$. Чтобы
избежать громоздких выражений мы написали эту оценку для негильбертовых норм. Она
справедлива и для соответствующих гильбертовых норм (порождающих скалярные
произведения в $H^{s,\varphi}(\Omega)$ и $\mathcal{H}_{s,\varphi}(\Omega,\Gamma)$),
поскольку они оцениваются через использованные нормы с постоянными, не зависящими от
$s$, $\varphi$ и $\lambda$. Дополнительное условие $|\lambda|\geq1$ вызвано тем, что
функция $\varphi(t)$ определена лишь при $t\geq1$. Заметим, что оценка \eqref{4.28}
представляет интерес лишь при $|\lambda|\gg1$.
\end{remark}


В соболевском случае $\varphi\equiv1$ и $s\geq2q$ теорема \ref{th4.9} доказана
М.~С.~Аграновичем и М.~И.~Вишиком [\ref{AgranovichVishik64}] (\S~4, 5),
[\ref{Agranovich97}, c.~25]. Двусторонняя априорная оценка \eqref{4.28} имеет ряд
приложений, в частности, в теории параболических задач [\ref{AgranovichVishik64}].
Отметим [\ref{AgranovichVishik64}, с.~74], что правая часть оценки \eqref{4.28}
справедлива и без предположения об эллиптичности с параметром задачи \eqref{4.24},
\eqref{4.25}.

Докажем отдельно утверждения (i) и (ii) теоремы \ref{th4.9}.

\medskip

\textbf{Доказательство утверждения (i) теоремы \ref{th4.9}.} Число $\lambda_{0}>0$
берем из формулировки этой теоремы в случае $\varphi\equiv1$. Пусть $\lambda\in K$,
$|\lambda|\geq\lambda_{0}$, и $s>2q$, $\varphi\in\mathcal{M}$. Положим
$\varepsilon:=(s-2q)/2>0$. Имеем гомеоморфизмы
$$
(L(\lambda),B(\lambda)):\,H^{s\mp\varepsilon}(\Omega)\leftrightarrow
\mathcal{H}_{s\mp\varepsilon}(\Omega,\Gamma).
$$
Применив здесь интерполяцию с функциональным параметром $\psi$ из теоремы
\ref{th2.14}, где $\varepsilon=\delta$, получим в силу интерполяционных теорем
\ref{th2.5}, \ref{th2.21} и \ref{th3.2} гомеоморфизм
$$
(L(\lambda),B(\lambda)):\,H^{s,\varphi}(\Omega)\leftrightarrow
\mathcal{H}_{s,\varphi}(\Omega,\Gamma).
$$
Утверждение (i) доказано.


Пред тем как доказывать утверждение (ii), полезно ввести некоторые пространства
Хермандера, нормы в которых зависят от дополнительного параметра $\varrho$. Нам
понадобится одно интерполяционное свойство этих пространств.

Пусть $\sigma>0$, $\varphi\in\mathcal{M}$ и $\varrho\geq1$. Предположим, что
$G\in\{\mathbb{R}^{k},\Omega,\Gamma\}$, где $k\in\nobreak\mathbb{N}$. Обозначим
через $H^{\sigma,\varphi}(G,\varrho)$ пространство $H^{\sigma,\varphi}(G)$,
наделенное нормой, зависящей от параметра $\varrho$ следующим образом:
\begin{equation}\label{4.29}
\|u\|_{H^{\sigma,\varphi}(G,\varrho)}:= \bigl(\,\|u\|_{H^{\sigma,\varphi}(G)}^{2}+
\varrho^{2\sigma}\varphi^{2}(\varrho)\|u\| _{L_{2}(G)}^{2}\,\bigl)^{1/2}.
\end{equation}
Эта норма эквивалентна норме в пространстве $H^{\sigma,\varphi}(G)$ при каждом
значении параметра $\varrho$. Поэтому пространство $H^{\sigma,\varphi}(G,\varrho)$
является полным. Оно гильбертово, поскольку норма \eqref{4.29} порождена скалярным
произведением
$$
(u_{1},u_{2})_{H^{\sigma,\varphi}(G,\varrho)}:=
(u_{1},u_{2})_{H^{\sigma,\varphi}(G)}+
\varrho^{2\sigma}\varphi^{2}(\varrho)(u_{1},u_{2})_{L_{2}(\Omega)}.
$$
Как обычно, в соболевском случае $\varphi\equiv1$ мы опускаем индекс $\varphi$ в
обозначениях.

В силу теорем \ref{th2.14}, \ref{th2.21} и \ref{th3.2} пространства
$$
[H^{\sigma-\varepsilon}(G,\varrho),H^{\sigma+\delta}(G,\varrho)]_{\psi}
\quad\mbox{и}\quad H^{\sigma,\varphi}(G,\varrho)
$$
равны с точностью до эквивалентности норм. Оказывается, в оценках норм этих
пространств можно выбрать постоянные так, чтобы они не зависели от параметра
$\varrho$.


\begin{lemma}\label{lem4.1}
Пусть заданы функция $\varphi\in\mathcal{M}$ и положительные числа $\sigma$,
$\varepsilon$, $\delta$, причем $\sigma-\varepsilon>0$. Тогда существует число
$c\geq1$ такое, что для произвольных параметра $\varrho\geq1$ и  функции $u\in
H^{\sigma,\varphi}(G)$ выполняется двусторонняя оценка норм
\begin{equation}\label{4.30}
c^{-1}\,\|u\|_{H^{\sigma,\varphi}(G,\varrho)}\leq\|u\|
_{[H^{\sigma-\varepsilon}(G,\varrho), H^{\sigma+\delta}(G,\varrho)]_{\psi}}\leq
c\,\|u\|_{H^{\sigma,\varphi}(G,\varrho)}.
\end{equation}
Здесь $G\in\{\mathbb{R}^{k},\Omega,\Gamma\}$, $\psi$ --- интерполяционный параметр
из теоремы $\ref{th2.14}$, а число $c$ не зависит как от $\varrho$, так и от $u$.
\end{lemma}


\textbf{Доказательство.} Сначала докажем лемму \ref{lem4.1} для $G=\mathbb{R}^{k}$,
где $k\in\mathbb{N}$. Отсюда мы выведим ее в нужных нам случаях $G=\Omega$ и
$G=\Gamma$.

Случай $G=\mathbb{R}^{k}$. Пусть $\varrho\geq1$ и $u\in
H^{\sigma,\varphi}(\mathbb{R}^{k})$. Воспользовавшись определением норм в
пространствах $H^{\sigma,\varphi}(\mathbb{R}^{k},\varrho)$ и
$H^{\sigma,\varphi}(\mathbb{R}^{k})$, запишем
\begin{gather}\notag
\|u\|_{H^{\sigma,\varphi}(\mathbb{R}^{k},\varrho)}=\\
\biggl(\;\int\limits_{\mathbb{R}^{k}}\,\bigl(\langle\xi\rangle^{2\sigma}
\varphi^{2}(\langle\xi\rangle)+\varrho^{2\sigma}\varphi^{2}(\varrho)\bigr)\,
|\widehat{u}(\xi)|^{2}\,d\xi\,\biggr)^{1/2}.\label{4.31}
\end{gather}
Наряду с \eqref{4.31} рассмотрим еще одну гильбертову норму функции $u$:
\begin{equation}\label{4.32}
\biggl(\;\int\limits_{\mathbb{R}^{k}}\,(\langle\xi\rangle+\varrho)^{2\sigma}\,
\varphi^{2}(\langle\xi\rangle+\varrho)\,
|\widehat{u}(\xi)|^{2}\,d\xi\,\biggr)^{1/2}.
\end{equation}
Нормы \eqref{4.31} и \eqref{4.32} эквивалентные, причем постоянные, с помощью
которых одна норма оценивается через другую, зависят лишь от $\sigma$, $\varphi$ и,
значит, не зависят от параметра $\varrho$. Это следует из леммы \ref{lem4.2},
которая будет установлена сразу после доказательства настоящей леммы. Обозначим
через $H^{\sigma,\varphi}(\mathbb{R}^{k},\varrho,1)$ гильбертово пространство
$H^{\sigma,\varphi}(\mathbb{R}^{k})$, наделенное нормой \eqref{4.32} и
соответствующим скалярным произведением
$$
(u_{1},u_{2})_{H^{\sigma,\varphi}(\mathbb{R}^{k},\varrho,1)}:=
\int\limits_{\mathbb{R}^{k}}(\langle\xi\rangle+\varrho)^{2\sigma}\,
\varphi^{2}(\langle\xi\rangle+\varrho)\,
\widehat{u}_{1}(\xi)\,\overline{\widehat{u}_{2}(\xi)}\,d\xi.
$$
Имеем эквивалентность норм
\begin{equation}\label{4.33}
c_{0}^{-1}\,\|u\|_{H^{\sigma,\varphi}(\mathbb{R}^{k},\varrho)}\leq
\|u\|_{H^{\sigma,\varphi}(\mathbb{R}^{k},\varrho,1)}\leq
c_{0}\,\|u\|_{H^{\sigma,\varphi}(\mathbb{R}^{k},\varrho)}.
\end{equation}
Здесь число $c_{0}=c_{0}(\sigma,\varphi)\geq1$ не зависит от $u$ и $\varrho$.

Проинтерполируем пару $[H^{\sigma-\varepsilon}(\mathbb{R}^{k},\varrho,1),
H^{\sigma+\delta}(\mathbb{R}^{k},\varrho,1)]$ с параметром $\psi$ из теоремы
\ref{th2.14}. Псевдодифференциальный оператор $J_{\varrho}$ с символом
$(\langle\xi\rangle+\varrho)^{\varepsilon+\delta}$ является порождающим для этой
пары. С помощью преобразования Фурье
$$
\mathcal{F}:H^{\sigma-\varepsilon}\bigl(\mathbb{R}^{k},\varrho,1\bigr)\leftrightarrow
L_{2}\bigl(\mathbb{R}^{k},(\langle\xi\rangle+\varrho)^{2(\sigma-\varepsilon)}
d\xi\bigr)
$$
оператор $\psi(J_{\varrho})$ приведен к виду умножения на функцию
$\psi((\langle\xi\rangle+\varrho)^{(\varepsilon+\delta)})=
(\langle\xi\rangle+\varrho)^{\varepsilon}\varphi(\langle\xi\rangle+\varrho)$
аргумента $\xi\in\mathbb{R}^{k}$. Следовательно,
\begin{gather*}
\|u\|^{2} _{[H^{\sigma-\varepsilon}(\mathbb{R}^{k},\varrho,1),
H^{\sigma+\delta}(\mathbb{R}^{k},\varrho,1)]_{\psi}}=
\|\psi(J_{\varrho})u\|_{H^{\sigma-\varepsilon}(\mathbb{R}^{k},\varrho,1)}^{2}= \\
\int\limits_{\mathbb{R}^{k}}(\langle\xi\rangle+\varrho)^{2(\sigma-\varepsilon)}\,
\psi^{2}((\langle\xi\rangle+\varrho)^{\varepsilon+\delta})\,
|\widehat{u}(\xi)|^{2}\,d\xi= \\
\int\limits_{\mathbb{R}^{k}}(\langle\xi\rangle+\varrho)^{2\sigma}\,
\varphi^{2}(\langle\xi\rangle+\varrho)\,|\widehat{u}(\xi)|^{2}\,d\xi=
\|u\|^{2}_{H^{\sigma,\varphi}(\mathbb{R}^{k},\varrho,1)}.
\end{gather*}
Таким образом,
\begin{equation}\label{4.34}
\|u\|_{[H^{\sigma-\varepsilon}(\mathbb{R}^{k},\varrho,1),
H^{\sigma+\delta}(\mathbb{R}^{k},\varrho,1)]_{\psi}}=
\|u\|_{H^{\sigma,\varphi}(\mathbb{R}^{k},\varrho,1)}.
\end{equation}

Отметим, что
\begin{gather}\notag
c_{1}^{-1}\,\|u\| _{[H^{\sigma-\varepsilon}(\mathbb{R}^{k},\varrho),
H^{\sigma+\delta}(\mathbb{R}^{k},\varrho)]_{\psi}}\leq \\ \notag
\|u\|_{[H^{\sigma-\varepsilon}(\mathbb{R}^{k},\varrho,1),
H^{\sigma+\delta}(\mathbb{R}^{k},\varrho,1)]_{\psi}}\leq \\
c_{1}\,\|u\|_{[H^{\sigma-\varepsilon}(\mathbb{R}^{k},\varrho),
H^{\sigma+\delta}(\mathbb{R}^{k},\varrho)]_{\psi}}, \label{4.35}
\end{gather}
где число $c_{1}\geq1$ не зависит от $u$ и $\varrho$. В самом деле, тождественный
оператор $I$ задает гомеоморфизмы
\begin{gather*}
I:\,H^{\sigma-\varepsilon}(\mathbb{R}^{k},\varrho)\leftrightarrow
H^{\sigma-\varepsilon}(\mathbb{R}^{k},\varrho,1), \\
I:\,H^{\sigma+\delta}(\mathbb{R}^{k},\varrho)\leftrightarrow
H^{\sigma+\delta}(\mathbb{R}^{k},\varrho,1).
\end{gather*}
Здесь нормы прямых и обратных операторов ограничены равномерно по параметру
$\varrho$. Отсюда в силу теоремы \ref{th2.8} получаем гомеоморфизм
\begin{gather*}
I:\,[H^{\sigma-\varepsilon}(\mathbb{R}^{k},\varrho),
H^{\sigma+\delta}(\mathbb{R}^{k},\varrho)]_{\psi}\leftrightarrow \\
[H^{\sigma-\varepsilon}(\mathbb{R}^{k},\varrho,1),
H^{\sigma+\delta}(\mathbb{R}^{k},\varrho,1)]_{\psi},
\end{gather*}
такой, что нормы прямого и обратного оператора ограничены равномерно по параметру
$\varrho$. (Заметим, что записанные здесь пары пространств нормальные.) Это означает
двустороннюю оценку норм \eqref{4.35}.

Теперь из формул \eqref{4.33} -- \eqref{4.35} следует требуемая оценка \eqref{4.30}
для $G=\mathbb{R}^{k}$, где число $c:=c_{0}c_{1}$ не зависит от параметра $\varrho$.

Случай $G=\Omega$. Мы выведем его из предыдущего случая, в котором берем $k:=n$.
Пусть $\varrho\geq1$. Обозначим через $R_{\Omega}$ линейный оператор сужения
распределения с пространства $\mathbb{R}^{n}$ в область $\Omega$. Имеем ограниченные
операторы
\begin{gather} \label{4.36}
R_{\Omega}:\,H^{\sigma,\varphi}(\mathbb{R}^{n},\varrho)\rightarrow
H^{\sigma,\varphi}(\Omega,\varrho), \\
R_{\Omega}:\,H^{\alpha}(\mathbb{R}^{n},\varrho)\rightarrow
H^{\alpha}(\Omega,\varrho),\quad\alpha>0. \label{4.37}
\end{gather}
Их нормы, очевидно, не превосходят числа 1. Применим к пространствам, в которых
действуют операторы \eqref{4.37}, где
$\alpha\in\{\sigma-\varepsilon,\sigma+\delta\}$, интерполяцию с параметром $\psi$. В
силу теоремы \ref{th2.8} заключаем, что норма оператора
$$
R_{\Omega}:\,[H^{\sigma-\varepsilon}(\mathbb{R}^{n},\varrho),
H^{\sigma+\delta}(\mathbb{R}^{n},\varrho)]_{\psi}\rightarrow
[H^{\sigma-\varepsilon}(\Omega,\varrho),H^{\sigma+\delta}(\Omega,\varrho)]_{\psi},
$$
ограничена равномерно по параметру $\varrho$. (Поскольку здесь левая пара
пространств нормальная, то и правая пара тоже нормальная.) Отсюда, воспользовавшись
неравенством \eqref{4.30} для случая $G=\mathbb{R}^{n}$, делаем вывод, что норма
оператора
\begin{equation}\label{4.38}
R_{\Omega}:\,H^{\sigma,\varphi}(\mathbb{R}^{n},\varrho)\rightarrow
[H^{\sigma-\varepsilon}(\Omega,\varrho),H^{\sigma+\delta}(\Omega,\varrho)]_{\psi}
\end{equation}
ограничена равномерно по параметру $\varrho$.

Нам понадобится линейный ограниченный оператор, правый обратный к \eqref{4.38}. В
монографии [\ref{Triebel80}, c. 386] для каждого $l\in\mathbb{N}$ построено линейное
отображение $T_{l}$, продолжающее произвольное распределение $u\in H^{-l}(\Omega)$ в
пространство $\mathbb{R}^{n}$ и являющееся ограниченным оператором
\begin{equation}\label{4.39}
T_{l}:\,H^{\alpha}(\Omega)\rightarrow
H^{\alpha}(\mathbb{R}^{n})\quad\mbox{при}\quad\alpha\in\mathbb{R},\;\;|\alpha|<l.
\end{equation}
(Мы этот оператор уже использовали в п.~3.2.1.) Выберем целое $l>\sigma+\delta$.
Применив интерполяцию с параметром $\psi$ к \eqref{4.39}, где
$\alpha\in\{\sigma-\varepsilon,\sigma+\delta\}$, получим в силу теорем \ref{th2.14}
и \ref{th3.2} ограниченный оператор
\begin{equation}\label{4.40}
T_{l}:\,H^{\sigma,\varphi}(\Omega)\rightarrow H^{\sigma,\varphi}(\mathbb{R}^{n}).
\end{equation}

Из ограниченности операторов \eqref{4.39} и \eqref{4.40} следует, что нормы
операторов
\begin{gather}\label{4.41}
T_{l}:\,H^{\alpha}(\Omega,\varrho)\rightarrow
H^{\alpha}(\mathbb{R}^{n},\varrho),\quad 0<\alpha<l, \\
T_{l}:\,H^{\sigma,\varphi}(\Omega,\varrho)\rightarrow
H^{\sigma,\varphi}(\mathbb{R}^{n},\varrho) \label{4.42}
\end{gather}
ограничены равномерно по параметру $\varrho$. Применим к пространствам, в которых
действуют операторы \eqref{4.41}, где
$\alpha\in\{\sigma-\nobreak\varepsilon,\allowbreak\sigma+\nobreak\delta\}$,
интерполяцию с параметром $\psi$. На основании теоремы \ref{th2.8} и неравенства
\eqref{4.30} для случая $G=\mathbb{R}^{n}$ заключаем, что норма оператора
\begin{equation}\label{4.43}
T_{l}:\,[H^{\sigma-\varepsilon}(\Omega,\varrho),
H^{\sigma+\delta}(\Omega,\varrho)]_{\psi}\rightarrow
H^{\sigma,\varphi}(\mathbb{R}^{n},\varrho)
\end{equation}
ограничена равномерно по $\varrho$.

Оператор $R_{\Omega}T_{l}=I$ тождественный. Поэтому из равномерной ограниченности по
параметру $\varrho$ норм операторов \eqref{4.43}, \eqref{4.36} и \eqref{4.42},
\eqref{4.38} следует равномерная ограниченность по $\varrho$ норм операторов
вложения
\begin{gather*}
I=R_{\Omega}T_{l}:\,[H^{\sigma-\varepsilon}(\Omega,\varrho),
H^{\sigma+\delta}(\Omega,\varrho)]_{\psi}\rightarrow
H^{\sigma,\varphi}(\Omega,\varrho), \\
I=R_{\Omega}T_{l}:\,H^{\sigma,\varphi}(\Omega,\varrho)\rightarrow
[H^{\sigma-\varepsilon}(\Omega,\varrho), H^{\sigma+\delta}(\Omega,\varrho)]_{\psi}.
\end{gather*}
Это сразу дает двустороннюю оценку \eqref{4.30} для $G=\Omega$.

Случай $G=\Gamma$. Мы выведем его из первого случая $G=\mathbb{R}^{k}$, в котором
берем $k:=n-1$. Пусть $\varrho\geq1$. Будет рассуждать подобно доказательству леммы
\ref{lem2.14}. Воспользуемся локальным определением \ref{def2.13} пространств
$H^{s,\varphi}(\Gamma)$, $s\in\mathbb{R}$, $\varphi\in\mathcal{M}$, для
фиксированных конечного атласа и разбиения единицы на $\Gamma$. Рассмотрим линейное
отображение „распрямления” многообразия $\Gamma$:
$$
T:\,u\mapsto((\chi_{1}u)\circ\alpha_{1},\ldots, (\chi_{r}u)\circ\alpha_{r}),\quad
u\in\mathcal{D}'(\Gamma).
$$
Непосредственно проверяется, что это отображение задает изометрические операторы
\begin{gather}\label{4.44}
T:\,H^{\sigma,\varphi}(\Gamma,\varrho)\rightarrow
(H^{\sigma,\varphi}(\mathbb{R}^{n-1},\varrho))^{r}, \\
T:\,H^{\alpha}(\Gamma,\varrho)\rightarrow
(H^{\alpha}(\mathbb{R}^{n-1},\varrho))^{r},\quad\alpha>0. \label{4.45}
\end{gather}

Применим к пространствам, в которых действует операторы \eqref{4.45}, где
$\alpha\in\{\sigma-\varepsilon,\sigma+\delta\}$, интерполяцию с параметром $\psi$.
На основании теоремы \ref{th2.8} заключаем, что норма оператора
\begin{gather*}
T:\,[H^{\sigma-\varepsilon}(\Gamma,\varrho),
H^{\sigma+\delta}(\Gamma,\varrho)]_{\psi}\rightarrow \\
\bigl[(H^{\sigma-\varepsilon}(\mathbb{R}^{n-1},\varrho))^{r},
(H^{\sigma+\delta}(\mathbb{R}^{n-1},\varrho))^{r}]_{\psi}
\end{gather*}
ограничена равномерно по параметру $\varrho$. (Записанные здесь пары пространств,
очевидно, нормальные.) Отсюда в силу теоремы \ref{th2.5} и формулы \eqref{4.30},
доказанной для $G=\mathbb{R}^{n-1}$, делаем вывод, что норма оператора
\begin{equation} \label{4.46}
T:\,[H^{\sigma-\varepsilon}(\Gamma,\varrho),
H^{\sigma+\delta}(\Gamma,\varrho)]_{\psi}\rightarrow
(H^{\sigma,\varphi}(\mathbb{R}^{n-1},\varrho))^{r}
\end{equation}
также ограничена равномерно по $\varrho$.

Наряду с $T$ рассмотрим линейное отображение „склейки”
$$
K:\,(w_{1},\ldots,w_{r})\mapsto\sum_{j=1}^{r}
\Theta_{j}((\eta_{j}w_{j})\circ\alpha_{j}^{-1}),
$$
где $w_{1},\ldots,w_{r}$ --- распределения в $\mathbb{R}^{n-1}$. Здесь функция
$\eta_{j}\in C^{\infty}(\mathbb{R}^{n-1})$ финитная и равная 1 на множестве
$\alpha_{j}^{-1}(\mathrm{supp}\,\chi_{j})$, а $\Theta_{j}$ --- оператор продолжения
нулем на $\Gamma$. В силу \eqref{2.70} имеем ограниченные операторы
\begin{gather*}
K:\,(H^{\sigma,\varphi}(\mathbb{R}^{n-1}))^{r}\rightarrow
H^{\sigma,\varphi}(\Gamma), \\
K:\,(H^{\alpha}(\mathbb{R}^{n-1}))^{r}\rightarrow H^{\alpha}(\Gamma),\quad
\alpha\in\mathbb{R}.
\end{gather*}
Непосредственно проверяется, что норма каждого из операторов
\begin{gather}\label{4.47}
K:\,(H^{\sigma,\varphi}(\mathbb{R}^{n-1},\varrho))^{r}\rightarrow
H^{\sigma,\varphi}(\Gamma,\varrho),\\
K:\,(H^{\alpha}(\mathbb{R}^{n-1},\varrho))^{r}\rightarrow
H^{\alpha}(\Gamma,\varrho),\quad\alpha>0, \label{4.48}
\end{gather}
ограничена равномерно по параметру $\varrho$.

Применим к пространствам, в которых действуют операторы \eqref{4.48}, где
$\alpha\in\{\sigma-\varepsilon,\sigma+\delta\}$, интерполяцию с параметром $\psi$.
На основании теоремы \ref{th2.8} заключаем, что норма оператора
\begin{gather*}
K:\,[(H^{\sigma-\varepsilon}(\mathbb{R}^{n-1},\varrho))^{r},
(H^{\sigma+\delta}(\mathbb{R}^{n-1},\varrho))^{r}]_{\psi}\rightarrow\\
[H^{\sigma-\varepsilon}(\Gamma,\varrho), H^{\sigma+\delta}(\Gamma,\varrho)]_{\psi}
\end{gather*}
ограничена равномерно по $\varrho$. Отсюда в силу теоремы \ref{th2.5} и формулы
\eqref{4.30}, доказанной для $G=\mathbb{R}^{n-1}$, делаем вывод, что норма оператора
\begin{equation}\label{4.49}
K:\,(H^{\sigma,\varphi}(\mathbb{R}^{n-1},\varrho))^{r}\rightarrow
[H^{\sigma-\varepsilon}(\Gamma,\varrho), H^{\sigma+\delta}(\Gamma,\varrho)]_{\psi}
\end{equation}
также ограничена равномерно по $\varrho$.

В силу \eqref{2.69}, произведение $KT=I$ --- тождественный оператор. Поэтому из
равномерной ограниченности по параметру $\varrho$ норм операторов \eqref{4.44},
\eqref{4.49} и \eqref{4.46}, \eqref{4.47} следует равномерная ограниченность по
$\varrho$ норм операторов вложения
\begin{gather*}
I=KT:\,H^{\sigma,\varphi}(\Gamma,\varrho)\rightarrow
[H^{\sigma-\varepsilon}(\Gamma,\varrho),
H^{\sigma+\delta}(\Gamma,\varrho)]_{\psi}, \\
I=KT:\,[H^{\sigma-\varepsilon}(\Gamma,\varrho),
H^{\sigma+\delta}(\Gamma,\varrho)]_{\psi}\rightarrow
H^{\sigma,\varphi}(\Gamma,\varrho).
\end{gather*}
Это сразу дает двустороннюю оценку \eqref{4.30} для $G=\Gamma$.

Лемма \ref{lem4.1} доказана.


В доказательстве леммы \ref{lem4.1} был использован следующий результат


\begin{lemma}\label{lem4.2}
Пусть $\sigma>0$, $\varphi\in\mathcal{M}$ и
$\varphi_{\sigma}(t):=t^{\sigma}\varphi(t)$ при $t\geq1$. Тогда существует число
$c=c(\sigma,\varphi)\geq1$ такое, что
\begin{equation}\label{4.50}
c^{-1}\varphi_{\sigma}(t_{1}+t_{2})\leq\varphi_{\sigma}(t_{1})+\varphi_{\sigma}(t_{2})
\leq c\,\varphi_{\sigma}(t_{1}+t_{2})
\end{equation}
для всех $t_{1},t_{2}\geq1$.
\end{lemma}


\textbf{Доказательство.} Так как $\varphi_{\sigma}\in\mathrm{RO}$ (см. п.~2.8.1), то
\begin{equation}\label{4.51}
\varphi_{\sigma}(2t)\asymp\varphi_{\sigma}(t)\quad\mbox{при}\quad t\geq1.
\end{equation}
Кроме того, поскольку функция $\varphi_{\sigma}$ имеет порядок изменения $\sigma>0$,
она (слабо) эквивалентна некоторой  возрастающей положительной функции $\psi$:
\begin{equation}\label{4.52}
\varphi_{\sigma}(t)\asymp\psi(t)\quad\mbox{при}\quad t\geq1.
\end{equation}

Теперь \eqref{4.50} вытекает из формул \eqref{4.51}, \eqref{4.52} и возрастания
функции $\psi$. В самом деле, для каждого номера $j=1,2$ можем записать
$$
\varphi_{\sigma}(t_{j})\asymp\psi(t_{j})\leq\psi(t_{1}+t_{2})\asymp
\varphi_{\sigma}(t_{1}+t_{2})\quad\mbox{при}\quad t_{1},t_{2}\geq1.
$$
Поэтому существует число $c_{1}>0$ такое, что
\begin{equation}\label{4.53}
\varphi_{\sigma}(t_{1})+\varphi_{\sigma}(t_{2})\leq
c_{1}\varphi_{\sigma}(t_{1}+t_{2})\quad\mbox{при}\quad t_{1},t_{2}\geq1.
\end{equation}
Обратно, предположив без потери общности, что $t_{1}\leq t_{2}$, имеем:
\begin{gather*}
\varphi_{\sigma}(t_{1}+t_{2})\asymp\psi(t_{1}+t_{2})\leq\psi(2t_{2})\asymp
\varphi_{\sigma}(2t_{2})\asymp \\
\varphi_{\sigma}(t_{2})\leq
\varphi_{\sigma}(t_{1})+\varphi_{\sigma}(t_{2})\quad\mbox{при}\quad
t_{1},t_{2}\geq1.
\end{gather*}
Следовательно, существует число $c_{2}>0$ такое, что
\begin{equation}\label{4.54}
\varphi_{\sigma}(t_{1}+t_{2})\leq
c_{2}\,(\varphi_{\sigma}(t_{1})+\varphi_{\sigma}(t_{2})) \quad\mbox{при}\quad
t_{1},t_{2}\geq1.
\end{equation}
Формулы \eqref{4.53} и \eqref{4.54} означают двойное неравенство \eqref{4.50}.

Лемма \ref{lem4.2} доказана.


Используя лемму \ref{lem4.1}, мы можем дать


\textbf{Доказательство утверждения (ii) теоремы \ref{th4.9}.} Пусть $s>2q$,
$\varphi\in\mathcal{M}$, а параметр $\lambda\in K$ такой, что
$|\lambda|\geq\max\{\lambda_{0},1\}$, где число $\lambda_{0}>0$ взято из утверждения
(i) этой теоремы. Положим $\varepsilon=\delta=(s-2q)/2>0$. Как отмечалось выше, в
соболевском случае $\varphi\equiv1$ теорема \ref{th4.9} верна. Значит, мы имеем
гомеоморфизмы
\begin{gather}\label{4.55}
(L(\lambda),B(\lambda)):\;H^{s\mp\varepsilon}(\Omega,|\lambda|)\leftrightarrow
\\ H^{s\mp\varepsilon-2q}(\Omega,|\lambda|)\oplus\bigoplus_{j=1}^{q}
H^{s\mp\varepsilon-m_{j}-1/2}(\Gamma,|\lambda|)
=:\mathcal{H}_{s\mp\varepsilon}(\Omega,\Gamma,|\lambda|)\notag
\end{gather}
такие, что нормы прямого и обратного операторов ограничены равномерно по параметру
$\lambda$. (Заметим, что мы перешли в оценке \eqref{4.28} к гильбертовым нормам.)
Пусть $\psi$ --- интерполяционный параметр из теоремы \ref{th2.14}. Применив
интерполяцию с этим параметром к пространствам, в которых действует гомеоморфизм
\eqref{4.55}, получим еще один гомеоморфизм
\begin{gather}\notag
(L(\lambda),B(\lambda)):\;[H^{s-\varepsilon}(\Omega,|\lambda|),
H^{s+\varepsilon}(\Omega,|\lambda|)]_{\psi}\leftrightarrow \\
[\mathcal{H}_{s-\varepsilon}(\Omega,\Gamma,|\lambda|),
\mathcal{H}_{s+\varepsilon}(\Omega,\Gamma,|\lambda|)]_{\psi}. \label{4.56}
\end{gather}
Здесь нормы прямого и обратного операторов ограничены равномерно по $\lambda$ в силу
теоремы \ref{th2.8}. (Заметим, что в \eqref{4.56} записаны нормальные пары
пространств.)

Далее, согласно теореме \ref{th2.5} об интерполяции ортогональных сумм пространств
можем записать
\begin{gather*}
[\mathcal{H}_{s-\varepsilon}(\Omega,\Gamma,|\lambda|),
\mathcal{H}_{s+\varepsilon}(\Omega,\Gamma,|\lambda|)]_{\psi}=\\
\bigl[H^{s-\varepsilon-2q}(\Omega,|\lambda|),
H^{s+\varepsilon-2q}(\Omega,|\lambda|)]_{\psi}\,\oplus \\
\bigoplus_{j=1}^{q}\;[H^{s-\varepsilon-m_{j}-1/2} (\Gamma,|\lambda|),
H^{s+\varepsilon-m_{j}-1/2} (\Gamma,|\lambda|)]_{\psi}
\end{gather*}
с равенством норм. Отсюда, в силу леммы \ref{lem4.1} и \eqref{4.56} следует
гомеоморфизм
\begin{gather*}
(L(\lambda),B(\lambda)):\;H^{s,\varphi}(\Omega,|\lambda|)\leftrightarrow\\
H^{s-2q,\varphi}(\Omega,|\lambda|)\oplus\bigoplus_{j=1}^{q}
H^{s-m_{j}-1/2,\varphi}(\Gamma,|\lambda|)
\end{gather*}
такой, что нормы прямого и обратного операторов ограничены равномерно по параметру
$\lambda$. Это немедленно дает двустороннюю оценку \eqref{4.28}.

Утверждение (ii) доказано.


Таким образом, теорема \ref{th4.9} доказана


\begin{corollary}\label{cor4.1}
Пусть краевая задача \eqref{4.24}, \eqref{4.25} эллиптическая с параметром на
некотором замкнутом луче $K:=\{\lambda\in\mathbb{C}:\arg\lambda =\mathrm{const}\}$.
Тогда ограниченный оператор \eqref{4.27} имеет нулевой индекс при любых $s>2q$,
$\varphi\in\mathcal{M}$ и $\lambda\in\mathbb{C}$.
\end{corollary}


\textbf{Доказательство.} Как отмечалось выше, оператор \eqref{4.27} нетеров, и его
индекс не зависит от указанных параметров $s$, $\varphi$ и $\lambda$. В~силу теоремы
\ref{th4.9} индекс оператора \eqref{4.27} равен $0$ при $|\lambda|\gg1$. Значит, он
равен $0$ для всех $\lambda\in\mathbb{C}$. Следствие \ref{cor4.1} доказано.

\subsection[Формально смешанная эллиптическая краевая задача]
{Формально смешанная \\ эллиптическая краевая задача}\label{sec4.1.5}

В этом пункте рассматривается эллиптическая краевая задача для линейного
дифференциального уравнения $Lu=f$ в многосвязной области $\Omega$. В отличие от
предыдущих пунктов предполагается, что порядки граничных выражений различны на
разных связных компонентах границы $\Gamma$. Например, для уравнения Лапласа в
кольце можно задавать граничное условие Дирихле на одной и граничное условие Неймана
на другой компонентах границы. Рассматриваемая задача относится к классу смешанных
задач [\ref{VishikEskin69}, \ref{Peetre61}, \ref{Schechter60b}, \ref{Simanca87}].
Они изучены существенно менее полно, чем несмешанные задачи. Это связано с тем, что
при сведении смешанной задачи к псевдодифференциальному оператору на границе
возникают определенные трудности (см., например, [\ref{Simanca87}]).
В~рассматриваемой нами задаче участки границы, на которых порядок граничного
выражения различен, не примыкают друг к другу. Такую задачу мы называем формально
смешанной. Ее с помощью локальных построений можно свести к эллиптической модельной
задаче в полупространстве [\ref{94Dop12}].

В этом пункте предполагается, что граница $\Gamma$ области $\Omega$ состоит из
$r\geq2$ непустых связных компонент $\Gamma_{1},\ldots,\Gamma_{r}$. Рассмотрим
формально смешанную краевую задачу в области $\Omega$
\begin{gather} \label{4.57}
Lu=f\quad\text{в}\quad\Omega \\
B_{k,j}u=g_{k,j}\quad\text{на}\quad\Gamma_{k}\quad\text{при}\;\;
j=1,\ldots,q\;\;\text{и}\;\;k=1,\ldots,r. \label{4.58}
\end{gather}
Здесь, как и прежде, $L$ --- линейное дифференциальное выражение на
$\overline{\Omega}$ четного порядка $2q$, а $\{B_{k,j}:j=1,\ldots,q\}$ --- система
граничных линейных дифференциальных выражений, заданных на компоненте $\Gamma_{k}$.
Предполагается, что все порядки $m_{k,j}:=\mathrm{ord}\,B_{k,j}\leq2q-1$.
Коэффициенты дифференциальных выражений $L$ и $B_{k,j}$ являются бесконечно гладкими
комплекснозначными функциями. Положим
\begin{gather*}
\Lambda:=(L,B_{1,1},\ldots,B_{1,q},\ldots,B_{r,1},\ldots,B_{r,q}), \\
N_{\Lambda}:=\{u\in C^{\infty}(\,\overline{\Omega}\,):\,\Lambda u=0\}.
\end{gather*}


\begin{definition}\label{def4.2}
Формально смешанная краевая задача \eqref{4.57}, \eqref{4.58} называется
\emph{эллиптической} в многосвязной области $\Omega$, если выполняются следующие
условия:
\begin{itemize}
\item [$\mathrm{(i)}$] Дифференциальное выражение $L$ правильно эллиптическое на
$\overline{\Omega}$.
\item [$\mathrm{(ii)}$] Для каждого $k=1,\ldots,r$ система граничных выражений $\{B_{k,j}:
j=1,\ldots,q\}$ удовлетворяет условию дополнительности по отношению к $L$ на
$\Gamma_{k}$.
\end{itemize}
\end{definition}


\begin{theorem}\label{th4.10}
Предположим, что краевая задача \eqref{4.57}, \eqref{4.58} эллиптическая в области
$\Omega$. Пусть $s>2q$ и $\varphi\in\mathcal{M}$. Тогда отображение $u\mapsto\Lambda
u$, где $u\in C^{\infty}(\,\overline{\Omega}\,)$, продолжается по непрерывности до
ограниченного нетерового оператора
\begin{equation}\label{4.59}
\Lambda:\,H^{s,\varphi}(\Omega)\rightarrow
H^{s-2q,\;\varphi}(\Omega)\oplus\bigoplus_{k=1}^{r}\bigoplus_{j=1}^{q}
H^{s-m_{k,j}-1/2,\varphi}(\Gamma_{k}) =:\mathcal{H}_{s,\varphi}.
\end{equation}
Ядро этого оператора совпадает с $N_{\Lambda}$, а область значений состоит из всех
векторов
$$
(f,g_{1,1},\ldots,g_{1,q},\ldots,g_{r,1},\ldots,g_{r,q})\in \mathcal{H}_{s,\varphi}
$$
таких, что
$$
(f,w_{0})_{\Omega}+\sum_{k=1}^{r}\,\sum_{j=1}^{q}\; (g_{k,j},w_{k,j})_{\Gamma_{k}}=0
$$
для любой вектор-функции
$$
(w_{0},w_{1,1},\ldots,w_{1,q},\ldots,w_{r,1},\ldots,w_{r,q})\in W_{\Lambda}.
$$
Здесь $W_{\Lambda}$ --- некоторое не зависящее  от $s$ и $\varphi$ конечномерное
подпространство в
$$
C^{\infty}(\,\overline{\Omega}\,)\times
\prod_{j=1}^{r}\,(C^{\infty}(\Gamma_{j}))^{q}.
$$
Индекс оператора \eqref{4.59} равен числу $\dim N_{\Lambda}-\dim W_{\Lambda}$ и
также не зависит от $s$, $\varphi$.
\end{theorem}


В соболевском случае $\varphi\equiv1$ эта теорема является частным случаем теоремы~1
работы [\ref{94Dop12}, с.~38]. (Там ограниченность и нетеровость оператора $\Lambda$
была установлена в шкале пространств Лизоркина-Трибеля, которая содержит в себе
соболевскую шкалу). Отсюда общий случай $\varphi\in\mathcal{M}$ выводится с помощью
интерполяции аналогично доказательству теоремы~\ref{th4.1}.




\newpage

\markright{\emph \ref{sec4.2}. Эллиптическая краевая задача в двусторонней шкале}

\section[Эллиптическая краевая задача в двусторонней\break шкале]
{Эллиптическая краевая задача \\ в двусторонней шкале}\label{sec4.2}

\markright{\emph \ref{sec4.2}. Эллиптическая краевая задача в двусторонней шкале}

В этом пункте мы изучим регулярную эллиптическую краевую задачу \eqref{4.1} в
двусторонней шкале пространств Хермандера. В отличие от предыдущего п.~\ref{sec4.1},
здесь числовой индекс, задающий основную гладкость, пробегает всю вещественную ось.

\subsection{Предварительные замечания}\label{sec4.2.1}

Теоремы о разрешимости эллиптических краевых задач, доказанные в п.~\ref{sec4.1},
вообще говоря, не верны для произвольного вещественного параметра $s$, задающего
основную гладкость решения задачи. Это обусловлено тем, что отображение $u\mapsto
u\!\upharpoonright\!\Gamma$, где $u\in C^{\infty}(\,\overline{\Omega}\,)$, не
продолжается до непрерывного оператора взятия следа
$R_{\Gamma}:H^{s,\varphi}(\Omega)\rightarrow\mathcal{D}'(\Gamma)$ при $s<1/2$ (см.
замечание \ref{rem3.5}). Следовательно, оператор \eqref{4.3}, соответствующий
краевой задаче \eqref{4.1}, не определен при $s<m+1/2$, где $m$~--- максимум
порядков граничных выражений $B_{j}$. Аналогично и для других краевых задач
рассмотренных в п.~\ref{sec4.1}.

Для того, чтобы получить ограниченный оператор $(L,B)$ при любом $s<2q$, надо вместо
$H^{s,\varphi}(\Omega)$ брать иное пространство в качестве области определения этого
оператора.  Известно два принципиально различных способа построения такой области
определения, предложенные Я.~А.~Ройтбергом [\ref{Roitberg64}, \ref{Roitberg65},
\ref{Roitberg96}] и Ж.-Л.~Лионсом, Е.~Мадженесом [\ref{LionsMagenes71},
\ref{Magenes66}, \ref{LionsMagenes62}, \ref{LionsMagenes63}] в соболевском случае.
Они приводят к различным типам теорем о разрешимости эллиптических краевых задач~---
общим и индивидуальным теоремам. В~общей теореме область определения оператора
$(L,B)$ не зависит от коэффициентов эллиптического выражения $L$ и является единой
для всех краевых задач одного порядка. В индивидуальных теоремах она зависит от
коэффициентов выражения $L$ (даже от коэффициентов при младших производных).

Заметим, что теоремы о разрешимости эллиптических краевых задач, доказанные в
п.~\ref{sec4.1}, являются общими (для соответствующих классов задач).

В настоящем п.~\ref{sec4.2} мы реализуем подход Я.~А.~Ройтберга применительно к
двусторонней уточненной шкале пространств Хермандера --- модифицируем по Ройтбергу
эту шкалу и получим общую теорему о разрешимости регулярной эллиптической краевой
задачи \eqref{4.1} в такой модифицированной шкале. (Для других эллиптических краевых
задач, рассмотренных в п.~\ref{sec4.1}, это несложно сделать, следуя развитой ниже
теории. Мы не будем на этом останавливаться.) Подход Ж.-Л.~Лионса и Е.~Мадженеса,
приводящий к индивидуальным теоремам о разрешимости, будет рассмотрен далее в
п~\ref{sec4.4}, \ref{sec4.5}.

Отметим, что неоднородная эллиптическая краевая задача \eqref{4.1} не сводится при
$s<m+1/2$ к двум полуоднородным краевым задачам, изученным в п. \ref{sec3.3} и
\ref{sec3.4}. В самом деле, если $s<-1/2$, то решения этих задач принадлежат к
пространствам распределений различной природы. А~именно, решения полуоднородной
задачи \eqref{3.36} принадлежат к пространству $K_{L}^{s,\varphi}(\Omega)\subset
H^{s,\varphi}(\Omega)$ и являются распределениями, заданными в открытой области
$\Omega$, а решения полуоднородной задачи \eqref{3.82}, \eqref{3.83} принадлежат к
пространству $H^{s,\varphi}(\mathrm{b.c.})\subset
H^{s,\varphi}_{\overline{\Omega}}(\mathbb{R}^{n})$ и являются распределениями,
сосредоточенными в замкнутой области $\overline{\Omega}$. Если же $-1/2<s<m+1/2$, то
решения обеих полуоднородных задач будут распределениями в области $\Omega$ (см.
формулу \eqref{3.32}), но оператор $(L,B)$ нельзя корректно определить на
пространстве $K_{L}^{s,\varphi}(\Omega)\cup H^{s,\varphi}(\mathrm{b.c.})$, поскольку
$(K_{L}^{s,\varphi}(\Omega)\cap H^{s,\varphi}(\mathrm{b.c.}))\setminus
N\neq\varnothing$. Например, левая часть этого неравенства содержит ввиду теоремы
\ref{th3.11} решение $u\in C^{\infty}(\,\overline{\Omega}\,)$ краевой задачи
$$
Lu=0\;\;\mbox{в}\;\;\Omega,\quad
B_{1}u=0,\ldots,B_{q-1}u=0,\,B_{q}u=h\;\;\mbox{на}\;\;\Gamma.
$$
Здесь мы считаем, что $\mathrm{ord}\,B_{q}=m$, а ненулевую функцию $h\in
C^{\infty}(\Gamma)$ выбираем так, чтобы $(h,C^{+}_{q}v)_{\Gamma}=0$ для всех $v\in
N^{+}$.

\subsection{Модифицированная уточненная шкала}\label{sec4.2.2}

Сначала введем понятие обобщенного решения по Ройтбергу [\ref{Roitberg64};
\ref{Roitberg65}; \ref{Roitberg96}, с.~80] краевой задачи \eqref{4.1}.

В окрестности границы $\Gamma$ запишем дифференциальные выражения $L$ и $B_{j}$ в
виде
\begin{equation}\label{4.63}
L=\sum_{k=0}^{2q}\;L_{k}\,D_{\nu}^{k}\quad\mbox{и}\quad
B_{j}=\sum_{k=0}^{m_{j}}\;B_{j,k}\,D_{\nu}^{k}.
\end{equation}
Здесь, напомним, $D_{\nu}:=i\,\partial/\partial\nu$, а $L_{k}$ и $B_{j,k}$ ---
некоторые тангенциальные (по отношению к границе $\Gamma$) дифференциальные
выражения. Проинтегрировав по частям, запишем следующую формулу Грина:
\begin{equation}\label{4.64}
(Lu,v)_{\Omega}=(u,L^{+}v)_{\Omega}-
i\sum_{k=1}^{2q}\;(D_{\nu}^{k-1}u,L^{(k)}v)_{\Gamma}
\end{equation}
для произвольных функций $u,v\in C^{\infty}(\,\overline{\Omega}\,)$. Здесь
$$
L^{(k)}:=\sum_{r=k}^{2q}D_{\nu}^{r-k}L_{r}^{+},
$$
где $L_{r}^{+}$~--- дифференциальное выражение, формально сопряженное к $L_{r}$. С
помощью предельного перехода получаем, что формула \eqref{4.64} верна для каждой
функции $u\in H^{2q}(\Omega)$. Обозначим:
\begin{equation}\label{4.65}
u_{0}:=u\quad\mbox{и}\quad
u_{k}:=(D_{\nu}^{k-1}u)\!\upharpoonright\!\Gamma\quad\mbox{при}\quad k=1,\ldots2q.
\end{equation}

В силу \eqref{4.63}, \eqref{4.64} краевая задача \eqref{4.1} относительного искомой
функции $u\in H^{2q}(\Omega)$ равносильна системе условий
\begin{gather} \label{4.66}
(u_{0},L^{+}v)_{\Omega}- i\sum_{k=1}^{2q}\;(u_{k},L^{(k)}v)_{\Gamma}=(f,v)_{\Omega}
\quad\mbox{для
всех}\quad v\in C^{\infty}(\,\overline{\Omega}\,), \\
\sum_{k=0}^{m_{j}}\;B_{j,k}\,u_{k+1}=g_{j}
\quad\mbox{на}\;\;\Gamma\quad\mbox{при}\quad j=1,\ldots,q. \label{4.67}
\end{gather}
Заметим, что эти условия имеют смысл в случае произвольных распределений
\begin{gather}\label{4.68}
u_{0}\in\mathcal{D}'(\mathbb{R}^{n}),\;\;
\mathrm{supp}\,u_{0}\subseteq\overline{\Omega},\quad
u_{1},\ldots,u_{2q}\in\mathcal{D}'(\Gamma),\\
f\in\mathcal{D}'(\mathbb{R}^{n}),\;\;\mathrm{supp}\,f\subseteq\overline{\Omega},\quad
g_{1},\ldots,g_{q}\in\mathcal{D}'(\Gamma).\quad\quad\; \notag
\end{gather}
Поэтому полезно ввести следующее определение.


\begin{definition}\label{def4.5}
Вектор $u=(u_{0},u_{1},\ldots,u_{2q})$, удовлетворяющий условию \eqref{4.68},
называем \emph{обобщенным по Ройтбергу решением} краевой задачи \eqref{4.1}, если
выполняются условия \eqref{4.66}, \eqref{4.67}.
\end{definition}


Введем гильбертовы пространства, элементы которых можно рассматривать как обобщенные
по Ройтбергу решения.

Пусть $r\in\mathbb{N}$, $s\in\mathbb{R}$ и $\varphi\in\mathcal{M}$. Положим
$$
E_{r}:=\{k-1/2:k=1,\ldots,r\}.
$$


\begin{definition}\label{def4.6}
В случае $s\in\mathbb{R}\setminus E_{r}$ линейное пространство
$H^{s,\varphi,(r)}(\Omega)$ является, по определению, пополнением линейного
многообразия $C^{\infty}(\,\overline{\Omega}\,)$ по гильбертовой норме
\begin{gather}\notag
\|u\|_{H^{s,\varphi,(r)}(\Omega)}:=\\
\biggl(\,\|u\|_{H^{s,\varphi,(0)}(\Omega)}^{2}+
\sum_{k=1}^{r}\;\|(D_{\nu}^{k-1}u)\!\upharpoonright\!\Gamma\|
_{H^{s-k+1/2,\varphi}(\Gamma)}^{2}\biggr)^{1/2}.\label{4.69}
\end{gather}
В случае $s\in E_{r}$ пространство $H^{s,\varphi,(r)}(\Omega)$ определяем
посредством интерполяции:
\begin{gather}\notag
H^{s,\varphi,(r)}(\Omega):=\\ [H^{s-\varepsilon,\varphi,(r)}(\Omega),\,
H^{s+\varepsilon,\varphi,(r)}(\Omega)]_{t^{1/2}}\quad\mbox{при}\quad
0<\varepsilon<1. \label{4.70}
\end{gather}
\end{definition}


\begin{remark}\label{rem4.6}
В формуле \eqref{4.70} использована интерполяция гильбертовых пространств со
степенным параметром $\psi(t)=t^{1/2}$. Ниже мы покажем, что записанная в правой
части этой формулы пара пространств допустимая (теорема \ref{th4.13} (i), (iv)), и
что результат интерполяции не зависит от $\varepsilon$ с точностью до
эквивалентности норм (п.~\ref{sec4.2.5}, теорема \ref{th4.22}).
\end{remark}


В соболевском случае $\varphi\equiv1$ пространство $H^{s,\varphi,(r)}(\Omega)$
введено Я.~А.~Ройтбергом [\ref{Roitberg64}; \ref{Roitberg65}; \ref{Roitberg96},
с.~69]. Как обычно, мы полагаем $H^{s,(r)}(\Omega):=H^{s,1,(r)}(\Omega)$.


\begin{definition}\label{def4.7}
Семейство гильбертовых пространств
\begin{equation}\label{4.71}
\{H^{s,\varphi,(r)}(\Omega):\,s\in\mathbb{R},\varphi\in\mathcal{M}\}
\end{equation}
называем \emph{модифицированной} по Ройтбергу уточненной шкалой. Число $r$ называем
порядком модификации.
\end{definition}


С точки зрения приложения к краевой задаче \eqref{4.1} нам интересен случай четного
порядка модификации $r=2q$. В силу определения \ref{def4.6}, отображение
$u\mapsto(D_{\nu}^{k-1}u)\!\upharpoonright\!\Gamma$, где $u\in
C^{\infty}(\,\overline{\Omega}\,)$ и $k\in\{0,1,\ldots,2q\}$, продолжается по
непрерывности до ограниченного оператора, действующего из пространства
$H^{s,\varphi,(r)}(\Omega)$ в пространство $H^{s-k+1/2,\varphi}(\Gamma)$ при любом
$s\in\mathbb{R}$. Следовательно, для произвольного элемента $u\in
H^{s,\varphi,(2q)}(\Omega)$ корректно определен по формулам \eqref{4.65} посредством
замыкания вектор
\begin{equation}\label{4.72}
(u_{0},u_{1},\ldots,u_{2q})\in
H^{s,\varphi,(0)}(\Omega)\oplus\bigoplus_{k=1}^{2q}\,H^{s-k+1/2,\varphi}(\Gamma).
\end{equation}
Поэтому можно рассматривать элемент $u$ как обобщенное (по Ройтбергу) решение
\eqref{4.72} краевой задачи \eqref{4.1}.

Изучим свойства модифицированной шкалы \eqref{4.71}. Для произвольных
$s\in\nobreak\mathbb{R}$ и $\varphi\in\mathcal{M}$ положим
$$
\Pi_{s,\varphi,(r)}(\Omega,\Gamma):=H^{s,\varphi,(0)}(\Omega)\oplus\bigoplus
_{k=1}^{r}\,H^{s-k+1/2,\varphi}(\Gamma).
$$
Кроме того, обозначим
\begin{gather*}
K_{s,\varphi,(r)}(\Omega,\Gamma):=
\{(u_{0},u_{1},\ldots,u_{r})\in\Pi_{s,\varphi,(r)}(\Omega,\Gamma):\\
u_{k}=(D_{\nu}^{k-1}u_{0})\!\upharpoonright\!\Gamma\;\mbox{для всех}\;k=1,\ldots
r\;\mbox{таких, что}\;s>k-1/2\}.
\end{gather*}
В силу теоремы \ref{th3.5}, $K_{s,\varphi,(r)}(\Omega,\Gamma)$ замкнуто в
$\Pi_{s,\varphi,(r)}(\Omega,\Gamma)$. Мы рассматриваем
$K_{s,\varphi,(r)}(\Omega,\Gamma)$ как гильбертово пространство относительно
скалярного произведения в $\Pi_{s,\varphi,(r)}(\Omega,\Gamma)$. В соболевском случае
$\varphi\equiv1$ мы опускаем индекс $\varphi$ в обозначениях пространств, вводимых в
этой главе.


\begin{theorem}\label{th4.12}
Пусть $r\in\mathbb{N}$, $s\in\mathbb{R}\setminus E_{r}$ и $\varphi\in\mathcal{M}$.
Справедливы следующие утверждения.
\begin{itemize}
\item [$\mathrm{(i)}$] Линейное отображение
\begin{equation}\label{4.73}
T_{r}:\,u\mapsto(u,u\!\upharpoonright\!\Gamma,\ldots,
(D_{\nu}^{r-1}u)\!\upharpoonright\!\Gamma),\quad u\in
C^{\infty}(\,\overline{\Omega}\,),
\end{equation}
продолжается по непрерывности до изометрического изоморфизма
\begin{equation}\label{4.74}
T_{r}:\,H^{s,\varphi,(r)}(\Omega)\leftrightarrow K_{s,\varphi,(r)}(\Omega,\Gamma).
\end{equation}
\item [$\mathrm{(ii)}$] Для произвольных положительных чисел $\varepsilon$, $\delta$
таких, что числа $s$, $s-\varepsilon$, $s+\delta$ принадлежат одному из интервалов
\begin{gather} \notag
\alpha_{0}:=(-\infty,1/2),\\ \notag
 \alpha_{k}:=(k-1/2,\,k+1/2),\;\;k=1,\ldots,r-1,\\
\alpha_{r}:=(r-1/2,+\infty),\label{4.75}
\end{gather}
справедливо равенство пространств с эквивалентностью норм в них:
\begin{equation}\label{4.76}
[H^{s-\varepsilon,(r)}(\Omega),H^{s+\delta,(r)}(\Omega)]_{\psi}=
H^{s,\varphi,(r)}(\Omega).
\end{equation}
Здесь $\psi$ --- интерполяционный параметр из теоремы $\ref{th2.14}$.
\end{itemize}
\end{theorem}


\textbf{Доказательство.} В случае $\varphi\equiv1$ утверждение (i) установлено
Я.~А.~Рой\-тберг\-ом [\ref{Roitberg96}, с.~71] (лемма 2.2.1). Мы выведем отсюда
утверждение (ii) для произвольного $\varphi\in\mathcal{M}$, а затем и
утверждение~(i).

Обозначим через $X_{\psi}$ левую часть равенства \eqref{4.76}. (Записанная там пара
пространств, очевидно, допустимая.) Рассмотрим изометрические операторы
$$
T_{r}:\,H^{\sigma,(r)}(\Omega)\rightarrow
\Pi_{\sigma,(r)}(\Omega,\Gamma),\quad\sigma\in\{s-\varepsilon,s+\delta\}.
$$
Применив интерполяцию с параметром $\psi$, получим ограниченный оператор
$$
T_{r}:\,X_{\psi}\rightarrow[\Pi_{s-\varepsilon,(r)}(\Omega,\Gamma),
\Pi_{s+\delta,(r)}(\Omega,\Gamma)]_{\psi}.
$$
В силу интерполяционных теорем \ref{th2.5}, \ref{th2.21} и \ref{th3.10} имеем
равенства пространств с эквивалентностью норм:
\begin{gather*}
[\Pi_{s-\varepsilon,(r)}(\Omega,\Gamma),
\Pi_{s+\delta,(r)}(\Omega,\Gamma)]_{\psi}= \\
[H^{s-\varepsilon,(0)}(\Omega),H^{s+\delta,(0)}(\Omega)]_{\psi}\oplus\\
\bigoplus_{k=1}^{r}\;[H^{s-\varepsilon-k+1/2}(\Gamma),
H^{s+\delta-k+1/2}(\Gamma)]_{\psi}= \\
H^{s,\varphi,(0)}(\Omega)\oplus \bigoplus_{k=1}^{r}\,H^{s-k+1/2,\varphi}(\Gamma)=
\Pi_{s,\varphi,(r)}(\Omega,\Gamma).
\end{gather*}
Следовательно, ограничен оператор
\begin{equation}\label{4.77}
T_{r}:\,X_{\psi}\rightarrow\Pi_{s,\varphi,(r)}(\Omega,\Gamma).
\end{equation}
Отсюда получаем оценку
\begin{equation}\label{4.78}
\|u\|_{H^{s,\varphi,(r)}(\Omega)}=\|T_{r}\,u\|_{\Pi_{s,\varphi,(r)}(\Omega,\Gamma)}\leq
c_{1}\,\|u\|_{X_{\psi}}
\end{equation}
для любого $u\in C^{\infty}(\,\overline{\Omega}\,)$. Здесь $c_{1}$ --- норма
оператора \eqref{4.77}.

Докажем неравенство, обратное к \eqref{4.78}. По условию, $s,\nobreak
s-\nobreak\varepsilon,\allowbreak s+\delta\in\alpha_{p}$ для некоторого номера
$p\in\{0,1,\ldots,r\}$. Рассмотрим линейное отображение
$$
T_{r,p}:\,u\mapsto\bigl(u,\,\{(D_{\nu}^{k-1}u)\!\upharpoonright\!\Gamma:p+1\leq
k\leq r\}\,\bigr),\quad u\in C^{\infty}(\,\overline{\Omega}\,).
$$
(Как и прежде, индекс $k$ целый.) Это отображение продолжается по непрерывности до
гомеоморфизма
\begin{gather}\notag
T_{r,p}:\,H^{\sigma,(r)}(\Omega)\,\leftrightarrow\\
H^{\sigma,(0)}(\Omega)\,\oplus\bigoplus_{p+1\leq k\leq
r}H^{\sigma-k+1/2}(\Gamma),\quad\sigma\in\{s-\varepsilon,s+\delta\}.\label{4.79}
\end{gather}

Действительно, ограниченность оператора \eqref{4.79} следует из определения
пространства $H^{\sigma,(r)}(\Omega)$. Покажем, что этот оператор биективный. Пусть
$u\in\nobreak H^{\sigma,(r)}(\Omega)$ и
$$
\bigl(u_{0},\,\{u_{k}:\,p+1\leq k\leq r\}\bigr)\in
H^{\sigma,(0)}(\Omega)\,\oplus\bigoplus_{p+1\leq k\leq r}\,H^{\sigma-k+1/2}(\Gamma).
$$
Положим $u_{k}:=(D_{\nu}^{k-1}u_{0})\!\upharpoonright\!\Gamma$ при $1\leq k\leq p$.
Распределение $u_{k}$ определено корректно в силу теоремы \ref{th3.5}, поскольку
$\sigma>k-1/2$ для указанных номеров~$k$. Заметим, что $\sigma<k-1/2$ при $p+1\leq
k\leq r$. Поэтому $(u_{0},u_{1},\ldots,u_{r})\in K_{\sigma,(r)}(\Omega,\Gamma)$. Как
отмечалось выше, утверждение (i) доказываемой теоремы верно в случае
$\varphi\equiv\nobreak 1$. Поэтому мы имеем гомеоморфизмы
$$
T_{r}:H^{\sigma,(r)}(\Omega)\leftrightarrow
K_{\sigma,(r)}(\Omega,\Gamma),\quad\sigma\in\{s-\varepsilon,s+\delta\}.
$$
Отсюда, поскольку
$$
T_{r}\,u=(u_{0},u_{1},\ldots,u_{r})\;\Leftrightarrow\;
T_{r,p}\,u=(\,u_{0},\,\{u_{k}:\,p+1\leq k\leq r\}\,),
$$
вытекает, что ограниченный оператор \eqref{4.79} биективный. Следовательно, он~---
гомеоморфизм (по теореме Банаха об обратном операторе).

Применим к \eqref{4.79} интерполяцию с параметром $\psi$. В силу теорем \ref{th2.5},
\ref{th2.21} и \ref{th3.10} получаем гомеоморфизм
\begin{equation}\label{4.80}
T_{r,p}:\,X_{\psi}\,\leftrightarrow\,
H^{s,\varphi,(0)}(\Omega)\,\oplus\bigoplus_{p+1\leq k\leq
r}\,H^{s-k+1/2,\varphi}(\Gamma).
\end{equation}
Отсюда вытекает неравенство, обратное к \eqref{4.78}:
\begin{gather*}
\|u\|_{X_{\psi}}\leq \\
c_{2} \biggl(\|u\|_{H^{s,\varphi,(0)}(\Omega)}^{2}+\sum_{p+1\leq k\leq
r}\|(D_{\nu}^{k-1}u)\!\upharpoonright\!\Gamma\|
_{H^{s-k+1/2,\varphi}(\Gamma)}^{2}\biggr)^{1/2}\leq \\
c_{2}\,\|u\|_{H^{s,\varphi,(r)}(\Omega)}
\end{gather*}
для всех $u\in C^{\infty}(\,\overline{\Omega}\,)$. Здесь $c_{2}$~--- норма
оператора, обратного к \eqref{4.80}. Таким образом, нормы в пространствах $X_{\psi}$
и $H^{s,\varphi,(r)}(\Omega)$ эквивалентны на множестве
$C^{\infty}(\,\overline{\Omega}\,)$. Оно плотно в $H^{s,\varphi,(r)}(\Omega)$ по
определению и в $X_{\psi}$ в силу теоремы~\ref{th2.1}. Следовательно,
$X_{\psi}=H^{s,\varphi,(r)}(\Omega)$ с точностью до эквивалентных норм. Утверждение
(ii) доказано.

Докажем утверждение (i). Согласно определению пространства
$H^{s,\varphi,(r)}(\Omega)$ отображение \eqref{4.73} продолжается по непрерывности
до изометрического оператора
\begin{equation}\label{4.81}
T_{r}:\,H^{s,\varphi,(r)}(\Omega)\rightarrow \Pi_{s,\varphi,(r)}(\Omega,\Gamma).
\end{equation}
На основании теоремы \ref{th3.5} имеем включение
$$
T_{r}(H^{s,\varphi,(r)}(\Omega))\subseteq K_{s,\varphi,(r)}(\Omega,\Gamma).
$$
Докажем обратное включение. Пусть
$$
(u_{0},u_{1},\ldots,u_{r})\in K_{s,\varphi,(r)}(\Omega,\Gamma).
$$
В силу \eqref{4.80} и равенства $X_{\psi}=H^{s,\varphi,(r)}(\Omega)$ мы имеем
гомеоморфизм
$$
T_{r,p}:\,H^{s,\varphi,(r)}(\Omega)\,\leftrightarrow\,
H^{s,\varphi,(0)}(\Omega)\,\oplus\bigoplus_{p+1\leq k\leq
r}\,H^{s-k+1/2,\varphi}(\Gamma).
$$
Поэтому найдется такое $u\in H^{s,\varphi,(r)}(\Omega)$, что
$$
T_{r,p}\,u=(u_{0},\,\{u_{k}:\,p+1\leq k\leq r\}\,).
$$
Отсюда на основании теоремы \ref{th3.5} получаем равенство
$T_{r}\,u=(u_{0},u_{1},\ldots,u_{r})$. Тем самым, доказано включение
$$
K_{s,\varphi,(r)}(\Omega,\Gamma)\subseteq T_{r}(H^{s,\varphi,(r)}(\Omega)).
$$
Таким образом,
$$
T_{r}(H^{s,\varphi,(r)}(\Omega))= K_{s,\varphi,(r)}(\Omega,\Gamma),
$$
что вместе с изометрическим оператором \eqref{4.81} влечет изометрический изоморфизм
\eqref{4.74}. Утверждение (i) доказано.

Теорема \ref{th4.12} доказана.

Утверждение (i) теоремы \ref{th4.12} дает полезное изометрическое описание
пространства $H^{s,\varphi,(r)}(\Omega)$ при $s\notin E_{r}$:
$$
T_{r}(H^{s,\varphi,(r)}(\Omega))=K_{s,\varphi,(r)}(\Omega,\Gamma).
$$


\begin{remark}\label{rem4.7}
Если $s\in E_{r}$, то отображение \eqref{4.73} продолжается по непрерывности до
ограниченного оператора \eqref{4.81}, но можно утверждать лишь только, что
$$
T_{r}(H^{s,\varphi,(r)}(\Omega))\subseteq K_{s,\varphi,(r)}(\Omega,\Gamma).
$$
Это следует из формулы \eqref{4.70}, теоремы \ref{th4.12} (i) и нижеследующей
интерполяционной леммы.
\end{remark}


\begin{lemma}\label{lem4.3}
Для произвольных $\sigma\in\mathbb{R}$, $\varepsilon>0$ и $\varphi\in\mathcal{M}$
справедливы следующие равенства пространств с точностью до эквивалентности норм в
них:
\begin{gather}\label{4.82}
[H^{\sigma-\varepsilon,\varphi,(0)}(\Omega),
H^{\sigma+\varepsilon,\varphi,(0)}(\Omega)]_{t^{1/2}}=
H^{\sigma,\varphi,(0)}(\Omega) \\
[H^{\sigma-\varepsilon,\varphi}(\Gamma),
H^{\sigma+\varepsilon,\varphi}(\Gamma)]_{t^{1/2}}=H^{\sigma,\varphi}(\Gamma).
\label{4.83}
\end{gather}
\end{lemma}


\textbf{Доказательство.} Мы выведем равенство \eqref{4.82} из интерполяционных
теорем \ref{th3.10} и \ref{th2.3}. Согласно первой из них
$$
H^{\sigma\mp\varepsilon,\varphi,(0)}(\Omega)= [H^{\sigma-2\varepsilon,(0)}(\Omega),
H^{\sigma+2\varepsilon,(0)}(\Omega)]_{\psi_{\mp}}.
$$
Здесь интерполяционные параметры $\psi_{\mp}$ определяются по формулам
$$
\psi_{-}(t):=t^{1/4} \varphi(t^{1/(4\varepsilon)}),\quad \psi_{+}(t):=t^{3/4}
\varphi(t^{1/(4\varepsilon)})\quad\mbox{при}\quad t\geq1
$$
и $\psi_{\mp}(t):=1$ при $0<t<1$. Отсюда в силу теоремы \ref{th2.3} о повторной
интерполяции получаем:
\begin{gather} \notag
[H^{\sigma-\varepsilon,\varphi,(0)}(\Omega),
H^{\sigma+\varepsilon,\varphi,(0)}(\Omega)]_{t^{1/2}}=\\ \notag
\bigl[\,[H^{\sigma-2\varepsilon,(0)}(\Omega),
H^{\sigma+2\varepsilon,(0)}(\Omega)]_{\psi_{-}},\\
 [H^{\sigma-2\varepsilon,(0)}(\Omega),
H^{\sigma+2\varepsilon,(0)}(\Omega)]_{\psi_{+}}\,\bigr]_{t^{1/2}}=\notag\\
[H^{\sigma-2\varepsilon,(0)}(\Omega), H^{\sigma+2\varepsilon,(0)}(\Omega)]_{\psi}.
\label{4.84}
\end{gather}
Здесь интерполяционный параметр $\psi$ определяется по формулам
$$
\psi(t):=\psi_{-}(t)\,(\psi_{+}(t)/\psi_{-}(t))^{1/2}= t^{1/2}
\varphi(t^{1/(4\varepsilon)})\quad\mbox{при}\quad t\geq1
$$
и $\psi(t)=1$ при $0<t<1$. Поэтому на основании теоремы \ref{th3.10}
\begin{equation}\label{4.85}
[H^{\sigma-2\varepsilon,(0)}(\Omega), H^{\sigma+2\varepsilon,(0)}(\Omega)]_{\psi}=
H^{\sigma,\varphi,(0)}(\Omega).
\end{equation}
Теперь равенства \eqref{4.84}, \eqref{4.85} влекут \eqref{4.82}. Так как полученные
в доказательстве равенства пространств выполняются с точностью до эквивалентности
норм, то это верно и для равенства \eqref{4.82}. Равенство \eqref{4.83} является
частным случаем теоремы \ref{th2.23}. Лемма \ref{lem4.3} доказана.

\medskip

Доказанная лемма пригодится нам и в дальнейшем (в тех местах, где будет использована
интерполяционная формула \eqref{4.70}).

Продолжим изучение свойств модифицированной уточненной шкалы.


\begin{theorem} \label{th4.13}
Пусть $r\in\mathbb{N}$, $s\in\mathbb{R}$ и $\varphi,\varphi_{1}\in\mathcal{M}$.
Справедливы следующие утверждения.
\begin{itemize}
\item [$\mathrm{(i)}$] Гильбертово пространство $H^{s,\varphi,(r)}(\Omega)$
сепарабельное.
\item [$\mathrm{(ii)}$] Множество $C^{\infty}(\,\overline{\Omega}\,)$ плотно в
пространстве $H^{s,\varphi,(r)}(\Omega)$
\item [$\mathrm{(iii)}$] Если $s>r-1/2$, то нормы в пространствах
$H^{s,\varphi,(r)}(\Omega)$ и $H^{s,\varphi}(\Omega)$ эквивалентны на плотном
множестве $C^{\infty}(\,\overline{\Omega}\,)$, и следовательно, эти пространства
равны с точностью до эквивалентности норм.
\item [$\mathrm{(iv)}$] Для произвольного числа $\varepsilon>0$ выполняется
непрерывное и плотное вложение
$H^{s+\varepsilon,\varphi_{1},(r)}(\Omega)\hookrightarrow
H^{s,\varphi,(r)}(\Omega)$. Это вложение компактно.
\item [$\mathrm{(v)}$] Если функция $\varphi/\varphi_{1}$ ограничена в окрестности
$\infty$, то выполняется непрерывное и плотное вложение
$H^{s,\varphi_{1},(r)}(\Omega)\hookrightarrow H^{s,\varphi,(r)}(\Omega)$. Оно
компактно, если $\varphi(t)/\varphi_{1}(t)\rightarrow0$ при
$t\rightarrow\nobreak\infty$.
\end{itemize}
\end{theorem}


\begin{remark} \label{rem4.8}
Непрерывное вложение $H^{s+\varepsilon,\varphi_{1},(r)}(\Omega)\hookrightarrow
H^{s,\varphi,(r)}(\Omega)$, фигурирующее в теореме \ref{th4.13} (iv), (v) при
$\varepsilon\geq0$, понимается в следующем смысле [\ref{BerezanskyUsSheftel90},
с.~504]:
\begin{itemize}
\item [$\mathrm{(i)}$] cуществует число $c>0$ такое, что
$$
\|u\|_{H^{s+\varepsilon,\varphi_{1},(r)}(\Omega)}\leq
c\,\|u\|_{H^{s,\varphi,(r)}(\Omega)}\quad\mbox{для всех}\quad u\in
C^{\infty}(\,\overline{\Omega}\,);
$$
\item [$\mathrm{(ii)}$] тождественное отображение на функциях $u\in
C^{\infty}(\,\overline{\Omega}\,)$ продолжается по непрерывности до
\emph{инъективного} оператора, действующего из
$H^{s+\varepsilon,\varphi_{1},(r)}(\Omega)$ в $H^{s,\varphi,(r)}(\Omega)$.
\end{itemize}
\end{remark}

\textbf{Доказательство теоремы \ref{th4.13}.} (i) Для параметра $s\notin E_{r}$
сепарабельность пространства $H^{s,\varphi,(r)}(\Omega)$ вытекает из теоремы
\ref{th4.12} (i) и сепарабельности пространства $K_{s,\varphi,(r)}(\Omega,\Gamma)$.
Если $s\in E_{r}$, то пространство $H^{s,\varphi,(r)}(\Omega)$ сепарабельно в силу
\eqref{4.70} как результат интерполяции сепарабельных гильбертовых пространств.

(ii) Утверждение (ii) в случае $s\notin E_{r}$ содержится в определении пространства
$H^{s,\varphi,(r)}(\Omega)$. Если $s\in E_{r}$, то в силу \eqref{4.70} и теоремы
\ref{th2.1} справедливо непрерывное плотное вложение
$H^{s+\varepsilon,\varphi,(r)}(\Omega)\hookrightarrow H^{s,\varphi,(r)}(\Omega)$ для
достаточно малого $\varepsilon>0$. Поскольку $s+\varepsilon\notin E_{r}$, множество
$C^{\infty}(\,\overline{\Omega}\,)$ плотно в пространстве
$H^{s+\varepsilon,\varphi,(r)}(\Omega)$. Следовательно, это множество плотно и в
пространстве $H^{s,\varphi,(r)}(\Omega)$.

(iii) Если $s>r-1/2$ и $k\in\{1,\ldots,r\}$, то в силу теоремы \ref{th3.5}
$$
\|(D_{\nu}^{k-1}u)\!\upharpoonright\!\Gamma\|_{H^{s-k+1/2,\varphi}(\Gamma)}\leq
c\,\|u\|_{H^{s,\varphi}(\Omega)}
$$
для любого $u\in C^{\infty}(\,\overline{\Omega}\,)$, где число $c>0$ не зависит от
$u$. Поэтому нормы в пространствах $H^{s,\varphi,(r)}(\Omega)$ и
$H^{s,\varphi,(0)}(\Omega)=H^{s,\varphi}(\Omega)$ эквивалентны на плотном линейном
многообразии $C^{\infty}(\,\overline{\Omega}\,)$. Значит, эти пространства равны.

(iv) В силу теорем \ref{th3.9} (iv) и \ref{th2.22} (iii) выполняется компактное
вложение
$$
\Pi_{s+\varepsilon,\varphi_{1},(r)}(\Omega,\Gamma)\hookrightarrow
\Pi_{s,\varphi,(r)}(\Omega,\Gamma).
$$
Оно влечет компактное вложение подпространств
$$
K_{s+\varepsilon,\varphi_{1},(r)}(\Omega,\Gamma)\hookrightarrow
K_{s,\varphi,(r)}(\Omega,\Gamma).
$$
Отсюда на основании теоремы \ref{th4.12} (i) получаем в случае
$s,s+\varepsilon\notin E_{r}$ компактное вложение пространства
$$
H^{s+\varepsilon,\varphi_{1},(r)}(\Omega)=
T_{r}^{-1}(K_{s+\varepsilon,\varphi_{1},(r)}(\Omega,\Gamma))
$$
в пространство
$$
H^{s,\varphi,(r)}(\Omega)=T_{r}^{-1}(K_{s,\varphi,(r)}(\Omega,\Gamma)).
$$
Если $\{s,s+\varepsilon\}\cap E_{r}\neq\varnothing$, то в силу \eqref{4.70} и
теоремы \ref{th2.1} справедливы непрерывные вложения для достаточно малого числа
$\varepsilon_{0}>0$:
\begin{gather*}
H^{s+\varepsilon,\varphi_{1},(r)}(\Omega)\hookrightarrow
H^{s+\varepsilon-\varepsilon_{0},\varphi_{1},(r)}(\Omega)\hookrightarrow \\
H^{s+\varepsilon_{0},\varphi,(r)}(\Omega)\hookrightarrow H^{s,\varphi,(r)}(\Omega).
\end{gather*}
Здесь число $\varepsilon_{0}$ должно удовлетворять условиям
$$
0<\varepsilon_{0}<1,\quad\varepsilon_{0}<\varepsilon/2\quad\mbox{и}\quad
s+\varepsilon-\varepsilon_{0},s+\varepsilon_{0}\notin E_{r}.
$$
По уже доказанному, среднее вложение компактно. Следовательно, вложение крайних
пространств также компактно. Это вложение плотно в силу утверждения (ii).

(v) Предположим, что функция $\varphi/\varphi_{1}$ ограничена в окрестности
$\infty$. Тогда в силу утверждения (iv) теорем \ref{th2.22}, \ref{th3.3} и
\ref{th3.8} имеем непрерывное вложение
$$
K_{s,\varphi_{1},(r)}(\Omega,\Gamma)\hookrightarrow
K_{s,\varphi,(r)}(\Omega,\Gamma).
$$ Отсюда на основании теоремы \ref{th4.12} (i) получаем требуемое непрерывное вложение в
случае $s\notin E_{r}$:
\begin{gather}\notag
H^{s,\varphi_{1},(r)}(\Omega)=
T_{r}^{-1}(K_{s,\varphi_{1},(r)}(\Omega,\Gamma))\hookrightarrow\\
T_{r}^{-1}(K_{s,\varphi,(r)}(\Omega,\Gamma))=H^{s,\varphi,(r)}(\Omega).\label{4.86}
\end{gather}
Если $s\in E_{r}$, то в силу интерполяционной формулы \eqref{4.70} имеем:
\begin{gather}\notag
H^{s\mp1/2,\varphi_{1},(r)}(\Omega)\hookrightarrow
H^{s\mp1/2,\varphi,(r)}(\Omega)\Rightarrow \\
H^{s,\varphi_{1},(r)}(\Omega)\hookrightarrow H^{s,\varphi,(r)}(\Omega).\label{4.87}
\end{gather}
Здесь левые непрерывные вложения уже доказаны, поскольку $s\mp1/2\notin E_{r}$.
Следовательно, выполняется правое непрерывное вложение. Он плотно при любом
$s\in\mathbb{R}$ согласно утверждению (ii).

Далее, если $\varphi(t)/\varphi_{1}(t)\rightarrow0$ при
$t\rightarrow\nobreak\infty$, то в силу утверждения (iv) теорем \ref{th2.22},
\ref{th3.3} и \ref{th3.8} вложение
$$
\Pi_{s,\varphi_{1},(r)}(\Omega,\Gamma)\hookrightarrow
\Pi_{s,\varphi,(r)}(\Omega,\Gamma)
$$
компактно. Отсюда следует компактность вложения подпространств
$$
K_{s,\varphi_{1},(r)}(\Omega,\Gamma)\hookrightarrow
K_{s,\varphi,(r)}(\Omega,\Gamma).
$$
Значит, компактно вложение \eqref{4.86} в случае $s\notin E_{r}$. Теперь, если $s\in
E_{r}$, то в импликации \eqref{4.87} левые вложения компактны, что влечет
компактность правого вложения в силу \eqref{4.70} и теоремы о том, что при
интерполяции со степенным параметром наследуется компактность операторов
[\ref{Triebel80}, с.~135]. Утверждение (v) доказано.

Теорема \eqref{th4.13} доказана.

\medskip

В заключение этого пункта изучим свойства дифференциальных выражений как операторов
на пространстве $H^{s,\varphi,(r)}(\Omega)$. Пусть $K=K(x,D)$~--- линейное
дифференциальное выражение, заданное на $\overline{\Omega}$, а $R=R(x,D)$~---
граничное линейное дифференциальное выражение, заданное на $\Gamma$. Коэффициенты
этих выражений --- бесконечно гладкие комплекснозначные функции, а порядки
произвольные.


\begin{theorem}\label{th4.14}
Пусть $r\in\mathbb{N}$, $s\in\mathbb{R}$ и $\varphi\in\mathcal{M}$. Справедливы
следующие утверждения.
\begin{itemize}
\item [$\mathrm{(i)}$] Если $\varkappa:=\mathrm{ord}\,K\leq r$, то
$$
\|Ku\|_{H^{s-\varkappa,\varphi,(0)}(\Omega)}\leq
c_{1}\,\|u\|_{H^{s,\varphi,(r)}(\Omega)}\quad\mbox{для всех}\quad u\in
C^{\infty}(\,\overline{\Omega}\,),
$$
где число $c_{1}>0$ не зависит от $u$. Поэтому отображение $u\mapsto Ku$, где $u\in
C^{\infty}(\,\overline{\Omega}\,)$, продолжается по непрерывности до ограниченного
линейного оператора
$$
K:\,H^{s,\varphi,(r)}(\Omega)\rightarrow H^{s-\varkappa,\varphi,(0)}(\Omega).
$$
\item [$\mathrm{(ii)}$] Если $\varrho:=\mathrm{ord}\,R\leq r-1$, то
$$
\|Ru\|_{H^{s-\varrho-1/2,\varphi}(\Gamma)}\leq
c_{2}\,\|u\|_{H^{s,\varphi,(r)}(\Omega)}\;\;\;\mbox{для всех}\;\;\;u\in
C^{\infty}(\,\overline{\Omega}\,),
$$
где число $c_{2}>0$ не зависит от $u$. Поэтому отображение $u\mapsto Ru$, где $u\in
C^{\infty}(\,\overline{\Omega}\,)$, продолжается по непрерывности до ограниченного
линейного оператора
$$
R:\,H^{s,\varphi,(r)}(\Omega)\rightarrow H^{s-\varrho-1/2,\varphi}(\Gamma).
$$
\end{itemize}
\end{theorem}


\textbf{Доказательство.} В соболевском случае $\varphi\equiv1$ эта теорема доказана
Я.~А.~Ройтбергом [\ref{Roitberg96}, с.~74] (лемма 2.3.1). Отсюда случай
произвольного $\varphi\in\mathcal{M}$ выводится с помощью интерполяции. Сделаем это,
например, для утверждения (i). Доказательство утверждения (ii) аналогично.

Предположим сначала, что $s\notin E_{r}$. Пусть положительное число
$\varepsilon=\nobreak\delta$ такое как в теореме \ref{th4.12} (ii). Отображение
$u\mapsto Ku$, где $u\in C^{\infty}(\,\overline{\Omega}\,)$, продолжается по
непрерывности до ограниченных линейных операторов
$$
K:\,H^{s\mp\varepsilon,(r)}(\Omega)\rightarrow
H^{s\mp\varepsilon-\varkappa,(0)}(\Omega).
$$
Применим теперь интерполяцию с функциональным параметром $\psi$ из теоремы
\ref{th2.14}. В~силу теорем \ref{th4.12} (ii) и \ref{th3.10} получаем утверждение
(i) в случае $s\notin E_{r}$.

Предположим теперь, что $s\in E_{r}$. Выберем произвольное число
$\varepsilon\in(0,1)$. Так как $s\mp\varepsilon\notin E_{r}$, то по уже доказанному
имеем ограниченные линейные операторы
$$
K:\,H^{s\mp\varepsilon,\varphi,(r)}(\Omega)\rightarrow
H^{s\mp\varepsilon-\varkappa,\varphi,(0)}(\Omega).
$$
Применив здесь интерполяцию со степенным параметром $t^{1/2}$, получим в силу
\eqref{4.70} и леммы \ref{lem4.3} утверждение (i) в случае $s\notin E_{r}$.

Теорема \ref{th4.14} доказана.

\medskip

Из теоремы \ref{th4.14} (i) в случае $\varkappa=0$ следует, что умножение на функцию
из класса $C^{\infty}(\,\overline{\Omega}\,)$ является ограниченным линейным
оператором в каждом из пространств $H^{s,\varphi,(r)}(\Omega)$.

\subsection[Теоремы о нетеровости и полном наборе гомеоморфизмов]
{Теоремы о нетеровости и  \\ полном наборегомеоморфизмов}\label{sec4.2.3}

Изучим регулярную эллиптическую краевую задачу \eqref{4.1} в модифицированной по
Ройтбергу уточненной шкале, где порядок модификации $r=2q$. Результаты этого пункта
и их доказательства полезно сравнить с результатами п.~\ref{sec4.2.1}.

Согласно теореме \ref{th4.14} отображение \eqref{4.2} продолжается по непрерывности
до ограниченного линейного оператора
\begin{gather}\notag
(L,B):H^{s,\varphi,(2q)}(\Omega)\rightarrow \\
H^{s-2q,\varphi,(0)}(\Omega)\oplus\bigoplus_{j=1}^{q}H^{s-m_{j}-1/2,\varphi}(\Gamma)=:
\mathcal{H}_{s,\varphi,(0)}(\Omega,\Gamma).\label{4.88}
\end{gather}
Следовательно, для произвольного элемента $u\in H^{s,\varphi,(2q)}(\Omega)$
определены по замыканию правые части $f\in H^{s-2q,\varphi,(0)}(\Omega)$ и $g_{j}\in
H^{s-m_{j}-1/2,\varphi}(\Gamma)$ краевой задачи \eqref{4.1}. В силу теоремы
\ref{th4.12} (i) равенство $(L,B)u=(f,g_{1},\ldots,g_{q})$ равносильно тому, что
вектор $(u_{0},u_{1},\ldots,u_{2q}):=T_{2q}u$ является обобщенным по Ройтбергу
решением этой краевой задачи. Элемент $u$ мы отождествляем с вектором
$(u_{0},u_{1},\ldots,u_{2q})$ и также называем обобщенным по Ройтбергу решением
краевой задачи \eqref{4.1}.

Изучим свойства оператора \eqref{4.88}.


\begin{theorem}\label{th4.15} Для произвольных параметров
$s\in\mathbb{R}$ и $\varphi\in\mathcal{M}$ ограниченный оператор \eqref{4.88}
нетеров. Его ядро равно $N$, а область значений состоит из всех векторов
$(f,g_{1},\ldots,g_{q})\in\mathcal{H}_{s,\varphi,(0)}(\Omega,\Gamma)$,
удовлетворяющих условию \eqref{4.4}. Индекс оператора \eqref{4.88} равен $\dim
N-\dim N^{+}$ и не зависит от $s$, $\varphi$.
\end{theorem}


\textbf{Доказательство.} В соболевском случае $\varphi\equiv1$ эта теорема
установлена Я.~А.~Ройтбергом [\ref{Roitberg64}, \ref{Roitberg65}]; см. также его
монографию [\ref{Roitberg96}, с.~96, 155] (теоремы 4.1.1 и 5.3.1). Отсюда мы выведем
общий случай $\varphi\in\mathcal{M}$ с помощью интерполяции.

Предположим сначала, что $s\notin E_{2q}$. Пусть положительное число
$\varepsilon=\delta$ такое как в теореме \ref{th4.12} (ii). Отображение \ref{th4.2}
продолжается по непрерывности до ограниченных нетеровых операторов
\begin{equation}\label{4.89}
(L,B):\,H^{s\mp\varepsilon,(2q)}(\Omega)\rightarrow
\mathcal{H}_{s\mp\varepsilon,(0)}(\Omega,\Gamma).
\end{equation}
Они имеют общее ядро $N$, одинаковый индекс $\varkappa:=\dim N-\dim N^{+}$ и области
значений
\begin{gather}\notag
(L,B)(H^{s\mp\varepsilon}(\Omega))=\\
\{(f,g_{1},\ldots,g_{q})\in\mathcal{H}_{s\mp\varepsilon,(0)}(\Omega,\Gamma):\,\mbox{верно
\eqref{4.4}}\}.\label{4.90}
\end{gather}

Применим к пространствам, в которых действуют операторы \eqref{4.89}, интерполяцию с
функциональным параметром $\psi$ из теоремы \ref{th2.14}. В силу теоремы~\ref{th2.7}
получим ограниченный нетеров оператор
\begin{gather}\notag
(L,B):\,[H^{s-\varepsilon,(2q)}(\Omega),
H^{s+\varepsilon,(2q)}(\Omega)]_{\psi}\rightarrow\\
[\mathcal{H}_{s-\varepsilon,(0)}(\Omega,\Gamma),
\mathcal{H}_{s+\varepsilon,(0)}(\Omega,\Gamma)]_{\psi},\label{4.91}
\end{gather}
продолжающий по непрерывности отображение \eqref{4.2}. Здесь в силу интерполяционных
теорем \ref{th2.5}, \ref{th2.21} и \ref{th3.10}
\begin{equation}\label{4.92}
[\mathcal{H}_{s-\varepsilon,(0)}(\Omega,\Gamma),
\mathcal{H}_{s+\varepsilon,(0)}(\Omega,\Gamma)]_{\psi}=
\mathcal{H}_{s,\varphi,(0)}(\Omega,\Gamma)
\end{equation}
с эквивалентностью норм. Отсюда и на основании теоремы \ref{th4.12} (ii) следует,
что отображение \eqref{4.2} продолжается по непрерывности до ограниченного оператора
\eqref{4.88}, равного \eqref{4.91}. Согласно теореме \ref{th2.7} нетеровость
операторов \eqref{4.89} влечет нетеровость оператора \eqref{4.88}, который наследует
их ядро $N$ и индекс $\varkappa$. Кроме того, область значений оператора
\eqref{4.88} равна
$$
\mathcal{H}_{s,\varphi,(0)}(\Omega,\Gamma)\cap(L,B)(H^{s-\varepsilon,(2q)}(\Omega)).
$$
Отсюда, в силу \eqref{4.90} получаем, что она такая как в формулировке настоящей
теоремы.

Предположим теперь, что $s\in E_{2q}$. Выберем произвольное число
$\varepsilon\in(0,1)$. Так как $s\mp\varepsilon\notin E_{2q}$, то по доказанному
отображение \eqref{4.2} продолжается по непрерывности до ограниченных нетеровых
операторов
$$
(L,B):\,H^{s\mp\varepsilon,\varphi,(2q)}(\Omega)\rightarrow
\mathcal{H}_{s\mp\varepsilon,\varphi,(0)}(\Omega,\Gamma).
$$
Они имеют общие ядро $N$ и индекс $\varkappa$. Применив здесь интерполяцию со
степенным параметром $t^{1/2}$, получим в силу теоремы \ref{th2.7}, равенства
\eqref{4.70} и леммы \ref{lem4.3} ограниченный нетеров оператор \eqref{4.88}. Он
продолжает по непрерывности отображение \eqref{4.2}, имеет те же ядро $N$ и индекс
$\varkappa$. Кроме того, область значений этого оператора равна
$$
\mathcal{H}_{s,\varphi,(0)}(\Omega,\Gamma)\cap(L,B)(H^{s-\varepsilon,\varphi,(2q)}(\Omega)).
$$
Отсюда и из доказанного выше следует, что она такая как в формулировке настоящей
теоремы.

Теорема \ref{th4.15} доказана.

\medskip

Теорема \ref{th4.15} является общей теоремой о разрешимости эллиптической краевой
задачи \eqref{4.1}, поскольку в ней область определения $H^{s,\varphi,(2q)}(\Omega)$
оператора $(L,B)$ не зависит от эллиптического выражения $L$.

Если $N=N^{+}=\{0\}$ (дефект краевой задачи отсутствует), то оператор \eqref{4.88}
является гомеоморфизмом пространства $H^{s,\varphi,(2q)}(\Omega)$ на
$\mathcal{H}_{s,\varphi,(0)}(\Omega,\Gamma)$. Это следует из теоремы \ref{th4.15} и
теоремы Банаха об обратном операторе. В~общем случае гомеоморфизм удобно задавать с
помощью следующих проекторов.


\begin{lemma}\label{lem4.4}
Для произвольных параметров $s\in\mathbb{R}$ и $\varphi\in\mathcal{M}$ справедливы
следующие разложения пространств $H^{s,\varphi,(2q)}(\Omega)$ и
$\mathcal{H}_{s,\varphi,(0)}(\Omega,\Gamma)$ в прямые суммы замкнутых
подпространств:
\begin{gather}\notag
H^{s,\varphi,(2q)}(\Omega)=\\ \label{4.93} N\dotplus\{u\in
H^{s,\varphi,(2q)}(\Omega):\,(u_{0},w)_{\Omega}=0\;\,\mbox{для всех}\,\;w\in N\},\\
\notag \mathcal{H}_{s,\varphi,(0)}(\Omega,\Gamma)=\\ \{(v,0,\ldots,0):\;v\in
N^{+}\}\dotplus (L,B)(H^{s,\varphi,(2q)}(\Omega)). \label{4.94}
\end{gather}
Здесь $u_{0}$ --- начальная компонента вектора
$(u_{0},u_{1},\ldots,u_{2q}):=T_{2q}u$. Обозначим через $P_{2q}$ косой проектор
пространства $H^{s,\varphi,(2q)}(\Omega)$ на второе слагаемое суммы \eqref{4.93}, а
через $Q^{+}_{0}$ косой проектор пространства
$\mathcal{H}_{s,\varphi,(0)}(\Omega,\Gamma)$ на второе слагаемое суммы \eqref{4.94}
параллельно первому слагаемому. Эти проекторы не зависят от $s$ и $\varphi$.
\end{lemma}


\textbf{Доказательство.} Сначала докажем равенство \eqref{4.93}. Из определения
пространства $H^{s,\varphi,(2q)}(\Omega)$ вытекает, что отображение $u\mapsto u_{0}$
является ограниченным оператором
$$
T_{0}:H^{s,\varphi,(2q)}(\Omega)\rightarrow H^{s,\varphi,(0)}(\Omega).
$$
Поэтому, второе слагаемое суммы \eqref{4.93} замкнуто в $H^{s,\varphi,(2q)}(\Omega)$
и имеет тривиальное пересечение с $N$. Следовательно, справедливо разложение
\begin{gather}\notag
H^{s,\varphi,(0)}(\Omega)=\\ N\dotplus\{u_{0}\in
H^{s,\varphi,(0)}(\Omega):\,(u_{0},w)_{\Omega}=0\;\,\mbox{для всех}\,\;w\in
N\}.\label{4.95}
\end{gather}
В самом деле, пространство, двойственное к подпространству $N$ пространства
$H^{-s,1/\varphi,(0)}(\Omega)$, будет гомеоморфно факторпространству пространства
$H^{s,\varphi,(0)}(\Omega)$ по второму слагаемому суммы \eqref{4.95}. Значит, это
слагаемое имеет коразмерность, равную $\dim N$, что влечет \eqref{4.95}.

Обозначим через $\Pi$ косой проектор пространства $H^{s,\varphi,(0)}(\Omega)$ на
первое слагаемое суммы  \eqref{4.95} параллельно второму слагаемому. Для
произвольного элемента $u\in H^{s,\varphi,(2q)}(\Omega)$ запишем $u=u'+u''$, где
$u':=\Pi u_{0}\in N$, а
$$
u'':=u-\Pi u_{0}\in H^{s,\varphi,(2q)}(\Omega)
$$
удовлетворяет условию
$$
(u''_{0},w)_{\Omega}=(u_{0}-\Pi u_{0},w)_{\Omega}=0\quad\mbox{для любого}\quad w\in
N.
$$
Равенство \eqref{4.93} доказано.

Равенство \eqref{4.94} вытекает из того, что в силу теоремы \ref{th4.15} линейные
многообразия, записанные в левой части этого равенства, замкнутые, имеют тривиальное
пересечение, и конечная размерность первого из них совпадает с коразмерностью
второго. Наконец, независимость проекторов $P_{2q}$ и $Q^{+}_{0}$ от параметров
$s,\varphi$ вытекает из включений $N,N^{+}\subset
C^{\infty}(\,\overline{\Omega}\,)$.

Лемма \ref{lem4.4} доказана.


\begin{theorem} \label{th4.16}
Для произвольных параметров $s\in\mathbb{R}$ и $\varphi\in\mathcal{M}$ сужение
отображения $\eqref{4.88}$ на подпространство $P_{2q}(H^{s,\varphi}(\Omega))$
является гомеоморфизмом
\begin{equation}\label{4.96}
(L,B):\,P_{2q}(H^{s,\varphi}(\Omega))\leftrightarrow
Q^{+}_{0}(\mathcal{H}_{s,\varphi,(0)}(\Omega,\Gamma)).
\end{equation}
\end{theorem}


\textbf{Доказательство.} Согласно теореме \ref{th4.15}, $N$~--- ядро, а
$Q^{+}_{0}(\mathcal{H}_{s,\varphi,(0)}(\Omega,\Gamma))$~--- область значений
оператора \eqref{4.88}. Следовательно, ограниченный оператор \eqref{4.96}~---
биекция. Поэтому он является гомеоморфизмом в силу теоремы Банаха об обратном
операторе. Теорема~\ref{th4.16} доказана.

\medskip

При каждом фиксированном $\varphi\in\mathcal{M}$ набор гомеоморфизмов \eqref{4.96}
полный, поскольку $s$ пробегает всю вещественную ось. Этим он отличается от
одностороннего набора гомеоморфизмов, фигурирующего в теореме \ref{th4.2}.

Из теоремы \ref{th4.16} вытекает следующая априорная оценка решения эллиптической
краевой задачи \eqref{4.1}.


\begin{theorem}\label{th4.17} Пусть $s\in\mathbb{R}$,
$\varphi\in\mathcal{M}$ и число $\varepsilon>0$. Предположим, что элемент $u\in
H^{s,\varphi,(2q)}(\Omega)$ является обобщенным по Ройтбергу решением краевой задачи
\eqref{4.1}, где
$(f,g_{1},\ldots,g_{q})\in\mathcal{H}_{s,\varphi,(0)}(\Omega,\Gamma)$. Тогда
выполняется оценка
\begin{gather}\notag
\|u\|_{H^{s,\varphi,(2q)}(\Omega)}\leq \\
c\,\bigl(\,\|(f,g_{1},\ldots,g_{q})\|_{\mathcal{H}_{s,\varphi,(0)}(\Omega,\Gamma)}+
\|u\|_{H^{s-\varepsilon,\varphi,(2q)}(\Omega)}\,\bigr),\label{4.97}
\end{gather}
в которой число $c=c(s,\varphi,\varepsilon)>0$ не зависит от $u$ и вектора
$(f,g_{1},\ldots,g_{q})$.
\end{theorem}


\textbf{Доказательство.} Воспользуемся разложением \eqref{4.93} и запишем элемент
$u\in H^{s,\varphi,(2q)}(\Omega)$ в виде $u=u'+u''$, где
$$
u':=(1-P_{2q})u\in N\quad\mbox{и}\quad u'':=P_{2q}u\in
P_{2q}(H^{s,\varphi}(\Omega)).
$$
В силу теоремы \ref{th4.16}
\begin{gather*}
\|u''\|_{H^{s,\varphi,(2q)}(\Omega)}\leq
c_{2}\|(L,B)u''\|_{\mathcal{H}_{s,\varphi,(0)}(\Omega,\Gamma)}=\\
c_{2}\|(L,B)u\|_{\mathcal{H}_{s,\varphi,(0)}(\Omega,\Gamma)}=
c_{2}\|(f,g_{1},\ldots,g_{q})\|_{\mathcal{H}_{s,\varphi,(0)}(\Omega,\Gamma)}.
\end{gather*}
Здесь $c_{2}$ --- норма оператора, обратного к \eqref{4.96}. Кроме того, так как
пространство $N$ конечномерно и $1-P_{2q}$ суть проектор на $N$ в пространстве
$H^{s-\varepsilon,\varphi,(2q)}(\Omega)$ (лемма \ref{lem4.4}), то
$$
\|u'\|_{H^{s,\varphi,(2q)}(\Omega)}\leq
c_{0}\|u'\|_{H^{s-\varepsilon,\varphi,(2q)}(\Omega)}\leq
c_{1}\|u\|_{H^{s-\varepsilon,\varphi,(2q)}(\Omega)}.
$$
Здесь положительные числа $c_{0}$ и $c_{1}$ не зависят от $u'$, $u$ и
$(f,g_{1},\ldots,g_{q})$. Сложив эти неравенства, получаем оценку \eqref{4.97}.
Теорема \ref{th4.17} доказана.

\medskip

Если $N=\{0\}$, т.~е.  краевая задача \eqref{4.1} имеет не более одного решения, то
величина $\|u\|_{H^{s-\varepsilon,\varphi,(2q)}(\Omega)}$ в правой части априорной
оценки \eqref{4.97} отсутствует.

В соболевском случае $\varphi\equiv1$ теоремы \ref{th4.15} -- \ref{th4.17} доказаны
Я.~А.~Ройтбергом [\ref{Roitberg64}, \ref{Roitberg65}], [\ref{Roitberg96}] (п.~3.3,
5.3); см. также монографию [\ref{Berezansky65}, с.~248] и обзор [\ref{Agranovich97},
с.~84]. В случае произвольного $\varphi\in\mathcal{M}$ эти теоремы уточняют
результаты Ройтберга применительно к модифицированной уточненной шкале \eqref{4.71}.

В силу теорем \ref{th4.13} (iii) и \ref{th3.9} (i) выполняются следующие равенства
пространств с точностью до эквивалентности норм в них:
\begin{gather}\notag
H^{s,\varphi,(2q)}(\Omega)=H^{s,\varphi}(\Omega)\quad\mbox{и}\\
H^{s-2q,\varphi,(0)}(\Omega)=H^{s-2q,\varphi}(\Omega)\quad\mbox{при}\quad
s>2q-1/2.\label{4.98}
\end{gather}
Поэтому теоремы \ref{th4.15} -- \ref{th4.17} содержат в себе теоремы \ref{th4.1} --
\ref{th4.3} о разрешимости эллиптической краевой задачи \eqref{4.1} в шкале
позитивных пространств Хермандера. Более того, справедлив следующий факт.


\begin{theorem}\label{th4.18}
Теоремы $\ref{th4.1}$ -- $\ref{th4.3}$ верны для произвольных параметров $s>2q-1/2$
и $\varphi\in\mathcal{M}$.
\end{theorem}

Это следует из теорем \ref{th4.15} -- \ref{th4.17}, формул \eqref{4.98} и равенства
проекторов $P_{2q}=P$ на $H^{s,\varphi,(2q)}(\Omega)$.

\subsection[Гладкость обобщенного решения вплоть до границы]
{Гладкость обобщенного решения \\ вплоть до границы}\label{sec4.2.4}

Изучим гладкость обобщенного (по Ройтбергу) решения краевой задачи \eqref{4.1} в
модифицированной уточненной шкале пространств. Результаты этого пункта полезно
сравнить с результатами п.~\ref{sec4.1.2}.

Для целого $r\geq0$ обозначим через $H^{-\infty,(r)}(\Omega)$ объединение всех
пространств $H^{s,\varphi,(r)}(\Omega)$, где $s\in\mathbb{R}$ и
$\varphi\in\mathcal{M}$. В пространстве $H^{-\infty,(r)}(\Omega)$ вводится топология
индуктивного предела.

Ограниченные операторы \eqref{4.88}, рассмотренные для всех параметров
$s\in\nobreak\mathbb{R}$ и $\varphi\in\mathcal{M}$, порождают линейный непрерывный
оператор
\begin{gather}\label{4.99}
(L,B):\,H^{-\infty,(2q)}(\Omega)\rightarrow
H^{-\infty,(0)}(\Omega)\times(\mathcal{D}'(\Gamma))^{q}=:\\
\mathcal{H}_{-\infty,(0)}(\Omega,\Gamma). \notag
\end{gather}
В силу теоремы \ref{th4.15} ядро этого оператора равно $N$, а область значений
состоит из всех векторов
$(f,g_{1},\ldots,g_{q})\in\mathcal{H}_{-\infty,(0)}(\Omega,\Gamma)$, удовлетворяющих
условию \eqref{4.4}.


\begin{theorem}\label{th4.19} Предположим, что элемент
$u\in H^{-\infty,(2q)}(\Omega)$ является обобщенным решением краевой задачи
\eqref{4.1}, где
$$
f\in H^{s-2q,\varphi,(0)}(\Omega)\quad\mbox{и}\quad g_{j}\in
H^{s-m_{j}-1/2,\varphi}(\Gamma)\quad\mbox{при}\quad j=1,\ldots,q
$$
для некоторых параметров $s\in\mathbb{R}$, $\varphi\in\mathcal{M}$. Тогда $u\in
H^{s,\varphi,(2q)}(\Omega)$.
\end{theorem}


\textbf{Доказательство.} В силу указанных выше свойств оператора \eqref{4.99} вектор
$F:=(f,g_{1},\ldots,g_{q})=(L,B)u$ удовлетворяет условию \eqref{4.4}. По условию,
$F\in\mathcal{H}_{s,\varphi}(\Omega,\Gamma)$; следовательно, согласно теореме
\ref{th4.15} верно включение $F\in(L,B)(H^{s,\varphi}(\Omega))$. Поэтому наряду с
условием $(L,B)u=F$ справедливо равенство $(L,B)v=F$ для некоторого элемента $v\in
H^{s,\varphi}(\Omega)$. Отсюда, $(L,B)(u-v)=0$, что влечет включение
$$
w:=u-v\in N\subset H^{s,\varphi}(\Omega).
$$
Таким образом, $u=v+w\in H^{s,\varphi}(\Omega)$. Теорема \ref{th4.19} доказана.

\medskip

Теорема \ref{th4.19} --- это утверждение о глобальной гладкости обобщенного решения
(т.е. во всей области замкнутой области $\overline{\Omega}$). В соболевском случае
$\varphi\equiv1$ она доказана  Я.~А.~Ройтбергом [\ref{Roitberg64},
\ref{Roitberg65}], [\ref{Roitberg96}, с.~214] (теорема 7.1.1).

Теперь рассмотрим случай локальной гладкости. Пусть $U$ --- открытое множество в
$\mathbb{R}^{n}$, имеющее непустое пересечение с областью $\Omega$. Как и в
п.~\ref{sec4.1.2} положим $\Omega_{0}:=U\cap\Omega$ и $\Gamma_{0}:=U\cap\Gamma$
(возможен случай $\Gamma_{0}=\varnothing$). Введем следующий локальный аналог
пространства $H^{\sigma,\varphi,(r)}(\Omega)$, где $\sigma\in\mathbb{R}$,
$\varphi\in\mathcal{M}$ и целое $r\geq0$. Положим
\begin{gather*}
H^{\sigma,\varphi,(r)}_{\mathrm{loc}}(\Omega_{0},\Gamma_{0}):= \{u\in
H^{-\infty,(r)} (\Omega):\\\chi\,u\in H^{\sigma,\varphi,(r)}(\Omega)\;\,\mbox{для
всех}\,\;\chi\in
C^{\infty}(\,\overline{\Omega}\,),\;\mathrm{supp}\,\chi\subset\Omega_{0}\cup\Gamma_{0}\}.
\end{gather*}
В силу теоремы \ref{th4.14} (i), умножение на функцию $\chi\in
C^{\infty}(\,\overline{\Omega}\,)$ является ограниченным оператором в пространстве
$H^{\sigma,\varphi,(r)}(\Omega)$. Поэтому для каждого элемента $u\in
H^{-\infty,(r)}(\Omega)$ корректно определено произведение $\chi u\in
H^{-\infty,(r)}(\Omega)$. Кроме того,
$$
H^{\sigma,\varphi,(r)}(\Omega)\subset
H^{\sigma,\varphi,(r)}_{\mathrm{loc}}(\Omega_{0},\Gamma_{0}).
$$

Нам понадобится также локальное пространство
$H^{\sigma,\varphi}_{\mathrm{loc}}(\Gamma_{0})$, введенное в п.~\ref{sec2.6.3b}.


\begin{theorem}\label{th4.20} Предположим, что элемент
$u\in H^{-\infty,(2q)}(\Omega)$ является обобщенным решением краевой задачи
\eqref{4.1}, где
\begin{equation}\label{4.100}
f\in H^{s-2q,\varphi,(0)}_{\mathrm{loc}}(\Omega_{0},\Gamma_{0}),\;\;g_{j}\in
H^{s-m_{j}-1/2,\varphi}_{\mathrm{loc}}(\Gamma_{0})\;\;\mbox{при}\;\;j=1,\ldots,q
\end{equation}
для некоторых параметров $s\in\mathbb{R}$ и $\varphi\in\mathcal{M}$. Тогда верно
включение $u\in H^{s,\varphi,(2q)}_{\mathrm{loc}}(\Omega_{0},\Gamma_{0})$.
\end{theorem}


\textbf{Доказательство.} В соболевском случае $\varphi\equiv1$ эта теорема доказана
Я.~А.~Ройтбергом [\ref{Roitberg64}, \ref{Roitberg65}], [\ref{Roitberg96}, с.~216]
(теорема 7.2.1). Следовательно, так как
\begin{gather*}
f\in H^{s-2q,\varphi,(0)}_{\mathrm{loc}}(\Omega_{0},\Gamma_{0})\subset
H^{s-1/2-2q,(0)}_{\mathrm{loc}}(\Omega_{0},\Gamma_{0}),\\
g_{j}\in H^{s-m_{j}-1/2,\varphi}_{\mathrm{loc}}(\Gamma_{0})\subset
H^{s-1/2-m_{j}-1/2}_{\mathrm{loc}}(\Gamma_{0}),
\end{gather*}
то
\begin{equation}\label{4.101}
u\in H^{s-1/2,(2q)}_{\mathrm{loc}}(\Omega_{0},\Gamma_{0})\subset
H^{s-1,\varphi,(2q)}_{\mathrm{loc}}(\Omega_{0},\Gamma_{0})
\end{equation}
Покажем, что отсюда и из условия \eqref{4.100} следует требуемое свойство решения:
$u\in H^{s,\varphi}_{\mathrm{loc}}(\Omega_{0},\Gamma_{0})$. Будем рассуждать подобно
доказательству теоремы~\ref{th4.5}.

Пусть функции $\chi,\eta\in C^{\infty}(\,\overline{\Omega}\,)$ такие, что
$\mathrm{supp}\,\chi,\,\mathrm{supp}\,\eta\subset\Omega_{0}\cup\Gamma_{0}$ и
$\eta=\nobreak1$ в окрестности $\mathrm{supp}\,\chi$. Переставив оператор умножения
на функцию $\chi$ с дифференциальными операторами $L$ и $B_{j}$, где $j=1,\ldots,q$,
получим равенства \eqref{4.18}, \eqref{4.19}, где дифференциальные выражения $L'$ и
$B'_{j}$ такие как в доказательстве теоремы \ref{th4.5}. Напомним, что
$\mathrm{ord}\,L'\leq2q-1$ и $\mathrm{ord}\,B'_{j}\leq m_{j}-1$.

Из включений \eqref{4.100}, \eqref{4.101} следует в силу теоремы \ref{th4.14}, что
правые части эллиптической краевой задачи \eqref{4.18}, \eqref{4.19} имеют следующую
гладкость:
$$
\chi\,f+L'(\eta u)\in H^{s-2q,\varphi,(0)}(\Omega),\quad \chi\,g_{j}+B'_{j}(\eta
u)\in H^{s-m_{j}-1/2,\varphi}(\Gamma).
$$
Поэтому решение $\chi u$ этой задачи принадлежит пространству
$H^{s,\varphi,(2q)}(\Omega)$ на основании теоремы \ref{th4.19}. Таким образом, $u\in
H^{s,\varphi,(2q)}_{\mathrm{loc}}(\Omega_{0},\Gamma_{0})$ ввиду произвольности
функции $\chi\in C^{\infty}(\,\overline{\Omega}\,)$, удовлетворяющей условию
$\mathrm{supp}\,\chi\subset\Omega_{0}\cup\Gamma_{0}$.

Теорема \ref{th4.20} доказана.

\medskip

Теоремы \ref{th4.19} и \ref{th4.20} уточняют применительно к шкале пространств
Хермандера теоремы Я.~А.~Ройтберга о повышении гладкости обобщенного решения
эллиптической краевой задачи в модифицированной шкале пространств Соболева
[\ref{Roitberg64}, \ref{Roitberg65}], [\ref{Roitberg96}] (п.~7.1, 7.2). Из равенств
\eqref{4.98} следует, что теоремы \ref{th4.19} и \ref{th4.20} содержат в себе
теоремы \ref{th4.4} и \ref{th4.5}. В теореме \ref{th4.20} отметим случай
$\Gamma_{0}=\varnothing$, который приводит к утверждению о повышении локальной
гладкости обобщенного решения в окрестностях внутренних точек области $\Omega$.

В качестве приложения теорем \ref{th4.19} и \ref{th4.20} установим достаточное
условие того, что обобщенное по Ройтбергу решение $u\in H^{-\infty,(2q)}(\Omega)$
эллиптической краевой задачи \eqref{4.1} является классическим, т.~е. удовлетворяет
условию
\begin{equation}\label{4.102}
u\in H^{\sigma+2q}(\Omega)\cap C^{2q}(\Omega)\cap C^{m}(\,\overline{\Omega}\,),
\end{equation}
где $\sigma>-1/2$ и $m:=\max\{m_{1},\ldots,m_{q}\}$. Это условие возникает следующим
образом. В силу теорем \ref{th4.13} (iii) и \ref{th3.5}, из включения
$$
u\in H^{\sigma+2q,(2q)}(\Omega)=H^{\sigma+2q}(\Omega)
$$
вытекает, что элемент $u$ является решением задачи \eqref{4.1} в смысле теории
распределений, заданных в области $\Omega$. Теперь корректно рассматривать включение
$u\in C^{2q}(\Omega)\cap C^{m}(\,\overline{\Omega}\,)$. Оно означает, что в
\eqref{4.1} функции $Lu$ и $B_{j}u$ вычисляются с помощью классических производных,
т. е. решение $u$ классическое.


\begin{theorem}\label{th4.21}
Предположим, что элемент $u\in H^{-\infty,(2q)}(\Omega)$ является обобщенным по
Ройтбергу решением задачи  \eqref{4.1}, где
\begin{gather}\label{4.103}
f\in H^{n/2,\varphi,(0)}_{\mathrm{loc}}(\Omega,\varnothing)\cap
H^{m-2q+n/2,\varphi,(0)}(\Omega)\cap H^{\sigma,(0)}(\Omega), \\
g_{j}\in H^{\,m-m_{j}+(n-1)/2,\varphi}(\Gamma)\cap
H^{\sigma+2q-m_{j}-1/2}(\Gamma)\;\;\mbox{при}\;\; j=1,\ldots,q \label{4.104}
\end{gather}
для некоторых числа $\sigma>-1/2$ и функционального параметра
$\varphi\in\mathcal{M}$, удовлетворяющего условию \eqref{2.37}. Тогда решение $u$
классическое, т.~е. удовлетворяет включению \eqref{4.102}.
\end{theorem}


\textbf{Доказательство.} В силу теорем \ref{th4.19} и \ref{th4.20}, из условий
\eqref{4.103} и \eqref{4.104} вытекает включение
$$
u\in H^{2q+n/2,\varphi,(2q)}_{\mathrm{loc}}(\Omega,\varnothing)\cap
H^{m+n/2,\varphi,(2q)}(\Omega)\cap H^{\sigma+2q,(2q)}(\Omega).
$$
Отсюда на основании теорем \ref{th4.13} (iii) и \ref{th3.4} имеем:
\begin{gather*}
u\in H^{m+n/2,\varphi,(2q)}(\Omega)\cap H^{\sigma+2q,(2q)}(\Omega)=\\
H^{m+n/2,\varphi}(\Omega)\cap H^{\sigma+2q}(\Omega)\subseteq
C^{m}(\,\overline{\Omega}\,)\cap H^{\sigma+2q}(\Omega).
\end{gather*}
(Последнее равенство становится ясным, если рассмотреть отдельно случаи
$m+n/2\geq\sigma+2q$ и $m+n/2<\sigma+2q$.) Кроме того,
$$
\chi u\in H^{2q+n/2,\varphi,(2q)}(\Omega)= H^{2q+n/2,\varphi}(\Omega)\subset
C^{2q}(\,\overline{\Omega}\,)
$$
для любой функции $\chi\in C^{\infty}_{0}(\Omega)$, что влечет включение $u\in
C^{2q}(\Omega)$. Таким образом, выполняется условие \eqref{4.102}, т.~е. $u$~---
классическое решение. Теорема \ref{th4.21} доказана.

\subsection[Интерполяция в модифицированной шкале]
{Интерполяция \\ в модифицированной шкале}\label{sec4.2.5}

Изучим два вопроса, связанных с интерполяцией в модифицированной уточненной шкале.
Во-первых, мы докажем, что правая часть равенства \eqref{4.70} не зависит от выбора
параметра $\varepsilon$. Во-вторых, мы установим, что интерполяционная формула
\eqref{4.76} остается верной при существенно более слабых условиях на входящие в нее
параметры.


\begin{theorem} \label{th4.22}
Пусть $r\in\mathbb{N}$, $s\in E_{r}$ и $\varphi\in\mathcal{M}$. Пространство
$$
H^{s,\varphi,(r)}(\Omega,\varepsilon):=
\bigl[\,H^{s-\varepsilon,\varphi,(r)}(\Omega),
H^{s+\varepsilon,\varphi,(r)}(\Omega)\,\bigr]_{1/2}
$$
не зависит с точностью до эквивалентности норм от параметра $\varepsilon\in(0,1)$.
\end{theorem}


\textbf{Доказательство.} Предположим сначала, что $r=2q$
--- четное число. Рассмотрим регулярную эллиптическую краевую задачу в области
$\Omega$:
\begin{gather}\label{4.105}
Lu\equiv(1-\Delta)^{q}\,u=f\quad\mbox{в}\quad\Omega, \\
B_{j}u\equiv D^{j-1}_{\nu}\,u=g_{j}\quad\text{на}\quad\Gamma,\quad j=1,\ldots,q.
\label{4.106}
\end{gather}
Как обычно, $\Delta$ --- оператор Лапласа. Эта задача формально самосопряженная и
для нее $N=N^{+}=\{0\}$. Согласно теореме \ref{th4.15} мы имеем гомеоморфизм
$$
(L,B):\,H^{s,\varphi,(2q)}(\Omega,\varepsilon)\leftrightarrow
H^{s-2q,\varphi,(0)}(\Omega)\oplus\bigoplus_{j=1}^{q}\,H^{s-j+1/2,\varphi}(\Gamma)
$$
при $0<\varepsilon<1$. Отсюда сразу следует теорема для четного $r=2q$.

Предположим далее, что число $r$ нечетное. В силу теоремы \ref{th4.12} (i), для
любого числа $\sigma\in(-\infty,r+1/2)\setminus E_{r}$ мы имеем изометрические
изоморфизмы
\begin{gather*}
T_{r}:\,H^{\sigma,\varphi,(r)}(\Omega)\leftrightarrow
K_{\sigma,\varphi,(r)}(\Omega,\Gamma),\\
T_{r+1}:\,H^{\sigma,\varphi,(r+1)}(\Omega)\leftrightarrow
K_{\sigma,\varphi,(r+1)}(\Omega,\Gamma)=\\
K_{\sigma,\varphi,(r)}(\Omega,\Gamma)\oplus H^{\sigma-r-1/2,\varphi}(\Gamma).
\end{gather*}
Поэтому композиция отображений
$$
u\mapsto T_{r+1}\,u=:(u_{0},u_{1},\ldots,u_{r},u_{r+1})\mapsto
(T_{r}^{-1}(u_{0},u_{1},\ldots,u_{r}),u_{r+1}),
$$
где $u\in H^{\sigma,\varphi,(r+1)}(\Omega)$, определяет изометрический изоморфизм
\begin{equation}\label{4.107}
T:\,H^{\sigma,\varphi,(r+1)}(\Omega)\leftrightarrow
H^{\sigma,\varphi,(r)}(\Omega)\oplus H^{\sigma-r-1/2,\varphi}(\Gamma).
\end{equation}
Возьмем здесь $\sigma=s\mp\varepsilon$, где $0<\varepsilon<1$, и применим
интерполяцию со степенным параметром $t^{1/2}$. Ввиду леммы \ref{lem4.3} получаем
гомеоморфизм
$$
T:\,H^{s,\varphi,(r+1)}(\Omega,\varepsilon)\leftrightarrow
H^{s,\varphi,(r)}(\Omega,\varepsilon)\oplus
H^{s-r-1/2,\varphi}(\Gamma)=:X(\varepsilon).
$$
Следовательно,
$$
\|u\|_{H^{s,\varphi,(r)}(\Omega,\varepsilon)}= \|(u,0)\|_{X(\varepsilon)}\asymp
\|T^{-1}(u,0)\|_{H^{s,\varphi,(r+1)}(\Omega,\varepsilon)}.
$$
Отсюда, поскольку параметр $r+1$ четный, вытекает, по доказанному, что нормы в
пространствах $H^{s,\varphi,(r)}(\Omega,\varepsilon)$, где $0<\varepsilon<1$,
эквивалентны. Значит, эти пространства равны, поскольку множество
$C^{\infty}(\,\overline{\Omega}\,)$ плотно в каждом из них согласно теореме
\ref{th4.13} (ii).

Теорема \ref{th4.22} доказана.


\begin{theorem} \label{th4.23}
Пусть заданы $r\in\mathbb{N}$, $s\in\mathbb{R}$, $\varphi\in\mathcal{M}$ и
положительные числа $\varepsilon,\delta$. Если число $r$ нечетное, дополнительно
предположим, что выполняется хотя бы одно из неравенств $s-\varepsilon>r-1/2$ и
$s+\delta<r+1/2$. Тогда верна интерполяционная формула \eqref{4.76}:
\begin{equation}\label{4.108}
[H^{s-\varepsilon,(r)}(\Omega),H^{s+\delta,(r)}(\Omega)]_{\psi}=
H^{s,\varphi,(r)}(\Omega)
\end{equation}
с эквивалентностью норм. Здесь $\psi$~--- интерполяционный параметр из теоремы
$\ref{th2.14}$.
\end{theorem}


\textbf{Доказательство.} Предположим сначала, что $r=2q$~--- четное число. В силу
теоремы \ref{th4.15} имеем для эллиптической краевой задачи \eqref{4.105},
\eqref{4.106} гомеоморфизмы
\begin{gather}\notag
(L,B):\,H^{s,\varphi,(2q)}(\Omega)\leftrightarrow\\ \label{4.109}
H^{s-2q,\varphi,(0)}(\Omega)\oplus\bigoplus_{j=1}^{q}\,H^{s-j+1/2,\varphi}(\Gamma),\\
\notag (L,B):\,H^{\sigma,(2q)}(\Omega)\leftrightarrow\\
H^{\sigma-2q,(0)}(\Omega)\oplus\bigoplus_{j=1}^{q}\,H^{\sigma-j+1/2}(\Gamma),
\quad\sigma\in\mathbb{R}. \label{4.110}
\end{gather}
Применим интерполяцию с параметром $\psi$ к пространствам, в которых действуют
гомеоморфизмы \eqref{4.110} для $\sigma\in\{s-\varepsilon,\,s+\delta\}$. Так как
$\psi$ --- интерполяционный параметр, то в силу теорем \ref{th3.10} и \ref{th2.21}
получим еще один гомеоморфизм
\begin{gather}\notag
(L,B):[H^{s-\varepsilon,(2q)}(\Omega), H^{s+\delta,(2q)}(\Omega)]_{\psi}\leftrightarrow\\
H^{s-2q,\varphi,(0)}(\Omega)\oplus\bigoplus_{j=1}^{q}\,H^{s-j+1/2,\varphi}(\Gamma).
\label{4.111}
\end{gather}
Теперь из гомеоморфизмов \eqref{4.109} и \eqref{4.111} следует равенство пространств
\eqref{4.108} вместе с эквивалентностью норм в них.

Предположим далее, что число $r$ нечетное. Рассмотрим отдельно случаи
$s-\varepsilon>r-1/2$ и $s+\delta<r+1/2$.

Если $s-\varepsilon>r-1/2$, то в силу теорем \ref{th4.13} (iii) и \ref{th3.2}
справедливы следующие равенства пространств с точностью до эквивалентности норм в
них:
\begin{gather*}
[H^{s-\varepsilon,(r)}(\Omega),H^{s+\delta,(r)}(\Omega)]_{\psi}=
[H^{s-\varepsilon}(\Omega),H^{s+\delta}(\Omega)]_{\psi}= \\
H^{s,\varphi}(\Omega)=H^{s,\varphi,(r)}(\Omega).
\end{gather*}
Формула \eqref{4.108} доказана в рассматриваемом случае.

Рассмотрим теперь случай $s+\delta<r+1/2$. Заметим, что гомеоморфизм \eqref{4.107}
верен для любого $\sigma<r+1/2$. Если $\sigma\notin E_{r}$, то этот гомеоморфизм
установлен в доказательстве теоремы \ref{th4.22}. Отсюда он следует для $\sigma\in
E_{r}$ посредством интерполяции в силу равенств \eqref{4.70} и \eqref{4.83}.
В~частности, мы имеем гомеоморфизмы
\begin{gather}\label{4.112}
T:\,H^{s,\varphi,(r+1)}(\Omega)\leftrightarrow H^{s,\varphi,(r)}(\Omega)\oplus
H^{s-r-1/2,\varphi}(\Gamma),\\
T:\,H^{\sigma,(r+1)}(\Omega)\leftrightarrow
H^{\sigma,(r)}(\Omega)\oplus H^{\sigma-r-1/2}(\Gamma),\label{4.113}\\
\sigma\in\{s-\varepsilon,\,s+\delta\}.\notag
\end{gather}

Применим к \eqref{4.113} интерполяцию с параметром $\psi$. В силу уже доказанного
(так как число $r+1$ четное) и теоремы \ref{th2.21} получим еще один гомеоморфизм
\begin{gather}\notag
T:\,H^{s,\varphi,(r+1)}(\Omega)\leftrightarrow\\
[H^{s-\varepsilon,(r)}(\Omega),H^{s+\delta,(r)}(\Omega)]_{\psi}\oplus
H^{s-r-1/2,\varphi}(\Gamma).\label{4.114}
\end{gather}
Теперь, воспользовавшись гомеоморфизмами \eqref{4.112} и \eqref{4.114}, можем
записать следующее:
\begin{gather*}
u\in[H^{s-\varepsilon,(r)}(\Omega),H^{s+\delta,(r)}(\Omega)]_{\psi}\,\Leftrightarrow\\
T^{-1}(u,0)\in H^{s,\varphi,(r+1)}(\Omega)\,\Leftrightarrow\, u\in
H^{s,\varphi,(r)}(\Omega).
\end{gather*}
Значит, верно равенство пространств \eqref{4.108} в случае $s+\delta<r+1/2$. Нормы в
этих пространствах эквивалентны:
\begin{gather*}
\|u\|_{[H^{s-\varepsilon,(r)}(\Omega),H^{s+\delta,(r)}(\Omega)]_{\psi}}\asymp\\
\|T^{-1}(u,0)\|_{H^{s,\varphi,(r+1)}(\Omega)}\asymp
\|u\|_{H^{s,\varphi,(r)}(\Omega)},
\end{gather*}
где $u\in H^{s,\varphi,(r)}(\Omega)$.

Теорема \ref{th4.23} доказана.





\section{Некоторые свойства модифицированной уточненной шкалы}\label{sec4.3}

\markright{\emph \ref{sec4.3}. Свойства модифицированной уточненной шкалы}

В этом пункте мы сформулируем и докажем два важных свойства модифицированной
уточненной шкалы пространств $H^{s,\varphi,(2q)}(\Omega)$, которая имеет четный
порядок модификации $2q$, где $q\in\mathbb{N}$. Они будут использованы в
п.~\ref{sec4.5} при доказательстве индивидуальных теорем о разрешимости
эллиптической краевой задачи.

\subsection{Формулировки результатов}\label{sec4.3.1}

Первое свойство дает эквивалентное альтернативное определение пространства
$H^{s,\varphi,(2q)}(\Omega)$. Оказывается, норма в нем эквивалентна некоторой норме,
содержащей правильно эллиптическое выражение $L$ и не включающей в себя граничные
значения функций. При этом, в отличие от \eqref{4.69}, не надо исключать полуцелые
значения $s\in E_{2q}$.


\begin{theorem}\label{th4.24}
Пусть $s\in\mathbb{R}$ и $\varphi\in\mathcal{M}$. Справедливы следующие утверждения.
\begin{itemize}
\item [$\mathrm{(i)}$] На множестве функций $u\in C^{\infty}(\,\overline{\Omega}\,)$
норма в пространстве $H^{s,\varphi,(2q)}(\Omega)$ эквивалентна норме графика
\begin{equation}\label{4.115}
\bigl(\,\|u\|_{H^{s,\varphi,(0)}(\Omega)}^{2}+
\|Lu\|_{H^{s-2q,\varphi,(0)}(\Omega)}^{2}\,\bigr)^{1/2}.
\end{equation}
Поэтому пространство $H^{s,\varphi,(2q)}(\Omega)$ совпадает с пополнением линейной
системы $C^{\infty}(\,\overline{\Omega}\,)$ по норме \eqref{4.115}.
\item [$\mathrm{(ii)}$] Отображение
\begin{equation}\label{4.116}
I_{L}:\,u\mapsto(u,Lu),\quad u\in C^{\infty}(\,\overline{\Omega}\,),
\end{equation}
продолжается по непрерывности до гомеоморфизма
\begin{equation}\label{4.117}
I_{L}:\,H^{s,\varphi,(2q)}(\Omega)\leftrightarrow K_{s,\varphi,L}(\Omega).
\end{equation}
Здесь
\begin{gather}\notag
K_{s,\varphi,L}(\Omega):=\{\,(u_{0},f):\,u_{0}\in
H^{s,\varphi,(0)}(\Omega),\;f\in H^{s-2q,\varphi,(0)}(\Omega),\\
(u_{0},L^{+}w)_{\Omega}=(f,w)_{\Omega}\;\,\mbox{для всех}\,\;w\in
C^{\infty}_{\nu,2q}(\,\overline{\Omega}\,)\,\} \label{4.118}
\end{gather}
--- замкнутое подпространство в
$$
H^{s,\varphi,(0)}(\Omega)\oplus H^{s-2q,\varphi,(0)}(\Omega).
$$
\end{itemize}
\end{theorem}


Напомним (см. п.~\ref{sec3.5.1}), что
$$
C^{\infty}_{\nu,2q}(\,\overline{\Omega}\,):=\{u\in
C^{\infty}(\,\overline{\Omega}\,):\,D^{j-1}_{\nu}u=0\;\;\mbox{на}\;\;\Gamma,
\;\;j=1,\ldots,2q\}.
$$

Если $s>-1/2$, то для распределений $u_{0}$, $f$ из теоремы \ref{th4.24} (ii)
условие \eqref{4.118} эквивалентно условию
\begin{equation}\label{4.119}
(u_{0},L^{+}w)_{\Omega}=(f,w)_{\Omega}\quad\mbox{для всех}\quad w\in
C^{\infty}_{0}(\Omega),
\end{equation}
которое означает, что $Lu_{0}=f$ в области $\Omega$. (Эквивалентность мы докажем
ниже в лемме \ref{lem4.5}.) Если $s<-1/2$, то это не так, и нам понадобится второе
свойство.


\begin{theorem}\label{th4.25}
Пусть $s<-1/2$, $s+1/2\notin\mathbb{Z}$, и $\varphi\in\mathcal{M}$. Тогда для
произвольных распределений
\begin{equation}\label{4.120}
u_{0}\in H^{s,\varphi,(0)}(\Omega),\quad f\in H^{s-2q,\varphi,(0)}(\Omega),
\end{equation}
удовлетворяющих условию \eqref{4.119}, существует единственная пара
$(u_{0}^{\ast},f)\in K_{s,\varphi,L}(\Omega)$ такая, что
\begin{equation}\label{4.121}
(u_{0},w)_{\Omega}=(u_{0}^{\ast},w)_{\Omega}\quad\mbox{для всех}\quad w\in
C^{\infty}_{0}(\Omega).
\end{equation}
При этом выполняется неравенство
\begin{equation}\label{4.122}
\|u_{0}^{\ast}\|_{H^{s,\varphi,(0)}(\Omega)}\leq c\,
\bigl(\,\|u_{0}\|_{H^{s,\varphi,(0)}(\Omega)}^{2}+
\|f\|_{H^{s-2q,\varphi,(0)}(\Omega)}^{2}\,\bigr)^{1/2},
\end{equation}
в котором число $c=c(s,\varphi)>0$ не зависит от $u_{0}$, $f$.
\end{theorem}


В соболевском случае $\varphi\equiv1$ теоремы \ref{th4.24} и \ref{th4.25} доказаны
Я.~А.~Ройтбергом [\ref{Roitberg71}], [\ref{Roitberg96}, с. 174, 187] (теоремы 6.1.1
и 6.2.1). Им использовались несколько иные, но равносильные нашим, формулировки.
Опишем их отличия.

В формулировке теоремы \ref{th4.24} Ройтберг [\ref{Roitberg96}, с. 174] (теорема
6.1.1) использовал вместо \eqref{4.118} условие (при $\varphi\equiv1$)
\begin{equation}\label{4.123}
(u_{0},L^{+}w)_{\Omega}=(f,w)_{\Omega}\;\;\mbox{для всех}\;\;w\in
H^{2q,1/\varphi}_{0}(\Omega)\cap H^{2q-s,1/\varphi,(0)}(\Omega).
\end{equation}
Это равносильная замена, как мы покажем ниже в лемме \ref{lem4.5}

Далее, если воспользоваться теоремой \ref{th4.24} и положить
$u^{\ast}:=I_{L}^{-1}(u_{0}^{\ast},f)$, то теорему \ref{th4.25} можно
переформулировать следующим эквивалентным образом. Для произвольных распределений
\eqref{4.120}, удовлетворяющих условию \eqref{4.121}, существует единственный
элемент $u^{\ast}\in H^{s,\varphi,(2q)}(\Omega)$ такой, что
$$
u_{0}^{\ast}=u_{0}\;\;\mbox{в}\;\;\Omega\quad\mbox{и}\quad
Lu^{\ast}=f\;\;\mbox{в}\;\;H^{s-2q,\varphi,(0)}(\Omega).
$$
При этом
$$
\|u^{\ast}\|_{H^{s,\varphi,(2q)}(\Omega)}\leq c\,
\bigl(\,\|u_{0}\|_{H^{s,\varphi,(0)}(\Omega)}^{2}+
\|f\|_{H^{s-2q,\varphi,(0)}(\Omega)}^{2}\bigr)^{1/2}.
$$
Здесь $u_{0}^{\ast}$ --- начальная компонента вектора $T_{2q}u^{\ast}$. Эта
равносильная формулировка теоремы  \ref{th4.25} и была использована Ройтбергом в
[\ref{Roitberg96}, с. 187] (теорема 6.2.1) при $\varphi\equiv1$.

\subsection{Доказательства результатов}\label{sec4.3.2}

Мы выведем теоремы \ref{th4.24} и \ref{th4.25} из соболевского случая
$\varphi\equiv1$ с помощью интерполяции. Перед этим докажем лемму \ref{lem4.5},
упоминавшуюся ранее.


\begin{lemma}\label{lem4.5}
Пусть $s\in\mathbb{R}$ и $\varphi\in\mathcal{M}$. Тогда для произвольно заданных
распределений \eqref{4.120} условия \eqref{4.118} и \eqref{4.123} равносильны. Если
$s>-1/2$, то равносильны условия \eqref{4.118} и \eqref{4.119}. Они равносильны
также при $s=-1/2$ в случае $\varphi\equiv1$.
\end{lemma}


\textbf{Доказательство.} Сначала покажем, что
$\eqref{4.118}\Leftrightarrow\eqref{4.123}$. В силу теоремы \ref{3.20} (i) имеем:
$\eqref{4.123}\Rightarrow\eqref{4.118}$. Докажем обратную импликацию. Предположим,
что выполняется условие \eqref{4.118}. Рассмотрим отдельно случаи $s\geq0$ и $s<0$.

Случай $s\geq0$. В силу теоремы \ref{3.20} (i) имеем:
\begin{equation}\label{4.124}
H^{2q,1/\varphi}_{0}(\Omega)\cap H^{2q-s,1/\varphi,(0)}(\Omega)=
H^{2q,1/\varphi}_{0}(\Omega)=H^{2q,1/\varphi}_{\nu,2q}(\Omega).
\end{equation}
Аппроксимируем произвольное распределение $w$ из пространства \eqref{4.124}
некоторой последовательностью функций $w_{j}\in
C^{\infty}_{\nu,2q}(\,\overline{\Omega}\,)$ по норме пространства
$H^{2q,1/\varphi}(\Omega)$. Согласно условию \eqref{4.118}
\begin{equation}\label{4.125}
(u_{0},L^{+}w_{j})_{\Omega}=(f,w_{j})_{\Omega}\quad\mbox{для всех}\quad
j\in\mathbb{N}.
\end{equation}
Отсюда с помощью предельного перехода при $j\rightarrow\infty$ получается условие
\eqref{4.123}. Это вытекает из включений \eqref{4.120} и справедливости следующих
пределов:
\begin{gather*}
\lim_{j\rightarrow\infty}L^{+}w_{j}=L^{+}w\quad\mbox{в}\;\;
H^{0,1/\varphi}(\Omega)\hookrightarrow H^{-s,1/\varphi,(0)}(\Omega), \\
\lim_{j\rightarrow\infty}w_{j}=w\quad\mbox{в}\;\;H^{2q,1/\varphi}(\Omega)
\hookrightarrow H^{2q-s,1/\varphi,(0)}(\Omega).
\end{gather*}

Случай $s<0$. В силу теоремы \ref{3.20} (i) и формулы \eqref{3.119}
\begin{equation}\label{4.126}
H^{2q,1/\varphi}_{0}(\Omega)\cap H^{2q-s,1/\varphi,(0)}(\Omega)=
H^{2q-s,1/\varphi}_{\nu,2q}(\Omega).
\end{equation}
Аппроксимируем произвольное распределение $w$ из пространства \eqref{4.126}
некоторой последовательностью функций $w_{j}\in
C^{\infty}_{\nu,2q}(\,\overline{\Omega}\,)$ по норме пространства
$H^{2q-s,1/\varphi}(\Omega)$. Согласно условию \eqref{4.118} имеем равенства
\eqref{4.125}. Переходя в них к пределу при $j\rightarrow\infty$, снова получаем
\eqref{4.123}. Это вытекает из \eqref{4.120} и следующих пределов:
\begin{gather*}
\lim_{j\rightarrow\infty}L^{+}w_{j}=L^{+}w\quad\mbox{в}\quad
H^{-s,1/\varphi}(\Omega)= H^{-s,1/\varphi,(0)}(\Omega), \\
\lim_{j\rightarrow\infty}w_{j}=w\quad\mbox{в}\quad
H^{2q-s,1/\varphi}(\Omega)=H^{2q-s,1/\varphi,(0)}(\Omega).
\end{gather*}

Мы доказали, что $\eqref{4.118}\Leftrightarrow\eqref{4.123}$ при любом
$s\in\mathbb{R}$.

Теперь предположим, что $s>-1/2$ и докажем эквивалентность условий \eqref{4.118} и
\eqref{4.119}. Импликация $\eqref{4.118}\Rightarrow\eqref{4.119}$ очевидна. Докажем
обратную импликацию. Предположим, что выполняется условие \eqref{4.119}. В силу уже
доказанного нам достаточно показать, что $\eqref{4.119}\Rightarrow\eqref{4.123}$. Из
теоремы \ref{th3.20} (i) и неравенства $s>-1/2$ следует равенство
\begin{equation}\label{4.127}
H^{2q,1/\varphi}_{0}(\Omega)\cap H^{2q-s,1/\varphi,(0)}(\Omega)=
H^{\lambda,1/\varphi}_{0}(\Omega),
\end{equation}
где $\lambda:=\max\{2q,\,2q-s\}$. Поэтому произвольно выбранное распределение $w$ из
пространства \eqref{4.127} можно аппроксимировать некоторой последовательностью
функций $w_{j}\in C^{\infty}_{0}(\Omega)$ по норме пространства
$H^{\lambda,1/\varphi}(\Omega)$. Согласно условию \eqref{4.119} выполняется
равенство \eqref{4.125}. Отсюда с помощью предельного перехода при
$j\rightarrow\infty$ получается условие \eqref{4.123}. Это вытекает из включений
\eqref{4.120} и следующих пределов:
\begin{gather*}
\lim_{j\rightarrow\infty}L^{+}w_{j}=L^{+}w\quad\mbox{в}\;\;
H^{\lambda-2q,1/\varphi}(\Omega)\hookrightarrow
H^{-s,1/\varphi,(0)}(\Omega), \\
\lim_{j\rightarrow\infty}w_{j}=w\quad\mbox{в}\;\;H^{\lambda,1/\varphi}(\Omega)
\hookrightarrow H^{2q-s,1/\varphi,(0)}(\Omega).
\end{gather*}
Таким образом, мы доказали, что $\eqref{4.118}\Leftrightarrow\eqref{4.119}$ при
$s>-1/2$.

Наконец, если $s=-1/2$ и $\varphi\equiv1$, то равенство \eqref{4.127} остается
верным в силу [\ref{Triebel80}, с.~412] (теорема 4.7.1 (а)). Повторив рассуждения
предыдущего абзаца, получим эквивалентность
$\eqref{4.118}\Leftrightarrow\eqref{4.119}$ и в этом случае.

Лемма \ref{lem4.5} доказана.


Следующая лемма о проекторах потребуется нам для доказательства теоремы
\ref{th4.24}.


\begin{lemma}\label{lem4.6}
Для каждого числа $\sigma\in\mathbb{R}$ cуществует проектор $\Pi_{\sigma}$
пространства $H^{\sigma,(0)}(\Omega)\times H^{\sigma-2q,(0)}(\Omega)$ на
подпространство $K_{\sigma,L}(\Omega)$ такой, что $\Pi_{\sigma}$~--- продолжение
отображения $\Pi_{\lambda}$ при $\sigma<\lambda$.
\end{lemma}


\textbf{Доказательство.} Пусть $\sigma\in\mathbb{R}$. Предварительно установим одно
полезное равенство. С~помощью предельного перехода и теоремы \ref{th4.14} выводим,
что формула Грина \eqref{3.6} остается верной для произвольных распределений $u\in
H^{\sigma,(2q)}(\Omega)$, $v\in H^{2q-\sigma,(2q)}(\Omega)$:
\begin{equation}\label{4.128}
(Lu,v_{0})_{\Omega}+\sum_{j=1}^{q}\;(B_{j}u,\,C_{j}^{+}v)_{\Gamma}
=(u_{0},L^{+}v)_{\Omega}+\sum_{j=1}^{q}\;(C_{j}u,\,B_{j}^{+}v)_{\Gamma}.
\end{equation}
Здесь компоненты полуторалинейных форм $(\cdot,\cdot)_{\Omega}$ и
$(\cdot,\cdot)_{\Gamma}$ принадлежат следующим  пространствам, взаимно двойственным
относительно этих форм:
\begin{gather*}
Lu\in H^{\sigma-2q,(0)}(\Omega),\quad (v_{0},v_{1},\ldots,v_{2q}):=T_{2q}v,\quad
v_{0}\in
H^{2q-\sigma,(0)}(\Omega),\\
(u_{0},u_{1},\ldots,u_{2q}):=T_{2q}u,\quad u_{0}\in
H^{\sigma,(0)}(\Omega),\quad L^{+}v\in H^{-\sigma,(0)}(\Omega),\\
B_{j}u\in H^{\sigma-m_{j}-1/2}(\Gamma),\quad C_{j}^{+}v\in
H^{-\sigma+m_{j}+1/2}(\Gamma),\\
C_{j}u\in H^{\sigma-2q+m^{+}_{j}+1/2}(\Gamma),\quad B_{j}^{+}v\in
H^{2q-\sigma-m^{+}_{j}-1/2}(\Gamma).
\end{gather*}
Напомним, что $\mathrm{ord}\,L=\mathrm{ord}\,L^{+}=2q$, а также ввиду \eqref{3.7}
\begin{gather*}
\mathrm{ord}\,B_{j}=m_{j}\leq2q-1,\;\;\mathrm{ord}\,C_{j}^{+}=2q-1-m_{j},\\
\mathrm{ord}\,B_{j}^{+}=:m_{j}^{+}\leq2q-1,\;\;\mathrm{ord}\,C_{j}=2q-1-m_{j}^{+}.
\end{gather*}

Далее, поскольку $\{B^{+}_{1},\ldots,B^{+}_{q},C^{+}_{1},\ldots,C^{+}_{q}\}$~---
система Дирихле порядка $2q$ (см. [\ref{LionsMagenes71}], ч.~2, теорема 2.1), для
нее справедливо следующее утверждение [\ref{Roitberg96}, с.~176] (лемма 6.1.2).
Ограниченный оператор
\begin{gather*}
\bigl(B^{+}_{1},\ldots,B^{+}_{q},C^{+}_{1},\ldots,C^{+}_{q}\bigr):\,
H^{2q-\sigma,(2q)}(\Omega)\rightarrow \\
\rightarrow\bigoplus_{j=1}^{q} \,H^{2q-\sigma-m^{+}_{j}-1/2}(\Gamma)\oplus
\bigoplus_{j=1}^{q}\,H^{-\sigma+m_{j}+1/2}(\Gamma)
\end{gather*}
имеет линейный ограниченный правый обратный оператор
\begin{gather}\notag
\Phi_{\sigma}:\,\bigoplus_{j=1}^{q}\,H^{2q-\sigma-m^{+}_{j}-1/2}(\Gamma)\oplus
\bigoplus_{j=1}^{q}\,H^{-\sigma+m_{j}+1/2}(\Gamma)\rightarrow\\
H^{2q-\sigma,(2q)}(\Omega)\label{4.129}
\end{gather}
такой, что $\Phi_{\sigma}$ --- сужение оператора $\Phi_{\lambda}$ при
$\sigma<\lambda$.

Для произвольного вектора
\begin{equation}\label{4.130}
h:=(0,\ldots,0,h_{1},\ldots,h_{q})
\in\{0\}^{q}\oplus\bigoplus_{j=1}^{q}\,H^{-\sigma+m_{j}+1/2}(\Gamma)
\end{equation}
положим $v:=\Phi_{\sigma}h\in H^{2q-\sigma,(2q)}(\Omega)$ в формуле Грина
\eqref{4.128}. Поскольку
$$
B^{+}_{j}v=0,\quad C^{+}_{j}v=h_{j}\quad\mbox{при}\quad j=1,\ldots,q,
$$
получим полезное равенство
\begin{equation}\label{4.131}
(u_{0},L^{+}\Phi_{\sigma}h)_{\Omega}-(Lu,(\Phi_{\sigma}h)_{0})_{\Omega}=
\sum_{j=1}^{q}\;(B_{j}u,h_{j})_{\Gamma}.
\end{equation}
Здесь элемент $u\in H^{\sigma,(2q)}(\Omega)$ и вектор $h$ вида \eqref{4.130}
произвольные. Как и прежде,
$$
u_{0}\in H^{\sigma,(0)}(\Omega)\quad\mbox{и}\quad(\Phi_{\sigma}h)_{0}=v_{0}\in
H^{2q-\sigma,(0)}(\Omega)
$$
--- начальные компоненты векторов $T_{2q}u$ и $T_{2q}\Phi_{\sigma}h=T_{2q}v$ соответственно.

Теперь построим проектор $\Pi_{\sigma}$. Сделаем это с помощью нижеследующих пяти
отображений. Пусть выбран произвольный вектор
\begin{equation}\label{4.132}
(u_{0},f)\in H^{\sigma,(0)}(\Omega)\oplus H^{\sigma-2q,(0)}(\Omega).
\end{equation}
Обозначим через $\|\cdot\|_{\sigma}$ норму в пространстве \eqref{4.132}.

Отображения 1 и 2. Рассмотрим следующее разложение пространства
$H^{\sigma,(0)}(\Omega)$ в прямую сумму замкнутых подпространств (см. доказательство
леммы \ref{lem4.4}):
\begin{gather}\notag
H^{\sigma,(0)}(\Omega)=\\N\dotplus\left\{u_{0}'\in
H^{\sigma,(0)}(\Omega):\;(u_{0}',w)_{\Omega}=0\;\,\mbox{для всех}\,\;w\in
N\right\}.\label{4.133}
\end{gather}
Оно существует, поскольку $N$ --- конечномерное подпространство в
$H^{\sigma,(0)}(\Omega)$. Обозначим через $\Psi_{\sigma}$ проектор пространства
$H^{\sigma(0)}(\Omega)$ на подпространство $N$ параллельно второму слагаемому суммы
\eqref{4.133}. Поскольку $N\subset C^{\infty}(\,\overline{\Omega}\,)$, то
$\Psi_{\sigma}$ --- продолжение проектора $\Psi_{\lambda}$ при $\sigma<\lambda$.
Линейные отображения 1 и 2 определяем соответственно по формулам:
\begin{gather}\label{4.134}
(u_{0},f)\mapsto(\Psi_{\sigma}u_{0},0)\in K_{\sigma,L}(\Omega), \\
(u_{0},f)\mapsto(u_{0}-\Psi_{\sigma}u_{0},f)=:(u_{0}',f). \label{4.135}
\end{gather}
Записанное в \eqref{4.134} включение вытекает из свойства $\Psi_{\sigma}u_{0}\in
N\subset C^{\infty}(\,\overline{\Omega}\,)$ и определения пространства
$K_{\sigma,L}(\Omega)$:
\begin{gather*}
\bigl((\Psi_{\sigma}u_{0},L^{+}w)_{\Omega}=(L\Psi_{\sigma}u_{0},w)_{\Omega}=
(0,w)_{\Omega}\;\mbox{для всех}\;w\in C^{\infty}_{\nu,2q}(\,\overline{\Omega})\bigr)
\Rightarrow\\(\Psi_{\sigma}u_{0},0)\in K_{\sigma,L}(\Omega).
\end{gather*}
Здесь первое равенство получено по формуле Грина \eqref{4.64}, при этом интегралы по
$\Gamma$ не появятся, так как $w\in C^{\infty}_{\nu,2q}(\,\overline{\Omega}\,)$.

Отображение 3. Обратимся к формуле \eqref{4.131}. По паре $(u_{0}',f)$ построим
функционал
\begin{equation}\label{4.136}
l_{\sigma}(h):=(u_{0}',L^{+}\Phi_{\sigma}h)_{\Omega}-(f,(\Phi_{\sigma}h)_{0})_{\Omega},\quad
h\mbox{ --- вектор \eqref{4.130}}.
\end{equation}
Этот функционал ограничен в силу следующей цепочки неравенств:
\begin{gather*}
|l_{\sigma}(h)|\leq|(u_{0}',L^{+}\Phi_{\sigma}h)_{\Omega}|+
|(f,(\Phi_{\sigma}h)_{0})_{\Omega}|\leq \\
\|u_{0}'\|_{H^{\sigma,(0)}(\Omega)}\,
\|L^{+}\Phi_{\sigma}h\|_{H^{-\sigma,(0)}(\Omega)}+\\
\|f\|_{H^{\sigma-2q,(0)}(\Omega)}\,
\|(\Phi_{\sigma}h)_{0}\|_{H^{-\sigma+2q,(0)}(\Omega)}\leq\\
\bigl(c_{1}\,\|u_{0}'\|_{H^{\sigma,(0)}(\Omega)}+
c_{2}\,\|f\|_{H^{\sigma-2q,(0)}(\Omega)}\bigr)\,
\|\Phi_{\sigma}h\|_{H^{-\sigma+2q,(2q)}(\Omega)}\leq\\
(c_{1}+c_{2})\,c_{3}\,\|(u_{0}',f)\|_{\sigma}\,
\Bigl(\,\sum_{j=1}^{q}\;\|h_{j}\|_{H^{-\sigma+m_{j}+1/2}(\Gamma)}^{2} \Bigr)^{1/2}.
\end{gather*}
Здесь $c_{1}$, $c_{2}$ и $c_{3}$ --- соответственно нормы операторов
\begin{gather*}
L^{+}:H^{-\sigma+2q,(2q)}(\Omega)\rightarrow H^{-\sigma,(0)}(\Omega), \\
T_{2q}:\,H^{-\sigma+2q,(2q)}(\Omega)\rightarrow H^{-\sigma+2q,(0)}(\Omega)\oplus
\bigoplus_{j=1}^{2q}\,H^{-\sigma+2q-j+1/2}(\Gamma)
\end{gather*}
и оператора \eqref{4.129}. Итак, $l_{\sigma}$ --- антилинейный ограниченный
функционал на пространстве $\bigoplus_{j=1}^{q}\,H^{-\sigma+m_{j}+1/2}(\Gamma)$,
причем его норма удовлетворяет неравенству $\|l_{\sigma}\|\leq
c_{4}\,\|(u_{0}',f)\|_{\sigma}$, где $c_{4}:=(c_{1}+c_{2})\,c_{3}$. В силу теоремы
\ref{th2.22} (v)
\begin{equation}\label{4.137}
\Bigl(\,\exists!\;
g=(g_{1},\ldots,g_{q})\in\bigoplus_{j=1}^{q}\,H^{\sigma-m_{j}-1/2}(\Gamma)\,\Bigr):\,
l_{\sigma}(h)=\sum_{j=1}^{q}\,(g_{j},h_{j})_{\Gamma}.
\end{equation}
При этом
\begin{equation}\label{4.138}
\Bigl(\,\sum_{j=1}^{q}\;\|g_{j}\|_{H^{\sigma-m_{j}-1/2}(\Gamma)}^{2}
\Bigr)^{1/2}\asymp\|l_{\sigma}\|\leq c_{4}\,\|(u_{0}',f)\|_{\sigma}.
\end{equation}

Линейное отображение 3 определяем по формуле
\begin{equation}\label{4.139}
R_{\sigma}:\,(u_{0}',f)\mapsto l_{\sigma}\mapsto g.
\end{equation}
Покажем, что $R_{\sigma}$ --- продолжение отображения $R_{\lambda}$ при
$\sigma<\lambda$. Если
$$
(u_{0}',f)\in H^{\lambda,(0)}(\Omega)\oplus H^{\lambda-2q,(0)}(\Omega),
$$
то наряду с $l_{\sigma}$ определен функционал $l_{\lambda}$. При этом $l_{\sigma}$
--- сужение функционала $l_{\lambda}$, поскольку $\Phi_{\sigma}$~--- сужение
оператора $\Phi_{\lambda}$. В частности, $l_{\sigma}(h)=l_{\lambda}(h)$ для любого
вектора $(h_{1},\ldots,h_{q})\in (C^{\infty}(\Gamma))^{q}$. Отсюда в силу
\eqref{4.137} и \eqref{4.139} получаем $R_{\sigma}(u_{0}',f)=R_{\lambda}(u_{0}',f)$,
т.~е. $R_{\sigma}$~--- продолжение отображения $R_{\lambda}$.

Отображение 4 определяем с помощью гомеоморфизма \eqref{4.96} (теорема \ref{th4.16})
так:
\begin{equation}\label{4.140}
(f,g)\mapsto(L,B)^{-1}Q^{+}(f,g)=:\omega\in P(H^{\sigma,(2q)}(\Omega)).
\end{equation}
Напомним, что $f\in H^{\sigma-2q,(0)}(\Omega)$, а вектор $g$ удовлетворяет условию
\eqref{4.137}. Линейное отображение \eqref{4.140} не зависит от $\sigma$ и
ограничено:
\begin{gather}\notag
\|\omega\|_{H^{\sigma,(2q)}(\Omega)}\leq c_{5}\,
\|Q^{+}(f,g)\|_{\mathcal{H}_{\sigma,(0)}(\Omega,\Gamma)}\leq \\
c_{5}\,c_{6}\,\|(f,g)\|_{\mathcal{H}_{\sigma,(0)}(\Omega,\Gamma)}.\label{4.141}
\end{gather}
Здесь $c_{5}$ --- норма оператора, обратного к \eqref{4.96}, а $c_{6}$~--- норма
проектора $Q^{+}$ в пространстве $\mathcal{H}_{\sigma,(0)}(\Omega,\Gamma)$.

Отображение 5 строим на основании частного случая $\varphi\equiv1$ теоремы
\ref{th4.24}, доказанного Я.~А.~Ройтбергом [\ref{Roitberg96}, с. 174] (теорема
6.1.1):
\begin{equation}\label{4.142}
\omega\mapsto I_{L}\,\omega\in K_{\sigma,L}(\Omega),\quad\omega\in
P(H^{\sigma,(2q)}(\Omega)).
\end{equation}
Это отображение не зависит от $\sigma$ и удовлетворяет двусторонней оценке
\begin{equation}\label{4.143}
\|I_{L}\,\omega\|_{\sigma}\asymp\|\omega\|_{H^{\sigma,(2q)}(\Omega)}.
\end{equation}

Теперь, используя отображения 1--5, определим на векторах \eqref{4.132} оператор
$\Pi_{\sigma}$ по формуле
\begin{equation}\label{4.144}
\Pi_{\sigma}:\,(u_{0},f)\mapsto(\Psi_{\sigma}u_{0},0)+I_{L}\,\omega.
\end{equation}
Здесь, напомним,
\begin{equation}\label{4.145}
\omega=(L,B)^{-1}Q^{+}(f,g),\;\; g=R_{\sigma}(u_{0}',f),\;\;
u_{0}'=u_{0}-\Psi_{\sigma}u_{0}.
\end{equation}
Оператор $\Pi_{\sigma}$ линейный, поскольку линейны отображения 1 -- 5. Он ограничен
в пространстве \eqref{4.132} в силу оценок \eqref{4.138}, \eqref{4.141},
\eqref{4.143} и ограниченности проектора $\Psi_{\sigma}$ в пространстве
$H^{\sigma,(0)}(\Omega)$. Кроме того, включения \eqref{4.134} и \eqref{4.142} влекут
за собой свойство $\Pi_{\sigma}(u_{0},f)\in K_{\sigma,L}(\Omega)$. Таким образом, мы
имеем линейный ограниченный оператор
$$
\Pi_{\sigma}:\,H^{\sigma,(0)}(\Omega)\oplus H^{\sigma-2q,(0)}(\Omega)\rightarrow
K_{\sigma,L}(\Omega).
$$
Так как с уменьшением параметра $\sigma$ отображения 1 -- 5 расширяются, то
$\Pi_{\sigma}$~--- продолжение оператора $\Pi_{\lambda}$ при $\sigma<\lambda$.

Остается показать, что $\Pi_{\sigma}$ --- проектор на подпространство
$K_{\sigma,L}(\Omega)$, т. е. $\Pi_{\sigma}(u_{0},f)=(u_{0},f)$ для любого вектора
$(u_{0},f)\in K_{\sigma,L}(\Omega)$. Выберем произвольно такой вектор. В силу
\eqref{4.134} и \eqref{4.135} имеем включение $(u_{0}',f)\in K_{\sigma,L}(\Omega)$.
Здесь, как и прежде, $u_{0}'=u_{0}-\Psi_{\sigma}u_{0}$. Поэтому согласно частному
случаю $\varphi\equiv1$ теоремы \ref{th4.24}
\begin{equation}\label{4.146}
\bigl(\,\exists!\;u'\in H^{\sigma,(2q)}(\Omega)\,\bigr):\,I_{L}\,u'=(u_{0}',f).
\end{equation}
Последнее равенство означает следующее:
\begin{gather}\label{4.147}
u_{0}'\;\,\mbox{--- начальная компонента вектора}\;\,
T_{2q}u'=(u_{0}',u_{1}',\ldots,u_{2q}'), \\
Lu'=f\quad\mbox{в}\quad H^{\sigma-2q,(0)}(\Omega).\label{4.148}
\end{gather}

Покажем, что $u'=\omega$, где элемент $\omega\in P(H^{\sigma,(2q)}(\Omega))$
определен согласно \eqref{4.145}. Напомним, что в силу \eqref{4.136}, \eqref{4.137}
справедливо следующее равенство для произвольного вектора $h$ вида \eqref{4.130}:
$$
l_{\sigma}(h):=(u_{0}',L^{+}\Phi_{\sigma}h)_{\Omega}-(f,(\Phi_{\sigma}h)_{0})_{\Omega}=
\sum_{j=1}^{q}\,(g_{j},h_{j})_{\Gamma}.
$$
С другой стороны, взяв в формуле \eqref{4.131} $u:=u'$ и подставив в нее соотношения
\eqref{4.147}, \eqref{4.148}, получим еще одно равенство
$$
(u_{0}',L^{+}\Phi_{\sigma}h)_{\Omega}-(f,(\Phi_{\sigma}h)_{0})_{\Omega}=
\sum_{j=1}^{q}\;(B_{j}u',h_{j})_{\Gamma}.
$$
Из этих равенств следует, что $B_{j}u'=g_{j}$ для всех $j=1,\ldots,q$. Это вместе с
\eqref{4.148} означает равенство $(L,B)\,u'=(f,g)$. Следовательно, в силу теоремы
\ref{th4.16} верно равенство $Q^{+}(f,g)=(f,g)$. Кроме того, из формул
\eqref{4.147}, \eqref{4.135} и определения проектора $\Psi_{\sigma}$ вытекает на
основании той же теоремы \ref{th4.16}, что $Pu'=u'$. Поэтому, воспользовавшись
изоморфизмом \eqref{4.96}, можем записать равенство $u'=(L,B)^{-1}Q^{+}(f,g)$. Оно
вместе с \eqref{4.145} влечет требуемое равенство $u'=\omega$.

Теперь в силу \eqref{4.146} справедливо $I_{L}\,\omega=I_{L}\,u'=(u_{0}',f)$. Отсюда
и из формул \eqref{4.144}, \eqref{4.145} следует, что
$$
\Pi_{\sigma}(u_{0},f)=(\Psi_{\sigma}u_{0},0)+I_{L}\,\omega=
(\Psi_{\sigma}u_{0},0)+(u_{0}',f)=(u_{0},f)
$$
для произвольного вектора $(u_{0},f)\in K_{\sigma,L}(\Omega)$. Таким образом,
искомый проектор $\Pi_{\sigma}$ построен.

Лемма \ref{lem4.6} доказана.

\medskip

Опираясь на эту лемму, докажем теорему \ref{th4.24}.


\textbf{Доказательство теоремы \ref{th4.24}.} Пусть $s\in\mathbb{R}$,
$\varphi\in\mathcal{M}$ и число $\varepsilon>0$. Воспользуемся частным случаем
$\varphi\equiv1$ этой теоремы \ref{th4.24}, который доказан Я.~А.~Ройтбергом
[\ref{Roitberg96}, с. 174] (теорема 6.1.1). Отображение \eqref{4.116} продолжается
по непрерывности до гомеоморфизмов
$$
I_{L}:\,H^{s\mp\varepsilon,(2q)}(\Omega)\leftrightarrow
K_{s\mp\varepsilon,L}(\Omega).
$$
Применим интерполяцию с функциональным параметром $\psi$ из теоремы \ref{th2.14},
где $\varepsilon=\delta$. Получим еще один гомеоморфизм
\begin{gather}\notag
I_{L}:\,[H^{s-\varepsilon,(2q)}(\Omega),
H^{s+\varepsilon,(2q)}(\Omega)]_{\psi}\leftrightarrow\\
[K_{s-\varepsilon,L}(\Omega),K_{s+\varepsilon,L}(\Omega)]_{\psi}.\label{4.149}
\end{gather}
Опишем пространства, в которых он действует. Согласно теореме \ref{th4.23}
\begin{equation}\label{4.150}
[H^{s-\varepsilon,(2q)}(\Omega),
H^{s+\varepsilon,(2q)}(\Omega)\bigr]_{\psi}=H^{s,\varphi,(2q)}(\Omega).
\end{equation}
Далее, в силу леммы \ref{lem4.6} и теоремы \ref{th2.6} можно проинтерполировать пару
подпространств $K_{s\mp\varepsilon,L}(\Omega)$ следующим образом:
\begin{gather*}
[K_{s-\varepsilon,L}(\Omega),K_{s+\varepsilon,L}(\Omega)\bigr]_{\psi}=\\
[H^{s-\varepsilon,(0)}(\Omega)\oplus H^{s-\varepsilon-2q,(0)}(\Omega),
H^{s+\varepsilon,(0)}(\Omega)\oplus H^{s+\varepsilon-2q,(0)}(\Omega)\bigr]_{\psi}\cap\\
K_{s-\varepsilon,L}(\Omega)=\\
H^{s,\varphi,(0)}(\Omega)\oplus H^{s-2q,\varphi,(0)}(\Omega)\cap
K_{s-\varepsilon,L}(\Omega)=K_{s,\varphi,L}(\Omega).
\end{gather*}
Здесь мы использовали также теоремы \ref{th2.5}, \ref{th3.10} и определение множеств
$K_{s-\varepsilon,L}(\Omega)$ и $K_{s,\varphi,L}(\Omega)$. В~силу теоремы
\ref{th2.6}, $K_{s,\varphi,L}(\Omega)$ является подпространством в
$H^{s,\varphi,(0)}(\Omega)\oplus H^{s-2q,\varphi,(0)}(\Omega)$, а равенства
\eqref{4.150} и
\begin{equation}\label{4.151}
[K_{s-\varepsilon,L}(\Omega),K_{s+\varepsilon,L}(\Omega)\bigr]_{\psi}=
K_{s,\varphi,L}(\Omega).
\end{equation}
выполняются с точностью до эквивалентности норм. Подставив эти равенства в
\eqref{4.149}, получим гомеоморфизм \eqref{4.117}. Утверждение (ii) теоремы
\ref{th4.24} доказано, а с ним~--- и пункт~(i).

Теорема \ref{th4.24} доказана.

\medskip

Для доказательства теоремы \ref{th4.25} нам понадобится еще одна лемма о проекторах.
Пусть $\sigma\in\mathbb{R}$ и $\varphi\in\mathcal{M}$. Положим
\begin{gather*}
K_{\sigma,\varphi,L}^{0}(\Omega):=\\ \{\,(u_{0},f)\in
H^{\sigma,\varphi,(0)}(\Omega)\oplus H^{\sigma-2q,\varphi,(0)}(\Omega):\,\mbox{верно
\eqref{4.119}}\,\}.
\end{gather*}
В силу теоремы \ref{th3.9} (iii), $K_{\sigma,\varphi,L}^{0}(\Omega)$ является
подпространством в $H^{\sigma,\varphi,(0)}(\Omega)\oplus
H^{\sigma-2q,\varphi,(0)}(\Omega)$. При этом $K_{\sigma,\varphi,L}(\Omega)\subseteq
K_{\sigma,\varphi,L}^{0}(\Omega)$.


\begin{lemma}\label{lem4.7}
Пусть числа $r\in\mathbb{N}$ и $\sigma\in\mathbb{R}$ удовлетворяют условию
\begin{equation}\label{4.152}
-r-1/2\leq\sigma<-r+1/2.
\end{equation}
Тогда существует проектор $\Pi^{(r)}_{\sigma}$ пространства
$H^{\sigma,(0)}(\Omega)\oplus H^{\sigma-2q,(0)}(\Omega)$ на подпространство
$K_{\sigma,L}^{0}(\Omega)$ такой, что $\Pi^{(r)}_{\sigma}$ --- продолжение
отображения $\Pi^{(r)}_{\lambda}$ при $\sigma<\lambda<-r+1/2$.
\end{lemma}


\textbf{Доказательство.} Предварительно установим одно полезное равенство. Пусть
$(u_{0},f)\in K_{\sigma,L}^{0}(\Omega)$. Воспользуемся тем, что в соболевском случае
$\varphi\equiv\nobreak1$ теорема \ref{th4.25} доказана Я.~А.~Ройтбергом
[\ref{Roitberg96}, с. 187] (теорема 6.2.1)  для любого параметра $s<-1/2$.
Следовательно, существует единственная пара $(u_{0}^{\ast},f)\in
K_{\sigma,L}(\Omega)$ такая, что $u_{0}=u_{0}^{\ast}$ в области $\Omega$. Поскольку
$u_{0},u_{0}^{\ast}\in
H^{\sigma,(0)}(\Omega)=H^{\sigma}_{\overline{\Omega}}(\mathbb{R}^{n})$ и
$\mathrm{supp}(u_{0}-u_{0}^{\ast})\subseteq\Gamma$, из неравенства \eqref{4.152}
вытекает следующее представление распределения $u_{0}-u_{0}^{\ast}$ (см., например,
[\ref{Roitberg96}, с. 185], лемма 6.2.2). Существует единственный вектор
\begin{equation}\label{4.153}
\omega^{\ast}=(\omega^{\ast}_{1},\ldots,\omega^{\ast}_{r})\in
\bigoplus_{j=1}^{r}\,H^{\sigma+j-1/2}(\Gamma)
\end{equation}
такой, что
\begin{equation}\label{4.154}
(u_{0}-u_{0}^{\ast},w)_{\Omega}=
\sum_{j=1}^{r}\,(\omega^{\ast}_{j},D_{\nu}^{j-1}w)_{\Gamma}.
\end{equation}
для любого $w\in H^{-\sigma,(0)}(\Omega)=H^{-\sigma}(\Omega)$.

Положим в формуле Грина \eqref{4.128}
$$
u:=u^{\ast}:=I_{L}^{-1}(u_{0}^{\ast},f)\in H^{\sigma,(2q)}(\Omega)
$$
и воспользуемся тем, что
$$
Lu=Lu^{\ast}=f\in
H^{\sigma-2q,(0)}(\Omega)=H^{\sigma-2q}_{\overline{\Omega}}(\mathbb{R}^{n}).
$$
Получим для любого $v\in H^{2q-\sigma,(2q)}(\Omega)=H^{2q-\sigma}(\Omega)$ (см.
теорему \ref{th4.13} (iii)) равенство
\begin{equation}\label{4.155}
(f,v)_{\Omega}+\sum_{j=1}^{q}\;(B_{j}u^{\ast},\,C_{j}^{+}v)_{\Gamma}
=(u_{0}^{\ast},L^{+}v)_{\Omega}+
\sum_{j=1}^{q}\;(C_{j}u^{\ast},\,B_{j}^{+}v)_{\Gamma}.
\end{equation}
Здесь $u_{0}^{\ast}$ и $v=v_{0}$~--- начальные компоненты векторов $T_{2q}u^{\ast}$
и $T_{2q}v$ соответственно.

Теперь из формул \eqref{4.155} и \eqref{4.154}, где берем $w:=L^{+}v\in
H^{-\sigma}(\Omega)$, следует равенство
\begin{gather}\notag
(f,v)_{\Omega}+\sum_{j=1}^{q}\;(B_{j}u^{\ast},\,C_{j}^{+}v)_{\Gamma}=\\
(u_{0},L^{+}v)_{\Omega}-
\sum_{j=1}^{r}\,(\omega^{\ast}_{j},D_{\nu}^{j-1}L^{+}v)_{\Gamma}+
\sum_{j=1}^{q}\;(C_{j}u^{\ast},\,B_{j}^{+}v)_{\Gamma}.\label{4.156}
\end{gather}
Граничные операторы $B_{j}^{+}$, $C_{j}^{+}$, $j=1,\ldots,q$, и
$(D_{\nu}^{j-1}L^{+}\cdot)\!\upharpoonright\!\Gamma$, $j=1,\ldots,r$, образуют
систему Дирихле порядка $2q+r$. Для нее в силу упомянутой выше леммы 6.1.2 из
монографии [\ref{Roitberg96}, с.~176] справедливо следующее утверждение.
Ограниченный линейный оператор
\begin{gather*}
(B^{+}_{1},\ldots,B^{+}_{q},C^{+}_{1},\ldots,C^{+}_{q},
(L^{+}\cdot)\upharpoonright\Gamma,\ldots,
(D_{\nu}^{r-1}L^{+}\cdot)\upharpoonright\Gamma):\\
H^{2q-\sigma,(2q+r)}(\Omega)\rightarrow \\
\bigoplus_{j=1}^{q} \,H^{2q-\sigma-m^{+}_{j}-1/2}(\Gamma)\oplus
\bigoplus_{j=1}^{q}\,H^{-\sigma+m_{j}+1/2}(\Gamma)\oplus
\bigoplus_{j=1}^{r}\,H^{-\sigma-j+1/2}(\Gamma)
\end{gather*}
имеет ограниченный правый обратный оператор
\begin{gather}\notag
\Phi^{(r)}_{\sigma}: \bigoplus_{j=1}^{q}\,H^{2q-\sigma-m^{+}_{j}-1/2}(\Gamma)\oplus
\notag \bigoplus_{j=1}^{q}\,H^{-\sigma+m_{j}+1/2}(\Gamma)\oplus\\
\bigoplus_{j=1}^{r}\,H^{-\sigma-j+1/2}(\Gamma)\rightarrow
H^{2q-\sigma,(2q+r)}(\Omega)=H^{2q-\sigma}(\Omega) \label{4.157}
\end{gather}
такой, что $\Phi^{(r)}_{\sigma}$ --- сужение оператора $\Phi^{(r)}_{\lambda}$ при
$\sigma<\lambda<-r+1/2$. Здесь $\mathrm{ord}\,C_{j}^{+}=2q-1-m_{j}$ и
$\mathrm{ord}\,B_{j}^{+}=:m_{j}^{+}$, а равенство пространств, записанное в
\eqref{4.157}, следует из теоремы \ref{th4.13} (iii) ввиду условия $\sigma<-r+1/2$.

Для произвольного вектора
\begin{equation}\label{4.158}
h:=(0,\ldots,0,h_{1},\ldots,h_{r})
\in\{0\}^{2q}\oplus\bigoplus_{j=1}^{r}\,H^{-\sigma-j+1/2}(\Gamma)
\end{equation}
положим $v:=\Phi^{(r)}_{\sigma}h\in H^{2q-\sigma}(\Omega)$ в формуле \eqref{4.156}.
Поскольку
\begin{gather*}
B^{+}_{j}v=0,\quad C^{+}_{j}v=0\quad\mbox{на}\quad\Gamma\quad\mbox{при}\quad
j=1,\ldots,q,\\
D_{\nu}^{j-1}L^{+}v=h_{j}\quad\mbox{на}\quad\Gamma\quad\mbox{при}\quad j=1,\ldots,r,
\end{gather*}
то получим равенство
\begin{equation}\label{4.159}
(u_{0},L^{+}\Phi^{(r)}_{\sigma}h)_{\Omega}-(f,\Phi^{(r)}_{\sigma}h)_{\Omega}=
\sum_{j=1}^{r}\,(\omega^{\ast}_{j},h_{j})_{\Gamma}.
\end{equation}
Напомним, что здесь пара $(u_{0},f)\in K_{\sigma,L}^{0}(\Omega)$ и вектор $h$ вида
\eqref{4.158} произвольные, а вектор $(\omega^{\ast}_{1},\ldots,\omega^{\ast}_{r})$
удовлетворяет условиям \eqref{4.153}, \eqref{4.154} и однозначно ими определяется по
паре $(u_{0},f)$.

Теперь перейдем к построению проектора $\Pi^{(r)}_{\sigma}$. Выберем произвольный
вектор $(u_{0},f)$, удовлетворяющий включению \eqref{4.132}. Как и прежде, обозначим
через $\|\cdot\|_{\sigma}$ норму в пространстве \eqref{4.132}. Рассмотрим функционал
\begin{equation}\label{4.160}
l^{(r)}_{\sigma}(h):=(u_{0},L^{+}\Phi^{(r)}_{\sigma}h)_{\Omega}-
(f,\Phi^{(r)}_{\sigma}h)_{\Omega},\quad h\mbox{ --- вектор \eqref{4.158}}.
\end{equation}
Он ограничен в силу следующей цепочки неравенств:
\begin{gather*}
|l^{(r)}_{\sigma}(h)|\leq|(u_{0},L^{+}\Phi^{(r)}_{\sigma}h)_{\Omega}|+
|(f,\Phi^{(r)}_{\sigma}h)_{\Omega}|\leq \\
\|u_{0}\|_{H^{\sigma,(0)}(\Omega)}\,
\|L^{+}\Phi^{(r)}_{\sigma}h\|_{H^{-\sigma,(0)}(\Omega)}+\\
\|f\|_{H^{\sigma-2q,(0)}(\Omega)}\,
\|\Phi^{(r)}_{\sigma}h\|_{H^{-\sigma+2q,(0)}(\Omega)}\leq \\
\bigl(c_{1}\,\|u_{0}\|_{H^{\sigma,(0)}(\Omega)}+
\|f\|_{H^{\sigma-2q,(0)}(\Omega)}\bigr)\,
\|\Phi^{(r)}_{\sigma}h\|_{H^{-\sigma+2q}(\Omega)}\leq \\
\leq (c_{1}+1)\,c_{2}\,\|(u_{0},f)\|_{\sigma}\,
\Bigl(\,\sum_{j=1}^{r}\;\|h_{j}\|_{H^{-\sigma-j+1/2}(\Gamma)}^{2}\Bigr)^{1/2}.
\end{gather*}
Здесь $c_{1}$ --- норма оператора
$$
L^{+}:H^{-\sigma+2q}(\Omega)\rightarrow H^{-\sigma}(\Omega)=H^{-\sigma,(0)}(\Omega),
$$
а $c_{2}$ --- норма оператора \eqref{4.157}. Итак, $l^{(r)}_{\sigma}$~---
антилинейный ограниченный функционал на пространстве
$\bigoplus_{j=1}^{r}\,H^{-\sigma-j+1/2}(\Gamma)$, причем его норма удовлетворяет
неравенству $\|l^{(r)}_{\sigma}\|\leq c_{3}\,\|(u_{0},f)\|_{\sigma}$, где
$c_{3}:=(c_{1}+1)\,c_{2}$. В силу теоремы \ref{th2.22} (v)
\begin{equation}\label{4.161}
\Bigl(\,\exists!\;(\omega_{1},\ldots,\omega_{r})\in
\bigoplus_{j=1}^{r}\,H^{\sigma+j-1/2}(\Gamma)\,\Bigr):\,
l^{(r)}_{\sigma}(h)=\sum_{j=1}^{r}\,(\omega_{j},h_{j})_{\Gamma}.
\end{equation}
При этом
\begin{equation}\label{4.162}
\Bigl(\,\sum_{j=1}^{r}\;\|\omega_{j}\|_{H^{\sigma+j-1/2}(\Gamma)}^{2}
\Bigr)^{1/2}\asymp\|l^{(r)}_{\sigma}\|\leq c_{3}\,\|(u_{0},f)\|_{\sigma}.
\end{equation}

По вектору $(\omega_{1},\ldots,\omega_{r})$ образуем распределение $\omega'\in
H^{\sigma,(0)}(\Omega)$ по формуле
\begin{equation}\label{4.163}
(\omega',w)_{\Omega}:=\sum_{j=1}^{r}\;(\omega_{j},D_{\nu}^{j-1}w)_{\Gamma},
\end{equation}
где $w\in H^{-\sigma,(0)}(\Omega)=H^{-\sigma}(\Omega)$ произвольно. Это определение
корректно в силу теорем \ref{th2.22} (v), \ref{th3.5} и \ref{th3.8} (iii).
Действительно, поскольку $-\sigma>r-1/2$ ввиду условия \eqref{4.152}, мы имеем
ограниченные операторы
$$
D_{\nu}^{j-1}:H^{-\sigma}(\Omega)\rightarrow
H^{-\sigma-j+1/2}(\Gamma)\quad\mbox{при}\quad j=1,\ldots,r.
$$
Следовательно, $(\omega',\cdot)_{\Omega}$~--- антилинейный ограниченный функционал
на пространстве $H^{-\sigma}(\Omega)$, что влечет за собой включение
$$
\omega'\in H^{\sigma,(0)}(\Omega)=H^{\sigma}_{\overline{\Omega}}(\mathbb{R}^{n}).
$$
При этом верно включение $\mathrm{supp}\,\omega'\subseteq\Gamma$ и двусторонняя
оценка
\begin{equation}\label{4.164}
\|\omega'\|_{H^{\sigma,(0)}(\Omega)}\asymp
\Bigl(\,\sum_{j=1}^{r}\;\|\omega_{j}\|_{H^{\sigma+j-1/2}(\Gamma)}^{2} \Bigr)^{1/2}
\end{equation}
(см., например, лемму~6.2.2 из монографии [\ref{Roitberg96}, с.~185]).

На векторах \eqref{4.132} определим линейное отображение $\Upsilon^{(r)}_{\sigma}$
по формуле
\begin{equation}\label{4.165}
\Upsilon^{(r)}_{\sigma}:\,(u_{0},f)\mapsto l^{(r)}_{\sigma}\mapsto
(\omega_{1},\ldots,\omega_{r})\mapsto\omega'.
\end{equation}
В силу оценок \eqref{4.162} и \eqref{4.164} оно является ограниченным оператором
\begin{equation}\label{4.166}
\Upsilon^{(r)}_{\sigma}:\,H^{\sigma,(0)}(\Omega)\oplus
H^{\sigma-2q,(0)}(\Omega)\rightarrow H^{\sigma,(0)}(\Omega).
\end{equation}
Покажем, что $\Upsilon^{(r)}_{\sigma}$ --- продолжение отображения
$\Upsilon^{(r)}_{\lambda}$ при $\sigma<\lambda<-r+1/2$. Если
$$
(u_{0},f)\in H^{\lambda,(0)}(\Omega)\oplus H^{\lambda-2q,(0)}(\Omega),
$$
то наряду с $l^{(r)}_{\sigma}$ определен функционал $l^{(r)}_{\lambda}$. При этом
$l^{(r)}_{\sigma}$ --- сужение функционала $l^{(r)}_{\lambda}$, поскольку
$\Phi^{(r)}_{\sigma}$ --- сужение оператора $\Phi^{(r)}_{\lambda}$. В частности,
$l^{(r)}_{\sigma}(h)=l^{(r)}_{\lambda}(h)$ для любого вектора
$(h_{1},\ldots,h_{q})\in (C^{\infty}(\Gamma))^{q}$. Отсюда в силу \eqref{4.161},
\eqref{4.163} получаем
$\Upsilon^{(r)}_{\sigma}(u_{0},f)=\Upsilon^{(r)}_{\lambda}(u_{0},f)$, т.~е.
$\Upsilon^{(r)}_{\sigma}$~--- продолжение отображения $\Upsilon^{(r)}_{\lambda}$.

Линейное отображение $\Pi^{(r)}_{\sigma}$ определяем с помощью оператора
$\Upsilon^{(r)}_{\sigma}$ и проектора $\Pi_{\sigma}$ из леммы \ref{lem4.6} следующим
образом:
\begin{equation}\label{4.167}
\Pi^{(r)}_{\sigma}(u_{0},f):=(\omega',0)+\Pi_{\sigma}(u_{0}-\omega',f),\;\;
\omega':=\Upsilon^{(r)}_{\sigma}(u_{0},f).
\end{equation}
Здесь $(u_{0},f)$ --- произвольный вектор \eqref{4.132}. Как указывалось выше, верны
включения $\omega'\in H^{\sigma,(0)}(\Omega)$ и
$\mathrm{supp}\,\omega'\subseteq\Gamma$. Следовательно, $(\omega',0)\in
K_{\sigma,L}^{0}(\Omega)$. Этот факт вместе с ограниченностью оператора
\eqref{4.166}, леммой \ref{lem4.6} и включением $K_{\sigma,L}(\Omega)\subseteq
K_{\sigma,L}^{0}(\Omega)$ влечет за собой ограниченность оператора
$$
\Pi^{(r)}_{\sigma}:\,H^{\sigma,(0)}(\Omega)\oplus
H^{\sigma-2q,(0)}(\Omega)\rightarrow K_{\sigma,L}^{0}(\Omega).
$$
При этом, так как с уменьшением параметра $\sigma$ операторы
$\Upsilon^{(r)}_{\sigma}$ и $\Pi_{\sigma}$ расширяются, то $\Pi^{(r)}_{\sigma}$
является продолжением оператора $\Pi^{(r)}_{\lambda}$ при $\sigma<\lambda<-r+1/2$.

Остается показать, что $\Pi^{(r)}_{\sigma}$ --- проектор на подпространство
$K_{\sigma,L}^{0}(\Omega)$, т.~е. $\Pi^{(r)}_{\sigma}(u_{0},f)=(u_{0},f)$ для любого
вектора $(u_{0},f)\in K_{\sigma,L}^{0}(\Omega)$. Выберем произвольно такой вектор. В
силу теоремы \ref{th4.25} в случае $\varphi\equiv1$, для $(u_{0},f)$ существует
единственная пара $(u_{0}^{\ast},f)\in K_{\sigma,L}(\Omega)$ такая, что выполняется
условие \eqref{4.121}. Покажем, что $u_{0}-u_{0}^{\ast}=\omega'$, где
$\omega':=\Upsilon^{(r)}_{\sigma}(u_{0},f)$. Для произвольного распределения $w\in
H^{-\sigma}(\Omega)$ определим вектор $h$ по формуле \eqref{4.158}, в которой
полагаем
\begin{equation}\label{4.168}
h_{j}:=(D_{\nu}^{j-1}w)\!\upharpoonright\!\Gamma\in
H^{-\sigma-j+1/2}(\Gamma)\quad\mbox{при}\quad j=1,\ldots,r.
\end{equation}
Как указывалось выше, следы $h_{j}$ существуют ввиду условия \eqref{4.152}. Напомним
также, что условие \eqref{4.121} означает включение
$\mathrm{supp}(u_{0}-u_{0}^{\ast})\subseteq\Gamma$, откуда следует равенство
\eqref{4.154} для некоторого (единственного) вектора \eqref{4.153}. Теперь из формул
\eqref{4.154}, \eqref{4.159} и \eqref{4.168} вытекает равенство
$$
(u_{0}-u_{0}^{\ast},w)_{\Omega}=
(u_{0},L^{+}\Phi^{(r)}_{\sigma}h)_{\Omega}-(f,\Phi^{(r)}_{\sigma}h)_{\Omega}.
$$
С другой стороны, в силу \eqref{4.160}, \eqref{4.161}, \eqref{4.163} и \eqref{4.168}
имеем
$$
(\omega',w)_{\Omega}=l^{(r)}_{\sigma}(h)=
(u_{0},L^{+}\Phi^{(r)}_{\sigma}h)_{\Omega}-(f,\Phi^{(r)}_{\sigma}h)_{\Omega}.
$$
Следовательно,
$$
(u_{0}-u_{0}^{\ast},w)_{\Omega}=(\omega',w)_{\Omega}\quad\mbox{для любого}\quad w\in
H^{-\sigma}(\Omega).
$$
В силу теоремы \ref{th3.8} (iii) это равносильно равенству
$$
u_{0}-u_{0}^{\ast}=\omega':=\Upsilon^{(r)}_{\sigma}(u_{0},f)\quad\mbox{в}\quad
H^{\sigma,(0)}(\Omega).
$$
Отсюда и из \eqref{4.167} получаем
\begin{gather*}
\Pi^{(r)}_{\sigma}(u_{0},f)=(\omega',0)+\Pi_{\sigma}(u_{0}-\omega',f)=\\
(\omega',0)+\Pi_{\sigma}(u_{0}^{\ast},f)=(\omega',0)+(u_{0}^{\ast},f)=(u_{0},f)
\end{gather*}
для произвольного вектора $(u_{0},f)\in K_{\sigma,L}^{0}(\Omega)$. Здесь мы также
воспользовались тем, что $(u_{0}^{\ast},f)\in K_{\sigma,L}(\Omega)$ и
$\Pi_{\sigma}$~--- проектор на подпространство $K_{\sigma,L}(\Omega)$. Таким
образом, искомый проектор $\Pi^{(r)}_{\sigma}$ построен.

Лемма \ref{lem4.7} доказана.

\medskip

Воспользовавшись этой леммой, докажем теорему \ref{th4.25}.

\medskip

\textbf{Доказательство теоремы \ref{th4.25}.} Пусть $s<-1/2$,
$s+1/2\notin\mathbb{Z}$, и $\varphi\in\mathcal{M}$. Выберем числа $r\in\mathbb{N}$ и
$\varepsilon\in(0,1/2)$ такие, что
\begin{equation}\label{4.169}
-r-1/2<s\mp\varepsilon<-r+1/2.
\end{equation}
Как уже отмечалось, в соболевском случае $\varphi\equiv1$ теорема \ref{th4.25}
доказана Я.~А.~Ройтбергом [\ref{Roitberg96}, с.~187]. Поэтому мы можем ввести
линейное отображение $G:(u_{0},f)\mapsto(u_{0}^{\ast},f)$, где $(u_{0},f)\in
K_{s-\varepsilon,L}^{0}(\Omega)$, а $(u_{0}^{\ast},f)\in
K_{s-\varepsilon,L}(\Omega)$ удовлетворяет \eqref{4.121}. Оно определяет
ограниченные операторы
\begin{equation}\label{4.170}
G:\,K_{s\mp\varepsilon,L}^{0}(\Omega)\rightarrow K_{s\mp\varepsilon,L}(\Omega).
\end{equation}

Мы докажем теорему  \ref{th4.25} для любого $\varphi\in\mathcal{M}$, если установим,
что $G$ является ограниченным оператором из $K_{s,\varphi,L}^{0}(\Omega)$ в
$K_{s,\varphi,L}(\Omega)$. Сделаем это с помощью интерполяции с параметром $\psi$ из
теоремы \ref{th2.14}, где берем $\varepsilon=\delta$. Применяя ее к \eqref{4.170},
получаем еще один ограниченный оператор
\begin{equation}\label{4.171}
G:\,[K_{s-\varepsilon,L}^{0}(\Omega),K_{s+\varepsilon,L}^{0}(\Omega)]_{\psi}\rightarrow
[K_{s-\varepsilon,L}(\Omega),K_{s+\varepsilon,L}(\Omega)]_{\psi}.
\end{equation}
Опишем пространства, в которых он действует. На основании леммы \ref{lem4.7},
неравенства \eqref{4.169} и теоремы \ref{th2.6} мы можем проинтерполировать пару
подпространств $K_{s\mp\varepsilon,L}^{0}(\Omega)$ следующим образом:
\begin{gather*}
[K_{s-\varepsilon,L}^{0}(\Omega),K_{s+\varepsilon,L}^{0}(\Omega)]_{\psi}=\\
[H^{s-\varepsilon,(0)}(\Omega)\oplus H^{s-\varepsilon-2q,(0)}(\Omega),
H^{s+\varepsilon,(0)}(\Omega)\oplus
H^{s+\varepsilon-2q,(0)}(\Omega)]_{\psi}\cap\\ K_{s-\varepsilon,L}^{0}(\Omega)=\\
H^{s,\varphi,(0)}(\Omega)\oplus H^{s-2q,\varphi,(0)}(\Omega)\cap
K_{s-\varepsilon,L}^{0}(\Omega)=K_{s,\varphi,L}^{0}(\Omega).
\end{gather*}
Здесь мы также воспользовались теоремами \ref{th2.5}, \ref{th3.10} и определением
пространств $K_{s-\varepsilon,L}^{0}(\Omega)$, $K_{s,\varphi,L}^{0}(\Omega)$.
Полученное равенство пространств
\begin{equation}\label{4.172}
[K_{s-\varepsilon,L}^{0}(\Omega),
K_{s+\varepsilon,L}^{0}(\Omega)]_{\psi}=K_{s,\varphi,L}^{0}(\Omega)
\end{equation}
выполняется с точностью до эквивалентности норм. Подставив равенства \eqref{4.172} и
\eqref{4.151} в формулу \eqref{4.170}, получаем ограниченный оператор
$$
G:\,K_{s,\varphi,L}^{0}(\Omega)\rightarrow K_{s,\varphi,L}(\Omega),
$$
что и надо было доказать.

Теорема \ref{th4.25} доказана.


\begin{remark}\label{rem4.9}
Теоремы \ref{th4.24}, \ref{th4.25} и леммы \ref{lem4.6}, \ref{lem4.7} доказаны в
предположении, что задача \eqref{4.1} регулярная эллиптическая. На самом деле, они
справедливы для произвольного дифференциального выражения $L$, правильно
эллиптического в области $\overline{\Omega}$. Это следует из того, что для него
существует регулярная эллиптическая краевая задача, например, задача Дирихле.
\end{remark}




\markright{\emph \ref{sec4.4}. Расширение теорем Лионса--Мадженеса}

\section[Расширение теорем Лионса--Мадженеса]
{Расширение теорем \\ Лионса--Мадженеса}\label{sec4.4}

\markright{\emph \ref{sec4.4}. Расширение теорем Лионса--Мадженеса}

В этом пункте мы установим индивидуальные теоремы о разрешимости регулярной
эллиптической краевой задачи \eqref{4.1} в шкалах пространств Соболева. В~отличие от
общих теорем \ref{th4.1}, \ref{th4.15}, в индивидуальных теоремах область
определения оператора $(L,B)$ зависит от коэффициентов эллиптического выражения $L$.
В~работах Ж.-Л.~Лионса и Е.~Мадженеса [\ref{LionsMagenes71}, \ref{Magenes66},
\ref{LionsMagenes62}, \ref{LionsMagenes63}] получен ряд индивидуальных теорем для
оператора $(L,B)$, действующего в пространствах Соболева, содержащих нерегулярные
распределения. Мы докажем некоторую общую форму теорем Лионса-Мадженеса --- будет
найдено общее условие на пространство правых частей эллиптического уравнения, при
котором оператор $(L,B)$ ограничен и нетеров в соответствующей паре гильбертовых
пространств. Мы укажем также широкие классы пространств, удовлетворяющих этому
условию. Они содержат пространства, использованные Лионсом и Мадженесом, и много
других важных пространств, среди них и некоторые пространства Хермандера. Отметим,
что в рассматриваемых нами индивидуальных теоремах решение и правая часть
эллиптического уравнения являются распределениями в области $\Omega$, в отличие от
общей теоремы~\ref{th4.15}.

Результаты этого пункта мотивируют формулировки и методы доказательства
индивидуальных теорем о разрешимости эллиптической краевой задачи в шкалах
пространств Хермандера. Эти теоремы будут установлены в п.~\ref{sec4.5}.

\subsection{Теоремы Лионса--Мадженеса}\label{sec4.4.1}

Предварительно укажем на одно важное различие в определении пространств Соболева
отрицательного порядка на евклидовой области у нас и в цитированных работах
Ж.-Л.~Лионса, Е.~Мадженеса. Напомним, что мы, следуя [\ref{Triebel80},
\ref{Triebel86}], использовали такое определение пространства Соболева произвольного
порядка $s\in\mathbb{R}$ в области $\Omega$:
\begin{gather} \label{4.173}
H^{s}(\Omega):=\{u:=w\!\upharpoonright\!\Omega:\,w\in
H^{s}(\mathbb{R}^{n})\},\\
\|u\|_{H^{s}(\Omega)}:= \inf\bigl\{\|w\|_{H^{s}(\mathbb{R}^{n})}:w\in
H^{s}(\mathbb{R}^{n}),\;w=u\;\mbox{на}\;\Omega\bigr\}. \label{4.174}
\end{gather}
(см. определение \ref{def3.2} в случае $\varphi\equiv1$). Лионс и Мадженес
используют это определение только при $s\geq0$, а для $s<0$ в качестве пространства
Соболева порядка $s$ в области $\Omega$ берут двойственное пространство
$(H^{-s}_{0}(\Omega))'$. Здесь, напомним, $H^{-s}_{0}(\Omega)$ --- замыкание
множества $C^{\infty}_{0}(\Omega)$ в пространстве $H^{-s}(\Omega)$, а двойственность
пространств рассматривается относительно скалярного произведения в $L_{2}(\Omega)$.
Известно [\ref{Triebel80}, с.~414] (теорема 4.8.2), что эти определения дают одно и
тоже пространство (с точностью до эквивалентности норм), если порядок $s<0$
неполуцелый. Для полуцелых $s<0$ получаются различные пространства.

В настоящем п.~\ref{sec4.4} мы будем использовать определение пространств Соболева в
области $\Omega$, данное Лионсом и Мадженесом:
\begin{equation}\label{4.175}
\mathrm{H}^{s}(\Omega):=
\begin{cases}
\;H^{s}(\Omega)&\;\text{при}\;s\geq0, \\
\;(H^{-s}_{0}(\Omega))'&\;\text{при}\;s<0.
\end{cases}
\end{equation}
Слева в обозначении используется прямое написание буквы $\mathrm{H}$, а не курсивное
как в определении \eqref{4.173}. Здесь двойственное пространство
$(H^{-s}_{0}(\Omega))'$ состоит из антилинейных функционалов.

Функционалы из пространства $\mathrm{H}^{s}(\Omega)$, где $s<0$, однозначно
определяются по своим значениям на пробных функциях из пространства
$C^{\infty}_{0}(\Omega)$. Поэтому корректно отождествлять эти функционалы с
распределениями в области $\Omega$. При этом для любого $s<0$
\begin{gather}\label{4.176}
\mathrm{H}^{s}(\Omega)=\bigl\{w\!\upharpoonright\!\Omega:\,w\in
H^{s}_{\overline{\Omega}}(\mathbb{R}^{n})\bigr\}, \\
\|u\|_{\mathrm{H}^{s}(\Omega)}=\inf\,\bigl\{\,\|w\|_{H^{s}(\mathbb{R}^{n})}:\, w\in
H^{s}_{\overline{\Omega}}(\mathbb{R}^{n}),\;w=u\;\mbox{в}\;\Omega\,\bigr\}\label{4.177}
\end{gather}
(см. [\ref{LionsMagenes71}], гл.~1, замечание 12.5). Отсюда и из формул
\eqref{4.173}, \eqref{4.174} следует непрерывное вложение
$\mathrm{H}^{s}(\Omega)\hookrightarrow H^{s}(\Omega)$ и плотность множества
$C^{\infty}_{0}(\Omega)$ в $\mathrm{H}^{s}(\Omega)$ при каждом $s<0$.

Как отмечалось выше,
\begin{equation}\label{4.178}
\mathrm{H}^{s}(\Omega)=H^{s}(\Omega)\;\Leftrightarrow\;
s\in\mathbb{R}\setminus\{-1/2,-3/2,-5/2,\ldots\}.
\end{equation}
Если параметр $s<0$ полуцелый, то пространство $\mathrm{H}^{s}(\Omega)$ \'уже
пространства $H^{s}(\Omega)$.

Отметим, что для произвольных $s\in\mathbb{R}$ и $\varepsilon>0$ выполняется
компактное и плотное вложение
$\mathrm{H}^{s+\varepsilon}(\Omega)\hookrightarrow\mathrm{H}^{s}(\Omega)$.

Теперь перейдем к теоремам Лионса-Мадженеса. Для этого сначала напомним классическую
общую теорему о разрешимости краевой задачи \eqref{4.1} в шкале позитивных
пространств Соболева, сформулированную выше в предложении~\ref{prop3.1}.

\medskip

\textbf{Теорема A.} \it Отображение
\begin{equation}\label{4.179}
u\mapsto(Lu,Bu),\quad u\in C^{\infty}(\,\overline{\Omega}\,),
\end{equation}
продолжается по непрерывности до ограниченного нетерового оператора
\begin{gather}\label{4.180}
(L,B):\,\mathrm{H}^{\sigma+2q}(\Omega)\rightarrow\mathrm{H}^{\sigma}(\Omega)\oplus
\bigoplus_{j=1}^{q}\,H^{\sigma+2q-m_{j}-1/2}(\Gamma)=:\\ \notag
\mathbf{H}_{\sigma}(\Omega,\Gamma)
\end{gather}
при любом вещественном $\sigma\geq0$. Ядро этого оператора совпадает с $N$, а
область значений состоит из всех векторов
$(f,g_{1},\ldots,g_{q})\in\mathbf{H}_{\sigma}(\Omega,\Gamma)$, удовлетворяющих
условию \eqref{4.4}. Индекс оператора \eqref{4.180} равен $\dim N-\dim N^{+}$ и не
зависит от~$\sigma$. \rm

\medskip

Здесь и далее, в отличие от предыдущих теорем о разрешимости краевой задачи
\eqref{4.1}, индекс $s$, задающий гладкость области определения оператора $(L,B)$,
удобно представлять в виде суммы $s=\sigma+2q$. Отметим, что теорема~A уже
распространена на уточненную шкалу пространств Хермандера в п.~\ref{sec4.1.1}
(теорема~\ref{th4.1}).

Как отмечалось выше в п.~\ref{sec4.2.1}, теорема~A не верна для произвольного
значения $\sigma<0$. Подход Ж.-Л.~Лионса и Э.~Мадженеса состоит в том, чтобы вместо
$\mathrm{H}^{\sigma+2q}(\Omega)$ использовать в качестве области определения
оператора $(L,B)$ более узкое пространство
\begin{equation}\label{4.181}
D^{\sigma+2q}_{L,X}(\Omega):=\{u\in\mathrm{H}^{\sigma+2q}(\Omega):\,Lu\in
X^{\sigma}(\Omega)\},
\end{equation}
где $X^{\sigma}(\Omega)$ --- некоторое гильбертово пространство, непрерывно
вложенное в $\mathrm{H}^{\sigma}(\Omega)$. Здесь и далее образ $Lu$ для
$u\in\mathcal{D}'(\Omega)$ понимается в смысле теории распределений. В~пространстве
\eqref{4.181} вводится скалярное произведение графика
\begin{equation}\label{4.182}
(u_{1},u_{2})_{D^{\sigma+2q}_{L,X}(\Omega)}:=
(u_{1},u_{2})_{\mathrm{H}^{\sigma+2q}(\Omega)}+(Lu_{1},Lu_{2})_{X^{\sigma}(\Omega)}
\end{equation}
и соответствующая ему норма.

Пространство $D^{\sigma+2q}_{L,X}(\Omega)$ со скалярным произведением \eqref{4.182}
является полным. В самом деле, если последовательность $(u_{k})$ фундаментальная в
$D^{\sigma+2q}_{L,X}(\Omega)$, то, поскольку пространства
$\mathrm{H}^{\sigma+2q}(\Omega)$ и $X^{\sigma}(\Omega)$ полные, существуют пределы
$u:=\lim u_{k}$ в
$\mathrm{H}^{\sigma+2q}(\Omega)\hookrightarrow\mathcal{D}'(\Omega)$ и $f:=\lim
Lu_{k}$ в $X^{\sigma}(\Omega)\hookrightarrow\mathcal{D}'(\Omega)$ (вложения
непрерывны). Так как дифференциальный оператор $L$ непрерывен в
$\mathcal{D}'(\Omega)$, то в силу первого предела имеем: $Lu=\lim Lu_{k}$ в
$\mathcal{D}'(\Omega)$. Отсюда на основании второго предела получаем равенство
$Lu=f\in X^{\sigma}(\Omega)$. Следовательно, $u\in D^{\sigma+2q}_{L,X}(\Omega)$ и
$\lim u_{k}=u$ в пространстве $D^{\sigma+2q}_{L,X}(\Omega)$, т.~е. это пространство
полное.

В работах Ж.-Л.~Лионса и Э.~Мадженеса [\ref{LionsMagenes71}, \ref{Magenes66},
\ref{LionsMagenes62}, \ref{LionsMagenes63}] найдены некоторые важные пространства
$X^{\sigma}(\Omega)$, такие, что отображение \eqref{4.179} продолжается по
непрерывности до ограниченного нетерового оператора
\begin{gather}\label{4.183}
(L,B):\,D^{\sigma+2q}_{L,X}(\Omega)\rightarrow
X^{\sigma}(\Omega)\oplus\bigoplus_{j=1}^{q}\,H^{\sigma+2q-m_{j}-1/2}(\Gamma) =:\\
\mathbf{X}_{\sigma}(\Omega,\Gamma),\notag
\end{gather}
если $\sigma<0$. В отличие от теоремы~A, область определения оператора \eqref{4.183}
зависит вместе с топологией от коэффициентов эллиптического выражения~$L$. Поэтому
теоремы о свойствах оператора \eqref{4.183} являются индивидуальными теоремами о
разрешимости краевой задачи \eqref{4.1}.

Сформулируем две индивидуальные теоремы, принадлежащие Ж.-Л.~Лионсу и Э.~Мадженесу.

\medskip

\textbf{Теорема LM$_{\mathbf{1}}$} [\ref{LionsMagenes62}, \ref{LionsMagenes63}]. \it
Пусть $\sigma<0$ и $X^{\sigma}(\Omega):=L_{2}(\Omega)$. Тогда отображение
\eqref{4.179} продолжается по непрерывности до ограниченного нетерового оператора
\eqref{4.183}. Его ядро совпадает с $N$, а область значений состоит из всех векторов
$(f,g_{1},\ldots,g_{q})\in\mathbf{X}_{\sigma}(\Omega,\Gamma)$, удовлетворяющих
условию \eqref{4.4}. Индекс оператора \eqref{4.183} равен $\dim N-\dim N^{+}$ и не
зависит от~$\sigma$. \rm

\medskip

Здесь следует выделить случай $\sigma=-2q$, который важен в спектральной теории
эллиптических операторов [\ref{Grubb68}, \ref{Grubb71}, \ref{Grubb96},
\ref{Mikhailets82}, \ref{Mikhailets89}]. В этом случае пространство
\begin{equation}\label{4.184}
D^{0}_{L,L_{2}}(\Omega)=\{u\in L_{2}(\Omega):\,Lu\in L_{2}(\Omega)\}
\end{equation}
является областью определения максимального оператора, отвечающего дифференциальному
выражению $L$ в пространстве $L_{2}(\Omega)$. Заметим, что даже если все
коэффициенты выражения $L$ постоянны, пространство \eqref{4.184} существенно зависит
от каждого из них. Это видно из следующего результата Л.~Хермандера
[\ref{Hermander59}, с.~78].

Пусть $L$ и $M$~--- два линейных дифференциальных выражения с постоянными
коэффициентами. Тогда, если $D^{0}_{L,L_{2}}(\Omega)\subseteq
D^{0}_{M,L_{2}}(\Omega)$, то либо $M=\alpha L+\beta$ для некоторых
$\alpha,\beta\in\mathbb{C}$, либо $L$ и $M$~--- многочлены относительно оператора
дифференцирования вдоль некоторого вектора $e$, причем
$\mathrm{ord}\,M\leq\mathrm{ord}\,L$. Заметим, что для эллиптических операторов
вторая возможность исключается.

Для формулировки второй теоремы Лионса-Мадженеса нам понадобится весовое
пространство
\begin{equation}\label{4.185}
\varrho\mathrm{H}^{\sigma}(\Omega):=\{f=\varrho
v:\,v\in\mathrm{H}^{\sigma}(\Omega)\},
\end{equation}
где $\sigma<0$ и функция $\varrho\in C^{\infty}(\Omega)$ положительна. Наделим это
пространство скалярным произведением
\begin{equation}\label{4.186}
\quad (f_{1},f_{2})_{\varrho\mathrm{H}^{\sigma}(\Omega)}:=
(\varrho^{-1}f_{1},\varrho^{-1}f_{2})_{\mathrm{H}^{\sigma}(\Omega)}
\end{equation}
и соответствующей ему нормой. Пространство $\varrho\mathrm{H}^{\sigma}(\Omega)$
полно и непрерывно вложено в $\mathcal{D}'(\Omega)$. Это следует из того, что
оператор умножения на функцию $\varrho$ непрерывен в $\mathcal{D}'(\Omega)$ и
устанавливает гомеоморфизм полного пространства $\mathrm{H}^{\sigma}(\Omega)$ на
$\varrho\mathrm{H}^{\sigma}(\Omega)$.

Рассмотрим весовые функции вида $\varrho:=\varrho_{1}^{-\sigma}$, где
\begin{equation}\label{4.187}
\varrho_{1}\in C^{\infty}(\,\overline{\Omega}\,),\;\;\varrho_{1}>0\;\,\mbox{в
$\Omega$},\;\;\varrho_{1}(x)=\mathrm{dist}(x,\Gamma)\;\,\mbox{в окрестности
$\Gamma$}.
\end{equation}

\medskip

\textbf{Теорема LM$_{\mathbf{2}}$} [\ref{LionsMagenes71}, с.~216]. \it Пусть
$\sigma<0$ и
\begin{equation}\label{4.188}
X^{\sigma}(\Omega):=
\begin{cases}
\;\varrho_{1}^{-\sigma}\mathrm{H}^{\sigma}(\Omega)
\quad\;\mbox{при}\;\;\sigma+1/2\notin\mathbb{Z},\\
\;[\,\varrho_{1}^{-\sigma+1/2}\,\mathrm{H}^{\sigma-1/2}(\Omega),\,
\varrho_{1}^{-\sigma-1/2}\,\mathrm{H}^{\sigma+1/2}(\Omega)]_{t^{1/2}}\\
\qquad\qquad\qquad \mbox{при}\;\;\sigma+1/2\in\mathbb{Z}.
\end{cases}
\end{equation}
Тогда остается в силе утверждение теоремы $\mathrm{LM}_{1}$. \rm


\begin{remark}\label{rem4.10}
Ж.-Л. Лионс и Э. Мадженес [\ref{LionsMagenes71}, с.~197--199] использовали в
качестве $X^{\sigma}(\Omega)$ некоторое гильбертово пространство
$\Xi^{\sigma}(\Omega)$. Для целых параметров $\sigma\geq0$ оно определяется по
формуле \eqref{3.53}, для дробных $\sigma>0$~--- посредством интерполяции со
степенным параметром:
$$
\Xi^{\sigma}(\Omega):=
\bigl[\Xi^{[\sigma]}(\Omega),\Xi^{[\sigma]+1}(\Omega)\bigr]_{t^{\{\sigma\}}},
$$
и, наконец, для отрицательных $\sigma<0$~--- посредством перехода к двойственному
пространству (относительно скалярного произведения в $L_{2}(\Omega)$):
$\Xi^{\sigma}(\Omega):=(\Xi^{-\sigma}(\Omega))'$. Для произвольного $\sigma<0$
пространство $\Xi^{\sigma}(\Omega)$ совпадает (с точностью до эквивалентности норм)
с правой частью формулы \eqref{4.188}. Это вытекает из результата Лионса и Мадженеса
[\ref{LionsMagenes71}, с.~212] (следствие~7.4).
\end{remark}

\subsection{Ключевая индивидуальная теорема}\label{sec4.4.2}

Здесь мы докажем ключевую индивидуальную теорему о разрешимости краевой задачи
\eqref{4.1}. Эта теорема утверждает, что оператор \eqref{4.183} корректно определен,
ограничен и нетеров при $\sigma<0$, если гильбертово пространство
$X^{\sigma}(\Omega)\hookrightarrow\mathcal{D}'(\Omega)$ удовлетворяет приводимому
ниже условию~I$_{\sigma}$. Она послужит ключом для доказательства других
индивидуальных теорем.

\medskip

\textbf{Условие I$_{\sigma}$.} Множество $X^{\infty}(\Omega):=X^{\sigma}(\Omega)\cap
C^{\infty}(\,\overline{\Omega}\,)$ плотно в $X^{\sigma}(\Omega)$, и существует число
$c>0$ такое, что
\begin{equation}\label{4.189}
\|\mathcal{O}f\|_{H^{\sigma}(\mathbb{R}^{n})}\leq
c\,\|f\|_{X^{\sigma}(\Omega)}\quad\mbox{для любого}\quad f\in X^{\infty}(\Omega).
\end{equation}
Здесь $\mathcal{O}f(x)=f(x)$ при $x\in\overline{\Omega}$ и  $\mathcal{O}f(x)=0$ при
$x\in\mathbb{R}^{n}\setminus\overline{\Omega}$.

\medskip

Заметим, что чем меньше $\sigma$ тем слабее условие~I$_{\sigma}$ для одного и того
же пространства $X^{\sigma}(\Omega)$.


\begin{remark}\label{rem4.11}
Я.~А.~Ройтберг [\ref{Roitberg71}] (п.~2.4) рассмотрел условие на пространство
$X^{\sigma}(\Omega)$ несколько более сильное, чем наше условие~I$_{\sigma}$. Им
дополнительно требовалось, чтобы $C^{\infty}(\,\overline{\Omega}\,)\subset
X^{\sigma}(\Omega)$. При таком условии Ройтберг [\ref{Roitberg71}, с.~261;
\ref{Roitberg96}, с.~190] доказал ограниченность оператора \eqref{4.183} для всех
$\sigma<0$. Отметим, что это условие не охватывает важный случай
$X^{\sigma}(\Omega)=\{0\}$, а также некоторые весовые пространства
$X^{\sigma}(\Omega)=\varrho\mathrm{H}^{\sigma}(\Omega)$, рассмотренные нами ниже.
\end{remark}

\medskip

Сформулируем ключевую индивидуальную теорему.


\begin{theorem}\label{th4.26}
Пусть число $\sigma<0$, а $X^{\sigma}(\Omega)$ --- произвольное гильбертово
пространство, непрерывно вложенное в $\mathcal{D}'(\Omega)$ и удовлетворяющее
условию~$\mathrm{I}_{\sigma}$. Тогда справедливы следующие утверждения.
\begin{itemize}
\item[$\mathrm{(i)}$] Множество
$$
D^{\infty}_{L,X}(\Omega):=\{u\in C^{\infty}(\,\overline{\Omega}\,):\,Lu\in
X^{\sigma}(\Omega)\}
$$
плотно в пространстве $D^{\sigma+2q}_{L,X}(\Omega)$.
\item[$\mathrm{(ii)}$] Отображение $u\rightarrow(Lu,Bu)$,
где $u\in D^{\infty}_{L,X}(\Omega)$, продолжается по непрерывности до ограниченного
линейного оператора \eqref{4.183}.
\item[$\mathrm{(iii)}$] Оператор \eqref{4.183} нетеров. Его ядро совпадает с $N$, а
область значений состоит из всех векторов
$(f,g_{1},\ldots,g_{q})\in\mathbf{X}_{\sigma}(\Omega,\Gamma)$, удовлетворяющих
условию \eqref{4.4}.
\item[$\mathrm{(iv)}$] Если множество
$\mathcal{O}(X^{\infty}(\Omega))$ плотно в пространстве
$H^{\sigma}_{\overline{\Omega}}(\mathbb{R}^{n})$, то индекс оператора \eqref{4.183}
равен $\dim N-\dim N^{+}$.
\end{itemize}
\end{theorem}


\textbf{Доказательство} будет опираться на теоремы \ref{th4.15}, \ref{th4.19},
\ref{th4.24} и \ref{th4.25} для соболевского случая $\varphi\equiv1$, в которых мы
берем $s=\sigma+2q$. Напомним, что в этом случае, они доказаны Я.~А.~Ройтбергом
[\ref{Roitberg96}]. При этом теорема \ref{th4.25} установлена им для произвольного
$s<-1/2$.

Из условия теоремы следует, что отображение $f\mapsto\mathcal{O}f$, где $f\in
X^{\infty}(\Omega)$, продолжается по непрерывности до ограниченного линейного
оператора
\begin{equation}\label{4.190}
\mathcal{O}:\,X^{\sigma}(\Omega)\rightarrow
H^{\sigma}_{\overline{\Omega}}(\mathbb{R}^{n})=H^{\sigma,(0)}(\Omega).
\end{equation}
Этот оператор инъективный. В самом деле, пусть $\mathcal{O}f=0$ для некоторого
распределения $f\in X^{\sigma}(\Omega)$. Выберем последовательность $(f_{k})\subset
X^{\infty}(\Omega)$ такую, что $f_{k}\rightarrow f$ в
$X^{\sigma}(\Omega)\hookrightarrow\mathcal{D}'(\Omega)$. Тогда
$\mathcal{O}f_{k}\rightarrow0$ в
$H^{\sigma}_{\overline{\Omega}}(\mathbb{R}^{n})\hookrightarrow\mathcal{S}'(\mathbb{R}^{n})$,
откуда
$$
(f,v)_{\Omega}=\lim_{k\rightarrow\infty}(f_{k},v)_{\Omega}=
\lim_{k\rightarrow\infty}(\mathcal{O}f_{k},v)_{\Omega}=0\;\;\mbox{для всех}\;\;v\in
C^{\infty}_{0}(\Omega).
$$
Таким образом, $f=0$ как распределение из пространства
$X^{\sigma}(\Omega)\hookrightarrow\mathcal{D}'(\Omega)$, т.~е. оператор
\eqref{4.190} инъективный. Он определяет непрерывное вложение
$X^{\sigma}(\Omega)\hookrightarrow H^{\sigma,(0)}(\Omega)$.

Согласно теореме~\ref{th4.14}, для произвольного $u\in H^{\sigma+2q,(2q)}(\Omega)$
определен посредством предельного перехода элемент $Lu\in H^{\sigma,(0)}(\Omega)$.
Положим
$$
D^{\sigma+2q,(2q)}_{L,X}(\Omega):=\{\,u\in H^{\sigma+2q,(2q)}(\Omega):\,Lu\in
X^{\sigma}(\Omega)\,\}.
$$
В пространстве $D^{\sigma+2q,(2q)}_{L,X}(\Omega)$ введем скалярное произведение
графика
$$
(u_{1},u_{2})_{D^{\sigma+2q,(2q)}_{L,X}(\Omega)}:=(u_{1},u_{2})_{H^{\sigma+2q,(2q)}(\Omega)}+
(Lu_{1},Lu_{2})_{X^{\sigma}(\Omega)}.
$$
Пространство $D^{\sigma+2q,(2q)}_{L,X}(\Omega)$ полное относительно этого скалярного
произведения. Действительно, пусть последовательность $(u_{k})$ фундаментальная в
$D^{\sigma+2q,(2q)}_{L,X}(\Omega)$. Тогда, поскольку пространства
$H^{\sigma+2q,(2q)}(\Omega)$ и $X^{\sigma}(\Omega)$ полные, существуют пределы
$u:=\lim u_{k}$ в $H^{\sigma+2q,(2q)}(\Omega)$ и $f:=\lim Lu_{k}$ в
$X^{\sigma}(\Omega)$. Первый из них влечет сходимость $\lim Lu_{k}=Lu$ в
$H^{\sigma,(0)}(\Omega)$. Отсюда в силу справедливости второго предела и
непрерывности вложения \eqref{4.190}, получаем: $Lu=f\in X^{\sigma}(\Omega)$.
Следовательно, $u\in D^{\sigma+2q,(2q)}_{L,X}(\Omega)$ и $\lim u_{k}=u$ в
$D^{\sigma+2q,(2q)}_{L,X}(\Omega)$, т.~е. пространство
$D^{\sigma+2q,(2q)}_{L,X}(\Omega)$ полное.

На основании теоремы~\ref{th4.15} делаем вывод, что сужение оператора
\begin{gather}\notag
(L,B):H^{\sigma+2q,(2q)}(\Omega)\rightarrow\\  H^{\sigma,(0)}(\Omega)\oplus
\bigoplus_{j=1}^{q}H^{\sigma+2q-m_{j}-1/2}(\Gamma)=:
\mathbf{H}_{\sigma,(0)}(\Omega,\Gamma).\label{4.191}
\end{gather}
на пространство $D^{\sigma+2q,(2q)}_{L,X}(\Omega)$ является ограниченным оператором
\begin{equation}\label{4.192}
(L,B):\,D^{\sigma+2q,(2q)}_{L,X}(\Omega)\rightarrow\mathbf{X}_{\sigma}(\Omega,\Gamma).
\end{equation}
Ядро оператора \eqref{4.192} равно $N$, а область значений состоит из всех векторов
$(f,g_{1},\ldots,g_{q})\in\mathbf{X}_{\sigma}(\Omega,\Gamma)$, удовлетворяющих
условию \eqref{4.4}. Следовательно, оператор  \eqref{4.192} нетеров, причем его
коядро имеет размерность $\beta\leq\dim N^{+}$.

Кроме того, если множество $\mathcal{O}(X^{\infty}(\Omega)$ плотно в пространстве
$H^{\sigma}_{\overline{\Omega}}(\mathbb{R}^{n})$, то $\beta=\nobreak\dim N^{+}$. В
самом деле, обозначив оператор \eqref{4.191} через $\Lambda$, а более узкий оператор
\eqref{4.192} через $\Lambda_{0}$, рассмотрим сопряженные к ним операторы
$\Lambda^{\ast}$ и $\Lambda_{0}^{\ast}$. Так как теперь непрерывное вложение
$\mathbf{X}_{\sigma}(\Omega,\Gamma)\hookrightarrow
\mathbf{H}_{\sigma,(0)}(\Omega,\Gamma)$ плотно, то
$\ker\Lambda_{0}^{\ast}\supseteq\ker\Lambda^{\ast}$. Следовательно,
\begin{gather*}
\beta=\dim\mathrm{coker}\,\Lambda_{0}=\dim\ker\Lambda_{0}^{\ast}\geq \\
\dim\ker\Lambda^{\ast}=\dim\mathrm{coker}\,\Lambda=\dim N^{+}.
\end{gather*}
Значит, $\beta=\dim N^{+}$ и индекс оператора \eqref{4.192} равен $\dim N-\dim
N^{+}$ в рассмотренной ситуации.

Покажем, что множество $D^{\infty}_{L,X}(\Omega)$ плотно в пространстве
$D^{\sigma+2q,(2q)}_{L,X}(\Omega)$. Поскольку множество
$X^{\infty}(\Omega)\times(C^{\infty}(\Gamma))^{q}$ плотно в пространстве
$\mathbf{X}_{\sigma}(\Omega,\Gamma)$, то согласно лемме Гохберга-Крейна
[\ref{GohbergKrein57}] (лемма 2.1) можно записать
\begin{equation}\label{4.193}
\mathbf{X}_{\sigma}(\Omega,\Gamma)=(L,B)\bigl(D^{\sigma+2q,(2q)}_{L,X}(\Omega)\bigr)\dotplus
\mathcal{Q}(\Omega,\Gamma),
\end{equation}
где $\mathcal{Q}(\Omega,\Gamma)$ --- некоторое конечномерное пространство,
удовлетворяющее условию
\begin{equation}\label{4.194}
\mathcal{Q}(\Omega,\Gamma)\subset X^{\infty}(\Omega)\times(C^{\infty}(\Gamma))^{q}.
\end{equation}
Обозначим через $\Pi$ оператор проектирования пространства
$\mathcal{X}_{\sigma}(\Omega,\Gamma)$ на первое слагаемое в \eqref{4.193}
параллельно второму слагаемому.

Пусть $u\in D^{\sigma+2q,(2q)}_{L,X}(\Omega)$. Аппроксимируем $F:=(L,B)u$ некоторой
последовательностью $(F_{k})\subset
X^{\infty}(\Omega)\times(C^{\infty}(\Gamma))^{q}$ в норме пространства
$\mathbf{X}_{\sigma}(\Omega,\Gamma)$. Тогда
\begin{equation}\label{4.195}
\lim_{k\rightarrow\infty}\Pi F_{k}=\Pi F=F=(L,B)u\quad
\mbox{в}\quad\mathbf{X}_{\sigma}(\Omega,\Gamma),
\end{equation}
где в силу \eqref{4.194}
\begin{equation}\label{4.196}
(\Pi F_{k})\subset X^{\infty}(\Omega)\times(C^{\infty}(\Gamma))^{q}
\end{equation}

Нетеров оператор \eqref{4.192} естественным образом порождает гомеоморфизм
$$
\Lambda_{0}:=(L,B):\,D^{\sigma+2q,(2q)}_{L,X}(\Omega)/N\leftrightarrow
\Pi(\mathbf{X}_{\sigma}(\Omega,\Gamma)).
$$
В силу \eqref{4.195}
$$
\lim_{k\rightarrow\infty}\Lambda_{0}^{-1}\,\Pi F_{k}=\{u+w:w\in
N\}\quad\mbox{в}\quad D^{\sigma+2q,(2q)}_{L,X}(\Omega)/N.
$$
Следовательно, существует последовательность представителей $u_{k}\in
D^{\sigma+2q,(2q)}_{L,X}(\Omega)$ классов смежности $\Lambda_{0}^{-1}\,\Pi F_{k}$
такая, что
\begin{equation}\label{4.197}
\lim_{k\rightarrow\infty}u_{k}=u\quad\mbox{в}\quad D^{\sigma+2q,(2q)}_{L,X}(\Omega).
\end{equation}
При этом в силу \eqref{4.196}
$$
(L,B)u_{k}=\Pi F_{k}\in
C^{\infty}(\,\overline{\Omega}\,)\times(C^{\infty}(\Gamma))^{q}.
$$
Отсюда на основании теоремы \ref{th4.19} и теоремы вложения Соболева получаем:
$$
u_{k}\in\bigcap_{s>2q}\,H^{s,(2q)}(\Omega)=
\bigcap_{s>2q}\,H^{s}(\Omega)=C^{\infty}(\,\overline{\Omega}\,).
$$
Таким образом, в \eqref{4.197} последовательность $(u_{k})\subset
D^{\infty}_{L,X}(\Omega)$. Тем самым доказано, что $D^{\infty}_{L,X}(\Omega)$ плотно
в $D^{\sigma+2q,(2q)}_{L,X}(\Omega)$.

Далее рассмотрим отдельно случаи $-2q-1/2\leq\sigma<0$ и $\sigma<-2q-1/2$.

Первый случай: $-2q-1/2\leq\sigma<0$. Тогда
$H^{\sigma+2q,(0)}(\Omega)=\mathrm{H}^{\sigma+2q}(\Omega)$. Действительно, так как
$H^{\lambda}_{0}(\Omega)=H^{\lambda}(\Omega)$ при $0\leq\lambda\leq1/2$ (см.
[\ref{LionsMagenes71}], гл.~1, теорема~11.1), то
$$
H^{s,(0)}(\Omega)=(H^{-s}(\Omega))'=(H^{-s}_{0}(\Omega))'=\mathrm{H}^{s}(\Omega)
\;\;\mbox{при}\;\;-1/2\leq s<0.
$$
Отсюда
\begin{equation}\label{4.198}
H^{s,(0)}(\Omega)=\mathrm{H}^{s}(\Omega)\quad\mbox{при}\quad s\geq-1/2
\end{equation}
с равенством норм.

Воспользуемся теоремой~\ref{th4.24} и рассмотрим отображение $I_{0}:u\mapsto u_{0}$,
где $u\in D^{\sigma+2q,(2q)}_{L,X}(\Omega)$ и $(u_{0},f):=I_{L}u$. Оно определяет
гомеоморфизм
\begin{equation}\label{4.199}
I_{0}:\,D^{\sigma+2q,(2q)}_{L,X}(\Omega)\leftrightarrow D^{\sigma+2q}_{L,X}(\Omega).
\end{equation}
В самом деле, для произвольных распределений
$$
u_{0}\in H^{\sigma+2q,(0)}(\Omega)=\mathrm{H}^{\sigma+2q}(\Omega)\quad\mbox{и}\quad
f\in X^{\sigma}(\Omega)\hookrightarrow H^{\sigma,(0)}(\Omega)
$$
условия \eqref{4.118} и \eqref{4.119} эквивалентны в силу леммы \ref{lem4.5}.
Условие \eqref{4.119} означает, что $Lu_{0}=f$ в смысле теории распределений в
области $\Omega$. Отсюда на основании теоремы~\ref{th4.24} получаем равенство
$I_{0}(D^{\sigma+2q,(2q)}_{L,X}(\Omega))=D^{\sigma+2q}_{L,X}(\Omega)$. Кроме того,
для $u\in D^{\sigma+2q,(2q)}_{L,X}(\Omega)$ и $(u_{0},f)=I_{L}u$ выполняется
эквивалентность норм:
\begin{gather*}
\|u\|_{D^{\sigma+2q,(2q)}_{L,X}(\Omega)}^{2}=\|u\|_{H^{\sigma+2q,(2q)}(\Omega)}^{2}+
\|f\|_{X^{\sigma}(\Omega)}^{2}\asymp \\
\|u_{0}\|_{H^{\sigma+2q,(0)}(\Omega)}^{2}+\|f\|_{H^{\sigma,(0)}(\Omega)}^{2}+
\|f\|_{X^{\sigma}(\Omega)}^{2}\asymp \\
\|u_{0}\|_{\mathrm{H}^{\sigma+2q}(\Omega)}^{2}+\|f\|_{X^{\sigma}(\Omega)}^{2}=
\|u_{0}\|_{D^{\sigma+2q}_{L,X}(\Omega)}^{2}.
\end{gather*}
Следовательно, отображение $I_{0}$ определяет гомеоморфизм \eqref{4.199}.

Из свойств операторов \eqref{4.199} и \eqref{4.192} (последний, напомним, обозначен
через $\Lambda_{0}$) следует, что оператор
\begin{equation}\label{4.200}
\Lambda_{0}I_{0}^{-1}:\,D^{\sigma+2q}_{L,X}(\Omega)\rightarrow
\mathbf{X}_{\sigma}(\Omega,\Gamma)
\end{equation}
ограничен и нетеров. При этом его ядро, область значений и индекс те же, что и у
оператора \eqref{4.192}. Так как оператор $I_{0}$ осуществляет взаимно однозначное
отображение множества $D^{\infty}_{L,X}(\Omega)$ на себя, то это множество плотно в
пространстве $D^{\sigma+2q}_{L,X}(\Omega)$, а оператор \eqref{4.200} является
продолжением по непрерывности отображения $u\rightarrow(Lu,Bu)$, где $u\in
D^{\infty}_{L,X}(\Omega)$. Теорема~\ref{th4.26} доказана в первом случае.

Второй случай: $\sigma<-2q-1/2$. Тогда
$H^{\sigma+2q,(0)}(\Omega)=H^{\sigma+2q}_{\overline{\Omega}}(\mathbb{R}^{n})$.
Положим $Rw:=w\!\upharpoonright\!\Omega$ для $w\in\mathcal{D}'(\mathbb{R}^{n})$.
Докажем, что отображение $I_{0}:u\mapsto Ru_{0}$, где $u\in
D^{\sigma+2q,(2q)}_{L,X}(\Omega)$ и $(u_{0},f):=I_{L}u$ определяет гомеоморфизм
\eqref{4.199} в рассматриваемом случае (в первом случае было $Ru_{0}=u_{0}$).
Воспользуемся теоремой~\ref{th4.24} и заметим, что
$\eqref{4.118}\Rightarrow\eqref{4.119}$. Для произвольного $u\in
D^{\sigma+2q,(2q)}_{L,X}(\Omega)$ имеем: $Ru_{0}\in\mathrm{H}^{\sigma+2q}(\Omega)$
(см. \eqref{4.176}), $f=Lu\in X^{\sigma}(\Omega)$ и
$$
(Ru_{0},L^{+}v)_{\Omega}=(u_{0},L^{+}v)_{\Omega}=(f,v)_{\Omega}\quad\mbox{для
всех}\quad v\in C^{\infty}_{0}(\Omega),
$$
т.~е. $LRu_{0}=f$ в смысле теории распределений в области $\Omega$. Следовательно,
$I_{0}u=Ru_{0}\in D^{\sigma+2q}_{L,X}(\Omega)$. Кроме того, в силу \eqref{4.177} и
определения пространства $H^{\sigma+2q,(2q)}(\Omega)$ справедлива оценка
\begin{gather*}
\|I_{0}u\|_{D^{\sigma+2q}_{L,X}(\Omega)}^{2}=
\|Ru_{0}\|_{\mathrm{H}^{\sigma+2q}(\Omega)}^{2}+\|f\|_{X^{\sigma}(\Omega)}^{2}\leq \\
\|u_{0}\|_{H^{\sigma+2q}(\mathbb{R}^{n})}^{2}+\|f\|_{X^{\sigma}(\Omega)}^{2}\leq
\|u\|_{D^{\sigma+2q,(2q)}_{L,X}(\Omega)}^{2}.
\end{gather*}
Следовательно, оператор $I_{0}:D^{\sigma+2q,(2q)}_{L,X}(\Omega)\rightarrow
D^{\sigma+2q}_{L,X}(\Omega)$ ограниченный.

Покажем, что этот оператор биективный. Пусть $\omega\in D^{\sigma+2q}_{L,X}(\Omega)$
и $f:=L\omega\in X^{\sigma}(\Omega)$. В~силу \eqref{4.176} найдется распределение
$u_{0}\in H^{\sigma+2q}_{\overline{\Omega}}(\mathbb{R}^{n})$ такое, что
$\omega=Ru_{0}$. При этом выполняется равенство \eqref{4.119}:
$$
(u_{0},L^{+}v)_{\Omega}=(\omega,L^{+}v)_{\Omega}=(f,v)_{\Omega}\quad\mbox{для
всех}\quad v\in C^{\infty}_{0}(\Omega).
$$
Согласно теореме~\ref{th4.25}, для распределений $u_{0}\in
H^{\sigma+2q,(0)}(\Omega)$ и $f\in X^{\sigma}(\Omega)\hookrightarrow
H^{\sigma,(0)}(\Omega)$ существует единственная пара $(u_{0}^{\ast},f)\in
K_{\sigma+2q,L}(\Omega)$, удовлетворяющая условию \eqref{4.121}. (В соболевском
случае $\varphi\equiv1$ эта теорема верна для $s\in\mathbb{R}$, как показал
Я.~А.~Ройтберг [\ref{Roitberg96}, с.~187], теорема 6.2.1.) Отсюда на основании
теоремы~\ref{th4.24} можем записать
$$
u^{\ast}:=I_{L}^{-1}(u_{0}^{\ast},f)\in
D^{\sigma+2q,(2q)}_{L,X}(\Omega)\quad\mbox{и}\quad
I_{0}u^{\ast}=Ru_{0}^{\ast}=Ru_{0}=\omega.
$$
Элемент $u^{\ast}$ --- единственный прообраз распределения $\omega$ при отображении
$I_{0}$. Действительно, если $I_{0}u'=\omega$ для некоторого $u'\in
D^{\sigma+2q,(2q)}_{L,X}(\Omega)$, то пара $(u'_{0},f'):=I_{L}u'\in
K_{\sigma+2q,L}(\Omega)$ удовлетворяет условиям:
\begin{gather*}
f'=LRu'_{0}=L\omega=f,\\
(u'_{0},v)_{\Omega}=(\omega,v)_{\Omega}=(u_{0},v)_{\Omega}\quad\mbox{для всех}\quad
v\in C^{\infty}_{0}(\Omega).
\end{gather*}
Поэтому в силу теоремы~\ref{th4.25} пары $(u'_{0},f')=(u'_{0},f)$ и
$(u_{0}^{\ast},f)$ равны, а значит равны их прообразы $u'$ и $u^{\ast}$ при
отображении $I_{L}$.

Таким образом, линейный ограниченный оператор \eqref{4.199} биективный в
рассматриваемом случае. Следовательно, по теореме Банаха об обратном операторе, он
является гомеоморфизмом. Теперь, воспользовавшись нетеровостью оператора
\eqref{4.192} и рассуждая также как и в первом случае, мы завершаем доказательство
во втором случае.

Теорема~\ref{th4.26} доказана.


\begin{remark}\label{rem4.12}
Утверждение, подобное теореме \ref{th4.26}, доказано в [\ref{Magenes66}, с.~202]
(теорема 6.16) в случае неполуцелого $\sigma\leq-2q$ и краевой задачи Дирихле. При
этом на пространство $X^{\sigma}(\Omega)$ накладывались некоторые иные условия,
зависящие от рассматриваемой краевой задачи. Наше условие I$_{\sigma}$ от нее не
зависит.
\end{remark}


\begin{remark}\label{rem4.13}
В теореме \ref{th4.26}, как и в индивидуальных теоремах LM$_{1}$, LM$_{2}$, решение
и правая часть эллиптического уравнения $Lu=f$ являются распределениями, заданными в
$\Omega$. Иные индивидуальные теоремы доказаны в работах Я.~А.~Ройтберга
[\ref{Roitberg68}] и Ю.~В.~Костарчука, Я.~А.~Ройтберга [\ref{KostarchukRoitberg73}].
В~этих теоремах $u$ и/или $f$ не являются распределениями в области $\Omega$.
\end{remark}


Пространство $X^{\sigma}(\Omega):=\{0\}$, очевидно, удовлетворяет
условию~I$_{\sigma}$. В~этом случае теорема~\ref{th4.26} описывает свойства
полуоднородной краевой задачи \eqref{4.1}, где $f=0$, и справедлива для любого
$\sigma\in\mathbb{R}$ ввиду теоремы~\ref{th3.11}. При этом, хотя и
$\mathrm{H}^{s}(\Omega)\neq H^{s}(\Omega)$ для полуцелого $s=\sigma+2q<0$, но
\begin{equation}\label{4.201}
\{u\in\mathrm{H}^{s}(\Omega):\,Lu=0\;\mbox{в}\;\Omega\}=\{u\in
H^{s}(\Omega):\,Lu=0\;\mbox{в}\;\Omega\}.
\end{equation}
Кроме того, нормы в пространствах $\mathrm{H}^{s}(\Omega)$ и $H^{s}(\Omega)$
эквивалентны на распределениях $u$, указанных в формуле \eqref{4.201}.

\medskip

Далее мы рассмотрим различные приложения теоремы \ref{th4.26}, обусловленные
конкретным выбором пространства $X^{\sigma}(\Omega)$. Мы изучим отдельно случаи,
когда в качестве $X^{\sigma}(\Omega)$ берутся гильбертовы пространства Соболева и их
весовые аналоги.

\subsection[Индивидуальная теорема для пространств Соболева]
{Индивидуальная теорема \\ для пространств Соболева}\label{sec4.4.3}

В~следующей теореме описаны все гильбертовы (невесовые) пространств Соболева,
удовлетворяющие условию I$_{\sigma}$.


\begin{theorem}\label{th4.27}
Пусть $\sigma<0$ и $\lambda\in\mathbb{R}$. Пространство
$X^{\sigma}(\Omega):=\mathrm{H}^{\lambda}(\Omega)$ удовлетворяет условию
$\mathrm{I}_{\sigma}$ тогда и только тогда, когда
\begin{equation}\label{4.202}
\lambda\geq\max\{\sigma,-1/2\}.
\end{equation}
\end{theorem}


\textbf{Доказательство.} Заметим, что так как
$X^{\sigma}(\Omega)=\mathrm{H}^{\lambda}(\Omega)$, то множество
$X^{\infty}(\Omega)=C^{\infty}(\,\overline{\Omega}\,)$ плотно в пространстве
$X^{\sigma}(\Omega)$.

\emph{Достаточность.}  Предположим, что выполняется неравенство \eqref{4.202}. Тогда
ввиду \eqref{4.198}
$$
\mathrm{H}^{\lambda}(\Omega)=H^{\lambda,(0)}(\Omega)\hookrightarrow
H^{\sigma,(0)}(\Omega)
$$
вместе с топологией. Напомним, что всякая функция $f\in
C^{\infty}(\,\overline{\Omega}\,)$ отождествляется с функционалом
$(f,\cdot\,)_{\Omega}$, а он в свою очередь --- с функцией $\mathcal{O}f$ как
элементом пространства
$H^{\sigma}_{\overline{\Omega}}(\mathbb{R}^{n})=H^{\sigma,(0)}(\Omega)$.
Следовательно,
$$
\|\mathcal{O}f\|_{H^{\sigma}(\mathbb{R}^{n})}\leq
c\,\|f\|_{\mathrm{H}^{\lambda}(\Omega)}\quad\mbox{для любого}\quad f\in
C^{\infty}(\,\overline{\Omega}\,),
$$
где число $c>0$ не зависит от $f$. Достаточность доказана.

\emph{Необходимость.} Предположим, что пространство
$X^{\sigma}(\Omega):=\mathrm{H}^{\lambda}(\Omega)$ удовлетворяет
условию~I$_{\sigma}$. Если $\lambda\geq0$, то неравенство \eqref{4.202} выполняется.
Поэтому ограничимся случаем $\lambda<0$. Оператор \eqref{4.190} определяет
непрерывное плотное вложение $\mathrm{H}^{\lambda}(\Omega)\hookrightarrow
H^{\sigma,(0)}(\Omega)$. Отсюда
$$
H^{-\sigma}(\Omega)=(H^{\sigma,(0)}(\Omega))'\subseteq(\mathrm{H}^{\lambda}(\Omega))'=
H^{-\lambda}_{0}(\Omega).
$$
Следовательно, $-\sigma\geq-\lambda$. Кроме того, $-\lambda\leq1/2$, поскольку если
$-\lambda>1/2$, то функция $f\equiv1\in H^{-\sigma}(\Omega)$ не принадлежит
пространству $H^{-\lambda}_{0}(\Omega)$ в силу теоремы \ref{th3.20} (i). Таким
образом, $\sigma$ удовлетворяет неравенству \eqref{4.202}. Необходимость доказана.

Теорема \ref{th4.27} доказана.

\medskip

Из теорем \ref{th4.26} и \ref{th4.27} вытекает следующая индивидуальная теорема о
разрешимости краевой задачи \eqref{4.1} в пространствах Соболева.


\begin{theorem}\label{th4.28}
Пусть $\sigma<0$ и $\lambda\geq\max\{\sigma,-1/2\}$. Тогда отображение \eqref{4.179}
продолжается по непрерывности до ограниченного нетерового оператора
\begin{gather}\notag
(L,B):\,\{u\in\mathrm{H}^{\sigma+2q}(\Omega):Lu\in
\mathrm{H}^{\lambda}(\Omega)\}\rightarrow \\
\mathrm{H}^{\lambda}(\Omega)\oplus\bigoplus_{j=1}^{q}\,H^{\sigma+2q-m_{j}-1/2}(\Gamma),
\label{4.203}
\end{gather}
область определения которого является гильбертовым пространством относительно нормы
графика
$$
\bigl(\,\|u\|_{\mathrm{H}^{\sigma+2q}(\Omega)}^{2}+
\|Lu\|_{\mathrm{H}^{\lambda}(\Omega)}^{2}\bigr)^{1/2}.
$$
Индекс оператора \eqref{4.203} равен $\dim\ N-\dim N^{+}$ и не зависит от $\sigma$ и
$\lambda$.
\end{theorem}


\textbf{Доказательство.}  Ограниченность и нетеровость оператора \eqref{4.203} сразу
следуют из теорем \ref{th4.26} и \ref{th4.27}, где
$X^{\sigma}(\Omega):=H^{\lambda}(\Omega)$ и
$D^{\infty}_{L,X}(\Omega)=C^{\infty}(\,\overline{\Omega}\,)$. Кроме того, так как
множество $\mathcal{O}(C^{\infty}(\,\overline{\Omega}\,))$, отождествляемое с
$C^{\infty}(\,\overline{\Omega}\,)$, плотно в пространстве
$H^{\sigma}_{\overline{\Omega}}(\mathbb{R}^{n})=H^{\sigma,(0)}(\Omega)$, то в силу
утверждения (iv) теоремы~\ref{th4.26} индекс оператора \eqref{4.203} равен
$\dim{N}-\dim{N}^{+}$ и не зависит от $\sigma$ и~$\lambda$. Теорема \ref{th4.28}
доказана.

\medskip

Теорема LM$_{1}$ являются частным случаем теоремы \ref{th4.28}, отвечающим случаю,
когда $\sigma=0$, т.~е. $X^{\sigma}(\Omega)=L_{2}(\Omega)$.

В теореме~\ref{th4.28} допускаются некоторые пространства $X^{\sigma}(\Omega)$,
содержащие в себе $L_{2}(\Omega)$. Среди них наиболее широким является пространство
$X^{\sigma}(\Omega)=\mathrm{H}^{-1/2}(\Omega)$ при $\sigma\leq-1/2$.

Отметим, что в случае $-1/2\leq\sigma=\lambda<0$ область определения оператора
\eqref{4.203} не зависит от $L$. В самом деле, если $-1/2<\sigma=\lambda<0$, то
поскольку число $\sigma$ неполуцелое, имеем ограниченный оператор $L:\nobreak
\mathrm{H}^{\sigma+2q}(\Omega)\rightarrow\mathrm{H}^{\sigma}(\Omega)$ (см.
[\ref{LionsMagenes71}, с.~107], гл.~1, предложение 12.1). Отсюда следует, что
область определения оператора \eqref{4.203} совпадает с
$\mathrm{H}^{\sigma+2q}(\Omega)$ и не зависит от $L$. В этом случае
теорема~\ref{th4.28} совпадает с теоремой~7.5 из монографии Ж.-Л.~Лионса и
Э.~Мадженеса [\ref{LionsMagenes71}, с.~218], доказанной в предположении, что
$N=N^{+}=\{0\}$.

Если $\sigma=\lambda=-1/2$, то поскольку пространство $\mathrm{H}^{-1/2}(\Omega)$
\'уже, чем $L(\mathrm{H}^{2q-1/2}(\Omega))$, предыдущие рассуждения не проходят.
Однако, в силу теоремы~\ref{th4.26} (i), равенства \eqref{4.198} и
теоремы~\ref{th4.24}, область определения оператора \eqref{4.203} является
пополнением множества функций $u\in C^{\infty}(\,\overline{\Omega}\,)$ по норме
\begin{gather*}
\|u\|_{\mathrm{H}^{2q-1/2}(\Omega)}^{2}+\|Lu\|_{\mathrm{H}^{-1/2}(\Omega)}^{2}=\\
\|u\|_{H^{2q-1/2,(0)}(\Omega)}^{2}+\|Lu\|_{H^{-1/2,(0)}(\Omega)}^{2}\asymp
\|u\|_{H^{2q-1/2,(2q)}(\Omega)}.
\end{gather*} Следовательно, она совпадает с
пространством $H^{2q-1/2,(2q)}(\Omega)$, которое не зависит от $L$.

В заключение этого пункта отметим, что условию I$_{\sigma}$ при $\sigma<\nobreak0$
удовлетворяют также пространства Хермандера $H^{\lambda,\varphi}(\Omega)$, где
$\lambda>\max\{\sigma,-1/2\}$ и $\varphi\in\mathcal{M}$. Индивидуальные теоремы,
связанные с пространствами Хермандера, будут рассмотрены ниже в п.~\ref{sec4.5}.

\subsection[Индивидуальная теорема для весовых пространств]
{Индивидуальная теорема \\ для весовых пространств}\label{sec4.4.4}

В~теореме \ref{th4.27} всегда
$X^{\sigma}(\Omega)\subseteq\mathrm{H}^{-1/2}(\Omega)$. Пространство
$X^{\sigma}(\Omega)$, содержащее широкий класс распределений $f\notin
\mathrm{H}^{-1/2}(\Omega)$ и удовлетворяющее условию I$_{\sigma}$, можно получить,
используя некоторые весовые пространства $\varrho\mathrm{H}^{\sigma}(\Omega)$.

Функцию $\varrho:\Omega\rightarrow\mathbb{C}$ называют \emph{мультипликатором} в
пространстве $\mathrm{H}^{\lambda}(\Omega)$, где $\lambda\geq0$, если оператор
умножения на $\varrho$ ограничен на этом пространстве. Этот оператор мы обозначаем
через $M_{\varrho}$. Описание класса всех мультипликаторов в пространстве
$\mathrm{H}^{\lambda}(\Omega)=H^{\lambda}(\Omega)$ при $\lambda\geq0$ дано в
монографии В.~Г.~Мазьи и Т.~О.~Шапошниковой [\ref{MazyaShaposhnikova86}] (п.~6.3.3).

Пусть $\sigma<-1/2$. Рассмотрим следующее условие на функцию~$\varrho$.

\medskip

\textbf{Условие II$_{\sigma}$.} Функция $\varrho$ является мультипликатором в
пространстве $\mathrm{H}^{-\sigma}(\Omega)$ и
\begin{equation}\label{4.204}
D_{\nu}^{j}\,\varrho=0\;\;\mbox{на}\;\;\Gamma\;\;\mbox{для всех}\;\;
j\in\mathbb{Z},\;\;0\leq j<-\sigma-1/2.
\end{equation}


Отметим, что если $\varrho$ является мультипликатором в пространстве
$\mathrm{H}^{-\sigma}(\Omega)$, то, очевидно,
$\varrho\in\mathrm{H}^{-\sigma}(\Omega)$. Следовательно, в силу теоремы~\ref{th3.5}
существует след $(D_{\nu}^{j}\varrho)\!\!\upharpoonright\!\Gamma\in
H^{-\sigma-j-1/2}(\Gamma)$ для любого целого $j\geq0$ такого, что
$-\sigma-j-1/2>\nobreak0$. Поэтому условие II$_{\sigma}$ сформулировано корректно.


\begin{theorem}\label{th4.29}
Пусть $\sigma<-1/2$ и функция $\varrho\in C^{\infty}(\Omega)$ положительна.
Пространство $X^{\sigma}(\Omega):=\varrho\mathrm{H}^{\sigma}(\Omega)$ удовлетворяет
условию~$\mathrm{I}_{\sigma}$ тогда и только тогда, когда функция $\varrho$
удовлетворяет условию $\mathrm{II}_{\sigma}$.
\end{theorem}


Для доказательства этой теоремы нам понадобится следующая лемма.


\begin{lemma}\label{lem4.8} Пусть $\sigma<-1/2$.
Умножение на функцию $\varrho$ является ограниченным оператором
\begin{equation}\label{4.205}
M_{\varrho}:\,H^{-\sigma}(\Omega)\rightarrow H^{-\sigma}_{0}(\Omega)
\end{equation}
тогда и только тогда, когда выполняется условие~$\mathrm{II}_{\sigma}$.
\end{lemma}


\textbf{Доказательство.} \emph{Необходимость.} Если умножение на функцию $\varrho$
определяет ограниченный оператор \eqref{4.205}, то эта функция является
мультипликатором в пространстве $H^{-\sigma}(\Omega)$ и, кроме того, $\varrho\in
H^{-\sigma}_{0}(\Omega)$. Как известно [\ref{Triebel80}, с.~412] (теорема 4.7.1
(а)),
\begin{equation}\label{4.206}
\varrho\in H^{-\sigma}_{0}(\Omega)\;\Leftrightarrow\; (\varrho\in
H^{-\sigma}(\Omega)\;\;\mbox{и верно \eqref{4.204}})
\end{equation}
Отсюда следует условие~II$_{\sigma}$. Необходимость доказана.

\emph{Достаточность.} Пусть функция $\varrho$ удовлетворяет условию~II$_{\sigma}$.
Надо доказать только лишь, что $\varrho u\in H^{-\sigma}_{0}(\Omega)$ для любого
$u\in H^{-\sigma}(\Omega)$. Условие~II$_{\sigma}$ влечет в силу~\eqref{4.206}
включение $\varrho\in H^{-\sigma}_{0}(\Omega)$. Пусть последовательности
$(u_{k})\subset C^{\infty}(\,\overline{\Omega}\,)$ и $(\varrho_{j})\subset
C^{\infty}_{0}(\Omega)$ такие, что $u_{k}\rightarrow u$ и
$\varrho_{j}\rightarrow\varrho$ в $H^{-\sigma}(\Omega)$. Поскольку функции $\varrho$
и $u_{k}$~--- мультипликаторы в пространстве $H^{-\sigma}(\Omega)$, то в нем
$$
\lim_{k\rightarrow\infty}(\varrho u_{k})=\varrho u\quad\mbox{и}\quad
\lim_{j\rightarrow\infty}(\varrho_{j}u_{k})=\varrho u_{k}\quad\mbox{для всех $k$}.
$$
Отсюда ввиду $\varrho_{j}u_{k}\in C^{\infty}_{0}(\Omega)$ следует, что $\varrho u\in
H^{-\sigma}_{0}(\Omega)$. Достаточность доказана.

Лемма \ref{lem4.8} доказана.

\medskip

\textbf{Доказательство теоремы \ref{th4.29}.} По условию, $\sigma<-1/2$ и функция
$\varrho\in C^{\infty}(\Omega)$ положительна. Обозначим через $M_{\varrho}$ и
$M_{\varrho^{-1}}$ операторы умножения на функции $\varrho$ и $\varrho^{-1}$
соответственно. Имеем гомеоморфизм
$M_{\varrho}:\mathrm{H}^{\sigma}(\Omega)\leftrightarrow\varrho\mathrm{H}^{\sigma}(\Omega)$.
Отсюда и из плотности множества $C^{\infty}_{0}(\Omega)$ в
$\mathrm{H}^{\sigma}(\Omega)$ следует, что $C^{\infty}_{0}(\Omega)$ плотно в
$X^{\sigma}(\Omega):=\varrho\mathrm{H}^{\sigma}(\Omega)$. Значит, более широко
множество $X^{\infty}(\Omega)$ плотно в $X^{\sigma}(\Omega)$.

Определим пространство со скалярным произведением:
\begin{gather*}
\varrho^{-1} H^{-\sigma}_{0}(\Omega):=\{f=\varrho^{-1}v:\,v\in
H^{-\sigma}_{0}(\Omega)\,\}, \\
(f_{1},f_{2})_{\varrho^{-1} H^{-\sigma}_{0}(\Omega)}:=(\varrho f_{1},\varrho
f_{2})_{H^{-\sigma}(\Omega)}.
\end{gather*}
Имеем гомеоморфизм
\begin{equation}\label{4.207}
M_{{\varrho}^{-1}}:H^{-\sigma}_{0}(\Omega)\leftrightarrow\varrho^{-1}H^{-\sigma}_{0}(\Omega).
\end{equation}
Следовательно, пространство $\varrho^{-1}H^{-\sigma}_{0}(\Omega)$ полно и множество
$C^{\infty}_{0}(\Omega)$ плотно в нем.

Отметим, что
\begin{equation}\label{4.208}
(\varrho^{-1}H^{-\sigma}_{0}(\Omega))'=\varrho\mathrm{H}^{\sigma}(\Omega)\quad\mbox{с
равенством норм}.
\end{equation}
В самом деле, перейдя в \eqref{4.207} к сопряженному оператору, получим гомеоморфизм
$$
M_{{\varrho}^{-1}}:(\varrho^{-1}H^{-\sigma}_{0}(\Omega))'\leftrightarrow
(H^{-\sigma}_{0}(\Omega))'=\mathrm{H}^{\sigma}(\Omega).
$$
Отсюда и из определения пространства $\varrho\mathrm{H}^{\sigma}(\Omega)$ следует
справедливость еще одного гомеоморфизма
$$
I=M_{\varrho}M_{{\varrho}^{-1}}:(\varrho^{-1}H^{-\sigma}_{0}(\Omega))'\leftrightarrow\varrho
\mathrm{H}^{\sigma}(\Omega),
$$
где $I$ --- тождественный оператор. Тем самым доказано соотношение \eqref{4.208},
поскольку использованные нами гомеоморфизмы являются изометрическими.

Теперь мы можем завершить наше доказательство следующими рассуждениями. В силу
леммы~\ref{lem4.8}, условие~II$_{\sigma}$ равносильно ограниченности оператора
\eqref{4.205}, что ввиду \eqref{4.207} эквивалентно непрерывному вложению
$H^{-\sigma}(\Omega)\hookrightarrow \varrho^{-1}H^{-\sigma}_{0}(\Omega)$. Это
вложение плотное. В силу \eqref{4.208}, оно равносильно плотному непрерывному
вложению
$$
\varrho\mathrm{H}^{\sigma}(\Omega)=(\varrho^{-1}H^{-\sigma}_{0}(\Omega))'
\hookrightarrow(H^{-\sigma}(\Omega))'=
H^{\sigma}_{\overline{\Omega}}(\mathbb{R}^{n}).
$$
Наконец, непрерывное вложение $\varrho\mathrm{H}^{\sigma}(\Omega)\hookrightarrow
H^{\sigma}_{\overline{\Omega}}(\mathbb{R}^{n})$ эквивалентно условию~I$_{\sigma}$.
Отметим, что последнее вложение плотно, поскольку множество $C^{\infty}_{0}(\Omega)$
плотно в пространстве $H^{\sigma}_{\overline{\Omega}}(\mathbb{R}^{n})$. Таким
образом, доказана эквивалентность условий II$_{\sigma}$ и~I$_{\sigma}$ для
пространства $X^{\sigma}(\Omega)=\varrho\mathrm{H}^{\sigma}(\Omega)$.

Теорема \ref{th4.29} доказана.

\medskip

Из теорем \ref{th4.26} и \ref{th4.29} вытекает следующая индивидуальная теорема о
разрешимости краевой задачи \eqref{4.1} в весовых пространствах Соболева.


\begin{theorem}\label{th4.30}
Пусть $\sigma<-1/2$ и положительная функция $\varrho\in C^{\infty}(\Omega)$
удовлетворяет условию~$\mathrm{II}_{\sigma}$. Тогда отображение
$u\rightarrow(Lu,Bu)$, где $u\in C^{\infty}(\,\overline{\Omega}\,)$,
$Lu\in\varrho\mathrm{H}^{\sigma}(\Omega)$, продолжается по непрерывности до
ограниченного нетерового оператора
\begin{gather}\notag
(L,B):\,\bigl\{u\in\mathrm{H}^{\sigma+2q}(\Omega):
Lu\in\varrho\mathrm{H}^{\sigma}(\Omega)\bigr\}\rightarrow\\
\varrho\mathrm{H}^{\sigma}(\Omega)\oplus\bigoplus_{j=1}^{q}\,H^{\sigma+2q-m_{j}-1/2}(\Gamma),
\label{4.209}
\end{gather}
область определения которого является гильбертовым пространством относительно нормы
графика
$$
\bigl(\,\|u\|_{\mathrm{H}^{\sigma+2q}(\Omega)}^{2}+
\|\varrho^{-1}Lu\|_{\mathrm{H}^{\sigma}(\Omega)}^{2}\bigr)^{1/2}.
$$
Индекс оператора \eqref{4.209} равен $\dim N-\dim N^{+}$ и не зависит от $\sigma$ и
$\varrho$.
\end{theorem}


\textbf{Доказательство.} Ограниченность и нетеровость оператора \eqref{4.209} сразу
следуют из теорем \ref{th4.26} и \ref{th4.29}, где
$X^{\sigma}(\Omega):=\varrho\mathrm{H}^{\sigma}(\Omega)$. Кроме того, так как
множество $\mathcal{O}(X^{\infty}(\Omega))$, содержащее $C^{\infty}_{0}(\Omega)$,
плотно в пространстве $H^{\sigma}_{\overline{\Omega}}(\mathbb{R}^{n})$, то в силу
теоремы~\ref{th4.26} (iv) индекс оператора \eqref{4.209} равен $\dim N-\dim N^{+}$ и
не зависит от $\sigma$ и $\varrho$. Теорема~\ref{th4.30} доказана.

\medskip

Важный пример положительной функции $\varrho\in C^{\infty}(\Omega)$, удовлетворяющей
условию II$_{\sigma}$ дает следующая теорема.


\begin{theorem}\label{th4.31}
Пусть заданы число $\sigma<-1/2$ и функция $\varrho_{1}$, удовлетворяющая условию
\eqref{4.187}. Предположим, что $\delta\geq-\sigma-1/2\in\mathbb{Z}$ либо
$\delta>-\sigma-1/2\notin\mathbb{Z}$. Тогда функция $\varrho:=\varrho_{1}^{\delta}$
удовлетворяет условию~$\mathrm{II}_{\sigma}$.
\end{theorem}


Отсюда следует, что теорема LM$_{2}$ для неполуцелых $\sigma<-1/2$ является частным
случаем индивидуальной теоремы~\ref{th4.30}. Напомним, что в теореме LM$_{2}$
использована весовая функция $\varrho:=\varrho_{1}^{-\sigma}$.

\textbf{Доказательство теоремы \ref{th4.31}.} Для функции
$\varrho=\varrho_{1}^{\delta}$ условие \eqref{4.204} выполняется, поскольку
$\varrho_{1}=0$ на $\Gamma$ и $\delta\geq-\sigma-1/2$. Поэтому нам остается
доказать, что $\varrho_{1}^{\delta}$~--- мультипликатор в пространстве
$\mathrm{H}^{-\sigma}(\Omega)=H^{-\sigma}(\Omega)$. Если положительное число
$\delta$ целое, то функция $\varrho_{1}^{\delta}$ принадлежит пространству
$C^{\infty}(\overline{\Omega})$ и поэтому является мультипликатором в
$H^{-\sigma}(\Omega)$. Предположим далее, что $\delta\notin\mathbb{Z}$. Тогда, по
условию, $\delta>-\sigma-1/2$.

Нетрудно проверить, что функция $\eta_{\delta}(t):=\nobreak t^{\delta}$, $0\leq
t\leq1$, принадлежит пространству $H^{-\sigma}((0,1))$ (это мы сделаем в следующем
абзаце). Тогда она имеет некоторое продолжение на $\mathbb{R}$, принадлежащее
$H^{-\sigma}(\mathbb{R})$. Сохраним за ним обозначение $\eta_{\delta}$. По теореме
Р.~Стрихартца [\ref{Strichartz67}; \ref{MazyaShaposhnikova86}, с.~90], всякая
функция из пространства $H^{-\sigma}(\mathbb{R})$ является мультипликатором в нем,
если $-\sigma>1/2$. Значит, $\eta_{\delta}$~--- мультипликатор в
$H^{-\sigma}(\mathbb{R})$.  Тогда функция
$\eta_{\delta,n}(t',t_{n}):=\eta_{\delta}(t_{n})$ аргументов
$t'\in\mathbb{R}^{n-1}$, $t_{n}\in\mathbb{R}$ является мультипликатором в
$H^{-\sigma}(\mathbb{R}^{n})$ [\ref{MazyaShaposhnikova86}, с.~110] (предложение~5).
Она совпадает с $\varrho_{1}^{\delta}$ в специальных локальных координатах
$(x',t_{n})$ вблизи границы $\Gamma$. Здесь $x'$~--- координаты точки в локальной
карте на поверхности $\Gamma$, а $t_{n}$~--- расстояние до $\Gamma$. Отсюда следует
, что $\varrho_{1}^{\delta}$~--- мультипликатор в каждом пространстве
$H^{-\sigma}(\Omega\cap V_{j})$, где $\{V_{j}:j=1,\ldots,r\}$~--- конечная система
шаров в $\mathbb{R}^{n}$ достаточно малого радиуса $\varepsilon$, покрывающая
границу $\Gamma$ [\ref{MazyaShaposhnikova86}, с.~339] (лемма~3). Добавив к этой
системе множество $V_{0}:=\{x\in\Omega:\mathrm{dist}(x,\Gamma)>\varepsilon/2\}$,
получим конечное открытое покрытие замкнутой области $\overline{\Omega}$. Пусть
функции $\chi_{j}\in C^{\infty}_{0}(V_{j})$, $j=0,1,\ldots,r$, образуют разбиение
единицы на $\overline{\Omega}$, подчиненное этому покрытию. Поскольку умножение на
функцию из класса $C^{\infty}_{0}(V_{j})$ является ограниченным оператором в
пространстве $H^{-\sigma}(\Omega\cap V_{j})$, то $\chi_{j}\varrho_{1}^{\delta}$
--- мультипликатор в этом пространстве. Следовательно, функция
$\varrho_{1}^{\delta}=\sum_{j=0}^{r}\chi_{j}\varrho_{1}^{\delta}$ является
мультипликатором в $H^{-\sigma}(\Omega)$.

Нам остается показать, что $\eta_{\delta}\in H^{-\sigma}((0,1))$. Воспользуемся
внутренним описанием пространства $H^{-\sigma}((0,1))$. Если $-\sigma\in\mathbb{Z}$,
то включение $\eta_{\delta}\in H^{-\sigma}((0,1))$ равносильно тому, что
$\eta_{\delta}\in L_{2}((0,1))$ и $\eta_{\delta}^{(-\sigma)}\in L_{2}((0,1))$.
Последние два включения очевидно выполняются при $\delta>-\sigma-1/2$. Значит,
$\eta_{\delta}\in H^{-\sigma}((0,1))$ в рассмотренном случае. Если
$-\sigma\notin\mathbb{Z}$, то включение $\eta_{\delta}\in H^{-\sigma}((0,1))$
равносильно тому, что $\eta_{\delta}\in H^{[-\sigma]}((0,1))$ и
\begin{equation}\label{4.210}
\int_{0}^{1}\int_{0}^{1}\frac{|D_{t}^{[-\sigma]}t^{\delta}-
D_{\tau}^{[-\sigma]}\tau^{\delta}|^{2}}
{|t-\tau|^{1+2\{-\sigma\}}}\,dt\,d\tau<\infty
\end{equation}
(см., например, [\ref{Adams75}, с.~214], теорема 7.48). Здесь, как обычно,
$[-\sigma]$ и $\{-\sigma\}$~--- соответственно целая и дробная части числа
$-\sigma$. Так как $\delta>[-\sigma]-1/2$, то $\eta_{\delta}\in
H^{[-\sigma]}((0,1))$ по доказанному. Кроме того, неравенство \eqref{4.210} верно в
силу следующей элементарной леммы, которую мы докажем в конце этого пункта.


\begin{lemma}\label{lem4.9}
Пусть $\alpha,\beta,\gamma\in\mathbb{R}$, причем $\alpha\neq0$ и $\gamma>0$. Тогда
\begin{equation}\label{4.211}
I(\alpha,\beta,\gamma):=\int_{0}^{1}\int_{0}^{1}\frac{|t^{\alpha}-\tau^{\alpha}|^{\gamma}}
{|t-\tau|^{\beta}}\,dt\,d\tau<\infty
\end{equation}
тогда и только тогда, когда выполняются все три неравенства
\begin{equation}\label{4.212}
\alpha\gamma-\beta>-2,\quad\gamma-\beta>-1,\quad\alpha\gamma>-1.
\end{equation}
\end{lemma}


В самом деле, двойной интеграл в \eqref{4.210} равен $c\,I(\alpha,\beta,\gamma)$,
где $c$~--- некоторое положительное число, а $\alpha=\delta-[-\sigma]$,
$\beta=1+2\{-\sigma\}$ и $\gamma=2$. Для чисел $\alpha$, $\beta$ и $\gamma$
неравенства \eqref{4.212} выполняются:
\begin{gather*}
\alpha\gamma-\beta=2(\delta+\sigma)-1>-2,\quad\gamma-\beta=1-2\{-\sigma\}>-1,\\
\alpha\gamma=2(\delta-[-\sigma])>-1
\end{gather*}
В первом и третьем неравенствах мы используем условие $\delta>-\sigma-1/2$. Таким
образом, $\eta_{\delta}\in H^{-\sigma}((0,1))$ и в случае нецелого $\sigma-1/2$.

Теорема \ref{th4.31} доказана.

\medskip

Нам остается доказать лемму \ref{lem4.9}

\medskip

\textbf{Доказательство леммы \ref{lem4.9}.} Сделав замену $\lambda:=\tau/t$ во
внутреннем интеграле, запишем вследствие очевидных преобразований:
\begin{gather*}
I(\alpha,\beta,\gamma)=\\
2\int_{0}^{1}dt\int_{0}^{t}\frac{|t^{\alpha}-\tau^{\alpha}|^{\gamma}}
{|t-\tau|^{\beta}}\,d\tau=
2\int_{0}^{1}t^{\alpha\gamma-\beta+1}dt\int_{0}^{1}\frac{|1-\lambda^{\alpha}|^{\gamma}}
{|1-\lambda|^{\beta}}\,d\lambda.
\end{gather*}
Здесь интеграл по переменной $t$ конечен тогда и только тогда, когда
$\alpha\gamma-\beta>-2$, а интеграл по $\tau$ конечен тогда и только тогда, когда
$\alpha\gamma>-1$ и $\gamma-\beta>-1$. Следовательно,
$\eqref{4.211}\Leftrightarrow\eqref{4.212}$. Лемма \ref{lem4.9} доказана.





\markright{\emph \ref{sec4.5}. Индивидуальные теоремы о разрешимости}

\section[Индивидуальные теоремы о разрешимости краевой задачи и пространства
Хермандера]{Индивидуальные теоремы \\ о разрешимости краевой задачи \\ и
пространства Хермандера}\label{sec4.5}

\markright{\emph \ref{sec4.5}. Индивидуальные теоремы о разрешимости}

Установим аналоги индивидуальных теорем \ref{th4.26}, \ref{th4.28} и \ref{th4.30}
для уточненной шкалы.

\subsection[Ключевая индивидуальная теорема и пространства Хермандера]
{Ключевая индивидуальная теорема \\ и пространства Хермандера}\label{sec4.5.1}

Пусть $\sigma<0$ и $\varphi\in\mathcal{M}$. Предположим, что задано гильбертово
пространство $X^{\sigma,\varphi}(\Omega)$, непрерывно вложенное в
$\mathcal{D}'(\Omega)$. Рассмотрим следующий аналог условия I$_{\sigma}$.

\medskip

\textbf{Условие I$_{\sigma,\varphi}$.} Множество
$X^{\infty}(\Omega)=X^{\sigma,\varphi}(\Omega)\cap
C^{\infty}(\,\overline{\Omega}\,)$ плотно в $X^{\sigma,\varphi}(\Omega)$ и
существует число $c>0$ такое, что
\begin{equation}\label{4.213}
\|\mathcal{O}f\|_{H^{\sigma,\varphi}(\mathbb{R}^{n})}\leq
c\,\|f\|_{X^{\sigma,\varphi}(\Omega)}\quad\mbox{для любого}\quad f\in
X^{\infty}(\Omega).
\end{equation}


Положим
$$
D^{\sigma+2q,\varphi}_{L,X}(\Omega):=\{u\in H^{\sigma+2q,\varphi}(\Omega):\,Lu\in
X^{\sigma,\varphi}(\Omega)\}.
$$
В~пространстве $D^{\sigma+2q,\varphi}_{L,X}(\Omega)$ вводим скалярное произведение
графика
$$
(u_{1},u_{2})_{D^{\sigma+2q,\varphi}_{L,X}(\Omega)}:=
(u_{1},u_{2})_{H^{\sigma+2q,\varphi}(\Omega)}+(Lu_{1},Lu_{2})_{X^{\sigma,\varphi}(\Omega)}.
$$
Это пространство полное, что доказывается аналогично случаю $\varphi\equiv1$ (см.
п.~\ref{sec4.4.1}).

Сформулируем ключевую индивидуальную теорему.


\begin{theorem}\label{th4.32}
Пусть $\varphi\in\mathcal{M}$, число $\sigma<0$ такое, что
\begin{equation}\label{4.214}
\sigma+2q\neq1/2-k\quad\mbox{для любого}\quad k\in\mathbb{N},
\end{equation}
а $X^{\sigma,\varphi}(\Omega)$~--- произвольное гильбертово пространство, непрерывно
вложенное в $\mathcal{D}'(\Omega)$ и удовлетворяющее
условию~$\mathrm{I}_{\sigma,\varphi}$. Тогда справедливы следующие утверждения.
\begin{itemize}
\item[$\mathrm{(i)}$] Множество
$$
D^{\infty}_{L,X}(\Omega):=\{u\in C^{\infty}(\,\overline{\Omega}\,):\,Lu\in
X^{\sigma,\varphi}(\Omega)\}
$$
плотно в пространстве $D^{\sigma+2q,\varphi}_{L,X}(\Omega)$.
\item[$\mathrm{(ii)}$] Отображение $u\rightarrow(Lu,Bu)$,
где $u\in D^{\infty}_{L,X}(\Omega)$, продолжается по непрерывности до ограниченного
линейного оператора
\begin{gather}\notag
(L,B):\,D^{\sigma+2q,\varphi}_{L,X}(\Omega)\rightarrow \\
X^{\sigma,\varphi}(\Omega)\oplus\bigoplus_{j=1}^{q}\,H^{\sigma+2q-m_{j}-1/2,\varphi}(\Gamma)
=:\mathbf{X}_{\sigma,\varphi}(\Omega,\Gamma).\label{4.215}
\end{gather}
\item[$\mathrm{(iii)}$] Оператор \eqref{4.215} нетеров. Его ядро совпадает с $N$, а
область значений состоит из всех векторов
$(f,g_{1},\ldots,g_{q})\in\mathbf{X}_{\sigma,\varphi}(\Omega,\Gamma)$,
удовлетворяющих условию \eqref{4.4}.
\item[$\mathrm{(iv)}$] Если множество
$\mathcal{O}(X^{\infty}(\Omega))$ плотно в пространстве
$H^{\sigma,\varphi}_{\overline{\Omega}}(\mathbb{R}^{n})$, то индекс оператора
\eqref{4.215} равен $\dim N-\dim N^{+}$.
\end{itemize}
\end{theorem}


Эта теорема доказывается аналогично теореме \ref{th4.26}. При этом необходимо
использовать теоремы \ref{th4.15}, \ref{th4.19}, \ref{th4.24} и \ref{th4.25} для
произвольного $\varphi\in\mathcal{M}$. Кроме того, вместо формул \eqref{4.176} и
\eqref{4.177} следует воспользоваться теоремой~\ref{th3.21}. Дополнительное условие
\eqref{4.214}, отсутствовавшее в теореме \ref{th4.26}, вызвано тем, что в теоремах
\ref{th3.21} и \ref{th4.25} параметр $s=\sigma+2q$ неполуцелый.

\subsection{Другие индивидуальные теоремы}\label{sec4.5.2}

Укажем важные приложения теоремы \ref{th4.32}, обусловленные конкретным выбором
пространства $X^{\sigma,\varphi}(\Omega)$. Пусть $\sigma<0$ и
$\varphi\in\mathcal{M}$

Пространство $X^{\sigma,\varphi}(\Omega):=\{0\}$ удовлетворяет условию
I$_{\sigma,\varphi}$. В~этом случае теорема~\ref{th4.32} описывает свойства
полуоднородной краевой задачи \eqref{4.1}, где $f=0$, и справедлива при любом
$s\in\mathbb{R}$, что установлено в теореме~\ref{th3.11}.

Пространство $X^{\sigma,\varphi}(\Omega):=L_{2}(\Omega)$ также удовлетворяет условию
I$_{\sigma,\varphi}$. Как мы упоминали в п.~\ref{sec4.4.1}, этот выбор пространства
$X^{\sigma,\varphi}(\Omega)$ важен в спектральной теории эллиптических операторов.

В качестве $X^{\sigma,\varphi}(\Omega)$ допустим выбор некоторых пространств
Хермандера. Ввиду теоремы \ref{th4.18} здесь содержательным является случай
$\sigma<-1/2$. Всякое пространство Хермандера
$X^{\sigma,\varphi}(\Omega):=H^{\lambda,\eta}(\Omega)$, где $\lambda>-1/2$ и
$\eta\in\mathcal{M}$, удовлетворяет условию I$_{\sigma,\varphi}$ при $\sigma<-1/2$.
Действительно, в силу теоремы \ref{th3.9} (i), (iv) пространство
$H^{\lambda,\eta}(\Omega)=H^{\lambda,\eta,(0)}(\Omega)$ непрерывно вложено в
$H^{\sigma,\varphi,(0)}(\Omega)=H^{\sigma,\varphi}_{\overline{\Omega}}(\mathbb{R}^{n})$.
Отсюда следует неравенство \eqref{4.213}, в котором $X^{\infty}(\Omega)=
C^{\infty}(\,\overline{\Omega}\,)$. Полагая
$X^{\sigma,\varphi}(\Omega):=H^{\lambda,\eta}(\Omega)$ в ключевой теореме
\ref{th4.32}, получаем следующую индивидуальную теорему о разрешимости краевой
задачи~\ref{4.1} в пространствах Хермандера.


\begin{theorem}\label{th4.33}
Пусть число $\sigma<-1/2$ удовлетворяет условию \eqref{4.214}, $\lambda>-1/2$ и
$\varphi,\eta\in\mathcal{M}$. Тогда отображение $u\mapsto(Lu,Bu)$, где $u\in
C^{\infty}(\,\overline{\Omega}\,)$, продолжается по непрерывности до ограниченного
нетерового оператора
\begin{gather}\notag
(L,B):\,\{u\in H^{\sigma+2q,\varphi}(\Omega):Lu\in
H^{\lambda,\eta}(\Omega)\}\rightarrow\\ \label{4.216}
H^{\lambda,\eta}(\Omega)\oplus\bigoplus_{j=1}^{q}\,H^{\sigma+2q-m_{j}-1/2,\varphi}(\Gamma),
\end{gather}
область определения которого является гильбертовым пространством относительно нормы
графика
$$
\bigl(\,\|u\|_{H^{\sigma+2q,\varphi}(\Omega)}^{2}+
\|Lu\|_{H^{\lambda,\eta}(\Omega)}^{2}\bigr)^{1/2}.
$$
Индекс оператора \eqref{4.216} равен $\dim\ N-\dim N^{+}$ и не зависит от параметров
$\sigma$, $\varphi$ и $\lambda$, $\eta$.
\end{theorem}


Отметим, что в теореме~\ref{th4.33} решение и правая часть эллиптического уравнения
$Lu=f$ могут иметь различные дополнительные гладкости $\varphi$ и~$\eta$.

В теореме~\ref{th4.33} пространство
$X^{\sigma,\varphi}(\Omega):=H^{\lambda,\eta}(\Omega)$ лежит в $H^{-1/2}(\Omega)$,
поскольку $\lambda>-1/2$. Пространство $X^{\sigma,\varphi}(\Omega)$, содержащее
широкий класс распределений $f\notin H^{-1/2}(\Omega)$ и удовлетворяющее условию
I$_{\sigma,\varphi}$, можно получить, используя некоторые весовые пространства
Хермандера.

Пусть $\sigma<-1/2$, $\varphi\in\mathcal{M}$, а функция $\varrho\in
C^{\infty}(\Omega)$ положительная. Определим пространство со скалярным произведением
\begin{gather*}
\varrho H^{\sigma,\varphi}(\Omega):=\{f=\varrho v:\,v\in
H^{\sigma,\varphi}(\Omega)\,\}, \\
(f_{1},f_{2})_{\varrho H^{\sigma,\varphi}(\Omega)}:=
(\varrho^{-1}f_{1},\varrho^{-1}f_{2})_{H^{\sigma,\varphi}(\Omega)}.
\end{gather*}
Умножение на функцию $\varrho$ определяет гомеоморфизм
$$
M_{\varrho}:H^{\sigma,\varphi}(\Omega)\leftrightarrow\varrho
H^{\sigma,\varphi}(\Omega).
$$
Следовательно, пространство $\varrho H^{\sigma,\varphi}(\Omega)$ полно (гильбертово)
и непрерывно вложено в $\mathcal{D}'(\Omega)$.

Если функция $\varrho$ является мультипликатором в пространстве
$H^{-\sigma,1/\varphi}(\Omega)$ и выполняется \eqref{4.204}, то пространство
$X^{\sigma,\varphi}(\Omega):=\varrho H^{\sigma,\varphi}(\Omega)$ удовлетворяет
условию I$_{\sigma,\varphi}$. Это доказывается аналогично теореме~\ref{th4.29}.
Отличие состоит лишь в том, что вместо формул \eqref{4.206} и \eqref{4.175}
необходимо использовать теорему \ref{th3.20} (i), (iii). Полагая
$X^{\sigma,\varphi}(\Omega):=\varrho H^{\sigma,\varphi}(\Omega)$ в ключевой теореме
\ref{th4.32}, получаем следующую индивидуальную теорему о разрешимости краевой
задачи~\eqref{4.1} в пространствах, связанных с весовыми пространствами Хермандера.


\begin{theorem}\label{th4.34}
Пусть число $\sigma<-1/2$ удовлетворяет условию \eqref{4.214},
$\varphi\in\mathcal{M}$, а функция $\varrho\in C^{\infty}(\Omega)$ положительна.
Предположим, что $\varrho$ является мультипликатором в пространстве
$H^{-\sigma,1/\varphi}(\Omega)$ и удовлетворяет условию \eqref{4.204}. Тогда
отображение $u\rightarrow(Lu,Bu)$, где $u\in C^{\infty}(\,\overline{\Omega}\,)$,
$Lu\in\varrho H^{\sigma,\varphi}(\Omega)$, продолжается по непрерывности до
ограниченного нетерового оператора
\begin{gather*}
(L,B):\,\bigl\{u\in H^{\sigma+2q,\varphi}(\Omega): Lu\in\varrho
H^{\sigma,\varphi}(\Omega)\bigr\}\rightarrow\\
\varrho H^{\sigma,\varphi}(\Omega)
\oplus\bigoplus_{j=1}^{q}\,H^{\sigma+2q-m_{j}-1/2,\varphi}(\Gamma),
\end{gather*}
область определения которого является гильбертовым пространством относительно нормы
графика
$$
\bigl(\,\|u\|_{H^{\sigma+2q,\varphi}(\Omega)}^{2}+
\|\varrho^{-1}Lu\|_{H^{\sigma,\varphi}(\Omega)}^{2}\bigr)^{1/2}.
$$
Индекс этого оператора равен $\dim N-\dim N^{+}$ и не зависит от $\sigma$, $\varphi$
и $\varrho$.
\end{theorem}


Важный пример функции $\varrho$, удовлетворяющей условию теоремы~\ref{th4.34},
получается, если положить $\varrho:=\varrho_{1}^{\delta}$ для произвольного
фиксированного числа $\delta>-\sigma-1/2$. Напомним, что $\varrho_{1}$~--- функция
расстояния до границы $\Gamma$, подчиненная условию \eqref{4.187}. В~самом деле,
функция $\varrho\in C^{\infty}(\Omega)$ положительна и удовлетворяет условию
\eqref{4.204}. В~силу теоремы \ref{th4.31}, эта функция является мультипликатором в
пространстве $H^{-\sigma\mp\varepsilon}(\Omega)$, где число $\varepsilon>0$ такое,
что $\delta>-\sigma\mp\varepsilon-1/2$. Отсюда на основании интерполяционной теоремы
\ref{th3.2}, следует, что $\varrho$~--- мультипликатор в пространстве
$H^{-\sigma,1/\varphi}(\Omega)$. Таким образом, функция
$\varrho=\varrho_{1}^{\delta}$ удовлетворяет условию теоремы~\ref{th4.34}.




\markright{\emph \ref{notes4}. Примечания и комментарии}

\section{Примечания и комментарии}\label{notes4}

\markright{\emph \ref{notes4}. Примечания и комментарии}

\small

\textbf{К п. 4.1.} В фундаментальной работе С.~Агмона, А.~Дуглиса и Л.~Ниренберга
[\ref{AgmonDouglisNirenberg62}] доказаны априорные оценки решений общей
эллиптической краевой задачи в соответствующих парах пространств Гельдера (нецелых
порядков) и позитивных пространств Соболева; при этом задача предполагается заданной
в ограниченной евклидовой области с гладкой границей. Для позитивных пространств
Соболева эти оценки доказаны независимо также Ф.~Е.~Браудером [\ref{Browder59}],
Л.~Н.~Слободецким [\ref{Slobodetsky58a}, \ref{Slobodetsky58b}, \ref{Slobodetsky60}],
М.~Шехтером [\ref{Schechter60a}] и другими (см. литературу, приведенную в обзоре
М.~С.~Аграновича [\ref{Agranovich97}]). Априорные оценки эквивалентны нетеровости
оператора, порожденного задачей, в указанных парах пространств. Если краевая задача
регулярная эллиптическая, то дефектное подпространство и область значений оператора
можно описать с помощью дифференциальных выражений, связанных с формально
сопряженной задачей. Доказано также, что наличие априорных оценок решений в парах
гильбертовых пространств Соболева влечет выполнение условия Лопатинского для набора
граничных выражений, т.~е. эллиптичность краевой задачи. Эти вопросы систематически
изложены, например, в монографиях Ю.~М.~Березанского [\ref{Berezansky65}],
Ж.-Л.~Лионса и Э.~Мадженеса [\ref{LionsMagenes71}], О.~И.~Панича [\ref{Panich86}],
Х.~Трибеля [\ref{Triebel80}], Л.~Хермандера [\ref{Hermander65}, \ref{Hermander87}],
М.~Шехтера [\ref{Schechter77}], в обзоре М.~С.~Аграновича [\ref{Agranovich97}].
Индекс эллиптической краевой задачи вычислен М.~Ф.~Атьёй и Р.~Боттом
[\ref{AtiyahBott64}], при этом использована фундаментальная формула индекса
эллиптического матричного ПДО, установленная М.~Ф.~Атьёй и И.~М.~Зингером
[\ref{AtiyahSinger63}]. Подробное изложение теории индекса эллиптических краевых
задач дано в монографии Ш.~Ремпеля и Б.-В.~Шульце [\ref{RempelSchulze86}].

Теоремы о разрешимости регулярных эллиптических краевых задач доказаны и для иных
шкал позитивных функциональных пространств: Г.~Шлензак [\ref{Shlenzak74}] для
некоторого класса гильбертовых пространств Хермандера, а Х.~Трибелем
[\ref{Triebel80}, \ref{Triebel86}] и Й.~Франке [\ref{Franke85}] для пространств
Лизоркина--Трибеля и Никольского--Бесова, как банаховых так и квазибанаховых.

Теоремы о разрешимости эллиптических краевых задач имеют различные и важные
приложения. Среди них отметим теоремы о повышении гладкости решения эллиптического
уравнения вплоть до границы области, приложения к исследованию функции Грина
эллиптической краевой задачи, к задачам оптимального управления, к нелокальным
эллиптическим краевым задачам, к некоторым классам нелинейных эллиптических краевых
задач и другие. См. монографии Ю.~М.~Березанского [\ref{Berezansky65}],
Ю.~М.~Березанского, Г.~Ф.~Уса и З.~Г.~Шефтеля [\ref{BerezanskyUsSheftel90}],
О.~А.~Ладыженской и Н.~Н.~Уральцевой [\ref{LadyszenskayaUraltseva64}], Ж.-Л.~Лионса
[\ref{Lions72a}, \ref{Lions72b}], Ж.-Л.~Лионса и Э.~Мадженеса
[\ref{LionsMagenes71}], Я.~А.~Ройтберга [\ref{Roitberg96}, \ref{Roitberg99}],
И.~В.~Скрыпника [\ref{Skrypnik73}], обзор М.~С.~Аграновича [\ref{Agranovich97}] и
приведенную там литературу.

В работах Ш.~Агмона и Л.~Ниренберга [\ref{Agmon62}, \ref{AgmonNirenberg63}],
М.~С.~Аграновича и М.~И.~Вишика [\ref{AgranovichVishik64}] выделен подкласс
эллиптических краевых задач~--- эллиптические задачи с параметром, обладающие
следующим важным свойством. При достаточно больших по модулю значениях комплексного
параметра оператор, соответствующий задаче, является гомеоморфизмом в подходящих
парах соболевских пространств, причем норма оператора допускает двустороннюю оценку
с постоянными, не зависящими от параметра. Индекс такой задачи равен нулю для всех
значений параметра. Эллиптические краевые задачи с параметром имеют важные
приложения в теории параболических задач и в спектральной теории дифференциальных
операторов. Различные более широкие классы эллиптических операторов и эллиптических
краевых задач с параметром изучены в работах М.~С.~Аграновича [\ref{Agranovich90b},
\ref{Agranovich92}], Г.~Грубб [\ref{Grubb96}] (гл.~2), Р.~Денка, Р.~Менникена и
Л.~Р.~Волевича [\ref{DenkMennickenVolevich98}, \ref{DenkMennickenVolevich01},
\ref{DenkVolevich02}], А.~Н.~Кожевникова [\ref{Kozhevnikov73}, \ref{Kozhevnikov96},
\ref{Kozhevnikov97}], О.~И.~Панича [\ref{Panich66}, \ref{Panich73}]; см. также
обзоры М.~С.~Аграновича [\ref{Agranovich90}] (\S~4) и [\ref{Agranovich97}] (\S~3,
п.~6.4) и приведенную там литературу.

Все теоремы п. 4.1 за исключением последней доказаны нами в [\ref{06UMJ3},
\ref{07UMJ5}], причем в более общей ситуации, когда эллиптическая краевая задача
задана на гладком компактном многообразии с краем. Последняя теорема 4.10
установлена в [\ref{07Dop4}].

\medskip

\textbf{К п. 4.2.} В работах Я.~А.~Ройтберга [\ref{Roitberg64}, \ref{Roitberg65},
\ref{Roitberg96}] и Ж.-Л.~Лионса, Э.~Мадженеса [\ref{LionsMagenes71},
\ref{Magenes66}, \ref{LionsMagenes62}, \ref{LionsMagenes63}] исследована
разрешимость эллиптической краевой задачи в двусторонних шкалах видоизмененных
соболевских пространств. Ими предложены принципиально различные способы построения
области определения оператора, соответствующего задаче, которые приводят к различным
типам теорем о разрешимости~--- общей и индивидуальным теоремам. В общей теореме
Ройтберга область определения оператора не зависит от коэффициентов эллиптического
уравнения и является общей для всех краевых задач одного порядка. В~индивидуальных
теоремах Лионса--Мадженеса она зависит от коэффициентов эллиптического уравнения.
Термины „общая” и „индивидуальная” теорема о разрешимости предложены в
[\ref{09OperatorTheory191}, \ref{09MFAT2}].

В общей теореме Ройтберга используются введенные им [\ref{Roitberg64},
\ref{Roitberg65}] двусторонняя шкала видоизмененных соболевских пространств и
понятие обобщенного решения в этих пространствах. Я.~А.~Ройтберг называет их
пространствами „типа соболевских”. Мы используем термин „модифицированная по
Ройтбергу соболевская шкала” для класса этих пространств.

Общая теорема о разрешимости, доказанная Я.~А.~Ройтбергом [\ref{Roitberg64},
\ref{Roitberg65}] для регулярных эллиптических краевых задач, была затем установлена
им для нерегулярных эллиптических краевых задач [\ref{Roitberg69}, \ref{Roitberg70}]
и краевых задач для общих эллиптических систем [\ref{Roitberg75}]. Эти результаты
известны в литературе как „теоремы о полном наборе гомеоморфизмов (или
изоморфизмов)”. Им посвящена монография Я.~А.~Ройтберга [\ref{Roitberg96}], их
изложение имеется также в книге Ю.~М.~Березанского [\ref{Berezansky65}, с.~248],
справочнике [\ref{FunctionalAnalysis72}, с.~170] и обзоре М.~С.~Аграновича
[\ref{Agranovich97}, с.~81]. В наиболее общей форме теорема о полном наборе
гомеоморфизмов доказана А.~Н.~Кожевниковым [\ref{Kozhevnikov01}] для общих
эллиптических псевдодифференциальных краевых задач.

Я.~А.~Ройтберг, З.~Г.~Шефтель и их ученики систематически применяли теоремы о полном
наборе гомеоморфизмов и модифицированную соболевскую шкалу к изучению различных
классов эллиптических краевых задач: задачи со степенными особенностями в правых
частях, задачи с сильным вырождением на границе, задачи трансмиссии, нелокальные
задачи, задача Соболева и другие; см. монографии Я.~А.~Ройтберга [\ref{Roitberg96},
\ref{Roitberg99}] и приведенную там литературу. Я.~А.~Ройтберг [\ref{Roitberg68}]
вывел из теоремы о полном наборе гомеоморфизмов ряд других утверждений о
разрешимости регулярной эллиптической краевой задачи. Они были распространены на
нерегулярные эллиптические краевые задачи Я.~А.~Ройтбергом и Ю.~В.~Костарчуком
[\ref{KostarchukRoitberg73}], а также на эллиптические краевые задачи для систем
И.~Я.~Ройтберг и Я.~А.~Ройтбергом [\ref{RoitbergRoitberg00}] (см. монографию
Я.~А.~Ройтберга [\ref{Roitberg99}], п.~1.3). А.~А.~Мурачем [\ref{94Dop12},
\ref{94UMJ12}] доказаны теоремы о полном наборе гомеоморфизмов для двусторонних шкал
пространств Лизоркина-Трибеля и Никольского-Бесова, модифицированных по Ройтбергу.

Концепция модифицированной по Ройтбергу соболевской шкалы и обобщенного решения в
ней использована В.~А.~Козловым, В.~Г.~Мазьей и Й.~Россманном
[\ref{KozlovMazyaRossmann97}] в теории эллиптических краевых задач в областях с
негладкой границей, Н.~В.~Житарашу и С.~Д.~Эйдельманом [\ref{ZhitarashuEidelman92}]
в теории параболических уравнений, Я.~А.~Ройтбергом [\ref{Roitberg99}] в теории
гиперболических уравнений.

Все теоремы п.~4.2 доказаны нами в [\ref{08UMJ4}], за исключением последней теоремы
4.22, публикуемой впервые.

\medskip

\textbf{К п. 4.3.} Основной результат этого пункта составляют теоремы 4.6 и 4.7.
Первая из них вводит эквивалентную норму в пространствах уточненной шкалы, связанную
с эллиптическим выражением, а вторая утверждает существование следов у обобщенного
решения эллиптического уравнения. В соболевском случае эти теоремы доказаны
Я.~А.~Ройтбергом [\ref{Roitberg71}], [\ref{Roitberg96}, с. 174, 187]. Для
пространств Хермандера они установлены нами в монографии впервые. При этом
используются леммы 4.6 и 4.7 о проекторах в некоторых соболевских пространствах,
представляющие самостоятельный интерес.

\medskip

\textbf{К п. 4.4.} Различные индивидуальные теоремы о разрешимости регулярной
эллиптической краевой в шкалах соболевских пространств, содержащих нерегулярные
распределения, установлены Ж.-Л.~Лионсом и Э.~Мадженесом в [\ref{LionsMagenes71},
\ref{Magenes66}, \ref{LionsMagenes62}, \ref{LionsMagenes63}]. В этих теоремах
областью определения оператора, соответствующего задаче, служат пространства
ограничений эллиптического оператора, построенные по заданным пространствам правых
частей эллиптического уравнения. В частности, она может совпадать с областью
определения максимального оператора, отвечающего эллиптическому дифференциальному
выражению. Как показал Л.~Хермандер [\ref{Hermander59}, с.~78], последняя
существенно зависит от всех коэффициентов эллиптического выражения, даже если они
постоянны. Эта индивидуальная теорема используется в спектральной теории
эллиптических операторов; см. работы Г.~Грубб [\ref{Grubb68}, \ref{Grubb71},
\ref{Grubb96}] и В.~А.~Михайлеца [\ref{Mikhailets82}, \ref{Mikhailets89},
\ref{Mikhailets90}].

Предложенное нами условие I$_{\sigma}$ на пространство правых частей эллиптического
уравнения является достаточно общим. Ему удовлетворяют пространства, использованные
Ж.-Л.~Лионсом и Э.~Мадженесом. Я.~А.~Ройтберг [\ref{Roitberg71}] (п.~2.4)
использовал условие более сильное, чем наше. При таком условии Ройтберг
[\ref{Roitberg71}, с.~261; \ref{Roitberg96}, с.~190] доказал ограниченность
оператора, соответствующего эллиптической краевой задаче, для всех значений
$\sigma<0$.

Все результаты п. 4.4.2 -- 4.4.4, в частности индивидуальные теоремы 4.25, 4.27 и
4.29, доказаны в [\ref{09MFAT2}]. Они содержат в себе, как частный случай,
индивидуальные теоремы Лионса-Мадженеса, сформулированные в п.~4.4.1. К ключевой
индивидуальной теореме 4.25 близки по духу (выбор широкого класса пространств правых
частей) две индивидуальные теоремы, сформулированные в обзорах Э.~Мадженеса
[\ref{Magenes66}, с.~202] и М.~С.~Аграновича [\ref{Agranovich97}, с.~85] (первая из
них принадлежит Ж.-Л.~Лионсу и Э.~Мадженесу, а вторая~--- Я.~А.~Ройтбергу).

Иные индивидуальные теоремы доказаны в работах Я.~А.~Ройтберга [\ref{Roitberg68},
с.~545] и Ю.~В.~Костарчука, Я.~А.~Ройтберга [\ref{KostarchukRoitberg73}, с.~275].
В~этих теоремах решения и/или правые части эллиптического уравнения не являются
распределениями, заданными в евклидовой области, что отличает их от индивидуальных
теорем п.~4.4.

\medskip

\textbf{К п. 4.5.} Результаты этого пункта анонсированы нами в
[\ref{09OperatorTheory191}] (п.~6). Они распространяют индивидуальные теоремы,
доказанные в п.~4.4, на широкие классы пространств, связанных с пространствами
Хермандера.

\normalsize



\chapter{\textbf{Эллиптические системы}}\label{ch5}

\chaptermark{\emph Гл. \ref{ch5}. Эллиптические системы}



\section[Равномерно эллиптические системы в уточненной шкале]
{Равномерно эллиптические \\ системы в уточненной шкале}\label{sec2.4}

\markright{\emph\ref{sec2.4}. Равномерно эллиптические системы в уточненной шкале}

Здесь мы изучим равномерно эллиптические системы псевдодифференциальных уравнений,
заданных в евклидовом пространстве. Мы установим априорную оценку их решений в
уточненной шкале пространств и исследуем внутреннюю гладкость решений. Эти
результаты полезно сравнить с теоремами п.~\ref{sec2.4b}, где исследован случай
одного уравнения. Их доказательства проводятся по той же схеме, что и в скалярном
случае. Для полноты изложения мы приводим их целиком.

\subsection{Равномерно эллиптические системы}\label{sec2.4.1}

В  пространстве $\mathbb{R}^{n}$ рассматривается линейная система
псевдодифференциальных уравнений
\begin{equation}\label{2.42}
\sum_{k=1}^{p}\:A_{j,k}\,u_{k}=f_{j},\quad j=1,\ldots,p.
\end{equation}
Здесь $n,p\in\mathbb{N}$, а $A_{j,k}\in\Psi^{\infty}(\mathbb{R}^{n})$, где
$j,k=1,\ldots,p$,~--- скалярные полиоднородные псевдодифференциальные операторы в
$\mathbb{R}^{n}$. Решения уравнений \eqref{2.42} рассматриваются в классе
распределений в $\mathbb{R}^{n}$.

Важным примером системы \eqref{2.42} является система линейных дифференциальных
уравнений с коэффициентами из $C^{\infty}_{\mathrm{b}}(\mathbb{R}^{n})$.

Мы исследуем следующий класс систем псевдодифференциальных уравнений
[\ref{Agranovich90}, с.~51].


\begin{definition}\label{def2.12}
Cистема \eqref{2.42} называется \emph{равномерно эллиптической} в $\mathbb{R}^{n}$
по Дуглису-Ниренбергу, если существуют наборы вещественных чисел
$\{l_{1},\ldots,l_{p}\}$ и $\{m_{1},\ldots,m_{p}\}$ такие, что:
\begin{itemize}
\item [(i)] $\mathrm{ord}\,A_{j,k}\leq l_{j}+m_{k}$ для всех
$j,k=1,\ldots,p$;
\item [(ii)] найдется число $c>0$ такое, что
$$
|\det\bigl(a^{(0)}_{j,k}(x,\xi)\bigr)_{j,k=1}^{p}|\geq c\quad\mbox{для любых}\;\;
x,\xi\in\mathbb{R}^{n},\;\;|\xi|=1,
$$
где $a^{(0)}_{j,k}(x,\xi)$ --- главный символ ПДО $A_{j,k}$ в случае
$\mathrm{ord}\,A_{j,k}=l_{j}+m_{k}$, либо $a^{(0)}_{j,k}(x,\xi)\equiv0$ в случае
$\mathrm{ord}\,A_{j,k}<l_{j}+m_{k}$.
\end{itemize}
\end{definition}


Это определение задает достаточно широкий класс эллиптических систем ПДО. В случае
дифференциальных уравнений он был введен А.~Дуглисом и Л.~Ниренбергом
[\ref{DouglisNirenberg55}] (ими рассматривалась эллиптичность в евклидовых
областях). Этот класс содержит в себе однородные эллиптические системы, для которых
все $l_{j}=0$ и $m_{k}=m\in\mathbb{R}$, а также~--- эллиптические по
И.~Г.~Петровскому системы [\ref{Petrovsky39}], для которых все $l_{j}=0$, а $m_{k}$
могут быть различными.

Далее мы предполагаем, что система \eqref{2.42} удовлетворяет определению
\ref{def2.12}.

Запишем эту систему в матричной форме: $Au=f$. Здесь $A:=(A_{j,k})_{j,k=1}^{p}$~---
матричный ПДО, а $u=\mathrm{col}\,(u_{1},\ldots,u_{p})$,
$f=\mathrm{col}\,(f_{1},\ldots,f_{p})$~--- функциональные столбцы. Как и сама
система, матричный ПДО $A$ называется равномерно эллиптическим в $\mathbb{R}^{n}$.
Для него существует параметрикс $B$, т. е. справедливо следующее утверждение
[\ref{Agranovich90}, с.~52].


\begin{proposition}\label{prop2.6}
Существует матричный классический ПДО $B=(B_{k,j})_{k,j=1}^{p}$ такой, что
$B_{k,j}\in\Psi^{-m_{k}-l_{j}}(\mathbb{R}^{n})$ и
\begin{equation}\label{2.46}
BA=I+T_{1},\quad AB=I+T_{2},
\end{equation}
где $T_{1}=(T_{1}^{j,k})_{j,k=1}^{p}$ и $T_{2}=(T_{2}^{k,j})_{k,j=1}^{p}$ ---
некоторые матричные ПДО, элементы которых принадлежат классу
$\Psi^{-\infty}(\mathbb{R}^{n})$, а $I$ --- тождественный оператор в
$\mathcal{S}'(\mathbb{R}^{n})$.
\end{proposition}

\subsection{Априорная оценка решения системы}\label{sec2.4.2}

Установим априорную оценку решения системы $Au=f$ в подходящих парах пространств
Хермандера. Поскольку $A_{j,k}\in\Psi^{l_{j}+m_{k}}(\mathbb{R}^{n})$, согласно
лемме~\ref{lem2.8} имеем линейный ограниченный оператор
\begin{equation}\label{2.45}
A:\,\bigoplus_{k=1}^{p}\,H^{s+m_{k},\,\varphi}(\mathbb{R}^{n})\rightarrow
\bigoplus_{j=1}^{p}\,H^{s-l_{j},\,\varphi}(\mathbb{R}^{n})
\end{equation}
для любых $s\in\mathbb{R}$ и $\varphi\in\mathcal{M}$.


\begin{theorem}\label{th2.16}
Пусть $s\in\mathbb{R}$, $\sigma>0$ и $\varphi\in\mathcal{M}$. Существует число
$c=c(s,\sigma,\varphi)>0$ такое, что для произвольных вектор-функций
\begin{gather}\label{2.47}
u=\mathrm{col}\,(u_{1},\ldots,u_{p})
\in\bigoplus_{k=1}^{p}\,H^{s+m_{k},\,\varphi}(\mathbb{R}^{n}), \\
f=\mathrm{col}\,(f_{1},\ldots,f_{p})\
\in\bigoplus_{j=1}^{p}\,H^{s-l_{j},\,\varphi}(\mathbb{R}^{n}),\label{2.48}
\end{gather}
удовлетворяющих уравнению $Au=f$ в $\mathbb{R}^{n}$, справедлива априорная оценка
\begin{gather}\notag
\Bigl(\sum_{k=1}^{p}\|u_{k}\|_{H^{s+m_{k},\,\varphi}(\mathbb{R}^{n})}^{2}\Bigr)^{1/2}\leq\\
c\,\Bigl(\sum_{j=1}^{p}\|f_{j}\|_{H^{s-l_{j},\,\varphi}(\mathbb{R}^{n})}^{2}\Bigr)^{1/2}+
c\,\Bigl(\sum_{k=1}^{p}\|u_{k}\|_{H^{s+m_{k}-\sigma,\,\varphi}(\mathbb{R}^{n})}^{2}\Bigr)^{1/2}.
\label{2.49}
\end{gather}
\end{theorem}


\textbf{Доказательство.} Обозначим через $\|\cdot\|'_{s,\varphi}$,
$\|\cdot\|''_{s,\varphi}$ и $\|\cdot\|'_{s-\sigma,\varphi}$ соответственно нормы в
пространствах
$$
\bigoplus_{k=1}^{p}\,H^{s+m_{k},\varphi}(\mathbb{R}^{n}),\;\;
\bigoplus_{j=1}^{p}\,H^{s-l_{j},\varphi}(\mathbb{R}^{n})\;\;\mbox{и}\;\;
\bigoplus_{k=1}^{p}\,H^{s+m_{k}-\sigma,\varphi}(\mathbb{R}^{n}).
$$
Пусть вектор-функции \eqref{2.47}, \eqref{2.48} удовлетворяют уравнению $Au=f$ в
$\mathbb{R}^{n}$. В силу первого равенства в \eqref{2.46} запишем $u=Bf-T_{1}u$.
Отсюда следует оценка \eqref{2.49}:
$$
\|u\|'_{s,\varphi}=\|Bf-T_{1}u\|'_{s,\varphi}\leq\|Bf\|'_{s,\varphi}+
\|T_{1}u\|'_{s,\varphi}\leq c\,\|f\|''_{s,\varphi}+c\,\|u\|'_{s-\sigma,\varphi}.
$$
Здесь $c$ --- максимум норм операторов
\begin{gather}\label{2.50}
B:\,\bigoplus_{j=1}^{p}\,H^{s-l_{j},\,\varphi}(\mathbb{R}^{n})\rightarrow
\bigoplus_{k=1}^{p}\,H^{s+m_{k},\,\varphi}(\mathbb{R}^{n}),\\
T_{1}:\,\bigoplus_{k=1}^{p}\,H^{s+m_{k}-\sigma,\,\varphi}(\mathbb{R}^{n})
\rightarrow\bigoplus_{k=1}^{p}\,H^{s+m_{k},\,\varphi}(\mathbb{R}^{n}).\label{2.51}
\end{gather}
Эти операторы ограниченные в силу леммы \ref{lem2.8} и предложения \ref{prop2.6}.
Теорема~\ref{th2.16} доказана.

\medskip

Теорема \ref{2.16} уточняет применительно к шкале пространств
$H^{s,\varphi}(\mathbb{R}^{n})$ априорную оценку, полученную Л.~Хермандером
[\ref{Hermander67}, с.~175] в соболевской шкале.

\subsection{Гладкость решения системы}\label{sec2.4.3}

Предположим, что правая часть уравнения  $Au=f$ имеет некоторую внутреннюю гладкость
в уточненной шкале на заданном открытом непустом множестве
$V\subseteq\mathbb{R}^{n}$. Изучим внутреннюю гладкости решения $u$ на этом
множестве. Рассмотрим сначала случай, когда $V=\mathbb{R}^{n}$. Напомним, что
$H^{-\infty}(\mathbb{R}^{n})$ объединение всех пространств
$H^{s,\varphi}(\mathbb{R}^{n})$, где $s\in\mathbb{R}$, $\varphi\in\mathcal{M}$.


\begin{theorem}\label{th2.17}
Предположим, что $u\in(H^{-\infty}(\mathbb{R}^{n}))^{p}$ является решением уравнения
$Au=f$ в $\mathbb{R}^{n}$, где
$$
f_{j}\in H^{s-l_{j},\,\varphi}(\mathbb{R}^{n})\quad\mbox{при}\quad j=1,\ldots,p
$$
для некоторых параметров $s\in\mathbb{R}$ и $\varphi\in\mathcal{M}$. Тогда
$$
u_{k}\in H^{s+m_{k},\,\varphi}(\mathbb{R}^{n})\quad\mbox{при}\quad k=1,\ldots,p.
$$
\end{theorem}


\textbf{Доказательство.} В силу теоремы \ref{th2.15} (i), для вектор-функции
$u\in(H^{-\infty}(\mathbb{R}^{n}))^{p}$ существует число $\sigma>0$ такое, что
\begin{equation}\label{2.52}
u\in\bigoplus_{k=1}^{p}\,H^{s+m_{k}-\sigma,\,\varphi}(\mathbb{R}^{n}).
\end{equation}
Отсюда и из условия теоремы получаем на основании формул \eqref{2.46}, \eqref{2.50},
\eqref{2.51} требуемое свойство:
$$
u=BAu-T_{1}u=Bf-T_{1}u
\in\bigoplus_{k=1}^{p}\,H^{s+m_{k},\,\varphi}(\mathbb{R}^{n}).
$$
Теорема~\ref{th2.17} доказана.

\medskip

Рассмотрим теперь общий случай, когда $V$ --- произвольное открытое непустое
подмножество пространства $\mathbb{R}^{n}$. Напомним, что пространство
$H^{\sigma,\varphi}_{\mathrm{int}}(V)$ распределений, имеющих заданную внутреннюю
гладкость на $V$ в уточненной шкале, определено в п.~\ref{sec2.4.3b}.


\begin{theorem}\label{th2.18}
Предположим, что $u\in(H^{-\infty}(\mathbb{R}^{n}))^{p}$ является решением уравнения
$Au=f$ на множестве $V$, где
\begin{equation}\label{2.53}
f_{j}\in H^{s-l_{j},\,\varphi}_{\mathrm{int}}(V)\quad\mbox{при}\quad j=1,\ldots,p
\end{equation}
для некоторых параметров $s\in\mathbb{R}$ и $\varphi\in\mathcal{M}$. Тогда
\begin{equation}\label{2.54}
u_{k}\in H^{s+m_{k},\,\varphi}_{\mathrm{int}}(V)\quad\mbox{при}\quad k=1,\ldots,p.
\end{equation}
\end{theorem}


\textbf{Доказательство.} Докажем сначала, что из условия \eqref{2.53} вытекает
следующее свойство повышения внутренней гладкости решения уравнения $Au=f$: для
каждого числа $r\geq1$ справедлива импликация
\begin{equation}\label{2.55}
u\in\bigoplus_{k=1}^{p}\,H^{s-r+m_{k},\,\varphi}_{\mathrm{int}}(V)\;\Rightarrow\;
u\in\bigoplus_{k=1}^{p}\,H^{s-r+1+m_{k},\,\varphi}_{\mathrm{int}}(V).
\end{equation}

Произвольно выберем функцию $\chi\in C^{\infty}_{\mathrm{b}}(\mathbb{R}^{n})$ такую,
что
\begin{equation}\label{2.56}
\mathrm{supp}\,\chi\subset V\quad\mbox{и}\quad
\mathrm{dist}(\mathrm{supp}\,\chi,\partial V)>0.
\end{equation}
Для нее существует функция $\eta\in C^{\infty}_{\mathrm{b}}(\mathbb{R}^{n})$ такая,
что
\begin{equation}\label{2.57}
\mathrm{supp}\,\eta\subset V,\;\;\mathrm{dist}(\mathrm{supp}\,\eta,\,\partial
V)>0,\;\;\eta=1\;\mbox{в окрестности}\;\mathrm{supp}\,\chi
\end{equation}
(это показано в доказательстве теоремы \ref{th2.18b}).

Переставив матричный ПДО $A$ и оператор умножения на функцию $\chi$, запишем
\begin{gather}\notag
A\chi u=A\chi\eta u=\chi\,A\eta u+A'\eta u=\\ \notag
\chi\,Au+\chi\,A(\eta-1)u+A'\eta u=\\
=\chi f+\chi\,A(\eta-1)u+A'\eta u\quad\mbox{в}\quad\mathbb{R}^{n}. \label{2.58}
\end{gather}
Здесь матричный ПДО $A'=(\,A'_{j,k}\,)_{j,k=1}^{p}$~--- коммутатор ПДО $A$ и
оператора умножения на функцию $\chi$. Поскольку
$A'_{j,k}\in\Psi^{l_{j}+m_{k}-1}(\mathbb{R}^{n})$, то в силу леммы~\ref{lem2.8}
имеем ограниченный оператор
$$
A':\,\bigoplus_{k=1}^{p}\,H^{s-r+m_{k},\,\varphi}(\mathbb{R}^{n})\,\rightarrow\,
\bigoplus_{j=1}^{p}\,H^{s-r+1-l_{j},\,\varphi}(\mathbb{R}^{n}).
$$
Следовательно,
\begin{equation}\label{2.59}
u\in\bigoplus_{k=1}^{p}\,H^{s-r+m_{k},\,\varphi}_{\mathrm{int}}(V)\;\Rightarrow\;
A'\eta u\in\bigoplus_{j=1}^{p}\,H^{s-r+1-l_{j},\,\varphi}(\mathbb{R}^{n}).
\end{equation}
Далее, согласно условию \eqref{2.53} и ввиду неравенства $r\geq1$ имеем
\begin{equation}\label{2.60}
\chi f\in\bigoplus_{j=1}^{p}\,H^{s-l_{j},\varphi}(\mathbb{R}^{n})\hookrightarrow
\bigoplus_{j=1}^{p}\,H^{s-r+1-l_{j},\varphi}(\mathbb{R}^{n}).
\end{equation}
Кроме того, так как носители функций $\chi$ и $\eta-1$ не пересекаются, то ПДО
$$
\chi A_{j,k}(\eta-1)\in\Psi^{-\infty}(\mathbb{R}^{n})
$$
для всех $j,k=1,\ldots,p$. Отсюда, поскольку для вектор-функции
$u\in(H^{-\infty}(\mathbb{R}^{n}))^{p}$ справедливо \eqref{2.52} при некотором
$\sigma>0$, мы получаем в силу леммы~\ref{lem2.8} включение
\begin{equation}\label{2.61}
\chi\,A(\eta-1)u\in\bigoplus_{j=1}^{p}\,H^{s-r+1-l_{j},\varphi}(\mathbb{R}^{n}).
\end{equation}

На основании формул \eqref{2.58} -- \eqref{2.61} и теоремы~\ref{th2.17} получаем,
что
\begin{gather*}
u\in\bigoplus_{k=1}^{p}\,H^{s-r+m_{k},\,\varphi}_{\mathrm{int}}(V)\;\Rightarrow\;
A\chi u\in\bigoplus_{j=1}^{p}\,H^{s-r+1-l_{j},\varphi}(\mathbb{R}^{n})\;\Rightarrow\\
\chi u\in\bigoplus_{k=1}^{p}\,H^{s-r+1+m_{k},\,\varphi}(\mathbb{R}^{n}).
\end{gather*}
Тем самым доказана импликация \eqref{2.55} ввиду произвольности выбора функции
$\chi\in C^{\infty}_{\mathrm{b}}(\mathbb{R}^{n})$, удовлетворяющей условию
\eqref{2.56}.

Теперь с помощью \eqref{2.55} легко вывести свойство \eqref{2.54}. Можно считать,
что в формуле \eqref{2.52} число $\sigma>0$ целое. Значит,
$$
u\in\bigoplus_{k=1}^{p}\,H^{s-\sigma+m_{k},\,\varphi}_{\mathrm{int}}(V).
$$
Применив \eqref{2.55} последовательно для
$r=\sigma,\,\sigma\nobreak-\nobreak1,\ldots,1$, выводим свойство \eqref{2.54}:
\begin{gather*}
u\in\bigoplus_{k=1}^{p}\,H^{s-\sigma+m_{k},\,\varphi}_{\mathrm{int}}(V)
\;\Rightarrow\;
u\in\bigoplus_{k=1}^{p}\,H^{s-\sigma+1+m_{k},\,\varphi}_{\mathrm{int}}(V)\;
\Rightarrow
\\ \ldots\;\Rightarrow\;
u\in\bigoplus_{k=1}^{p}\,H^{s+m_{k},\,\varphi}_{\mathrm{int}}(V).
\end{gather*}

Теорема \ref{th2.18} доказана.

\medskip

Теорема~\ref{th2.18} позволяет установить наличие непрерывных производных у
выбранной компоненты $u_{k}$ решения системы \eqref{2.42}. При этом используется
теорема~\ref{th2.15} (iii).


\begin{theorem}\label{th2.19}
Пусть заданы целые числа $k\in\{1,\ldots,p\}$, $r\geq0$ и функция
$\varphi\in\mathcal{M}$, удовлетворяющая условию \eqref{2.37}. Предположим, что
$u\in(H^{-\infty}(\mathbb{R}^{n}))^{p}$ является решением уравнения $Au=f$ на
открытом множестве $V\subseteq\mathbb{R}^{n}$, где
\begin{equation}\label{2.62}
f_{j}\in H^{r-m_{k}-l_{j}+n/2,\,\varphi}_{\mathrm{int}}(V)\quad\mbox{для всех}\quad
j=1,\ldots,p.
\end{equation}
Тогда компонента $u_{k}$ решения имеет на множестве $V$ непрерывные частные
производные до порядка $r$ включительно, причем эти производные ограничены на каждом
множестве $V_{0}\subset V$ таком, что $\mathrm{dist}(V_{0},\partial V)>0$.
В~частности, если $V=\mathbb{R}^{n}$, то $u_{k}\in
C^{\,r}_{\mathrm{b}}(\mathbb{R}^{n})$.
\end{theorem}


\textbf{Доказательство.} В силу теоремы~\ref{th2.18}, где полагаем $s:=r-m_{k}+n/2$,
справедливо включение $u_{k}\in H^{r+n/2,\,\varphi}_{\mathrm{int}}(V)$. Пусть
функция $\eta\in C^{\infty}_{\mathrm{b}}(\mathbb{R}^{n})$ удовлетворяет условиям
$$
\mathrm{supp}\,\eta\subset V,\;\;\mathrm{dist}(\mathrm{supp}\,\eta,\,\partial
V)>0,\;\;\eta=1\;\mbox{в окрестности}\;V_{0}.
$$
Для распределения $\eta u_{k}$ в силу теоремы~\ref{th2.15} (iii) имеем
$$
\eta u_{k}\in H^{r+n/2,\,\varphi}(\mathbb{R}^{n})\hookrightarrow
C^{\,r}_{\mathrm{b}}(\mathbb{R}^{n}).
$$
Отсюда вытекает, что все частные производные функции $u_{k}$ до порядка $r$
включительно непрерывны и ограничены в некоторой окрестности множества $V_{0}$.
Тогда эти производные непрерывны и на множестве $V$, поскольку можно взять
$V_{0}:=\{x_{0}\}$ для любой точки $x_{0}\in V$. Теорема~\ref{th2.19} доказана.


В качестве приложения теоремы~\ref{th2.19} мы дадим следующее достаточное условие
классичности решения системы $Au=f$ в случае, когда все $A_{j,k}$~---  линейные
дифференциальные операторы с коэффициентами класса
$C^{\infty}_{\mathrm{b}}(\mathbb{R}^{n})$, и все $l_{j}=0$. Иными словами, когда
$Au=f$ является равномерно эллиптической в $\mathbb{R}^{n}$ по Петровскому системой
дифференциальных уравнений.


\begin{corollary}
Предположим, что $u\in(H^{-\infty}(\mathbb{R}^{n}))^{p}$ является решением уравнения
$Au=f$ на открытом множестве $V\subseteq\mathbb{R}^{n}$, где $f_{j}\in
H^{n/2,\varphi}_{\mathrm{int}}(V)$ для каждого номера $j=1,\ldots,p$ и некоторого
параметра $\varphi\in\mathcal{M}$, удовлетворяющего условию \eqref{2.37}. Тогда
решение $u$ классическое на множестве $V$, т.~е. $u_{k}\in C^{\,m_{k}}(V)$ для всех
$k=1,\ldots,p$.
\end{corollary}


Это сразу следует из теоремы~\ref{th2.19} при $r=m_{k}$. Отметим, что для
классического решения $u$ системы \eqref{2.42} ее левые части вычисляются с помощью
классических производных (а не обобщенных, как в теории распределений), и эти
производные непрерывны на множестве $V$.




\markright{\emph \ref{sec2.6}. Эллиптические системы на замкнутом многообразии}

\section[Эллиптические системы на замкнутом многообра-\break зии]
{Эллиптические системы \\ на замкнутом многообразии}\label{sec2.6}

\markright{\emph \ref{sec2.6}. Эллиптические системы на замкнутом многообразии}

Здесь мы изучим эллиптические по Дуглису-Ниренбергу системы псевдодифференциальных
уравнений, заданных на бесконечно гладком замкнутом (компактном) ориентированном
многообразии $\Gamma$. Мы докажем, что оператор, соответствующий этим системам,
ограничен и нетеров в подходящих парах пространств уточненной шкалы. Кроме того,
будет исследован один класс эллиптических систем (эллиптические системы с
параметром), для которых указанный оператор является гомеоморфизмом при больших по
модулю значениях комплексного параметра. Результаты, которые мы получим, полезно
сравнить с теоремами п.~\ref{sec2.6b}, где исследован случай одного уравнения. Их
доказательства проводятся по той же схеме, что и в скалярном случае.

\subsection{Эллиптические системы}\label{sec2.6.1}

Рассмотрим систему линейных уравнений
\begin{equation}\label{2.106}
\sum_{k=1}^{p}\:A_{j,k}\,u_{k}=f_{j}\;\;\mbox{на}\;\;\Gamma,\;\;\mbox{где}\;\;
j=1,\ldots,p.
\end{equation}
Здесь $p\in\mathbb{N}$, а $A_{j,k}\in\Psi^{\infty}(\Gamma)$, где
$j,k=1,\ldots,p$,~--- скалярные классические ПДО, заданные на многообразии $\Gamma$.
Уравнения \eqref{2.106} понимаются в смысле теории распределений.


\begin{definition}\label{def2.17}
Cистема \eqref{2.106} называется \emph{эллиптической} на $\Gamma$ по
Дуглису-Ниренбергу, если существуют наборы вещественных чисел
$\{l_{1},\ldots,l_{p}\}$ и $\{m_{1},\ldots,m_{p}\}$ такие, что:
\begin{itemize}
\item [(i)] $\mathrm{ord}\,A_{j,k}\leq{l_{j}+m_{k}}$
для всех $j,k=1,\ldots,p$;
\item [(ii)] $\det\bigl(\,a^{(0)}_{j,k}(x,\xi)\,\bigr)_{j,k=1}^{p}\neq0$
для произвольных точки $x\in\Gamma$ и ковектора $\xi\in
T^{\ast}_{x}\Gamma\setminus\{0\}$; здесь $a^{(0)}_{j,k}(x,\xi)$ --- главный символ
ПДО $A_{j,k}$ в случае $\mathrm{ord}\,A_{j,k}=l_{j}+m_{k}$, либо
$a^{(0)}_{j,k}(x,\xi)\equiv0$ в случае $\mathrm{ord}\,A_{j,k}<l_{j}+m_{k}$.
\end{itemize}
\end{definition}


Далее в п.~\ref{sec2.6} предполагается, что система \eqref{2.106} удовлетворяет
определению~\ref{def2.17}.

Запишем эту систему в матричной форме: $Au=f$ на $\Gamma$. Здесь
$A:=(A_{j,k})_{j,k=1}^{p}$~--- матричный ПДО на $\Gamma$, а
$u=\mathrm{col}\,(u_{1},\ldots,u_{p})$, $f=\mathrm{col}\,(f_{1},\ldots,f_{p})$~---
функциональные столбцы. Как и сама система, соответствующий ей матричный ПДО $A$
называется эллиптическим на $\Gamma$.

Обозначим через $A^{+}=(A_{j,k}^{+})_{j,k=1}^{p}$ матричный ПДО, формально
сопряженный к оператору $A$ относительно $C^{\infty}$-плотности $dx$ на $\Gamma$;
здесь каждый ПДО $A_{j,k}^{+}$ формально сопряженный к $A_{j,k}$. Эллиптичность
системы $Au=f$ равносильна эллиптичности сопряженной системы $A^{+}v=g$ (по Дуглису
-- Ниренбергу). Положим
\begin{gather*}
\mathcal{N}:=\{\,u\in
(C^{\infty}(\Gamma))^{p}:\,Au=0\;\;\mbox{на}\;\;\Gamma\,\},\\
\mathcal{N}^{+}:=\{\,v\in
(C^{\infty}(\Gamma))^{p}:\,A^{+}v=0\;\;\mbox{на}\;\;\Gamma\,\}.
\end{gather*}
Поскольку системы $Au=f$ и $A^{+}v=g$ эллиптичны, то [\ref{Agranovich90}, с.~52]
пространства $\mathcal{N}$ и $\mathcal{N}^{+}$ конечномерны.

\subsection[Оператор эллиптической системы в уточненной шкале]
{Оператор эллиптической системы \\ в уточненной шкале}\label{sec2.6.2}

Изучим свойства матричного ПДО $A$ в уточненной шкале на многообразии~$\Gamma$.
Поскольку $A_{j,k}\in\Psi^{l_{j}+m_{k}}(\Gamma)$, то в силу леммы \ref{lem2.13} мы
имеем линейный ограниченны оператор
\begin{equation}\label{2.108}
A:\,\bigoplus_{k=1}^{p}\,H^{s+m_{k},\varphi}(\Gamma)\,\rightarrow\,
\bigoplus_{j=1}^{p}\,H^{s-l_{j},\varphi}(\Gamma)
\end{equation}
для любых $s\in\mathbb{R}$ и $\varphi\in\mathcal{M}$. Изучим его свойства.


\begin{theorem}\label{th2.28}
Для произвольных параметров $s\in\mathbb{R}$ и $\varphi\in\nobreak\mathcal{M}$
ограниченный оператор \eqref{2.108} нетеров. Его ядро совпадает с пространством
$\mathcal{N}$, а область значений состоит из всех таких вектор-функций
\begin{equation}\label{2.109}
f=\mathrm{col}\,(f_{1},\ldots,f_{p})\in
\bigoplus_{j=1}^{p}\,H^{s-l_{j},\varphi}(\Gamma),
\end{equation}
которые удовлетворяют условию
\begin{equation}\label{2.110}
\sum_{j=1}^{p}\,(f_{j},w_{j})_{\Gamma}=0\;\;\mbox{для любого}\;\;
w=(w_{1},\ldots,w_{p})\in \mathcal{N}^{+}.
\end{equation}
Индекс оператора \eqref{2.108} равен $\dim\mathcal{N}-\dim\mathcal{N}^{+}$ и не
зависит от $s$, $\varphi$.
\end{theorem}


\textbf{Доказательство.} В случае $\varphi\equiv1$ (соболевская шкала) эта теорема
известна [\ref{Agranovich90}, с.~52]. Отсюда общий случай $\varphi\in\mathcal{M}$
получается с помощью интерполяции с функциональным параметром. А именно, пусть
$s\in\mathbb{R}$. Имеем ограниченные нетеровы операторы
\begin{equation}\label{2.111}
A:\,\bigoplus_{k=1}^{p}\,H^{s\mp1+m_{k}}(\Gamma)\,\rightarrow\,
\bigoplus_{j=1}^{p}\,H^{s\mp1-l_{j}}(\Gamma)
\end{equation}
с общим ядром $\mathcal{N}$ и одинаковым индексом
$\varkappa:=\dim\mathcal{N}-\dim\mathcal{N}^{+}$. При этом
\begin{gather}\notag
A\Bigl(\,\bigoplus_{k=1}^{p}\,H^{s\mp1+m_{k}}(\Gamma)\,\Bigr)=\\
\Bigl\{\,f=\mathrm{col}\,(f_{1},\ldots,f_{p})\in
\bigoplus_{j=1}^{p}\,H^{s\mp1-l_{j}}(\Gamma):\,\mbox{верно \eqref{2.110}}\Bigr\}.
\label{2.112}
\end{gather}
Применим к \eqref{2.111} интерполяцию с функциональным параметром $\psi$ из теоремы
\ref{th2.21}, где $\varepsilon=\delta=1$. Получим ограниченный оператор
\begin{gather*}
A:\,\Bigl[\;\bigoplus_{k=1}^{p}\,H^{s-1+m_{k}}(\Gamma),\,
\bigoplus_{k=1}^{p}\,H^{s+1+m_{k}}(\Gamma)\,\Bigr]_{\psi}\rightarrow\\
\Bigl[\;\bigoplus_{j=1}^{p}\,H^{s-1-l_{j}}(\Gamma),\,
\bigoplus_{j=1}^{p}\,H^{s+1-l_{j}}(\Gamma)\,\Bigr]_{\psi},
\end{gather*}
который в силу теорем \ref{th2.5} и \ref{th2.21} совпадает с оператором
\eqref{2.108}. Следовательно, согласно теореме \ref{th2.7} оператор \eqref{2.108}
нетеров с ядром $\mathcal{N}$ и индексом
$\varkappa=\dim\mathcal{N}-\dim\mathcal{N}^{+}$. Область значений этого оператора
равна
$$
\Bigl(\,\bigoplus_{j=1}^{p}\,H^{s-l_{j},\varphi}(\Gamma)\,\Bigr)\,\cap\,
A\Bigl(\,\bigoplus_{k=1}^{p}\,H^{s-1+m_{k}}(\Gamma)\,\Bigr).
$$
Отсюда, в силу \eqref{2.112} получаем, что она такая как в формулировке настоящей
теоремы. Теорема \ref{th2.28} доказана.

\medskip

Согласно этой теореме, $\mathcal{N}^{+}$~--- дефектное подпространство оператора
\eqref{2.108}. Заметим, что в виду теоремы \ref{th2.22} (v) оператор
\begin{equation}\label{2.113}
A^{+}:\,\bigoplus_{j=1}^{p}\,H^{-s+l_{j},1/\varphi}(\Gamma)\,\rightarrow\,
\bigoplus_{k=1}^{p}\,H^{-s-m_{k},1/\varphi}(\Gamma)
\end{equation}
сопряженный к оператору \eqref{2.108}. Поскольку сопряженная система $A^{+}v=g$
эллиптическая, то в силу теоремы \ref{th2.28} ограниченный оператор \eqref{2.113}
нетеров и имеет ядро $\mathcal{N}^{+}$ и дефектное подпространство $\mathcal{N}$.

Если пространства $\mathcal{N}$ и $\mathcal{N}^{+}$ тривиальны, то из теоремы
\ref{th2.28} и теоремы Банаха об обратном операторе вытекает, что оператор
\eqref{2.108} является гомеоморфизмом. В общем случае гомеоморфизм удобно задавать с
помощью следующих проекторов.

Представим пространства, в которых действует оператор \eqref{2.108}, в виде прямых
сумм замкнутых подпространств:
\begin{gather*}
\bigoplus_{k=1}^{p}\,H^{s+m_{k},\varphi}(\Gamma)=\\ \mathcal{N}\dotplus
\Bigl\{u\in\bigoplus_{k=1}^{p}\,H^{s+m_{k},\varphi}(\Gamma):\;
\sum_{k=1}^{p}\,(u_{k},v_{k})_{\Gamma}=0\;\,\mbox{для
всех}\,\;v\in\mathcal{N}\Bigr\},\\
\bigoplus_{j=1}^{p}\,H^{s-l_{j},\varphi}(\Gamma)=\\ \mathcal{N}^{+}\dotplus
\Bigl\{f\in\bigoplus_{j=1}^{p}\,H^{s-l_{j},\varphi}(\Gamma):\;
\sum_{j=1}^{p}\,(f_{j},w_{j})_{\Gamma}=0\;\,\mbox{для
всех}\,\;w\in\mathcal{N}^{+}\Bigr\}.
\end{gather*}
Здесь, напомним, $u=\mathrm{col}\,(u_{1},\ldots,u_{p})$,
$f=\mathrm{col}\,(f_{1},\ldots,f_{p})$, а также $v=(v_{1},\ldots,v_{p})$,
$w=(w_{1},\ldots,w_{p})$. Указанные разложения в прямые суммы существуют, поскольку
в них слагаемые имеют тривиальное пересечение, и конечная размерность первого из них
равна коразмерности второго. (Действительно, например, в первой сумме
факторпространство пространства $\bigoplus_{k=1}^{p}\,H^{s+m_{k},\varphi}(\Gamma)$
по второму слагаемому является двойственным пространством к подпространству
$\mathcal{N}$ пространства $\bigoplus_{k=1}^{p}\,H^{-s-m_{k},1/\varphi}(\Gamma)$).
Обозначим через $\mathcal{P}$ и $\mathcal{P}^{+}$ косые проекторы соответственно
пространств
$$
\bigoplus_{k=1}^{p}\,H^{s+m_{k},\varphi}(\Gamma)\quad\mbox{и}\quad
\bigoplus_{j=1}^{p}\,H^{s-l_{j},\varphi}(\Gamma)
$$
на вторые слагаемые в указанных суммах параллельно первым слагаемым. Эти проекторы
не зависят от $s$ и $\varphi$.


\begin{theorem}\label{th2.29}
Для произвольных $s\in\mathbb{R}$ и $\varphi\in\mathcal{M}$ сужение оператора
\eqref{2.108} на подпространство
$\mathcal{P}\bigl(\,\bigoplus_{k=1}^{p}\,H^{s+m_{k},\varphi}(\Gamma)\,\bigr)$
является гомеоморфизмом
\begin{equation}\label{2.114}
A:\mathcal{P}\Bigl(\,\bigoplus_{k=1}^{p}\,H^{s+m_{k},\varphi}(\Gamma)\,\Bigr)
\leftrightarrow
\mathcal{P}^{+}\Bigl(\,\bigoplus_{j=1}^{p}\,H^{s-l_{j},\varphi}(\Gamma)\,\Bigr).
\end{equation}
\end{theorem}


\textbf{Доказательство.} По теореме \ref{th2.28}, $\mathcal{N}$ --- ядро, а
$$
\mathcal{P}^{+}\Bigl(\,\bigoplus_{j=1}^{p}\,H^{s-l_{j},\varphi}(\Gamma)\Bigr)
$$
--- область значений оператора \eqref{2.108}. Следовательно, оператор \eqref{2.114}~---
биекция. Кроме того, этот оператор ограничен. Значит, он является гомеоморфизмом в
силу теоремы Банаха об обратном операторе. Теорема~\ref{th2.29} доказана.

\medskip

Из теоремы \ref{th2.29} вытекает следующая априорная оценка решения уравнения
$Au=f$.


\begin{theorem}\label{th2.30}
Пусть $s\in\mathbb{R}$ и $\varphi\in\mathcal{M}$. Предположим, что вектор-функция
\begin{equation}\label{2.115}
u=\mathrm{col}\,(u_{1},\ldots,u_{p})
\in\bigoplus_{k=1}^{p}\,H^{s+m_{k},\varphi}(\Gamma)
\end{equation}
является решением уравнения $Au=f$ на $\Gamma$, правая часть которого удовлетворяет
условию \eqref{2.109}. Тогда для выбранных параметров $s$, $\varphi$ и произвольного
числа $\sigma>0$ существует такое число $c>0$, не зависящее от $u$ и $f$, что
\begin{gather}\notag
\Bigl(\,\sum_{k=1}^{p}\|u_{k}\|_{H^{s+m_{k},\varphi}(\Gamma)}^{2}\Bigr)^{1/2}\leq\\
c\,\Bigl(\,\sum_{j=1}^{p}\|f_{j}\|_{H^{s-l_{j},\varphi}(\Gamma)}^{2}\Bigr)^{1/2}+
c\,\Bigl(\,\sum_{k=1}^{p}\|u_{k}\|_{H^{s+m_{k}-\sigma,\varphi}(\Gamma)}^{2}\Bigr)^{1/2}.
\label{2.116}
\end{gather}
\end{theorem}


\textbf{Доказательство.} Использованные в \eqref{2.116} нормы в пространствах
$$
\bigoplus_{k=1}^{p}\,H^{s+m_{k},\varphi}(\Gamma),\quad
\bigoplus_{j=1}^{p}\,H^{s-l_{j},\varphi}(\Gamma)\quad\mbox{и}\quad
\bigoplus_{k=1}^{p}\,H^{s-\sigma+m_{k},\varphi}(\Gamma)
$$
обозначим для краткости через $\|\cdot\|^{'}_{s,\varphi}$\,,
$\|\cdot\|^{''}_{s,\varphi}$ и $\|\cdot\|^{'}_{s-\sigma,\varphi}$ соответственно.
Так как $\mathcal{N}$~--- конечномерное подпространство в этих пространствах, то
указанные нормы эквивалентны на $\mathcal{N}$. В частности, для вектор-функции
$u-\mathcal{P}u\in\mathcal{N}$ справедливо
$$
\|u-\mathcal{P}u\|^{'}_{s,\varphi}\leq\,
c_{1}\,\|u-\mathcal{P}u\|^{'}_{s-\sigma,\varphi}
$$
с постоянной $c_{1}>0$, не зависящей от $u$. Отсюда получаем
\begin{gather*}
\|u\|^{'}_{s,\varphi}\leq\,
\|u-\mathcal{P}u\|^{'}_{s,\varphi}+\|\mathcal{P}u\|^{'}_{s,\varphi}\leq \\
c_{1}\,\|u-\mathcal{P}u\|^{'}_{s-\sigma,\varphi}+\|\mathcal{P}u\|^{'}_{s,\varphi}\leq\,
c_{1}\,c_{2}\,\|u\|^{'}_{s-\sigma,\varphi}+\|\mathcal{P}u\|^{'}_{s,\varphi},
\end{gather*}
где $c_{2}$ --- норма проектора $1-\mathcal{P}$, действующего в пространстве
$$
\bigoplus_{k=1}^{p}\,H^{s-\sigma+m_{k},\varphi}(\Gamma).
$$
Итак,
\begin{equation}\label{2.117}
\|u\|^{'}_{s,\varphi}\leq\,\|\mathcal{P}u\|^{'}_{s,\varphi}+
c_{1}\,c_{2}\,\|u\|^{'}_{s-\sigma,\varphi}.
\end{equation}
Воспользуемся теперь условием $Au=f$. Поскольку $\mathcal{N}$~--- ядро оператора
\eqref{2.108} и $u-\mathcal{P}u\in\mathcal{N}$, то $A\mathcal{P}u=f$. Таким образом,
$\mathcal{P}u$~--- прообраз вектор-функции $f$ относительно гомеоморфизма
\eqref{2.114}. Следовательно,
$$
\|\mathcal{P}u\|^{'}_{s,\varphi}\leq\, c_{3}\,\|f\|^{''}_{s,\varphi},
$$
где $c_{3}$ --- норма оператора, обратного к \eqref{2.114}. Отсюда и из неравенства
\eqref{2.117} немедленно вытекает оценка \eqref{2.116}. Теорема \ref{th2.30}
доказана.

\medskip

Отметим, что, если $\mathcal{N}=\{0\}$, т. е.  уравнение $Au=f$ имеет не более
одного решения, то величину
$$
\sum_{k=1}^{p}\;\|u_{k}\|_{H^{s-\sigma+m_{k},\varphi}(\Gamma)}
$$
в правой части оценки \eqref{2.116} можно опустить. Если же $\mathcal{N}\neq\{0\}$,
то для каждой вектор-функции $u$ эту величину можно сделать как угодно малой за счет
выбора достаточно большого числа $\sigma$.

\subsection{Локальная гладкость решения}\label{sec2.6.3}

Пусть $\Gamma_{0}$ --- открытое непустое подмножество многообразия $\Gamma$.
Исследуем локальную гладкость решения эллиптического уравнения $Au=f$ на
$\Gamma_{0}$ в уточненной шкале. Рассмотрим сначала случай, когда
$\Gamma_{0}=\Gamma$.


\begin{theorem}\label{th2.31}
Предположим, что вектор-функция $u\in(\mathcal{D}'(\Gamma))^{p}$ является решением
уравнения $Au=f$ на многообразии $\Gamma$, где
\begin{equation}\label{2.118}
f_{j}\in H^{s-l_{j},\varphi}(\Gamma)\quad\mbox{при}\quad j=1,\ldots,p
\end{equation}
для некоторых параметров $s\in\mathbb{R}$ и $\varphi\in\mathcal{M}$. Тогда
\begin{equation}\label{2.119}
u_{k}\in H^{s+m_{k},\varphi}(\Gamma)\quad\mbox{при}\quad k=1,\ldots,p.
\end{equation}
\end{theorem}


\textbf{Доказательство.} Поскольку многообразие $\Gamma$ компактно, то пространство
$\mathcal{D}'(\Gamma)$ является объединением соболевских пространств
$H^{\sigma}(\Gamma)$, где $\sigma\in\mathbb{R}$. Следовательно, для вектор-функции
$u\in(\mathcal{D}'(\Gamma))^{p}$ существует число $\sigma<s$, такое что
$$
u\in\bigoplus_{k=1}^{p}\,H^{\sigma+m_{k}}(\Gamma).
$$
В~силу теорем \ref{th2.28} и \ref{th2.22} (iii) справедливо равенство
$$
\Bigl(\,\bigoplus_{j=1}^{p}\,H^{s-l_{j},\varphi}(\Gamma)\Bigr)\cap
A\Bigl(\,\bigoplus_{k=1}^{p}\,H^{\sigma+m_{k}}(\Gamma)\Bigr)=
A\Bigl(\,\bigoplus_{k=1}^{p}\,H^{s+m_{k},\varphi}(\Gamma)\Bigr).
$$
Поэтому, из условия \eqref{2.118} вытекает, что
$$
f=Au\in A\Bigl(\,\bigoplus_{k=1}^{p}\,H^{s+m_{k},\varphi}(\Gamma)\,\Bigr).
$$
Таким образом, на многообразии $\Gamma$ наряду с равенством $Au=f$ выполняется также
равенство $Av=f$ для некоторой вектор-функции
$$
v\in\bigoplus_{k=1}^{p}\,H^{s+m_{k},\varphi}(\Gamma).
$$
Следовательно, $A(u-v)=0$ на $\Gamma$ и, согласно теореме \ref{th2.28} справедливо
$w:=u-v\in\mathcal{N}$. Однако,
$$
\mathcal{N}\subset(C^{\infty}(\Gamma))^{p}\subset
\bigoplus_{k=1}^{p}\,H^{s+m_{k},\varphi}(\Gamma).
$$
Значит,
$$
u=v+w\in\bigoplus_{k=1}^{p}\,H^{s+m_{k},\varphi}(\Gamma),
$$
т. е. выполняется свойство \eqref{2.119}. Теорема \ref{th2.31} доказана.

\medskip

Рассмотрим теперь случай произвольного $\Gamma_{0}$. Напомним, что пространство
$H^{\sigma,\varphi}_{\mathrm{loc}}(\Gamma_{0})$ распределений, имеющих заданную
локальную гладкость на $\Gamma_{0}$ в уточненной шкале, определено в
п.~\ref{sec2.6.3b}.


\begin{theorem}\label{th2.32}
Предположим, что $u\in(\mathcal{D}'(\Gamma))^{p}$ является решением уравнения $Au=f$
на множестве $\Gamma_{0}$, где
$$
f_{j}\in H^{s-l_{j},\varphi}_{\mathrm{loc}}(\Gamma_{0})\quad\mbox{при}\quad
j=1,\ldots,p
$$
для некоторых параметров $s\in\mathbb{R}$ и $\varphi\in\mathcal{M}$. Тогда
$$
u_{k}\in H^{s+m_{k},\varphi}_{\mathrm{loc}}(\Gamma_{0})\quad\mbox{при}\quad
k=1,\ldots,p.
$$
\end{theorem}


Теорема \ref{th2.32} доказывается аналогично теореме \ref{th2.18}. При этом надо
вместо теоремы \ref{th2.17} использовать теорему \ref{th2.31}.

Из теорем \ref{th2.32} и \ref{th2.27} немедленно вытекает следующее достаточное
условие непрерывности производных у выбранной компоненты $u_{k}$ решения системы
$Au=f$.


\begin{corollary}\label{cor2.4}
Пусть заданы целые числа $k\in\{1,\ldots,p\}$, $r\geq0$ и функция
$\varphi\in\mathcal{M}$, удовлетворяющая неравенству \eqref{2.37}. Предположим, что
$u\in(\mathcal{D}'(\Gamma))^{p}$ является решением уравнения $Au=f$ на множестве
$\Gamma_{0}$, где
$$
f_{j}\in H^{r-m_{k}-l_{j}+n/2,\,\varphi}_{\mathrm{loc}}(\Gamma_{0})\quad\mbox{для
всех}\quad j=1,\ldots,p.
$$
Тогда $u_{k}\in C^{\,r}(\Gamma_{0})$.
\end{corollary}

\subsection{Эллиптическая система с параметром}\label{sec2.6.4}

Следуя обзору [\ref{Agranovich90}, с.~64], мы рассмотрим достаточно широкий класс
эллиптических систем с параметром на многообразии $\Gamma$. Произвольно зафиксируем
числа $p,q\in\mathbb{N}$, $m>0$ и $m_{1},\ldots,m_{p}\in\mathbb{R}$. Рассмотрим
матричный ПДО $A(\lambda)$, зависящий от комплексного параметра $\lambda$ следующим
образом:
\begin{equation}\label{2.120}
A(\lambda):=\,\sum_{r\,=\,0}^{q}\,\lambda^{q-r}\,A^{(r)}.
\end{equation}
Здесь $A^{(r)}:=\bigl(\,A^{(r)}_{j,k}\,\bigr)_{j,k=1}^{p}$ --- квадратная матрица,
образованная скалярными полиоднородными ПДО $A^{(r)}_{j,k}$ на $\Gamma$ такими, что
$$
\mathrm{ord}\,A^{(r)}_{j,k}\leq mr+m_{k}-m_{j}.
$$
При этом будем считать, что $A^{(0)}=-I$, где $I$~--- единичная матрица.

Рассмотрим зависящую от параметра $\lambda\in\mathbb{C}$ систему линейных уравнений
\begin{equation}\label{2.121}
A(\lambda)\,u=f\quad\mbox{на}\;\;\Gamma.
\end{equation}
Здесь, как и прежде, $u=\mathrm{col}\,(u_{1},\ldots,u_{p})$,
$f=\mathrm{col}\,(f_{1},\ldots,f_{p})$ --- функциональные столбцы, компоненты
которых --- распределения на многообразии $\Gamma$.

Пусть $K$ --- фиксированный замкнутый угол на комплексной плоскости с вершиной в
начале координат (не исключается случай, когда $K$ вырождается в луч).


\begin{definition}\label{def2.18}
Система \eqref{2.121} называется \emph{эллиптической с параметром} в угле $K$, если
\begin{equation}\label{2.122}
\det\sum_{r\,=\,0}^{q}\,\lambda^{q-r}\,a^{r,0}(x,\xi)\neq0
\end{equation}
для любых $x\in\Gamma$, $\xi\in T_{x}^{\ast}\Gamma$ и $\lambda\in K$ таких, что
$(\xi,\lambda)\neq0$. Здесь
$a^{r,0}(x,\xi):=\bigl(\,a^{r,0}_{j,k}(x,\xi)\,\bigr)_{j,k=1}^{p}$~--- квадратная
матрица, произвольный элемент $a^{r,0}_{j,k}(x,\xi)$ которой определен так: он~---
либо главный символ ПДО $A^{(r)}_{j,k}$ в случае $\mathrm{ord}\,A^{(r)}_{j,k}=
mr+m_{k}-m_{j}$, либо нулевая функция в противном случае. При этом в случае $r\geq1$
функции $a^{r,0}_{j,k}(x,\xi)$ считаются равными 0 при $\xi=0$ (такое допущение
связано с тем, что главные символы изначально не определены при $\xi=0$).
\end{definition}


Далее в п. \ref{sec2.6.4} предполагается, что система \eqref{2.121} удовлетворяет
определению \ref{def2.18}.

Из определению \ref{def2.18} следует, что система \eqref{2.121} эллиптическая по
Дуглису--Ниренбергу при каждом фиксированном $\lambda\in\mathbb{C}$. В~самом деле, в
силу \eqref{2.120} матрица $A(\lambda)$ образована элементами
\begin{equation}\label{2.123}
\sum_{r\,=\,0}^{q}\,\lambda^{q-r}\,A^{(r)}_{j,k},\quad j,k=1,\ldots,p,
\end{equation}
которые являются классическими скалярными ПДО порядка не выше, чем ${l_{j}+m'_{k}}$,
где $l_{j}:=-m_{j}$ и $m'_{k}:=mq+m_{k}$. Главный символ ПДО \eqref{2.123} равен
$a^{q,0}_{j,k}(x,\xi)$ для каждого фиксированного $\lambda$. Согласно условию
\eqref{2.122}, в котором берем $\lambda:=0$, справедливо неравенство
$$
\det\,\bigl(\,a^{q,0}_{j,k}(x,\xi)\,\bigr)_{j,k=1}^{p}\neq0\;\;\mbox{для
любых}\;\;x\in\Gamma,\;\xi\in T_{x}^{\ast}\Gamma\setminus\{0\}.
$$
Последнее означает, что система \eqref{2.121} эллиптическая по
Дуглису--Ни\-рен\-бер\-гу при каждом $\lambda\in\mathbb{C}$.

Следовательно, для эллиптической системы \eqref{2.121} справедлива теорема
\ref{th2.28}, согласно которой ограниченный оператор
\begin{equation}\label{2.124}
A(\lambda):\,\bigoplus_{k=1}^{p}\,H^{s+mq+m_{k},\varphi}(\Gamma)\,\rightarrow\,
\bigoplus_{j=1}^{p}\,H^{s+m_{j},\varphi}(\Gamma)
\end{equation}
нетеров для произвольных $\lambda\in\mathbb{C}$, $s\in\mathbb{R}$ и
$\varphi\in\mathcal{M}$. Более того, поскольку система \eqref{2.121} эллиптическая с
параметром в угле $K$, этот оператор имеет следующие дополнительные свойства.


\begin{theorem}\label{th2.33}
Справедливы следующие утверждения.
\begin{itemize}
\item [$\mathrm{(i)}$] Существует число $\lambda_{0}>0$ такое, что для
каждого значения параметра $\lambda\in K$, удовлетворяющего условию $|\lambda|\geq
\lambda_{0}$, матричный ПДО $A(\lambda)$ устанавливает при любых $s\in\mathbb{R}$ и
$\varphi\in\mathcal{M}$ гомеоморфизм
\begin{equation}\label{2.125}
A(\lambda):\,\bigoplus_{k=1}^{p}\,H^{s+mq+m_{k},\varphi}(\Gamma)\,\leftrightarrow\,
\bigoplus_{j=1}^{p}\,H^{s+m_{j},\varphi}(\Gamma).
\end{equation}
\item [$\mathrm{(ii)}$] Для произвольно выбранных параметров
$s\in\mathbb{R}$ и $\varphi\in\nobreak\mathcal{M}$ найдется число $c\geq1$ такое,
что для каждого значения параметра $\lambda\in K$, $|\lambda|\geq\lambda_{0}$, и
любых вектор-функций
\begin{gather}\label{2.126}
u=\mathrm{col}\,(u_{1},\ldots,u_{p})\in
\bigoplus_{k=1}^{p}\,H^{s+mq+m_{k},\varphi}(\Gamma),\\
f=\mathrm{col}\,(f_{1},\ldots,f_{p})\in
\bigoplus_{j=1}^{p}\,H^{s+m_{j},\varphi}(\Gamma),\label{2.127}
\end{gather}
удовлетворяющих уравнению \eqref{2.121}, верна двусторонняя оценка
\begin{gather}\notag
c^{-1}\sum_{j=1}^{p}\;\|f_{j}\|_{H^{s+m_{j},\varphi}(\Gamma)}\leq \\ \notag
\sum_{k=1}^{p}\;\|u_{k}\|_{H^{s+mq+m_{k},\varphi}(\Gamma)}+
|\lambda|^{q}\,\sum_{k=1}^{p}\;\|u_{k}\|_{H^{s+m_{k},\varphi}(\Gamma)}\leq\\
c\,\sum_{j=1}^{p}\;\|f_{j}\|_{H^{s+m_{j},\varphi}(\Gamma)}.\label{2.128}
\end{gather}
Здесь число $c$ не зависит от параметра $\lambda$ и вектор-функций $u$,~$f$.
\end{itemize}
\end{theorem}


В случае $\varphi\equiv1$ (шкала соболевских пространств) эта теорема известна
[\ref{Agranovich90}, с.~64]. Отметим, что левое неравенство в двусторонней оценке
\eqref{2.128} справедливо без предположения об эллиптичности с параметром системы
\eqref{2.121} (сравнить с [\ref{AgranovichVishik64}, с.~65]). Чтобы избежать
громоздкой формулы, мы используем в \eqref{2.128} эквивалентные негильбертовы нормы
в пространствах \eqref{2.126} и \eqref{2.127}.

Докажем отдельно утверждения (i) и (ii) теоремы \ref{th2.33}. Мы выведем общий
случай $\varphi\in\mathcal{M}$ из соболевского случая $\varphi\equiv1$.

\medskip

\textbf{Доказательство теоремы~\ref{th2.33}.} Пусть $s\in\mathbb{R}$ и
$\varphi\in\mathcal{M}$.

(i) Поскольку для каждого $\lambda\in\mathbb{C}$ система \eqref{2.121} эллиптическая
по Дуглису--Ниренбергу, то в силу теоремы \ref{th2.28} ограниченный нетеров оператор
\eqref{2.124} имеет независящие от $s$ и $\varphi$ конечномерные ядро
$\mathcal{N}(\lambda)$ и дефектное подпространство $\mathcal{N}^{+}(\lambda)$. Далее
воспользуемся тем, что теорема \ref{th2.33} верна в случае $\varphi\equiv1$.
Существует число $\lambda_{0}>0$ такое, что для каждого $\lambda\in K$,
удовлетворяющего условию $|\lambda|\geq\lambda_{0}$, матричный ПДО $A(\lambda)$
устанавливает гомеоморфизм
$$
A(\lambda):\,\bigoplus_{k=1}^{p}\,H^{s+mq+m_{k}}(\Gamma)\,\leftrightarrow\,
\bigoplus_{j=1}^{p}\,H^{s+m_{j}}(\Gamma).
$$
Следовательно, для указанного $\lambda$ пространства $\mathcal{N}(\lambda)$ и
$\mathcal{N}^{+}(\lambda)$ тривиальны, т.~е. линейный ограниченный оператор
\eqref{2.124} является биекцией. Отсюда по теореме Банаха об обратном операторе
получаем гомеоморфизм \eqref{2.125}. Утверждение (i) теоремы~\ref{th2.33} доказано.

(ii) Поскольку теорема \ref{th2.33} справедлива в соболевском случае
$\varphi\equiv1$, существует такое число $\lambda_{0}>0$, что для каждого значения
параметра $\lambda\in K$, удовлетворяющего условию $|\lambda|\geq \lambda_{0}$,
справедливы гомеоморфизмы
\begin{equation}\label{2.142}
A(\lambda):\;\bigoplus_{k=1}^{p}\,H^{s\mp1+mq+m_{k}}(\Gamma,|\lambda|^{q},mq)
\,\leftrightarrow\,\bigoplus_{j=1}^{p}\,H^{s\mp1+m_{j}}(\Gamma),
\end{equation}
причем норма оператора \eqref{2.142} вместе с нормой обратного оператора ограничены
равномерно по параметру $\lambda$. Мы используем гильбертовы пространства
$H^{\sigma,\varphi}(\Gamma,\varrho,\theta)$, зависящие от дополнительных параметров
$\varrho=|\lambda|^{q}$ и $\theta=mq$, определенные в п. \ref{sec2.6.4b}. Пусть
$\psi$~--- интерполяционный параметр из предложения \ref{th2.21}, где берем
$\varepsilon=\delta=1$. Применив интерполяцию с этим параметром к \eqref{2.142},
получим гомеоморфизм
\begin{gather}\notag
A(\lambda):\\ \notag
\Bigl[\,\bigoplus_{k=1}^{p}H^{s-1+mq+m_{k}}(\Gamma,|\lambda|^{q},mq),\,
\bigoplus_{k=1}^{p}H^{s+1+mq+m_{k}}(\Gamma,|\lambda|^{q},mq)\Bigr]_{\psi}\leftrightarrow \\
\Bigl[\,\bigoplus_{j=1}^{p}\,H^{s-1+m_{j}}(\Gamma),\,
\bigoplus_{j=1}^{p}H^{s+1+m_{j}}(\Gamma)\Bigr]_{\psi}. \label{2.143}
\end{gather}
При этом в силу теоремы \ref{th2.8} норма оператора \eqref{2.143} и норма обратного
к нему оператора ограничены равномерно по параметру~$\lambda$. (Записанные в формуле
\eqref{2.143} допустимые пары пространств нормальные.) Отсюда на основании теоремы
\ref{th2.5} об интерполяции прямых сумм пространств получаем гомеоморфизм
\begin{gather}\notag
A(\lambda):\\ \notag \bigoplus_{k=1}^{p}
\left[\,H^{s-1+mq+m_{k}}(\Gamma,|\lambda|^{q},mq),
\,H^{s+1+mq+m_{k}}(\Gamma,|\lambda|^{q},mq)\right]_{\psi}\leftrightarrow\\
\bigoplus_{j=1}^{p}
\left[\,\,H^{s-1+m_{j}}(\Gamma),\,H^{s+1+m_{j}}(\Gamma)\,\right]_{\psi}.
\label{2.144}
\end{gather}
 При этом нормы операторов \eqref{2.143} и \eqref{2.144} равны, а также равны нормы
обратных к ним операторов. Теперь воспользуемся леммой \ref{lem2.14}, где полагаем
$$
\sigma:=s+mq+m_{k},\;\;\varrho:=|\lambda|^{q},\;\;\theta:=mq,\;\;\varepsilon=\delta=1,
$$
и теоремой \ref{th2.21}. Согласно ним \eqref{2.144} влечет гомеоморфизм
\begin{equation}\label{2.145}
A(\lambda):\;\bigoplus_{k=1}^{p}\,H^{s+mq+m_{k},\varphi}(\Gamma,|\lambda|^{q},mq)\,
\leftrightarrow\,\bigoplus_{j=1}^{p}\,H^{s+m_{j},\varphi}(\Gamma)
\end{equation}
такой, что норма оператора \eqref{2.145} вместе с нормой обратного оператора
ограничены равномерно по параметру $\lambda$. Это означает двустороннюю оценку
\eqref{2.128}, где число $c$ не зависит от параметра $\lambda$ и от вектор-- функций
\eqref{2.126}, \eqref{2.127}. Утверждение (ii) теоремы~\ref{th2.33} доказано.

\medskip

Из теоремы \ref{th2.33} (i) вытекает следующее утверждение об индексе оператора,
соответствующего эллиптической системе с параметром.


\begin{corollary}\label{cor2.5}
Предположим, что система \eqref{2.121} эллиптическая с параметром на некотором
замкнутом луче $K:=\{\lambda\in\mathbb{C}:\,\arg\lambda =\mathrm{const}\}$. Тогда
оператор \eqref{2.124} имеет нулевой индекс при любом $\lambda\in\mathbb{C}$.
\end{corollary}


\textbf{Доказательство.} Для каждого фиксированного $\lambda\in\mathbb{C}$ система
\eqref{2.121} эллиптическая по Дуглису--Ниренбергу. Поэтому в силу теоремы
\ref{th2.28} оператор \eqref{2.124} имеет конечный индекс, не зависящий от
$s\in\mathbb{R}$ и $\varphi\in\mathcal{M}$. Кроме того, этот индекс не зависит и от
параметра $\lambda$. В самом деле, в силу \eqref{2.120} параметр $\lambda$ влияет
лишь на младшие члены элементов матричного ПДО $A(\lambda)$:
$$
A(\lambda)-A(0)=\sum_{r=0}^{q-1}\,\lambda^{q-r}A^{(r)}=
\Bigl(\,\sum_{r=0}^{q-1}\,\lambda^{q-r}A^{(r)}_{j,k}\,\Bigr)_{j,k=1}^{p},
$$
где
$$
\mathrm{ord}\;\sum_{r=0}^{q-1}\,\lambda^{q-r}A^{(r)}_{j,k}\;\leq\,
m(q-1)+m_{k}-m_{j}.
$$
Поэтому ввиду леммы \ref{lem2.13} мы имеем ограниченный оператор
$$
A(\lambda)-A(0):\,
\bigoplus_{k=1}^{p}\,H^{s+mq+m_{k},\varphi}(\Gamma)\,\rightarrow\,
\bigoplus_{j=1}^{p}\,H^{s+m+m_{j},\varphi}(\Gamma).
$$
Но в силу теоремы \ref{th2.22} (iii) и условия $m>0$ справедливо компактное вложение
пространства $H^{s+m+m_{j},\varphi}(\Gamma)$ в пространство
$H^{s+m_{j},\varphi}(\Gamma)$. Следовательно, компактен оператор
$$
A(\lambda)-A(0):\,
\bigoplus_{k=1}^{p}\,H^{s+mq+m_{k},\varphi}(\Gamma)\,\rightarrow\,
\bigoplus_{j=1}^{p}\,H^{s+m_{j},\varphi}(\Gamma).
$$
Отсюда вытекает (см., например, [\ref{Hermander87}, с.~249]), что операторы
$A(\lambda)$ и $A(0)$ имеют одинаковый индекс, т.~е. последний не зависит от
параметра $\lambda$. Далее, согласно теореме \ref{th2.33} (i) при достаточно больших
по модулю значениях параметра $\lambda\in K$, справедлив гомеоморфизм \eqref{2.125}.
Следовательно, индекс оператора $A(\lambda)$ равен нулю при $\lambda\in K$,
$|\lambda|\gg1$, а значит, и при любом $\lambda\in\mathbb{C}$. Следствие
\ref{cor2.5} доказано.




\markright{\emph\ref{sec4.1.6}. Эллиптические краевые задачи для систем уравнений}

\section[Эллиптические краевые задачи для систем уравне-\break ний]
{Эллиптические краевые задачи \\ для систем уравнений}\label{sec4.1.6}

\markright{\emph\ref{sec4.1.6}. Эллиптические краевые задачи для систем уравнений}

Здесь мы изучим эллиптические краевые задачи для систем дифференциальных уравнений в
уточненной шкале пространств. Известно, что они порождают ограниченные и нетеровы
операторы в соответствующих парах позитивных пространств Соболева; см. монографию
[\ref{WlokaRowleyLawruk95}] (п.~9.4), обзор [\ref{Agranovich97}] (\S~6) и
приведенную там литературу. Мы уточним этот результат для шкалы пространств
Хермандера применительно к системам, эллиптическим по Петровскому.

\subsection{Постановка задачи}\label{sec4.1.6.1}

Рассмотрим линейную систему дифференциальных уравнений
\begin{equation}\label{4.60}
\sum_{k=1}^{p}\,L_{j,k}\,u_{k}=f_{j}\quad\mbox{в}\quad\Omega,\quad j=1,\ldots,p.
\end{equation}
Здесь $p\in\mathbb{N}$, а $L_{j,k}=L_{j,k}(x,D)$, где $j,k=1,\ldots,p$, ---
скалярные линейные дифференциальные выражения, заданные на $\overline{\Omega}$.
Выражение $L_{j,k}$ имеет произвольный конечный порядок и бесконечно гладкие на
$\overline{\Omega}$ комплекснозначные коэффициенты. Для каждого номера
$k=1,\ldots,p$ положим
$$
m_{k}:=\max\{\mathrm{ord}\,L_{j,k}:\,j=1,\ldots,p\}.
$$
Таким образом, $m_{k}$ --- максимальный порядок дифференцирования искомой функции
$u_{k}$. Предполагается, что все $m_{k}\geq1$.

Свяжем с системой \eqref{4.60} квадратную матрицу порядка $p$:
$$
\mathbf{L}^{(0)}(x,\xi):=\bigl(L_{j,k}^{(0)}(x,\xi)\bigr)_{j,k=1}^{p},\quad
x\in\overline{\Omega},\;\;\xi\in\mathbb{C}^{n}.
$$
Здесь $L_{j,k}^{(0)}(x,\xi)$ --- главный символ дифференциального выражения
$L_{j,k}$ в случае, когда $\mathrm{ord}\,L_{j,k}=m_{k}$, либо
$L_{j,k}^{(0)}(x,\xi)\equiv0$ в случае, когда $\mathrm{ord}\,L_{j,k}<m_{k}$.

\medskip

\begin{definition}\label{def4.3}
Система \eqref{4.60} называется \emph{эллиптической по Петровскому} на
$\overline{\Omega}$, если $\det\mathbf{L}^{(0)}(x,\xi)\neq0$ для любой точки
$x\in\overline{\Omega}$ и произвольного ненулевого вектора $\xi\in\mathbb{R}^{n}$.
\end{definition}

\medskip

Если система \eqref{4.60} эллиптическая по Петровскому и $n\geq3$, то однородный
полином $\det\mathbf{L}^{(0)}(x,\xi)$ имеет четный порядок [\ref{Agranovich97},
с.~53]:
$$
\sum_{k=1}^{p}\,m_{k}=2q\quad\mbox{для некторого}\quad q\in\mathbb{N}.
$$
Мы предполагаем, что это условие выполняется при всяком целом $n\geq2$.

Будем рассматривать решения системы \eqref{4.60}, удовлетворяющие краевым условиям
\begin{equation}\label{4.61}
\sum_{k=1}^{p}\,B_{j,k}\,u_{k}=g_{j}\quad\mbox{на}\quad\Gamma,\quad j=1,\ldots,q.
\end{equation}
Здесь $B_{j,k}=B_{j,k}(x,D)$, где $j=1,\ldots,q$ и $k=1,\ldots,p$, --- скалярные
линейные граничные дифференциальные выражения, заданные на $\Gamma$. Выражение
$B_{j,k}$ имеет порядок $\mathrm{ord}\,B_{j,k}\leq m_{k}-1$ и бесконечно гладкие
комплекснозначные коэффициенты. Для каждого номера $j=1,\ldots,q$ положим
$$
r_{j}:=\min\{m_{k}-\mathrm{ord}\,B_{j,k}:\,k=1,\ldots,p\}\geq1.
$$
Здесь мы принимаем $\mathrm{ord}\,B_{j,k}:=-\infty$, если $B_{j,k}\equiv0$. Теперь
$\mathrm{ord}\,B_{j,k}\leq m_{k}-r_{j}$ для всех $j=1,\ldots,q$ \,и
\,$k=1,\ldots,p$.

Свяжем с краевыми условиями \eqref{4.61} матрицу размера $q\times p$:
$$
\mathbf{B}^{(0)}(x,\xi):= \bigl(B_{j,k}^{(0)}(x,\xi)\bigr)_{\substack{j=1,\ldots,q
\\ k=1,\ldots,p}}\,,\quad x\in\Gamma,\;\;\xi\in\mathbb{C}^{n}.
$$
Здесь $B_{j,k}^{(0)}(x,\xi)$ --- главный символ дифференциального выражения
$B_{j,k}$ в случае, когда $\mathrm{ord}\,B_{j,k}=m_{k}-r_{j}$, либо
$B_{j,k}^{(0)}(x,\xi)\equiv0$ в случае, когда $\mathrm{ord}\,B_{j,k}<m_{k}-r_{j}$.

\medskip

\begin{definition}\label{def4.4}
Краевая задача \eqref{4.60}, \eqref{4.61} называется эллиптической по Петровскому в
области $\Omega$, если выполняются следующие условия:
\begin{itemize}
\item [$\mathrm{(i)}$] Система \eqref{4.60} является правильно эллиптической
на $\overline{\Omega}$, т.~е. для произвольной точки $x\in\Gamma$ и любых линейно
независимых векторов $\xi',\xi''\in\mathbb{R}^n$ многочлен
$\det\mathbf{L}^{(0)}(x,\xi'+\tau\xi'')$ переменного $\tau$ имеет ровно $q$ корней
$\tau^{+}_{j}(x;\xi',\xi'')$, $j=1,\ldots,q$, с положительной мнимой частью и
столько же корней с отрицательной мнимой частью (выписанных с учетом их кратности).
\item [$\mathrm{(ii)}$] Краевые условия \eqref{4.61} удовлетворяют условию
дополнительности по отношению к системе \eqref{4.60} на $\Gamma$, т.~е. для
произвольной точки $x\in\Gamma$ и любого вектора $\xi\neq0$ касательного к $\Gamma$
в точке $x$, строки матрицы
$$
\mathbf{B}^{(0)}(x,\xi+\tau\nu(x))\cdot
\mathbf{L}^{(0)}_{\mathrm{c}}(x,\xi+\tau\nu(x)),
$$
элементы которой рассматриваются как многочлены от $\tau$, линейно независимы по
модулю многочлена $\prod_{j=1}^{q}(\tau-\tau_{j}^{+}(x;\xi,\nu(x)))$. Здесь
$\mathbf{L}^{(0)}_{\mathrm{c}}(x,\xi)$~--- транспониро\-ван\-ная матрица
алгебраических дополнений элементов матрицы $\mathbf{L}^{(0)}(x,\xi)$.
\end{itemize}
\end{definition}

\medskip

\begin{remark}\label{rem4.5}
Если система \eqref{4.60} удовлетворяет условию (i) определения \ref{def4.4}, то она
эллиптическая по Петровскому на $\overline{\Omega}$. Обратное верно при $n\geq3$
[\ref{Agranovich97}, с.~53].
\end{remark}

\medskip

Приведем примеры эллиптических краевых задач для систем уравнений.

\medskip

\begin{example}\label{ex4.4}
Эллиптическая краевая задача для системы Коши-Римана:
\begin{gather*}
\frac{\partial u_{1}}{\partial x_{1}}-\frac{\partial u_{2}}{\partial
x_{2}}=f_{1},\quad \frac{\partial u_{1}}{\partial x_{2}}+\frac{\partial
u_{2}}{\partial x_{1}}=f_{2}\quad\mbox{в}\quad\Omega,\\
u_{1}+u_{2}=g\quad\mbox{на}\quad\Gamma.
\end{gather*}
Здесь $n=p=2$, $m_{1}=m_{2}=1$ и, значит, $q=1$. Система Коши-Римана --- пример
однородной эллиптической системы. Такие системы удовлетворяют определению
\ref{def4.3}, где $m_{1}=\ldots=m_{p}$.
\end{example}

\medskip

\begin{example}\label{ex4.5}
Эллиптическая по Петровскому краевая задача:
\begin{gather*}
\frac{\partial u_{1}}{\partial x_{1}}-\frac{\partial^{3}u_{2}}{\partial
x_{2}^{3}}=f_{1},\quad \frac{\partial u_{1}}{\partial x_{2}}+\frac{\partial^{3}
u_{2}}{\partial x_{1}^{3}}=f_{2}\quad\mbox{в}\quad\Omega, \\
u_{1}=g_{1},\quad u_{2}\,\biggl(\mbox{или}\;\frac{\partial
u_{2}}{\partial\nu},\;\mbox{или}\;\frac{\partial^{2}
u_{2}}{\partial\nu^{2}}\biggr)=g_{2}\quad\mbox{на}\quad\Gamma.
\end{gather*}
Здесь $n=p=2$, $m_{1}=1$, $m_{2}=3$ и, значит, $q=2$. Отметим, что рассмотренная
здесь система не является однородной эллиптической [\ref{Volevich65}, с.~376].
\end{example}

\medskip

Другие содержательные примеры приведены, например, в [\ref{Agranovich97}] (\S~6.2).

\subsection{Теорема о разрешимости}\label{sec4.1.6.2}

Запишем краевую задачу \eqref{4.60}, \eqref{4.61} в матричной форме:
$$
\mathbf{L}u=f\;\;\mbox{в}\;\;\Omega,\quad\mathbf{B}u=g\;\;\mbox{на}\;\;\Gamma.
$$
Здесь $\mathbf{L}:=(L_{j,k})_{j,k=1}^{p}$ и
$\mathbf{B}:=(B_{j,k})_{\substack{j=1,\ldots,q
\\ k=1,\ldots,p}}$~--- матричные дифференциальные выражения, а
$u:=\mathrm{col}\,(u_{1},\ldots,u_{p})$, $f:=\mathrm{col}\,(f_{1},\ldots,f_{p})$ и
$g:=\mathrm{col}\,(g_{1},\ldots,g_{q})$~--- функциональные столбцы.

\medskip

\begin{theorem}\label{th4.11}
Предположим, что краевая задача \eqref{4.60}, \eqref{4.61} эллиптическая по
Петровскому в области $\Omega$. Пусть $s>0$ и $\varphi\in\mathcal{M}$. Тогда
отображение $u\mapsto(\mathbf{L}u,\mathbf{B}u)$, где
$u\in(C^{\infty}(\,\overline{\Omega}\,))^{p}$, продолжается по непрерывности до
ограниченного нетерового оператора
\begin{gather}\notag
(\mathbf{L},\mathbf{B}):\,\bigoplus_{k=1}^{p}H^{s+m_{k},\varphi}(\Omega)\rightarrow\\
(H^{s,\varphi}(\Omega))^{p}\oplus\bigoplus_{j=1}^{q}
H^{s+r_{j}-1/2,\varphi}(\Gamma)=: \mathbb{H}_{s,\varphi}(\Omega,\Gamma).\label{4.62}
\end{gather}
Ядро $\mathbf{N}$ оператора \eqref{4.62} лежит в
$(C^{\infty}(\,\overline{\Omega}\,))^{p}$ и не зависит от $s$, $\varphi$. Область
значений этого оператора состоит из всех вектор-функций
$$
(f_{1},\ldots,f_{p};g_{1},\ldots,g_{q})\in \mathbb{H}_{s,\varphi}(\Omega,\Gamma)
$$
таких, что
$$
\sum_{j=1}^{p}\,(f_{j},w_{j})_{\Omega}+\sum_{j=1}^{q}\, (g_{j},h_{j})_{\Gamma}=0
$$
для произвольной вектор-функции
$$
(w_{1},\ldots,w_{p};\,h_{1},\ldots,h_{q})\in\mathbf{W}.
$$
Здесь $\mathbf{W}$ --- некоторое не зависящее от $s$ и $\varphi$ конечномерное
подпространство в $(C^{\infty}(\,\overline{\Omega}\,))^{p}\times
(C^{\infty}(\Gamma))^{q}$. Индекс оператора \eqref{4.62} равен
$\dim\mathbf{N}-\dim\mathbf{W}$ и также не зависит от $s$,~$\varphi$.
\end{theorem}

\medskip

В соболевском случае $\varphi\equiv1$ эта теорема является частным случаем известной
теоремы о разрешимости эллиптических краевых задач для общих систем смешанного
порядка; см. например монографии [\ref{FunctionalAnalysis72}, с.~176],
[\ref{WlokaRowleyLawruk95}, с.~393], обзор [\ref{Agranovich97}, с.~58]. Отсюда общий
случай $\varphi\in\mathcal{M}$ выводится с помощью интерполяции аналогично
доказательству теоремы~\ref{th4.1}.




\markright{\emph \ref{notes5}. Примечания и комментарии}

\section{Примечания и комментарии}\label{notes5}

\markright{\emph \ref{notes5}. Примечания и комментарии}

\small

\textbf{К п. 5.1.} Широкие классы эллиптических систем дифференциальных уравнений
введены и изучены И.~Г.~Петровским [\ref{Petrovsky39}; \ref{Petrovsky61}, с.~328;
\ref{Petrovsky86}, с.~235] и А.~Дуглисом, Л.~Ниренбергом [\ref{DouglisNirenberg55}].
Для этих систем установлены теоремы о гладкости решений и получены внутренние
априорные оценки решений в подходящих парах пространств Гельдера (с нецелыми
индексами) и пространств Соболева. См. также монографию Л.~Хермандера
[\ref{Hermander65}] (п.~10.6), где рассмотрен случай соболевской шкалы. Им же
[\ref{Hermander67}, с.~175] установлены априорные оценки решений эллиптических
систем псевдодифференциальных уравнений.

Все теоремы этого пункта доказаны в [\ref{08UMB3}]. Они распространяют результаты
п.~1.4 на равномерно эллиптические по Дуглису--Ниренбергу системы
псевдодифференциальных уравнений, заданные в евклидовом пространстве. Случай систем,
равномерно эллиптических по Петровскому, рассмотрен отдельно в [\ref{09Collection1},
\ref{09UMJ3}].

\medskip

\textbf{К п. 5.2.} Теория эллиптических систем на замкнутых гладких многообразиях
изложена, например, в монографиях Л.~Хермандера [\ref{Hermander87}] (гл.~19) и
обзоре М.~С.~Аграновича [\ref{Agranovich90}] (п.~3.2). Для таких систем априорные
оценки их решений равносильны тому, что ограниченный оператор, отвечающий системе,
является нетеровым в соответствующих парах соболевских пространств. Его индекс
вычислен М.~Ф.~Атьёй и И.~М.~Зингером [\ref{AtiyahSinger63}].

Для систем, эллиптических с параметром, индекс равен нулю и соответствующий оператор
осуществляет гомеоморфизмы в подходящих парах соболевских пространств для больших по
модулю значений параметра. См., например, обзор М.~С.~Аграновича
[\ref{Agranovich90}] (п.~4.3). Различные классы эллиптических систем с параметром
изучены в работах А.~Н.~Кожевникова [\ref{Kozhevnikov73}, \ref{Kozhevnikov96}] и
М.~С.~Аграновича [\ref{Agranovich90b}, \ref{Agranovich92}].

Между эллиптическим краевыми задачами и матричными эллиптическими ПДО имеется тесная
связь~--- каждая эллиптическая краевая задача эквивалентна некоторой эллиптической
системе псевдодифференциальных уравнений на границе области. Это позволяет применить
теорию эллиптических систем для доказательства теорем о разрешимости эллиптических
краевых задач; см., например, монографии Й.~Т.~Влоки, Б.~Ровли и Б.~Лорака
[\ref{WlokaRowleyLawruk95}] (ч.~IV), Ю.~В.~Егорова [\ref{Egorov84}] (гл.~III, \S~3),
Л.~Хермандера [\ref{Hermander87}] (гл.~20).

Системы, близкие по своим свойствам к эллиптическим, изучены в работах
Б.~Р.~Вайнберга, В.~В.~Грушина [\ref{VainbergGrushin67a}] и Р.~С.~Сакса
[\ref{Saks78}, \ref{Saks79}].

Все теоремы этого пункта доказаны в [\ref{08MFAT2}]. Они распространяют результаты
п.~2.2 на эллиптические по Дуглису--Ниренбергу системы псевдодифференциальных
уравнений. Случай систем, эллиптических по Петровскому, рассмотрен отдельно в
[\ref{07Dop5}, \ref{08BPAS3}].

\medskip

\textbf{К п. 5.3.} Эллиптические краевые задачи для различных систем
дифференциальных уравнений исследованы в работах С.~Агмона, А.~Дуглиса и
Л.~Ниренберга [\ref{AgmonDouglisNirenberg64}], М.~C.~Аграновича и А.~С.~Дынина
[\ref{AgranovichDynin62}], Л.~Р.~Волевича [\ref{Volevich63}, \ref{Volevich65}],
Л.~Н.~Слободецкого [\ref{Slobodetsky58a}, \ref{Slobodetsky60}], В.~А.~Солонникова
[\ref{Solonnikov63}, \ref{Solonnikov64}, \ref{Solonnikov66}, \ref{Solonnikov67}],
Л.~Хермандера [\ref{Hermander65}] (п.~10.6), [\ref{Hermander87}] (п.~19.5). См.
также монографию Й.~Т.~Влоки, Б.~Ровли и Б.~Лорака [\ref{WlokaRowleyLawruk95}],
посвященную эллиптическим краевым задачам для систем, обзор М.~С.~Аграновича
[\ref{Agranovich97}] (\S~6) и приведенную там литературу. Установлено, что оператор
задачи ограничен и нетеров в соответствующих парах позитивных пространств Соболева.

В двусторонней модифицированной соболевской шкале эти задачи изучены в работах
Я.~А.~Ройтберга и З.~Г.~Шефтеля [\ref{RoitbergSheftel69}, \ref{RoitbergSheftel72}],
И.~А.~Коваленко [\ref{Kovalenko73}], Я.~А.~Ройтберга [\ref{Roitberg75}],
И.~Я.~Ройтберг и Я.~А.~Ройтберга [\ref{RoitbergRoitberg00}]. Ими установлены теоремы
о полном наборе гомеоморфизмов, который порождает оператор, отвечающий задаче.
Изложение этих результатов дано в монографиях Я.~А.~Ройтберга [\ref{Roitberg96}]
(гл.~10) и [\ref{Roitberg99}] (п.~1.3) применительно к системам, эллиптическим по
Дуглису--Ниренбергу.

Краевые задачи для систем, близким по своим свойствам к эллиптическим, изучались в
работах Б.~Р.~Вайнберга, В.~В.~Грушина [\ref{VainbergGrushin67b}] и Р.~С.~Сакса
[\ref{Saks77}, \ref{Saks80}, \ref{Saks81}]. Эти задачи исследовались в парах
позитивных соболевских пространств.

Основной результат этого пункта~--- теорема 5.10 о разрешимости в уточненной шкале
краевых задач для эллиптических по Петровскому систем установлен в [\ref{07Dop6}].
Он обобщает теорему 4.7 на эллиптические системы уравнений.

\normalsize


\chaptermark{}



\chapter*{\textbf{Литература}}

\chaptermark{\emph Литература}



\addcontentsline{toc}{chapter}{\textbf{Литература}}
\thispagestyle{empty}

\small

\begin{enumerate}

\label{part_bibl}

\item\label{AgmonDouglisNirenberg62}
Агмон С., Дуглис А., Ниренберг Л. Оценки вблизи границы решений эллиптических
уравнений в частных производных при общих граничных условиях.~I.~--- Москва: Изд.
иностр. лит., 1962.~--- 206~с. (Перевод с английского статьи в Commun. Pure Appl.
Math.~--- 1959.~--- \textbf{12}, no.~4.~--- P. 623--727.)

\item\label{Agranovich65}
Агранович М. С. Эллиптические сингулярные интегро-дифферен\-ци\-аль\-ные
операторы~// Успехи матем. наук.~--- 1965.~--- \textbf{20}, №~5.~--- С. 3--120.

\item\label{Agranovich90}
Агранович М. С. Эллиптические операторы на замкнутых многообразиях~// Итоги науки и
техн. ВИНИТИ. Совр. пробл. матем. Фунд. напр.~--- 1990.~--- \textbf{63}.~--- С.
5--129.

\item\label{Agranovich90b}
Агранович М. С. О несамосопряженных задачах с параметром, эллиптических по
Агмону--Дуглису--Ниренбергу~// Функц. анализ и его прил.~--- 1990.~--- \textbf{24},
№~1.~--- С. 59--61.

\item\label{Agranovich92}
Агранович М. С. О модулях собственных значений несамосопряженных задач с параметром,
эллиптических по Агмону--Дуглису--Ниренбергу~// Функц. анализ и его прил.~---
1992.~--- \textbf{26}, №~2.~--- С. 51--55.

\item\label{AgranovichVishik64}
Агранович М. С., Вишик М. И. Эллиптические задачи с параметром и параболические
задачи общего вида~// Успехи матем. наук.~--- 1964.~--- \textbf{19}, №~3.~--- С.
53--161.

\item\label{AgranovichDynin62}
Агранович М. С., Дынин А. С. Общие краевые задачи для эллиптических систем в
многомерной области~// Доклады АН СССР.~--- 1962.~--- \textbf{146}, №~3.~--- С.
511--514.

\item\label{Alexits63}
Алексич Г. Проблемы сходимости ортогональных рядов.~--- Москва: Иностр. лит.,
1963.~--- 360~с.

\item\label{BergLefstrem80}
Берг Й., Лёфстрём Й. Интерполяционные пространства. Введение.~--- Москва: Мир,
1980.~--- 264~с.

\item\label{BerezanskyKreinRoitberg63}
Березанский Ю. М., Крейн С.~Г., Ройтберг Я.~А. Теорема о гомеоморфизмах и локальное
повышение гладкости вплоть до границы решений эллиптических уравнений~// Доклады АН
СССР.~--- 1963.~--- \textbf{148}, №~4.~--- С. 745--748.

\item\label{Berezansky65}
Березанский Ю. М. Разложение по собственным функциям самосопряженных операторов.~---
Киев: Наукова думка, 1965.~--- 800~с.

\item\label{BerezanskyUsSheftel90}
Березанский Ю. М., Ус Г.~Ф., Шефтель З. Г.  Функциональный анализ.~--- Киев: Выща
школа, 1990.~--- 600~с.

\item\label{BesovIlinNikolsky75}
Бесов О. В., Ильин В.~П. Интегральные представления функций и теоремы вложения.~---
Москва: Наука, 1975.~--- 480~с.

\item\label{VainbergGrushin67a}
Вайнберг Б. Р., Грушин В. В. О равномерно неэллиптических задачах.~I~//
Математический сборник.~--- 1967.~--- \textbf{72(114)}, №~1.~--- С. 602--636.

\item\label{VainbergGrushin67b}
Вайнберг Б. Р., Грушин В. В. О равномерно неэллиптических задачах.~II~//
Математический сборник.~--- 1967. --- \textbf{73(115)}, №~4.~--- С. 126--154.

\item\label{VishikEskin69}
Вишик М. И., Эскин Г.~И. Смешанные краевые задачи для эллиптических систем
дифференциальных уравнений~// Труды Ин-та прикладной матем. Тбилисского гос.
ун-та.~--- 1969.~--- \textbf{11}.~--- С. 31--48.

\item\label{Vladimirov79}
Владимиров В. С. Обобщенные функции в математической физике.~--- Москва: Наука,
1979.~--- 320~с.

\item\label{Volevich63}
Волевич Л. Р. К теории краевых задач для общих эллиптических систем~// Доклады АН
СССР.~--- 1963.~--- \textbf{148}, №~3.~--- С. 489--492.

\item\label{Volevich65}
Волевич Л. Р. Разрешимость краевых задач для общих эллиптических систем~// Матем.
сборник.~--- 1965.~--- \textbf{68}, №~3.~--- С. 373--416.

\item\label{VolevichPaneah65}
Волевич Л. Р., Панеях Б. П. Некоторые пространства обобщенных функций и теоремы
вложения~// Успехи матем. наук.~--- 1965.~--- \textbf{20}, №~1.~--- С. 3--74.

\item\label{GelfandShilov58}
Гельфанд И. М., Шилов Г.~Е. Пространства основных и обобщенных функций.~--- Москва:
Физматгиз, 1958.~--- 307~с.

\item\label{GelfandShilov59}
Гельфанд И. М., Шилов Г.~Е. Обобщенные функции и действия над ними.~--- Москва:
Физматгиз, 1959.~--- 470~с.

\item\label{GilbargTrudinger89}
Гилбарг Д., Трудингер Н. Эллиптические дифференциальные уравнения с частными
производными второго порядка.~--- Москва: Наука, 1989.~--- 464~с.

\item\label{GohbergKrein57}
Гохберг И. Ц., Крейн М.~Г. Основные положения о дефектных числах, корневых векторах
и индексах линейных операторов~// Успехи матем. наук.~--- 1957.~--- \textbf{12},
№~2.~--- С. 43--118.

\item\label{DjakovMityagin06}
Джаков П. Б., Митягин Б. С. Зоны неустойчивости одномерных периодических операторов
Шрёдингера и Дирака~// Успехи матем. наук.~--- 2006.~--- \textbf{61}, №~4.~--- С.
77--182.

\item\label{Egorov84}
Егоров Ю. В. Линейные дифференциальные уравнения главного типа.~--- Москва: Наука,
1984.~--- 360~с.

\item\label{ZhitarashuEidelman92}
Житарашу Н. В., Эйдельман С. Д. Параболические граничные задачи.~--- Кишинев,
Штиинца, 1992.~--- 328~с.

\item\label{Kalugina75}
Калугина Т. Ф. Интерполяция банаховых пространств с функциональным параметром.
Теорема реитерации~// Вестник Московского ун-та. Сер. 1, матем., механика.~---
1975.~--- \textbf{30}, №~6.~--- С. 68--77.

\item\label{KashinSaakyan84}
Кашин Б. С., Саакян А. А. Ортогональные ряды.~--- Москва: Наука, 1984.~--- 496~с.

\item\label{Kovalenko73}
Коваленко И. А. Теоремы об изоморфизмах для эллиптических систем с граничными
условиями, не являющимися нормальными~// Укр. матем. журн.~--- 1973.~---
\textbf{25}, №~3.~--- С. 373--379.

\item\label{Kozhevnikov73}
Кожевников А. Н. Спектральные задачи для псевдодифференциальных систем,
эллиптических по Дуглису~-- Ниренбергу и их приложения~// Матем. сборник.~---
1973.~--- \textbf{92}, №~1.~--- С. 60--88.

\item\label{KostarchukRoitberg73}
Костарчук Ю. В., Ройтберг Я. А. Теореми про ізоморфізми для еліптичних граничних
задач з граничними умовами, які не є нормальними~// Укр. матем. журн.~--- 1973.~---
\textbf{25}, №~2.~--- С. 271--277.

\item\label{Krein60a}
Крейн С.~Г. Об одной интерполяционной теореме в теории операторов~// Доклады АН
СССР.~--- 1960.~--- \textbf{130}, №~3.~--- С. 491--494.

\item\label{Krein60b}
Крейн С.~Г. О понятии нормальной шкалы пространств~// Доклады АН СССР.~--- 1960.~---
\textbf{132}.~--- С. 510--513.

\item\label{KreinPetunin66}
Крейн С. Г., Петунин Ю. И. Шкалы банаховых пространств~// Успехи матем. наук.~---
1966.~--- \textbf{21}, вып.~2 (128).~--- С. 89--168.

\item\label{KreinPetuninSemenov78}
Крейн С. Г., Петунин Ю. И., Семенов Е. М. Интерполяция линейных операторов.~---
Москва: Наука, 1978.~--- 400~с.

\item\label{LadyszenskayaUraltseva64}
Ладыженская О. А., Уральцева Н. Н. Линейные и квазилинейные уравнения эллиптического
типа.~--- Москва: Наука, 1964.~--- 538~с.

\item\label{Lizorkin86}
Лизоркин П. И. Пространства обобщенной гладкости~// Х.~Трибель. Теория
функциональных пространств.~--- Москва: Мир, 1986.~--- С. 381--415.

\item\label{Lions72a}
Лионс Ж.-Л. Оптимальное управление системами, описываемыми уравнениями с частными
производными.~--- Москва: Мир, 1972.~--- 414~с.

\item\label{Lions72b}
Лионс Ж.-Л. Некоторые методы решения нелинейных задач.~--- Москва: Мир, 1972.~---
588~с.

\item\label{LionsMagenes71}
Лионс Ж.-Л., Мадженес Э. Неоднородные граничные задачи и их приложения.~--- Москва:
Мир, 1971.~--- 372~с.

\item\label{Lopatinsky53}
Лопатинский Я. Б. Об одном способе сведения краевых задач для систем
дифференциальных уравнений эллиптического типа к регулярным интегральным
уравнениям~// Укр. матем. журн.~--- 1953.~--- \textbf{5}, №~2.~--- С. 123--151.

\item\label{Lopatinsky84}
Лопатинский Я. Б. Теория общих граничных задач. Избранные труды.~--- Киев: Наукова
думка, 1984.~--- 316~с.

\item\label{Magenes66}
Мадженес Э. Интерполяционные пространства и уравнения в частных производных~//
Успехи матем. наук.~--- 1966.~--- \textbf{21}, №~2.~--- С. 169--218.

\item\label{MazyaShaposhnikova86}
Мазья В. Г., Шапошникова Т. О. Мультипликаторы в пространствах дифференцируемых
функций.~--- Ленинград: Изд. Ленинградского ун-та., 1986.~--- 404~с.

\item\label{Mikhailets82}
Михайлец В. А. Асимптотика спектра эллиптических операторов и краевые условия~//
Доклады АН СССР.~--- 1982.~--- \textbf{266}, №~5.~--- С. 1059--1062.

\item\label{Mikhailets89}
Михайлец В. А. Точная оценка остаточного члена в спектральной асимптотике общих
эллиптических краевых задач~// Функц. анализ и его прил.~--- 1989.~--- \textbf{23},
№~2.~--- С. 65--66.

\item\label{Mikhailets90}
Михайлец В. А. К теории общих граничных задач для эллиптических уравнений~//
Нелинейные граничные задачи. Т. 2.~--- Киев: Наукова думка, 1990.~--- С. 66--70.

\item\label{05UMJ5}
Михайлец В. А., Мурач А.~А. Эллиптические операторы в уточненной шкале
функциональных пространств~// Укр. матем. журн.~--- 2005.~--- \textbf{57}, №~5~---
С. 689--696.

\item\label{06UMJ2}
Михайлец В. А., Мурач А.~А. Уточненные шкалы пространств и эллиптические краевые
задачи. I~// Укр. матем. журн.~--- 2006.~--- \textbf{58}, №~2~--- С. 217--235.

\item\label{06UMJ3}
Михайлец В. А., Мурач А.~А. Уточненные шкалы пространств и эллиптические краевые
задачи. II~// Укр. матем. журн.~--- 2006.~--- \textbf{58}, №~3~--- С. 352--370.

\item\label{06Dop6}
Михайлец В. А., Мурач А.~А. Интерполяция пространств с функциональным параметром и
пространства дифференцируемых функций~// Доповіді НАН України. Матем. Природозн.
Технічні науки.~--- 2006.~--- №~6.~--- С. 13--18.

\item\label{06Dop10}
Михайлец В. А., Мурач А.~А. Эллиптический оператор в уточненной шкале пространств на
замкнутом многообразии~// Доповіді НАН України. Матем. Природозн. Технічні
науки.~--- 2006.~--- №~10.~--- С. 27--33.

\item\label{06UMJ11}
Михайлец В. А., Мурач А.~А. Регулярная эллиптическая граничная задача для
однородного уравнения в двусторонней уточненной шкале пространств~// Укр. матем.
журн.~--- 2006.~--- \textbf{58}, №~11.~--- С. 1536--1555.

\item\label{06UMB4}
Михайлец В. А., Мурач А.~А. Эллиптический оператор с однородными регулярными
граничными условиями в двусторонней уточненной шкале пространств~// Укр. матем.
вісник.~--- 2006.~--- \textbf{3}, №~4.~--- С. 547--580.

\item\label{07UMJ5}
Михайлец В. А., Мурач А.~А. Уточненные шкалы пространств и эллиптические краевые
задачи.~III~// Укр. матем. журн.~--- 2007.~--- \textbf{59}, №~5.~--- С. 679--701.

\item\label{08UMJ4}
Михайлец В. А., Мурач А.~А. Эллиптическая краевая задача в двусторонней уточненной
шкале пространств~// Укр. матем. журн.~--- 2008.~--- \textbf{60}, №~4.~--- С.
497--520.

\item\label{08Collection1}
Михайлец В. А., Мурач А.~А. Интерполяционные пространства Хермандера и эллиптические
операторы~// Збірник праць Ін-ту математики НАН України.~--- \textbf{5}, №~1.~---
С.~205--226.

\item\label{09Dop3}
Михайлец В. А., Мурач А.~А. Об эллиптических операторах на замкнутом многообразии~//
Доповіді НАН України. Матем. Природозн. Технічні науки.~--- 2009.~--- №~3.~--- С.
13--19.

\item\label{Mihlin68}
Михлин С. Г. Курс математической физики~--- Москва: Наука, 1968.~--- 576~с.

\item\label{94Dop12}
Мурач А. А. Эллиптические краевые задачи в полных шкалах пространств типа
Лизоркина--Трибеля~// Доклады АН Украины.~--- 1994.~--- №~12.~--- С. 36--39.

\item\label{94UMJ12}
Мурач А. А. Эллиптические краевые задачи в полных шкалах пространств типа
Никольского~// Укр. матем. журн.~--- 1994.~--- \textbf{46}.~--- №~12.~--- С.
1647--1654.

\item\label{07UMJ6}
Мурач А. А. Эллиптические псевдодифференциальные операторы в уточненной шкале
пространств на замкнутом многообразии~// Укр. матем. журн.~--- 2007.~---
\textbf{59}, №~6.~--- С. 798--814.

\item\label{07Dop4}
Мурач А. А. Эллиптические краевые задачи в многосвязных областях в уточненной шкале
пространств~// Доповіді НАН України. Матем. Природозн. Технічні науки.~--- 2007.~---
№~4.~--- С. 29--35.

\item\label{07Dop5}
Мурач А. А. Эллиптические по Петровскому системы дифференциальных уравнений в
уточненной шкале пространств на замкнутом многообразии~// Доповіді НАН України.
Матем. Природозн. Технічні науки.~--- 2007.~--- №~5.~--- С. 29--35.

\item\label{08UMB3}
Мурач А. А. Эллиптические по Дуглису-Ниренбергу системы в пространствах обобщенной
гладкости~// Укр. матем. вісник.~--- 2008.~--- \textbf{5}, №~3.~--- С. 350--365.

\item\label{09Collection1}
Мурач А. А. Эллиптические системы в двусторонней шкале пространств Хермандера~//
Збірник праць Ін-ту математики НАН України.~--- \textbf{6}, №~1.~--- С. 126--141.

\item\label{09UMJ3}
Мурач А. А. Об эллиптических системах в пространствах Хермандера~// Укр. матем.
журн.~--- 2009.~--- \textbf{61}, №~3.~--- С. 391--399.

\item\label{07Dop6}
Мурач О. О. Крайова задача для еліптичної за Петровським системи диференціальних
рівнянь в уточненій шкалі просторів~// Доповіді НАН України. Матем. Природозн.
Технічні науки.~--- 2007.~--- №~6.~--- С. 24--31.

\item\label{Nikolsky69}
Никольский С. М. Приближение функций многих переменных и теоремы вложения.~---
Москва: Наука, 1969.~--- 480~с.

\item\label{Panich66}
Панич О. И. Эллиптические краевые задачи с параметром только в краевых условиях~//
Доклады АН СССР.~--- 1966.~--- \textbf{170}, №~5.~--- С. 1020--1023.

\item\label{Panich73}
Панич О. И. О краевых задачах с параметром в краевых условиях~// Математическая
физика, вып. 14.~--- Киев: Наукова думка.--- С. 140--146.

\item\label{Panich86}
Панич О. И. Введение в общую теорию эллиптических краевых задач.~--- Киев: Выща
школа, 1986.~--- 128~с.

\item\label{Petrovsky39}
Петровский И. Г. Об аналитичности решений систем уравнений с частными
производными~// Матем. сборник.~--- 1939.~--- \textbf{5}, №~1.~--- С. 3--70.

\item\label{Petrovsky61}
Петровский И. Г. Лекции об уравнениях с частными производными.~--- Москва:
Физматгиз, 1961.~--- 400~с.

\item\label{Petrovsky86}
Петровский И. Г. Системы уравнений с частными производными. Алгебраическая
геометрия. Избр. тр.~--- Москва: Наука, 1986.~--- 499~с.

\item\label{PsDO67}
Псевдодифференциальные операторы. Сборник статей.~--- Москва: Мир, 1967.~--- 368~с.

\item\label{Pustylnik82}
Пустыльник Е. И. О перестановочно-интерполяционных гильбертовых пространствах~//
Известия вузов. Матем.~--- 1982.~--- №~5 (240).~--- С. 43--46.

\item\label{RempelSchulze86}
Ремпель Ш., Шульце Б.-В. Теория индекса эллиптических краевых задач.~--- Москва:
Мир, 1986.~--- 575~с.

\item\label{Roitberg63}
Ройтберг Я. А. Локальное повышение гладкости вплоть до границы решений эллиптических
уравнений~// Укр. матем. журн.~--- 1963.~--- \textbf{15}, №~4.~--- С. 444--448.

\item\label{Roitberg64}
Ройтберг Я. А. Эллиптические задачи с неоднородными граничными условиями и локальное
повышение гладкости вплоть до границы обобщенных решений~// Доклады АН СССР.~---
1964.~--- \textbf{157}, №~4.~--- С. 798--801.

\item\label{Roitberg65}
Ройтберг Я. А. Теорема о гомеоморфизмах, осуществляемых в $L_{p}$ эллиптическими
операторами, и  локальное повышение гладкости обобщенных решений~// Укр. матем.
журн.~--- 1965.~--- \textbf{17}, №~5.~--- С. 122--129.

\item\label{Roitberg68}
Ройтберг Я. А. Теоремы о гомеоморфизмах, осуществляемых эллиптическими
операторами~// Доклады АН СССР.~--- 1968.~--- \textbf{180}, №~3.~--- С. 542--545.

\item\label{Roitberg69}
Ройтберг Я. А. Формула Грина и теорема о гомеоморфизмах для общих эллиптических
граничных задач с граничными условиями, не являющимися нормальными~// Укр. матем.
журн.~--- 1969.~--- \textbf{21}, №~3.~--- С. 406--413.

\item\label{Roitberg70}
Ройтберг Я. А. Теоремы о гомеоморфизмах и формула Грина для общих эллиптических
граничных задач с граничными условиями, не являющимися нормальными~// Матем.
сборник.~--- 1970.~--- \textbf{83}, №~2.~--- С. 181--213.

\item\label{Roitberg71}
Ройтберг Я. А. О значениях на границе области обобщенных решений эллиптических
уравнений~// Матем. сборник.~--- 1971.~--- \textbf{86}, №~2.~--- С. 248--267.

\item\label{Roitberg75}
Ройтберг Я. А. Теорема о полном наборе изоморфизмов для эллиптических по
Дуглису-Ниренбергу систем~// Укр. матем. журн.~--- 1975.~--- \textbf{27}, №~4.~---
С. 544--548.

\item\label{RoitbergSheftel69}
Ройтберг Я. А., Шефтель З. Г. Теоремы о гомеоморфизмах для эллиптических систем и ее
приложения.~// Матем. сборник.~--- 1969.~--- \textbf{78}, №~3.~--- С. 446--472.

\item\label{RoitbergSheftel72}
Ройтберг Я. А., Шефтель З. Г. Теорема о гомеоморфизмах для эллиптических систем с
граничными условиями, не являющимися нормальными.~// Матем. исследования,
Кишинев.~--- 1972.~--- \textbf{7}, №~2.~--- С. 143--157.

\item\label{Saks77}
Сакс Р. С. Краевые задачи для слабо эллиптических систем дифференциальных
уравнений~// Докл. АН СССР.~--- 1977.~--- \textbf{236}, №~6.~--- С. 1312--1316.

\item\label{Saks78}
Сакс Р. С. Слабоэллиптические системы псевдодифференциальных уравнений на
многообразии без края~// Труды семинара С.~Л.~Соболева. Дифференциальные уравнения с
частными производными.~--- Новосибирск: Ин-т математики СО АН СССР, 1978.~---
№~2.~--- С.~103--126.

\item\label{Saks79}
Сакс Р. С. Слабоэллиптические системы дифференциальных уравнений и их свойства~//
Труды семинара С.~Л.~Соболева. Теория кубатурных формул и приложения функционального
анализа к задачам математической физики.~--- Новосибирск: Ин-т математики СО АН
СССР, 1979.~--- №~1.~--- С. 91--118.

\item\label{Saks80}
Сакс Р. С. Слабоэллиптические краевые задачи~// Труды семинара С.~Л.~Соболева.
Дифференциальные уравнения в частных производных.~--- Новосибирск: Ин-т математики
СО АН СССР, 1980.~--- №~2.~--- С. 57--78.

\item\label{Saks81}
Сакс Р. С. Краевые задачи для обобщенных эллиптических систем дифференциальных
уравнений~// Труды семинара С.~Л.~Соболева. Дифференциальные уравнения в частных
производных.~--- Новосибирск: Ин-т математики СО АН СССР, 1981.~--- №~2.~--- С.
86--108.

\item\label{Seneta85}
Сенета Е. Правильно меняющиеся функции.~--- Москва: Наука, 1985.~--- 142~с.

\item\label{Skrypnik73}
Скрыпник И. В. Нелинейные эллиптические уравнения высшего порядка.~--- Киев: Наукова
думка, 1973.~--- 220~с.

\item\label{Slobodetsky58a}
Слободецкий Л. Н. Оценка решений эллиптических и параболических систем~// Доклады АН
СССР.~--- 1958.~--- \textbf{120}, №~3.~--- С. 468--471.

\item\label{Slobodetsky58b}
Слободецкий Л. Н. Оценки в $L_{p}$ решений эллиптических систем~// Доклады АН
СССР.~--- 1958.~--- \textbf{123}, №~4.~--- С. 616--619.

\item\label{Slobodetsky58c}
Слободецкий Л. Н. Обобщенные пространства С.~Л.~Соболева и их приложения к краевым
задачам для дифференциальных уравнений в частных производных~// Ученые записки
Ленинградского гос. пед. ин-та.~--- 1958.~--- \textbf{187}.~--- С. 54--112.

\item\label{Slobodetsky60}
Слободецкий Л. Н. Оценки в $L_{2}$ решений линейных эллиптических и параболических
систем, I~// Вестник ЛГУ, сер. матем., мех. и астр.~--- 1960.~--- №~7.~--- С.
28--47.

\item\label{Sobolev50}
Соболев С. Л. Некоторые применения функционального анализа в математической
физике.~--- Ленинград: Изд. Ленинград. ун-та, 1950.~--- 255~с.

\item\label{Sobolev74}
Соболев С. Л. Введение в теорию кубатурных формул.~--- Москва: Наука, 1974.~---
808~с.

\item\label{Solonnikov63}
Солонников В. А. Оценки решений общих краевых задач для эллиптических систем~//
Доклады АН СССР.~--- 1963.~--- \textbf{151}, №~4.~--- С. 783--785.

\item\label{Solonnikov64}
Солонников В. А. Об общих краевых задачах, эллиптических в смысле
А.~Дуглиса-Л.~Ниренберга. I~// Известия АН СССР, сер. матем.~--- 1964.~---
\textbf{28}, №~3.~--- С. 665--706.

\item\label{Solonnikov66}
Солонников В. А. Об общих краевых задачах, эллиптических в смысле
А.~Дуглиса-Л.~Ниренберга. II~// Труды Матем. ин-та АН СССР.~--- 1966.~---
\textbf{92}.~--- С. 233--297.

\item\label{Solonnikov67}
Солонников В. А. Об оценках в $L_{p}$ решений эллиптических и параболических
систем~// Труды Матем. ин-та АН СССР.~--- 1967.~--- \textbf{102}.~--- С. 137--160.

\item\label{Stepanets87}
Степанец А. И. Классификация и приближение периодических функций.~--- Киев: Наукова
думка, 1987.~--- 268~с.

\item\label{Stepanets02}
Степанец А. И. Методы теории приближений. В 2-х томах.~--- Киев: Институт математики
НАН Украины, 2002.~--- 468~с., 427~с.

\item\label{Taylor85}
Тейлор М. Псевдодифференциальные операторы.~--- Москва: Мир, 1985.~--- 472~с.

\item\label{Treves84a}
Трев Ф. Введение в теорию псевдодифференциальных операторов и интегральных
операторов Фурье: В 2-х томах. Т.~1. Псевдодифференциальные операторы.~--- Москва:
Мир, 1984.~--- 360~с.

\item\label{Treves84b}
Трев Ф. Введение в теорию псевдодифференциальных операторов и интегральных
операторов Фурье: В 2-х томах. Т.~2. Интегральные операторы Фурье.~--- Москва: Мир,
1984.~--- 400~с.

\item\label{Triebel80}
Трибель Х. Теория интерполяции, функциональные пространства, дифференциальные
операторы.~--- Москва: Мир, 1980.~--- 664~с.

\item\label{Triebel86}
Трибель Х. Теория функциональных пространств. --- Москва: Мир, 1986.~--- 447~с.

\item\label{Ulyanov63}
Ульянов П. Л. О множителях Вейля для безусловной сходимости~// Матем. сборник.~---
1963.~--- \textbf{60}, №~1.~--- С. 39--62.

\item\label{Ulyanov64}
Ульянов П. Л. Решенные и нерешенные проблемы теории тригонометрических и
ортогональных рядов~// Успехи матем. наук.~--- 1964.~--- \textbf{19}, №~1.~--- С.
3--69.

\item\label{FunctionalAnalysis72}
Функциональный анализ / Под общ. ред. С.~Г.~Крейна.~--- Москва: Наука, 1972.~---
544~с.

\item\label{Hermander59}
Хермандер Л. К теории общих дифференциальных операторов в частных производных.~---
Москва: Изд. иностр. лит., 1959.~--- 132~с.

\item\label{Hermander65}
Хермандер Л. Линейные дифференциальные операторы с частными производными.~---
Москва: Мир, 1965.~--- 380~с. (Перевод англоязычного издания: Berlin,
Springer--Verlag, 1963.)

\item\label{Hermander67b}
Хермандер Л. Псевдодифференциальные операторы~// Псевдодифференциальные
операторы.~--- Москва: Мир, 1967.~--- С. 63--87.

\item\label{Hermander67}
Хермандер Л. Псевдодифференциальные операторы и неэллиптические краевые задачи~//
Псевдодифференциальные операторы.~--- Москва: Мир, 1967.~--- С. 166--296.

\item\label{Hermander86}
Хермандер Л. Анализ линейных дифференциальных операторов с частными производными: В
4-х т. Т. 2. Дифференциальные операторы с постоянными коэффициентами.~--- Москва:
Мир, 1986.~--- 456~с.

\item\label{Hermander87}
Хермандер Л. Анализ линейных дифференциальных операторов с частными производными: В
4-х т. Т. 3. Псевдодифференциальные операторы.~--- Москва: Мир, 1987.~--- 696~с.

\item\label{Shapiro53}
Шапиро З. Я. Об общих краевых задачах для эллиптических уравнений// Известия АН
СССР, сер. матем.~--- 1953.~--- \textbf{17}.~--- С. 539--562.

\item\label{Schechter60a}
Шехтер М. Общие граничные задачи для эллиптических уравнений в частных
производных~// Математика~: сб-к переводов.~--- 1960.~--- № 4--5.~--- С. 93--122.

\item\label{Shlenzak74}
Шлензак Г. Эллиптические задачи в уточненной шкале пространств~// Вестник
Московского ун-та. Сер. 1, матем., механика.~--- 1974.~--- \textbf{29}, №~4.~--- С.
48--58.

\item\label{Shubin78}
Шубин М. А. Псевдодифференциальные операторы и спектральная теория.~--- Москва:
Наука, 1978.~--- 280~с.

\item\label{Adams75}
Adams R. A. Sobolev Spaces.~--- New York: Academic Press, 1975.~--- xiii+268~p.

\item\label{Agmon62}
Agmon S. On the eigenfunctions and on the eigenvalues of general elliptic boundary
value problems~// Commun. Pure Appl. Math.~--- 1962.~--- \textbf{15}, no.~2.~--- P.
119--147.

\item\label{AgmonDouglisNirenberg64}
Agmon S., Douglis A., Nirenberg L. Estimates near the boundary for solutions of
elliptic partial differential equations satisfying general boundary
conditions.~II~// Commun. Pure Appl. Math.~--- 1964.~--- \textbf{17}, no.~1.~--- P.
35--92.

\item\label{AgmonNirenberg63}
Agmon S., Nirenberg L. Properties of solutions of ordinary differential equations in
Banach space~// Commun. Pure Appl. Math.~--- 1963.~--- \textbf{16}, no.~2.~--- P.
121--239.

\item\label{Agranovich97}
Agranovich M. S. Elliptic boundary problems~// Partial differential equa\-tions IX.
Encyclopedia of Mathematical Sciences. Vol. 79.~--- Berlin: Springer, 1997.~--- P.
1--144.

\item\label{AronszajnMilgram52}
Aronszajn N., Milgram A.~N. Differential equations on Riemannian mani\-folds~//
Rend. Circ. Mat. Palermo.~--- 1952.~--- \textbf{2}.~--- P. 1--61.

\item\label{AtiyahBott64}
Atiyah M. F., Bott R. The index theorem for manifolds with boundary~// Proc. Symp.
on Differentional Analysis. Bombay Coll.~--- London: Oxford University Press,
1964.~--- P. 175--186.

\item\label{AtiyahSinger63}
Atiyah M. F., Singer M.~F. The index of elliptic operators on compact manifolds~//
Bull. Amer. Math. Soc.~--- 1963.~--- \textbf{69}, no.~3.~--- P. 422--433.

\item\label{Avakumovic36}
Avakumovi\'c V. G. O jednom O-inverznom stavu~// Rad Jugoslovenske Akad. Znatn.
Umjetnosti.~--- 1936.~--- \textbf{254}.~--- P. 167--186.

\item\label{BennetSharpley88}
Bennet K, Sharpley R. Interpolation of operators.~--- Boston: Academic Press,
1988.~--- xiv+469~p.

\item\label{BinghamGoldieTeugels89}
Bingham N. H., Goldie C.~M., Teugels J.~L. Regular variation.~--- Cambridge:
Cambridge Univ. Press, 1989.~--- 512~p.

\item\label{Boyd67}
Boyd D. W. The Hilbert transform on rearrangement-invariant spaces~// Canadian J.
Math.~--- 1967.~--- \textbf{19}.~--- P. 599--616.

\item\label{Browder56}
Browder F. E. On the regularity properties of solutions of elliptic differential
equations~// Commun. Pure Appl. Math.~--- 1956.~--- \textbf{9}, no.~3.~--- P.
351--361.

\item\label{Browder59}
Browder F. E. Estimates and existence theorems for elliptic boundary-value
problems~// Proc. Nat. Acad. Sci.~--- 1959.~--- \textbf{45}, no.~3.~--- P. 365--372.

\item\label{BrudnyiKrugljak91}
Brudny\"{\i} Yu. A., Krugljak N. Ya. Interpolation functors and interpolation
spaces.~--- Amsterdam: North-Holland, 1991.~--- xvi+718~p.

\item\label{Burenkov99}
Burenkov V. Extension theorems for Sobolev spaces~// Oper. Theory Adv. Appl.~---
1999.~--- \textbf{109}.~--- P. 187--200.

\item\label{Calderon64}
Calderon A. P. Intermediate spaces and interpolation, the compex method~// Studia
Math.~--- 1964.~--- \textbf{24}.~--- P. 113--190.

\item\label{CarroCerda90}
Carro M. J., Cerd\`a J. On complex interpolation with an analytic functional~//
Math. Scand.~--- 1990.~--- \textbf{66}, no.~2.~--- P. 264--274.

\item\label{CobosFernandez88}
Cobos F., Fernandez D.~L. Hardy-Sobolev spaces and Besov spaces with a function
parameter~// Lecture Notes in Math.~--- 1988.~--- \textbf{1302}.~--- P. 158--170.

\item\label{DenkMennickenVolevich98}
Denk R., Mennicken R., Volevich L. The Newton polygon and elliptic problems with
parameter~// Math. Nachr.~--- 1998.~-- \textbf{192}, no.~1.~--- P. 125--157.

\item\label{DenkMennickenVolevich01}
Denk R., Mennicken R., Volevich L. On elliptic operator pensils with general
boundary conditions~// Integral Equations Operator Theory.~--- 2001.~---
\textbf{39}, no.~1.~--- P. 15--40.

\item\label{DenkVolevich02}
Denk R., Volevich L. Parameter-elliptic boundary value problems connected with the
Newton polygon~// Differential Integral Equations.~--- 2002.~--- \textbf{15},
no.~3.~--- P. 289--326.

\item\label{DjakovMityagin09}
Djakov P., Mityagin B. Spectral gaps of Schr\"odinger operators with periodic
singular potentials~// Dyn. Partial Differ. Equ.~--- 2009.~--- \textbf{6},
no.~2.~--- P. 95--165.

\item\label{Donoghue67}
Donoghue W. F. The interpolation of quadratic norms~// Acta Math.~--- 1967.~---
\textbf{118}, no. 3--4.~--- P. 251--270.

\item\label{DouglisNirenberg55}
Douglis A., Nirenberg L. Interior estimates for elliptic systems of partial
differential equations~// Comm. Pure Appl. Math.~--- 1955.~--- \textbf{8},
no.~4.~--- P. 503--538.

\item\label{EdmundsGurkaOpic97}
Edmunds D. E., Gurka P., Opic P. On embeddings of logarithmic Bessel potential
spaces~// J. Funct. Anal.~--- 1997.~--- \textbf{146}.~--- P. 116--150.

\item\label{EdmundsHaroske99}
Edmunds D. E., Haroske D. Spaces of Lipschitz type, embeddings and entropy
numbers~// Diss. Math.~--- 1999.~--- \textbf{380}.~--- P. 1--43.



\item\label{FarkasJacobScilling01a}
Farkas W., Jacob N., Schilling R.~L. Feller semigroups, $L_{p}$-sub-Markovian
semigroups, and app\-li\-ca\-tions to pseudo-differential operators with negative
definite symbols~// Forum Math.~--- 2001.~--- \textbf{13}.~--- P. 59--90.

\item\label{FarkasJacobScilling01b}
Farkas W., Jacob N., Schilling R.~L. Function spaces related to continuous negative
definite functions: $\psi$-Bessel potential spaces~// Diss. Math.~--- 2001.~---
\textbf{393}.~--- P. 1--62.

\item\label{FarkasLeopold06}
Farkas W., Leopold H.-G. Characterisations of function spaces of genera\-lized
smoothness~// Ann. Mat. Pura Appl.~--- 2006.~--- \textbf{185}, no.~1.~--- P. 1--62.

\item\label{FoiasLions61}
Foia\c{s} C., Lions J.-L. Sur certains th\'eor\`emes d'interpolation~// Acta Scient.
Math. Szeged.~--- 1961.~--- \textbf{22}, no. 3--4.~--- P. 269--282.

\item\label{Franke85}
Franke J. Besov--Triebel--Lizorkin spaces and boundary value problems~// Seminar
analysis 1984/85.~--- Berlin: Akad. Wiss. DDR, 1985.~--- P. 89--104.

\item\label{Gagliardo68}
Gagliardo E. Caratterizzazione construttiva di tutti gli spazi di interpolazione tra
spazi di Banach~// Symposia Math. Vol. II, INDAM, Rome, 1968.~--- London: Academic
Press, 1969.~--- P. 95--106.

\item\label{GelukHaan87}
Geluk J. L., Haan L.~de. Regular variation, extensions and Tauberian theorems.~---
Amsterdam: Stichting Mathematisch Centrum, 1987.~--- iv+132~p.

\item\label{Geymonat65}
Geymonat G. Sui problemi ai limiti per i sistemi lineari ellittici~// Annali di
Matematica Pura ed Applicata. (4).~--- 1965.~--- \textbf{69}.~--- P. 207--284.

\item\label{Grisvard67}
Grisvard P. Caract\'erisation de quelques espaces d'inter\-po\-la\-tion~// Arch.
Rat. Mech. Anal.~--- 1967.~--- \textbf{25}, no.~1.~--- P. 40--63.

\item\label{Grubb68}
Grubb G. A characterization of non-local boundary value problems associ\-ated with
an elliptic operator~// Ann. Scuola Norm. Sup. Pisa (3).~--- 1968.~--- \textbf{22},
no.~3.~--- P. 425--513.

\item\label{Grubb71}
Grubb G. On coerciveness and semiboundedness of general boundary problems~/
G.~Grubb~// Isr. J. Math.~--- 1971.~--- \textbf{10}.~--- P. 32--95.

\item\label{Grubb96}
Grubb G. Functional Calculas of Pseudo-Differential Boundary Problems.~--- 2-nd
ed.~--- Boston: Birkh\"aser, 1996.~--- 522~p.

\item\label{Gustavsson78}
Gustavsson J. A function parameter in connection with interpolation of Banach
spaces~// Math. Scand.~--- 1978.~--- \textbf{42}, no.~2.~--- P. 289--305.

\item\label{Haan70}
Haan L.~de. On Regular variation and its application to the weak convergence of
sample extremes.~--- Amsterdam: Mathematisch Centrum, 1970.~--- v+124~p.

\item\label{HaroskeMoura04}
Haroske D. D., Moura S.~D. Continuity envelopes of spaces of generalised smoothness,
entropy and approximation numbers~// J. Approximation Theory.~--- 2004.~---
\textbf{128}.~--- P. 151--174.

\item\label{Hegland95}
Hegland M. Variable Hilbert scales and their interpolation inequalities with
applications to Tikhonov regularization~// Appl. Anal.~--- 1995.~--- \textbf{59},
no.~1--4.~--- P. 207--223.

\item\label{Hegland10}
Hegland M. Error bounds for spectral enhancement which are based on variable Hilbert
scale inequalities~// J. Integral Equations Appl.~--- 2010.~--- \textbf{22},
no.~2.~--- P. 285--312.

\item\label{Jacob010205}
Jacob N. Pseudodifferential operators and Markov processes (in 3 volumes).~--
London: Imperial College Press, 2001, 2002, 2005.~--- xxii+493~p., xxii+453~p.,
xxviii+474~p.

\item\label{Janson81}
Janson S. Minimal and maximal methods of interpolation~// J. Funct. Anal.~---
1981.~--- \textbf{44}, no.~1.~--- P. 50--73.

\item\label{KalyabinLizorkin87}
Kalyabin G. A., Lizorkin P.~I. Spaces of functions of generalized smoothness~//
Math. Nachr.~--- 1987.~--- \textbf{133}, no.~1.~--- P. 7--32.

\item\label{Karamata30a}
Karamata J. Sur certains "Tauberian theorems" \;de M.~M.~Hardy et Litt\-lewood~//
Mathematica (Cluj).~--- 1930.~--- \textbf{3}.~--- P. 33--48.

\item\label{Karamata30b}
Karamata J. Sur un mode de croissance r\'eguli\`ere des fonctions~//
Ma\-the\-ma\-ti\-ca (Cluj).~--- 1930.~--- \textbf{4}.~--- P. 38--53.

\item\label{Karamata33}
Karamata J. Sur un mode de croissance r\'eguli\`ere. Th\'eor\`ems
foun\-d\-amen\-taux~// Bull. Soc. Math. France.~--- 1933.~--- \textbf{61}.~--- P.
55--62.

\item\label{Kozhevnikov96}
Kozhevnikov A. Asymptotics of the spectrum of the Douglis-Nirenberg elliptic
operators on a compact manifold~// Math. Nachr.~--- 1996.~--- \textbf{182},
no.~1.~--- P. 261--293.

\item\label{Kozhevnikov97}
Kozhevnikov A. Parameter-ellipticity for mixed-order elliptic boundary problems~//
C.~R. Math. Acad. Sci. Paris.~--- 1997.~--- \textbf{324}, no.~12.~--- P. 1361--1366.

\item\label{Kozhevnikov01}
Kozhevnikov A. Complete scale of isomorphisms for elliptic pseudodifferen\-tial
boundary-value problems~// J. London Math. Soc. (2).~--- 2001.~--- \textbf{64},
no~2.~--- P.~409--422.

\item\label{KozlovMazyaRossmann97}
Kozlov V. A., Maz'ya V. G., Rossmann J. Elliptic boundary value problems in domains
with point singularities.~--- Providence: American Math. Soc., 1997.~--- ix+414~p.

\item\label{Krugljak93}
Krugljak N. Ja. On the reiteration property of $\overrightarrow{X}_{\varphi,q}$
spaces~// Math. Scand.~--- 1993.~--- \textbf{73}, no.~1.~--- P. 65--80.

\item\label{Leopold98}
Leopold H.-G. Embeddings and entropy numbers in Besov spaces of generalized
smoothness~// Lect. Notes Pure Appl. Math.~--- 2000.~--- \textbf{213}.~--- P.
323--336.

\item\label{Lewy57}
Lewy H. An example of a smooth linear partial differential equation without
solution~// Ann. of Math. (2).~--- 1957.~--- \textbf{66}.~--- P. 155--158.

\item\label{Lions58}
Lions J.-L. Espaces interm\'ediaires entre espaces hilbertiens et applications~//
Bull. Math. Soc. Sci. Math. Phys. R. P. Roumanie.~--- 1958.~--- \textbf{50},
no.~4.~--- P. 419--432.

\item\label{Lions59}
Lions J.-L. Th\'eor\`emes de trace et d'interpolation.~I~// Ann. Scuola Norm. Sup.
Pisa (3).~--- 1959.~--- \textbf{13}.~--- P. 389--403.

\item\label{Lions60}
Lions J.-L. Une construction d'espaces d'interpolation~// C.~R. Math. Acad. Sci.
Paris.~--- 1961.~--- \textbf{251}.~--- P. 1853--1855.

\item\label{LionsMagenes62}
Lions J.-L., Magenes E. Probl\'emes aux limites  non homog\'enes, V~// Ann. Scuola
Norm. Sup. Pisa (3).~--- 1962.~--- \textbf{16}.~--- P. 1--44.

\item\label{LionsMagenes63}
Lions J.-L., Magenes E. Probl\'emes aux limites non homog\'enes, VI~// J. Anal.
Math.~--- 1963.~--- \textbf{11}.~--- P. 165--188.

\item\label{LionsPeetre61}
Lions J.-L., Peetre J. Propri\'et\'es d'espaces d'interpolation~// C.~R. Math. Acad.
Sci. Paris.~--- 1961.~--- \textbf{253}.~--- P. 1747--1749.

\item\label{Lefstrem92}
L\"ofstr\"om J. Interpolation of boundary value problems of Neumann type on smooth
domains~// J. London Math. Soc (2).~--- 1992.~--- \textbf{46}, no.~3.~--- P.
499--516.

\item\label{Malgrange57}
Malgrange B. Sur une classe d'op\'eratuers diff\'erentiels hypoelliptiques~// Bull.
Soc. Math. France.~--- 1957.~--- \textbf{85}.~--- P. 283--306.

\item\label{Maric00}
Maric V. Regular variation and differential equations.~--- New York: Springer
Verlag, 2000.~--- 127~p.

\item\label{MatheTautenhahn06}
Math\'e P., Tautenhahn U. Interpolation in variable Hilbert scales with application
to innverse problems~// Inverse Problems.~--- 2006.~--- \textbf{22}, no.~6.~--- P.
2271--2297.

\item\label{Meaney82}
Meaney C. On almost-everywhere convergent eigenfunction expansions of the
Laplace-Beltrami operator~// Math. Proc. Camb. Phil. Soc.~--- 1982.~---
\textbf{92}.~--- P. 129--131.

\item\label{Menschoff23}
Menschoff D. Sur les series de fonctions orthogonales I~// Fund. Math.~--- 1923.~---
\textbf{4}.~--- P. 82--105.

\item\label{Merucci82}
Merucci C. Interpolation r\'eelle avec fonction param\`etre: r\'eit\'eration et
applications aux espaces $\Lambda^{\varrho}(\varphi)$~// C.~R. Math. Acad. Sci.
Paris.~--- 1982.~--- \textbf{295}, no.~6.~--- P. 427--430.

\item\label{Merucci84}
Merucci C. Application of interpolation with a function parameter to Lorentz,
Sobolev and Besov spaces~// Lecture Notes in Math.~--- 1984.~--- \textbf{1070}.~---
P. 183--201.

\item\label{MikhailetsMolyboga09a}
Mikhailets V., Molyboga V. Spectral gaps of the one-dimensional Schr\"odinger
operators with singular periodic potentials~// Methods Funct. Anal. Topology.~---
2009.~--- \textbf{15}, no.~1.~--- P. 31--40.

\item\label{MikhailetsMolyboga10}
Mikhailets V., Molyboga V. Smoothness of Hill's potential and lengths of spectral
gaps~// arXiv: math.SP/1003.5000.~--- March 2010.~--- 8~p.

\item\label{08MFAT1}
Mikhailets V. A., Murach A.~A. Interpolation with a function parameter and refined
scale of spaces~// Methods Funct. Anal. Topology.~--- 2008.~--- \textbf{14},
no.~1.~--- P. 81--100.

\item\label{08BPAS3}
Mikhailets V. A., Murach A.~A. Elliptic systems of pseudodifferential equations in a
refined scale on a closed mani\-fold~// Bull. Pol. Acad. Sci. Math.~--- 2008.~---
\textbf{56}, no. 3--4.~--- P. 213--224.

\item\label{09OperatorTheory191}
Mikhailets V. A., Murach A.~A. Elliptic problems and H\"ormander spaces~// Oper.
Theory Adv. Appl.~--- 2009.~--- \textbf{191}.~--- P. 447--470.

\item\label{Moura01}
Moura S. Function spaces of generalised smoothness~// Diss. Math.~--- 2001.~---
\textbf{398}.~--- P. 1--87.

\item\label{08MFAT2}
Murach A. A. Douglis-Nirenberg elliptic systems in the refined scale of spaces on a
closed manifold~// Methods Funct. Anal. Topo\-lo\-gy.~--- 2008.~--- \textbf{14},
no.~2.~--- P. 142--158.

\item\label{09MFAT2}
Murach A. A. Extension of some Lions-Magenes theorems~// Methods Funct. Anal.
Topology.~--- 2009.~--- \textbf{15}, no.~2.~--- P. 152--167.

\item\label{OpicTrebels00}
Opic B., Trebels W. Bessel potentials with logarithmic components and Sobolev-type
embeddings~// Anal. Math.~--- 2000.~--- \textbf{26}.~--- P. 299--319.

\item\label{Orlicz27}
Orlicz W. Zur Theorie der Orthogonalreihen~// Bull. Acad. Sci. Polon. Cracovie.~---
1927.~--- P. 81--115.

\item\label{Ovchinnikov84}
Ovchinnikov V. I. The methods of orbits in interpolation theory~// Mathematical
Reports.~--- 1984.~--- \textbf{1}, no.~2.~--- P. 349--515.

\item\label{Paneah00}
Paneah B. The oblique derivative problem. The Poincar\'e problem.--- Berlin:
Wiley--VCH, 2000.~--- 348~p.

\item\label{Peetre61}
Peetre J. Mixed problems for higher order elliptic equations in two variables.~I~//
Ann. Scuola Norm. Sup. Pisa (3).~--- 1961.~--- \textbf{15}.~--- P. 337--353.

\item\label{Peetre63}
Peetre J. Sur le nombre de param\`etres dans la d\'efinition de certain espaces
d'interpolation~// Ric. Mat.~--- 1963.~--- \textbf{15}.~--- P. 248--261.

\item\label{Peetre68}
Peetre J. On interpolation functions. II~/ J.~Peetre~// Acta Sci. Math.~---
1968.~--- \textbf{29}, no.~1.~--- P. 91--92.

\item\label{Persson86}
Persson L.-E. Interpolation with a function parameter~// Math. Scand.~--- 1986.~---
\textbf{59}, no.~2.~--- P. 199--222.

\item\label{Plis54}
Pli\'s A. The problem of uniqueness for the solution of a system of partial
differential equations~// Bull. Pol. Acad. Sci. Math.~--- 1954.~-- \textbf{2}.~---
P. 55--57.

\item\label{Poschel08}
P\"oschel J. Hill's potensials in weighted Sobolev spaces and their spectral gaps~//
Hamiltonian dynamical systems and applications.~--- Dordrecht: Springer--Verlag,
2008.~--- P. 421--430.

\item\label{Rademacher22}
Rademacher H. Einige S\"{a}tze \"{u}ber Reihen von allgemeinen
Orthogonal\-functionen~// Math. Annalen.~--- 1922.~--- \textbf{87}.~--- S. 111--138.

\item\label{Reshnick87}
Reshnick S. I. Extreme values, regular variation and point processes.~-- New York:
Springer Verlag, 1987.~--- 320~p.

\item\label{Roitberg96}
Roitberg Ya. A. Elliptic boundary value problems in the spaces of
distri\-bu\-tions.~--- Dordrecht: Kluwer Acad. Publishers, 1996.~--- xii+415~p.

\item\label{Roitberg99}
Roitberg Ya. A. Elliptic boundary value problems in the spaces of
distri\-bu\-tions.~--- Dordrecht: Kluwer Acad. Publishers, 1999.~--- x+276~p.

\item\label{RoitbergRoitberg00}
Roitberg I., Roitberg Ya. Green's formula and theorems on isomorphisms for general
elliptic problems for Douglis-Nirenberg elliptic systems~// Oper. Theory Adv.
Appl.~--- 2000.~--- \textbf{117}.~--- P. 281--299.

\item\label{Schechter60b}
Schechter M. Mixed boundary value problems for general elliptic equa\-ti\-ons~//
Comm. Pure Appl. Math.~--- 1960.~--- \textbf{13}, no.~2.~--- P. 183--201.

\item\label{Schechter61}
Schechter M. A local regularity theorem~// J. Math. Mech.~--- 1961.~--- \textbf{10},
no.~2.~--- P. 279--288.

\item\label{Schechter63}
Schechter M. On $L_{p}$ estimates and regularity,~I~// Amer. J. Math.~--- 1963.~---
\textbf{85}, no.~1.~--- P. 1--13.

\item\label{Schechter67}
Schechter M. Complex interpolation~// Compos. Math.~--- 1967.~--- \textbf{18}, no.
1,~2.~--- P. 117--147.

\item\label{Schechter77}
Schechter M. Modern methods in partial differential equations.~--- New York:
McGraw-Hill Inc, 1977.~--- xv+245~p.

\item\label{Schwartz50}
Schwartz L. Th\'eorie des distributions. Vol.~I.~--- Paris: Herman, 1950.~--- 148~p.

\item\label{Schwartz51}
Schwartz L. Th\'eorie des distributions. Vol.~II.~--- Paris: Herman, 1951.~---
169~p.

\item\label{Seeley66}
Seeley R. T. Singular integrals and boundary value problems~// Amer. J. Math.~---
1966.~--- \textbf{88}, no.~4.~--- P. 781--809.

\item\label{Seeley72}
Seeley R. Interpolation in $L_{p}$ with boundary conditions~// Studia Math.~---
1972.~--- \textbf{44}.~--- P. 47--60.

\item\label{Simanca87}
Simanca S. R. Mixed elliptic boundary value problems~// Comm. Partial Differential
Equations.~--- 1987.~--- \textbf{12}, no.~2.~--- P. 123--200.

\item\label{Strichartz67}
Strichartz R. S. Multipliers on fractional Sobolev spaces~// J. Math. Mech.~---
1967.~--- \textbf{16}, no.~9.~--- P. 1031--1060.

\item\label{Tandori62}
Tandori K. \"{U}ber die orthogonalen Functionen X (unbedingte Kover\-genz)~// Acta
Scient. Math.~--- 1962.~--- \textbf{23}, no. 3--4.~--- P. 185--221.

\item\label{Tartar07}
Tartar L. An introduction to Sobolev spaces and interpolation spaces.~--- Berlin:
Springer-Verlag, 2007.~--- xxv+219~p.

\item\label{Triebel92}
Triebel H. Theory of function spaces. II.~--- Basel: Birkh\"aser, 1992.~---
viii+370~p.

\item\label{Triebel01}
Triebel H. The structure of functions.~--- Basel: Birkh\"aser, 2001.~--- xii+425~p.

\item\label{Triebel06}
Triebel H. Theory of function spaces. III.~--- Basel: Birkh\"aser, 2006.~---
xii+426~p.

\item\label{WlokaRowleyLawruk95}
Wloka J. T., Rowley B., Lawruk B. Boundary value problems for elliptic systems.~---
Cambridge: Cambridge University Press, 1995.~--- xiv+641~p.

\end{enumerate}

\normalsize




\newpage

\thispagestyle{empty}

\Large

\vspace*{2cm}

\begin{center}

\textbf{Михайлец Владимир Андреевич \\ Мурач Александр Александрович}

\vspace{2cm}

\LARGE

\textbf{ПРОСТРАНСТВА ХЕРМАНДЕРА, \\ ИНТЕРПОЛЯЦИЯ И \\ ЭЛЛИПТИЧЕСКИЕ ЗАДАЧИ}

\end{center}

\normalsize

\newpage

\thispagestyle{empty}

\Large

\vspace*{2cm}

\begin{center}

\textbf{Mikhailets Vladimir Andreevich \\ Murach Aleksandr Aleksandrovich}

\vspace{2cm}

\LARGE

\textbf{H\"ORMANDER SPACES, \\ INTERPOLATION, \\
AND ELLIPTIC PROBLEMS}

\vspace{0.5cm}

\Large

(\textit{in Russian})

\end{center}

\normalsize


\end{document}